\documentclass[graybox,envcountchap,sectrefs,]{svmult}
\usepackage[a4paper,total={12cm,20cm},centering,includehead]{geometry}
\usepackage{latexsym,amssymb}
\usepackage{mathpazo}
\usepackage{fancyhdr}
\renewcommand{\chaptermark}[1]{\markboth{\thechapter.\ #1}{}}

\usepackage[Lenny]{fncychap}
\usepackage{mathptmx}
\usepackage{courier}
\usepackage{type1cm}         
\usepackage{makeidx}         
 \usepackage{graphicx}
\usepackage{multicol}        
\usepackage[bottom]{footmisc}

  \newcommand{\ds}{\displaystyle}

\makeindex             
  \newcommand{\be}{\begin{equation}}
\newcommand{\ee}{\end{equation}}

\newcommand{\ba}{\begin{eqnarray}}
\newcommand{\ea}{\end{eqnarray}}
\newcommand{\bal}{\begin{align}}
\newcommand{\eal}{\end{align}}
\usepackage[francais]{babel}
\usepackage[T1]{fontenc}
\usepackage[latin1]{inputenc}
  \newcommand{\un}{\underline}
	\newcommand{\ov}{\overline} 
	\def\BbR{{\mathbb R}}

	\usepackage{amsfonts}				
\newcommand{\lb}{\label}
\def\BbC{\mathbb C}
\def\BbS{\mathbb S}
\def\BbT{\mathbb T}

\def\m{\mathcal}
\newtheorem{Theorem}{Théorème}[chapter]
\newtheorem{Definition}{Définition}[chapter]
\newtheorem{Lemma}{Lemme}[chapter]
\newtheorem{Proposition}{Proposition}[chapter]
\newtheorem{Problem}{Problème}[chapter]
\newtheorem{Exercice}{Exercice}[chapter]
\newcommand{\bthm}{\begin{Theorem}}
\newcommand{\ethm}{\end{Theorem}}
\newcommand{\bpr}{\begin{Proposition}}
\newcommand{\epr}{\end{Proposition}}
\newcommand{\blm}{\begin{Lemma}}
\newcommand{\elm}{\end{Lemma}}
\newcommand{\bdf}{\begin{Definition}}
\newcommand{\edf}{\end{Definition}}
\newcommand{\bpb}{\begin{Problem}}
\newcommand{\epb}{\end{Problem}}
\newcommand{\bex}{\begin{Exercice}}
\newcommand{\eex}{\end{Exercice}}
\newcommand{\fii}{\varphi}
\newcommand{\lm}{\lambda}
\begin{document}

\frontmatter

\pagenumbering{roman}
\author{
\begin{center}
\textbf{\textsc{{Abdelmajid BEN HADJ SALEM}} \\
\textsc{{Ingénieur Général Géographe}}}
\end{center}}
\vspace{1cm}
\title{
\begin{center}
\textbf{\textit{{Eléments de Géodésie et de la Théorie des Moindres Carrés} \\ \\ \\  }}
\vspace{8cm}
\Large
\textsc{\textbf{Version provisoire - 15 décembre 2016}}
\end{center}
}
\cleardoublepage
\begin{center}
\textbf{\textsc{\textit{Eléments de Géodésie et de la Théorie des Moindres Carrés}}} \\  
\end{center}
\vspace{0.5cm}
\begin{center}
\Large
\textbf{Par} 

\vspace{1cm}
\textbf{\textsc{Abdelmajid BEN HADJ SALEM} \\
\textsc{Ingénieur Général Géographe}}
\end{center}
\newpage
\normalsize
Abdelmajid BEN HADJ SALEM\\
6, rue du Nil, Cité Soliman Erriadh\\
8020 Soliman, Tunisia

e-mail: abenhadjsalem@gmail.com

\vglue10cm
\copyright\ 2016 Abdelmajid BEN HADJ SALEM
         \vfill\eject

\frontmatter
\begin{dedication}
\begin{center}
{\huge
\texttt{O my Lord! Increase me further in knowledge}.
}
\\[5ex]

{\large
(Holy Quran, Surah Ta Ha, 20:114.)
}
\end{center}
\end{dedication}

\begin{dedication}
\Large
\vspace{1cm}
$\m A$ \texttt{mes chers parents, à ma femme, à mes enfants, à mes professeurs et à tous ceux qui m'ont apporté leur soutien.}
\vspace{1.5cm}

\noindent $\m A$\texttt{ux martyrs et les blessés de la Révolution Tunisienne.}
 \vspace{1cm}

\end{dedication}
\normalsize









\chapter*{\textit{\textbf{Pr\'eface}}}

\addcontentsline{toc}{section}{Pr\'eface}
     C'est un grand bonheur de remettre ce livre contenant un cours d'introduction à la géodésie destiné à la formation d'ingénieurs en sciences géographiques. Ce cours est le fruit de l'enseignement de la géodésie que j'ai pu donner depuis le début des années quatre vingt dix du dernier siècle aux étudiants des Instituts Supérieurs des Etudes Technologiques ou à ceux du Diplôme des Etudes Supérieures Spécialisées de Géomatique à l'Ecole Nationale des Ingénieurs de Tunis, aux étudiants de la Faculté des Sciences de Tunis sans oublier les ingénieurs et techniciens de l'Office de la Topographie et du Cadastre et tout récemment les élèves ingénieurs de l'option topographie et géomatique de l'Ecole Supérieure Privée d'Aéronautique et des Technologies de Tunis. 
     
  Cet ouvrage constitue ma modeste participation à enrichir la documentation nationale en matière des sciences géographiques ou sous l'appellation de nos jours la géomatique et en particulier concernant son pilier fondamental à savoir la géodésie. En plus de son aspect pédagogique, l'ouvrage collecte des informations que le géodésien et le géomètre pratiquant ont besoin et que souvent, elles sont dispersées.   
	\\
		
		Après un chapitre d'introduction, ce cours comprend deux parties:

     \textbf{- Partie I}: 
		
		On présente l'essentiel de la géodésie géométrique et spatiale avec un chapitre consacré à la géodésie tunisienne et son évolution depuis un siècle de sa mise en place. L'organisation de cette première partie de l'ouvrage est comme suit.

Dans le deuxième chapitre, on démontre les principales formules de la trigonométrie sphérique.

Le troisième chapitre présente les différents éléments de l'astronomie de position liés à la géodésie et en particulier les différents systèmes de coordonnées utilisés en astronomie de position.

Le quatrième chapitre est un rappel de la géométrie des courbes, le repère de Frenêt, la théorie des surfaces, la première forme fondamentale, et les théorèmes liés aux rayons principaux de courbure d'une surface de $\BbR^3$.\index{\textbf{Frenêt J.F.}}

La géométrie de l'ellipse et de l'ellipsoïde est l'objet du cinquième chapitre où on définit les formules des coordonnées tridimensionnelles d'un point, relatives à un ellipsoïde donné. On traite aussi les lignes géodésiques de l'ellipsoïde en présentant une méthode itérative de l'intégration de leurs équations. 
 
Dans le sixième chapitre, on donne les définitions des systèmes et des coordonnées géodésiques ainsi que du géoïde. On présente aussi les principaux systèmes géodésiques des pays de l'Afrique du Nord.    

Le septième chapitre traite les réseaux géodésiques terrestres et spatiaux. On présente leurs conceptions et réalisations. De même, les opérations de densification des réseaux terrestres et spatiaux par la technologie GPS sont traitées en donnant les principales phases.

Quant au huitième chapitre, il est consacré aux différentes corrections apportées aux distances observées lors de leurs réductions à la surface de l'ellipsoïde de référence et aux plans des représentations planes utilisées.

Le neuvième chapitre est l'un des chapitres importants de cet ouvrage où on traite les représentations planes et principalement celles qui sont conformes. Dans ce chapitre, on donne une démonstration de la condition de conformité d'une représentation plane. On présente aussi ce qu'on appelle en langage mathématique les représentations quasi-conformes en présentant un exemple.

Les chapitres dixième  et onzième sont consacrés respectivement à étudier en détail les représentations planes Lambert et l'UTM en démontrant pour chacune, les différentes formules des expressions des coordonnées rectangulaires $(X,Y)$ et du module linéaire.

Le datum altimétrique ainsi que les différentes définitions des systèmes d'altitudes font l'objet du douzième chapitre de l'ouvrage.

Parmi les modèles de passage entre les systèmes géodésiques, on traite en détail, dans le treizième chapitre, les modèles tridimensionnels de Bur$\breve{s}$a-Wolf, de Molodensky et le modèle bidimensionnel de Helmert.\index{\textbf{Helmert F.R.}}
\index{\textbf{Molodensky M.S.}}\index{\textbf{Wolf H.}}\index{\textbf{Bur$\breve{s}$a M.}} On présente une méthode de détermination directe des paramètres du modèle de Bur$\breve{s}$a-Wolf.

Des éléments historiques de la géodésie tunisienne sont présentés dans le quatorzième chapitre. On parlera des différents systèmes géodésiques tunisiens avant l'établissement du système NTT(Nouvelle Triangulation Tunisienne) le système géodésique terrestre officiel de la Tunisie. Ce chapitre n'a pas l'intention en tout cas de décrire l'historique de la géodésie tunisienne depuis les premiers travaux de rattachement des points géodésiques tunisiens à la géodésie italienne (\textit{C. Fezzani, 1979})\index{\textbf{Fezzani C.}}.

Dans le quinzième chapitre, on présente des notions sur le mouvement d'un satellite artificiel autour de la Terre avant d'entamer le positionnement par les satellites GPS objet du seizième chapitre.

La bibliographie relative à la \textbf{Partie I} est l'objet du dix-septième chapitre. 
\\
\\

- \textbf{Partie II:} 

Elle concerne une introduction à la théorie des moindres carrés pour les modèles linéaires avec une première présentation, dans un cours de géodésie destiné aux ingénieurs, de l'aspect non-linéaire de la méthode des moindres carrés. Cette deuxième partie comprend quatre chapitres en plus de la bibliographie.

En poursuivant la numérotation précédente des chapitres, le dix-huitième chapitre, le premier de la \textbf{Partie II}, traite les différentes définitions et théorèmes mathématiques qui seront utiles pour la théorie des moindres carrés.

Le chapitre important de la \textbf{Partie II} de l'ouvrage est le dix-neuvième chapitre où sont présentés les éléments fondamentaux de la théorie des moindres carrés des modèles linéaires. L'auteur a adopté la notation de P. Hottier utilisée dans son cours \textit{La Théorie des Erreurs} (\textit{P. Hottier, 1980}).

Le vingtième chapitre est consacré, pour la première d'un cours de la théorie des erreurs pour les ingénieurs, à l'aspect théorique de la géométrie de compensation d'un modèle non-linéaire par les moindres carrés. On rappelle les définitions nécessaires et on présente la méthode de Gauss-Newton pour la résolution du système donnant le minimum de la fonction \textit{objectif} ou encore dite fonction \textit{énergie}.

On termine avec le vingt-unième chapitre où on traite l'aspect géométrique des conditions obtenues pour la solution de la compensation par les moindres carrés des modèles non-linéaires.

Enfin, le chapitre vingt-deuxième constitue la bibliographie de la \textbf{Partie II} de l'ouvrage, suivi d'un index pour les noms propres et les mots clés cités dans tout l'ouvrage. 
 
Quant à l'aspect pratique, des exercices et des problèmes ont été ajoutés à la fin de la plupart des chapitres. De plus, des éléments historiques ont été formulés sous forme de notes historiques pour certains chapitres. 
 
 Enfin, pour signaler toute correction à cette publication, prière de nous écrire à l'adresse:  abenhadjsalem$@$gmail.com, et merci d'avance.   
\vspace{\baselineskip}
\begin{flushright}\noindent
Tunis,\hfill {\it \textbf{Abdelmajid}}\\
Décembre 2016\hfill {\it \textbf{Ben Hadj Salem, Dipl.-Ing.}}\\
\hfill{\it \textbf{Ingénieur Général Géographe}}
\end{flushright}

\vfill


\chapter*{\textit{\textbf{Remerciements}}}

\addcontentsline{toc}{section}{Remerciements}
Pour la réalisation de ce livre, nous avons profité des documents et de publications 
que nous citons ci-dessous à savoir (§ [\ref{biblio1},\ref{biblio2}]):

\textbf{Partie I:}
 
- La thèse de C. Fezzani (1979).\index{\textbf{Fezzani C.}}

- \textit{Cours de Cartographie Mathématique} de J. Commiot (1979).\index{\textbf{Commiot J.}}

- \textit{Cours de Géodésie Elémentaire} de J. Lemenestrel (1980).\index{\textbf{Lemenestrel J.}}

- \textit{Geodesy: the Concepts} de P. Vani$\vec{\breve{c}}$ek et E.J. Krakiwsky (1986).\index{\textbf{Vani$\vec{\breve{c}}$ek P.}}\index{\textbf{Krakiwsky E.J.}}
\\

\textbf{Partie II:}

- \textit{La Théorie des Erreurs} de P. Hottier (1980).\index{\textbf{Hottier P.}}

- La thèse de P.J.G. Teunissen (1985).\index{\textbf{Teunissen P.J.G.}}

- Les publications de E.W. Grafarend et B. Schaffrin (1989).\index{\textbf{Grafarend E.W.}} \index{\textbf{Schaffrin B.}}

- La publication \textit{Nonlinear Systems} de P.J. Olver (2013) de l'Université de Minnesota.\index{\textbf{Olver P.J.}}

Que leurs auteurs, dont certains étaient mes professeurs, trouvent ici l'expression de ma sincère gratitude.
 
Je voudrai aussi remercier mes collègues de l'OTC et particulièrement Messieurs les ingénieurs avec lesquels j'avait travaillé ou collaboré sans oublier mes étudiants qui avaient souvent souffert de mes cours exigeant un certain niveau mathématique.
\newpage
Je suis aussi très reconnaissant à Messieurs les ingénieurs  A. Boudabous, M. Charfi, C. Fezzani, M. Ben Mahmoud, K. Naouali, M. Rezgui, J. Zaiem, J. Zid, S. Kahlouche de l'Algérie et N. Rebai, Maître-assistant à l'Ecole Nationale des Ingénieurs de Tunis, pour avoir lu l'ouvrage et donné leurs commentaires et suggestions avant l'édition finale.

\vspace{\baselineskip}
\begin{flushright}\noindent
Tunis,\hfill {\it \textbf{Abdelmajid}}\\
Décembre 2016\hfill {\it \textbf{Ben Hadj Salem, Dipl.-Ing.}}\\
\hfill{\it \textbf{Ingénieur Général Géographe}}
\end{flushright}

\tableofcontents 







\clearpage



\setcounter{page}{1}
\pagenumbering{arabic}

%
%


\clearpage

%
%
    
 \normalsize                        

\mainmatter
\setcounter{page}{1}
\pagenumbering{arabic}
\part{\textit{Eléments de Géodésie}}

\chapter{\textit{\textbf{Introduction}}}
 \normalsize
\begin{tinted}
'' Consciente des résultats extraordinaires obtenus par les institutions
cartographiques nationales et les agences spatiales, les commissions géodésiques,
les instituts de recherche et les universités, et d'autres organisations internationales 
comme la Fédération internationale des géomètres, en faisant fond sur les initiatives
de l'Association internationale de géodésie, qui représente la communauté géodésique
mondiale, pour ce qui est d'évaluer et de suivre au mieux les changements du
système terrestre, notamment la mise au point du Repère de référence terrestre
international, qui a été adopté,...''\index{\textbf{UN-GGIM}}
\end{tinted}
\begin{flushright}
\footnote{Extrait de la Résolution A/RES/69/266, du Repère de référence géodésique mondial pour le développement durable, adoptée par l'Assemblée Générale des Nations-Unis le 26 février 2015(UN-GGIM, 2015).}
\end{flushright}
%
\section{\textsc{Définitions de la Géodésie}}
Suivant l'étymologie grecque, le mot géodésie \index{Géodésie}veut dire divise la Terre. Le grand géodésien Allemand \textbf{F.R. Helmert} (\textit{F.R. Helmert}, 1884)\footnote{\textbf{Friedrich Robert Helmert} (1843-1917): géodésien allemand.} définissait la Géodésie comme suit " \textbf{\emph{la Géodésie est la science de la mesure et de la représentation de la surface terrestre}}".\index{\textbf{Helmert F.R.}}
 
Une définition contemporaine de la Géodésie est donnée par le Comité Associé Canadien de Géodésie et de Géophysique (\textit{C.A.C.G.G.},1973) à savoir : \textbf{\emph{la Géodésie est la discipline qui concerne la mesure et la représentation de la Terre, incluant son champ de gravité, dans un espace tridimensionnel variant avec le temps}}.

Une autre définition récente (2002) est :''\textbf{\emph{ Geodesy is an interdisciplinary science which uses spaceborne and airborne remotely sensed, and ground-based measurements to study the shape and size of the Earth, the planets and their satellites, and their changes; to precisely determine position and velocity of points or objects at the surface or orbiting the planet, within a realized terrestrial reference system, and to apply these knowledge to a variety of scientific and engineering applications, using mathematics, physics, astronomy, and computer science.}}''(\textit{M. Lemmens}, 2011).\index{\textbf{Lemmens M.}}

La Géodésie a ainsi deux aspects :

*	un aspect scientifique et de recherches :

-	la mesure des  dimensions de la Terre et la détermination de sa forme géométrique.
\\

*	un aspect pratique :

-	l'établissement et la maintenance des réseaux géodésiques tridimensionnels nationaux et globaux et en tenant compte des variations de ces réseaux  en fonction du temps;

-	la mesure et la représentation  des phénomènes géodynamiques comme le mouvement des pôles, les marées terrestres et le mouvement de la croûte terrestre.

Dans cette première partie du livre, on s'intéresse aux réseaux géodésiques et à leurs établissements. 

Un \textbf{réseau géodésique} \index{Réseau géodésique} est un ensemble de points dont les coordonnées sont connues avec précision dans un système de référence donné. Ces points vont servir par la suite comme  points de référence pour tous les travaux topographiques et cartographiques.  

Cette première partie de l'ouvrage comprendra en plus de l'introduction les chapitres suivants:

2.	la trigonométrie sphérique;

3.	notions d'astronomie de position;

4.	courbes et surfaces;

5.	géométrie de l'ellipse et de l'ellipsoïde;

6.	les systèmes géodésiques;

7.	les réseaux géodésiques;

8.	réduction des distances;

9.	les représentations planes;

10. la représentation Lambert Tunisie;

11.	la représentation UTM;

12.	les transformations entre les systèmes géodésiques;

13.	les systèmes d'altitudes;

14. la géodésie tunisienne;

15. notions sur le mouvement d'un satellite artificiel autour de la Terre;

16.	le système GPS. 

\chapter {\textit{\textbf{La  Trigonométrie Sphérique}}}
      La trigonométrie sphérique établit les relations liant les grandeurs caractéristiques d'un triangle sphérique.
\section{\textsc{Le Triangle Sphérique}}\index{Triangle sphérique}

On considère une sphère de centre un point $O$ et de rayon l'unité et trois points sur la sphère $A,B,$ et $C$. 
\bdf
On appelle triangle sphérique la figure formée par les 3 arcs de grands cercles $AB,AC,$ et $CB$ inférieurs à $200$ grades (\textbf{Fig. \ref{doc11a}}).
\edf

\begin{figure}[ht]
\centering
\includegraphics[width=0.40\textwidth]{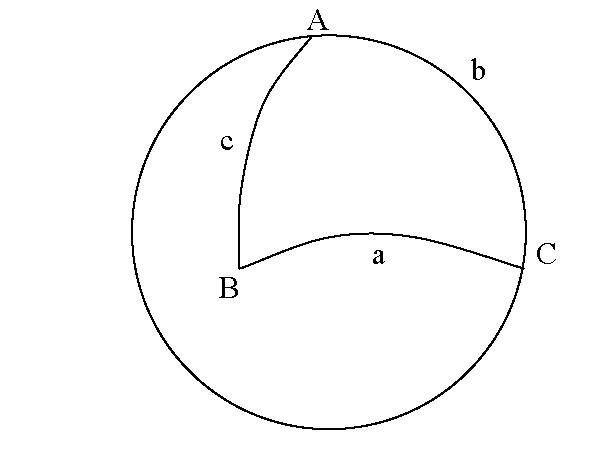}\\
\caption{Le triangle sphérique}
\label{fig:doc11a}
\end{figure}

Les grandeurs qui caractérisent le triangle sphérique $ABC$ sont :

-	les 3 côtés notés respectivement $a, b, c$, équivalents aux angles au centre des directions $OA,OB,OC$ soit 
$a = \widehat{(OB,OC)}, b =\widehat{(OA,OC)}, c =\widehat{(OA,OB)}$.

-	les 3 angles dièdres des faces du trièdre $OA,OB,OC$ notés $A,B,C$. 
\\

On remarque que les angles et côtés du triangle sont des grandeurs mesurables par des angles.
\section{\textsc{Le Trièdre Supplémentaire - Le Triangle Sphérique Polaire}}
      Au trièdre $OA,OB,OC$ on associe le trièdre supplémentaire dont les arrêtes $OA',OB',OC'$ sont respectivement orthogonales aux faces $OBC,OAC,OAB$. Le point $A'$ est choisi tel que $A$ et $A'$ soient dans la même demie sphère limitée par $BC$. Soit le point $C''$ diamétralement opposé au point $C$ (\textbf{Fig. \ref{fig:trianglepolaire}}). On a donc:

\begin{figure}
	\centering
		\includegraphics[width=0.50\textwidth]{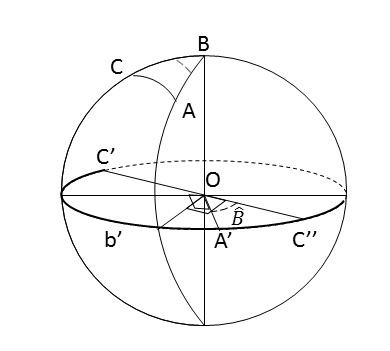}
	\caption{Le triangle sphérique polaire}
	\label{fig:trianglepolaire}
\end{figure}

$$\widehat{(OA,OC'')} = \pi-\widehat{(OBC,OAB)} = \pi-B=\widehat{(OA',OC')}$$
D'où les relations:
$$ \begin{array}{l}
\widehat{(OB',OC')}= a'= \pi - A \\
\widehat{(OA',OC')}= b'= \pi  - B  \\
\widehat{(OA',OB') }= c'=\pi -C  
  \end{array} $$
\bdf
      Le triangle sphérique $A',B',C'$ est dit triangle polaire du triangle $ABC$. 
			\edf
Comme le triangle $ABC$ est le triangle polaire de $A'B'C'$, on a :
$$ \begin{array}{l}
                           a = \pi  - A'   \\
                           b =  \pi - B'   \\   
													c = \pi  - C' 
  \end{array} $$
\section{\textsc{Les Formules de la Trigonométrie Sphérique}}
   Un triangle sphérique est entièrement défini par la donnée de 3 de ses 6 éléments. Alors entre 4 éléments quelconques, il y a :
	\[C^{4}_{6}=\frac{6!}{4!2!}=15 
\]
 relations non indépendantes comme suit :
 
-	3 côtés, 1  angle : 3 relations;

-	3 angles, 1 côté : 3 relations;

-	2 côtés, 2 angles(opposés aux côtés) : 3 relations;  
                                              
-	2 côtés, 2 angles (adjacents aux côtés) : 6 relations.
\subsection*{2.3.1. Etablissement de la Formule Fondamentale}
\begin{figure}
	\centering
		\includegraphics[width=0.50\textwidth]{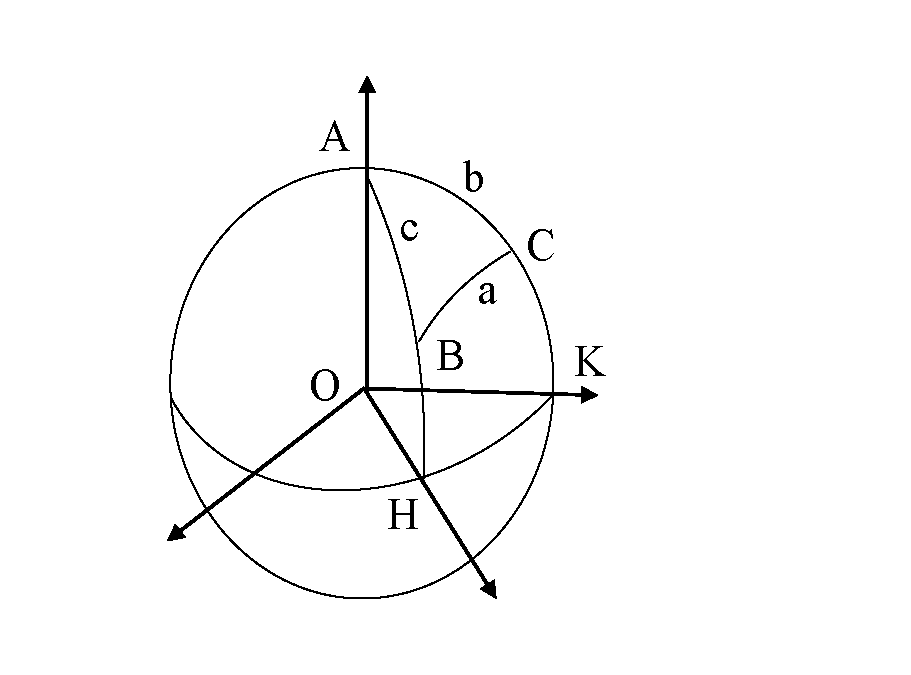}
	\caption{Calcul de la formule fondamentale}
	\label{fig:ffondamentale}
\end{figure}

      Soit un triangle sphérique $ABC$, en calculant le produit scalaire $\textbf{\textit{OB}}.\textbf{\textit{OC}}$ de 2 manières (\textbf{Fig. \ref{fig:ffondamentale}}), on arrive à la formule fondamentale :
\ba
\textbf{\textit{OB}} = cos(\pi /2 -c).\textbf{\textit{OH}} + sin( \pi /2 -c).\textbf{\textit{OA}} = sin c.\textbf{\textit{OH}} + cosc.\textbf{\textit{OA}}\nonumber  \\
 \textbf{\textit{OC}} = sinb.\textbf{\textit{OK}} + cosb.\textbf{\textit{OA}} \nonumber 
\ea
D'où:
	$$\textbf{\textit{OB}.\textit{OC}} = sinc.sinb.\textbf{\textit{OH}.\textit{OK}} + cosb.cosc $$
Or:
$$ \textbf{\textit{OH}.\textit{OK}} = ||\textbf{\textit{OH}}||.||\textbf{\textit{OK}}||.cos(\textbf{\textit{OH}},\textbf{\textit{OK}}) = 1.1.cosA=cosA$$
Et:
	\[\textbf{\textbf{\textit{OB}}.\textbf{\textit{OC}}} = ||\textbf{\textit{OB}}||.||\textbf{\textit{OC}}||.cos(\textbf{\textit{OB}},\textbf{\textit{OC}}) = 1.1.cosa = cosa 
\]
D'où : 
\be
    \fbox{$    cosa = cosb.cosc + sinb.sinc.cosA $}    \lb{p15}
		\ee
En utilisant le triangle polaire, on a :
$$	cosa' = cosb'.cosc' + sinb'.sinc'.cosA' $$
Or $ a'= \pi  - A,\,\,\,\, b' = \pi  - B,\,\,\,\mbox{et}\,\,\,\,c' = \pi  - C,\,\,\,\, a = \pi  - A'$, d'où :
\be
            \fbox{ $      cosA = - cosB.cosC + sinB.sinC.cosa   $}  \label{p16a}                       
\ee

\subsection*{2.3.2. La Formule des Sinus}
  De (\ref{p15}), on a : 
	\[cosA=\frac{cosa-cosbcocc}{sinbsinc}
\]
 Soit $sin^2A = 1 - cos^2A$, on arrive à :      
	$$ \frac{sin^2A}{sin^2a}=\frac{sin^2B}{sin^2b}=\frac{sin^2C}{sin^2c} $$ 
D'où:                                                                                                   
\be
\fbox{$ \ds 	\frac{sinA}{sina}=\frac{sinB}{sinb}=\frac{sinC}{sinc} $} \label{p18}
\ee
\subsection*{2.3.3. Formules des Sinus Cosinus}
      En utilisant la formule fondamentale, on a :
\ba
cosa = cosb.cosc + sinb.sinc.cosA \nonumber \\
cosb = cosa.cosc + sina.sinc.cosB \nonumber
\ea
Et en remplaçant dans la deuxième formule, l'expression de $cosa$, on obtient $sinc.cosb = sinb.cosc.cosA + sina.cosB$, d'où :
\be
      \fbox{$  sina.cosB = cosb.sinc- cosc.sinb.cosA $}                   \label{p19}
\ee
\subsection*{2.3.4. Formule des Cotangentes}
      En remplaçant dans (\ref{p19}) $sina$ par $sinA.sinb/sinB$, on obtient :
\be         
            \fbox{$   sinA.cotgB = cotgb.sinc - cosc.cosA   $}                     \label{p20}      
\ee
\subsection*{2.3.5. Cas d'un Triangle Rectangle}
\begin{figure}
	\centering
		\includegraphics[width=0.40\textwidth]{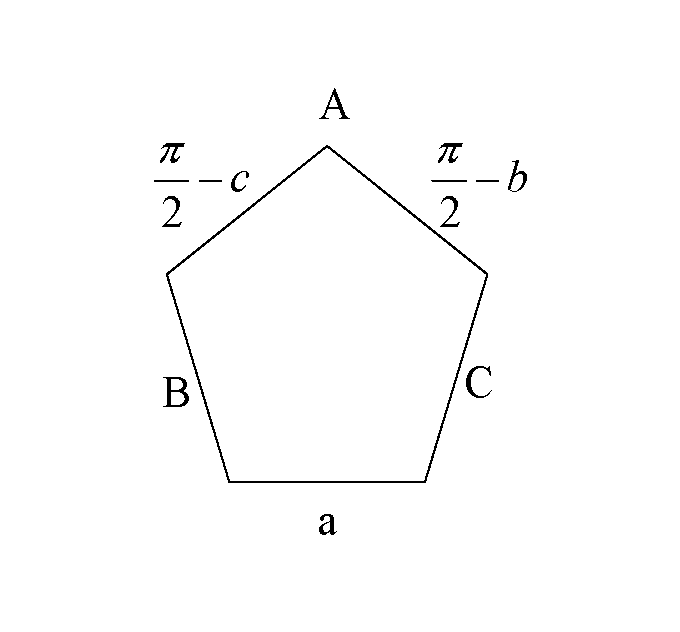}
			\caption{La règle de Nepier}
	\label{fig:neper}
\end{figure}

      Pour un triangle sphérique rectangle, un des angles vaut $\pi/2$ = 100 gr = 90°. Les formules se simplifient, leur nombre est :
 	\[C^{3}_{5}=\frac{5!}{3!2!}=10 
\]
Supposons que A = $\pi/2$, on fait le schéma ci-dessous (\textbf{Fig. \ref{fig:neper}}).
\\

On trouve les relations en appliquant la règle mnémonique de Nepier\footnote{\textbf{John Nepier} (1550 -1617): mathématicien écossais. }: \index{\textbf{Nepier J.}}

\textit{Le cosinus d'un élément quelconque est égal à :}

-	\textit{au produit des cotangentes des éléments adjacents;}

-	\textit{au produit des sinus des éléments non adjacents}.

Exemple : 

- $ cos a = cotgB.cotgC,$

- $cosa = sin(\pi/2-c).sin(\pi/2-b) = cosc.cosb$.
\section{\textsc{L'Excès Sphérique}}
\begin{figure}
	\centering
		\includegraphics[width=0.50\textwidth]{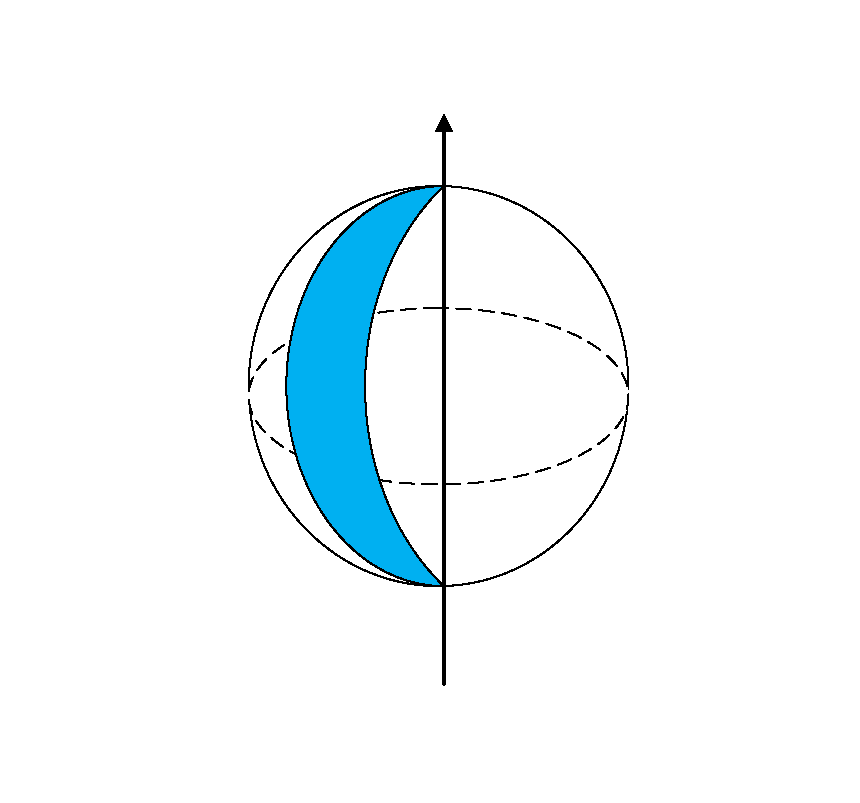}
	\caption{Un fuseau sphérique}
	\label{fig:fuseau}
\end{figure}

\index{Excès sphérique}
\bdf
On appelle fuseau sphérique la portion de la demi sphère limitée entre deux grands cercles (\textbf{Fig. \ref{fig:fuseau}}).
\edf
 La surface d'un fuseau sphérique d'un angle $A$ est proportionnelle à $AR^2$ où $R$ est le rayon de la sphère, soit $S = kAR^2$, pour $A = 2\pi$ on a $S = 4\pi R^2 = k2\pi R^2$ d'où $k = 2$, on obtient : 
	\[S = 2AR^2
\]

On considère maintenant un triangle sphérique ABC : 

- le fuseau $(AB,AC)$  donne $S_{1}= 2AR^2$;      

- le fuseau $(CA,CB)$  donne $S_{2}= 2CR^2$;      

- le fuseau $(BC,BA)$  donne $S_{3}= 2BR^2$;      

d'où : 
	\[S_{1} +  S_{2} + S_{3} = 2R^2(A+B+C)
\]
Or $S_{1} +  S_{2} + S_{3}$ = la surface de la demi-sphère + 2 fois la surface du triangle sphérique $ABC$. On note $T$ la surface du triangle sphérique ABC, on a alors:
	\[2R^2( A+B+C) = 2\pi R^2 + 2T
\]
ou encore : 
$$                                            A + B+ C =  \ds \pi  + \frac{T}{R^2}=\pi+\epsilon $$
Soit:
\be
                 \fbox{$      \epsilon (rd) =\ds  \frac{T}{R^2}=\frac{Aire\, ABC}{R^2}=\textrm{excès sphérique} $} \label{p22} 
\ee

  \section{\textsc{Exercices et Problèmes }}
\bex
 Calculer l'azimut d'une étoile de déclinaison $\delta  =+5^{\circ}$ quand sa distance zénithale est de $80^{\circ}$ pour un observateur situé à la latitude $\varphi= 56^{\circ}$. 
\eex
\bex
 En appliquant au triangle de position les formules de trigonométrie sphérique montrer que l'on peut calculer l'angle horaire $AH_c$ du coucher d'un astre par : $cos AH_c = -tg\varphi.tg\delta $.   
\eex
\bex
 Soit un triangle sphérique $ABC$. On donne les éléments suivants:

- $\hat{A}=80.1643\,3\,gr$;

- $\hat{B}=55.7735\,1\,gr$;

- $\hat{C}=64.0626\,1\,gr$;

- $AC= 20.1357\,km$;

- $AB= 22.1435\,km$.

1. Calculer $\alpha=\hat{A}+\hat{B}+\hat{C}$.

2. Déterminer $\epsilon$ l'excès sphérique de ce triangle.

3. Calculer la fermeture du triangle $ABC$, donnée par:
$$ f=\alpha -200.0000\,0\,gr-\epsilon$$
\eex
\bex
 Soit $(\BbS^2)$ une sphère de rayon égal à 1. Soit un carré sphérique $ABCD$ de côté $a$ (arc de grand cercle). On note $\alpha=\hat{A}=\hat{B}=\hat{C}=\hat{D}$. 

\begin{figure}
	\centering
		\includegraphics[width=0.50\textwidth]{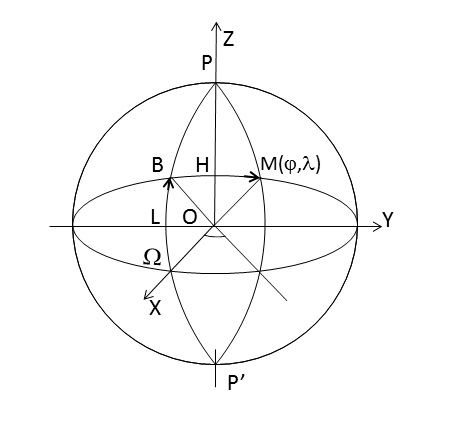}
	\caption{Les coordonnées de Cassini-Soldner}
	\label{fig:soldner}
\end{figure}

1. Montrer que:
$$ cos a=cotg^2\frac{\alpha}{2}$$
2. Donner l'expression de la diagonale $d=l'arc\,AC$.  
\eex
\bpb
Soit $(\BbS^2)$ une sphère de rayon égal à 1 et de centre le point $O$. Un point $M$ de 
$(\BbS^2)$ a pour coordonnées $(\varphi,\lambda)$. On appelle les coordonnées de Cassini-Soldner\index{Coordonnées de Cassini-Soldner}\index{\textbf{Soldner J.G}} \footnote{\textbf{César-François Cassini} (1714-1784): astronome et géodésien français.}\footnote{\textbf{Dr Johann Georg von Soldner} (1776-1833): mathématicien et astronome bavarois.} de $M$ les angles (\textbf{Fig. \ref{fig:soldner}}):\index{\textbf{Cassini C.F.}}

- $L=\widehat{O\Omega,OB}=\,Arc\,\,\Omega B$;

- $H=\widehat{OB,OM}=\,Arc\,\,BM$.

1. Déterminer les relations liant $L,H$ à $\varphi,\lambda$.

2. Inversement, donner les relations liant  $\varphi,\lambda$ à $L,H$.

\epb
\bpb
 Au lieu $M$ de latitude $\varphi= 38^{\circ}$  Nord, on observe l'étoile polaire  $A$ de déclinaison $\delta = + 89^{\circ}$ et d'ascension droite $\alpha$  = $+2\,h\, 13\,mn\, 52.90\,s$.
 
1.	Donner sur un graphique, les éléments du triangle sphérique $PAM$ où $P$ est le pôle Nord. 

2.	Sachant que l'heure sidérale locale $HSL$ est égale au moment de l'observation à  $6\,h\, 37\, mn\, 19.72\, s$,  calculer l'angle horaire $AH$.  

3.	En appliquant la formule des cotangentes, montrer  que l'azimut $Az$ de l'étoile est donné par la formule:
$$tgAz=\frac{sinAH}{cosAHsin\varphi-cos\varphi tg\delta}$$
4.	Calculer alors l'azimut $A$z.

5.	Calculer la distance zénithale $z$ de l'étoile.
\epb
\bpb
Soit la sphère unité $(\BbS^2)$ de $\BbR^3$. On considère le triangle sphérique $ABC$ avec:
$$ A=\begin{pmatrix}{
0\cr
0\cr
1}
\end{pmatrix},\quad B=\begin{pmatrix}{
cos\fii_B\cr
0\cr
sin\fii_B}
\end{pmatrix},\quad C=\begin{pmatrix}{
cos\fii_C cos\lm_C\cr
cos\fii_C sin\lm_C\cr
sin\fii_C}
\end{pmatrix}$$
et $M$ un point quelconque de $(\BbS^2)$: $ M=\begin{pmatrix}{
cos\fii cos\lm\cr
cos\fii sin\lm\cr
sin\fii}
\end{pmatrix}$. On appelle $(\Gamma)$ le grand cercle de  $(\BbS^2)$ passant par les points $B$ et $C$.

1. Montrer que si $M(\fii,\lm) \in (\Gamma)$, alors $\fii=\Phi(\lm)$ avec:
\be
cos\fii_Csin\fii_Csin\lm_Csin\Phi=cos\Phi(sin\fii_Bcos\fii_Csin(\lm_C-\lm)+cos\fii_Bsin\fii_Csin\lm) \lb{tsf}
\ee
2. Exprimer la surface du triangle sphérique $ABC$ en utilisant l'intégrale de la fonction $sin\Phi$ entre deux bornes que l'on déterminera.
\\ 

3. On introduit l'angle $\omega$ que fait le vecteur $T_{\Gamma}$ tangent en $M$ au grand cercle $(\Gamma)$ avec le vecteur $T_m$ tangent en $M$ au méridien passant par $M$. Vérifier que $\ds \frac{d\omega}{d\lm}=sin\fii$ (aide: on peut dériver, par rapport à $\lm$, deux fois l'équation (\ref{tsf})).
\\

4. En déduire la valeur de l'aire $S$ du triangle sphérique $ABC$.  

\epb
\chapter{\textit{\textbf{Notions d'Astronomie de Position}}}
\section{\textsc{Rappels Historiques}}
Il n'est pas possible de déterminer la période où l'homme avait commencé à s'intéresser aux étoiles et aux astres. Cependant, cet intérêt à l'astronomie s'était développé peu à peu (\textit{P. Vanic}$\breve{e}$\textit{k} \& \textit{E.J. Krakiwsky}, 1986). Ainsi, les astronomes chinois s'étaient occupés de l'observation des astres et des étoiles (observation d'une éclipse solaire 2667 avant J.C.), de même pour les peuples habitants la région des rivières le Tigre et l'Euphrate, les peuples d'Egypte et les Grecs. On trouve par exemple \emph{Thales de Miletus}\index{\textbf{Thales D.M}} qui prédit l'éclipse solaire en mai 685 avant J.C. \emph{Eratosthenes} (276-194 avant J.C.)\index{\textbf{Eratosthenes}} calcula la circonférence de la Terre (39000 km), l'astronome et mathématicien \emph{Hipparque} (190-120 avant J.C.)\index{\textbf{Hipparque}} avait trouvé les plus importantes irrégularités du mouvement apparent du soleil et de la lune et il découvrit la précession (mouvement du pôle). 
\\

L'astronome \emph{Claude Ptolémée} (80-168)\index{\textbf{Ptolémée C.}} avait compilé toutes les théories d'astronomie de l'époque que les Arabes avaient traduit en un ouvrage appelé "Al-Megiste". La civilisation islamique a apporté aussi sa contribution dans les connaissances astronomiques surtout sous les règnes d'\emph{Abu Jaafar Al Mansour} (754-775)\index{\textbf{Al Mansour A.J.}}, \emph{Haroun Al Rachid} (786-809)\index{\textbf{Al Rachid H.}} et \emph{Abu Jaafar Al Mamoun} (812-833)\index{\textbf{Al Mamoun A.J.}}. Parmi les astronomes arabes, on cite : \emph{Mohamed Ibn Jabir Al Battani} (mort en 929)\index{\textbf{Al Battani M.J.}} qui avait fait des catalogues astronomiques des mouvements des planètes. Ses travaux furent traduits en latin et imprimés à Nuremberg (Allemagne) en 1537. Il calcula la durée de l'année solaire à 2 mn près.
\\

\emph{Ibn Al Haitam}\index{\textbf{Ibn Al Haitam H.}} (965-1039), dit \emph{Al Hazan} en Europe, avait étudié la réfraction de l'atmosphère et calcula l'altitude de l'atmosphère pour la première fois. On cite aussi \emph{Abu Rayhan Al Biruni} (973-1055)\index{\textbf{Al Biruni A.R.}} qui était aussi un éminent astronome en plus de ses travaux en mathématiques. \emph{Nassireddine Al Tusi} (1201-1274) \index{\textbf{Al Tusi N.}}construisit le premier observatoire moderne et de plus, il avait développé la trigonométrie sphérique. 
\\

Au 15ème siècle, c'est le commencement de l'intérêt en Europe à l'astronomie grâce au monde Arabe. \emph{Nicolas Copernic} (1473-1543) \index{\textbf{Copernic N.}}a démenti la théorie qui dit que la Terre est le centre de l'Univers et a démontré que le soleil est le centre du système solaire. \emph{Tycho-Brahé} (1546-1609)\index{\textbf{Brahé T.}} l'éminent astronome danois a adhéré à la théorie de Plotémé avant l'invention du télescope. \emph{Johannes Kepler} \index{\textbf{Kepler J.}}(1571-1630) découvrit la forme elliptique des orbites des planètes. \emph{Galileo Galilée} (1564-1642)\index{\textbf{Galilée G.}} astronome italien, était le premier qui a fait usage du télescope. \emph{Isaac Newton} (1642-1727)\index{\textbf{Newton I.}} mathématicien et physicien anglais découvrit la loi de la Gravitation (1660).
\\

L'astronomie a connu un grand essor au 18ème siècle grâce aux mathématiciens de cette période comme  \emph{Leonhard Euler} (1707-1783)\index{\textbf{Euler L.}}, \emph{Joseph-Louis Lagrange} (1736-1813)\index{\textbf{Lagrange J.L.}} et aussi à la création des observatoires astronomiques comme ceux de  Paris en 1667 par \emph{Jean Dominique Cassini} (1625-1712)\index{\textbf{Cassini J.D.}} et de  Greenwich en 1675. \emph{Edmont Halley} (1656-1742) \index{\textbf{Halley E.}}avait compilé un catalogue des positions de l'hémisphère Sud et détecta une comète en 1682 que porta son nom après sa mort.
\\

Actuellement l'astronomie s'est avancée grâce à l'introduction de nouvelles méthodes telles que la photographie ou l'analyse spectrale, et l'envoi des satellites artificiels dont le premier a été lancé en octobre 1957.
\section{\textsc{Objectifs de l'Astronomie}}
Pour le géodésien, l'astronomie est un moyen de détermination de certaines inconnues du point stationné à partir d'observations sur les astres ou des étoiles. Les observations astronomiques effectuées dans ce cadre déterminent la verticale physique du point de l'observation, celle-ci étant matérialisée par l'axe de rotation de l'instrument. L'astronomie physique fournit comme résultat la distribution des verticales aux différents points stationnés. Si on assimile la verticale à la normale à la surface modèle de référence, on peut alors localiser ces points. On parlera alors d'\textit{astronomie de position}. En géodésie tridimensionnelle, l'astronomie donne la direction de la tangente à la ligne de force du champ de pesanteur au point considéré.
\\

Cependant, la géodésie ne peut se détacher de l'astronomie. En effet, il a toujours fallu, pour placer les points sur la sphère ou l'ellipsoïde de référence ou dans un trièdre trirectangulaire,  fixer les axes des coordonnées. Alors un des axes privilégiés est l'axe de rotation de la Terre. Ce dernier n'est pas matérialisé sur la surface topographique, mais il apparaît dans l'observation du mouvement de la Terre ou dans l'observation des étoiles. Donc, le géodésien est nécessairement astronome. Alors, les observations astronomiques permettent en géodésie de déterminer :

-	les 2 inconnues fixant la direction de la verticale physique du lieu ($\phi$,$\lambda$);

-	l'orientation d'une direction (l'azimut);

-	les coordonnées absolues d'un premier point d'un réseau géodésique appelé aussi point fondamental.
\section{\textsc{Les Systèmes de Référence}}
Le principe fondamental des déterminations astronomiques repose sur le fait que dans le repère lié aux étoiles, celles-ci occupent des positions pratiquement fixes, qu'il est  possible de calculer et de les regrouper en catalogues d'étoiles. Un catalogue d'étoiles comprend les coordonnées équatoriales célestes ($\alpha$,$\delta$) des étoiles observées, réduites à une époque moyenne conventionnelle. Un catalogue fondamental est issu de la compilation de plusieurs catalogues provenant de préférence d'observations absolues.  
\\

		Le système pratique de référence est défini par le catalogue fondamental, adopté à l'échelle internationale. Le catalogue fondamental FK4 était publié en 1963. Le catalogue le plus récent est le FK6 daté de 2000. Il est imprimé en deux volumes regroupant les données de 4150 étoiles.
\section{\textsc{Notions d'Astronomie de Position}}
\subsection*{3.4.1. Sphère Céleste - Mouvement Diurne}
\bdf
La sphère céleste \index{Sphère céleste} est une sphère de rayon infiniment grand sur laquelle sont projetées les perspectives des étoiles (\textbf{Fig. \ref{Sceleste2}}).

On appelle constellation la figure formée par les étoiles.
\edf

En regardant les étoiles, on s'aperçoit que les étoiles se déplacent dans leur ensemble : c'est le mouvement diurne (\textit{F. Tisserand \& H. Andoyer}, 1912).\index{\textbf{Tisserand F.}}\index{\textbf{Andoyer H.}}

Le mouvement diurne \index{Mouvement diurne}  obéit à 3 lois:

-	la sphère céleste tourne autour d'un de ses diamètres;

-	le mouvement s'effectue dans le sens rétrograde (non direct);

-	le mouvement est uniforme et sa période est voisine de $24\,h\, (23\,h\, 56\, mn)$. 
\subsection*{3.4.2. Définitions:}
\bdf
L'axe du monde est le diamètre autour duquel la sphère céleste effectue son mouvement.
\edf
\textbf{Pôles célestes} : P, P'; P pôle nord, P' pôle sud.
\begin{figure} [htp]
\centering
\includegraphics[width=0.50\textwidth]{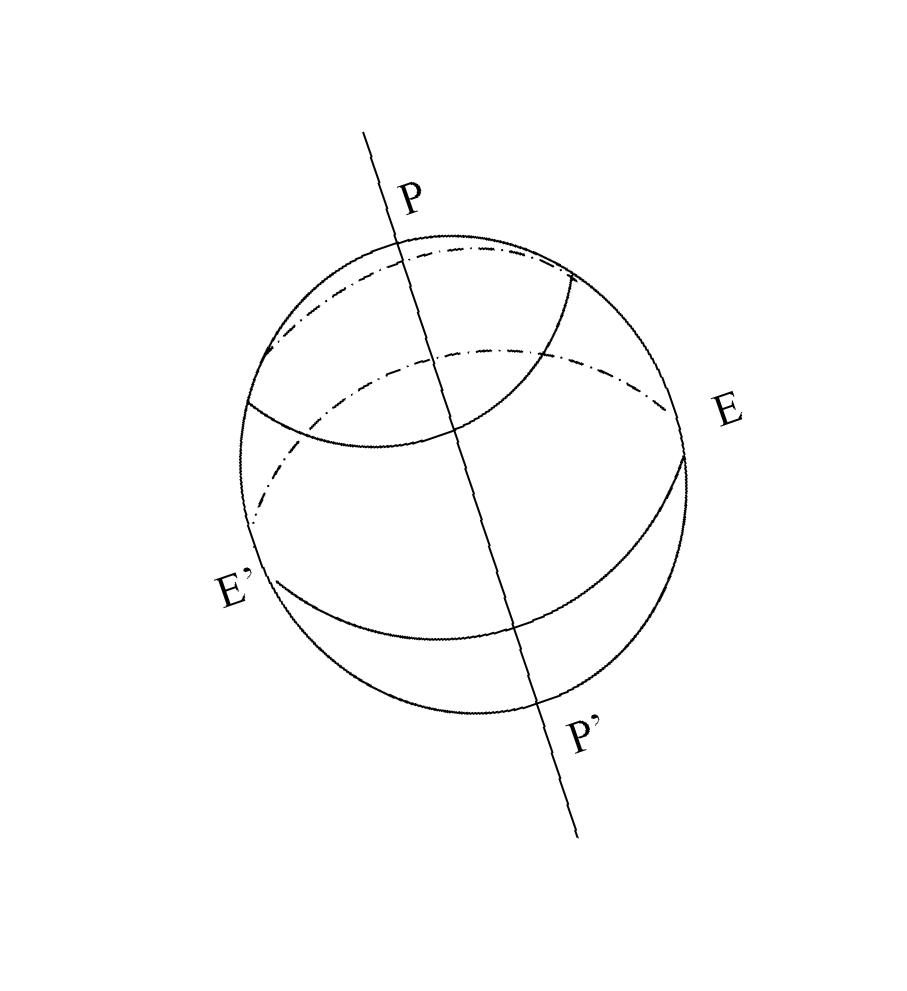}\\
\caption{La sphère céleste}
\label{Sceleste2}
\end{figure}

\begin{figure}
	\centering
		\includegraphics[width=0.40\textwidth]{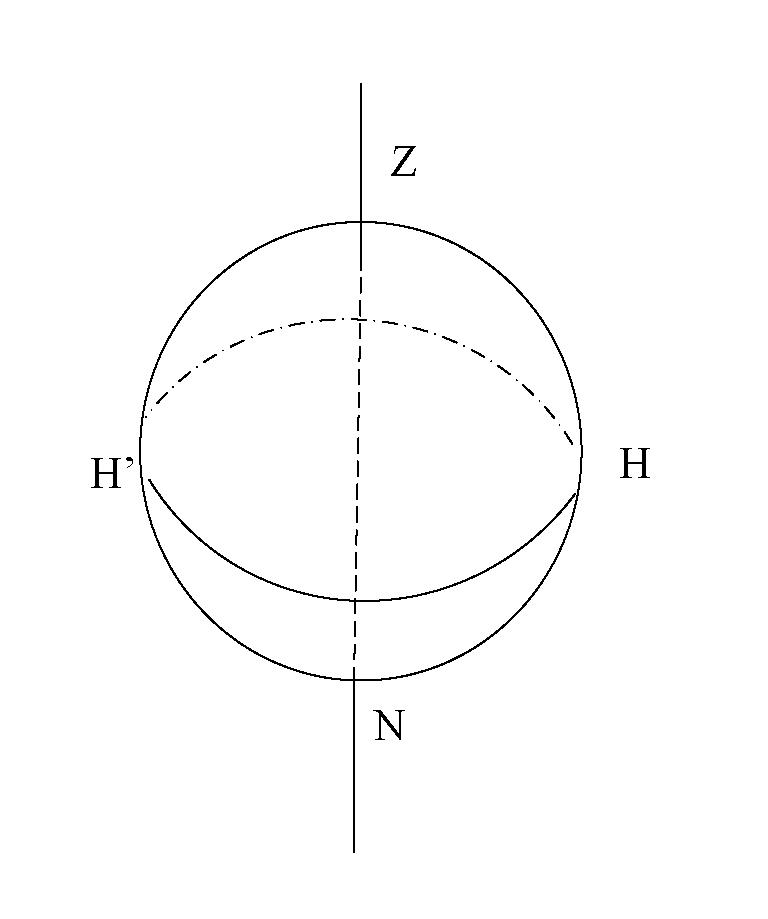}
	\caption{La verticale d'un lieu}
	\label{vertical}
\end{figure}

\begin{figure}
	\centering
		\includegraphics[width=0.40\textwidth]{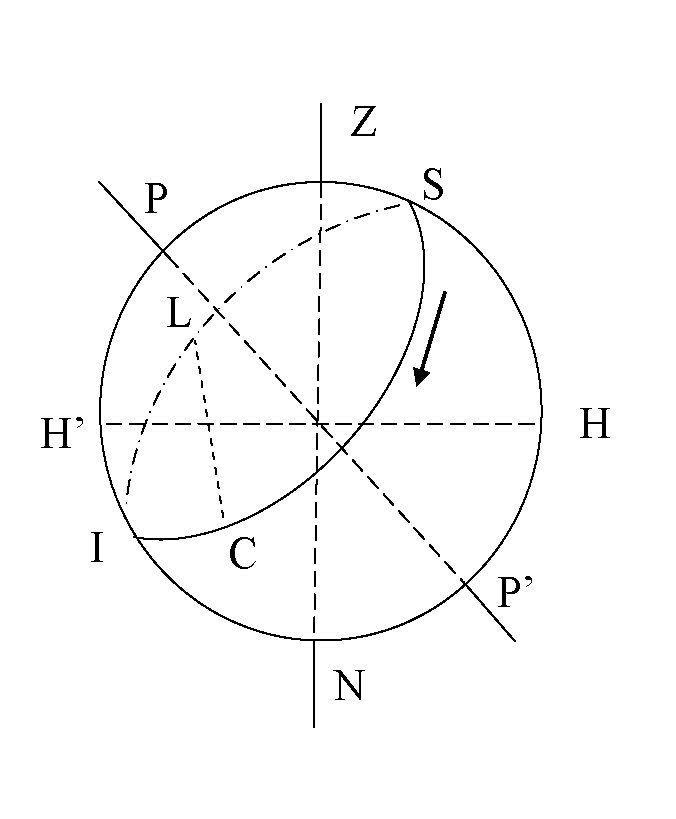}
	\caption{Le plan méridien}
	\label{fig:planmeridien}
\end{figure}

\bdf
L'équateur céleste est le grand cercle perpendiculaire à l'axe du monde.
L'intersection de la sphère céleste par un plan parallèle à EE' est un parallèle céleste. Un parallèle céleste est la trajectoire des étoiles en mouvement diurne. (\textbf{Fig. \ref{Sceleste2}})
\edf
 \bdf
La verticale d'un lieu est la direction donnée par un fil à plomb: Z c'est le zénith, N c'est le nadir. L'horizon est le grand cercle dont le plan est perpendiculaire à ZN (\textbf{Fig. \ref{vertical}}).
\edf
\bdf
Le plan méridien d'un lieu est le plan  défini par la verticale et l'axe du monde (\textbf{Fig. \ref{fig:planmeridien}}). Le méridien d'un lieu c'est un grand cercle intersection du plan méridien et de la sphère céleste. Le méridien est local.

Le demi-méridien supérieur : c'est le plan passant par PP' et contenant le zénith.
\edf

La trajectoire de l'étoile rencontre le méridien en deux points S et I:

-	S est \underline{le passage supérieur} (le plus près du zénith) ou culmination;

-	I est \underline{le passage inférieur}.
\\

\begin{figure}
	\centering
		\includegraphics[width=0.70\textwidth]{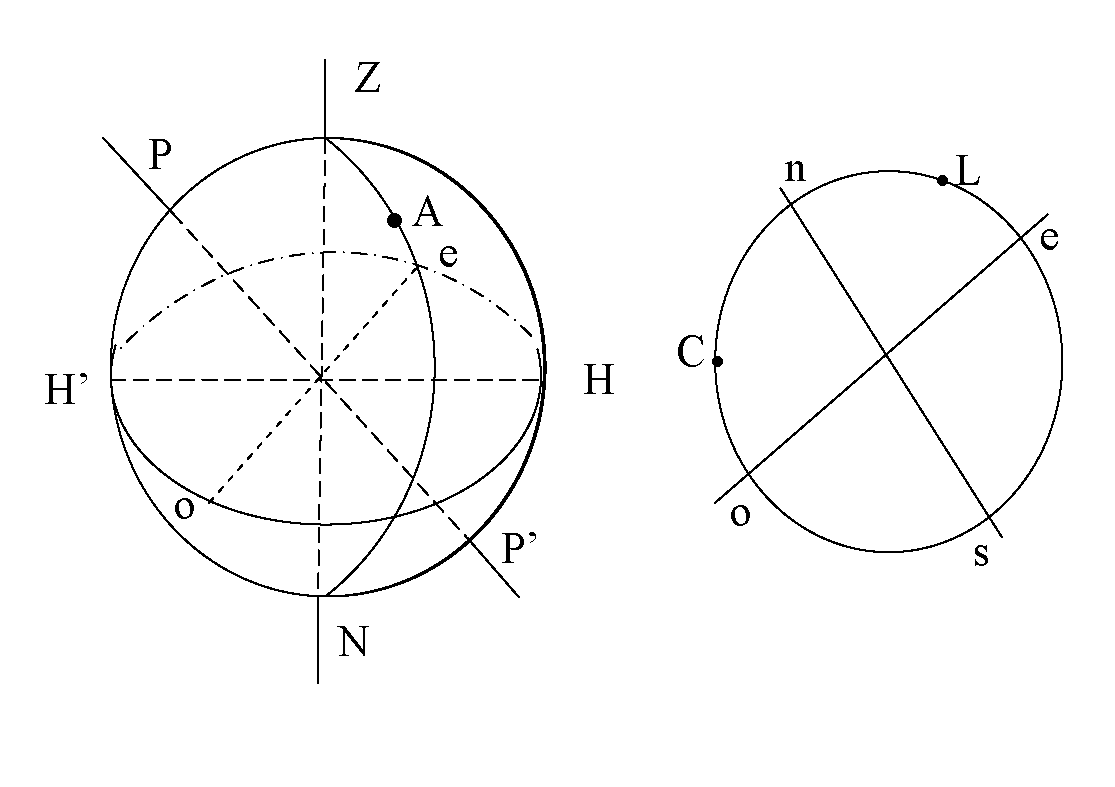}
	\caption{Vertical d'un astre}
	\label{fig:verticalastre}
\end{figure}

La trajectoire de l'étoile rencontre en général l'horizon en deux positions: 
\begin{itemize}
	\item L: est le lever, où l'étoile devient visible;
	\item C: est le coucher où l'étoile disparaît.
\end{itemize}
Les étoiles qui n'ont ni coucher, ni lever sont appelées étoiles circumpolaires.
\bdf
Le plan vertical est un plan contenant la verticale ZN.
\edf
\bdf
On appelle  vertical d'un astre  le plan vertical passant par l'astre (\textbf{Fig. \ref{fig:verticalastre}}).
\edf
\bdf
Le méridien est le vertical passant par le pôle; il rencontre l'horizon en un point n: c'est le Nord géographique, le point opposé au Nord c'est le Sud. Dans la direction perpendiculaire, on a l'Est et l'Ouest. L'Est se trouve à droite de la ligne Sud-Nord.
\edf
Pour un astre : le lever dans l'Est, à partir de la culmination il passe dans la région Ouest c'est le coucher.

\subsection*{3.4.3. Rappels d'Unités de Mesures des Arcs}
\begin{figure}
	\centering
		\includegraphics[width=0.40\textwidth]{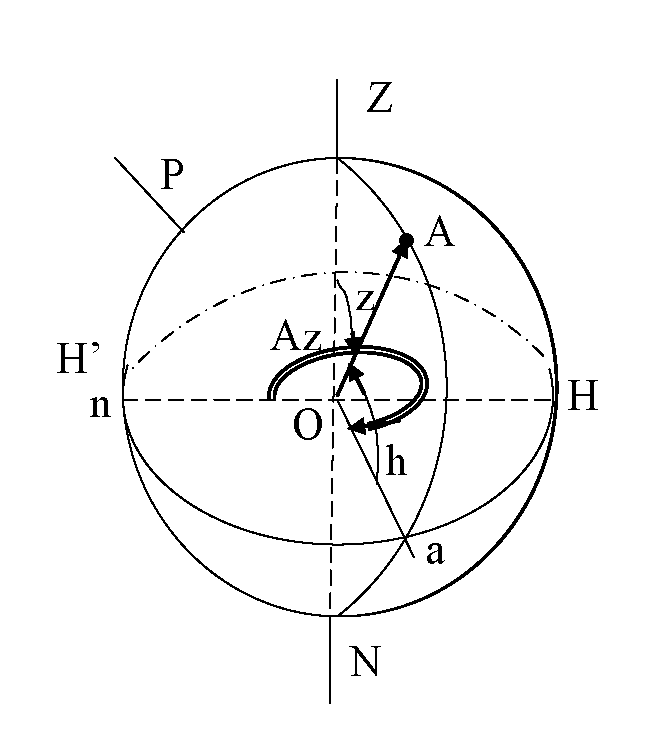}
	\caption{Les coordonnées locales horizontales}
	\label{fig:azimut1}
\end{figure}

On rappelle que le \textit{radian} $(rd)$ est l'unité internationale des mesures des angles. On donne ci-dessous les autres unités utilisées.
\begin{table}[htpb]
	\centering
	
\begin{center}
\[\begin{array}{ccc}
\hline
   \mbox{Le Système Centésimal}& \mbox{Le Système Sexagésimal}&\mbox{Le Système Horaire}\\ \hline
                400 gr       &              360^{\circ}           &               24 h     \\   \hline
          1\ grade =  1/64 rd&       1^{\circ} = 1/57 rd          &       1\ heure = 15^{\circ}     \\   \hline
          1\ cg =  1/6400 rd &      1\ minute = 1'=1^{\circ}/60   &        1\ mn = 1h/60 = 15'\\  \hline 
     1\ dcmg = 1/640000 rd   &    1\ seconde=1"=1'/60      &     1\ s = 1mn/60 = 15"   \\ \hline
\end{array}
\]
\caption{Table des Unités}
	\label{tab: Tableau des Unités}
\end{center}
\end{table}
\subsection*{3.4.4. Systèmes de Coordonnées Locales Horizontales (\textbf{Fig. \ref{fig:azimut1}})}\index{Coordonnées locales horizontales}
\bdf
L'azimut d'un astre \index{Azimut d'un astre} est l'angle formé par le vertical du l'astre et le plan méridien. Il est compté à partir du Nord dans le sens rétrograde (en grades).   
\edf
\be 
	 Az=\widehat{nOa}
\ee
\bdf
 La hauteur d'un astre \index{Hauteur d'un astre} est l'angle entre la direction de l'astre et l'horizon, compté à partir de ce dernier, positivement vers le zénith et négativement vers le nadir (en degrés). On le note par $h$.
\edf
\be 
	 h=\widehat{aOA}
\ee
\bdf
La distance zénithale \index{Distance zénithale} z est l'angle de la verticale avec la direction de l'astre, se compte du zénith vers le nadir. 
\edf
On a :
\ba
							z=\widehat{ZOA} \\
				 0\leq  z \leq 200\,gr \,\,\,\,\,\mbox{et}\,\,\,\,\,z + h = 100\,gr
\ea
\subsection*{3.4.5. Coordonnées Géographiques}\index{Coordonnées géographiques}
Soit M un point de la surface de la Terre.
\\
\bdf 
La latitude géographique $\phi$ est l'angle du plan de l'équateur avec la verticale du lieu, positivement vers le pôle Nord, négativement vers le pôle Sud.
\edf
\bdf
La longitude géographique $\lambda$ est l'angle formé par le méridien origine avec le méridien du lieu. Le méridien origine est le méridien passant par l'observatoire de Greenwich. $\lambda$ est comptée positivement vers l'Est en grades, degrés ou en heures. 
\edf 
\begin{figure}[htpb]
	\centering
		\includegraphics[width=0.60\textwidth]{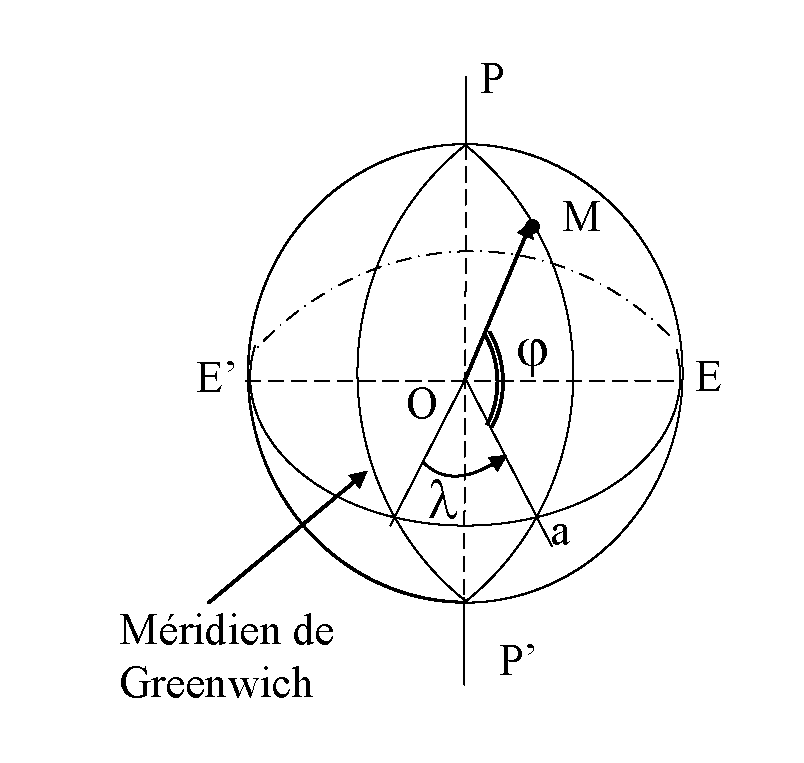}
	\caption{Latitude et longitude Géographiques}
	\label{fig:sphere1}
\end{figure}
\subsection*{3.4.6. Angle Horaire - Coordonnées Equatoriales Célestes - Temps Sidéral -}
\begin{figure}
	\centering
		\includegraphics[width=0.80\textwidth]{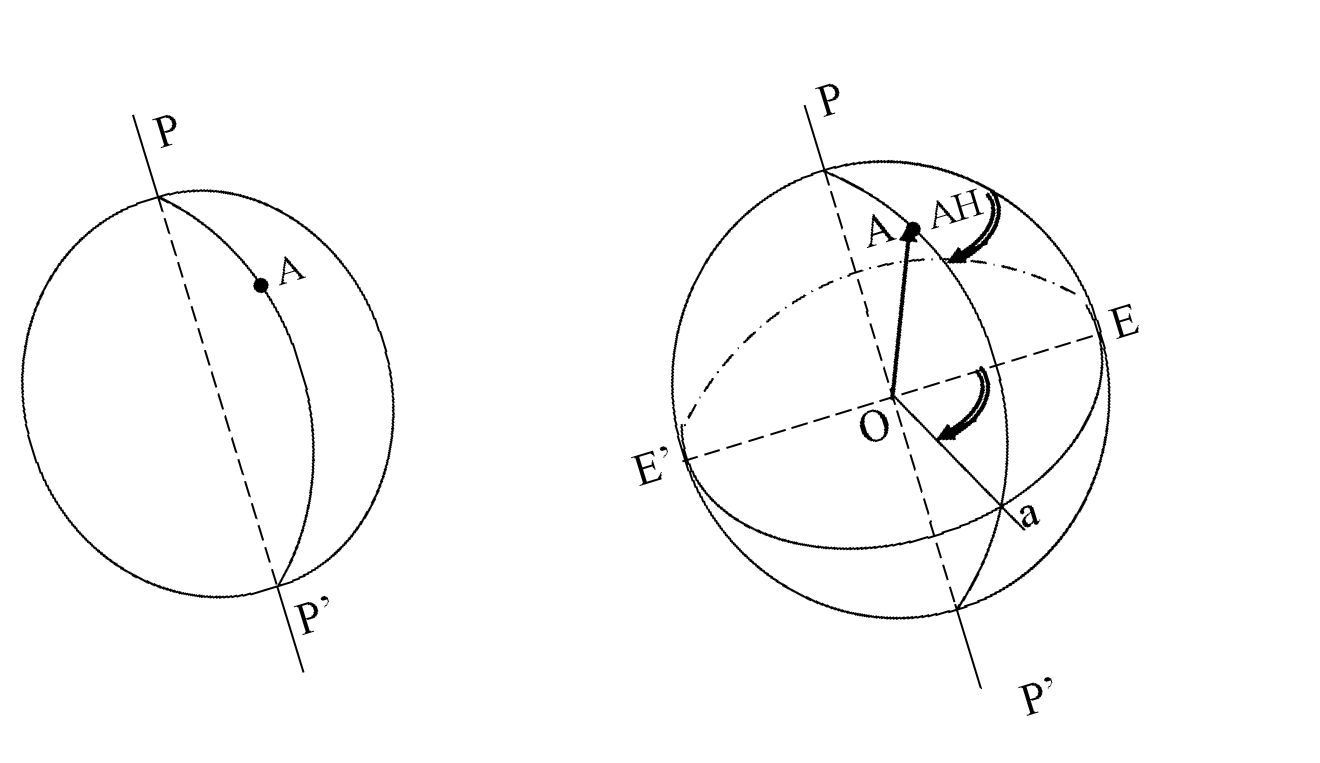}
	\caption{Cercle et angle horaires}
	\label{fig:anglehoraire}
\end{figure}

\bdf
Le cercle horaire d'un astre  \index{Cercle horaire de l'astre} est le demi-grand cercle passant par le pôle et l'astre (\textbf{Fig. \ref{fig:anglehoraire}}).
\edf
Le cercle horaire passant par le zénith c'est le méridien supérieur.                 
\bdf
L'angle horaire  \index{Angle horaire} $AH$ d'un astre $A$ en un lieu donné est l'angle formé par le méridien supérieur du lieu et  le cercle horaire de l'astre. Il se compte en heures dans le sens rétrograde à partir de la culmination (\textbf{Fig. \ref{fig:anglehoraire}}).
\edf
 \be
	AH= \widehat{EOa}  \label{p1}
\ee

\newpage

\textbf{Coordonnées Equatoriales  Célestes}\index{Coordonnées équatoriales célestes}

\begin{figure}
	\centering
		\includegraphics[width=0.50\textwidth]{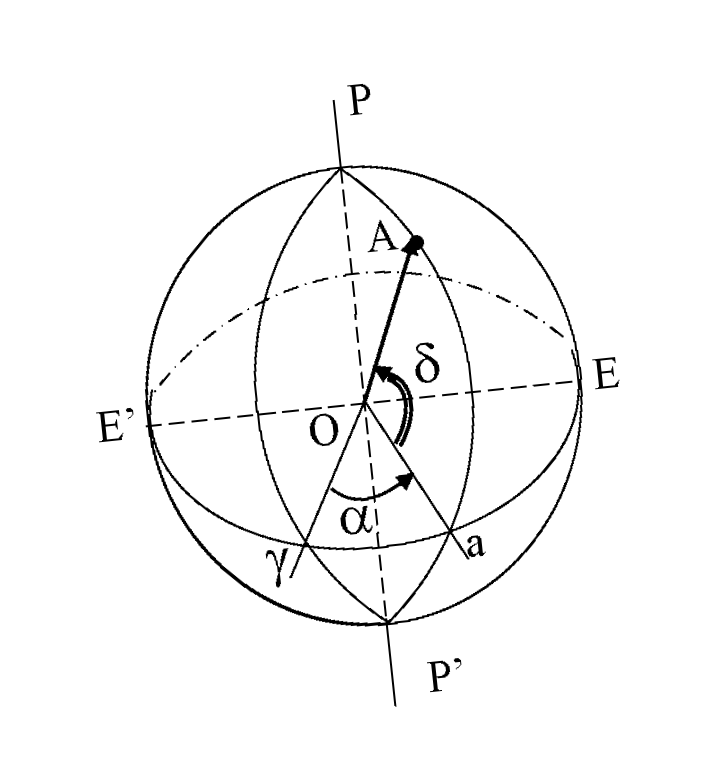}
	\caption{Les coordonnées équatoriales célestes}
	\label{fig:coordequceleste}
\end{figure}

Pour positionner le cercle horaire de l'astre A, on considère un astre fictif $\gamma$  situé sur l'équateur. On prendra comme origine le cercle horaire de $\gamma$ (\textbf{Fig. \ref{fig:coordequceleste}}).
\bdf                   
L'ascension droite $\alpha$ est l'angle entre le cercle horaire de $(\gamma)$ et le cercle horaire de l'astre, compté dans le sens direct, en heures, il mesure l'arc $\gamma a$ (\textbf{Fig. \ref{fig:coordequceleste}}).
\edf
\bdf
La déclinaison $\delta$ est l'angle du plan de l'équateur avec la direction de l'astre, compté à partir de l'équateur positivement vers P, négativement vers P'. $\delta$  mesure $aA$ (\textbf{Fig. \ref{fig:coordequceleste}}).
\edf

($\alpha$,$\delta$) constituent \textbf{les coordonnées équatoriales  célestes}. Elles sont indépendantes du temps.

\begin{figure}
	\centering
		\includegraphics[width=0.60\textwidth]{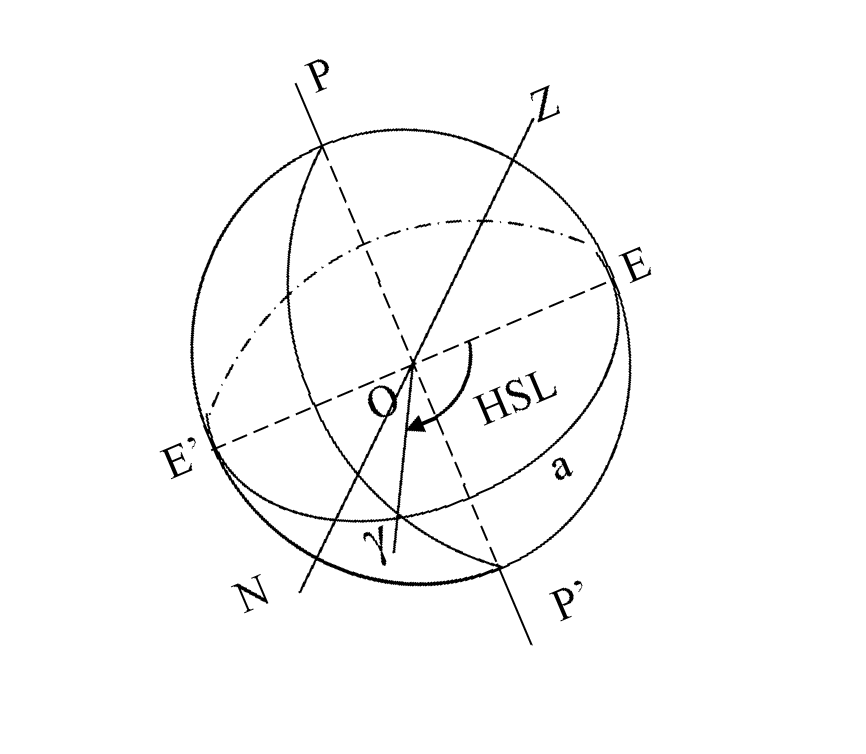}
	\caption{Heure sidérale locale}
	\label{fig:hsl}
\end{figure}

\textbf{Heure Sidérale Locale (HSL)}\index{Heure sidérale locale}
\bdf
On appelle $HSL$ en un lieu donné et à un instant donné l'angle horaire de l'astre $\gamma $. C'est l'angle formé par le méridien supérieur et le cercle horaire de $\gamma$ (\textbf{Fig. \ref{fig:hsl}}).
\edf
\be
HSL=\widehat{EO\gamma}
\ee
\bdf
Un jour sidéral: c'est l'intervalle de temps qui sépare 2 culminations successives du point $\gamma$. Il se divise en 24 heures sidérales.
\edf
\newpage

\textbf{Relation Fondamentale de l'Astronomie de Position} \index{Relation fondamentale de l'astronomie de position}
 
\begin{figure}
	\centering
		\includegraphics[width=0.60\textwidth]{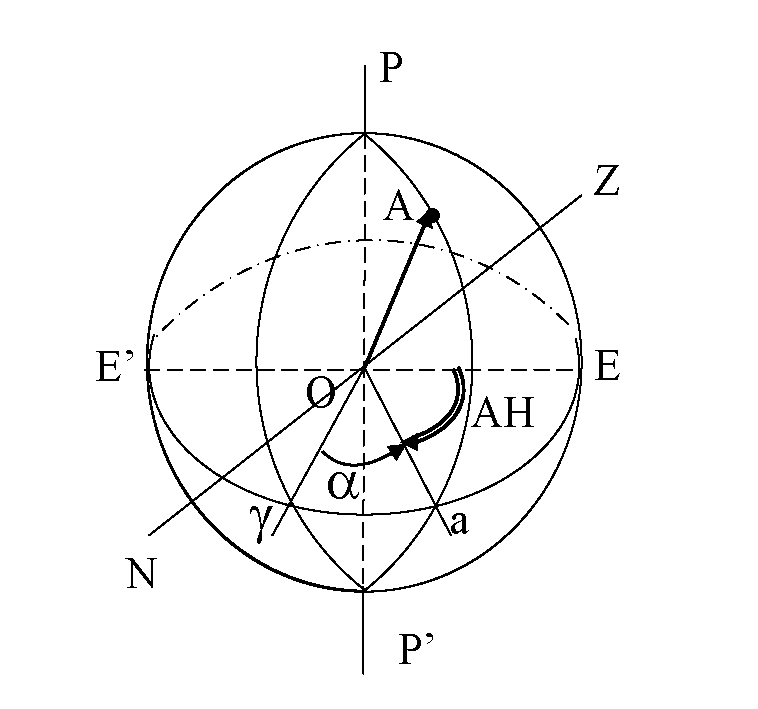}
	\caption{Relation fondamentale de l'astronomie de position}
	\label{fig:rfastronomie}
\end{figure}

A partir de la figure (\textbf{\ref{fig:rfastronomie}}), on a:
$$ \widehat{EO\gamma}=\widehat{EOa}+\widehat{aO\gamma} $$
Or:
\ba
	\widehat{EO\gamma}=HSL \nonumber  \\
	\widehat{EOa}= AH  \nonumber  
	\ea
Comme $	\widehat{aO\gamma}\geq 0$, on a donc:
$$ \widehat{aO\gamma}=\alpha $$
Par suite, on obtient la relation fondamentale de l'astronomie de position:
	\be 
 \fbox{$  HSL=AH+ \alpha $} \label{p4}
	\ee 

Au moment de \textbf{la culmination}, on a:\index{Culmination} 
\be
\fbox{ $ AH=0\quad \mbox{et} \quad HSL = \alpha $}
\ee  
\subsection*{3.4.7. Calcul de l'heure sidérale locale}
\begin{figure}
	\centering
		\includegraphics[width=0.60\textwidth]{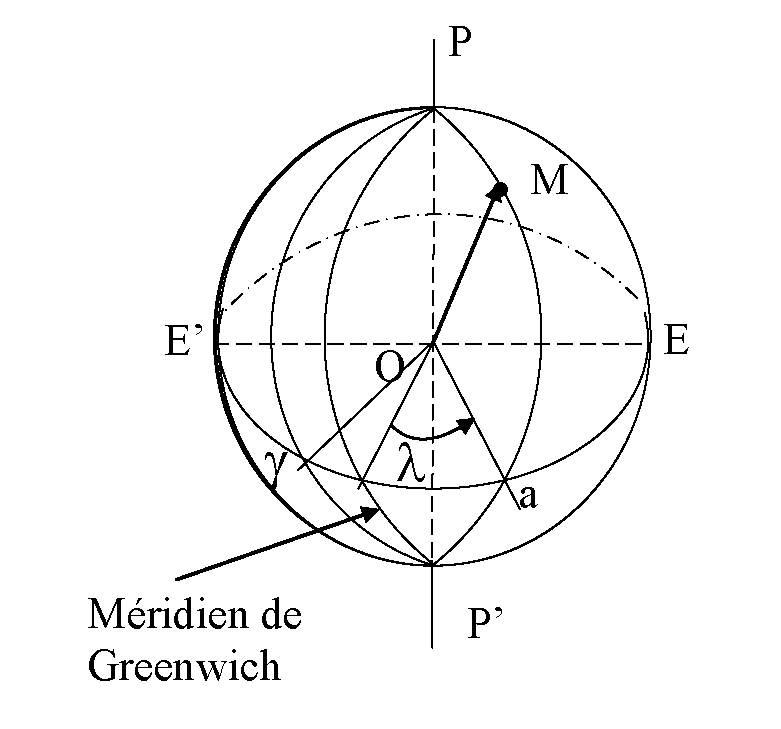}
	\caption{Relation entre $HSL_M$ et $HSG$}
	\label{fig:sphere1a}
\end{figure}

Soit un point $M$ de la Terre de longitude $\lambda$. Soit $HSL_M$ l'heure sidérale locale du lieu de $M$. Si on fait intervenir l'heure sidérale locale de Greenwich qu'on note $HSG$, on a la relation (\textbf{Fig. \ref{fig:sphere1a}}):
\be
	\fbox{ $ HSL_M=HSG+\lambda $}
\ee

\subsection*{3.4.8. Les Principales échelles de temps}
Parmi les variables primordiales de l'astronomie de position figure la variable temps. Le temps définit une échelle continue à partir d'une origine qu'on définit par convention. La mesure de l'unité de la variable temps varie d'une définition à une autre. On présente ci-après les principales échelles de temps.

\subsubsection*{Le Temps Atomique International $(TAI)$:} C'est l'étalon de temps fourni par les horloges atomiques (temps uniforme par définition). Il n'est plus rattaché à un mouvement de rotation uniforme, mais plutôt à la période de radiation de l'atome de césium 133. La seconde: c'est l'unité du Système International dont voici la définition:

\bdf
La seconde est la durée de $9\,192\,631\,770$ périodes de la radiation correspondant à la transition entre les niveaux hyperfins de l'état fondamental de l'atome de césium 133. 
\edf
\textbf{Temps Terrestre $(TT)$:} Il découle du mouvement de révolution de la Terre (temps quasi-uniforme), il s'est substitué depuis 1991 au Temps Dynamique Terrestre $(TDT)$. $TT$ assure la continuité en 1984 avec le Temps des Ephémérides $(ET)$. On a la relation suivante à la précision de la milliseconde :
\be
	\fbox{ $  TT = TDT = ET = TAI + 32,184\,s  $}
\ee
\textbf{Le Temps Universel $(UT)$:} On définit : 

*	$UT1$ (Temps Universel "Primaire") découle du mouvement de rotation de la Terre autour de son axe instantané de rotation (temps non uniforme connu a posteriori);
 
*	$UT$C (Temps Universel Coordonné) est une approximation linéaire continue par morceaux de $UT1$ dont il s'écarte au maximum de $0,9\,s$ (temps uniforme par paliers);
 
*	$UTC$ diffère de $TAI$ d'un nombre entier de secondes, ce jour\footnote{En vigueur à compter du 1er juillet 2012 à $0\,h\,\,UTC$ jusqu'à nouvel avis (Cf. bulletin C 43 de l'IERS). } : 
\be
	\fbox{ $ UTC = TAI - 35,000\,s $} 
\ee

*	A la précision considérée $(1\,s)$, $UT$ (Temps Universel) désigne indifféremment $UT1$ ou $UTC$. 
\subsection*{3.4.9. Le passage du  temps $TU$ au temps sidéral et vice versa}
Notons $TS$ et $TU$ respectivement le temps sidéral et le temps $TU,$ on a les formules suivantes (\textit{A. Danjon}, 1980):\index{\textbf{Danjon A.}}
\be
	\fbox{ $ \begin{array}{l}
TU=\ds TS\left(1-\frac{1}{366.2422}\right)= TS\times\frac{365.2422}{366.2422}\\
\\
TS=\ds TU\left(1+\frac{1}{365.2422}\right)=TU\times\frac{366.2422}{365.2422}
\end{array} $}
\ee
Application: calcul de $HSL$ à une certaine heure $TU,$ à un lieu de longitude $\lambda$, on a:
\be
\fbox{$ HSL=HSG_{0hTU}+TU\times \ds \frac{366.2422}{365.2422}+\lambda $}
\ee

\section{\textsc{Exercices et Problèmes}}
 \bex
 Au lieu de latitude $ \varphi= 36^{\circ} 54'$ Nord, on veut calculer les hauteurs $h_1$ et $h_2$ de l'étoile polaire de déclinaison $\delta= + 89^{\circ}$ respectivement à son passage supérieur et à son passage inférieur au méridien du lieu. Déterminer $h_1$ et $h_2$. 
\eex
%
 %
 %
 %
\bpb
1. En un lieu de latitude $\varphi$ quelles sont les étoiles : 

- qui ne se couchent pas ( qui sont toujours visibles);

- qui ne sont jamais visibles. 

Traiter le cas : lieu dans l'hémisphère nord. 

2. Quelle est la condition pour qu'une étoile culmine au zénith ? 

3. Cas particulier du soleil: la déclinaison du soleil varie de $-23^{\circ}27'$ à $+23^{\circ}27'$ au cours de l'année. On appelle jour le moment pendant lequel le soleil est au-dessus de l'horizon, nuit lorsque le soleil est au-dessous de l'horizon, midi l'instant de la culmination, minuit l'instant du passage inférieur. 

     a) Montrer qu'au moment des équinoxes le jour et la nuit sont d'égale durée quel que soit le lieu. 

     b) Montrer qu'à l'équateur, quelle que soit la date le jour et la nuit sont d'égale durée. 

     c) Au moment du solstice d'hiver quels sont les lieux :

 - où il fait constamment jour;

 - où il fait constamment nuit. 

Mêmes questions au moment du solstice d'été. 
   
		d) Quels sont les lieux de la Terre où le soleil culmine au zénith au moment du solstice d'hiver. Même question au moment du solstice d'été. 
    
		e) Quels sont les lieux de la Terre où au moins une fois dans l'année le soleil culmine au zénith. 
\epb
\bpb
Une station astronomique est située en un lieu de coordonnées géographiques : $\varphi =\, + 45^{\circ}\,00';\,\, \lambda =\, + 7\, h\, 20\, mn$.

En ce lieu, on observe une étoile $A$ de coordonnées équatoriales: 
 
$\alpha = \,+11\, h\, 13\, mn;\,\,\delta =\, 30^{\circ}\,00'.$

L'observation se fait le jour de l'équinoxe de printemps le 21 mars à $0$ heure $TU.$ L'heure sidérale de Greenwich est $11\, h\, 52\, mn$. 

1. Calculer l'heure sidérale locale du lever et du coucher de l'étoile $A$ au lieu considéré.
 
2. En déduire l'heure $TU$ du lever et du coucher de l'étoile au lieu considéré. 

Remarque: on choisira le coucher qui a lieu après le lever. 
\epb
\bpb
 En un lieu de latitude $43^{\circ},521$ et de longitude $+ 0\,h\,20\,mn\,57\,s$, on cherche à pointer la galaxie d'Andromède de coordonnées équatoriales $\alpha = 0\,h\,40\,mn,\, \delta = 41^{\circ}\,00'$ le 31 juillet 1992 à $21\, h\, TU.$ 

On donne l'heure sidérale de Greenwich à $0\,h\,TU$ le 31/07/1992: $HSG_{0hTU} = 20\,h\,35\,mn\,28\,s.$ 

1. Calculer l'heure sidérale locale à $21\, h\, TU.$

2. En déduire l'angle horaire de la galaxie.

3. Calculer la distance zénithale de la galaxie à $21\, h\, TU.$

4. Calculer son azimut à cette même heure.
\epb
\bpb
En un lieu de l'hémisphère Nord de latitude $\varphi$, on mesure la longueur de l'ombre portée $HC$, à midi vrai (passage du soleil au méridien), par une tige verticale $HA$ dont l'extrémité $H$ est sur le sol supposé horizontal. 

1. Donner l'expression $HC$ en fonction de $HA$ et de la distance zénithale $Dz$ du soleil. 

2. Donner l'expression de $HC$ en fonction de $HA$ et de $\varphi$: 

- aux équinoxes;

- aux solstices.
	
3. Quelle doit être la déclinaison du soleil et en quels lieux, pour que l'on ait $HC = HA$? 

4. En un lieu de latitude $\varphi = 47^{\circ}$ en quelles saisons peut on avoir $HC = HA$.
 
5. Si on déplace $HA$ le long d'un méridien, en restant dans l'hémisphère Nord, existe-t- il au cours de l'année des lieux où $HC= 0$, ou $HC$ devient infiniment grand.
\epb

\chapter{\textit{\textbf{Courbes et Surfaces}}}
\begin{svgraybox}
\begin{quotation}
	He who understands geometry may understand anything in this world.
	\end{quotation}
	\end{svgraybox}
		\begin{flushright}
	\textbf{Galileo Galilée} (1564 - 1642)
	\end{flushright}
\section{\textsc{Courbes Planes - Courbure}}
\bdf
Une courbe plane ($\mu$) est une application de $\BbR \Longrightarrow \BbR^2$  entièrement déterminée par la donnée d'une fonction vectorielle $\textbf{\textit{M}}(t)$ d'un paramètre $t$ :

$t\in \BbR \longrightarrow (x,y)\in  \BbR^2   /   \textbf{\textit{OM}}(x,y) =  x(t).\textbf{\textit{i}} + y(t).\textbf{\textit{j}} $

où $(\textbf{\textit{i}},\textbf{\textit{j}})$ la base orthonormée du plan $XOY$.
\edf
\subsection*{4.1.1. Longueur d'un arc de la courbe}
L'élément élémentaire de longueur d'un arc est la quantité $ds$ telle que: 
$$ 	ds^2 = dx^2 + dy^2 = (x'^2 + y'^2).dt^2 $$
avec $x'$ et $y'$ désignent les dérivées de $x(t)$ et $y(t)$ par rapport à la variable $t$, d'où :
$$ 	ds =  \sqrt{(x'^2 + y'^2)}.dt $$
Soit pour $t = t_0$, $M_0$  le point origine de l'arc, d'où en intégrant $s$, on obtient :
\be
	s=\int^{t}_{t_0}\sqrt{(x'^2 + y'^2)}.dt=F(t,t_0) \label{tu11}
\ee
De l'équation (\ref{tu11}), on peut exprimer $t$ en fonction de $s$. On peut alors adopter comme paramètre la longueur d'un arc de $(\mu)$ d'origine $M_0$  c'est-à-dire $s$ (l'abscisse curviligne) et de considérer la courbe définie par $M(s)$.

\begin{figure} [htp]
\centering
\includegraphics{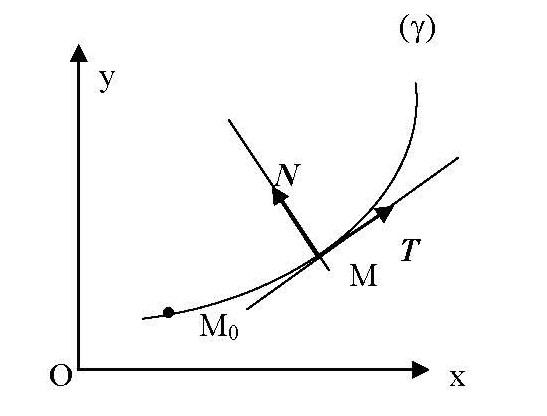}\\
\caption{Courbe plane}
\label{fig:cours3a20}
\end{figure}

\subsection*{4.1.2. La Tangente}
     Au point $M$, la courbe admet une tangente définie par le vecteur unitaire \textbf{\textit{T}}:
      \be
	\textit{\textbf{T}}=\frac{d\textit{\textbf{M}}}{ds}=\begin{pmatrix}{
\displaystyle {	\frac{dx}{ds}} \cr
\\ \cr
\displaystyle {	\frac{dy}{ds}}}
\end{pmatrix}  \label{tu12}
\ee
\subsection*{4.1.3. Normale et Courbure}
\bdf
La dérivée  de $\textbf{\textit{T}}$ par rapport à $s$  (lorsqu'elle existe et n'est pas nulle) définit une direction orthogonale à la tangente portant le vecteur unitaire $\textbf{\textit{N}}$ dite la normale au point $M$.
On a alors :
      \be
	\textit{\textbf{N}}=\frac{1}{\alpha}\frac{d \textit{\textbf{T}}}{ds} \label{tu13}
\ee
avec: 
\be
	\fbox{ $\alpha=\left\|\ds \frac{d \textit{\textbf{T}}}{ds}\right\|=\ds \frac{1}{R} $} \label{tu14}
\ee
                                               
$R$ est appelé rayon de courbure \index{Rayon de courbure} au point $M$.
\edf

\section{\textsc{Courbes Gauches}}

\subsection*{4.2.1. Trièdre de Frenêt\footnote{\textbf{Jean Frédéric Frenêt} (1816-1900): mathématicien, astronome et météorologue français.}-Courbure-Torsion}
\index{Trièdre de Frenêt}\index{\textbf{Frenêt J.F.}}
\bdf
Une courbe gauche $(\mu)$ est une application de $\BbR \Longrightarrow \BbR^3$  entièrement déterminée par la donnée d'une fonction vectorielle $\textbf{\textit{OM}(t)}$ d'un paramètre $t$ :
\be
	 t \in \BbR  \Longrightarrow  (x,y,z)\in  \BbR^3 \ /   \textbf{\textit{OM}}(x,y,z) =\begin{pmatrix}{
		x=x(t) \cr y=y(t) \cr z=z(t) } 
	\end{pmatrix}   \label{tu15}                       
\ee
\edf
\subsection*{4.2.2. Longueur d'un arc de la courbe}
L'élément élémentaire de longueur d'un arc est la quantité $ds$ telle que 
$$    ds^2 = dx^2 + dy^2 + dz^2 = (x'^2 + y'^2 +z'^2).dt^2    $$          
avec $x', y'$ et $z'$ désignent les dérivées de $x(t), y(t)$ et $z(t)$ par rapport à la variable $t$, d'où :
$$	s = \int^{t}_{t_{0}} \sqrt{x'^2 + y'^2 +z'^2}dt $$
Soit pour $t = t_0$, $M_0$ le point origine de l'arc, d'où en intégrant $s$, on obtient :
\be          
                s = F(t, t_0)            \label{tu18}                  
\ee          
De l'équation (\ref{tu18}), on peut exprimer $t$ en fonction de $s$. On peut alors adopter comme paramètre la longueur d'un arc de $(\mu)$ d'origine $M_0$  c'est-à-dire $s$ (l'abscisse curviligne) et de considérer la courbe gauche définie par $M(s)$.
\subsection*{4.2.3. La Tangente}
Au point $M$, la courbe admet une tangente définie par le vecteur unitaire \textbf{\textit{T}}. 
\be
	\textit{\textbf{T}}=\frac{d\textit{\textbf{M}}}{ds}=\begin{pmatrix}{
\displaystyle {	\frac{dx}{ds}} \cr
\\ \cr
\displaystyle {	 \frac{dy}{ds}}\cr
	\\ \cr 
	\displaystyle {  \frac{dz}{ds}}}
\end{pmatrix}  \label{tu19}
\ee                                                                         
\subsection*{4.2.4. La Normale - Courbure}
\bdf
La dérivée  de $\textbf{\textit{T}}$ par rapport à $s$, lorsqu'elle existe et n'est pas nulle, définit une direction orthogonale à la tangente portant le vecteur unitaire $\textbf{\textit{N}}$ dite la normale au point $M$. On a alors :
\be
	\textit{\textbf{N}}=\frac{1}{\alpha}\frac{\textit{d\textbf{T}}}{ds}  \label{tu20}
\ee
\newpage
avec: 
\be
	\fbox{ $\alpha=\left\|\ds \frac{ \textit{d\textbf{T}}}{ds}\right\|=\ds \frac{1}{R} $}\label{tu21}
\ee
$R$ est appelé rayon de courbure \index{Rayon de courbure}.
\edf
En effet, $\|\textbf{\textit{T}}\|=1   \Rightarrow \textbf{\textit{T}}.\textbf{\textit{T}} =1 \Rightarrow  \displaystyle 2\textbf{\textit{T}}. \frac{\textit{d\textbf{T}}}{ds} =0 $.
Donc: le vecteur $\textbf{\textit{T}}$ est orthogonal à $ \displaystyle \frac{\textit{d\textbf{T}}}{ds}$.
\begin{figure} [htp]
\centering
\includegraphics{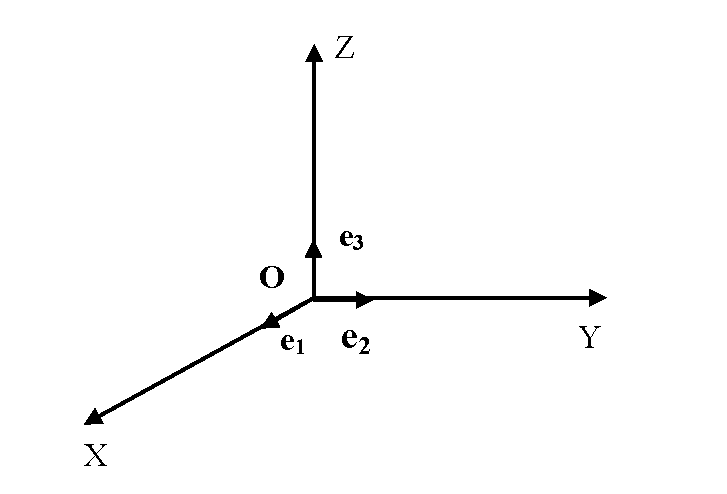}\\
\caption{Le trièdre de Frenêt}
\label{fig:repere}
\end{figure}
\subsection*{4.2.5. Binormale}
\bdf
La binormale est la droite passant par le point $M$ et de direction le vecteur $\textbf{\textit{B}}$ défini par :
\be
     \textbf{\textit{B}} = \textbf{\textit{T}} \wedge \textbf{\textit{N}}       \label{tu22} 
        \ee 
				\edf
On a évidemment :  $ \| \textbf{\textit{B}}\|$ = 1. Le triplet $(\textbf{\textit{T,N,B}})$ est direct et forme un trièdre dénommé le trièdre de Frenêt.
\bdf
Les plans définis par les vecteurs $(\textbf{\textit{T}},\textbf{\textit{N}})$, $(\textbf{\textit{N}},\textit{\textbf{B}})$ et $(\textit{\textit{B}},\textit{\textbf{T}})$ sont appelés respectivement plan osculateur, plan normal et plan rectifiant.
\edf

\subsection*{4.2.6. Torsion}
On calcule la dérivée du vecteur \textbf{\textit{B}} par rapport à $s$, on obtient: 
$$ \frac{d\textbf{\textit{B}}}{ds}=\frac{d\textbf{\textit{T}}}{ds}\wedge \textbf{\textit{N}}+\textbf{\textit{T}}\wedge\frac{d\textbf{\textit{N}}}{ds} $$
car $\textbf{\textit{T}}$ et $\displaystyle \frac{d\textbf{\textit{N}}}{ds}$ sont colinéaires, par conséquent $\displaystyle \frac{d\textbf{\textit{B}}}{ds}$ est orthogonal à $\textbf{\textit{T}}$. Comme $\textbf{\textit{B}}$ est unitaire, $\displaystyle \frac{d\textbf{\textit{B}}}{ds}$ est aussi orthogonal à $\textbf{\textit{B}}$, donc $\displaystyle \frac{d\textbf{\textit{B}}}{ds}$ est colinéaire à \textbf{\textit{N}}. On pose:
\be
	\fbox{ $ \ds \frac{d\textbf{\textit{B}}}{ds}=\frac{-1}{\tau(s)}\textbf{\textit{N}} $}  \label{tu24}
\ee
\bdf
Le réel $1/\tau(s)$ est appelé torsion\index{Torsion} de $(\mu)$ au point $M(s)$.
\edf
On calcule la dérivée du vecteur $\textbf{\textit{N}}$. Comme $\textbf{\textit{N}} = \textbf{\textit{B}}\wedge \textbf{\textit{T}}$, on obtient: 
$$ \frac{d\textbf{\textit{N}}}{ds}=\frac{d\textbf{\textit{B}}}{ds}\wedge \textbf{\textit{T}}+\textbf{\textit{B}}\wedge\frac{d\textbf{\textit{T}}}{ds}=\frac{-1}{\tau(s)}\textbf{\textit{N}}\wedge \textbf{\textit{T}}+\textbf{\textit{B}}\wedge\frac{\textbf{\textit{N}}}{R}$$
donc:
\be
\fbox{ $ \ds \frac{d\textbf{\textit{N}}}{ds}=\ds \frac{-\textbf{\textit{T}}}{R}+\frac{\textbf{\textit{B}}}{\tau(s)} $} \label{tu26}
\ee
Les trois relations exprimant les dérivées premières des vecteurs du repère de Frenêt peuvent s'écrire sous forme matricielle: 
\be
\fbox{ $
\ds \frac{d}{ds}\begin{pmatrix}{
\textbf{\textit{T}} \cr
 \textbf{\textit{N}} \cr
  \textbf{\textit{B}}}	 
\end{pmatrix}=\begin{pmatrix}{
	0     &     \ds \frac{1}{R}      &     0 \cr
\ds	\frac{-1}{R}  &     0   & \ds \frac{1}{\tau} \cr
	0  &     \ds   \frac{-1}{\tau}  &    0}
\end{pmatrix}.\begin{pmatrix}{
\textbf{\textit{T}} \cr \textbf{\textit{N}} \cr \textbf{\textit{B}} } 
\end{pmatrix} $}\label{tu27}
\ee
\section{\textsc{Surfaces}}
 \bdf
Une surface $(\sigma)$ de $\BbR^3$ est une application d'un domaine $\mathcal D \subset \BbR^2 \Rightarrow \BbR^3$ à $(u,v)\in \mathcal D \subset \BbR^2$ fait correspond un triplet $(x,y,z) \in \BbR^3$ où x,y,z sont des fonctions continues des deux paramètres $(u,v)$:
\be
	(u,v)\in \mathcal D \Rightarrow (x,y,z) \in \BbR^3/ \textbf{\textit{OM}}(u,v)=\begin{pmatrix}{
	x=x(u,v)\cr y=y(u,v)\cr z=z(u,v)}
\end{pmatrix}    \label{tu28}               
\ee
\edf
Donc $(u,v)\in \mathcal D \Rightarrow (x,y,z)\in(\sigma)$.

Si la fonction $\textbf{\textit{OM}}(u,v)$ est dérivable dans le domaine $\mathcal D$, on peut définir en tout point de $(\sigma)$ un plan tangent et une normale. 
\\

Soient $\textbf{\textit{M}}'_u$ et $\textbf{\textit{M}}'_v$ les deux vecteurs dérivées au point $M$ avec: 
\be
	    \frac{\partial \textbf{\textit{OM}}}{\partial u}(u,v)=\textbf{\textit{M'}}_u=\begin{pmatrix}{
\displaystyle {	\frac{\partial x}{\partial u}} \cr
\\ \cr
\displaystyle {	\frac{\partial y}{\partial u}}\cr
	\\ \cr
	\displaystyle { \frac{\partial z}{\partial u}}}
\end{pmatrix};\quad                        
	    \frac{\partial \textbf{\textit{OM}}}{\partial v}(u,v)=\textbf{\textit{M'}}_v=\begin{pmatrix}{
\displaystyle {	\frac{\partial x}{\partial v}} \cr
\\ \cr 
\displaystyle {	\frac{\partial y}{\partial v}}\cr
\\ \cr
	\displaystyle { \frac{\partial z}{\partial v}}}
\end{pmatrix}                 
\ee
Alors l'équation du plan tangent\index{Plan tangent} est définie par : 
$$ \textbf{\textit{MP}}.( \textbf{\textit{M}}'_u \wedge \textbf{\textit{M}}'_v )=0 $$
$P$ est un point courant du plan tangent. On pose: 
\be
	\textbf{\textit{n}} = \frac{\textbf{\textit{M}}'_u \wedge \textbf{\textit{M}}'_v}{\|\textbf{\textit{M}}'_u \wedge \textbf{\textit{M}}'_v \|} \label{tu31}
\ee
un vecteur unitaire porté par la normale à la surface $(\sigma)$ au point $M$.
 
 Les paramètres $(u,v)$ sont dits les \textit{coordonnées curvilignes}\index{Coordonnées curvilignes} sur la surface $(\sigma)$. Une courbe tracée sur la surface est définie par une relation $g(u,v) = 0$ ou par $u = u(t)$ ; $v = v(t)$ avec $t$ un paramètre. En particulier, les courbes $u =$ constante et $v =$ constante sont dites les \textit{courbes coordonnées}\index{Courbes coordonnées}.

\section{\textsc{La Première Forme Fondamentale }}
L'élément linéaire $ds$ sur la surface $(\sigma)$ est la distance de deux points infiniment voisins, le carré de $ds$ est le carré scalaire de $d\textbf{\textit{M}}$ soit : 
\be
	ds^2 = d\textbf{\textit{M}}.d\textbf{\textit{M}} = d\textbf{\textit{M}}^2 = \|d\textbf{\textit{M}} \|^2   \label{tu32}
\ee
\newpage
Or :
$$ 	\textbf{\textit{OM}}(u,v)\,\,\left(\begin{matrix}{
	x(u,v)\cr y(u,v)\cr z(u,v)}
\end{matrix}\right)\, \Longrightarrow d\textbf{\textit{M}}=\textbf{\textit{M}}'_udu+\textbf{\textit{M}}'_vdv=\left(\begin{matrix}{
	dx=x'_udu+x'_vdv \cr
	 dy=y'_udu+y'_vdv \cr 
	 dz=z'_udu+z'_vdv}
\end{matrix} \right) $$ 
 Par suite : 
$$ 	 d\textbf{\textit{M}}^2 = ds^2 = \textbf{\textit{M}}'_u.\textbf{\textit{M}}'_u du^2 + 2\textbf{\textit{M}}'_u.\textbf{\textit{M}}'_vdudv + \textbf{\textit{M}}'_v.\textbf{\textit{M}}'_vdv^2  $$
On pose :
\be
\left\{\begin{array}{l}
 E = \textbf{\textit{M}}'_u.\textbf{\textit{M}}'_u   \\ 
 F = \textbf{\textit{M}}'_u.\textbf{\textit{M}}'_v    \\ 
 G = \textbf{\textit{M}}'_v.\textbf{\textit{M}}'_v  
\end{array}\right. 
\ee
alors $ds^2$ s'écrit:
\be
  \fbox{ $ds^2 = E.du^2 + 2.F.dudv + G.dv^2  $}\label{tu38} 
\ee
(\ref{tu38}) est dite la première forme fondamentale,\index{Première forme fondamentale} elle définit la métrique de la surface $(\sigma)$. 

\subsection*{4.4.1. Ecriture matricielle de la première forme fondamentale}
On appelle $g=(g_{ij})$ la matrice carrée $2\times2$ telle que:
\ba
g_{11}=E \nonumber \\
g_{12}=g_{21}=F\nonumber \\
g_{22}=G \nonumber
\ea
Soit:
\be
g=\begin{pmatrix}{
g_{11} & g_{12} \cr
g_{21} & g_{22} }
\end{pmatrix}=\begin{pmatrix}{
E & F \cr
F & G}
\end{pmatrix} \lb{tu38a}
\ee
Alors l'équation (\ref{tu38}) s'écrit sous la forme:
\be
ds^2=(du,dv).g.\begin{pmatrix}{
du \cr
dv}
\end{pmatrix}=(du,dv).\begin{pmatrix}{
E & F \cr
F & G}
\end{pmatrix}.\begin{pmatrix}{
du \cr
dv}
\end{pmatrix}
\ee
La matrice $g$ s'appelle la matrice du tenseur métrique.\index{Tenseur métrique}
\subsection*{4.4.2. Angles de deux courbes coordonnées et Elément d'aire}
* On a : 
$F = \textbf{\textit{M}}'_u.\textbf{\textit{M}}'_v= \|\textbf{\textit{M}}'_u\|.\|\textbf{\textit{M}}'_v\|cos\alpha$, d'où:
$$ cos\alpha = \frac{F}{\|\textbf{\textit{M}}'_u\|.\|\textbf{\textit{M}}'_v\|}=\frac{F}{\sqrt{E}\sqrt{G}}=\frac{F}{\sqrt{EG}}  $$
et en considérant $\alpha \in [0,\pi]$:
$$ sin\alpha =\sqrt{1-cos\alpha^2}=\sqrt{1-\frac{F^2}{E.G}}=\sqrt{\frac{E.G-F^2}{E.G}}  $$
On pose parfois:
\be
	H=\sqrt{E.G-F^2}=h^2(u,v) \label{tu41}
\ee
soit:
\be
	H = \|\textbf{\textit{M}}'_u\wedge \textbf{\textit{M}}'_v\|=\|\textbf{\textit{M}}'_u\|.\|\textbf{\textit{M}}'_v\|sin\alpha  \label{tu42}
\ee
Par suite, le vecteur unitaire normal $\textbf{\textit{n}}$ a pour expression : 
\be
	\textbf{\textit{n}} = \frac{\textbf{\textit{M}}'_u\wedge\textbf{\textit{M}}'_v}{\|\textbf{\textit{M}}'_u\wedge \textbf{\textit{M}}'_v\|}= 
\frac{\textbf{\textit{M}}'_u\wedge\textbf{\textit{M}}'_v}{H} \label{tu43}
\ee
* En considérant maintenant le parallélogramme curviligne de sommet $M(u,v)$ et de côtés les vecteurs $\textbf{\textit{M}}'_udu$ et $\textbf{\textit{M}}'_vdv$, alors l'élément infinitésimal d'aire $d\m A$ a pour expression:
$$d\m A=||\textbf{\textit{M}}'_udu\wedge\textbf{\textit{M}}'_vdv||=||\textbf{\textit{M}}'_u||.||\textbf{\textit{M}}'_v||du.dv.sin\alpha=\sqrt{E.G-F^2}dudv=Hdudv $$
On le note aussi:
\be
\fbox{ $ d\m A=\ds \sqrt{E.G-F^2}du\wedge dv=Hdu\wedge dv $}
\ee
\subsection*{4.4.3. Coordonnées Orthogonales et Coordonnées Symétriques}
Les coordonnées $(u,v)$ sont dites \textbf{orthogonales} si $F = \textbf{\textit{M}}'_u .\textbf{\textit{M}}'_v =0$, soit $cos\alpha = 0$, donc $\alpha$ est un angle droit.

Les coordonnées orthogonales sont dites \textbf{coordonnées symétriques} si de plus $E = G$. Alors la première forme quadratique s'écrit: 
$$	ds^2 = Edu^2 + Gdv^2 = E(du^2 + dv^2) = H(du^2 + dv^2) = h^2(u,v)(du^2 + dv^2) $$
\textbf{Exemple:}

On considère une sphère de rayon $R$ qu'on note $(\sigma)$, elle est paramétrée par:
$$ \textbf{\textit{OM}}\left|\begin{array}{l}
Rcos\varphi cos\lambda \\
Rcos\varphi sin \lambda \\
Rsin\varphi 
\end{array} \right.$$
avec $\varphi \in [-\frac{\pi}{2},+\frac{\pi}{2}], \lambda \in [0,+2\pi[$. Les courbes coordonnées de $(\sigma)$ sont les méridiens $\lambda=constante$ et les parallèles $\varphi=constante$. On remarque qu'elles se coupent en un angle droit. On calcule la première forme fondamentale de la sphère:
\be
\textbf{\textit{OM}}'_{\varphi}=\left|\begin{array}{l} 
-Rsin\varphi cos\lambda \\
-Rsin\varphi sin\lambda \\
Rcos\varphi 
\end{array}\right.,\quad \textbf{\textit{OM}}'_{\lambda}=\left|\begin{array}{l} 
-Rcos\varphi sin\lambda \\
Rcos\varphi cos\lambda \\
0 
\end{array}\right. 
\ee
D'où:
\ba
E=\textbf{\textit{OM}}'_{\varphi}.\textbf{\textit{OM}}'_{\varphi}=R^2 \nonumber \\
F=\textbf{\textit{OM}}'_{\varphi}.\textbf{\textit{OM}}'_{\lambda}=0 \nonumber \\
G=\textbf{\textit{OM}}'_{\lambda}.\textbf{\textit{OM}}'_{\lambda}=R^2cos^2\varphi \nonumber
\ea
$F=0$ justifie ce qui a été dit ci-dessus sur l'orthogonalité des courbes coordonnées. Ces dernières sont orthogonales mais non symétriques. En effet : 
$$ 	ds^2 = R^2d \varphi^2+R^2cos^2\varphi d\lambda^2=R^2 cos^2\varphi \left(\frac{d\varphi^2}{cos^2\varphi}+d\lambda^2\right)   $$
La variable $L$ telle que : 
\be
	  dL  = \ds \frac{d\varphi}{cos\varphi} \lb{merc}
\ee
forme avec $\lambda$ un couple de coordonnées symétriques, car : 
\be
	ds^2 = R^2 cos^2\varphi(dL^2 + d\lambda^2)
\ee
$L$ est appelée \textbf{latitude croissante}\index{\textbf{Mercator G.}}\index{Latitude croissante} ou \textbf{latitude ou variable de Marcator}\footnote{\textbf{Gerhardus Mercator} (1512-1594): cartographe, astronome et ingénieur belge. Son nom était donné à la représentation cylindrique conforme proposée par lui-même. }\index{Variable de Mercator}.\index{Latitude de Mercator} On pose:
$$ t=tg\frac{\varphi}{2}\Longrightarrow cos\varphi=\frac{1-t^2}{1+t^2}$$
d'où:
$$ dt=(1+t^2)\frac{d\varphi}{2} \Longrightarrow d\varphi=\frac{2dt}{1+t^2}$$
De (\ref{merc}), on obtient $L$ vérifiant $L(0)=0$:
\be
L=\int_0^{\varphi}\frac{dw}{cosw}=\int_0^{tg\frac{\varphi}{2}}\frac{2dt}{1+t^2}.\frac{1+t^2}{1-t^2}=\int_0^{tg\frac{\varphi}{2}}\frac{2dt}{1-t^2}\lb{merc1}
\ee
On se restreint à $\varphi \in [-\varphi_1,+\varphi_1]$ où $0<\varphi_1<\pi/2$. L'équation (\ref{merc1}) s'écrit:
\ba
L= \ds \int_0^{\varphi}\frac{dw}{cosw}=\ds \int_0^{tg\frac{\varphi}{2}}\frac{dt}{1+t}+\int_0^{tg\frac{\varphi}{2}}\frac{dt}{1-t}=\ds \left[Log\left|\ds \frac{1+t}{1-t}\right|\right]_0^{tg\frac{\varphi}{2}}=\nonumber \\
Log\left|tg\left(\ds \frac{\pi}{4}+\frac{\varphi}{2}\right)\right|=Log tg\left(\ds \frac{\pi}{4}+\frac{\varphi}{2}\right)
\ea
car $tg\left(\ds \frac{\pi}{4}+\frac{\varphi}{2}\right) >0$, donc:
\be
	\fbox{ $ L = Log tg \left(\ds \frac{\varphi}{2}+\frac{\pi}{4}\right) $}
\ee
D'où :
$$ tg\left(\frac{\varphi}{2}+\frac{\pi}{4}\right)=exp(L)=\ds e^{\ds L}$$

Soit l'expression de la latitude $\varphi$ en fonction de $L$: 
$$ 	\varphi= \ds 2Arctg(e^{\ds L})-\frac{\pi}{2}  $$
\section{\textsc{La Deuxième Forme Fondamentale}}\index{Deuxième forme fondamentale}
  On calcule maintenant le vecteur $d^2\textbf{\textit{M}} $ la différentielle seconde de $\textit{\textbf{OM}}$. On a alors :
\ba
	   d^2\textbf{\textit{M}}= d(d\textbf{\textit{M}})=d(\textbf{\textit{OM}}'_u.du+\textbf{\textit{OM}}'_v.dv)= d\textbf{\textit{OM}}'_u.du +d\textbf{\textit{OM}}'_v.dv \nonumber \\	  
	          = (\textbf{\textit{OM}}''_{uu} .du + \textbf{\textit{OM}}''_{uv} .dv).du + (\textbf{\textit{OM}}''_{uv} .du + \textbf{\textit{OM}}''_{vv} .dv)dv \nonumber
\ea
soit : 
$$ d^2\textbf{\textit{M}}= \textbf{\textit{OM}}''_{uu} .(du)^2 + 2. \textbf{\textit{OM}}''_{uv} .du.dv + \textbf{\textit{OM}}''_{vv} .(dv)^2 $$
car $\textbf{\textit{OM}}''_{uv}  = \textbf{\textit{OM}}''_{vu}$. On peut écrire l'équation précédente sous la forme :
\be
	d^2\textbf{\textit{M}}=\frac{\partial ^2 \textbf{\textit{M}}}{\partial u^2}du^2+2\frac{\partial ^2 \textbf{\textit{M}}}{\partial u \partial v}dudv+\frac{\partial ^2 \textbf{\textit{M}}}{\partial v^2}dv^2 \label{dffeq1}
\ee
Soit $(\gamma)$ une courbe tracée sur la surface $(\sigma)$, définie par $u=u(s), v= v(s)$ où $s$ désigne l'abscisse curviligne. Soit $\textbf{\textit{n}}$ le vecteur normal à la surface et $\textbf{\textit{N}}$ le vecteur unitaire porté par la normale principale à la courbe $(\gamma)$. Si $\textbf{\textit{T}}$ est le vecteur porté par la  tangente à $(\gamma)$ au point $M(u,v)$, d'après les formules de Frenêt, on a :
$$ \ds 	\frac{d\textbf{\textit{T}}}{ds}=\frac{\textbf{\textit{N}}}{R} $$
où $\ds \frac{1}{R}$ est \textit{la  courbure} de $(\gamma)$ au point $M$. Or $d\textbf{\textit{M}} = \textbf{\textit{T}}ds$, par suite :
\be
	d(d\textbf{\textit{M}})=d^2\textbf{\textit{M}}=d(\textbf{\textit{T}}ds)=ds.d\textbf{\textit{T}}=ds.\left(\textbf{\textit{N}}\frac{ds}{R}\right) \label{dffeq2}
\ee
On multiplie vectoriellement l'équation (\ref{dffeq1}) par le vecteur normal unitaire $\textbf{\textit{n}}$, on obtient :
\be
	\textbf{\textit{n}}.d^2\textbf{\textit{M}}=\textbf{\textit{n}}.\frac{\partial ^2 \textbf{\textit{M}}}{\partial u^2}du^2+2\textbf{\textit{n}}.\frac{\partial ^2 \textbf{\textit{M}}}{\partial u \partial v}dudv+\textbf{\textit{n}}.\frac{\partial ^2 \textbf{\textit{M}}}{\partial v^2}dv^2 \label{dffeq3}
\ee
On pose :
\be
	\fbox{ $ \left\{\begin{array}{lll}
	L=\displaystyle \textbf{\textit{n}}.\frac{\partial ^2 \textbf{\textit{M}}}{\partial u^2} \\
	\\
	M=\displaystyle \textbf{\textit{n}}.\frac{\partial ^2 \textbf{\textit{M}}}{\partial u \partial v} \\
	\\
	N=\displaystyle \textbf{\textit{n}}.\frac{\partial ^2 \textbf{\textit{M}}}{\partial v^2}
	\end{array}\right. $}
\ee
(\ref{dffeq3}) s'écrit alors :
\be
	\textbf{\textit{n}}.d^2\textbf{\textit{M}}=Ldu^2+2Mdudv+Ndv^2 	\label{dffeq4}
\ee
On multiplie aussi l'équation (\ref{dffeq2}) par le vecteur $\textbf{\textit{n}}$ d'où :
$$  \textbf{\textit{n}}.d^2\textbf{\textit{M}}=\textbf{\textit{n}}.ds.\left(\textbf{\textit{N}}\frac{ds}{R}\right)=\textbf{\textit{n.N}}\frac{ds^2}{R}$$
Soit $\theta$  l'angle formé par $\textbf{\textit{n}}$ et $\textbf{\textit{N}}$, d'où :
\be
\textbf{\textit{n}}.d^2\textbf{\textit{M}}=\textbf{\textit{n}}.\textbf{\textit{N}}\frac{ds^2}{R}=\frac{cos\theta}{R}ds^2\label{dffeq5}
\ee
Comme (\ref{dffeq4}) est égal à (\ref{dffeq5}), on obtient :             
$$ 	\frac{cos\theta}{R}ds^2=Ldu^2+2Mdudv+Ndv^2=\Phi(u,v) $$
soit :         
$$ \ds 	\frac{cos\theta}{R}=\frac{Ldu^2+2Mdudv+Ndv^2}{ds^2} $$
ou encore :                                                       
   \be
	\fbox{ $ \ds \frac{cos\theta}{R}=\frac{Ldu^2+2Mdudv+Ndv^2}{Edu^2+2Fdudv+Gdv^2}=\frac{II(u,v)}{I(u,v)} $} \label{dffeq10}
\ee
avec $I(u,v)$ la première forme fondamentale et l'expression :                    
\be 
	\fbox{ $  II(u,v) = Ldu^2 + 2Mdudv + Ndv^2 =  \Phi(u,v) $}   \label{dffeq1044}                               
\ee
est appelée \textbf{\textit{la deuxième forme fondamentale}}.
D'où:
\bthm
Le produit de la courbure en un point donné d'une courbe tracée sur une surface dans l'espace à trois dimensions par le cosinus de l'angle entre la normale à la surface et la normale principale à la courbe est égale au rapport de la deuxième et la première formes fondamentales du vecteur tangent à la courbe en ce point.
\ethm
\bdf
La quantité $\displaystyle \frac{cos\theta}{R}$ invariante pour toutes les courbes ayant même vecteur tangent $\textbf{\textit{T}}$ en un point donné est dite la courbure normale de la surface en ce point.\index{Courbure normale}
\edf
\bpr
Si la courbe est la section d'une surface par un plan normal, on a :
\be
	\fbox{ $ \theta =0\,\, ou\,\;\theta=\pi \Rightarrow cos\theta=\pm1 \Rightarrow \,	\ds \pm \frac{1}{R}=\frac{Ldu^2+2Mdudv+Ndv^2}{Edu^2+2Fdudv+Gdv^2} $} \label{corr}
\ee
\epr         
\subsection*{4.5.1. Trièdre de Darboux\footnote{\textbf{Jean Gaston Darboux} (1842-1917): mathématicien français.} - Ribaucour\footnote{\textbf{Albert Ribaucour} (1845-1893): ingénieur et mathématicien français.}}\index{Trièdre de Darboux - Ribaucour}\index{\textbf{Ribaucour A.}} \index{\textbf{Darboux J.G.}}Soit $(\gamma)$ une courbe tracée sur une surface $(\sigma)$ pour laquelle on sait définir en un point donné $M$ le repère de Frénet $(\textbf{\textit{T}}, \textbf{\textit{N}}, \textbf{\textit{B}})$.
\bdf 
On appelle repère de Darboux - Ribaucour $(\textbf{\textit{T}}, \textbf{\textit{n}} , \textbf{\textit{g}})$ le repère orthonormé formé par les vecteurs $\textbf{\textit{T}}$, $\textbf{\textit{n}}$ et le vecteur $\textbf{\textit{g}}=\textbf{\textit{T}}\wedge \textbf{\textit{n}}$. 
\edf
La position relative des deux repères est donnée par l'angle :  
\be
	\theta=\widehat{\textbf{\textit{N}},\textbf{\textit{n}}}
\ee
\subsection*{4.5.2. Section Normale}\index{Section normale}
\bdf
Soit la courbe $(\gamma)$ tracée sur $(\sigma)$  et définie comme intersection de  $(\sigma)$  et du plan passant par le point $M$ et de directions $\textbf{\textit{n}}$ et $\textbf{\textit{T}}$, alors $(\gamma)$   est appelée section normale de   $(\sigma)$  en $M$ dans la direction $\textbf{\textit{T}}$.
\edf
La  normale principale de $(\gamma)$  est la droite portée par le vecteur $\textbf{\textit{n}}$. Si $R_n$  est le rayon de courbure de $(\gamma)$ au point $M$, on a par définition :
$$	\frac{d\textbf{\textit{T}}}{ds}=\frac{\textbf{\textit{n}}}{R_n} $$
par suite :
$$ \textbf{\textit{n}}.\frac{d\textbf{\textit{T}}}{ds}=\frac{1}{R_n} $$ 
or (\ref{dffeq2}) donne :       
$$ 	d\textbf{\textit{T}}=\frac{d^2\textbf{\textit{M}}}{ds} $$
d'où :   
 $$ \textbf{\textit{n}}.\frac{d\textbf{\textit{T}}}{ds}=\textbf{\textit{n}}.\frac{d^2\textbf{\textit{M}}}{ds^2}=\frac{1}{R_n}=\frac{II(u,v)}{I(u,v)} $$ 
soit :  
 \be 
\frac{1}{R_n}=\frac{Ldu^2+2Mdudv+Ndv^2}{Edu^2+2Fdudv+Gdv^2}=\frac{II(u,v)}{I(u,v)}
\ee 
En comparant l'équation ci-dessus avec l'équation (\ref{dffeq10}), on obtient :                       
\be 
\fbox{ $ 	R= R_n.cos\theta $}
\ee
D'où le deuxième théorème de Meusnier\footnote{\textbf{Jean Baptiste Meusnier} (1754-1793): militaire, géomètre et mathématicien français.} : \index{Deuxième théorème de Meusnier}\index{\textbf{Meusnier J.B.}}
\bthm
Le rayon de courbure $R$ d'une courbe $(\gamma)$ tracée sur une surface $(\sigma)$ et ayant même tangente de direction $\textbf{\textit{T}}$ est égal au produit  de $R_n$ rayon de courbure de la section normale par le cosinus de l'angle $\theta$ entre les vecteurs $\textbf{\textit{n}}$ et $\textbf{\textit{N}}$. 
\ethm
\subsection*{4.5.3. Indicatrice de Dupin\footnote{\textbf{Charles Dupin} (1784-1873): ingénieur et mathématicien français.}}\index{Indicatrice de Dupin}\index{\textbf{Dupin C.}}
On considère le repère orthonormé $\mathcal R$ au point $M$ défini par les vecteurs :
$$ 	\frac{1}{\sqrt{E}}\frac{\partial \textbf{\textit{M}}}{\partial u} \quad \mbox{et}\quad  	\frac{1}{\sqrt{G}}\frac{\partial \textbf{\textit{M}}}{\partial v} $$
 \bdf
L'indicatrice de Dupin est l'ensemble des points $P$ du plan tangent en $M$ vérifiant :
\be
	\textbf{\textit{MP}}=\sqrt{R_n}\textbf{\textit{T}}
\ee
quand  $\textbf{\textit{T}}$  varie autour de $M$.
\edf
Soit un point $P(\alpha,\beta)$ dans $\mathcal R$, on a alors:
\ba
	\textbf{\textit{MP}}=	\frac{\alpha}{\sqrt{E}}\frac{\partial \textbf{\textit{M}}}{\partial u}+	\frac{\beta}{\sqrt{G}}\frac{\partial \textbf{\textit{M}}}{\partial v} \nonumber \\
		\textbf{\textit{MP}}=\sqrt{R_n}\textbf{\textit{T}}=	\sqrt{R_n}\left(	\frac{du}{ds}\frac{\partial \textbf{\textit{M}}}{\partial u}+	\frac{dv}{ds}\frac{\partial \textbf{\textit{M}}}{\partial v}\right) \nonumber
\ea
d'où :
$$ 	\frac{du}{ds}=\frac{\alpha}{\sqrt{ER_n}} \quad \mbox{et} \quad \frac{dv}{ds}=\frac{\beta}{\sqrt{GR_n}} $$
En utilisant la deuxième forme quadratique, on a :
\ba
		\frac{1}{R_n}=L\left(\frac{du}{ds}\right)^2+2M\frac{du}{ds}\frac{dv}{ds}+N\left(\frac{dv}{ds}\right)^2 \Rightarrow \nonumber \\
			\frac{1}{R_n}=L\frac{\alpha^2}{ER_n}+2M\frac{\alpha \beta}{R_n\sqrt{EG}}+N\frac{\beta^2}{GR_n} \nonumber
\ea 
ou encore :
   \be
\fbox{ $ \ds \frac{L}{E}\alpha^2+2\frac{M}{\sqrt{EG}}\alpha \beta +\frac{N}{G}\beta^2=1 $}
\ee                                                                                                      
C'est l'équation d'une conique (ellipse, parabole, hyperbole) suivant le signe du discriminant $\displaystyle \frac{M^2-LN}{EG}$ ou $M^2-LN$ respectivement (négatif, nul ou positif).  
\subsection*{4.5.4. Les Directions principales}\index{Directions principales}
On suppose que $M^2 - LN <0$. 

\bdf
On appelle directions principales les directions des axes de symétrie de l'indicatrice de Dupin.
\edf
\bdf
On appelle les rayons  de courbure principaux $R_1$ et $R_2$ les rayons de courbure normale dans les deux directions principales.
\edf

Les directions principales sont orthogonales.
\subsection*{4.5.5. Formule d'Euler\footnote{\textbf{Leonhard Euler} (1707-1783): mathématicien et physicien suisse.}}\index{La formule d'Euler}\index{\textbf{Euler L.}}
En supposant que l'indicatrice de Dupin est une ellipse d'équation :
\be
 \ds 	\frac{\alpha^2}{a^2}+\frac{\beta^2}{b^2}=1 \lb{ind}
\ee 
où $a, b$ sont les 2 rayons de courbure normale principaux, on peut écrire :                                   
$$ \textbf{\textit{MP}}=\sqrt{R_n}\textbf{\textit{T}}=\textbf{\textit{i}}\sqrt{R_n}cos\psi+\textbf{\textit{j}}\sqrt{R_n}sin\psi $$
avec $ \psi= \widehat{\textbf{\textit{T}},\textbf{\textit{i}}}$, or:    
$ \textbf{\textit{MP}} = \alpha \textbf{\textit{i}}+\beta \textbf{\textit{j}}$, avec l'équation (\ref{ind}), on obtient alors la formule d'Euler :
$$ 	\frac{R_ncos^2\psi}{a^2}+\frac{R_nsin^2\psi}{b^2}=1 \Rightarrow \frac{1}{R_n}=	\frac{cos^2\psi}{a^2}+\frac{sin^2\psi}{b^2} \label{ind1} $$
D'où:
\bthm
(\textbf{Formule d'Euler}):La courbure de la section normale $ \displaystyle \frac{1}{R_n}$ en un point donné est égale à:
\be
\fbox{ $ \ds \frac{1}{R_n}=\frac{cos^2\psi}{a^2}+\frac{sin^2\psi}{b^2} $}
\ee
où $\displaystyle \frac{1}{a^2},\frac{1}{b^2}$ sont les courbures principales au point considéré et $\psi$ l'angle sur la surface entre le vecteur tangent à la section normale et la direction principale correspondante.
\ethm

\bdf
Le produit des courbures principales est la courbure de Gauss\footnote{\textbf{Carl Friedrich Gauss} (1777-1855): mathématicien et géomètre prussien, fondateur de la théorie des surfaces. } \index{Courbure de Gauss} ou courbure totale de la surface et la courbure moyenne \un{la somme} des courbures principales.\index{Courbure totale}\index{Courbure moyenne} 
\edf
\index{\textbf{Gauss C.F.}}
%
Pour la première forme fondamentale $ds^2$, on a déjà noté (\ref{tu38a}): 
$$g=\begin{pmatrix}{
E & F \cr
F & G}
\end{pmatrix}$$
et concernant la deuxième forme fondamentale  $\Phi(u,v)$ donnée par l'équation (\ref{dffeq1044}), elle peut s'écrire sous la forme:
 \be
\Phi=\Phi(u,v)= (du,dv).\begin{pmatrix}{
L & M \cr
M & N}
\end{pmatrix}.\begin{pmatrix}{
du \cr
dv}
\end{pmatrix}
\ee  
où par abus de notation, on a noté par $\Phi$ la matrice ci-dessus. 

Alors, on annonce les deux théorèmes suivants sans les démontrer (\textit{B. Doubrovine - S. Novikov - A. Fomenko}, 1982):\index{\textbf{Doubrovine B.}}\index{\textbf{Novikov S.}}\index{\textbf{Fomenko A.}}
\bthm
La courbure totale $K$ en un point d'une surface est égale au rapport des déterminants de ses deuxième et première formes fondamentales:
\be
\fbox{$ K= \ds \frac{\mbox{Dét} \Phi}{\mbox{Dét} g}=\frac{LN-M^2}{EG-F^2} $}
\ee
\ethm
et:
\bthm
La courbure moyenne $H$ en un point d'une surface est égale à la trace de la matrice $g^{-1}.\Phi$:
\be
\fbox{$ H=Tr(g^{-1}.\Phi) $}
\ee
\ethm
\vspace{1cm}
\textbf{\underline{Note historique:}} \textsl{La théorie des surfaces élaborée par Gauss était surtout influencée essentiellement par son travail comme géomètre topographe dans le Royaume de Hannover au Nord de l'Allemagne durant la période 1821-1825. En 1822, il présenta son mémoire intitulé " General solution of the problem of mapping parts of a given surface onto another given surface in such a way that image and pre-image become similar in their smallest parts", à la Société Royale des Sciences à Copenhague (Danemark) où il recevait un prix officiel.}\index{\textbf{Gauss C.F.}}

\textsl{Où se réside donc l'importance de son mémoire?  Ce dernier concernait l'étude du problème de cartographier une surface sur une autre en satisfaisant certaines propriétés. C'est le problème de base de la cartographie. Parmi les représentations planes dites abusivement projections sont celles qui conservent les angles ou représentations conformes. Elles ont un aspect pratique pour la navigation maritime. Ainsi, Gauss avait réussi à trouver une procédure pour déterminer toutes les représentations conformes localement pour les surfaces analytiques. Il ajouta dans le titre de son mémoire cette phrase en latin:}
\begin{center}
 \textsl{Ab his  via sterniture ad maiora}.
\end{center}
\textsl{soit} " \textbf{De là, le chemin de quelque chose plus importante est préparé} ". \textsl{En effet, Gauss présentait en octobre 1827 une théorie générale des surfaces à travers son papier " Disquisitiones generales circa superficies curvas\footnote{Voir aussi (\textbf{P. Dombrowski}, 1979).} " (Investigations about curved surfaces). L'important résultat de son papier est  \underline{le théorème egregium} dit encore \underline{le théorème merveilleux}. Ce dernier dit que la courbure de Gauss est une propriété intrinsèque pour les surfaces de dimension 2. La courbure de Gauss dépend des composantes $g_{ij}$ du tenseur métrique et de ses dérivées partielles premières et secondes par rapport aux coordonnées locales.} (\textit{E. Zeidler}, 2011)\index{\textbf{Zeidler E.}}\index{\textbf{Dombrowski P.}}
\section{\textsc{Exercices et Problèmes}}
\bex
Soit $(\Gamma)$ la surface paramétrée par $(u,v)$ dans $\BbR^2$ telle que:
\[M(u,v) \left\{ \begin{array}{l}
	     X =u(1-u^2)cosv \\
	    Y = u(1-u^2)sinv\\
	    Z = 1-u^2
	    \end{array} \right.
\]
1. Calculer l'expression de $ds^2$.

2. Montrer que l'équation cartésienne de $(\Gamma)$ est: 
$$ x^2+y^2=(1-z)z^2 $$
\eex
\bex
 Soit la surface d'Enneper\footnote{Surface paramétrée par le mathématicien allemand \textbf{Alfred Enneper} (1830-1885).}:
	\[M(u,v) \left\{ \begin{array}{l}
	     X =\ds u-\frac{u^3}{3}+uv^2 \\
	    Y = v-\ds \frac{v^3}{3}+vu^2 \\
	    Z = u^2-v^2
	    \end{array} \right.
\]
1. Montrer que:    $$ ds^2 = (1+u^2+v^2)^2.(du^2+dv^2)$$

2. Calculer un vecteur unitaire normal à la surface.

3. Montrer que la surface d'Enneper est de courbure moyenne nulle en chaque point.
\eex
\bex
 On suppose que la métrique d'une surface donnée est:
$$ds^2=A^2du^2+B^2dv^2, \quad A=A(u,v),\quad B=B(u,v)$$
1. Montrer alors que l'expression de la courbure totale est:
$$ K=-\frac{1}{AB}\left[\left(\frac{A'_v}{B}\right)'_v+\left(\frac{B'_u}{A}\right)'_u\right]$$
$'$ désigne la dérivation partielle. 
\eex  
\bpb
 Soit l'ellipse $(E)$ définie par les équations paramétriques:
$$ M \left\{\begin{array}{l}
x=acos u \\
y= bsinu \\
\mbox{avec}\quad a>b>0
\end{array}\right. $$
On  pose:
$$e^2=\frac{a^2-b^2}{a^2};\quad e'^2= \frac{a^2-b^2}{b^2}$$
1. Calculer la position sur l'axe des abscisses des deux points $F$ et $F'$ appelés foyers tels que $MF+MF'=2a$.

2. Montrer que le produit des distances des foyers à la tangente à l'ellipse en M est indépendant de $u$.

3. Donner l'expression de $ds$.

4. Déterminer les expressions des vecteurs unitaires $\textbf{\textit{T}}$ et $\textbf{\textit{N}}$ et en déduire le rayon de coubure de l'ellipse.  

5. Montrer qu'il passe par $M$ deux cercles tangents en ce point à la courbe et centrés sur $Ox,Oy$ respectivement (appelés cercles surosculateurs).

6. Que deviennent ces cercles lorsque $M$ est un sommet de l'ellipse.
\epb
\bpb 
 Soit la courbe $(C)$ définie par les formules:
$$M \left\{\begin{array}{l}
x=at^2 \\
y= at^3 \\
z=\ds \frac{9}{16}at^4\quad \mbox{avec}\quad a>0
\end{array}\right. $$
1. Calculer l'abscisse curviligne $s$ d'un point $M$ quelconque de cette courbe lorsqu'on prend pour origine des arcs l'origine des coordonnées et qu'on prend pour sens des arcs croissants celui des $y$ croissants.

2. Déterminer au point $M$ les vecteurs unitaires du trièdre de Frenêt.

3. Calculer le rayon de courbure et les coordonnées du centre de courbure.

4. Evaluer la torsion en $M$.
\epb
\bpb
On définit une surface $(S)$ par les équations:
	\[M(u,v) \left\{ \begin{array}{l}
	     X = u^2 + v \\
	    Y = u + v^2 \\
	    Z = uv
	    \end{array} \right.
\]
1.	Calculer les composantes des vecteurs $\textit{\textbf{OM}}'_u$ et $\textit{\textbf{OM}}'_v$.

2.	Calculer les coefficients $E,F,G$ de la première forme fondamentale de la surface $(S)$.

3.	En déduire l'expression de $ds^2$.

4.	Les coordonnées $(u,v )$ sont-elles orthogonales? symétriques?

5.	Calculer un vecteur normal de $(S)$. 
\epb
\bpb
On définit une surface $(\Sigma)$ par les équations:
 	\[M(u,v) \left\{ \begin{array}{l}
	           X =  a.cosu.cosv \\
	         Y = a.cosu.sinv \\
	          Z = b.sinu
	              \end{array} \right.
\] 
avec $a,b$ deux constantes positives. 

1.	Calculer les composantes des vecteurs $\textit{\textbf{OM}}'_u$ et $\textit{\textbf{OM}}'_v$.

2.	Calculer les coefficients $E,F,G$ de la première forme fondamentale de la surface $(\Sigma)$.

3.	En déduire l'expression de $ds^2$.

4.	Les coordonnées $( u,v )$ sont-elles orthogonales? symétriques?

5.	Calculer un vecteur unitaire normal $\textit{\textbf{n}}$ de $(\Sigma)$.

6.	Calculer les vecteurs :
$$   \textit{\textbf{OM}}''_{uu},\quad \textit{\textbf{OM}}''_{uv},\quad \textit{\textbf{OM}}''_{vv}$$
On pose:
$$L = \textit{\textbf{n}}.\textit{\textbf{OM}}''_{uu},\quad  M = \textit{\textbf{n}}.\textit{\textbf{OM}}''_{uv},\quad N = \textbf{\textit{n}}.\textit{\textbf{OM}}''_{vv} $$
7. Calculer les coefficients $L,M$ et $N$.
\epb
\bpb
 On considère la surface $(\Gamma)$ définie par les équations:
 	\[M(u,v)\left\{ \begin{array}{lll}
	           X =  sinu.cosv \\
	         Y = sinu.sinv \\
	          Z = cosu+Logtg\ds\frac{u}{2} +\psi(v)
	              \end{array} \right.
\] 
avec $\psi(v)$ est une fonction définie de classe $C^1$ de $v$. 

1. Donner le domaine de définition de la surface  $(\Gamma)$.

2.	Montrer que les courbes coordonnées $v=constante$ constituent une famille de courbes planes de $(\Gamma)$ et que leur plan coupe $(\Gamma)$ sous un angle constant.

3.  Calculer les composantes des vecteurs $\textit{\textbf{OM}}'_u$ et $\textit{\textbf{OM}}'_v$.

4.	Calculer les coefficients $E,F,G$ de la première forme fondamentale de la surface $(\Gamma)$.

5.	En déduire l'expression de $ds^2$.

6.	Les coordonnées $( u,v )$   sont-elles orthogonales? symétriques?

7.	On suppose pour la suite que $\psi(v)=0$, calculer un vecteur unitaire normal $\textit{\textbf{n}}$ de $\Gamma$.

8.	Calculer les vecteurs :
$$   \textit{\textbf{OM}}''_{uu},\quad \textit{\textbf{OM}}''_{uv},\quad \textit{\textbf{OM}}''_{vv}$$
On pose:
$$L = \textit{\textbf{n}}.\textit{\textbf{OM}}''_{uu},\quad  M = \textit{\textbf{n}}.\textit{\textbf{OM}}''_{uv},\quad N = \textbf{\textit{n}}.\textit{\textbf{OM}}''_{vv} $$
9. Calculer les coefficients $L,M$ et $N$.

10. En déduire l'expression des courbures moyenne et totale.
\epb
\bpb
 Soit la surface $(\Gamma)$ définie paramétriquement par:
 	\[M(u,v)\left\{ \begin{array}{lll}
	           X =  thu.cosv \\
	         Y = thu.sinv \\
	          Z = \ds \frac{1}{chu}+Logth\ds\frac{u}{2} 
	              \end{array} \right.
\] 
avec $chu$ et $thu$ sont respectivement le cosinus et la tangente hyperboliques définies par:
$$ chu=\frac{e^u+e^{-u}}{2},\quad thu=\frac{e^u+e^{-u}}{e^u-e^{-u}}$$ 
1.  Donner le domaine de définition de la surface  $(\Gamma)$.

2.  Calculer les composantes des vecteurs $\textit{\textbf{OM}}'_u$ et $\textit{\textbf{OM}}'_v$.

3.	Calculer les coefficients $E,F,G$ de la première forme fondamentale de la surface $(\Gamma)$.

4.	En déduire l'expression de $ds^2$.

5.	Les coordonnées $( u,v )$   sont-elles orthogonales? symétriques?

6.	Calculer un vecteur unitaire normal $\textit{\textbf{n}}$ de $(\Gamma)$.

7.	Calculer les vecteurs :
$$   \textit{\textbf{OM}}''_{uu},\quad \textit{\textbf{OM}}''_{uv},\quad \textit{\textbf{OM}}''_{vv}$$

On pose:
$$L = \textit{\textbf{n}}.\textit{\textbf{OM}}''_{uu},\quad  M = \textit{\textbf{n}}.\textit{\textbf{OM}}''_{uv},\quad N = \textbf{\textit{n}}.\textit{\textbf{OM}}''_{vv} $$
8. Calculer les coefficients $L,M$ et $N$.

9. Déterminer les coubures moyenne et totale.
\epb
  \bpb          
 Montrer que les courbures totale $K$ et moyenne $H$ en un point $M(x,y,z)$ d'une surface paramétrée par $z=f(x,y)$, où $f$ est une fonction lisse, sont données par:
$$K=\ds \frac{f''_{xx}f''_{yy}-f''^2_{xy}}{(1+f'^2_x+f'^2_y)^2}$$
 et:
$$ \ds H=\frac{(1+f'^2_x)f''_{xx}-2f'_xf'_yf''_{xy}+(1+f'^2_x)f''_{yy}}{(1+f'^2_x+f'^2_y)^{\frac{3}{2}}} $$
\epb
\bpb
 Soit $(\Sigma)$ une surface de $\BbR^3$ paramétrée par $OM(u,v)$ telle que sa première forme fondamentale s'écrit: $ds^2=Edu^2+2Fdudv+Gdv^2$
		
1. Montrer que les conditions suivantes sont équivalentes:

\quad i) - $\ds \frac{\partial E}{\partial v}=\frac{\partial G}{\partial u}=0$,
	
\quad	ii) - le vecteur $\ds \frac{\partial^2 OM }{\partial u\partial v}$ est parallèle au vecteur normal $N$ à la surface,

\quad iii) - les côtés opposés de tout quadrilatère curviligne formés par les courbes coordonnées $(u,v)$ ont même longueurs.

2. Quand ces conditions sont satisfaites, on dit que les courbes coordonnées de $(\Sigma)$ forment un réseau de \textit{Tchebychev}.\footnote{\textbf{Pafnouti Tchebychev} (1821 - 1894): mathématicien russe.} Montrer que dans ce cas, on peut paramétrer la surface par $(\tilde{u},\tilde{v})$ telle que $ds^2$ s'écrit:\index{\textbf{Tchebychev P.}}
$$ ds^2=d\tilde{u}^2+2cos\theta d\tilde{u}d\tilde{v}+d\tilde{v}^2$$
où $\theta$ est une fonction de $(\tilde{u},\tilde{v})$. Montrer que $\theta$ est l'angle entre les courbes coordonnées $\tilde{u},\tilde{v}$.

3. Montrer que l'expression de la courbure totale est donnée par:
$$ K= \ds \frac{1}{sin\theta}.\ds \frac{\partial^2\theta}{\partial \tilde{u} \partial \tilde{v}}$$
4. On pose :
$$\begin{array}{l}
\hat{u}=\tilde{u}+\tilde{v} \\
\hat{v}=\tilde{u}-\tilde{v}
\end{array}$$
Montrer que $ds^2$ s'écrit avec les nouvelles variables $(\hat{u},\hat{v})$:
$$ds^2=cos^2\omega d\hat{u}^2+sin^2\omega d\hat{v}^2$$ avec $\omega=\theta/2$. (\textit{A.N. Pressley}, 2010)\index{\textbf{Pressley A.N.}}
\epb
\bpb
    Soit $(\m F)$ une surface définie dans $\BbR^3$, paramétrée par la fonction vectorielle $\textbf{\textit{OM}}=F(u,v)$ telle que:
		$$
		F(u,v)\left|\begin{array}{l}
		x=F_1(u,v) \\
		y=F_2(u,v)\\
		z=F_3(u,v)
		\end{array}\right.
		$$
		$F$ est dite une paramétrisation conforme\index{Paramétrisation conforme} de $(\cal F)$ si on a les deux conditions suivantes:
		$$\frac{\partial F}{\partial u}.\frac{\partial F}{\partial u}=\frac{\partial F}{\partial v}.\frac{\partial F}{\partial v}=e^{\Phi(u,v)}\quad et \quad \frac{\partial F}{\partial u}.\frac{\partial F}{\partial v}=0$$
	1. Ecrire la première forme fondamentale de $(\cal F)$.
		
	2. Soit $n$	Le vecteur normal unitaire.
		$$ 		n=\ds \frac{\ds \frac{\partial F}{\partial u} \wedge \frac{\partial F}{\partial v}}{\left\| \ds \frac{\partial F}{\partial u} \wedge \frac{\partial F}{\partial v} \right\|} $$
		Quand le point $M$ varie sur la surface $(\cal F)$, le repère $ \ds (\frac{\partial F}{\partial u},\frac{\partial F}{\partial v},n)$ est un repère mobile. La deuxième forme fondamentale de $(\cal F)$ est définie par:
		$$ 		n.d^2F=Ldu^2+2Mdudv+Ndv^2 $$
	Si cette deuxième forme fondamentale s'écrit sous la forme :
	$$ 	-n.d^2F=e^{\Phi(u,v)}\left(\frac{du^2}{\rho_1}+\frac{dv^2}{\rho_2} \right)$$
		alors, la paramétrisation de $(\cal F)$ est dite isotherme. Dans ce cas, $\rho_1,\rho_2$ sont les rayons de courbure principaux de la surface $(\cal F)$. Une surface qui admet des coordonnées isothermes est dite isotherme.
	\index{Surface isotherme}	
				
	3. On considère que $(\cal F)$ est la sphère définie par:$ 		M=\left|\begin{array}{l}
		x=Rcos\varphi cos\lambda  \\
		y=Rcos\varphi sin\lambda \quad R>0\\
				z=Rsin\varphi              
				\end{array}\right. $
				
Soit $\m L_M$ la variable de Mercator. Montrer que la sphère paramétrée par $(\m L_M,\lambda)$ est une surface isotherme.

4. On considère $\m B$ la base du repère mobile  $ \ds (\frac{\partial F}{\partial u},\frac{\partial F}{\partial v},n)$. Exprimer les dérivées partielles $\ds \frac{\partial}{\partial u}$ et  $\ds \frac{\partial}{\partial u}$ des vecteurs de $\m B$ dans $\m B$, en tenant compte que la surface est isotherme c'est-à-dire qu'on a l'équation:
$$ 		-n.d^2F=e^{\Phi(u,v)}\left(\frac{du^2}{\rho_1}+\frac{dv^2}{\rho_2} \right)= -(L.du^2+2Mdu.dv+N.dv^2)$$
5. Montrer qu'on peut écrire les résultats de 4. sous la forme matricielle suivante:
$$ \frac{\partial}{\partial u}\begin{pmatrix}{
F'_u\cr
\cr
F'_v\cr
\cr
n }
\end{pmatrix}=\ds \begin{pmatrix} {
\ds \frac{\Phi'_u}{2} & \ds \frac{-\Phi'_v}{2} & \ds -\frac{e^\Phi}{\rho_1} \cr
\cr
\ds \frac{\Phi'_v}{2} & \ds \frac{\Phi'_u}{2} & 0 \cr
\cr
\ds \frac{1}{\rho_1} & 0 & 0}
\end{pmatrix}.\begin{pmatrix}{
F'_u\cr
\cr
F'_v\cr
\cr
n}
\end{pmatrix}$$
et:
$$\frac{\partial}{\partial v}\begin{pmatrix}{
F'_u\cr
\cr
F'_v\cr
\cr
n}
\end{pmatrix}=\ds \begin{pmatrix}{
\ds \frac{\Phi'_v}{2} & \ds \frac{\Phi'_u}{2} & 0 \cr
\cr
\ds-\frac{\Phi'_u}{2} & \ds \frac{\Phi'_v}{2} & \ds -\frac{e^\Phi}{\rho_2} \cr
\cr
0 & \ds \frac{1}{\rho_2} & 0 }
\end{pmatrix}.\begin{pmatrix}{
F'_u\cr
\cr
F'_v\cr
\cr
n}
\end{pmatrix}$$
Les deux dernières expressions ci-dessus sont appellées les équations de Gauss-Weingarten\footnote{\textbf{Julius Weingarten} (1836 - 1910) : mathématicien allemand.} de la surface $(\m F)$.
\epb
\index{Equations de Gauss-Weingarten}\index{\textbf{Gauss C.F.}}\index{\textbf{J. Weingarten}}
\chapter{\textit{\textbf{Géométrie de l'Ellipse et de l'Ellipsoïde}}}

\section{\textsc{Géométrie de l'Ellipse}}
\subsection*{5.1.1. Définitions}
\begin{figure}[htbp]
	\centering
		\includegraphics[width=1.00\textwidth]{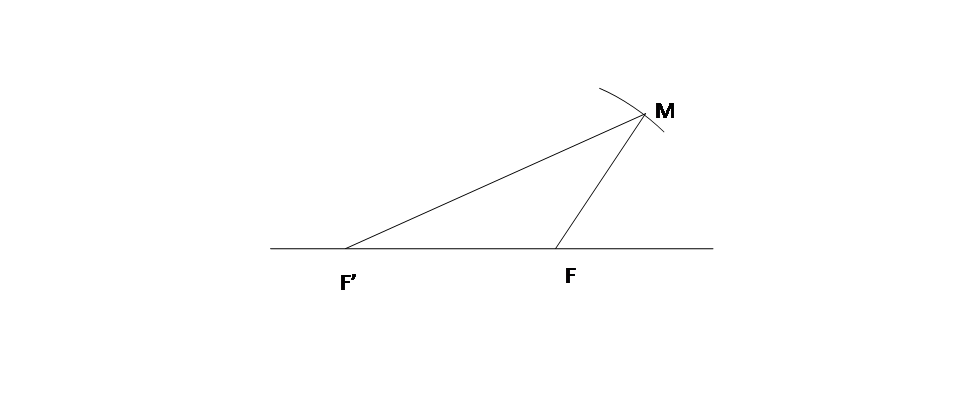}
	\caption{Définition de l'ellipse}
	\label{fig:defellipse1}
\end{figure}

\bdf
L'ellipse est le lieu des points dont la somme des distances à deux points fixes ou foyers est constante:
\be
	\fbox{ $  MF + MF' = constante=2a $}\label{p69}
\ee
où $a$ est dit le demi-grand axe de l'ellipse (\textbf{Fig. \ref{fig:defellipse1}}).
\edf

\begin{figure}[htbp]
	\centering
		\includegraphics[width=1.00\textwidth]{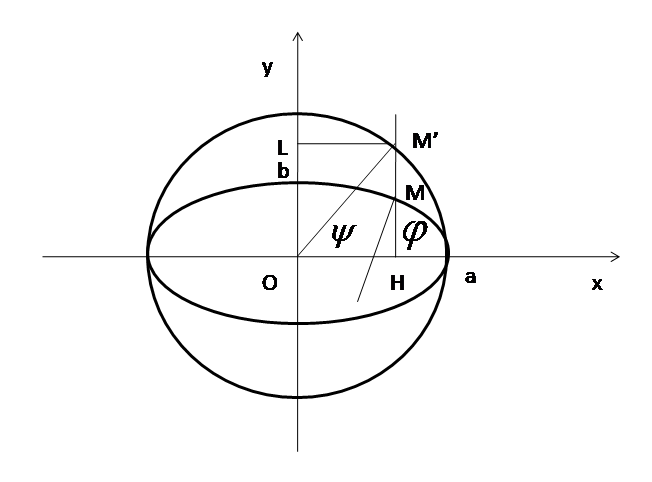}
	\caption{L'affinité}
	\label{fig:defellipse2}
\end{figure}

\bdf
Une ellipse est la transformée par affinité\index{Affinité} d'un cercle dans le rapport $b/a$ où $b$ est le demi-petit axe (\textbf{Fig. \ref{fig:defellipse2}}).
\edf
Au point $M'\in$  cercle   $\Longrightarrow   M \in$  ellipse  avec :
\be
	HM=\frac{b}{a}HM' \label{p70}
\ee
\newpage
Soit $\psi$ l'angle $\widehat{HOM'}$, $\psi$ est dite latitude paramétrique\index{Latitude paramétrique}  ou réduite\index{Latitude réduite}, d'où les coordonnées de $M'$:
\ba
	x=OH=OM'cos\psi  \nonumber \\
	y=OL=OM'sin\psi \nonumber
\ea
Par suite, les coordonnées de $M$ sur l'ellipse sont:
\be
\left\{\begin{array}{l}
	x=OH=acos\psi   \\
	y=OL=\ds \frac{b}{a}HM'=\frac{b}{a}asin\psi=bsin\psi  
	\end{array} \right. \label{p73}
\ee
Dans le système d'axes $Ox, Oy$, l'équation de l'ellipse s'écrit:
$$ \ds 	\frac{x^2}{a^2}+\frac{y^2}{b^2}=1 $$
On appelle respectivement aplatissement, le carré de la première excentricité et le carré de la deuxième excentricité les quantités:   
\be
	\fbox{ $	\alpha=\ds \frac{a-b}{a}, \quad e^2 = \ds \frac{a^2-b^2}{a^2}, \quad e'^2=\ds \frac{a^2-b^2}{b^2} $} \label{p78}
\ee

\section{\textsc{Equations paramétriques de l'ellipse}}
Les  équations (\ref{p73}) représentent les équations paramétriques\index{Equations paramétriques} de l'ellipse en fonction de la latitude $\psi$. On va exprimer ces équations en fonction de l'angle $\varphi$  de la normale en $M$ avec l'axe $Ox$.
\\

\begin{figure}[htp]
	\centering
		\includegraphics[width=1.00\textwidth]{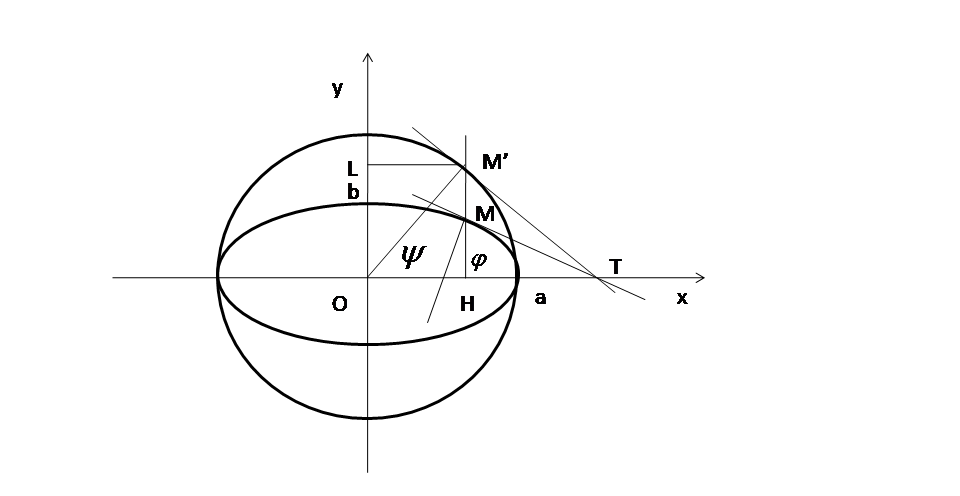}
	\caption{La relation entre $\varphi$ \,\,\, et\,\,\, $\psi$}
	\label{fig:relation}
\end{figure}

Soit $TM'$ la tangente en $M$' au cercle de rayon $a$, le point $T$ est l'intersection de cette tangente avec l'axe $Ox$. La transformée de cette tangente par affinité de rapport $b/a$ de cette tangente est la droite tangente à l'ellipse au point $M$ et elle passe par $T$ (\textbf{Fig. \ref{fig:relation}}).
\\

Dans le triangle $MHT$, on a: 
$$ tg\varphi  =\frac{HT}{MH} $$
et dans le triangle $M'HT$, on a:
	$$ tg\psi  = \frac{HT}{M'H} $$
d'où:
$$ 	\frac{tg\psi}{tg\varphi}=\frac{HT}{M'H}\frac{MH}{HT}=\frac{MH}{M'H}= \mbox{rapport de l'affinité}=\frac{b}{a} $$ 
Soit:          
\be
	\fbox{ $ tg\psi=\ds \frac{b}{a}tg\varphi $} \label{p80}
\ee
De (\ref{p80}), on exprime $cos\psi$  et $sin\psi$  en fonction de l'angle $\varphi $, d'où:
	$$\frac{1}{cos^2\psi}=1+tg^2\psi=1+(b/a)^2tg^2\varphi=\frac{a^2cos^2\varphi+b^2sin^2\varphi}{a^2cos^2\varphi}$$
D'où:                                                             
	$$cos^2\psi=\frac{a^2cos^2\varphi}{a^2cos^2\varphi+b^2sin^2\varphi} $$
On pose:          
\be
	W^2=\frac{a^2cos^2\varphi+b^2sin^2\varphi}{a^2}=1-e^2sin^2\varphi \label{p81}
\ee
d'où:                   
\be
	W=\frac{cos\varphi}{cos\psi} \label{p82}
\ee
\newpage
On calcule de même $sin\psi$: 
	$$  sin^2\psi= 1 - cos^2\psi= 1 - \frac{a^2cos^2\varphi}{a^2cos^2\varphi+b^2sin^2\varphi} $$
soit: 
\be
sin^2\psi=	\frac{b^2sin^2\varphi}{a^2cos^2\varphi+b^2sin^2\varphi} \label{p83}
\ee
On pose :                   
\be
	V^2=\frac{a^2cos^2\varphi+b^2sin^2\varphi}{b^2}=1-e'^2cos^2\varphi \label{p84}
\ee
avec $e'$ la 2ème excentricité, d'où:
\be
	V=\frac{sin\varphi}{sin\psi}=\frac{a}{b}W \label{p85}
\ee
  Alors les équations paramétriques de l'ellipse en fonction de $\varphi$  sont:
\ba
&		X= \ds acos\psi=a\frac{cos\varphi}{W}=a\frac{cos\varphi}{\sqrt{1-e^2sin^2\varphi}} & \label{p86} \\
& Y=bsin\psi=\ds \frac{bsin\varphi}{V}=\frac{b^2sin\varphi}{a\sqrt{1-e^2sin^2\varphi}}=a(1-e^2)\frac{sin\varphi}{\sqrt{1-e^2sin^2\varphi}} &\label{p87}
\ea
Soit:
\be
\fbox{ $ \begin{array}{l}
	X=\ds a\frac{cos\varphi}{\sqrt{1-e^2sin^2\varphi}} \\
	Y=\ds a(1-e^2)\frac{sin\varphi}{\sqrt{1-e^2sin^2\varphi}} 
	\end{array} $}
\ee
\subsection*{5.2.1. Relations différentielles entre $\varphi$ et $\psi$}
A partir de la relation (\ref{p80}), on obtient:
\be
\frac{d\psi}{cos^2\psi}=\frac{bd\varphi}{acos^2\varphi}\Longrightarrow\frac{d\psi}{d\varphi}=\frac{bcos^2\psi}{acos^2\varphi} \label{p90}
\ee
et en utilisant (\ref{p82}), on a:
\be
\fbox{ $ \ds \frac{d\psi}{d\varphi}=\frac{b}{aW^2} =\frac{b}{a(1-e^2sin^2\varphi)} $} \label{p91}
\ee
\section{\textsc{Calcul de la grande normale }}
\begin{figure}[htbp]
	\centering
		\includegraphics[width=1.00\textwidth]{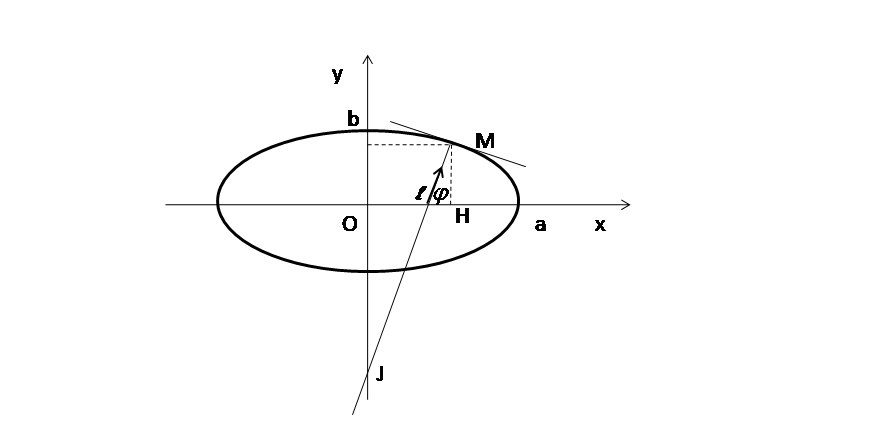}
	\caption{La grande normale}
	\label{fig:lagrandenormale}
\end{figure}

\bdf
On appelle la grande normale \index{Grande normale} la longueur de $JM$. $JM$ est porté par la normale à l'ellipse au point $M$. La normale a pour vecteur de direction, le vecteur $\textbf{\textit{l}}$ de composantes $(cos\varphi , sin\varphi )$ (\textbf{Fig. \ref{fig:lagrandenormale}}).
\edf
D'où l'équation cartésienne de la normale:
\be
	\frac{X-X_M}{cos\varphi}=\frac{Y-Y_M}{sin\varphi} \label{p92}
\ee

On obtient l'ordonnée de $J$ en faisant $X = 0$ dans (\ref{p92}), d'où:
$$ \ds	\frac{-X_M}{cos\varphi}=\frac{Y-Y_M}{sin\varphi}\Longrightarrow Y_J=Y_M-X_Mtg\varphi $$
Par suite, la distance $MJ$ est égale à:
$$	MJ = \sqrt{(Y_J-Y_M)^2+X^2_M}=\sqrt{X^2_Mtg^2\varphi+X^2_M}=X_M\sqrt{1+tg^2\varphi} $$
Soit:
$$ 	MJ=\ds \frac{X_M}{cos\varphi} $$
\newpage
Or:      
	$$ X_M =\ds \frac{acos\varphi}{W}\Longrightarrow MJ=  \frac{acos\varphi}{Wcos\varphi}=\frac{a}{W}=\frac{a}{\sqrt{1-e^2sin^2\varphi}} $$
On pose:          
\begin{equation}
\fbox{ $ 	N(\varphi)=MJ=\ds \frac{a}{\sqrt{1-e^2sin^2\varphi}} $} \label{p97}
\end{equation}
$N$ est appelé \textit{la grande normale}.
\\

Les équations paramétriques de l'ellipse (\ref{p73})  deviennent:
\ba
X=acos\psi=a\frac{cos\varphi}{W}=a\frac{cos\varphi}{\sqrt{1-e^2sin^2\varphi}}=N(\varphi)cos\varphi \nonumber \\
Y=bsin\psi=b\frac{sin\varphi}{V}=a(1-e^2)\frac{sin\varphi}{\sqrt{1-e^2sin^2\varphi}}=(1-e^2)N(\varphi)sin\varphi \nonumber
\ea
Soit:
\be
\fbox{ $ 
\begin{array}{l}
X=\ds a\frac{cos\varphi}{\sqrt{1-e^2sin^2\varphi}}=N(\varphi)cos\varphi  \\
Y=\ds a(1-e^2)\frac{sin\varphi}{\sqrt{1-e^2sin^2\varphi}}=(1-e^2)N(\varphi)sin\varphi 
\end{array} $} \label{p99}
\ee
\subsection*{5.3.1. Arc élémentaire $ds$ et rayon de courbure $\rho$ de l'ellipse}
L'arc élémentaire $ds$ se calcule à partir des équations paramétriques de l'ellipse par:
\ba
	                  ds^2 = dX^2 + dY^2 = a^2sin^2\psi d\psi^2  + b^2.cos^2\psi d\psi^2 \nonumber \\ 
\mbox{ou } ds^2 = (a^2sin^2\psi +b^2.cos^2\psi )d\psi^2  \nonumber 
\ea
En utilisant les équations (\ref{p82}) et (\ref{p85}), on obtient:
	$$ ds = \frac{b}{W}.d\psi  $$
Et en remplaçant $d\psi$  en utilisant (\ref{p91}), on trouve:
	$$	ds= a(1-e^2)\frac{d\varphi}{(1-e^2sin^2\varphi )^{3/2}}$$
\newpage
La longueur de l'arc de méridien comptée depuis l'équateur est:
\begin{equation}
	\fbox{ $ s(\varphi) =\ds  \int^{\varphi}_{0}ds =\ds a(1-e^2)\int^{\varphi}_{0}\ds  \frac{dt}{(1-e^2sin^2t)^{3/2}} $} \label{p103} 
\end{equation}
L'intégration se fait à  partir d'un développement limité de $(1-e^2sin^2t )^{-3/2}$ (Voir plus loin). Le rayon de courbure\index{Rayon de courbure de l'ellipse} $\rho$ de l'ellipse s'obtient à partir de $ds$ par:
\begin{equation}
	\fbox{ $\rho=\ds \frac{ds}{d\varphi}=\frac{b^2}{aW^3}=\frac{a(1-e^2)}{(1-e^2sin^2\varphi)^{3/2}} $} \label{p104} 
\end{equation}

\section{\textsc{Géométrie de l'Ellipsoïde de Révolution}}
On va étudier les propriétés de l'ellipsoïde de révolution\index{Ellipsoïde de révolution} obtenu par la rotation d'une ellipse autour du demi-petit axe comme le montre la figure ci-dessous (\textbf{Fig. \ref{fig:ellip}}):
\begin{figure}
\begin{center}
	\includegraphics[width=0.85\textwidth]{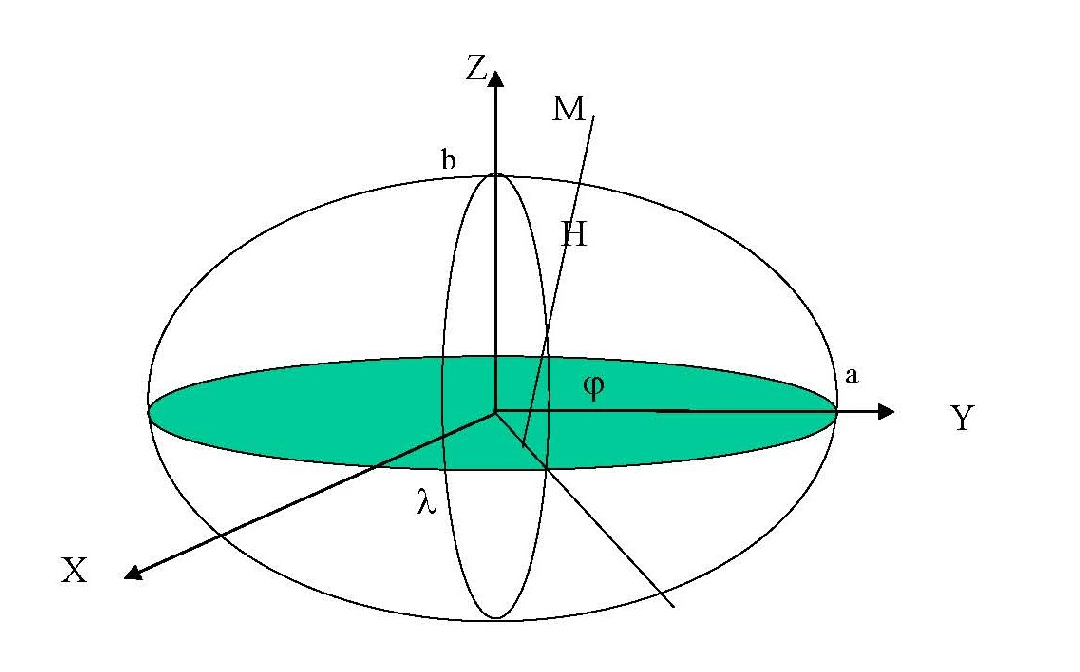}
	\caption{L'ellipsoïde de révolution: ellipsoïde de référence}
	\label{fig:ellip}
\end{center}
\end{figure}
\subsection*{5.4.1. Les Coordonnées Géographiques}
Les coordonnées géographiques définies sur l'ellipsoïde de révolution sont:

-	la longitude $\lambda$: angle du plan méridien du point M avec le plan méridien origine, dans notre cas, le plan origine est le plan $XOZ$,

-	la latitude $\varphi$: angle de la direction de la normale au point M avec le plan équatorial;

-	l'altitude ellipsoïdique $he$: si le point est sur l'ellipsoïde $he=0$.
\\

Dans le plan $ROZ$ (\textbf{Fig. \ref{fig:ellips3}}) avec  $\textbf{\textit{r}}$ et $\textbf{\textit{k}}$  les vecteurs unitaires des axes $OR$ et $OZ$, on peut écrire:
\begin{figure}
	\centering
		\includegraphics[width=0.70\textwidth]{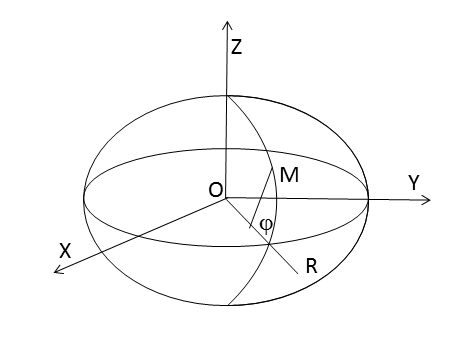}
	\caption{Calcul des coordonnées géodésiques}
	\label{fig:ellips3}
\end{figure}
\ba
	             \textbf{\textit{OM}} = acos\psi  \textbf{\textit{r}} + bsin\psi  \textbf{\textit{k}}  \nonumber  \\
\mbox{et} \, \, \,      \textbf{\textit{r}} = cos\lambda  \textbf{\textit{i}} + sin\lambda  \textbf{\textit{j}} \nonumber 
\ea
D'où:   
$$ \textbf{\textit{OM}} = acos\psi cos\lambda \textbf{\textit{i}} + acos\psi    sin\lambda \textbf{\textit{j}} +            bsin\psi    \textbf{\textit{k}}  $$
Donc, les équations paramétriques du point M sont:
\ba
	X=acos\psi cos\lambda  \nonumber\\
	Y = acos\psi sin\lambda \nonumber \\
	Z=bsin\psi \nonumber 
\ea
Et en exprimant $\psi$ en fonction de $\varphi$, on a:
\ba
	X=acos\psi cos\lambda=\frac{acos\varphi}{W}cos\lambda=Ncos\varphi cos\lambda \nonumber \\
	Y = acos\psi sin\lambda = Ncos\varphi sin\lambda \nonumber   \\
	Z=bsin\psi = \frac{b^2sin\varphi}{a\sqrt{1-e^2sin^2\varphi}}=N(1-e^2)sin\varphi \nonumber  
\ea
Soit:
 \be
\fbox{ $ 
\begin{array}{l}
	X=Ncos\varphi cos\lambda               \\
	Y = Ncos\varphi sin\lambda                \\
	Z=N(1-e^2)sin\varphi 
	\end{array} $} \label{p113a}  
\ee

Si $he \neq  0$, alors les coordonnées de $M$ sont:
\be
\fbox{ $ \begin{array}{l}
	X=(N+he)cos\varphi cos\lambda   \\
	Y = (N+he)cos\varphi sin\lambda   \\
	Z=(N(1-e^2)+he)sin\varphi 
		\end{array} $}  \label{p116}  
\ee
\subsection*{5.4.2. Passage des coordonnées tridimensionnelles $(X,Y,Z)$ aux coordonnées $(\varphi ,\lambda ,he)$}
Des deux premières équations de  (\ref{p116}) et ne pas tenir compte du cas particulier $(X=0)$, on obtient:                                         
\begin{equation}
	\fbox{ $ tg\lambda=\ds \frac{Y}{X}\Longrightarrow \,\lambda=Arctg\frac{Y}{X} $} \label{p117} 
\end{equation}
De même, on pose:
	$$ r=\sqrt{X^2+Y^2}=(N+he)cos\varphi $$
De (\ref{p116}), on peut écrire: 
\begin{equation}
	Z=(N+he)sin\varphi-Ne^2sin\varphi \label{p119} 
\end{equation}
soit:
	$$ Z=Z'-Ne^2sin\varphi $$
en posant:
\begin{equation}
	Z'=(N+he)sin\varphi \label{p121} 
\end{equation}
                                        \newpage
Le calcul de $\varphi$ se fait par itérations:  
1ère itération:                                     
\begin{equation}
	Z'=Z\Rightarrow tg\varphi=\frac{Z'}{r}\Rightarrow \varphi_1
=Arctg\frac{Z'}{r} \label{p122} 
\end{equation}
2ème itération:    
\ba
	N = a(1-e^2sin^2\varphi_1)^{-1/2} \nonumber \\
   Z' = Z + Ne^2.sin\varphi_1 \nonumber \\
    \varphi_2 = Arctg(Z'/r) \label{p123} 
\ea
3ème itération :   
\ba
	N = a(1-e^2sin^2\varphi_2)^{-1/2} \nonumber \\
                         Z' = Z + Ne^2.sin\varphi_2 \nonumber \\
                                  \varphi_3 = Arctg( Z'/r)  \label{p124}     
	\ea
En général, 3 à 4 itérations suffisent et on obtient:     
\begin{equation}
		  \fbox{ $  \varphi  =  \varphi_3 $} \label{p125}
	\end{equation}
Par suite, on peut déterminer l'altitude géodésique $he$ :                                                   
\begin{equation}
	\fbox{ $ he=\ds \frac{r}{cos\varphi}-N(\varphi) $}\label{p126} 
\end{equation}
\newpage
\section{\textsc{Calcul des Lignes Géodésiques de L'Ellipsoïde de Révolution}} 
\begin{verse}
\begin{svgraybox}
	" A côté de la difficulté principale, de celle qui tient au fond même des choses, il y a une foule de difficultés secondaires qui viennent compliquer encore la tâche du chercheur. Il y aurait donc intérêt à étudier d'abord un problème où l'on rencontrerait cette difficulté principale, mais où l'on serait affranchis de toutes les difficultés secondaires. Ce problème est tout trouvé, c'est celui des \textbf{lignes géodésiques} d'une surface; c'est encore un problème de dynamique, de sorte que la difficulté principale subsiste; mais c'est le plus simple de tous les problèmes de dynamique. " 
\end{svgraybox}
\begin{flushright}
(\textbf{H. Poincaré}\footnote{\textbf{Henri Poincaré} (1854-1912): mathématicien français, parmi les plus grands du XIXème siècle.}, 1905)
\end{flushright}
\end{verse}
\index{\textbf{Poincaré H.}}
Après avoir défini les lignes géodésiques \index{Ligne géodésique} d'une surface, on établit  les équations des géodésiques pour une surface donnée. Comme application, on détaille celles de l'ellipsoïde de révolution. On fera l'intégration de ces équations. 
\normalsize
\subsection*{5.5.1. Introduction et Notations}
Soit  $(S)$ une surface définie par les paramètres   $(u,v)$  avec  $(u,v)\in \m D$ un domaine  $ \subset \BbR^2$. Un point  $M\in  (S)$ vérifie :
\begin{equation}
	 \textbf{\textit{OM}} =\textbf{\textit{OM}}(u,v)\left|\begin{array}{l}
		x(u,v) \\
		y(u,v) \\
		z(u,v)
		\end{array}\right. \label{k1}
\end{equation}
On introduit les notations usuelles :
$$ E=\ds \frac{\partial \textbf{\textit{M}}}{\partial u}.\frac{\partial \textbf{\textit{M}}}{\partial u}=\left\|\frac{\partial \textbf{\textit{M}}}{\partial u}\right\|^2  $$
\be
\begin{array}{l}
			F=\ds \frac{\partial \textbf{\textit{M}}}{\partial u}.\frac{\partial \textbf{\textit{M}}}{\partial v} \\
		\\
	G=\ds \frac{\partial \textbf{\textit{M}}}{\partial v}.\frac{\partial \textbf{\textit{M}}}{\partial v} =\left\|\frac{\partial \textbf{\textit{M}}}{\partial v}\right\|^2 
	\end{array} \label{k1c}   
\ee
 Des équations  (\ref{k1c}), on obtient les équations :                     
\be
\begin{array}{l}
	\ds \frac{\partial E}{\partial u}=2\frac{\partial \textbf{\textit{M}}}{\partial u}.\frac{\partial ^2\textbf{\textit{M}}}{\partial u^2}  \\
		\\
		\ds \frac{\partial E}{\partial v}=2\frac{\partial \textbf{\textit{M}}}{\partial u}.\frac{\partial ^2\textbf{\textit{M}}}{\partial u \partial v}  \\
			\\
			\ds \frac{\partial F}{\partial u}=\frac{\partial^2 \textbf{\textit{M}}}{\partial u^2}.\frac{\partial \textbf{\textit{M}}}{\partial v}+ \frac{\partial \textbf{\textit{M}}}{\partial u}.\frac{\partial^2\textbf{\textit{M}}}{\partial u \partial v} \\
		\\
		\ds \frac{\partial F}{\partial v}=\frac{\partial^2 \textbf{\textit{M}}}{\partial v^2}.\frac{\partial \textbf{\textit{M}}}{\partial u}+ \frac{\partial \textbf{\textit{M}}}{\partial v}.\frac{\partial^2\textbf{\textit{M}}}{\partial u \partial v} \\
	\\
	\ds \frac{\partial G}{\partial u}=2\frac{\partial \textbf{\textit{M}}}{\partial v}.\frac{\partial ^2\textbf{\textit{M}}}{\partial u \partial v}  \\
		\\
		\ds \frac{\partial G}{\partial v}=2\frac{\partial \textbf{\textit{M}}}{\partial v}.\frac{\partial ^2\textbf{\textit{M}}}{\partial v^2} 
		\end{array}  \label{k2f} 
\ee
 Soit  $\textbf{\textit{n}}$ le vecteur unitaire normal en  $M(u,v)$ à la surface  $(S)$, $\textbf{\textit{n}}$ est donné par  :  
   \begin{equation}
\textbf{\textit{n}}=\frac{\ds \frac{\partial \textbf{\textit{M}}}{\partial u}\wedge \ds  \frac{\partial \textbf{\textit{M}}}{\partial u}}{H}	\label{k3}
\end{equation}
 avec:
\begin{equation}
	H=\left\|\frac{\partial \textbf{\textit{M}}}{\partial u}\wedge \frac{\partial \textbf{\textit{M}}}{\partial u} \right\|  \label{k4}
\end{equation}
D'où :   
\begin{equation}
	ds^2= E.du^2 + 2.F.du.dv + G.dv^2  \label{k5}
\end{equation}
  L'équation (\ref{k5}) représente le carré infinitésimal de la longueur de l'arc.
    
Soit une courbe  $(\Gamma)$  tracée sur  $(S)$ et  $\textbf{\textit{N}}$ est le vecteur unitaire de la normale  principale le long de $(\Gamma)$.                            
\\

 \bdf
Une courbe $(\Gamma)$ est dite ligne géodésique de la surface  $(S)$ si et seulement si les  vecteurs  $\textbf{\textit{n}}$ et  $\textbf{\textit{N}}$ sont colinéaires.
\edf  
On démontre par le calcul des variations (\textit{P. Petersen}, 1998) que la ligne géodésique entre deux points d'une surface $(S)$ lorsqu'elle existe est la courbe de longueur minimale joignant les deux points.  \index{\textbf{Petersen P.}}    
\subsection*{5.5.2. Les Equations Différentielles des Lignes Géodésiques }
On calcule l'expression de  $\textbf{\textit{N}}$, on obtient :
	$$ \textbf{\textit{N}}=R\frac{d\textbf{\textit{T}}}{ds}  $$
 or:
$$	\textbf{\textit{T}}=\frac{d\textbf{\textit{M}}}{ds}=\frac{\partial \textbf{\textit{M}}}{\partial u}\frac{du}{ds}+\frac{\partial \textbf{\textit{M}}}{\partial v}\frac{dv}{ds}  $$
 d'où:
	$$ \frac{d\textbf{\textit{T}}}{ds}=\frac{\partial^2 \textbf{\textit{M}}}{\partial u^2}\left(\frac{du}{ds}\right)^2+2\frac{\partial^2 \textbf{\textit{M}}}{\partial u \partial v}\frac{du}{ds} \frac{dv}{ds}+\frac{\partial \textbf{\textit{M}}}{\partial u}\frac{d^2u}{ds^2}+\frac{\partial \textbf{\textit{M}}}{\partial v}\frac{d^2u}{ds^2}+\frac{\partial^2 \textbf{\textit{M}}}{\partial v^2}\left(\frac{dv}{ds}\right)^2 $$

La condition $\textbf{\textit{n}}$ // $\textbf{\textit{N}}$ peut être écrite:
	$$ \textbf{\textit{N}}\wedge \textbf{\textit{n}}=0 $$
soit:
\begin{equation}
	R\frac{d\textbf{\textit{T}}}{ds}\wedge \left(\ds \frac{\ds \frac{\partial \textbf{\textit{M}}}{\partial u}\wedge \ds  \frac{\partial \textbf{\textit{M}}}{\partial u}}{H}\right)=0 \label{k9-2}
\end{equation}
Utilisant la  formule du produit vectoriel:  
\begin{equation}
	 \textbf{\textit{A}} \wedge ( \textbf{\textit{B}} \wedge \textbf{\textit{C}})= ( \textbf{\textit{A}}.\textbf{\textit{C}}) \textbf{\textit{B}}  - ( \textbf{\textit{A}}.\textbf{\textit{B}}) \textbf{\textit{C}} \label{k10}
\end{equation}
on obtient:
	$$ \left(\frac{d\textbf{\textit{T}}}{ds}.\frac{\partial \textbf{\textit{M}}}{\partial v}\right)\frac{\partial \textbf{\textit{M}}}{\partial u}-\left(\frac{d\textbf{\textit{T}}}{ds}.\frac{\partial \textbf{\textit{M}}}{\partial u}\right)\frac{\partial \textbf{\textit{M}}}{\partial v}=0 $$
Or $\displaystyle \frac{\partial \textbf{\textit{M}}}{\partial u}$ et  $\displaystyle \frac{\partial \textbf{\textit{M}}}{\partial v}$ forment une base du plan tangent en M, d'où les deux conditions:
\begin{equation}
	\frac{d\textbf{\textit{T}}}{ds}.\frac{\partial \textbf{\textit{M}}}{\partial v}=0 \quad \mbox{et} \quad \frac{d\textbf{\textit{T}}}{ds}.\frac{\partial \textbf{\textit{M}}}{\partial u}=0 \label{k12}
\end{equation}
Ce qui donne deux équations différentielles du second ordre:
\begin{equation}
 \frac{\partial^2 \textbf{\textit{M}}}{\partial u^2}.\frac{\partial \textbf{\textit{M}}}{\partial v}\left(\frac{du}{ds}\right)^2+F\frac{d^2u}{ds^2}+2\frac{\partial^2 \textbf{\textit{M}}}{\partial u \partial v}.\frac{\partial \textbf{\textit{M}}}{\partial v}\frac{du}{ds} \frac{dv}{ds}+\frac{\partial^2 \textbf{\textit{M}}}{\partial v^2}.\frac{\partial \textbf{\textit{M}}}{\partial v}\left(\frac{dv}{ds}\right)^2+G\frac{d^2 v}{ds^2}=0 \label{k13-1}	
\end{equation}
\newpage
et:
\begin{equation}
  \frac{\partial^2 \textbf{\textit{M}}}{\partial v^2}.\frac{\partial \textbf{\textit{M}}}{\partial u}\left(\frac{dv}{ds}\right)^2+F\frac{d^2v}{ds^2}+2\frac{\partial^2 \textbf{\textit{M}}}{\partial u \partial v}.\frac{\partial \textbf{\textit{M}}}{\partial u}\frac{du}{ds} \frac{dv}{ds}+\frac{\partial^2 \textbf{\textit{M}}}{\partial u^2}.\frac{\partial \textbf{\textit{M}}}{\partial u}\left(\frac{du}{ds}\right)^2+E\frac{d^2 u}{ds^2}=0 \label{k13-2}	
\end{equation}
On pose:
\ba
	E'_u=\frac{\partial E}{\partial u}; \quad 	E'_v=\frac{\partial E}{\partial v}; \quad 
	F'_u=\frac{\partial F}{\partial u} \nonumber \\
	F'_v=\frac{\partial F}{\partial v}; \quad	G'_u=\frac{\partial G}{\partial u}; \quad
	G'_v=\frac{\partial G}{\partial v} \label{k14} 
\ea
et on utilise les équations (\ref{k2f}), (\ref{k13-1}) et \ref{k13-2}), ces 2 dernières équations peuvent être écrites :
\ba
	\fbox{ $\ds (F'_u-\frac{E'_v}{2})\left(\frac{du}{ds}\right)^2+F\frac{d^2u}{ds^2}+G'_u\frac{du}{ds} \frac{dv}{ds}+\frac{G'_v}{2}\left(\frac{dv}{ds}\right)^2+G\frac{d^2 v}{ds^2}=0 $}\label{k15-1}	\\
	\nonumber \\
	\fbox{ $ \ds (F'_v-\frac{G'_u}{2})\left(\frac{dv}{ds}\right)^2+F\frac{d^2v}{ds^2}+E'_v\frac{dv}{ds} \frac{du}{ds}+\frac{E'_u}{2}\left(\frac{du}{ds}\right)^2+E\frac{d^2 u}{ds^2}=0 $}\label{k15-2}	
\ea
\subsection*{5.5.3. Détermination des Lignes Géodésiques de l'Ellipsoïde de révolution}
On considère maintenant comme surface l'ellipsoïde de  révolution qu'on paramètre comme suit: 
\ba
	X = N.cos\varphi cos\lambda \nonumber \\
  Y = Ncos\varphi sin\lambda  \label{k16} \\    
        Z = N(1-e^2)sin\varphi \nonumber 
\ea
où: 
	$$N=\frac{a}{\sqrt{1-e^2sin^2\varphi}}=aW^{-1/2} $$ 
est le rayon de courbure de la grande normale avec:         
	$$ W= 1-e^2sin^2\varphi  $$                                
On note:
	$$ r=Ncos\varphi $$
le rayon du parallèle de latitude $\varphi$ et $\rho$  le rayon de courbure de la méridienne donné par : 
	$$\rho=\frac{a(1-e^2)}{(1-e^2sin^2\varphi)\sqrt{1-e^2sin^2\varphi}}=a(1-e^2)W^{-3/2} $$
Alors la première forme fondamentale s'écrit :   
\begin{equation}
	ds^2=\rho^2d\varphi^2+r^2d\lambda^2 \label{k20}
\end{equation}
En prenant comme variables  $ u = \varphi$ et  $v =\lambda$, on obtient: 
\ba
	&E = E(\varphi)= \rho^2, \quad       F = 0, \quad           G = r^2 \label{k21} &\\&
	E'_{\varphi} = 2\rho \rho',\,E'_{\lambda} = 0,\,F'_{\varphi}=F'_{\lambda}=0,\,G'_{\varphi}= 2rr' =-2r\rho sin\varphi,\,G'_{\lambda} = 0     \label{k22} &  
\ea
Alors les équations (\ref{k15-1}) et (\ref{k15-2}) deviennent :
\ba
	-2r\rho sin\varphi \frac{d \varphi}{ds}\frac{d\lambda}{ds}+r^2\frac{d^2\lambda}{ds^2}=0 \label{k23} \\
		r\rho sin\varphi \left(\frac{d\lambda}{ds}\right)^2+\rho \rho'\left(\frac{d\varphi}{ds}\right)^2+\rho^2\frac{d^2\varphi}{ds^2}=0 \label{k24} 
\ea
La première équation s'écrit:
\begin{equation}
	\frac{d}{ds}\left(r^2\frac{d\lambda}{ds}\right)=0 \label{k25} 
\end{equation}
dont l'intégration donne  :
  \begin{equation}
	  \ds r^2\frac{d\lambda}{ds}=C=constante   \label{k26}
\end{equation}
On retrouve alors la  relation de Clairaut (\textit{J. Lemenestrel}, 1980):\footnote{\textbf{Alexis Claude de Clairaut} (1713-1765): mathématicien, astronome et géophysicien français.}\index{Relation de Clairaut}\index{\textbf{Clairaut A.C.}}\index{\textbf{Lemenestrel J.}}
\begin{equation}
	\fbox{$  r.sinAz = constante = C=asinAze   $}         \label{k27}                     
\end{equation}
où $Az$ est l'azimut de la géodésique au point $M$ et $Aze$ son azimut initial au point $M_0$ à l'équateur.

L'équation (\ref{k24}) s'écrit:
\begin{equation}
		\rho\left(	r sin\varphi \left(\frac{d\lambda}{ds}\right)^2+\rho'\left(\frac{d\varphi}{ds}\right)^2+\rho \frac{d^2\varphi}{ds^2}\right)=0 \label{k27a} 
\end{equation}
Ce qui donne:

- $\rho=0$ le point $M$ est sur l'équateur: $\varphi=0$ et $r=a$ le demi-grand axe de l'ellipsoïde et l'équation (\ref{k23}) devient:
\begin{equation}
	\frac{d^2\lambda}{ds^2}=0 \label{k27b}
\end{equation}
dont l'intégration donne:
\begin{equation}
	\lambda-\lambda_0=l(s-s_0)  \label{k27c}
\end{equation}
le point $M$ décrit l'équateur et la géodésique est le grand cercle de rayon $a$.

- $\rho \neq 0$, le point $M$ n'est pas sur l'équateur, l'équation  (\ref{k24}) s'écrit comme  suit:
\begin{equation}
	\rho\frac{d^2\varphi}{ds^2}	+\rho'\left(\frac{d\varphi}{ds}\right)^2+	rsin\varphi \left(\frac{d\lambda}{ds}\right)^2=0 \label{k28} 
\end{equation}
 Pour intégrer (\ref{k28}), on utilise une nouvelle fonction, soit :   
   \begin{equation}
	Z=\frac{d\lambda}{d\varphi}  \label{k29} 
\end{equation}
 De (\ref{k26}), on obtient :                                                         
$$
\frac{d\varphi}{ds}=\frac{d\varphi}{d\lambda}\frac{d\lambda}{ds}=\frac{C}{r^2}\frac{d\varphi}{d\lambda}=\frac{C}{r^2Z} $$
soit:
\begin{equation}
 \frac{d\varphi}{ds}=\frac{C}{r^2Z} \label{k30} 
\end{equation}
On exprime maintenant la dérivée seconde $d^2\varphi/ds^2$ :
\begin{equation}
 \frac{d^2\varphi}{ds^2}=\frac{d}{ds}\left(\frac{d\varphi}{ds} \right)= \frac{d}{d\varphi}\left(\frac{d\varphi}{ds} \right) \frac{d\varphi}{ds}=\frac{1}{2}\frac{d}{d\varphi}\left(\frac{d\varphi}{ds} \right)^2 \label{k31} 
\end{equation}
L'équation (\ref{k28}) s'écrit en utilisant (\ref{k26}) et (\ref{k31}) : 
\begin{equation}
	\frac{\rho}{2}\frac{d}{d\varphi}\left[\left( \frac{d\varphi}{ds}\right)^2\right]+\rho'\left( \frac{d\varphi}{ds}\right)^2+sin\varphi\left(\frac{C^2}{r^3} \right)=0 \label{k32}
\end{equation}
On pose:
\begin{equation}
	U=\left( \frac{d\varphi}{ds}\right)^2  \label{k33}
\end{equation}
L'équation (\ref{k32}) devient:
\begin{equation}
		\frac{\rho}{2}\frac{dU}{d\varphi}+\rho'U=-\frac{C^2sin\varphi}{r^3} \label{k34}
\end{equation}
L'équation (\ref{k34}) est une équation différentielle linéaire du premier ordre avec second membre. Sa résolution sans second membre donne :
\begin{equation}
	U=\frac{k}{\rho^2} \label{k35}
\end{equation}
En utilisant le second membre de (\ref{k34}), on considère que $k$ est une fonction de  $\varphi$, on a alors :
\begin{equation}
	U=\frac{1}{\rho^2}\left(k_0-\frac{C^2}{r^2}\right)=\frac{k_0r^2-C^2}{\rho^2r^2} \label{k36}
\end{equation}
avec $k_0$ la constante d'intégration. $U$ étant une fonction positive, on doit avoir :
\begin{equation}
	     k_0r^2 -C^2 > 0 \label{k37}
\end{equation}
En revenant à l'équation (\ref{k33}), on obtient :
\begin{equation}
	U=\left(\frac{d\varphi}{ds}\right)^2=\frac{k_0r^2-C^2}{\rho^2r^2} \label{k38}
\end{equation}
On utilise les équations (\ref{k30}) et (\ref{k38}), on obtient :
\begin{equation}
\left(\frac{d\varphi}{ds}\right)^2=\frac{k_0r^2-C^2}{\rho^2r^2}=	\left(\frac{C}{r^2Z}\right)^2=\frac{C^2}{r^4Z^2}=\frac{C^2}{r^4} \left(\frac{d\varphi}{d\lambda}\right)^2 \label{k39}
\end{equation}
ce qui donne :
\begin{equation}
\left(\frac{d\lambda}{d\varphi}\right)^2=\frac{\rho^2}{r^2}\frac{C^2}{k_0r^2-C^2} \label{k40}
\end{equation}
Pour déterminer la valeur de $k_0$, on exprime $\displaystyle \frac{d\lambda}{ds}$ en utilisant les équations (\ref{k26}) et (\ref{k40}). On écrit $ds^2$ :
$$ ds^2=\rho^2d\varphi^2+r^2d\lambda^2=\frac{r^2(k_0r^2-C^2)}{C^2}d\lambda^2+r^2d\lambda^2$$
soit:
\begin{equation}
	ds^2=\frac{r^4k_0}{C^2}d\lambda^2\Rightarrow \left(\frac{d\lambda}{ds}\right)^2=\frac{C^2}{k_0r^4}  \label{k41}
\end{equation}
Or d'après (\ref{k26}) : 
$$ \left(\frac{d\lambda}{ds}\right)^2=\frac{C^2}{r^4}  $$
d'où alors $k_0 = 1$ et par suite:
\begin{equation}
	\left(\frac{d\lambda}{d\varphi} \right)^2=\frac{\rho^2}{r^2}\frac{C^2}{r^2-C^2}  \label{k42}
\end{equation}
Pour pouvoir intégrer l'équation précédente, on exprime $r^2- C^2$, d'où :
\ba
r^2-C^2=N^2cos^2\varphi-C^2=\frac{a^2cos^2\varphi}{1-e^2sin^2\varphi}-C^2= \nonumber \\
\frac{(a^2-C^2)\left(1-\frac{a^2-C^2e^2}{a^2-C^2}sin^2\varphi\right)}{W}  \label{k43}	
\ea
On pose:
\begin{equation}
	k^2=\frac{a^2-C^2e^2}{a^2-C^2} \label{k44}
\end{equation}
D'où:
\begin{equation}
	r^2- C^2=(a^2-C^2)(1-k^2sin^2\varphi)/W  \label{k45}
\end{equation}
On remarque que  le coefficient $k$ est supérieur à 1, donc la latitude géodésique $\varphi$ reste inférieure à la latitude $\varphi_1$ définie par $sin\varphi_1 = 1/k$.

Alors l'équation (\ref{k42}) s'écrit :
\begin{equation}
	\left(\frac{d\lambda}{d\varphi}\right)^2=\frac{(1-e^2)^2C^2}{(a^2-C^2)cos^2\varphi(1-e^2sin^2\varphi)(1-k^2sin^2\varphi)}  \label{k46}
\end{equation}
D'où en remplaçant $C$ par $a.sin(Aze)$ et comme $tg(Aze)$ est de même signe que  $(d\lambda/d\varphi)$, on peut écrire alors :
\begin{equation}
	\frac{d\lambda}{d\varphi}=\frac{(1-e^2)tg(Aze)}{cos\varphi\sqrt{(1-e^2sin^2\varphi)(1-k^2sin^2\varphi)}}  \label{k47}
\end{equation}
Soit en intégrant entre 0 et $\varphi$: 
$$\lambda-\lambda_e=\int^{\varphi}_0\frac{(1-e^2)tg(Aze)}{cost\sqrt{(1-e^2sin^2t)(1-k^2sin^2t)}}dt= $$
$$(1-e^2)tg(Aze)\int^{\varphi}_0\frac{dt}{cost\sqrt{(1-e^2sin^2t)(1-k^2sin^2t)}} $$
ou encore : 
\begin{equation}
\lambda-\lambda_e=(1-e^2)tg(Aze)\int^{\varphi}_0\frac{dt}{cost\sqrt{(1-e^2sin^2t)(1-k^2sin^2t)}} \label{k48}
\end{equation}
En prenant comme variable $w=sint$, l'intégrale (\ref{k48}) devient:
\begin{equation}
\fbox{$	\lambda-\lambda_e=(1-e^2)tg(Aze)\ds \int^{sin \varphi}_0 \frac{dw}{(1-w^2)\sqrt{(1-e^2w^2)(1-k^2w^2)}} $}\label{k49}
\end{equation}
On cherche maintenant à exprimer l'abscisse curviligne $s$ en fonction de $\varphi$. Or l'expression de $ds^2$ est égale à :
$$ds^2=\rho^2d\varphi^2+r^2d\lambda^2=\rho^2d\varphi^2+\frac{C^2}{r^2}ds^2$$ 
soit:
\begin{equation}
	ds^2=\frac{r^2\rho^2d\varphi^2}{r^2-C^2}=\frac{a^2(1-e^2)^2cos^2\varphi d\varphi^2}{cos^2(Aze)(1-e^2sin^2\varphi)^3(1-k^2sin^2\varphi)}  \label{k50}
\end{equation}
D'où:
\begin{equation}
	ds=\frac{a(1-e^2)cos\varphi d\varphi}{cos(Aze)(1-e^2sin^2\varphi)\sqrt{(1-k^2sin^2\varphi)(1-e^2sin^2\varphi)}}  \label{k51}
\end{equation}
En prenant $t=sin\varphi$ comme nouvelle variable, l'intégrale de (\ref{k51}) donne en prenant comme origine de l'abscisse curviligne $s$ un point de l'équateur:
\begin{equation}
	\fbox{ $ s=\ds \frac{a(1-e^2)}{cosAze}\int^{sin\varphi}_0	\frac{dt}{(1-e^2t^2)\sqrt{(1-k^2t^2)(1-e^2t^2)}} $}  \label{k52}
\end{equation}
 Les intégrales (\ref{k49}) et (\ref{k52}) sont dites des intégrales elliptiques de troisième espèce.\index{Intégrale elliptique} 
\section{\textsc{Applications aux Problèmes Direct et Inverse  du Calcul des Lignes Géodésiques}}
Dans cette deuxième partie, on va traiter numériquement l'application des formules précédentes dans la résolution des problèmes dits respectivement direct et inverse du calcul des lignes géodésiques. 
\subsection*{5.6.1. Le Problème Direct}\index{Le problème direct}
On donne :

-	$(\varphi_1, \lambda_1)$ d'un point $M_1$;

-	la longueur $s$ de la géodésique de $M_1$ à $M_2$;

-	l'azimut géodésique $Az_1$ de la ligne  géodésique de  $M_1$ à  $M_2$.
\\

On demande de calculer :

-	les coordonnées géodésiques 	$(\varphi_2, \lambda_2)$ de $M_2$;

-	l'azimut géodésique $Az_2$ en $M_2$.
\\

\textbf{Solution:}
1. Calcul de la constante $C$, $C = N(\varphi_1).cos\varphi_1. sinAz_1 = a.sin(Aze)$ d'où $Aze$ et $k$.

2. Détermination de $\varphi_2$  à partir de :
$$ 	\Delta s=\frac{a(1-e^2)}{cosAze}	\frac{cos\varphi_1 \Delta \varphi}{(1-e^2sin^2\varphi_1)\sqrt{(1-k^2sin^2\varphi_1)(1-e^2sin^2\varphi_1)}}  
$$
avec $\Delta \varphi=\varphi_2-\varphi_1$.

3. Ayant $\varphi_2$, on calcule $\lambda_2$  par :
$$	\lambda_2-\lambda_1=(1-e^2)tg(Aze)\int^{sin \varphi_2}_{sin\varphi_1}\frac{dw}{(1-w^2)\sqrt{(1-e^2w^2)(1-k^2w^2)}} $$
4. Le calcul de $Az_2$ se fait par $sin(Az_2) = C/r(\varphi_2)$.
\subsection*{5.6.2. Le Problème Inverse}\index{Le problème inverse}
On donne les coordonnées $(\varphi_1, \lambda_1)$ et $(\varphi_2, \lambda_2)$ de deux points $M_1$ et $M_2$. On demande de calculer:

- la longueur $s$ de la ligne  géodésique de  $M_1$ à  $M_2$;

- l'azimut $Az_1$ en $M_1$;

-	l'azimut géodésique $Az_2$ en $M_2$.
\\

\textbf{Solution:}

1.	On doit calculer la constante $C$. A partir de l'équation (\ref{k42}), on peut écrire que: 
$$\left(\frac{\Delta \lambda}{\Delta \varphi} \right)^2=\frac{\rho^2(\varphi_1)}{r^2(\varphi_1)}\frac{C^2}{(r^2(\varphi_1)-C^2)}=\frac{(\lambda_2-\lambda_1)^2}{(\varphi_2-\varphi_1)^2}$$ ce qui donne $C$:
$$ C=\frac{\frac{r^2}{\rho}\frac{\Delta \lambda}{\Delta \varphi}}{\sqrt{1+\frac{r^2}{\rho^2}\left(\frac{\Delta \lambda}{\Delta \varphi}\right)^2}}$$
En considérant l'azimut compris entre 0 et $\pi$, donc $Az$ est positif, $C$ est positif. En le calculant pour  $\varphi_1$  et $\varphi_2$, on obtient $C$ par la valeur moyenne:$$ C=\frac{C_1(\varphi_1)+C_2(\varphi_2)}{2}$$

2. Par suite, on obtient la valeur de $k$ par (\ref{k44}):
$$k=\frac{a^2-C^2e^2}{a^2-C^2}$$
3. Ayant $C$, on a par (\ref{k27}), $Az_1$ et $Az_2$ :
$$ sin Az_1=\frac{C}{r(\varphi_1)}\quad  \mbox{et}\quad sinAz_2=\frac{C}{r(\varphi_2)} $$
4. Par suite, on a aussi $Az_e$:$$ sinAz_e=\frac{C}{a}$$
5. Enfin, l'équation (\ref{k52}) détermine $s$.

On itère le processus.
\subsection*{5.6.3. Calcul de l'Expression (\ref{k52})}
Dans ce paragraphe, on calcule en détail:
$$ 	s=\frac{a(1-e^2)}{cosAz_e}\int^{sin\varphi}_0	\frac{dt}{(1-e^2t^2)\sqrt{(1-k^2t^2)(1-e^2t^2)}} $$
Pour $|x| < 1$, on a les développements limités suivants:
\ba
	\frac{1}{(1+x)^{3/2}}=1-\frac{3}{2}x+\frac{15}{8}x^2-\frac{35}{16}x^3+\frac{315}{128}x^4+... \label{k53} \\
		\frac{1}{\sqrt{1-x}}=1+\frac{x}{2}+\frac{3x^2}{8}+\frac{5x^3}{16}+\frac{35x^4}{128}+... \label{k54} 
\ea
En prenant $x=-e^2t^2$ et $x=k^2t^2$, on obtient:
\ba
\frac{1}{(1-e^2t^2)^{3/2}}=1+\frac{3}{2}e^2t^2+\frac{15}{8}e^4t^4+\frac{35}{16}e^6t^6+\frac{315}{128}e^8t^8+... \label{k55} \nonumber \\
\frac{1}{\sqrt{1-k^2t^2}}=1+\frac{k^2t^2}{2}+\frac{3k^4t^4}{8}+\frac{5k^6t^6}{16}+\frac{35k^8t^8}{128}+... \label{k56} 
\ea
Par suite:
\ba
& \ds \frac{1}{(1-e^2t^2)\sqrt{(1-k^2t^2)(1-e^2t^2)}}=1+\frac{k^2+3e^2}{2}t^2+\frac{3k^4+6e^2k^2+15e^4}{8}t^4+ \nonumber & \\ & \ds \frac{5k^6+9k^4e^2+15k^2e^4+35e^6}{16}t^6+ \nonumber & \\ & \ds \frac{35k^8+60k^6e^2+90k^4e^4+140k^2e^6+315e^8}{128}t^8+...& \label{k57} 
\ea
ou encore à l'ordre 4 :
\be \ds \frac{1}{(1-e^2t^2)\sqrt{(1-k^2t^2)(1-e^2t^2)}}=1+mt^2+nt^4+... \lb{cgeod}
\ee
avec:
	\[m=\frac{k^2+3e^2}{2}; \quad n=\frac{3k^4+6e^2k^2+15e^4}{8}
\]
\subsection*{5.6.4. Calcul de l'expression (\ref{k49})}
On a:
$$  	\lambda-\lambda_e=(1-e^2)tg(Aze)\int^{sin \varphi}_0\frac{dw}{(1-w^2)\sqrt{(1-e^2w^2)(1-k^2w^2)}} $$ soit dans notre cas:
$$	\lambda_2-\lambda_1=(1-e^2)tg(Aze)\int^{sin \varphi_2}_{sin\varphi_1}\frac{dt}{(1-t^2)\sqrt{(1-e^2t^2)(1-k^2t^2)}} $$
Or d'après (\ref{k54}):
$$ \ds \frac{1}{\sqrt{1-e^2t^2}}=1+\frac{1}{2}e^2t^2+\frac{3}{8}e^4t^4+\frac{5}{16}e^6t^6+\frac{35}{128}e^8t^8+... $$
et :
$$\ds \frac{1}{\sqrt{1-k^2t^2}}=1+\frac{k^2t^2}{2}+\frac{3k^4t^4}{8}+\frac{5k^6t^6}{16}+\frac{35k^8t^8}{128}+... $$
et pour $(1-t^2)^{-1}$, on obtient:
	$$	\frac{1}{1-t^2}=1+t^2+t^4+t^6+t^8+...$$ 
D'où:
\ba
& \ds \frac{1}{(1-t^2)\sqrt{(1-e^2t^2)(1-k^2t^2)}}=1+\frac{2+k^2+e^2}{2}t^2+ \nonumber & \\ & \ds \frac{8+4k^2+4e^2+3k^4+2e^2k^2+3e^4}{8}t^4+\nonumber & \\ & \ds \frac{16+8k^2+8e^2+6k^4+4e^2k^2+6e^4+5k^6+3k^4e^2+3k^2e^4+5e^6}{16}t^6+...\nonumber &
\ea
Qu'on écrit sous la forme:
\begin{equation}
\frac{1}{(1-t^2)\sqrt{(1-e^2t^2)(1-k^2t^2)}}=1+\alpha t^2+ \beta t^4+\gamma t^6+...\label{k60} 
\end{equation}
avec:
\be
\left\{\begin{array}{l}
\alpha=\ds \frac{2+k^2+e^2}{2} \\
\beta=\ds \frac{8+4k^2+4e^2+3k^4+2e^2k^2+3e^4}{8} \\
\gamma=\ds \frac{16+8k^2+8e^2+6k^4+4e^2k^2+6e^4+5k^6+3k^4e^2+3k^2e^4+5e^6}{16} 
\end{array}\right.
\ee
\subsection*{5.6.5. Traitement d'un exemple}
\subsubsection*{Le Problème direct}
Soit le point $M_1$ avec:

- $\varphi_1=10.4549\,8299\,gr$;

- $\lambda_1=9.5954\, 2429\,gr$;

- $Az_1=249.3101\, 68\,gr$;

- $s=16255.206\,m$.

\textbf{Solution:}

- $C = N(\varphi_1).cos\varphi_1. sinAz_1 = -3^,594\,478.080\,m$;

- $Az_e=238.1131\, 471\,gr$;

- $k=\sqrt{\ds \frac{a^2-C^2e^2}{a^2-C^2}}=1.209227584$; 

- pour calculer $\varphi_2$, on pose $\Delta \varphi=\varphi_2 - \varphi_1$, et $s=\Delta s$, on a alors l'équation en utilisant (\ref{cgeod}):
	\[\frac{\Delta s.cosAz_e}{a(1-e^2)}=\ds \int_{sin\varphi_1}^{sin\varphi_2}\frac{dt}{(1-e^2t^2)\sqrt{(1-k^2t^2)(1-e^2t^2)}}=\int_{sin\varphi_1}^{sin\varphi_2}(1+mt+nt^2)dt
\]
A l'ordre 1, on a: $\ds \frac{\Delta s.cosAz_e}{a(1-e^2)}=sin\varphi_2-sin\varphi_1$.
\section{\textsc{Exercices et Problèmes}}
\bex
A partir de la définition géométrique de l'ellipse (\ref{p69}) donnée par: $$ MF + MF' = constante=2a $$ retrouver l'expression de l'équation cartésienne de l'ellipse.
\eex
\bex
 Montrer la formule très utilisée en géodésie:
$$ \frac{d(Ncos\varphi)}{d\varphi}=-\rho sin\varphi $$
 avec $N$ et $\rho$ les deux rayons de courbures principaux de l'ellipsoïde de révolution donnés respectivement par (\ref{p97}) et (\ref{p104}).
\eex
\bpb
A partir des équations de l'ellipsoïde de révolution:
$$M=\left\{\begin{array}{lll}
	X=Ncos\varphi cos\lambda  \\
	Y = Ncos\varphi sin\lambda  \\
	Z=N(1-e^2)sin\varphi   \end{array}\right.
	$$
1. Calculer les vecteurs:
$$ \displaystyle \frac{\partial \textbf{\textit{M}}}{\partial \lambda},\frac{\partial \textbf{\textit{M}}}{\partial \varphi}$$
2. Calculer les coefficients:
$$ E=\displaystyle \frac{\partial \textbf{\textit{M}}}{\partial \lambda}.\frac{\partial \textbf{\textit{M}}}{\partial \lambda}, \quad
	F=\displaystyle \frac{\partial \textbf{\textit{M}}}{\partial \lambda}.\frac{\partial \textbf{\textit{M}}}{\partial \varphi}, \quad
		G=\displaystyle \frac{\partial \textbf{\textit{M}}}{\partial \varphi}.\frac{\partial \textbf{\textit{M}}}{\partial \varphi}
   	$$
	Démontrer que l'expression de la première forme fondamentale s'écrit:
$$ ds^2=\rho^2d\varphi^2+N^2cos^2\varphi d\lambda^2 $$
\newpage
3. Calculer le vecteur normal $\textbf{\textit{n}}$ :
\begin{displaymath}
\textbf{\textit{n}}=\displaystyle \frac{\partial \textbf{\textit{M}}}{\partial \lambda} \wedge \frac{\partial \textbf{\textit{M}}}{\partial \varphi}\frac{1}{\left\|\displaystyle \frac{\partial \textbf{\textit{M}}}{\partial \lambda} \wedge \frac{\partial \textbf{\textit{M}}}{\partial \varphi}\right\|}
\end{displaymath}
4. Calculer les vecteurs:
$$ 	\displaystyle \frac{\partial^2 \textbf{\textit{M}}}{\partial \lambda^2},  \quad	\displaystyle \frac{\partial^2 \textbf{\textit{M}}}{\partial \lambda \partial \varphi}, \quad  \displaystyle \frac{\partial^2 \textbf{\textit{M}}}{\partial \varphi^2} $$
5. Déterminer les coefficients:
$$  	L=\displaystyle n.\frac{\partial^2 \textbf{\textit{M}}}{\partial \lambda^2}, \quad 
	M=\displaystyle n.\frac{\partial^2 \textbf{\textit{M}}}{\partial \lambda \partial \varphi}, \quad
		N=\displaystyle n.\frac{\partial^2 \textbf{\textit{M}}}{\partial ^2\varphi} 	$$
6. Ecrire la deuxième forme fondamentale $\Phi(\lambda,\varphi)$.

7. En appliquant la formule (\ref{corr}), Montrer que :
\begin{displaymath}
  N(\varphi)=\frac{a}{\sqrt{1-e^2sin^2\varphi}}
\end{displaymath}
est le rayon de courbure de la section normale au point $M$ perpendiculaire au plan de la méridienne de l'ellipsoïde de révolution.

8. En posant:
\begin{displaymath}
  d\m L= \frac{\rho d\varphi}{Ncos\varphi}
\end{displaymath}
En déduire que $ds^2$ s'écrit:
\begin{displaymath}
  ds^2=N^2cos^2\varphi (d\m L^2+d\lambda^2)
  \end{displaymath}
9. Montrer que $\mathcal L$ est donnée par:
\begin{displaymath}
\m L(\varphi)=Logtg\left(\ds \frac{\pi}{4}+\frac{\varphi}{2}\right)-\frac{e}{2}Log\left(\frac{1+esin\varphi}{1-esin\varphi}\right)
\end{displaymath}
\epb
\bpb
 Sur l'ellipsoïde, on note $\varphi$  la latitude géodésique et  $\psi$   la latitude réduite.

1. Calculer $\rho$  le  rayon de  courbure de l'ellipse méridienne en fonction de $\psi$.

2. Exprimer l'aplatissement de l'ellipsoïde en fonction des valeurs de $\rho$  au pôle et à l'équateur. 

3. On mesure la longueur d'un arc de méridien d'un degré à la fois au pôle et à l'équateur. On trouve respectivement $111\,695\, m$ et $110\,573\, m$. En déduire l'aplatissement.
\epb
\bpb
 On donne les coordonnées tridimensionnelles suivantes d'un point $M$:
	\[	M=(X,Y,Z)=(4\,300\,244.860\, m,1\,062\,094.681\, m,4\,574\,775.629\, m)
\]
Les paramètres de l'ellipsoïde de référence sont $a  = 6\,378\,137.00\, m,\quad  e^2 = 0.006\,694\,38$.
 
1.	Calculer le demi-petit axe $b$.

2.	Calculer l'aplatissement. 

3.	Calculer les coordonnées géodésiques $(\varphi ,\lambda , he)$ du point $M$. $\varphi$   et $\lambda$    seront calculées en grades avec  cinq chiffres après la virgule.
\epb
\bpb
Soit $\m E(a,e)$ un ellipsoïde de révolution où $a,e$ sont respectivement le demi-grand axe et la première excentricité. $(g)$ une géodésique partant d'un point $E(\varphi=0,\lambda_E)$ sur l'équateur et d'azimut $Az_E$. A cette géodésique, on lui fait correspondre une géodésique $(g')$ sur la sphère $\m S^2$ dite de Jacobi\footnote{\textbf{Carl Gustav Jacob Jacobi} (1804-1851): mathématicien allemand.} \index{\textbf{Jacobi C.G.J.}} de rayon $a$, ayant le même azimut $Az_E$ au point $E'(\varphi'=0,\lambda_E)$. De même au point $M(\varphi,\lambda)$ de la géodésique $(g)$ de l'ellipsoïde, on lui fait correspondre le point $M'(\varphi',\lambda')$ de $(g')$ de $\m S^2$ tel qu'il y a conservation des azimuts. 
\begin{figure}[h]
	\centering
		\includegraphics[width=0.70\textwidth]{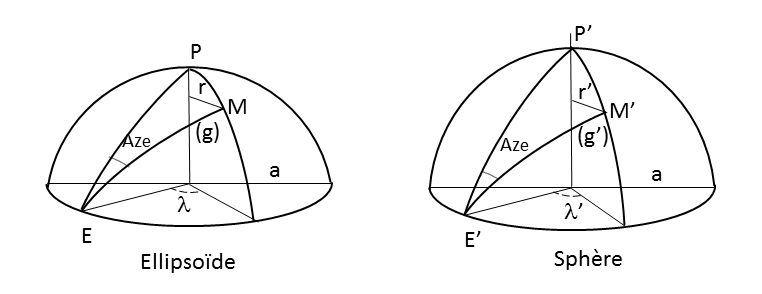}
	\caption{La correspondance de la sphère de Jacobi}
	\label{fig:sjacobi}
\end{figure}

1. Ecrire l'équation de Clairaut pour la géodésique $(g)$. 

2. On note $r'$ le rayon du parallèle passant par $M'$ de la géodésique $(g')$. Ecrire de même l'équation de Clairaut pour la géodésique $(g')$. 

3. Montrer que $\varphi$ et $\varphi'$ vérifient:
$$ Ncos\varphi=acos\varphi'$$
et en déduire que $\varphi'$ est la latitude paramétrique de $M$.

4. Ecrire les expressions de $tgAz_g$ et $tgAz_{g'}$ respectivement sur $(g)$ et $(g')$.

5. Montrer que:
$$  d\lambda=\frac{\rho d\varphi}{ad\varphi'}d\lambda'$$ 
En déduire que : 
$$ 	d\lambda=\sqrt{1-e^2cos^2\varphi'}d\lambda' $$
6. En intégrant l'équation précédente, montrer qu'on obtient:
$$	\lambda - \lambda_E=\int_{\lambda_E}^{\lambda'+\lambda_E}\sqrt{1-e^2cos^2\varphi'}d\lambda' $$
avec $\lambda>\lambda_E$ et $\lambda'$ est comptée à partir de $\lambda_E$.

7. En écrivant $ \sqrt{1-e^2cos^2\varphi'}=1-\frac{e^2}{2}cos^2\varphi'+o(e^4)$ où $o(e^4)$ est un infiniment petit d'ordre 4 en $e$ dont on néglige, écrire l'intégrale précédente entre $\lambda_E$ et $\lambda_E+\lambda$.

8. Comme $(g')$ est une géodésique de la sphère, on démontre que:
$$ 	cos^2\varphi'd\lambda'=\frac{sinAz_E}{a}ds'  $$
où $ds'$ est l'élément différentiel de l'abscisse curviligne sur la géodésique (un grand cercle). Alors en posant $s'=0$ au point $E'$, montrer que l'équation précédente s'écrit sous la forme:
$$ 		\lambda =\lambda_E+\lambda'-\frac{e^2sinAz_e}{2a}\int_0^{s'}ds' $$ 
9. On suppose que la géodésique $(g')$ coupe une première fois le plan de l'équateur en un point $F'$, montrer qu'on obtient:
\ba
	\lambda'_F=\pi \label{27}\nonumber \\
	s'=\pi a \label{28}  \nonumber \\
	\lambda_F=\lambda_E+\pi-\frac{e^2\pi sinAz_E}{2}  \nonumber
\ea
10. La géodésique $(g')$ partant de $F'$ a pour azimut $\pi-Az_E$, elle coupe une deuxième fois l'équateur au point $E'$, mais la géodésique $(g)$ sur l'ellipsoïde coupe une deuxième fois le plan de l'équateur au point correspondant à $H$ dont la longitude est $\lambda_H$. Montrer que $\lambda_H$ est donnée par:
$$	\lambda_H=\lambda_E+2\pi-\frac{e^2\pi sinAz_E}{2}-\frac{e^2\pi sin(\pi-Az_E)}{2}=\lambda_E+2\pi-e^2\pi sinAz_E $$
Quelle conclusion a-t-on sur les lignes géodésiques de l'ellipsoïde de révolution.
  \epb
\bpb
Un point $M$ de la surface d'une sphère $(S)$ de rayon $R$, a pour coordonnées $(X,Y,Z)$ dans un repère orthonormé: 
$$	  M=(X,Y,Z)=(Rcos\varphi .cos\lambda ,Rcos\varphi .sin \lambda ,Rsin\varphi)$$
1. Montrer qu'un vecteur normal unitaire $n$ à $(S)$ en $M$ est:
$$ n=	 (cos\varphi .cos\lambda,	 cos\varphi .sin \lambda,	 sin\varphi)^T$$
2. Soit $(C)$ le grand cercle passant par le point $A(R,0,0)$ et d'azimut $Az_E$. Le point $M$ peut être décrit par son abscisse curviligne $s$ mesurant l'arc $AM$. On note par $\omega$ représente l'angle au centre de l'arc $AM$. Utilisant la trigonométrie sphérique, montrer que:
	$$cos\varphi.sin\lambda= sin\omega.sinAz_E $$
3. En utilisant la formule fondamentale de la trigonométrie sphérique dans le triangle APM, montrer qu'on a les deux  relations :
$$\begin{array}{l}
	      cos\omega =  cos\varphi.cos\lambda \\
	                  sin\varphi = sin\omega.cosAz_E
										\end{array}$$
4. En déduire que les coordonnées de $M$ s'écrivent en fonction de $s$ comme suit:
$$ M\left\{
\begin{array}{l}
	 X= R.cos(s/R)\\
	 Y=RsinAz_Esin(s/R)\\
	 Z=RcosAz_Esin(s/R)
\end{array}\right.$$
5. Calculer les vecteurs $T$ et $N$ du repère de Frenêt. En déduire les composantes de $N$ en fonction de $\omega$.

6. Montrer que les vecteurs $N$ et $n$ sont parallèles.

7. Justifier que les géodésiques de la sphère sont les grands cercles. 
\epb
\bpb
Soit le tore $\BbT$ défini par les équations suivantes:
$$M(\fii,\lm)=\left\{\begin{array}{l}
x=(a+Rcos\fii)cos\lm \\
y=(a+Rcos\fii)sin\lm \\
z=Rsin\fii 
\end{array}\right.
$$
où $a,R$ deux constantes positives avec $a>R$, $(\fii,\lm)\in \,[0,2\pi]\times [0,2\pi]$. 

1. Calculer la première forme fondamentale $ds^2$.

2. Avec les notations usuelles, on pose: 
$$\frac{\partial E}{\partial \fii}=E'_\fii,\quad	\frac{\partial E}{\partial \lm}  =E'_\lm,\quad 	\frac{\partial F}{\partial \fii} =F'_\fii$$
 $$\quad \frac{\partial F}{\partial \lm} =F'_\lm,\quad \frac{\partial G}{\partial \fii}=G'_\fii,\quad	\frac{\partial G}{\partial \lm}= G'_\lm $$
Utilisant les équations des géodésiques (\ref{k15-1}) et (\ref{k15-2}) du cours, montrer que les équations des géodésiques du tore sont:
\[-2Rsin\fii(a+Rcos\fii)	\frac{d\fii}{ds} \frac{d\lm}{ds}+(a+Rcos\fii)^2\frac{d^2 \lm}{ds^2}=0 \]
\[Rsin\fii(a+Rcos\fii)\left(\frac{d\lm}{ds}\right)^2+R^2\frac{d^2 \fii}{ds^2}=0 \]	
3.  Montrer que la première équation ci-dessus donne: $$ \ds (a+Rcos\fii)^2\frac{d\lm}{ds}=C=cte$$
Montrer qu'on retrouve l'équation de Clairaut avec $C=(a+R)sinAze$ où $Aze$ est l'azimut de départ au point $M_0(\fii=0,\lm_0)$.

4.  On suppose au point $M_0$, la géodésique a pour azimut $Aze$ tel que:
 $$ 0 < Aze< \ds \frac{\pi}{2}$$
 Montrer que la deuxième équation des géodésiques s'écrit en utilisant le résultat précédent:
	\[  \frac{d^2 \fii}{ds^2}=-\frac{C^2}{R}\frac{sin\fii}{(a+Rcos\fii)^3}
\]
5.  Montrer qu'on arrive à: 	$$  \ds \left(\frac{d\fii}{ds}\right)^2=l-\frac{C^2}{R^2(a+Rcos\fii)^2}\geq 0$$
 où $l$ est une constante d'intégration.
\epb 

\chapter{\textit{\textbf{Les Systèmes Géodésiques}}}

Parmi les buts de la géodésie, on trouve la définition et la mise en place des systèmes géodésiques.
\\

A un système géodésique\index{Système géodésique}, on lui associe le réseau géodésique de base. On verra par la suite, l'établissement et le calcul des réseaux géodésiques.
\section{\textsc{Définition d'un système géodésique}}
\bdf
Un système géodésique donné est un système de coordonnées où sont représentés les points géodésiques. Il est défini par:

a - son origine;

b - son orientation;

c - le type de coordonnées choisies pour localiser les points. 
\edf
Le système le plus utilisé est le système cartésien\index{Système cartésien} formé par un repère $(OX,OY,OZ)$ tel que l'axe $OZ$ soit parallèle  à l'axe de rotation de la Terre, et le plan $OXZ$ parallèle au méridien de Greenwich origine des longitudes, l'axe $OY$  est tel que le trièdre $(OX,OY,OZ)$ soit orthogonal et direct (\textbf{Fig. \ref{fig:repere}}). A ce système, on lui associe une base orthonormée $(\textbf{\textit{e}}_1, \textbf{\textit{e}}_2, \textbf{\textit{e}}_3)$ c'est-à-dire :
	\[               \|\textbf{\textit{e}}_1\|=\|\textbf{\textit{e}}_2\|=\|\textbf{\textit{e}}_3\|= 1\, \mbox{mètre (l'unité des longueurs)}
\]

\begin{figure}[h]
	\centering
		\includegraphics[width=0.70\textwidth]{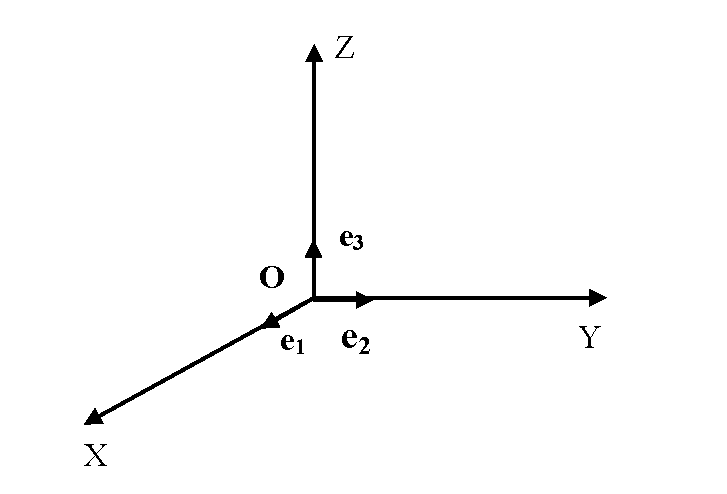}
	\caption{Le repère cartésien}
	\label{fig:repere}
\end{figure}

Pour les systèmes géodésiques classiques (terrestres), la position de l'origine est à 500 m environ du centre des masses de la Terre.
\\

Pour les systèmes géodésiques établis par la géodésie spatiale actuelle (comme le GPS - \textbf{G}lobal \textbf{P}ositioning \textbf{S}ystem), l'origine est presque confondue avec le centre des masses de la Terre (<2 m).
\\

L'orientation du système géodésique classique est faite à partir des observations astronomiques sur les étoiles. Ces observations vont orienter l'axe $OZ$ et le plan $OXZ$ du système à être respectivement parallèle à l'axe moyen de la rotation de la Terre et au méridien de Greenwich.
\\

Un système géodésique ou référentiel géodésique ou datum géodésique \index{Datum géodésique} obéit à certaines conditions à savoir:

- pas de déformation d'échelle;

- une meilleure distribution des points;

- la qualité homogène des coordonnées des points.
\\

En général, les référentiels géodésiques nationaux ne remplissent pas toujours ces conditions.
\\

Pour la mise en place d'un système géodésique, on adopte un modèle proche de la Terre. Un premier modèle est une sphère de rayon $R$  le rayon moyen de la terre. Dans ce cas, les coordoonées tridimensionnelles $(X,Y,Z)$ d'un point $M$ s'expriment par :
\be
M\left\{\begin{array}{l}
X=Rcos\varphi.cos\lambda \\
Y=Rcos\varphi.sin\lambda\\
Z=Rsin\varphi
\end{array}\right.
\ee
avec respectivement:

- $\varphi$ la latitude;

- $\lambda$ la longitude comptée à partir d'un méridien origine, positive vers l'Est;

- $h$ l'altitude au dessus de la surface.

Le deuxième modèle le plus approprié pour la Terre, après avoir fait de mesures, est l'ellipsoïde de révolution. Ainsi à chaque système géodésique est associé son ellipsoïde de révolution dit l'ellipsoïde de référence.  
\\

Un autre problème avec les systèmes géodésiques classiques est qu'il y a deux systèmes indépendants: l'un pour les coordonnées horizontales et un autre pour la composante verticale.
\\

Les réseaux planimétriques ou horizontaux sont déterminés à partir des observations de triangulation\index{Triangulation} (mesures angulaires) en général ou de trilatération\index{Trilatération} (mesures des distances) réduites à l'ellipsoïde adopté.
\\

Par contre, le système altimétrique\index{Système altimétrique} est observé par le nivellement de précision et la référence des altitudes est déterminée à partir des observations du niveau moyen des mers à l'aide des marégraphes\index{Marégraphe}.
\\
 
A un système donné de coordonnées planes, par exemple les coordonnées Lambert\footnote{La représentation conique conforme présentée par \textbf{Johann Heinrich Lambert} (1728-1777): mathématicien, physicien et astronome suisse.}\index{\textbf{Lambert J.H.}} $(x,y)$, ou des coordonnées tridimensionnelles $(X,Y,Z)$ géodésiques (par géodésie classique ou par les techniques spatiales), elles sont associées alors à un référentiel ou datum géodésique.
\\

 
\section{\textsc{Le Géoïde}}
Soit le repère $OXYZ$ et une masse ponctuelle $m'$ au point $O$ et soit un point $M (X,Y,Z)$ de masse ponctuelle $m$ (\textbf{Fig. \ref{p11}}).
\begin{figure} [htp]
\centering
\includegraphics{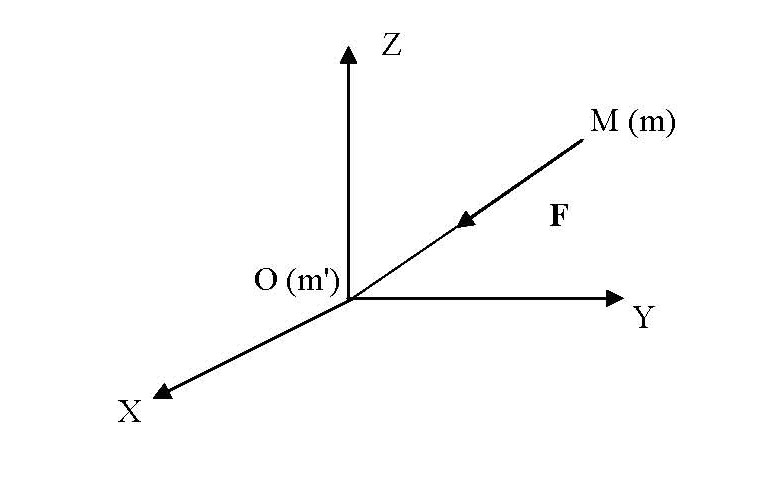}
\caption{Le repère 3D}
\label{fig:p11}
\end{figure}

Alors, le point $M$ est soumis  à une force \textbf{\textit{F}} due à l'attraction de la masse $m'$ au point $O$. Le module de cette force est : 
 \be
	\textit{F}=\frac{Gmm'}{r^2}= \textit{F}\left(X,Y,Z\right) \label{p127} 
\ee
où $G$ est la constante universelle de gravitation et $r$ est la distance $OM$.

\subsection*{6.2.1. Le Champ du Potentiel}
\bdf
On appelle champ du potentiel\index{Champ de potentiel} la fonction scalaire  $\textit{V}$ définie par \index{Datum géodésique}:
\be
V=\frac{Gmm'}{r}=V\left(X,Y,Z\right) \label{p128} 
\ee
\edf
   \subsection*{6.2.2. Gradient}
\bdf
On appelle gradient \index{Gradient} d'une fonction scalaire $U(X,Y,Z)$ le vecteur noté $\textbf{\textit{grad}}U$ et  de composantes :  
 \be
	\textbf{\textit{grad}}U = \begin{pmatrix}{
\displaystyle	\frac{\partial U}{\partial x} \cr
\nonumber \cr
\displaystyle \frac{\partial U}{\partial y}\cr
\nonumber \cr
\displaystyle \frac{\partial U}{\partial z} }
\end{pmatrix} \label{p129} 
\ee
\edf
\textit{\textbf{Exemple 1:}}     $U = X^2+Y^2+Z^2$, $\textbf{\textit{grad}}\textit{U}$ est le vecteur de composantes : 
\be
	\textbf{\textit{grad}}U = (2X,2Y,2Z)^T \label{p130} 
\ee
 où \textit{T} désigne transposé. 
                                                       
\textit{\textbf{Exemple 2:} } 
\be
U = \frac{1}{r}	\label{p131} 
\ee
 comme  $r^2  = X^2 + Y^2 + Z^2\,\Longrightarrow \,2rdr = 2XdX + 2YdY + 2ZdZ$
\be
	\textbf{\textit{grad}}U = \left(\frac{-X}{r{^3}},\frac{-Y}{r{^3}},\frac{-Z}{r{^3}}\right) ^T \label{p132} 
\ee
Si on pose :                                                       
\be
	\textbf{\textit{r}} = \textbf{\textit{OM}} = X\textbf{\textit{i}}  + Y\textbf{\textit{j}} + Z\textbf{\textit{k}} \label{p133} 
\ee
alors :   
\be
	\textbf{\textit{F}}=-\frac{\textbf{\textit{r}}}{r^{3}} \label{p134} 
\ee
On calcule maintenant le gradient de la fonction scalaire donnée par l'équation (\ref{p128}) c'est-à-dire le champ du potentiel. En utilisant l'exemple 2, on a :
\be
	\textbf{\textit{grad}}V=\textbf{\textit{grad}}\left(\frac{Gmm'}{r}\right)=Gmm'\textbf{\textit{grad}}\left(\frac{1}{r} \right) = -Gmm'\frac{\textbf{\textit{r}}}{r^{3}} \label{p135} 
\ee
On remarque si on pose :  
\be
	\textbf{\textit{n}}=\frac{\textbf{\textit{r}}}{r} \label{p136} 
\ee
$ \textbf{\textit{n}}$ est alors un vecteur unitaire porté par $\textbf{\textit{OM}}$ et dans la direction $\textbf{\textit{OM}}$.                                                                               
L'expression de la force $\textbf{\textit{F}}$ s'écrit :     
\be
	\textbf{\textit{F}}=-F\textbf{\textsl{n}} = - \frac{Gmm'}{r^2}\textbf{\textit{n}} \label{p137} 
\ee
\be
	\textbf{\textit{grad}} V=- \frac{Gmm'}{r^2}\textbf{\textit{n}} \label{p138} 
\ee
D'où :
\be
	\fbox{ $ \textbf{\textit{F}} = \textbf{\textit{grad}}V $} \label{p139}
\ee
On dit que la force \textbf{\textit{F}} dérive du champ de potentiel $V$.

\subsection*{6.2.3. Le Laplacien\footnote{Nommé à l'honneur de \textbf{Pierre Simon de Laplace} (1749-1827): mathématicien français.}}
\bdf
On \index{\textbf{Laplace P.S.}}appelle laplacien \index{Le laplacien} d'une fonction scalaire $U(X,Y,Z)$ la fonction scalaire notée $\Delta U$ définie par:
\be
\fbox{ $ \Delta U(X,Y,Z)=\ds \frac{\partial ^2U}{\partial X^2}+\frac{\partial ^2U}{\partial Y^2}+\frac{\partial ^2U}{\partial Z^2} $}
\ee
\edf
\subsection*{6.2.4. Le Champ Réel ou Champ du Potentiel de la Pesanteur}
Soit le repère $OXYZ$ tel que $O$ soit le centre de gravité de la Terre et $OZ$ son axe de rotation. Le plan $OXY$ contient le méridien de Greenwich (\textbf{Fig. \ref{p2}}). Un point $M(X,Y,Z)$ de masse unité est soumis au potentiel $V$ de gravitation et au potentiel $\Phi $ de la force centrifuge due à la rotation de la Terre.

L'expression de $V$ est :   
\be
	V=G\int \!\! \int \!\! \int_{Terre}\frac{dm}{r}  \label{p140}
\ee

\begin{figure} 
\centering
\includegraphics{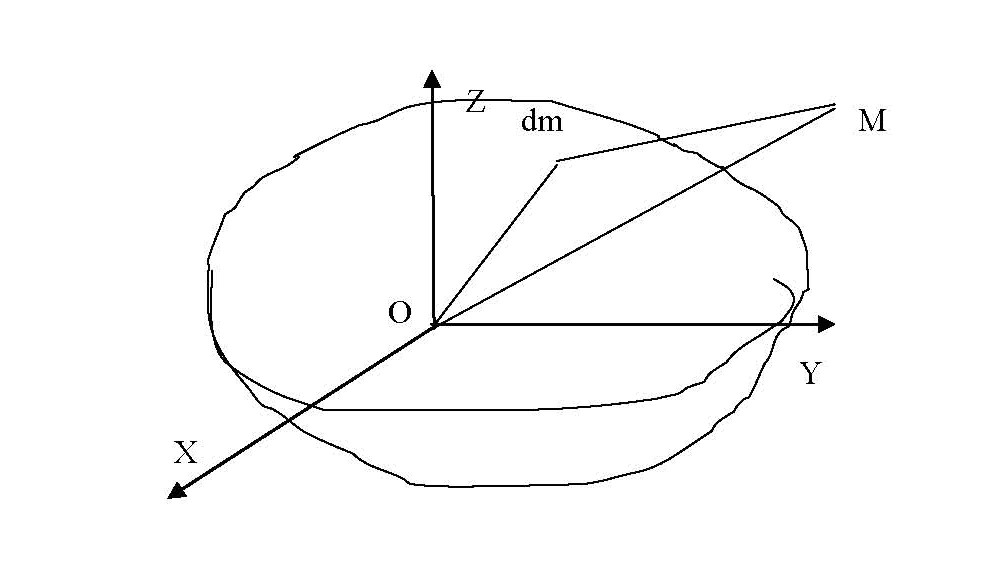}
\caption{Le potentiel}
\label{fig:p2}
\end{figure}

Malheureusement, cette expression n'est pas calculable car nous ignorons la distribution des masses à l'intérieur de la Terre.
\\

L'expression du potentiel $ \Phi $ de la force centrifuge est donnée par :
 \be
	\fbox{$\Phi=\ds \frac{1}{2}\left(x^2+y^2\right)\Omega ^2 $} \label{p141}
\ee
où $\Omega$  est la vitesse de la rotation de la Terre.
\\

\bdf
On appelle Potentiel du champ réel $W$ ou potentiel de la pesanteur\index{Potentiel de la pesanteur} la somme du potentiel $V$ et $\Phi$:
\be
	\fbox{$ W  = V + \Phi $}   \label{p142}
\ee
\edf
\bdf
On appelle vecteur de gravité\index{Vecteur de gravité} le vecteur \textbf{\textit{g}} tel que : 
 \be
		\fbox{$ \textbf{\textit{g}}  = \textbf{\textit{grad}}\textit{W} $} \label{p143}
\ee
\edf
  $ g $  mesure la gravité ou la pesanteur, a la dimension d'une accélération et exprimée en $m/s^2$ (Unité Système International) ou en $cm/s^2$ ($1cm/s^2 = 1 gal$ \index{Gal} en hommage à Galilée\index{\textbf{Galilée G.}}). $g$ mesure 978 $gals$ à l'équateur et 983 $gals$ aux pôles.\\
	
\bdf
Les surfaces $W(X,Y,Z) = W_0$ = constante, sur lesquelles le potentiel $W$ est constant sont appelées surfaces équipotentielles\index{Surface équipotentielle} ou surfaces de niveau\index{Surface de niveau}.
\edf
En différentiant le potentiel $W = W(X,Y,Z)$, on obtient :
\be
	d W =\frac{\partial W}{\partial x}dX+\frac{\partial W}{\partial y}dY+\frac{\partial W}{\partial z}dZ  \label{p144}
\ee
ou en notation vectorielle:   
\be
	dW = \textbf{\textit{grad}}W.\textbf{\textit{dM}} = \textbf{\textit{g}}.\textbf{\textit{dM}}  \label{p145}
\ee
Si $\textbf{\textit{dM}}$ est pris sur la surface équipotentielle $W = W_0$:
\be
	dW = 0 \Longrightarrow   \textbf{\textit{g.dM}} =0 \Longrightarrow       \textbf{\textit{g}}  \bot  \textbf{\textit{dM}}  \label{p146}
\ee
alors (\ref{p146}) exprime que le vecteur \textbf{\textit{g}} est normal à la surface équipotentielle passant par le même point.
\subsection*{6.2.5. Le Géoïde}
\index{Géoïde}
Une première définition du \textbf{géoïde} est due à C.F. Gauss \index{\textbf{Gauss C.F.}} (\textit{W.A. Heiskanen \& H. Moritz}, 1967): 
\\\index{\textbf{Heiskanen W.A.}}\index{\textbf{Moritz H.}}

'' \textit{Ce que nous appelons dans le sens géométrique la surface de la terre ce n'est que la surface qui coupe les lignes de la pesanteur sous un angle droit et qui fait partie de la surface des océans}''.
\\

Le terme géoïde fut introduit pour la première fois par J. Listing\footnote{\textbf{Johann Benedict Listing} (1808-1882): mathématicien prussien (élève de C.F. Gauss).} en 1872 (\textit{M. Bur$\breve{s}$a \& K. P$\breve{e}\breve{c}$}, 1986):'' \textit{nous appellerons la surface mathématique de la terre définie précédemment la surface à laquelle les océans font partie, surface géodale de la terre ou géoïde}''.\index{\textbf{Listing J.B.}}\index{\textbf{Bur$\breve{s}$a M.}} \index{\textbf{P$\breve{e}\breve{c}$ K.}}
\\

La surface des liquides et des fluides se met en équilibre perpendiculairement à la verticale. Si on considérait un ensemble fluide recouvrant toute la Terre, il définissait donc une surface de niveau de la pesanteur.
\\

La surface moyenne du niveau des mers, abstraction faite des marées et corrigée des variations, définissait une surface de niveau unique pour le monde entier. On peut définir d'autres surfaces de niveau de proche en proche à partir d'un point quelconque pris comme origine sur une verticale donnée, par la condition que cette surface soit en tout point perpendiculaire à toutes les autres verticales.
\\

Ces surfaces de niveau successives que l'on peut numéroter sont des surfaces fermées qui s'enveloppent les unes les autres, un autre exemple ce sont les surfaces de même pression atmosphérique (isobares) sont théoriquement des surfaces de niveau.
\\

Comme par définition, une surface de niveau est normale aux lignes de forces et que ces lignes sont des courbes gauches, les surfaces de niveau ne sont pas parallèles entr'elles, c'est-à-dire que la distance de deux surfaces de niveau n'est pas constante, ce qui reste constant c'est le travail qu'il faut accomplir contre la pesanteur pour déplacer une masse d'un point donné de ces surfaces à un point quelconque d'une autre de celles-ci.
\bdf \lb{defgeoid}
On appelle géoïde la surface de niveau qui coïncide avec la surface moyenne des mers et qui se prolonge sous les continents par la condition d'y rester normale à toutes les lignes de forces.
\edf
On peut dans ces conditions considérer que la Terre est consituée par le géoïde, surmontée du relief dont l'altitude au dessus du niveau moyen de la mer sera par définition égale à la distance qui le sépare du géoïde. L'expérience prouve que le géoïde s'écarte très peu d'un ellipsoïde de révolution: parceque le géoïde a une expression mathématique très compliquée, alors nous utilisons en géodésie comme surface mathématique du géoïde celle de l'ellipsoïde de révolution.  
\\

Le géoïde est donc une surface physique dont la modélisation mathématique est assez compliquée.
\\

L'un des buts de la géodésie est de déterminer la position de cette surface par rapport à la surface topographique.
\\

Généralement, l'origine des réseaux du nivellement de précision est déterminée à partir des mesures du niveau moyen des mers enregistrées par un marégraphe\index{Marégraphe}. Alors, on a la relation suivante entre l'altitude du nivellement\index{Altitude de nivellement} $(H)$ et l'altitude ellipsoïdale $(he)$: 
\be
\fbox{ $	he = H + N $} \label{p150}
\ee
où $N$ désigne la hauteur du géoïde\index{Hauteur du géoïde} par rapport à l'ellipsoïde de référence ou ondulation du géoïde\index{Ondulation du géoïde} (à ne pas confondre avec la grande normale donnée par l'équation (\ref{p97})).
\section{\textsc{Les Systèmes de Coordonnées}}
A chaque système géodésique, on lui associe un système de coordonnées\index{Système de coordonnées} avec lesquelles seront exprimées les positions des points géodésiques.
\subsection*{6.3.1. Les Coordonnées Sphériques}\index{Coordonnées sphériques}
Un point $M$ est défini par le triplet $(r, \lambda, \theta)$. Les coordonnées de $M$ s'expriment ainsi:
\be
M=\left\{\begin{array}{lll}
	X =rsin\theta cos\lambda \\ 
	 Y=rsin\theta sin\lambda	\\
	  Z=rcos\theta 
\end{array}\right.  \label{p151}
\ee
avec:

$r$: la distance géocentrique,
	
$\lambda$: la longitude,
	
$\theta$: le complément à la latitude géocentrique. 
\subsection*{6.3.2. Les Coordonnées Géodésiques}
Pour le modèle ellipsoïdique de la Terre, les coordonnées $(X, Y, Z)$ d'un point $M$ s'expriment par les formules:
\be
M=\left\{\begin{array}{lll}
	X =(N+he)cos\varphi cos\lambda \\
	 Y=(N+he)cos\varphi sin\lambda	\\ 
	 Z=(N(1-e^2)+he)sin\varphi
\end{array}\right.   \label{p152}
\ee
avec:
$$	N=\frac{a}{\sqrt{1-e^2sin^2\varphi}} $$
$a$, $e$ et $he$ sont respectivement le demi-grand axe, la première excentricité de l'ellipsoïde de référence et l'altitude ellipsoïdique au point concerné.
\\

Ces coordonnées sont dites des coordonnées géométriques.
\subsection*{6.3.3. Les Coordonnées Naturelles ou Géoidales ou physiques}
\bdf
On appelle coordonnées naturelles, géoidales ou physiques le triplet $(\Phi, \Lambda, H)$ avec $(\Phi, \Lambda)$ sont respectivement la latitude et la longitude astronomiques. \index{Coordonnées géoidales}
\edf
Ces coordonnées se rapportent à la verticale du lieu\index{Verticale du lieu} (définie par la direction donnée par un fil à plomb). La notion de la verticale est une notion fondamentale car elle correspond à une mesure physique. Elle n'est pas à confondre avec celle de la normale à l'ellipsoïde.
\subsection*{6.3.4. La D\'eviation de la verticale}
Si en un point donné, on a mesuré $(\Phi, \Lambda, H)$ et on a déterminé $(\varphi, \lambda, he)$, généralement on a:
\ba
	\Phi \neq \varphi  \label{p153}  \\
	\Lambda \neq \lambda  \label{p154}\\
	H \neq he  \label{p155}
\ea
On dit qu'on a une déviation de la verticale  en ce point.
\bdf
On définit les composantes de la déviation de la verticale par\index{D\'eviation de la verticale}:
\be
\fbox{ $
\begin{array}{l}
 \zeta=\Phi-\varphi  \label{p156} \\
 \eta = (\Lambda-\lambda)cos\varphi  
\end{array} $} \label{p157}
\ee
$\zeta$ et $\eta$ sont respectivement les composantes Nord et Est de la déviation de la verticale.
\edf

L'objet principal de la Géodésie est la détermination en chaque point de la Terre de $\zeta$, $\eta$ et $N=he - H$.
\newpage
\section{\textsc{Les Systèmes Géodésiques en Afrique du Nord}}
Les systèmes géodésiques terrestres de la Tunisie seront présentés au chapitre consacré à la Géodésie Tunisienne [§ \ref{geotun}].
\subsection*{6.4.1. Les Systèmes Géodésiques en Algérie}
\textbf{Le Système Voirol}\index{Système Voirol}

        C'était  le premier  système défini en Algérie (1875):
				
- le point fondamental (point de départ) : Voirol (près d'Alger) créé en 1875;

- la surface  de référence : c'est-à-dire le modèle choisi pour la Terre est l'ellipsoïde de Clarke\footnote{Voir note historique ci-après.} Français 1880 $(a\,=\,6\,378\,249.20\, m,\, e^2\,=\,0.0068034877)$ avec $a$ le demi-grand axe de l'ellipsoïde et $e$ la première excentricité;\index{Ellipsoïde de Clarke Français}

- l'orientation de départ est l'azimut astronomique de la direction Voirol-Meleb El Kora mesuré en 1874;

- la mise à l'échelle ou la qualité métrique de réseau : la mesure d'une distance ou base à Blida en Algérie mesurée en 1854.

\textbf{Le Système Europe 50}\index{Système Europe 50}
      
			Le système était mis en place par les Américains pour le besoin de l'OTAN(Organisation du Traité de l'Atlantique Nord)en faisant une compensation générale des réseaux de l'Europe Occidentale en associant les chaînes de triangles de l'Afrique du Nord. Il est défini par les éléments suivants:
			
- le point fondamental: Potsdam (Allemagne);

-	l'azimut d'orientation : l'azimut astronomique de la direction Potsdam-Golmberg;

-	l'ellipsoïde de référence : l'ellipsoïde international de Hayford:\footnote{\textbf{John Fillmore Hayford} (1868-1925): géodésien américain.} $$ a\,=\,  6\,378\,388.000\, m,\quad f\,=\,1/297.000 $$\index{\textbf{Hayford J.F.}}
\textbf{Le Système BT56}\index{Système BT56}

Les éléments de définition de ce système géodésique terrestre sont:

-	le point fondamental : point Bouzaréah;

-	l'ellipsoïde de référence : l'ellipsoïde de Clarke 1880 Anglais: $$a\,=\, 6\,378\,249.145\,m,\quad e^2\, =\, 0.0068\,0351\,128$$

-	l'orientation est assurée par trois gisements astronomiques traités en équations de condition.

\textbf{Le Système Nord Sahara}\index{Système Nord Sahara}
       
			Ce système a été défini par une transformation conforme des coordonnées du système Europe 50 de l'ellipsoïde international de Hayford sur l'ellipsoïde de Clarke 1880 Anglais.\index{Ellipsoïde de Clarke Anglais}
\subsection*{6.4.2. Les Systèmes Géodésiques en Libye}
\textbf{Le Système Europe 50}
        
					Ce système était utilisé dans la fin des années 50 et le début des années 60, il était mis en place par l'AMS (Army Map Service) des Etats-Unis. Le réseau géodésique comprenait une chaîne de triangles le long de la côte de la frontière Tuniso-Libyenne à la frontière Egypto-Libyenne et rattaché au système Europe 50 précédemment défini ci-dessus.

\textbf{Le Système LYB79}\index{Système LYB79}
        
					Ce système  est défini par une compensation générale des points du système Europe 50 avec 45 points observés par  la méthode Doppler\footnote{\textbf{Christian Doppler} (1803-1853): mathématicien et physicien autrichien.}. L'ellipsoïde de référence est l'ellipsoïde international de Hayford.\index{\textbf{Doppler C.}}\index{Ellipsoïde international de Hayford}
\subsection*{6.4.3. Les Systèmes Géodésiques au Maroc}
\textbf{Le Système Merchich}\index{Système Merchich}
        
				Il a été observé en 1922 ayant les éléments suivants :
				
- le point fondamental : le point Merchich (près de Casablanca au Maroc);

-	l'ellipsoïde de référence est l'ellipsoïde de Clarke 1880 Français;

-	l'azimut d'orientation : l'azimut astronomique de la direction Merchich-Berchid. 
\section{\textsc{Caractéristiques des Ellipsoïdes Géodésiques}}\index{Ellipsoïdes géodésiques}
\begin{center}
\[\begin{array}{cccccc}
\hline
\hline
     \mbox{Nom de } &\mbox{Demi-grand axe}& 1/f & e^2 & \mbox{Paramètres de}\\
 \mbox{l'Ellipsoïde}& \mbox{Demi-petit axe}&    &     & \mbox{définition}\\  \hline \hline
   \mbox{Clarke Français} & 6\,378\,249.200 &   293.46602 & 0.0068034877 & a,\,b \\       
      \mbox{1880}       &  6\,356\,515.000  &             &             &   \\ \hline
  \mbox{Clarke Anglais} & 6\,378\,249.145 &  293.46500 & 0.00680351128& a,\,1/f \\  
   \mbox{1880        } & 6\,356\,514.8696&        &              &         \\ \hline
      \mbox{Hayford 1909 ou} & 6\,378\,388.000& 297.00000 & 0.0067226700 & a \\ 
   \mbox{International 1924} &  6\,356\,911.940 &          &              & 1/f \\    \hline 
         \mbox{Krassovsky}$\footnote{Feodosy Nikolaevich Krassovsky (1878 - 1948) : astronome et géodésien russe.}$   & 6\,378\,245.000 &  298.30000 & 0.00669342162& a,\,1/f  \\ 
                             &6\,356\,863.0188 &             &             &          \\ \hline 
        \mbox{GRS 1967 (AIG)} & 6\,378\,160.000 &  298.24717 & 0.0066946053 & a,e^2 \\ 
                               & 6\,356\,774.516 &           &              &       \\ \hline
      \mbox{NWL 8} & 6\,378\,145.000 &  298.25000 & 0.0066945419 & a,\,1/f \\  
                   & 6\,356\,759.770 &            &              &         \\ \hline
     \mbox{WGS72 } & 6\,378\,135.000 &  298.26000 & 0.0066943178 & a,\,1/f \\  
                   &  6\,356\,750.520 &             &              &       \\ \hline
    \mbox{AIG 1975 } & 6\,378\,140.000 &  298.25700 & 0.0066943850 & a,\,1/f \\  
                     & 6\,356\,755.288 &            &              &           \\ \hline
 \mbox{APL Navigation} & 6\,378\,144.000 &  298.23000 & 0.0066949901 & a,\,1/f \\  
                       & 6\,356\,757.339 &             &             &         \\ \hline
 \mbox{GRS80 (AIG) } & 6\,378\,137.000 &  298.257222101& 0.0066943800229& a,b \\  
                      & 6\,356\,752.3141&              &                &      \\ \hline
  \mbox{WGS84 } & 6\,378\,137.000 &  298.257223563 & 0.0066943799 & a,\,1/f  \\  
                & 6\,356\,752.3142&                &              &          \\ \hline
 \hline 
     \end{array}
\]
\end{center}

\begin{table}[h]
	\centering
	\caption{Caractéristiques des ellipsoïdes géodésiques (\textit{H. Moritz}, 2000; \textit{C. Boucher}, 1979a; \textit{DMA}, 1987b)}
	\label{tab: Caractellipsoïdes}
\end{table}
\index{\textbf{Moritz H.}}\index{\textbf{Boucher C.}}\index{\textbf{Hayford J.F.}}\index{\textbf{Krassovsky N.F.}}
\newpage
\textbf{\un{Note historique:}} \textsl{L'ellipsoïde de Clarke 1880 Français a été déterminé par le géodésien anglais Clarke A. R. \footnote{\textbf{Clarke Alexander Ross} (1828-1914): géodésien anglais. \\ \,\,\, 7. \textbf{Feodosy Nikolaevich Krassovsky} (1878 - 1948) : astronome et géodésien russe.}entre les années 1878 et 1880. La méthode de détermination se repose sur des observations obtenues par les mesures d'arcs de méridiens de la Terre ainsi que les observations des latitudes géodésiques des points d'extrémités appuyés par la détermination de latitudes astronomiques.}

\textsl{A partir de la formule donnant la longueur d'un arc de la méridienne compris entre les latitudes géodésiques $\varphi$ et $\varphi'$, on écrit l'équation d'observation relative à la méthode des moindres carrés en partant d'un ellipsoïde approché de paramètres, dans la notation usuelle, $(a_0,b_0,e_0)$ soit:}
$$\left\{ \begin{array}{l}
ds +s_{cal}-s_{obs}=v \\
avec:\quad  s_{cal}=a_0(1-e^2_0)\ds \int_{\varphi}^{\varphi'}( 1-e_0^2sin^2t)^{-\frac{3}{2}}dt
\end{array}\right. $$
\textsl{Clarke A.R. avait introduit le coefficient $n=\ds \frac{a_0-b_0}{a_0+b_0}$, alors $s_{cal}$ s'écrit (à vérifier à titre d'exercice):}
$$  s_{cal}=s_{cal}(b_0,n)=\ds b_0(1+n)(1-n^2)\int_{\varphi}^{\varphi'}( 1+2ncos2t+n^2)^{-\frac{3}{2}}dt $$
\textsl{Il exprimait $ds$ en fonction de $db_0$ et $dn$ et il avait utilisé 56 observations d'arcs ou portions d'arcs. Les résultats obtenus sont:}

\textsl{- pour l'Angleterre: $a= 6\,378\,249.14533\,m$ arrondi à $6\,378\,249.145\,m$ et l'aplatissement $f=1/293.465$, c'est l'ellipsoïde de Clarke 1880 Anglais. Il ne sera pas utilisé pour la cartographie de la Grande-Bretagne mais c'est plutôt l'ellipsoïde d'Airy\footnote{\textbf{George Biddell Airy} (1801-1892): mathématicien et astronome anglais. }.}

- \textsl{pour la France: $a= 6\,378\,249.20\,m$, $b=6\,356\,515.00\,m$ et $f=1/293.466\,0208$. C'est l'ellipsoïde de Clarke 1880 Français. Il a été adopté en 1922 pour la carte de France et celles de l'Afrique du Nord.} (H. Monge, 1979) \index{Ellipsoïde de Clarke Français}\index{Ellipsoïde de Clarke Anglais}\index{\textbf{Clarke A.R.}}\index{Ellipsoïde d'Airy}\index{\textbf{Monge H.}}
\section{\textsc{Exercices et Problèmes}}
\bex
 Donner l'expression des composantes du gradient en coordonnées cylindriques.
\eex
\newpage
\bex
On donne l'expression scalaire d'une fonction $V(x,y,z)$  par : 
 	\[V(x,y,z)=\frac{ax^2+y^2}{z^2}+\frac{1}{2}\omega^2(x^2+y^2)
\]
1. Calculer les composantes du  vecteur $\textit{\textbf{grad}}V$ dans un domaine de $\BbR^3$ où $z\neq 0$. 
\eex
\bpb
Soit un point $ A (\varphi ,\lambda)$ sur un ellipsoïde de r\'evolution associ\'e \`a un r\'ef\'erentiel g\'eocentrique donn\'e $ \m R$. On consid\`ere le rep\`ere orthonorm\'e local en A $(e_{\lambda},e_{\varphi},e_n)$ d\'efini dans la base orthonorm\'ee $(i,j,k)$ de $ \m R$ où $e_{\lambda}$ est tangent au parallèle passant par $A$ et dirigé vers l'Est, $e_{\varphi}$ tangent à la méridienne, dirigé vers le nord et $e_n$ porté par la normale à l'ellipsoïde dirigé vers le zénith.

1. Exprimer les vecteurs de la base $(e_{\lambda},e_{\varphi},e_n)$ dans la base $(i,j,k)$ de $ \m R$.

2. Exprimer les vecteurs $i,j$ et $k$ dans la base $(e_{\lambda},e_{\varphi},e_n)$.

3. Calculer $de_{\lambda},de_{\varphi}$ et $de_n$ dans la base $(i,j,k)$.

4. En adoptant une écriture matricielle, montrer que :
$$	\begin{pmatrix}{
	de_{\lambda} \cr
de_{\varphi}  \cr
	de_n }
\end{pmatrix}=\begin{pmatrix}{
0 &	sin\varphi d\lambda & -cos\varphi d\lambda \cr
	-sin\varphi d\lambda & 0 & -d\varphi  \cr
	cos\varphi d\lambda & d\varphi & 0}
\end{pmatrix}\begin{pmatrix}{
	e_{\lambda} \cr
	e_{\varphi} \cr
	e_n }
\end{pmatrix}$$
\epb
\bpb
On considère les notations du précédent problème. Soit un point $M$. On pose:
$$ \Delta X=(X_M-X_A,Y_M-Y_A,Z_M-Z_A)^T,\, \Delta x=(x_M,y_M,z_M)^T$$
où $X$ et $x$ sont respectivement les composantes du vecteur $AM$ dans les repère $\m R$ et le repère local en $A$.

1. Montrer qu'on a la relation: $\Delta X=J.\Delta x$ avec $J$ une matrice orthogonale $(J^{-1}=J^T)$ qu'on déterminera. 

2. On suppose maintenant que $\m R$ est le repère GPS et que le passage du repère $\m R$ vers le repère terrestre est donné par le modèle dit à 7 paramètres (§ \ref{transfor}): 
\ba  
& \begin{pmatrix}{
X'\cr
Y'\cr
Z'}
	\end{pmatrix}=F(X)=\begin{pmatrix}{
Tx\cr
Ty\cr
Tz}
	\end{pmatrix}+(1+m).R.\begin{pmatrix}{
X\cr
Y\cr
Z}
	\end{pmatrix}& \nonumber \\ &=\begin{pmatrix}{
Tx\cr
Ty\cr
Tz}
	\end{pmatrix}+(1+m).\begin{pmatrix}{
1 & rz & -ry \cr
-rz & 1 & rx\cr
ry & -rx & 1}
	\end{pmatrix}
.\begin{pmatrix}{
X\cr
Y\cr
Z}
\end{pmatrix}
	\nonumber & 
	\ea
	 On note: $\delta X=F(\Delta X)-\Delta X$, que représente $\delta X$. Montrer qu'au premier ordre que:
$$\delta X\cong m\Delta X+(R-I)\Delta X$$
3. On appelle $\omega=(rx,ry,rz)^T$, montrer que:
$$ ||\delta X||=(m^2+\omega^2sin^2\theta)^{\frac{1}{2}}||\Delta X||$$
où $\theta$ est l'angle entre les vecteurs $\Delta X$ et $\omega$.	

4. En déduire que:
 $$|m|.||\Delta X||\leq ||\delta x||\leq (m^2+\omega^2)^{\frac{1}{2}}||\Delta X||$$
5. En utilisant la relation liant $\Delta X$ et $\Delta x$; montrer que:
$$\delta x=m\Delta x+J^T(R-I)J\Delta x$$
\epb
\bpb
  On définit dans $\BbR^3$ un point $M$ par ses coordonnées ellipsoïdiques de Jacobi $(\phi,\lambda,u)$ comme suit:
$$ M \left\{ \begin{array}{l}
x=\sqrt{u^2+\epsilon^2}.cos\phi cos\lambda \\
y=\sqrt{u^2+\epsilon^2}.cos\phi sin\lambda \\
z=u.sin\phi 
\end{array} \right. $$
avec: $\epsilon^2=\sqrt{a^2-b^2},\,\phi \in [-\pi/2,\pi/2 ],\, \lambda \in [0,2\pi]$ et $u\in ]0,+\infty[$, $a,b$ deux constantes réelles telles que $a>b>0$. 

1. Montrer que le point $M$ appartient à un ellipsoïde de révolution en précisant ses demi-axes.

2. Calculer $ds^2$ et montrer qu'il s'écrit sous la forme:
$$ ds^2=(d\phi,d\lambda,du).G.\begin{pmatrix}{
d\phi \cr 
d\lambda \cr
du }
\end{pmatrix}$$
avec $G$ donnée par :
$$ G=(g_{ij})=\begin{pmatrix}{
u^2+\epsilon^2sin^2\phi & 0 & 0 \cr 
0 & (u^2+\epsilon^2)cos^2\phi & 0  \cr
0 & 0 & \ds \frac{u^2+\epsilon^2sin^2\phi}{u^2+\epsilon^2} }
\end{pmatrix}$$

3. Sachant que l'expression du laplacien d'une fonction scalaire $V$ en coordonnées de Jacobi est exprimée par:
$$ 
\Delta V=\frac{1}{\sqrt{g}}\left\{\frac{\partial}{\partial \phi}\left(\frac{\sqrt{g}}{g_{11}}.\frac{\partial V}{\partial \phi}\right)+\frac{\partial}{\partial \lambda}\left(\frac{\sqrt{g}}{g_{22}}.\frac{\partial V}{\partial \lambda}\right)+\frac{\partial}{\partial u}\left(\frac{\sqrt{g}}{g_{33}}.\frac{\partial V}{\partial u}\right)\right\} $$
où $g$ est le déterminant de la matrice $G$, donner l'expression de $\Delta V$. 

4. Calculer $\Delta V$ sachant que $V$ est donnée par:
$$ V(\phi,u)=\frac{GM}{\epsilon}Arctg\frac{\epsilon}{u}+ \frac{1}{3}a^2\omega^2\frac{q}{q_0}\left(1-\frac{3}{2}cos^2\phi\right)+\frac{1}{2}\omega^2(u^2+\epsilon^2)cos^2\phi $$
avec $G,M$ et $\omega$ des constantes et:
\ba
q=q(u)=\frac{1}{2}\left[\left(1+3\frac{u^2}{\epsilon^2}\right)Arctg\frac{\epsilon}{u}-3\frac{u}{\epsilon}\right] \nonumber \\
q_0=q(u=b)=\frac{1}{2}\left[\left(1+3\frac{b^2}{\epsilon^2}\right)Arctg\frac{\epsilon}{b}-3\frac{b}{\epsilon}\right] \nonumber
\ea
\epb
\bpb
Avec les notations usuelles, un potentiel est donné avec les constantes $(GM, a, J_2,\omega)$ de GRS80 par :
$$  W(r,\theta,\lambda)=\ds \frac{GM}{r}\left(1-J_2\frac{a^2}{r^2}P_2(cos\theta)\right)+\frac{\omega^2}{2}r^2sin^2\theta $$
où $(r,\theta,\lambda)$ sont les coordonnées sphériques du point de calcul et :

$a = 6378\,137.00\, m\, ,b = 6356\, 752.31\, m$;

$GM = 3\, 986\, 005\times 10^8\, m^3 s^{-2}$;

$\omega = 7\, 292\, 115\times 10^{-11} rad. s^{-1}$;

$J_2 = 108\, 263\times 10^{-8}.$

1. Calculer le potentiel $W$ aux points suivants situés sur l'ellipsoïde GRS80:

* le Pôle Nord;

* sur l'equateur;

* au point $A$ sur l'ellipsoïde avec $\theta_A = 34^{\circ}$.

2. Calculer les variations de $W$ entre les 3 points.

3. Exprimer le potentiel $W$ en coordonnées cartésiennes $(X,Y,Z)$.

On rappelle : $P_2(cos\theta)=\ds \frac{1}{2}(3cos^2\theta -1)$.
\epb
\chapter{\textit{\textbf{Les Réseaux Géodésiques}}}
\section{\textsc{Introduction}}
L'un des buts de la Géodésie est l'établissement des réseaux géodésiques dans un territoire donné. Ces réseaux géodésiques vont constituer l'ossature des travaux cartographiques et topographiques.
\\

Généralement, la région à cartographier est une zone vaste très étendue. Les procédés topographiques comme les polygonales ne peuvent pas être utilisés et c'est dû :

-	premièrement, la surface topographique n'est pas un plan mais plutôt un sphéroïde. Ainsi la sphéricité de la Terre est négligée dans les travaux topographiques. De plus, les corrections de la représentation plane ne sont pas prises en compte;

-	deuxièmement, vu l'étendue de la zone, les levés topographiques ne peuvent pas être faits à partir d'une seule polygonale. On est amené à faire plusieurs polygonales, celles-ci sont déterminées les unes indépendantes des autres ainsi que leurs compensations ou ajustements.
\\

Le groupement de ces polygonales va cumuler les erreurs dès qu'on s'éloigne de la polygonale choisie polygonale de départ. Aussi, on ne peut laisser les erreurs, tolérées pour une polygonale, se cumuler et  falsifier la position des points.
\\

Le but de la géodésie est donc la détermination avec précision des coordonnées de ces points dispersés sur tout le territoire, objet de la carte. 
\section{\textsc{Les Réseaux Géodésiques Classiques}}
\subsection*{7.2.1. Conception d'un Réseau Géodésique}
Un réseau géodésique est généralement constitué par une chaîne de triangles où les sommets représentent les points géodésiques souvent matérialisés aux sommets des montagnes et des constructions et bâtiments élevés (châteaux d'eau, phares,...). Ce réseau de triangles couvre l'ensemble du territoire.
\subsection*{7.2.2. Point Fondamental}\index{Point fondamental}
Le réseau sera déterminé de proche en proche à partir d'un premier point. Ce point est appelé Point Fondamental du réseau géodésique ou du système géodésique associé. En ce point, on détermine $(\Phi, \Lambda)$  respectivement la latitude et la longitude astronomiques à partir d'observations sur des étoiles. La méthode la plus utilisée était celle des hauteurs égales. On observe aussi l'azimut astronomique de la direction vers un deuxième point. Cet azimut sera déterminé à partir de l'observation astronomique de l'étoile polaire. L'orientation du réseau est définie à partir de l'azimut géodésique $Azg$ de la direction citée précédemment. Au point fondamental, on prendra par convention:
\be
\fbox{ $ \begin{array}{l}
	\varphi_{\mbox{géodésique}} = \Phi     \\
    \lambda_{\mbox{géodésique}} = \Lambda  \\
			\quad Azg = Aza
	\end{array} 	$}\label{p162}
\ee
Ce choix permet de positionner le réseau par rapport à l'ellipsoïde de référence, ce qui implique qu'au point fondamental la normale à l'ellipsoïde est confondue avec la verticale.
\\

Pour déterminer les coordonnées d'un deuxième point $M_2$, il faut mesurer la distance $M_1M_2$. A partir des calculs géodésiques sur l'ellipsoïde, on peut déterminer $(\varphi_2 , \lambda_2)$ du point $M_2$  connaissant $(\varphi_1,\lambda_1)$ du point $M_1$ , $Azg$ et  $M_1M_2$.
\subsection*{7.2.3. Mise à l'Echelle }
La mesure de la distance définit la mise à l'échelle du réseau. En général, on mesure une base de quelques kms (exemple 10 kms), par des opérations de mesures angulaires, cette base est intégrée au réseau géodésique ce qu'on appelle '\textbf{amplification de la base}'.
\subsection*{7.2.4. Processus de Réalisation}
La méthode la plus utilisée dans le processus de réalisation d'un réseau géodésique terrestre était la triangulation : soit la mesure des angles des triangles. Ces angles sont obtenus par différences de lectures entre les points et mesurés à l'aide des théodolites de haute précision comme le TH3 par la méthode des tours d'horizon. Un premier contrôle de la qualité des observations s'obtient par le calcul de la fermeture des triangles géodésiques. Soit $\epsilon $ l'excès sphérique d'un triangle en grades, on détermine $\hat{F}$, la fermeture du triangle, par :
\be
	\fbox{$ \hat{F} = \hat{A}+\hat{B}+\hat{C} - 200\,gr - \epsilon $} \label{p163}
\ee
où $\hat{A},\hat{B},\hat{C}$ sont les angles bruts mesurés du triangle. 
\\

En fonction des longueurs des côtés des triangles, la fermeture des triangles doit être inférieure à la tolérance requise.
\subsection*{7.2.5. Equation de Laplace et Contrôle de la qualité du réseau géodésique}
Pour contrôler la qualité du réseau géodésique, on observe généralement tous les 200 kms un point de Laplace, soit la détermination de $\varphi_a$, $\lambda_a$ et $Aza$. On montre qu'en un point $M$ de coordonnées géodésiques $(\varphi_g,\lambda_g)$ le passage d'un azimut astronomique $Aza$ d'une direction observée en $M$ à l'azimut géodésique de la même direction est donné par la formule dite équation de Laplace: 
\be
\fbox{ $ Azg=Aza-(\lambda_g-\lambda_a)sin\varphi $} \label{p163a}
\ee
où $\lambda_a$ est la longitude astronomique de $M$ et $\varphi$ est la latitude géodésique $(\varphi_g)$ ou la latitude astronomique $(\Phi_a)$.
\\

A l'aide de l'équation de Laplace\index{Equation de Laplace} citée ci-dessus, on contrôle la fermeture Laplacienne par :\index{\textbf{Laplace P.S.}}
\be
	\fbox{ $ E = Azg - Aza - (\lambda_g - \lambda_a )sin\varphi $} \label{p164}
\ee
De même, on observe aussi une base pour améliorer la qualité de l'échelle du réseau géodésique.    
\section{\textsc{La Géodésie Spatiale}}
La géodésie spatiale a connu un essor dans la mise en place des réseaux géodésiques à partir des années 80 par l'utilisation de la méthode Doppler\index{La méthode Doppler}. Cette méthode repose sur l'effet dit Doppler sur les signaux émis par les satellites en mouvement et qui sont reçus par un récepteur centré en position sur le point géodésique. Le traitement des données enregistrées  permet la détermination des coordonnées géocentriques du point géodésique $(X,Y,Z)$ ou $(\varphi,\lambda,he)$ dans un système géodésique bien défini.
\\

La technique Doppler a été remplacée par la méthode de positionnement par les satellites GPS (\textbf{G}lobal \textbf{P}ositioning \textbf{S}ystem). Ce système permet d'observer en tout point de la Terre et en tout instant au moins quatre satellites de la constellation des satellites GPS (\textit{B. Hofmann-Wellenhof - H. Lichtenegger - J. Collins}, 1992).\index{\textbf{Hofmann-Wellenhof B.}}\index{\textbf{Lichtenegger H.}}\index{\textbf{Collins J.}}
\subsection*{7.3.1. Le Système de Référence $WGS84$}\index{Système de référence $WGS84$}
Le système $WGS84$ (\textbf{W}orld \textbf{G}eodetic \textbf{S}ystem 84) est le système lié aux satellites GPS. C'est un système géocentrique $(O,X,Y,Z)$, c'est-à-dire :

-	son origine $O$ est le centre des masses de la Terre;

-	l'axe $OZ$ est confondu avec l'axe de la rotation de la Terre défini par le BIH (Bureau International de l'Heure) à l'instant 1984.0;

-	le plan $OXY$ contient le méridien origine des longitudes;

-	l'axe $OY$ est perpendiculaire à $OXY$ tel que $(O,X,Y,Z)$ soit un système direct.
\\

A ce système est associé l'ellipsoïde de référence appelé ellipsoïde $WGS84$ avec $a = 6\,378\,137.00\,m$  et $1/f = 298.257\,223\,563\,40697$ (DMA, 1987a - 1987b - 1987c).
\subsection*{7.3.2. Le Système de Référence $ITRF$}\index{Système de référence $ITRF$}
Le système $WGS84$ a une précision absolue de $1\, m$. Pour avoir une précision centimétrique, la communauté géodésique internationale représentée par l'Association Internationale de Géodésie (AIG)\footnote{AIG: l'Association Internationale de Géodésie est l'une des huit associations de l'Union Géodésique et Géophysique Internationale (UGGI) (www.iag-aig.org).} a défini un système géocentrique très proche du système $WGS84$ dont la précision est centimétrique. Ce système est appelé le système $ITRF$ (\textbf{I}nternational \textbf{T}errestrial \textbf{R}eference \textbf{F}rame) associé à l'ellipsoïde $GRS80$\index{Ellipsoïde $GRS80$} défini par l'AIG avec ses paramètres $a = 6\,378\,137.00\, m$ et $1/f = 298.257\,222\,100\,88378$ ou $e^2 = 0.006\,694\,380\,023$ (\textit{H. Moritz}, 2000). Comme le système $WGS84$, le  système  $ITRF$ est géocentrique. Ce système est défini à partir de 1989 c'est pourquoi, on parle du système $ITRF89$ comme système utilisé par la communauté géodésique internationale.\index{\textbf{Moritz H.}}
\subsection*{7.3.3. La Détermination d'un Nouveau Réseau Géodésique}
La méthode de positionnement par GPS a révolutionné la détermination des réseaux géodésiques. On décrit ci-dessous le cas de la mise en place d'un nouveau réseau géodésique dans un pays.
\subsubsection*{Observations du Réseau Géodésique GPS de Base}
La technique GPS (\textit{B. Hofmann-Wellenhof - H. Lichtenegger - J. Collins}, 1992) ne nécessite plus le stationnement des points sur les sommets des montagnes pour faciliter le passage des visées entre les points. On choisit alors un ensemble de points d'accès facile, formant le Réseau Géodésique GPS de Base. Après la phase de la construction de ces points, on passe à la phase des observations par l'emploi de récepteurs GPS bi-fréquences avec au moins 3 récepteurs pour avoir une base commune. La durée des observations dépendra de la distance moyenne entre les points. Supposons que cette moyenne soit de l'ordre de 120-150 kms, alors la durée des observations est au minimum 3 sessions de 8 heures (une session est une durée d'observations sans interruption). On paramètre les récepteurs de façon à enregistrer tous les $30\,s$ les signaux émis  par les satellites GPS qui se trouvent à $10^{\circ}$ au-dessus de l'horizon des stations. 
\\

\begin{figure}[htp]
	\centering
		\includegraphics[width=0.80\textwidth]{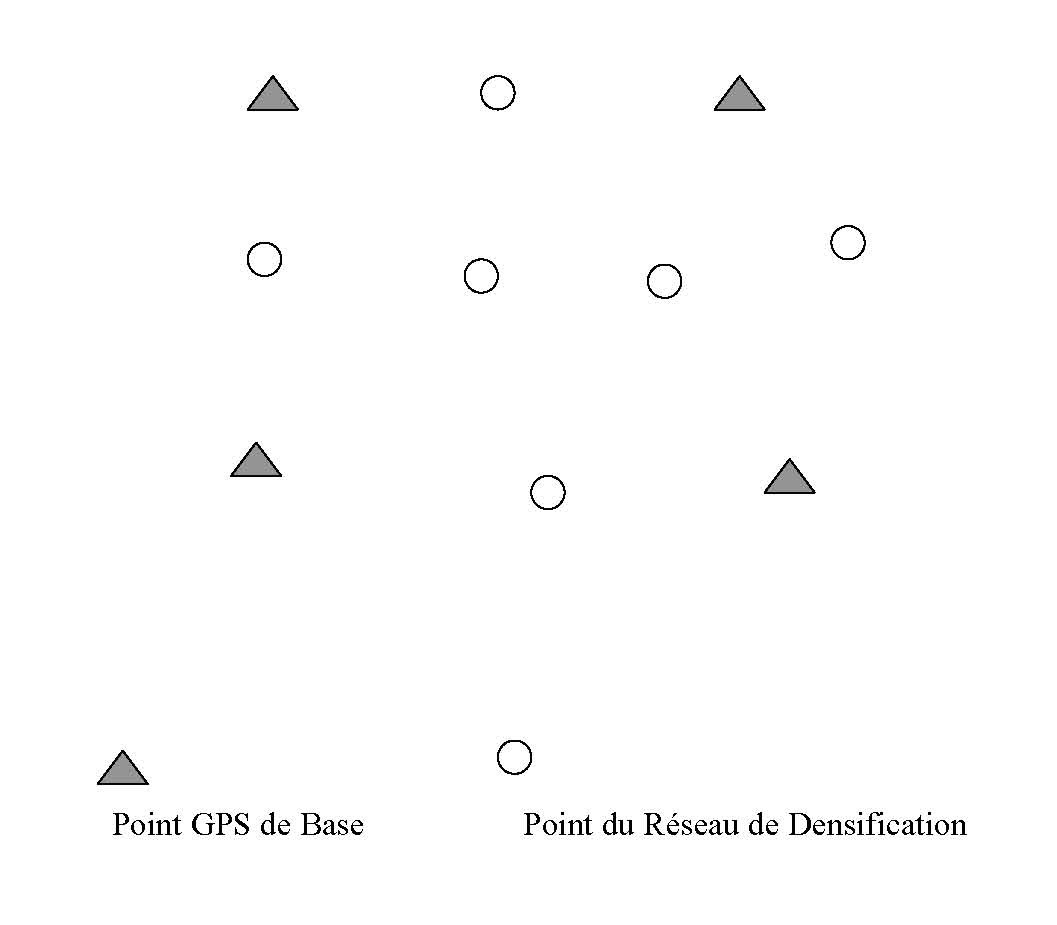}
	\caption{Canevas des points}
	\label{fig:canevas1}
\end{figure}

\textbf{Exemple :} on a à déterminer 6 points n° 1, 2, 3, 4, 5 et 6 avec 3 récepteurs (\textbf{Fig. \ref{fig:canevas1}}).
 Les observations seront faites comme suit :
\begin{center}
\[\begin{array}{ccc}
\hline
   \mbox{Jour n°}&\mbox{Points observés}\\ \hline
                1 - 2 - 3  &        1 - 2 - 3    \\   \hline
          4 - 5 - 6       &       2 - 3 - 4      \\   \hline
          7 - 8 - 9   &   3 - 6 - 4 \\  \hline 
     10 - 11 - 12   &   5 - 6 - 4 \\ \hline
     13 - 14 -15    &   5 - 4 - 2  \\ \hline
\end{array}
\]
\end{center}
\begin{table}[htp]
	\centering
	\caption{Table des observations GPS}
	\label{tab: Tableau des Observations}
\end{table}
\subsection*{7.3.4. Calculs des points du Réseau Géodésique GPS de Base}
Le système $ITRF$ est défini par un ensemble de stations dont les coordonnées sont connues avec précision. Ces points serviront comme points fixes dans le calcul des coordonnées des points du réseau géodésique GPS de base d'un pays. En effet, les observations de ces points sont accessibles via l'Internet et on récupère les données observées dans la même période. Parmi les points les plus proches, par exemple, des pays de l'Afrique du Nord, on cite les stations Matera et Cagliari(Italie), Grasse (France) et Madrid (Espagne).
\\

Pour calculer les coordonnées des points nouveaux, on utilise les éphémérides précises c'est-à-dire les positions précises des satellites au moment des observations. Ces éphémérides seront déterminées deux semaines après les observations par les centres de calculs de l'IGS (International GPS Service)\index{IGPS Service} service dépendant de l'AIG.
\\

A la fin des calculs, on obtient les coordonnées géodésiques $(\varphi , \lambda, he )$  dans le système $ITRF$, ellipsoïde $GRS80$. A partir des formules suivantes :
\be
M=\left\{ \begin{array}{lll}
	X =(N+he)cos\varphi cos\lambda \\
	Y=(N+he)cos\varphi sin\lambda	\\
	 Z=(N(1-e^2)+he)sin\varphi
\end{array}\right.  \nonumber
\ee
on obtient les coordonnées 3D $(X,Y,Z)$. Si on dispose d'un modèle de géoïde, on détermine les altitudes orthométriques des points à partir de la relation :
\be
	he =  N + H \nonumber
\ee
avec:
 
- $N$ désigne ici l'ondulation du géoïde\index{Ondulation du géoïde} au dessus/au dessous de l'ellipsoïde $GRS80$;
   
- $he$ l'altitude ellipsoïdique GPS\index{Altitude ellipsoïdique GPS};
              
- $H$ l'altitude orthométrique\index{Altitude orthométrique}.
\section{\textsc{Densification du Réseau Géodésique GPS de Base}}
Soit le schéma suivant : 
la phase des observations se fait comme suit pour 4 récepteurs (\textbf{Fig. \ref{fig:image1}}):

-	un récepteur sur un point connu du Réseau GPS de Base (exemple point 1);

-	les autres récepteurs sur les points nouveaux, par exemple: 

   *	1ère session (2 à 3 heures) : 1 récepteur sur le point 101, un autre sur le point 102 et un autre sur le point 103;
   
   *	2ème session :  1 récepteur sur le point 104, un autre sur le point 105 et un autre sur le point 106.

On aura le schéma suivant pour les observations :

\begin{figure}[htp]
	\centering
		\includegraphics[width=0.60\textwidth]{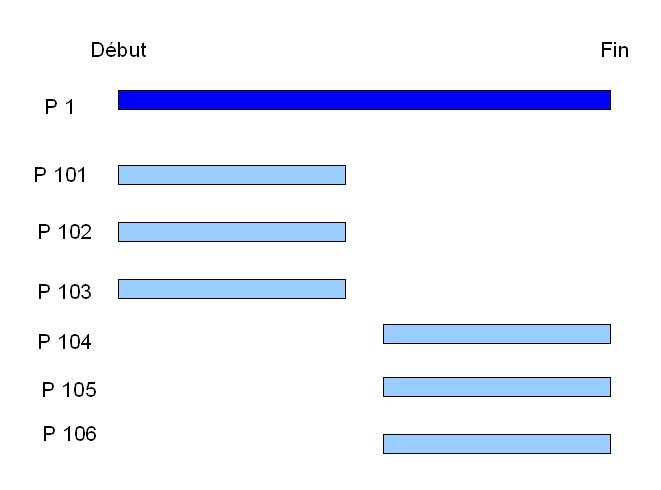}
	\caption{Sessions de densification d'observations GPS}
	\label{fig:image1}
\end{figure}
Les observations se font en groupe de points afin d'obtenir une meilleure détermination des points lors des calculs.
\\

La phase des calculs repose sur la fixation du point 1 dont les coordonnées sont connues et on utilise les éphémérides radiodiffusées\index{Ephémérides radiodiffusées} portées par les signaux envoyés par les satellites GPS. Nous obtenons les coordonnées géodésiques $(\varphi, \lambda, he)$  des points nouveaux.
\section{\textsc{La Densification d'un Réseau Géodésique Terrestre }}
A l'heure actuelle, la méthode de positionnement par GPS est devenue la méthode la plus utilisée pour densifier un réseau géodésique terrestre par de nouveaux points.
\\

 La méthodologie est la suivante:
 
\begin{enumerate}
	\item observations par GPS d'un ensemble de points anciens délimitant la zone à densifier, le nombre des points anciens est au minimum 3 points. Ces observations permettent de déterminer les paramètres de transformation entre le système GPS et le système géodésique terrestre en vigueur dans un pays;
  \item observations des points nouveaux : la durée des observations dépend de l'ordre des points nouveaux. La méthode des observations des points est la même que celle décrite pour la densification d'un nouveau réseau géodésique;
  \item calculs des observations GPS des points anciens et détermination des paramètres de transformation du système GPS  au système géodésique terrestre officiel;
  \item	calculs des observations GPS des points nouveaux en fixant les coordonnées GPS des points anciens;
  \item	application des paramètres de transformation aux coordonnées GPS pour avoir les coordonnées géodésiques des points nouveaux dans le système géodésique terrestre en vigueur;
  \item	calculs des coordonnées planimétriques des points nouveaux dans la représentation plane officielle du pays.
\end{enumerate}
\section{\textsc{Exercices et Problèmes}}
\bex
Soit un triangle plan $ABC$ de côtés respectivement $a,b$ et $c$ (\textbf{Fig. \ref{fig:triangle}}):
\begin{figure}[htbp]
	\centering
		\includegraphics[width=0.5\textwidth]{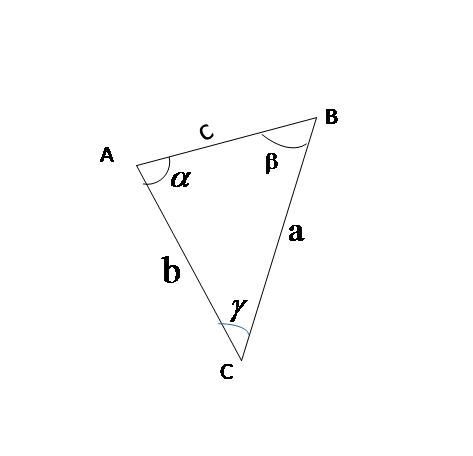}
		\caption{Le triangle plan}
		\label{fig:triangle}
\end{figure}

1. A partir de la formule:
$$cos\alpha=\frac{b^2+c^2-a^2}{2bc}$$
exprimer $d\alpha$ en fonction de $da,db$ et $dc$.

2. Donner l'expression de la variance $\sigma^2_{\alpha}$ en fonction des variances des côtés. 

3. On considère maintenant que le triangle $ABC$ est équilatéral et on note $L=a=b=c$. Montrer que :
$$ \sigma_{\alpha}=\ds \sqrt{2}\frac{\sigma_L}{L}$$ 
4. En déduire l'expression de la tolérance de la fermeture du triangle.
\eex

\chapter{\textit{\textbf{Réduction des Distances}}}
\section{\textsc{Introduction}}
Souvent dans les travaux géodésiques, on est amené à mesurer des distances entre des points géodésiques à l'aide des distancemètres. Pour pouvoir utiliser ces distances, il faut y apporter des corrections  pour réduire ces distances au plan de la représentation plane utilisée.
\section{\textsc{Corrections des Distances}}
 On donne une distance spatiale $D$ mesurée par un distancemètre entre deux points $A$ d'altitude  $H_A$ et $B$ d'altitude $H_B$.
\\

Vue la présence de l'atmosphère, la distance $D$ est une distance courbe de rayon de courbure $\rho$ qui dépend du type de l'appareil utilisé.
\\

Soit $\alpha$ l'angle au centre de rayon de courbure, on a alors : 
\be
	D  =  \rho \alpha  \label{da1}
\ee
De plus, on a :
\be
	sin\frac{\alpha}{2}=\frac{D_P}{2\rho} \label{da2}
\ee
\subsection*{8.2.1. Correction de courbure}\index{Correction de courbure}
La correction de courbure est $C_1$  telle que : 
\be
C_1 = D_P -  D = 2\rho sin\frac{\alpha}{2} -\rho \alpha =-\rho(\alpha-2sin\frac{\alpha}{2}) \label{da3}
\ee
 Comme  $ \alpha $ est un angle très petit < 3°, on peut écrire :
\[
	sin(\frac{\alpha}{2}) =\frac{\alpha}{2}- \frac{(\alpha/2)^3}{6}+... =\frac{\alpha}{2}- \frac{\alpha^3}{48}+... 
\]
soit:
\[	C_1 =   - \rho ( \alpha -\alpha+ \frac{\alpha^3}{24}+..  )
\]
D'où : 
\be
\fbox{ $ 	C_1 =D_P-D=  \ds  - \rho \frac{\alpha^3}{24} $}\label{da4}
\ee
Pour les distancemètres à ondes lumineuses, on a $\rho = 8R$ où $R$ désigne le rayon moyen de la Terre = 6378 $km$. Pour les distancemètres à ondes centimétriques on a $\rho = 4R$. La correction $C_1$ vaut en fonction de la distance mesurée :
\be
	C_1 =  -\frac{D^3}{24\rho^2} \label{da5}
\ee
\textbf{Application}:
 $ D = 10\, km$, $\rho   = 4R\Longrightarrow \,C_1 = -0.06\, mm $ quantité négligeable.
\subsection*{8.2.2. Réduction à l'horizontale }\index{Réduction à l'horizontale}
Soit $H_m =  (H_A + H_B )/2$ l'altitude moyenne de la visée $AB$, $i$ l'angle que fait $AB$ avec le plan horizontal, $D_H$ la distance horizontale à l'altitude moyenne. 

La correction à l'horizontale est donnée par $C_2$ telle que :           
\be
	 C_2 = D_H - D_P  \label{da6}
\ee
Or  $D_H= D_P.cosi$, d'où :                     
	\[ C_2 = D_P cosi - D_P  = D_P (cosi - 1) = -2D_P .sin^2\frac{i}{2} 
\]
En posant:           
\be
	 \Delta H    = H_B - H_A  \label{da7} 
\ee
alors, on peut écrire que :      
	\[         sin(\frac{i}{2})=\frac{\Delta H}{2D_P}  \]
Par suite :
\be
	 \fbox{$ C_2 = D_H-D_P=-\ds  \frac{\Delta H^2}{2D_P} $}  \label{da8}                       
\ee
    \subsection*{8.2.3. Réduction à la surface de référence ou au niveau zéro}\index{Réduction au niveau zéro}
 Soit $H_m = (H_A + H_B)/2 $ l'altitude moyenne au dessus de l'ellipsoïde. $R$ le rayon de courbure dans la direction $AB$. On a alors:
\be
	\frac{D_H}{R+H_m}=\frac{D_0}{R} \label{da9}
\ee
La correction $C_3$  est telle que :  
\be
	                          C_3 = D_0 - D_H   \label{da10}                           
\ee
Or:                                                                                                   
\be
	D_0=\frac{R}{R+H_m}D_H \label{da11}
\ee
Par suite :                                             
	\[C_3=\frac{R}{R+H_m}D_H -D_H= -D_H \frac{H_m}{R+H_m}
\]
Comme $H_m$  étant négligeable par rapport au rayon de courbure de l'ellipsoïde de référence, on obtient :
\be
	 \fbox{ $ C_3 = D_0 - D_H =\ds -D_H \frac{H_m}{R} $}\label{da12}                                        
\ee
\textbf{Exemple Numérique} :  
$D_H = 10\, km$, $H_m = 800\, m$, $R = 6378\, km$ d'où : $C_3 = D_0 - D_H = -10\times 0.8/6378 = -1.25 \, m$.
\subsection*{8.2.4. Correction de courbure au niveau zéro : Passage de $D_0$ à $D_e$}
Soit $D_e$ la distance suivant l'ellipsoïde en l'assimillant en un arc de cercle de rayon $R$, on a :                                                                                                        
\be
	D_e = R.\beta  \label{da13}
\ee
\begin{figure}
	\centering
		\includegraphics[width=0.50\textwidth]{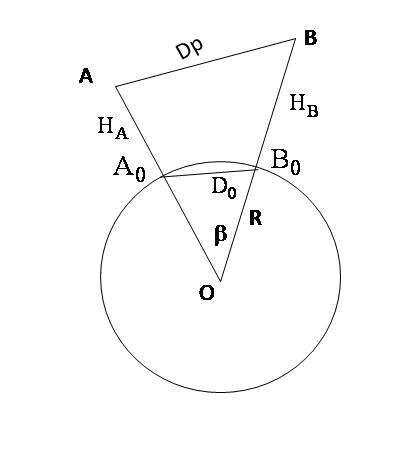}
	\caption{Réduction de la distance suivant la pente}
	\label{fig:dist1}
\end{figure}

et : 
\be
	sin\frac{\beta}{2}=\frac{D_0}{2R}\Longrightarrow \,D_0=2Rsin\frac{\beta}{2} \label{da14}
\ee
d'où :                            
\be
	C_4 =   D_e - D_0 =    R.\beta - 2.Rsin(\frac{\beta}{2}) = R.\beta - 2sin(\frac{\beta}{2}) \label{da15}
\ee
Comme $\beta$ est petit, on peut écrire que :
	\[sin\frac{\beta}{2}=\frac{\beta}{2}-\frac{1}{6}\frac{\beta^3}{2^3}=\frac{\beta}{2}-\frac{\beta^3}{48}+.... \]
d'où:
	\[C_4 = R(\beta-\beta + \frac{\beta^3}{24}+...)= \frac{R\beta^3}{24}= \frac{RD_e^3}{24R^3} \]
soit:
\be
		C_4 = \frac{D_{e}^3}{24R^2} \label{da16}
\ee
Dans (\ref{da16}), on peut remplacer $D_e$ par $D_0$ ce qui donne : 
\be
	\fbox{$ C_4 = D_e-D_0=\ds \frac{D_{0}^3}{24R^2} $} \label{da17}
\ee
\textbf{Exemple numérique} : 
$D_0  = 10\, km$ $\Longrightarrow \,C_4  = 0.001\,2\,m = 1.2\, mm$.
\subsection*{8.2.5. Réduction de la distance au plan de la représentation  plane}\index{Réduction au plan de la représentation plane}
\subsubsection*{8.2.5.1. Le module linéaire}\index{Module linéaire}
\bdf
On appelle le module linéaire en un point donné pour une représentation plane, le rapport $m$ tel que :
\be
\fbox{ $	m=\ds \frac{\mbox{distance plan}}{\mbox{distance ellipsoïde}}=\frac{D_r}{D_e} $} \label{da18}
\ee
\edf
\subsubsection*{8.2.5.2. L'altération linéaire}\index{Altération linéaire}
\bdf
On appelle l'altération linéaire en un point donné pour une représentation plane, le coefficient $\epsilon$ tel que :
\be
\fbox{$ \epsilon      =  m - 1  $}    \label{da19}                                
\ee
\edf
\subsection*{8.2.6. Correction de la distance au passage au plan}
 On a :   
\be
	 C_5 = D_r - D_e  \label{da20}
\ee
Soit : 
\be
 \fbox{ $ C_5 =  D_r - D_e = mD_e - D_e = ( m - 1)D_e=   \epsilon D_e $}    \label{da21}                                
\ee
\textbf{Exemple} :
On donne $D_e = 10\,000\,m$, $ \epsilon = +12\,cm/km$, alors $C_5 = 12\times10=120\,cm$ = $1.20\,m$ et $D_r = D_e + C_5 = 10\,000 + 1.20 = 10\,001.20\,m$.
\section{\textsc{Formule rigoureuse de passage de} $D_P$ \textsc{à} $D_0$}
Dans le triangle $OAB$ (\textbf{Fig. \ref{fig:dist1}}), on a :
\be
	D^2_P = (R+H_A)^2 +(R + H_B)^2 - 2(R +H_A)(R + H_B)cos\beta \label{da22}
\ee
D'où :
\be
	cos\beta=\frac{(R+H_A)^2 +(R + H_B)^2 -D_P^2}{2(R+H_A)(R + H_B)} \label{da23}
\ee
Dans le triangle $OA_0B_0$, on a :
\begin{equation}
	          D_0^2 = R^2 + R^2 -2R^2cos \beta    \label{da24}                                                  
\end{equation}
et : 
\be
	cos\beta=\frac{2R^2-D_0^2}{2R^2} \label{da25}
\ee
En écrivant que les équations (\ref{da23}) et (\ref{da25}) sont égales et en posant  $\Delta H = H_A -H_B$, on obtient la formule rigoureuse :
\be
	\fbox{ $ D_0=D_P\sqrt{\frac{\displaystyle 1-\frac{\Delta H^2}{D_P^2}}{\left(\displaystyle 1+\frac{H_A}{R}\right)\left(\displaystyle 1+\frac{H_B}{R}\right)}} $}\label{da26}
\ee

\section{\textsc{Exercices et Problèmes}}
\bex
 On a mesuré une distance suivant la pente $D_P = 20130.858\, m $ entre deux points $A$ et $B$ avec  $H_A = 235.07\, m,\, H_B = 507.75\, m$, on prendra comme rayon terrestre  $R = 6378\, km$.

1.	Calculer la distance suivant l'ellipsoïde :

-	en utilisant les différentes corrections;

-	en utilisant la formule rigoureuse.

2.	En prenant la valeur de la formule rigoureuse et sachant que le module linéaire $m$ vaut $0.999\,850\,371$, calculer la distance réduite au plan de la représentation plane utilisée.
\eex
\bex
Entre 2 points $A$ ( $H_A = 128.26\, m$ ) et  $B$  ( $H_B = 231.84\, m$), la distance $D_P$ suivant la pente est égale à $15\,498.823\, m$. Soit  $D_0$ la distance corde au niveau de la surface de référence. L'angle de site observé en $A$ en direction de $B$ est  $i = 0.3523\, gr$.

1.	Calculer la valeur de $D_0$ en utilisant la formule rigoureuse.

2.	Calculer $D_0$ par les corrections.

3.	En adoptant la moyenne des deux méthodes, calculer la distance  $D_e$ réduite à la surface de référence.

4.	Le module linéaire de la représentation plane utilisée est de $0.999\,648\,744$. Calculer alors la distance $D_r$ réduite au plan de la représentation.
\eex
\bex
On a mesuré une distance suivant la pente entre les points $A\, ( H_A =  1\,319.79\, m)$ et  $B\,( H_B = 1\,025.34\, m)$ avec $D_P = 16\,483.873\, m$.

1. Calculer la distance $D_e$ distance réduite à l'ellipsoïde de référence par la formule rigoureuse, on prendra le rayon de la Terre $R= 6378\, km$.

2. Calculer la distance  $D_r$  réduite à la représentation plane utlisée si l'altération linéaire de la zone est de $- 14\,cm/km$.
\eex
\chapter{\textit{\textbf{Les  Représentations  Planes}}}
\index{Représentations planes}
\section{\textsc{Introduction}}
On a vu dans les chapitres précédents qu'un point géodésique est représenté par ses coordonnées géodésiques $(\varphi,\lambda, he)$ dans un système géodésique donné relatif à un ellipsoïde donné.
\\

Les calculs géodésiques sur l'ellipsoïde étant compliqués d'une part et que l'homme a cherché toujours à visualiser le monde extérieur où il vit par des graphiques et des plans représentés sur des surfaces planes d'autre part.
\\

Le géodésien, par le moyen des représentations planes appelées incorrectement projections, donne une représentation du modèle terrestre (sphère ou ellipsoïde) sur le plan où il est plus facile d'effectuer les calculs d'angles, de distances et de gisements.
\\

Cependant, ces représentations subissent des déformations dues aux propriétés géométriques des surfaces modèles et images. Le problème fondamental ici est de définir des représentations minimisant ces déformations compte tenu d'un objectif déterminé.
\\

Dans la suite du cours, on étudie les représentations dites conformes en général et plus en détail: la représentation conique Lambert et la représentation UTM (Universal Transverse Mercator).
\\

On va établir une correspondance entre les points d'une surface modèle $\sigma $ et les points d'une surface image $\Sigma$, dans le cas particulier où:

-	la surface $\sigma$ est sphérique ou ellipsoïdique;

-	la surface $\Sigma$ est plane.
\section{\textsc{Eléments correspondants}}
\begin{figure}[hp]
	\centering
		\includegraphics[width=0.90\textwidth]{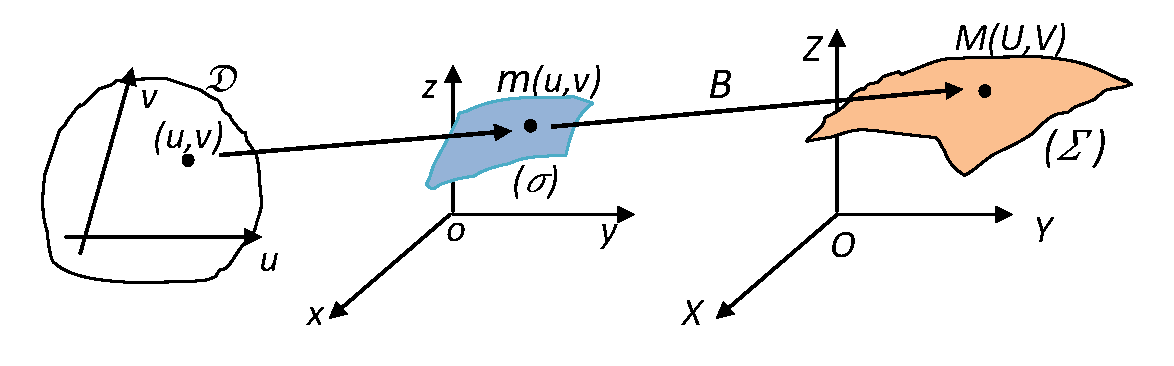}
	\caption{Représentation plane}
	\label{fig:replane}
\end{figure}

 \index{Eléments correspondants}
Représenter la surface $ \sigma $ sur $\Sigma$ consiste à définir une bijection \textit{B} de $\sigma \Longrightarrow \Sigma$ (\textbf{Fig. \ref{fig:replane}}):
 
à $m(u,v) \in (\sigma) \Longrightarrow \,  M(U,V) \in (\Sigma)$
avec:

   $(u,v) \in \mathcal D\subset \BbR^2,\,\, U= U(u,v)$, $V= V(u,v)$ et $\textbf{\textit{OM}} = B(\textbf{\textit{om}}) $                                                                                     

$(u,v)$ les paramètres qui définissent la surface $(\sigma)$ et $U,V$ sont ceux de la surface $(\Sigma)$.
\\

Les points $m(u,v)$ et $M(U,V)$ sont appelés points correspondants\index{Points correspondants}. Si le point $m$ décrit une courbe $(\gamma)$  sur $(\sigma)$, son image $M$ décrit une courbe $(\Gamma)$, on dit que les courbes $(\gamma)$  et $(\Gamma)$  sont dites courbes correspondantes\index{Courbes correspondantes}.
\\

\begin{figure}
	\centering
		\includegraphics[width=0.90\textwidth]{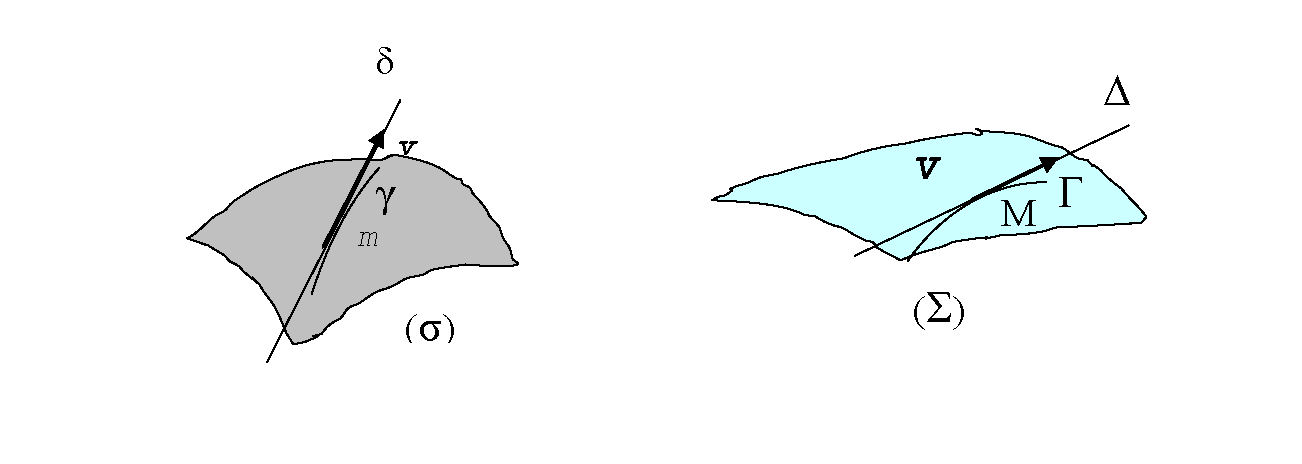}
	\caption{Tangentes correspondantes}
	\label{fig:tangcorresp}
\end{figure}

De même, on appelle tangentes correspondantes\index{Tangentes correspondantes}, les tangentes à deux courbes correspondantes en deux points correspondants (\textbf{Fig. \ref{fig:tangcorresp}}).

\begin{figure}[htbp]
	\centering
		\includegraphics[width=0.90\textwidth]{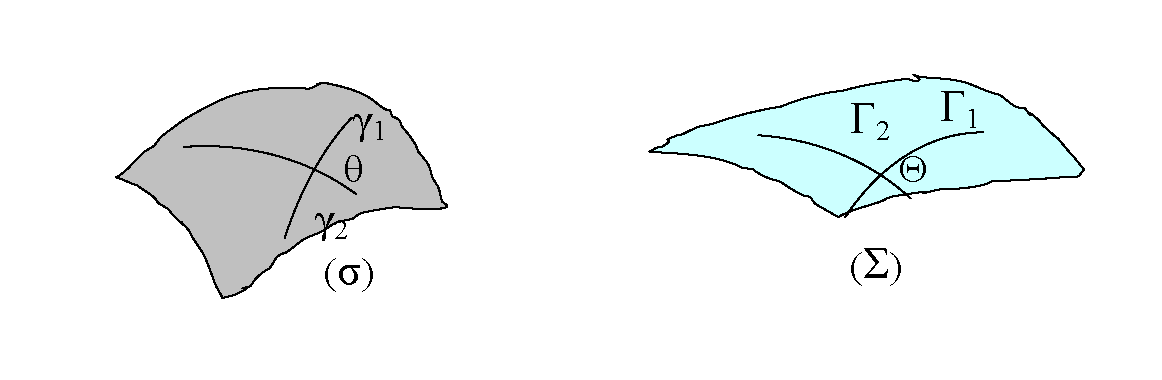}
	\caption{Angles correspondants}
	\label{fig:anglecorresp}
\end{figure}

L'angle de deux tangentes à deux courbes sur $(\sigma)$ et l'angle des tangentes correspondantes sont dites angles correspondants\index{Angles correspondants} (\textbf{Fig. \ref{fig:anglecorresp}}).

\section{\textsc{Canevas}}\index{Canevas}
Les représentations sont différenciées par deux aspects qui sont :

-	la nature des courbes coordonnées du modèle et celles de l'image qui définissent le caractère du canevas;

-	le type de l'altération: longueurs et/ou angles et/ou surfaces.
\bdf
On appelle canevas les images des courbes coordonnées du modèle.
\edf
Pour passer au plan, on peut considérer le passage du modèle ellipsoïdique au plan ou celui du modèle ellipsoïdique au plan via le modèle sphérique :
$$ \left\{\begin{array}{ll}
	 \mbox{\textbf{l'ellipsoïde}} \Longrightarrow   \mbox{\textbf{au plan}} \\
\\
 \mbox{\textbf{l'ellipsoïde}} \Longrightarrow \mbox{\textbf{à la sphère}} \Longrightarrow \mbox{\textbf{au plan}}
 \end{array}\right. $$
 Les représentations peuvent être classées selon la nature des courbes coordonnées $(u,v)$ et $(U,V)$. Pour le modèle sphérique, les courbes coordonnées $(u,v)$ déterminent toujours deux familles de courbes orthogonales, méridiens et parallèles ou pseudo-méridiens et pseudo-parallèles.
\\

\begin{figure}
	\centering
		\includegraphics[width=0.90\textwidth]{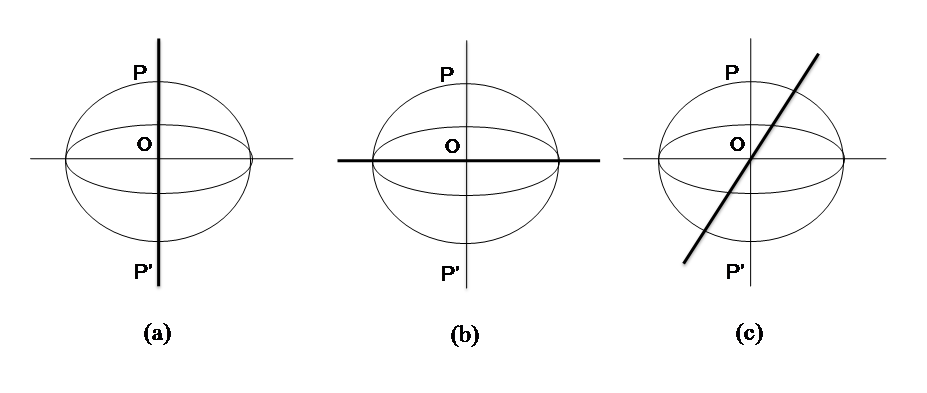}
	\caption{Types de représentation}
	\label{fig:pivotd}
\end{figure}

Soit $DD'$ le diamètre de référence du modèle, le point $D$ est le pivot\index{Pivot} de la représentation (\textbf{Fig. \ref{fig:pivotd}}).
\\

La représentation est dite :

-	directe\index{Représentation directe}, si le diamètre de référence est choisi sur la ligne des pôles $PP'$;

-	transverse\index{Représentation transverse} si le diamètre de référence est perpendiculaire à $PP'$;

-	oblique\index{Représentation oblique}: les autres cas, le pivot n'est ni pôles, ni sur l'équateur.
Quant aux courbes coordonnées\index{Courbes coordonnées} de l'image plane, elles sont :

-	soit deux familles de droites perpendiculaires; $U,V$ sont alors les coordonnées cartésiennes $X$ et $Y$ du plan;

-	soit une famille de droites concourantes $(\Omega)$ et la famille de cercles orthogonaux  $(R)$ définissant des coordonnées polaires.
\section{\textsc{Les Représentations Cylindriques}}
Les représentations cylindriques\index{Représentation cylindrique} sont définies par les représentations ayant comme courbes coordonnées images les coordonnées cartésiennes $X,Y$ correspondantes aux courbes coordonnées du modèle.
\\

Leurs équations sont  pour les représentations cylindriques quelconques :
\be
\begin{array} {l}
	X=X(u) \\
	Y=Y(v)
\end{array} \Longleftrightarrow \begin{array}{l}
 	u = u(X) \\
	v = v(Y)
\end{array}
\ee
Cas des représentations cylindriques directes : les paramètres du modèle sont $\varphi,\lambda$  respectivement la latitude et la longitude et les équations de la représentations sont de la forme :
\be
\begin{array} {l}
	X=X(\varphi) \\
	Y=Y(\lambda)
\end{array} \Longleftrightarrow \begin{array} {l}
 	\varphi = \varphi(X) \\
	\lambda = \lambda(Y)
\end{array}\quad ou \quad  \begin{array} {l}
	X=X(\lambda) \\
	Y=Y(\varphi)
\end{array} \Longleftrightarrow \begin{array} {l}
 	\lambda = \lambda(X) \\
	\varphi = \varphi(Y)
\end{array}
\ee
Cas des représentations cylindriques transverses : les paramètres du modèle sont les coordonnées de Cassini-Soldner\index{Coordonnées de Cassini-Soldner} $(L,H)$ et les équations sont de la forme :
\be
\begin{array} {l}
	X=X(L) \\
	Y=Y(H)
\end{array} \Longleftrightarrow \begin{array} {l}
 	L = L(X) \\
	H = H(Y) 
\end{array}\quad ou \quad \begin{array} {l}
	Y=Y(L) \\
	X=X(H)
\end{array} \Longleftrightarrow \begin{array} {l}
 	L = L(Y) \\
	H = H(X) 
\end{array}
\ee
En remplaçant $L$ et $H$ en fonction  de $\varphi$ et de $\lambda$, on obtient :
\be
\begin{array} {l}
	X=X(\varphi,\lambda) \\
	Y=Y(\varphi,\lambda) 
\end{array} \Longleftrightarrow \begin{array} {l}
 	\varphi = \varphi(X,Y) \\
	\lambda = \lambda(X,Y) 
\end{array}
\ee
\section{\textsc{Les Représentations Coniques et Azimutales}}
Ce sont les représentations planes telles que les courbes coordonnées images sont définies par les coordonnées polaires $R, \Omega$. Les courbes $R =$ Constante et $\Omega = $ constante sont les courbes correspondantes des courbes coordonnées $u$ et $v$ du modèle.
\\

Les équations générales de ces représentations sont de la forme :
\be
\begin{array} {l}
	\Omega=\Omega(u) \\
	R=R(v)
\end{array} \Longleftrightarrow \begin{array} {l}
 	u = u(\Omega) \\
	v = v(R)
\end{array}
\ee
 $u =$ cte $\Longrightarrow \,\Omega =$ cte  $\Longrightarrow$ les images de $u =$ cte sont des droites concourantes.\\
                 
$ v =$ cte $\Longrightarrow \, R =$ cte $\Longrightarrow$ les images de $v =$ cte sont des arcs de cercles concentriques.
\\

Parmi les représentations coniques\index{Représentation conique}, on trouve un groupe particulier de représentations où l'angle $\Omega$ est égal à l'angle du cercle méridien correspondant soit l'azimut  de la tangente au méridien au pôle $D$ de la représentation, ces représentations sont dites \textbf{représentations azimutales}\index{Représentation azimutale} :
\be
\begin{array} {l}
	\Omega=Az \\
	R=R(v)
\end{array} \Longleftrightarrow \begin{array} {l}
 	Az = \Omega \\
	v = v(R)
\end{array}
\ee
Les représentations \textbf{coniques directes} ont leurs équations comme suit :
\be
	\begin{array} {l}
R=R(\varphi) \\
\Omega=\Omega(\lambda) 
\end{array} \Longleftrightarrow \begin{array} {l}
\varphi=\varphi(R) \\
\lambda =\lambda(\Omega)
\end{array}
\ee
\section{\textsc{Les Altérations}}
\subsection*{9.6.1. L'Altération Angulaire}\index{Altération angulaire}
\bdf
On appelle altération angulaire la différence des deux angles correspondants soit:
\be 
	\fbox{$ \Theta -\theta $}
\ee
\edf
\subsection*{9.6.2. Le Module Linéaire dans une direction $\delta$}\index{Module linéaire}
Soit $\delta$ la direction de la tangente en un point donné $m$ du modèle $(\sigma)$,  $s$ et $S$ les abscisses curvilignes sur les 2 courbes correspondantes $(\gamma)$  et $(\Gamma)$.
\bdf
On appelle  module linéaire dans la direction $\delta$ le rapport :
\be
	\fbox{ $ m_\delta (u,v)=\ds \frac{dS}{ds}=\frac{\|\textbf{\textit{dM}}(U,V)\|}{\|\textbf{\textit{dm}}(u,v)\|}=\frac{\|\textbf{\textit{V}}\|}{\|\textbf{\textit{v}}\|} $}
\ee
où $\textbf{\textit{V}}$ est l'image du vecteur $\textbf{\textit{v}}$ unitaire dans la direction $\delta$.
\edf
En utilisant les éléments de la 1ère forme fondamentale\index{Première forme fondamentale} des surfaces $(\sigma)$  et $(\Omega)$, on a alors:
$$ 	m^2_\delta=\left(\frac{dS}{ds}\right)^2=\frac{dS^2}{ds^2}=\frac{EdU^2+2FdUdV+GdV^2}{edu^2+2fdudv+gdv^2} $$
soit:
\be
	m_\delta=\sqrt{\frac{EdU^2+2FdUdV+GdV^2}{edu^2+2fdudv+gdv^2}}
\ee
\subsection*{9.6.3. L'Altération Linéaire}\index{Altération linéaire}
\bdf
On appelle altération linéaire dans la direction $\delta$ la quantité sans unité:
\be 
	 \fbox{ $ \epsilon= m_\delta - 1 $}
\ee
\edf
\subsection*{9.6.4. Le Module aréolaire}\index{Module aréolaire}
Soient $da(\sigma)$  et $dA(\Sigma)$ des aires de domaines limités par des contours correspondants, alors:

\bdf
Le module aréolaire ou rapport des aires est donné par :
\be
	\fbox{ $ m_a=\ds \frac{dA(\Sigma)}{da(\sigma)}=\sqrt{\frac{EG-F^2}{eg-f^2}} $}
\ee
\edf
\section{\textsc{Indicatrice de Tissot}}\index{Indicatrice de Tissot}\index{\textbf{Tissot N.A.}}
\subsection*{9.7.1. Le Lemme de Tissot \footnote{\textbf{Nicolas Auguste Tissot} (1824 - 1897): cartographe français.}}
\begin{figure}
	\centering
		\includegraphics[width=0.90\textwidth]{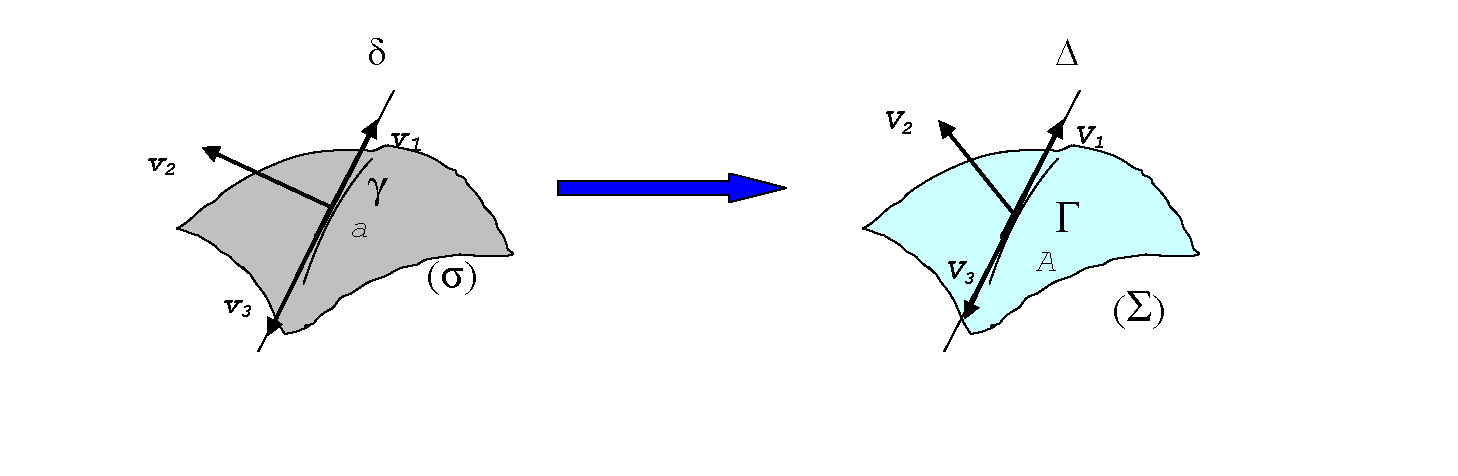}
	\caption{Les directions principales}
	\label{fig:lemmetissot}
\end{figure}

\blm
\index{Lemme de Tissot}
En 2 points correspondants, il existe au moins deux 2 vecteurs $\textbf{\textit{v}}_1$ et $\textbf{\textit{v}}_2$  orthogonaux sur $(\sigma)$ admettant deux vecteurs $\textbf{\textit{V}}_1$ et $\textbf{\textit{V}}_2$ correspondants  orthogonaux sur $(\Sigma)$ (\textbf{Fig. \ref{fig:lemmetissot}}).
\elm
\textit{Le couple de directions correspondantes orthogonales à la fois sur la surface image $(\Sigma)$  et sur la surface modèle $(\sigma)$   sont appelées directions principales \index{Directions principales} (au sens de Tissot).}

\begin{figure}
	\centering
		\includegraphics[width=0.95\textwidth]{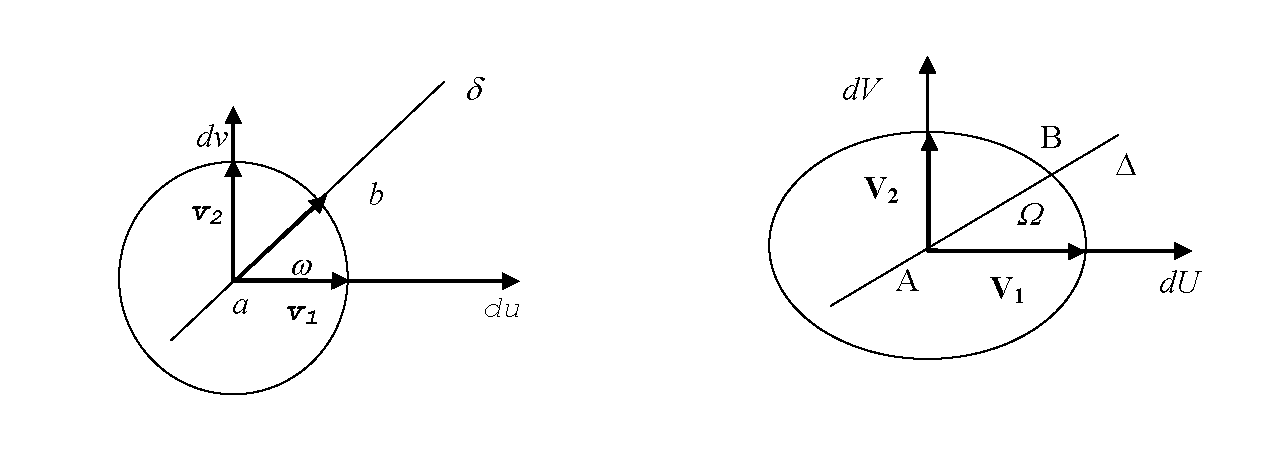}
	\caption{Indicatrice de Tissot}
	\label{fig:indicatricetissot}
\end{figure}

Soient $\delta_1$ et $\delta_2$ les directions  principales sur  $(\sigma)$. Les modules linéaires dans les directions $\delta_1$ et $\delta_2$ sont dits modules principaux:
\ba
	            m_{\delta1} = m_1  \nonumber \\
	           m_{\delta2}  = m_2 \nonumber
\ea
\textbf{Indicatrice de Tissot:}
Soient $a$ un point de $(\sigma)$, dans le plan tangent à $(\sigma)$ au point $a$ et le repère orthonormé $(a,du,dv)$ de vecteurs unitaires $\textbf{v}_1$ et $\textbf{v}_2$.
Soit $b$ un point voisin de $a$ tel que $\|\textbf{\textit{ab}}\| = 1$. Au repère $(a,du,dv)$ correspond le repère orthonormé $(A,dU,dV)$ sur la surface $\Sigma$, et au point $b$ correspond le point $B$ (\textbf{Fig. \ref{fig:indicatricetissot}}).

Par définition, on a :  
$$	m_\delta =\frac{\|\textbf{\textit{AB}}\|}{\|\textbf{\textit{ab}}\|}=\|\textbf{\textit{AB}}\| $$
\newpage
et:
\ba
	 m_1 =\frac{\| \textbf{\textit{V}}_1 \|}{\| \textbf{\textit{v}}_1 \|}=\| \textbf{\textit{V}}_1 \| \nonumber	\\
  m_2=\frac{\|\textbf{\textit{V}}_2 \|}{\|\textbf{\textit{v}}_2 \|}=\|\textbf{\textit{V}}_2 \| \nonumber
\ea
Par les définitions des modules linéaires $m_1$ et $m_2$, on peut écrire:
\be
\begin{array}{l}
 m_1=\frac{\|\textbf{\textit{AB}}\|cos\Omega}{\|\textbf{\textit{ab}}\|cos\omega}=\ds \frac{m_\delta cos\Omega}{1.cos\omega}\Rightarrow m_1 cos\omega=m_\delta cos\Omega   \\
\\
 m_2 =\frac{\|\textbf{\textit{AB}}\|sin\Omega}{\|\textbf{\textit{ab}}\|sin\omega}=\ds \frac{m_\delta sin\Omega}{1.sin\omega}\Rightarrow m_2 sin\omega=m_\delta sin\Omega 
\end{array} \label{fd2}
\ee
Les coordonnées du point $B$ dans $(A,dU,dV)$ sont donc:
\be
\begin{array}{l}
	dU = ABcos\Omega  = m_\delta.cos\Omega = m_1.cos\omega  \\
  dV = ABsin\Omega  = m_\delta.sin\Omega = m_2.sin\omega 
	\end{array}
\ee
\bthm
(\textbf{Indicatrice de Tissot}) Quand le point $b$ varie c'est-à-dire $\omega$ varie, le point $B$ image de $b$ décrit une ellipse d'équation:
\be
	\fbox{ $ \ds \frac{dU^2}{m_1^2}+\frac{dV^2}{m_2^2}=1 $}
\ee
Cette ellipse est appelée \textbf{indicatrice de Tissot}.\index{Indicatrice de Tissot}
\ethm
Elle est l'image du cercle de rayon unité dans le plan tangent au point $a$ de $(\sigma)$. Les longueurs des demis grands et petits axes sont les modules principaux\index{Modules principaux} $m_1$ et $m_2$. La longueur d'un demi-diamètre est le module linéaire dans la direction $\delta$ soit  $m_\delta$.

Dans le cas général, il existe un seul couple de vecteurs orthogonaux correspondants.

\textbf{Corollaire 9.1} \textit{S'il y'a une infinité de couples de vecteurs orthogonaux correspondants, l'indicatrice de Tissot est un cercle quelque soit la direction $\delta$ et le module linéaire est indépendant de la direction:}
\be
	\fbox{ $ m_\delta = m_1 = m_2 = m $}
\ee
\subsection*{9.7.2. Altération Angulaire}
L'altération angulaire est donnée par $\Omega - \omega$. Or d'après les coordonnées de B données par les équations (\ref{fd2}), on a :
	\[tg\Omega=\frac{m_2}{m_1}tg\omega
\]
On calcule:
$$\ds 	\frac{tg\Omega -tg\omega}{tg\Omega+tg\omega}=\frac{\left(\ds \frac{m_2}{m_1}-1\right)tg\omega}{\left(\ds \frac{m_2}{m_1}+1\right)tg\omega}=\frac{m_2 - m_1}{m_2 +m_1} $$
Or:
$$ \ds	\frac{tg\Omega -tg\omega}{tg\Omega+tg\omega}=\frac{sin(\Omega-\omega)}{sin(\Omega+\omega)}	=\frac{m_2 - m_1}{m_2 +m_1} $$
Si l'altération angulaire est nulle $\Longrightarrow \, \Omega-\omega=0\Longrightarrow\,\Omega=\omega$. D'où $m_2 = m_1$ et l'indicatrice de Tissot est un cercle. Par suite, on a l'équivalence :
\be
	\fbox{ $ \mbox{\textbf{Altération angulaire nulle}}  \Leftrightarrow  m_2 = m_1 \, \mbox{et pour toute}\, \delta \,\, m_\delta = m = \mbox{cte} $}
\ee
La représentation est dite dans ce cas \textbf{conforme}\index{Représentation conforme}.
\section{\textsc{Les Représentations Planes et les Fonctions Analytiques}}
\subsection*{9.8.1. Rappels Mathématiques}
On considère le plan complexe tel que à un point de coordonnées $(x,y)\in \BbR^2$, on associe le nombre complexe $z=x+iy \in \BbC$ et on peut écrire:
\be 
	z=\left|z\right|e^{i.arg(z)}=\left|z\right|(cos(arg(z))+isin(arg(z)))
\ee
où  $\left|z\right|$ est le module \index{Module} et $arg(z)$ est l'argument \index{Argument} du nombre complexe $z$ défini à $2k\pi$ avec:
\be
\begin{array}{l}
	\left|z\right|=\sqrt{x^2+y^2}\\
	tg(arg(z))=\ds \frac{y}{x},\quad x \neq 0
\end{array}
\ee
\subsubsection*{9.8.1.1. Logarithme Complexe}
Soit $t$ un nombre complexe donné, on cherche tous les nombres complexes $z$ tels que $e^z=t$. Il n'en existe que si $t\neq 0$. Supposons $t\neq 0$, on a alors:
	\[z=x+iy\Longrightarrow e^z=t\Longrightarrow e^{x}e^{iy}=t=\left|t\right|e^{i.arg(t)}\Rightarrow x=Log\left|t\right|\,\,\,\mbox{et}\,\, y=arg(t)\,\mbox{à}\,\, 2k\pi \,\,\mbox{près}
\]
On a donc:
$$ 	z=x+iy=Log\left|t\right| +i.arg(t) $$
Par définition, on pose:
\be 
\fbox{ $ 	z=Logt=Log\left|t\right| +i.arg(t) $}
\ee
\bdf
On appelle détermination de $Logt$ dans un ouvert connexe $\m D$ du plan complexe toute fonction $g$ continue de $t$ définie dans $\m D$ et telle que:
\be 
\fbox{ $	\forall t \in \m D ,\, \, e^{g(t)}=t $}
\ee
\edf
 

\subsection*{9.8.2. Fonction analytique}
A tout nombre complexe $z = x + iy$ on peut lui faire correspondre un nombre complexe $Z = X + iY$ par l'intermédiaire d'une fonction $f$. On note, en appelant $P$ et $Q$ les parties réelle et imaginaire de $f$ :
\be
	                                           Z = f(z) = P(x,y) + iQ(x,y)   \lb{fona}                                             
\ee
Cette correspondance entre $z$ et $Z$ définit une représentation d'un plan $(p)$ sur un plan $(P)$ dans laquelle le point $A$ d'affixe \index{Affixe} $Z$ du plan $(P$) est l'image d'un point $a$ d'affixe $z$ du plan $(p)$.
\begin{figure}[h]
	\centering
	\includegraphics{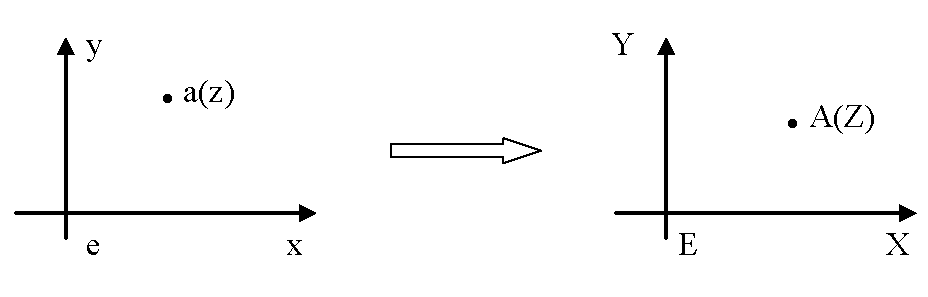}
\caption{Correspondance}
\label{fig:affixe}
\end{figure}
L'extension des propriétés concernant les limites et la continuité, pour la fonction $f$, se déduit immédiatement des propriétés analogues concernant les fonctions $P$ et $Q$  des deux variables $(x,y)$.

Pour étendre à la fonction $f$ la notion de dérivée, il faut étudier la limite, lorsque $z\longrightarrow   0$, du rapport $ \displaystyle {\frac{dZ}{dz}}$, qui s'écrit :
$$\ds \frac{dZ}{dz}=\frac{dX+idY}{dx+idy}=\frac{P'_xdx+P'_ydy+i(Q'_xdx+Q'_ydy)}{dx+idy}=\frac{(P'_x+iQ'_x)dx+(P'_y+iQ'_y)dy}{dx+idy}   $$            
Ce rapport dépend en général de $\displaystyle \frac{dy}{dx}$, sa limite dépend de la manière dont $dz$ tend vers zéro, ou encore de la façon dont le point $a'$ d'affixe $z + dz$ tend vers le point $a$ d'affixe $z$ : si $a'$ tend vers $a$ en décrivant une spirale dont $a$ est le point asymptotique, par exemple, la limite n'existe pas.

Mais le rapport $\displaystyle \frac{dZ}{dz}$ est une fonction homographique \index{Fonction homographique} de $(dx, dy)$. La limite quand $dz \longrightarrow  0$, est indépendante  de la façon dont $dz\longrightarrow   0$, c'est-à-dire  dont $dx$ et $dy$ tendent (indépendamment) vers 0, si :  
\be
P'_x+iQ'_x=\frac{P'_y+iQ'_y}{i}
\ee
C'est-à-dire si :       
\be
\fbox{$ P'_x=Q'_y \quad \mbox{et} \quad P'_y=-Q'_x $ }\label{eqcau}
\ee
ou en revenant à $X$ et $Y$:
\be
\fbox{ $ \mbox{Conditions de Cauchy}\left\{ \begin{array}{l}
 \displaystyle \frac{\partial X}{\partial x}=\frac{\partial Y}{\partial y}  \\
\\
  \displaystyle \frac{\partial X}{\partial y}=-\frac{\partial Y}{\partial x}  
  \end{array}\right. $} \lb{926}
\ee
relations connues sous le nom de \underline{\textbf{conditions de Cauchy}}\footnote{\textbf{Augustin-Louis Cauchy} (1789-1857): mathématicien français.}.\index{Conditions de Cauchy}
Lorsque ces conditions sont satisfaites, la fonction $f$ admet, en tout point de son domaine de définition, une dérivée notée :
\be 
	\fbox{ $   f'(z)= \ds \frac{df}{dz}=\frac{dZ}{dz}=P'_x+iQ'_x=Q'_y-iP'_y $ }
\ee
La fonction $f$ est dite \textbf{analytique}\index{Fonction analytique}.
\\

\subsubsection*{9.8.2.1. Autre définition de la fonction analytique}
Soit $F(x,y)$ une fonction complexe donc une application de $\BbR \times \BbR \longmapsto \BbC$. On suppose que $F$ soit différentiable entraîne que:
\be 
	dF(x,y)=\frac{\partial F(x,y)}{\partial x}dx+\frac{\partial F(x,y)}{\partial y}dy  \label{eqs14}
\ee
On pose:
\ba
	z=x+iy \nonumber \\
	\bar{z} = x-iy \nonumber
\ea
qui sont 2 fonctions différentiables en x et y, d'où:
\ba
	dz=dx+idy \nonumber \\
	d\bar{z}=dx-idy \nonumber
\ea
Par suite:
\ba
	dx=\frac{1}{2}(dz+d\bar{z}) \\
	dy=\frac{1}{2i}(dz-d\bar{z})
\ea
On les remplace dans l'équation (\ref{eqs14}), on obtient:
$$ 	dF(x,y)=\ds \frac{1}{2}\left(\frac{\partial F(x,y)}{\partial x} - i \frac{\partial F(x,y)}{\partial y}\right)dz+\frac{1}{2}\left(\frac{\partial F(x,y)}{\partial x} + i \frac{\partial F(x,y)}{\partial y}\right)d\bar{z} $$
En remplaçant $x$ par $(z+\bar{z})/2$ et $y$ par $(z-\bar{z})/2i$, $F(x,y)$ devient une fonction $\emph{G}(z,\bar{z})$, ce qui donne en posant:
\ba
	\ds \frac{\partial}{\partial z}=\frac{1}{2}\left(\frac{\partial}{\partial x}-i\frac{\partial}{\partial y}\right) \nonumber \\
	\ds 	\frac{\partial}{\partial \bar{z}}=\frac{1}{2}\left(\frac{\partial}{\partial x}+i\frac{\partial}{\partial y}\right) \nonumber \\
		dF(x,y)=d\emph{G}(z,\bar{z})=\ds \frac{\partial \emph{G}(z,\bar{z})}{\partial z}dz+\frac{\partial \emph{G}(z,\bar{z})}{\partial \bar{z}}d\bar{z}
\ea
On considère maintenant le cas où la fonction $F$ est analytique. De l'équation précédente, on a:
$$ dF(x,y)=d\emph{G}(z,\bar{z})=\ds \frac{\partial \emph{G}(z,\bar{z})}{\partial z}dz+\frac{\partial \emph{G}(z,\bar{z})}{\partial \bar{z}}d\bar{z}=F'(z)dz $$
C'est-à-dire:
\be 
	\frac{\partial \emph{G}(z,\bar{z})}{\partial \bar{z}}=0
\ee
On vérifie bien si ce terme là est nul:
\be 
	\frac{\partial \emph{G}(z,\bar{z})}{\partial \bar{z}}= \frac{1}{2}\left(\frac{\partial}{\partial x}+i\frac{\partial}{\partial y}\right)\emph{G}=\frac{1}{2}\left(\frac{\partial \emph{G}}{\partial x}+i\frac{\partial \emph{G}}{\partial y}\right) \label{eq425}
\ee
On introduit la notation suivante:
$$ F(x,y)=G(z,\bar{z})=P(x,y)+iQ(x,y)=P+iQ$$
On a donc de (\ref{eq425}):
\ba 
	\frac{\partial \emph{G}(z,\bar{z})}{\partial \bar{z}}= \frac{1}{2}\left(\frac{\partial \emph{G}}{\partial x}+i\frac{\partial \emph{G}}{\partial y}\right)=\frac{1}{2}\left(\frac{\partial (P+iQ)}{\partial x}+i\frac{\partial (P+iQ)}{\partial y}\right)=\nonumber \\
	\frac{1}{2}\left(\frac{\partial P}{\partial x}+i\frac{\partial Q}{\partial x}+i\frac{\partial P}{\partial y}-\frac{\partial Q}{\partial y}\right) \nonumber
\ea
Comme la fonction $F$ ou $G$ est analytique, et en utilisant les conditions de Cauchy (\ref{eqcau}), on a finalement:
\be 
	\frac{\partial \emph{G}(z,\bar{z})}{\partial \bar{z}}=0
\ee
\bdf
Une fonction $f(z,\bar{z})$ est analytique en $z$ si elle ne dépend que de $z$ soit:
\be 
	\fbox{$ \displaystyle \frac{\partial f(z,\bar{z})}{\partial \bar{z}}=0 $ }
\ee
\edf
Avec les notations de l'équation (\ref{fona}), au point $a$ du plan $(p)$, quelle que soit la direction du vecteur $\textbf{\textit{a}}_0\textbf{\textit{a}}$, d'affixe $dz$, on peut écrire :
$$ dZ=f'(z_0)dz=| f'(z_0)|e^{i.argf'(z_0)}dz $$
relation qui exprime que le vecteur $\textbf{\textit{A}}_0\textbf{\textit{A}}$ se déduit du vecteur $\textbf{\textit{a}}_0\textbf{\textit{a}}$ par une similitude, dont le rapport est $| f'(z_0)|$ et l'angle $argf'(z_0)$.
\\

\begin{flushleft}
		\begin{figure}[htp]
	\includegraphics{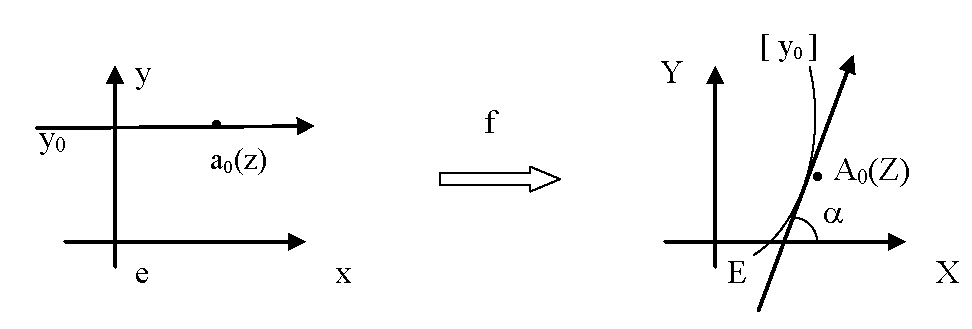}
	\caption{Image de $y=y_0$}
	\label{imagemeridien}
		\end{figure}
   		\end{flushleft}

La représentation du plan $(p)$ sur le plan $(P)$ est donc conforme. On peut écrire :
 \be
 dS=|dZ| \quad \mbox{et} \quad ds=  |dz|                                                                                        \ee
Le module linéaire de la représentation est :
\be
\begin{array}{l}
m=\left| \ds \frac{dZ}{dz}\right| \\
\mbox{et} \quad arg\left(\ds \frac{dZ}{dz}\right)  =Arctg\left(\ds \frac{\displaystyle \frac{\partial Y}{\partial x}}{\displaystyle \frac{\partial X}{\partial x}}\right)=\alpha 
\end{array} 
\ee                                                 
 $\alpha$ angle de la tangente en $A_0$ à l'image $y_0$  de la droite $y = y_0$ du plan $(p)$ (\textbf{Fig. \ref{imagemeridien}}). En effet, cette image est définie en fonction du paramètre $x$, par :
\be
\left\{\begin{array}{l} 
	      X = X(x, y_0)   \\   
				Y= Y(x, y_0) 
				\end{array} \right. 
  \ee  
Les conditions de Cauchy données par (\ref{926}) se traduisent aussi en disant que les fonctions $P$ et $Q$ sont des fonctions \underline{harmoniques}, c'est-à-dire qu'elles satisfont chacune \underline{l'équation de Laplace} :\index{Fonction harmonique}
\be
\fbox{ $ \begin{array}{ll}
 \displaystyle \frac{\partial^2 X}{\partial x^2}+\frac{\partial^2 X}{\partial y^2}=0 \\
 \displaystyle \frac{\partial^2 Y}{\partial x^2}+\frac{\partial^2 Y}{\partial y^2}=0 
  \end{array} $} \label{dfff}
\ee 
On démontre les propriétés suivantes (\textit{J. Dieudonné\footnote{\textbf{Jean Dieudonné} (1906 - 1992): mathématicien français.}}, 1970):\index{\textbf{Dieudonné J.}}

- si $f$ est une fonction analytique, elle admet des dérivées de tous les ordres : elle est donc développable en série entière en tout point de son domaine de définition;

- une fonction analytique est déterminée dans tout son domaine d'existence, si elle est définie dans une région, aussi petite qu'on la suppose, entourant un point $a$; ou même tout le long d'un arc de courbe, aussi petit qu'on le suppose, aboutissant au point $a$.

En effet, la connaissance de $f$ au voisinage de $a$ permet (théoriquement tout au moins) de former la suite des dérivées de $f$ au point $a$, donc, d'écrire son développement en série de Taylor. Si $b$ est un point intérieur au cercle de convergence de cette série, on peut alors calculer les dérivées successives de $f$ au point $b$, et ainsi de suite. L'opération est dite prolongement analytique de $f$.
\subsection*{9.8.3. Représentation Conforme d'une Surface sur une autre}
\begin{quotation}
\begin{svgraybox}
				Conformal mappings play also a fundamental role in modern physics, namely, in string theory and conformal quantum field theory.
\end{svgraybox}
\begin{flushright}
(\textbf{E. Zeidler\footnote{\textbf{E. Zeidler} (1905 - 2016): éminent physicien allemand.}}, 2011)\index{\textbf{Zeidler E.}}
	\end{flushright}
		\end{quotation}

La représentation d'une surface modèle $(\sigma )$ sur une surface image $(\Sigma )$ est définie en établissant une correspondance entre les coordonnées curvilignes $a(u,v)$ de $(\sigma )$ et $A(U,V)$ de $(\Sigma )$ :
\ba
   U = U(u,v) \label{eqU} \\ 
	V = V(u,v)  \label{eqV}     
\ea
Le module linéaire est alors :
$$ m(u,v)=\ds \frac{dS}{ds}\Rightarrow m^2(u,v)=\frac{dS^2}{ds^2}=\frac{EdU^2+2FdU.dV +GdV^2}{edu^2+2fdu.dv+gdv^2} $$
En remplaçant $dU$ et $dV$ en fonction de $du$ et $dv$ en utilisant (\ref{eqU}) et (\ref{eqV}), on peut écrire $dS^2$ sous la forme:
\be 
	dS^2=\m E du^2+2\m Fdu.dv +\m Gdv^2 \label{eqUV}
\ee
Soient le plan tangent à $(\sigma)$ au point $a(u,v)$, et deux courbes $(\gamma_1)$ et $(\gamma_2)$ passant par $a$ dont les tangentes respectivement à  $(\gamma_1)$ et $(\gamma_2)$ appartiennent au plan tangent. $ \displaystyle (\frac{\partial \textbf{a}}{\partial u},\frac{\partial \textbf{a}}{\partial v} )$ est une base du plan tangent.

Soit la direction de la tangente à $(\gamma_1)$ de direction:
$$ \ds 	\frac{\partial \textbf{a}}{\partial u}du+\frac{\partial \textbf{a}}{\partial v}dv $$
De même, soit la direction de la tangente à $(\gamma_2)$ de direction:
$$ \ds 	\frac{\partial \textbf{a}}{\partial u}\delta u+\frac{\partial \textbf{a}}{\partial v}\delta v $$
La forme fondamentale de $(\sigma)$ est :
$$ ds^2=edu^2+2fdu.dv+gdv^2 $$
En notant $\Omega$ l'angle des deux tangentes en $a$, on a la relation:
$$ (\ds		\frac{\partial \textbf{a}}{\partial u}du+\frac{\partial \textbf{a}}{\partial v}dv).(	\frac{\partial \textbf{a}}{\partial u}\delta u+\frac{\partial \textbf{a}}{\partial v}\delta v)=\|	\ds \frac{\partial \textbf{a}}{\partial u}du+\frac{\partial \textbf{a}}{\partial v}dv\|.\|	\frac{\partial \textbf{a}}{\partial u}\delta u+\frac{\partial \textbf{a}}{\partial v}\delta v\|.cos\Omega $$
En posant:
$$ ds=\|	\frac{\partial \textbf{a}}{\partial u}du+\frac{\partial \textbf{a}}{\partial v}dv\|\Rightarrow ds^2=edu^2+2fdu.dv+gdv^2 $$
et:
$$ \delta s=\|	\frac{\partial \textbf{a}}{\partial u}\delta u+\frac{\partial \textbf{a}}{\partial v}\delta v\|\Rightarrow \delta s^2=e\delta u^2+2f\delta u.\delta v+g\delta v^2 $$
Ce qui donne:
$$ cos\Omega=\frac{edu \delta u+f(du\delta v+\delta u dv)+gdv\delta v}{ds\delta s} $$
On pose:
\be
\left\{\begin{array}{l} 
p=\ds \frac{dv}{du}\Rightarrow dv=pdu \\
\\
q=\ds \frac{\delta v}{\delta u}\Rightarrow \delta v= q \delta u
\end{array}\right.
\ee
D'où:
$$ cos\Omega=\frac{edu \delta u+f(du.q\delta u+\delta u.p du)+g.pdu.q.\delta u}{ds\delta s} $$
ou encore:
$$ 	 cos\Omega=\frac{edu \delta u+f(p+q)du\delta u+gpqdu \delta u}{ds\delta s} $$
En utilisant les notations $p$ et $q$, on a:
\ba
ds^2=(e+2fp+gp^2)du^2 \nonumber \\
\delta s^2=(e+2fq+gq^2)\delta u^2 \nonumber 
\ea
On obtient finalement:
\be
\fbox{ $	 cos\Omega=\ds \frac{e+f(p+q)+gpq}{\sqrt{e+2fp+gp^2}\sqrt{e+2fq+gq^2}} $}
\ee
Sur la surface image $(\Sigma)$, l'expression de $dS^2$ donnée par l'équation (\ref{eqUV}):
$$	dS^2=\m E du^2+2\m Fdu.dv +\m Gdv^2 $$
Soit $\Omega'$ l'angle des tangentes correspondantes à $(\gamma_1)$ et $(\gamma_2)$, on a a aussi:
\be
\fbox{ $	 cos\Omega'=\ds \frac{\m E+\m F(p+q)+\m Gpq}{\sqrt{\m E+2\m Fp+\m Gp^2}\sqrt{\m E+2\m Fq+\m Gq^2}} $}
\ee
Il y'a conservation des angles si $\forall \,p,q$ on a:
\be
	\fbox{ $ cos\Omega=cos\Omega' $} \label{eqcos}
\ee
En particulier si :
$$ q=\ds \frac{\delta v}{\delta u}=0 $$
c'est-à-dire, on prend $\Omega$ l'angle d'une tangente avec la courbe coordonnée $v=constante\Rightarrow \delta v=0$ et cela suffit. Alors l'équation $(\ref{eqcos})$ devient $\forall \,p=\displaystyle \frac{dv}{du}$:
$$ \ds 	\frac{e+fp}{\sqrt{e^2+2efp+egp^2}}=\ds \frac{\m E+\m Fp}{\sqrt{\m E^2+2\m E\m Fp+\m E\m Gp^2}} $$
Elevant au carré, les deux membres de l'équation précédente s'écrivent:
$$ 1+	\frac{f^2p^2-egp^2}{e^2+2efp+egp^2}=1+\frac{\m F^2p^2-\m E \m G p^2}{\m E^2+2\m E\m Fp+\m E\m Gp^2} $$
En éliminant le 1 et simplifiant par $p\neq 0$, on obtient:
$$ \ds 	\frac{f^2-eg}{e^2+2efp+egp^2}=\ds \frac{\m F^2-\m E \m G }{\m E^2+2\m E\m Fp+\m E\m Gp^2} $$
Si:
$$ f^2-eg=0 \Rightarrow \m F^2-\m E \m G=0 $$
On a donc:
$$	cos\Omega=\cos\Omega'=\pm 1 $$
Maintenant on suppose que :
$$ f^2-eg \neq 0 \Rightarrow \m F^2-\m E \m G \neq 0 $$
On doit avoir  $\forall \,p$ le rapport :
$$ \ds 	\frac{e^2+2efp+egp^2 }{\m E^2+2\m E\m Fp+\m E\m Gp^2} $$
égal à $\displaystyle \frac{f^2-eg }{\m F^2-\m E \m G}$.

Pour cela, il faudra donc $\forall \,p$:
\be 
\fbox{ $	\ds \frac{e^2+2efp+egp^2 }{\m E^2+2\m E\m Fp+\m E\m Gp^2}= \ds \frac{f^2-eg }{\m F^2-\m E \m G}=\mu^2(u,v) $} \label{eqconf}
\ee
Ce qui donne (avec $\mu >0$) :
\ba 
& e^2+2efp+egp^2=\mu^2 \m E^2+2\mu^2 \m E\m Fp+\mu^2 \m E\m Gp^2\Rightarrow e^2-\mu^2 \m E^2+2p(ef-\mu^2 \m E \m F)\nonumber & \\ &
+p^2(eg-\mu^2 \m E \m G)=0 & \nonumber
\ea
Soit:
$$ 
\left\{\begin{array}{lll}
e^2-\mu^2 \m E^2=0\Rightarrow \mu^2=\ds \frac{e^2}{\m E^2} \\
ef-\mu^2 \m E \m F=0 \Longrightarrow \mu^2=\ds \frac{ef}{\m E \m F} \\
eg-\mu^2 \m E \m G=0\Rightarrow \mu^2=\ds \frac{eg}{\m E \m G}
\end{array}\right.
$$ 
Comme $\displaystyle \frac{e}{\m E} \neq 0$, on a:
$$ \ds  \frac{e}{\m E}=\frac{f}{\m F}=\frac{g}{\m G}\Rightarrow \mu^2 =\frac{f^2-eg }{\m F^2-\m E \m G} $$
La condition (\ref{eqconf}) est vérifiée si et seulement si :
\be
 \fbox{ $ \ds \frac{e}{\m E}=\frac{f}{\m F}=\frac{g}{\m G}\Rightarrow cos\Omega=cos\Omega' \Rightarrow \left\{\begin{array}{ll}
 \mbox{la représentation est conforme}\\
\mbox{ et le module linéaire} \,\, m=m(u,v)= \ds \frac{1}{\sqrt{\mu(u,v)}} \\
\mbox{est indépendant de la direction}
 \end{array}\right. $}
\ee
\index{Représentation conforme}
- Si les coordonnées $(u,v)$ sont orthogonales $(f=0)$, les deux conditions précédentes à gauche deviennent :
 \be 
\m F=0 \quad \mbox{et} \quad \frac{\m E}{e}=\frac{\m G}{g}
\ee
-	Si les coordonnées $(u,v)$ sont symétriques, alors $e = g$, les conditions de conformité s'écrivent :
\be 
     \m F = 0,  \quad      \m E = \m G    \Longrightarrow    \mbox{les coordonnées}\, U \,\mbox{et}\, V \, \mbox{sont symétriques}
\ee
Or :        
\be 
 \left\{\begin{array}{ll}
	 \m F = 0          \Longrightarrow \,               U'_u U'_v  + V'_u V'_v  = 0   \\                                           \m E =  \m G  \Longrightarrow \,             U'^2_u + V'^2_u = U'^2_v + V'^2_v 
	 \end{array}\right.         \label{156}                      
\ee
De la première relation de (\ref{156}), on tire en supposant $V'_v\neq0$: 
 \be
 V'_u=-\frac{U'_uU'_v}{V'_v} \lb{eqlapp}
\ee
D'où : 
 $$ U'^2_u\left(1+\frac{U'^2_v}{V'^2_v}\right)=U'^2_v+V'^2_v\Rightarrow (U'^2_u-V'^2_v)(U'^2_v+V'^2_v)=0 $$
Soit :      
\be 
	     \fbox{ $  U'_u = \pm V'_v $}       \label{eqlap}                  
\ee
L'équation (\ref{eqlapp}) donne: 
\be 
	   \fbox{ $   U'_v = \mp V'_u  $}     \label{eqlap1}                  
\ee

Les équations (\ref{eqlap}), (\ref{eqlap1}) sont, pour les fonctions $U$ et $V$, les conditions de Cauchy.

La correspondance entre le plan des $(u,v)$ et celui des $(U,V)$ est une représentation conforme (\textit{G. Julia}, 1955)\index{\textbf{Julia G.}} et on peut poser :  
                                $$ Z  = U + iV      ,\quad            z  = u + iv   $$                                  
                                
\underline{Toute fonction analytique $f$ définit une représentation conforme de $(\sigma )$  sur $(\Sigma )$.}

$(u,v)$ et $(U,V)$ étant des coordonnées symétriques, les éléments linéaires $ds$ et $dS$ s'écrivent :
$$ 	ds^2=h^2(u,v)(du^2+dv^2) \quad \mbox{ou}\quad ds=|h(u,v)| |dz| $$
et:
$$ 	dS^2=H^2(U,V)(dU^2+dV^2) \quad \mbox{ou}\quad dS=|H(U,V)| |dZ| $$
 Alors le module linéaire est donné par :                                                    
\be 
	\fbox{ $ m=\ds \left|\frac{H(U,V)}{h(u,v)}\right| \left|\frac{dZ}{dz}\right| $}
\ee

L'argument de  $ \displaystyle {\frac{dZ}{dz}}$  s'interprète d'une manière analogue à ce qui se passe dans la représentation d'un plan sur un plan.
\\

Si en particulier, la surface image est un plan :         
     $$     dS^2= dX^2+ dY^2                            $$
et:  
 \be 
	m=\left|\frac{1}{h(u,v)}\right| \left|\frac{dZ}{dz}\right|
\ee
Si la surface modèle est un modèle de la terre, on a :                           
\be 
	ds^2=r^2(d\lambda^2+dL^2)
\ee
- pour un ellipsoïde :  $r = Ncos\varphi$     et   $L$ est la latitude isométrique;

- pour une sphère :  $ r = acos\varphi$   et  $L$ est la latitude croissante ou latitude de Mercator.
\\

Pour une représentation plane conforme, le modèle linéaire est :
   \be 
	m=\frac{1}{r} \left|\frac{dZ}{dz}\right|
\ee
                                                                                                                
avec $z =   \lambda+ iL$  et  $Z  = X + iY$   ou $z = L + i\lambda$ et $Z  = Y  +  iX$.  Dans ce dernier cas, on a :
\be 
	\frac{dZ}{dz}=\frac{\partial Y}{\partial L}+i\frac{\partial X}{\partial L}
\ee
Donc :   
\be
 arg\left( \frac{dZ}{dz}\right) =Arctg\left( \frac{\displaystyle\frac{\partial X}{\partial L}}{\displaystyle\frac{\partial Y}{\partial L}} \right)=Arctg\left(\frac{dX}{dY}\right) 
 \ee                                                         
qui n'est autre que le gisement de l'image du méridien (\textit{J. Commiot}, 1980).\index{\textbf{Commiot J.}}
\section{\textsc{Les Représentations ou Transformations Quasi-Conformes}}\index{Représentation quasi-conforme}

Dans les paragraphes précédents, on a étudié les représentations de la sphère avec les variables $(L_M,\lambda)$ ou celles de l'ellipsoïde de révolution avec les coordonnées $(L,\lambda)$ vers le plan $(X,Y)$ avec:
\ba
\left\{
\begin{array}{ll}
		X=X(L_M,\lambda) \\
		Y=Y(L_M,\lambda)
\end{array}\right. \label{d1} \\
\mbox{ou} \quad \left\{ 
\begin{array}{ll}
		X=X(L,\lambda) \\
		Y=Y(L,\lambda)
\end{array}  \right. \label{d2} 
\ea
en notant:
\ba
	L_M=Log tg \left( \frac{\pi}{4}+\frac{\varphi}{2} \right) \quad \mbox{la latitude de Mercator} \nonumber  \\
L=Log tg \left( \frac{\pi}{4}+\frac{\varphi}{2} \right)-\frac{e}{2}Log \frac{1+esin\varphi}{1-esin\varphi} \quad \mbox{la latitude isométrique} \nonumber
\ea
En posant:
\ba
	z=L_M+i\lambda\quad (ou \quad z=L+i\lambda) \label{d6} \\
	Z=X+iY \label{d7}
\ea           
on a considéré les représentations conformes (c'est-à-dire qui conservent les angles) ou encore définies par :
\be
	Z=Z(z) \label{d8}
\ee          
avec $Z(z)$ une fonction dite holomorphe de $z$ soit:\index{Fonction holomorphe}
$$ \ds \frac{\partial Z}{\partial \bar{z}}=0 $$
où $\bar{z}$ est le conjugué de $z$ soit $\bar{z}=L_M-i\lambda$ (ou $\bar z=L-i\lambda$).
\bdf
Une fonction $f(z)=Z=Z(z)$ définie et dérivable sur un domaine $\mathcal{D} \subset \BbC$ (l'ensemble des nombres complexes) est dite quasi-conforme si elle vérifie (\textit{L. Bers}, 1977):\index{\textbf{Bers L.}}
	\be
\fbox{ $ \ds \frac{	\partial Z}{\partial \bar{z}}=\mu(z).\frac{	\partial Z}{\partial z} \quad avec \quad 	 |\mu(z)|<1 $} \label{d10}
\ee  
 Le coefficient $\mu$ s'appelle coefficient de Beltrami\footnote{\textbf{Eugeno Beltrami} (1835-1899): mathématicien italien.}.\index{\textbf{Beltrami E.}}\index{Coefficient de Beltrami}
\edf
\subsection*{9.9.1. Développement d'une fonction en un point $z_0$}
 Soit $f$ une fonction quasi-conforme et un point $z_0\in \mathcal{D}$. En écrivant un développement de $f$ au point $z_0$, on a alors:
	$$ f(z)=f(z_0)+(z-z_0)\frac{\partial f}{\partial z}(z_0)+(\bar{z}-\bar{z}_0)\frac{\partial f}{\partial \bar{z}}(\bar{z}_0)+... $$
Par un changement de variables, on peut prendre $z_0=0$, d'où:
$$ f(z)=f(z_0)+\ds z\frac{\partial f}{\partial z}(z_0)+\bar{z}\frac{\partial f}{\partial \bar{z}}(\bar{z}_0)+...  $$
Utilisant (\ref{d10}), l'équation précédente s'écrit en négligeant les termes du deuxième degré:
$$		f(z)=f(z_0)+\ds z\frac{\partial f}{\partial z}(z_0)+\bar{z}\mu(z_0).\frac{\partial f}{\partial z}(z_0)  $$
Donc $f(z)$ s'écrit localement:
\ba
	f(z)=\alpha+\beta z+\gamma \bar{z} \label{d14} \\
	\mbox{où}\quad \alpha,\beta,\gamma \,\, \mbox{des constantes complexes avec}\quad \left|\frac{\gamma}{\beta}\right|<1 \label{d15}
\ea
\subsection*{9.9.2. Etude de la Transformée d'un cercle}
On sait que pour une transformation conforme, l'image d'un cercle autour d'un point est un cercle (ou encore l'indicatrice de Tissot est un cercle). Soit un point $z_0$ qu'on peut prendre égal à 0. Par un changement de l'origine des axes, la fonction $f$ s'écrit:
	$$ f(z)=\beta z +\mu \beta \bar{z} $$
Par abus, on garde la même notation. On considère autour de l'origine $z_0=0$ un point $M(x=a.cos\theta,y=a.sin\theta)$ qui décrit un cercle infiniment petit de rayon $a$. On étudie ci-après l'image du point $M$ par $f$.
\\

De l'équation précédente, on a:
\ba
	z=acos\theta +iasin\theta =ae^{i\theta} \nonumber\\
	\mu=|\mu|e^{ik} \nonumber \\
	\beta=|\beta|e^{il} \nonumber \\
	f(z)=a|\beta|e^{il}(e^{i\theta}+|\mu|e^{i(k-\theta)})\label{d15f}
\ea
Si $\theta_1=\displaystyle \frac{k}{2}= \displaystyle \frac{arg(\mu)}{2}$, on a $z_1=ae^{ik/2}$ et :
\ba
	f(z_1)=	a|\beta|e^{il}e^{ik/2}(1+|\mu|) \nonumber \\
	|f(z_1)|=a|\beta|(1+|\mu|) \label{d15h}
\ea
Maintenant on prend $\theta_2=\displaystyle \theta_1+\frac{\pi}{2}=\displaystyle \frac{k}{2}+\frac{\pi}{2}$, alors $z_2=ae^{i\theta_2}=ae^{ik/2}e^{i\pi/2}=iae^{ik/2}$ et on obtient:
\ba
		f(z_2)=	ia|\beta|e^{il}e^{ik/2}(1-|\mu|)  \nonumber \\
	|f(z_2)|=a|\beta|(1-|\mu|) \label{d15j}
\ea
en tenant compte que $|\mu|<1$.
\\

\begin{figure}
	\centering
		\includegraphics{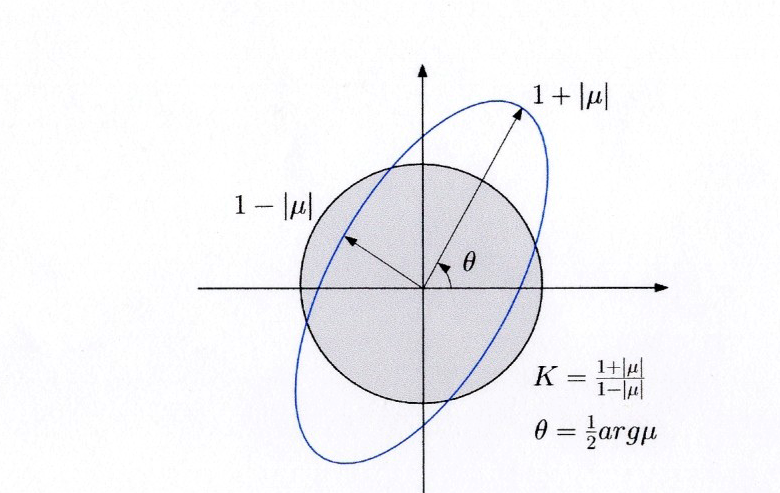}
	\caption{L'image d'un cercle}
	\label{fig:quasiconf}
\end{figure}

Des équations (\ref{d15f}),(\ref{d15h}) et (\ref{d15j}), on déduit que l'image de $M$ décrit une ellipse de demi-grand axe et demi- petit axe respectivement (\textbf{Fig. \ref{fig:quasiconf}}) :
\be
\fbox{ $ \begin{array}{l}
 a'=a|\beta|(1+|\mu|)   \\
 b'=a|\beta|(1-|\mu|) 
\end{array} $} 
\ee
On appelle:
\be
\fbox{ $ 	K=\ds \frac{1+|\mu|}{1-|\mu|} $} \label{d15m}
\ee
\textit{coefficient de distortion ou de dilatation}.\index{Coefficient de distorsion}

\subsection*{9.9.3. Calcul d'un élément de longueur sur le Plan}
Un élément de longueur sur le plan est donné par:
	$$ dS^2=dX^2+dY^2=|df|^2=df.\overline{df} $$
Comme  $df=\beta dz+\gamma d\bar{z}$ et $\overline{df}=\bar{\beta} d\bar{z}+\bar{\gamma} dz$, on a alors:
\ba
	dS^2=dX^2+dY^2=|df|^2=df.\overline{df}=(\beta dz+\gamma d\bar{z})( \bar{\beta} d\bar{z}+\bar{\gamma} dz) \nonumber \\
	=\beta\bar{\beta}dzd\bar{z}+\gamma\bar{\gamma}dzd \bar{z}+ dzd\bar{z} \left(\beta\bar{\gamma}\frac{dz}{d\bar{z}}+\gamma \bar{\beta}\frac{d\bar{z}}{dz} \right) \nonumber 
	\ea
 On pose:
$$ 	ds^2=dz.d\bar{z} $$
Le carré du module linéaire de la transformation quasi-conforme s'écrit:
\be
	m^2=\frac{dS^2}{ds^2}=|\beta|^2+|\gamma|^2+\left(\beta\bar{\gamma}\frac{dz}{d\bar{z}}+\gamma \bar{\beta}\frac{d\bar{z}}{dz} \right) \label{d19}
\ee
Dans l'équation (\ref{d19}), on considère $z=ae^{i\theta}$ varie le long d'un cercle de rayon $a$ infiniment petit et on fait tendre $\theta \longrightarrow 2\pi$. Alors, on obtient :
\ba
	\frac{dz}{d\bar{z}}=\frac{aie^{i\theta} d\theta}{-aie^{-i\theta}d\theta}=-e^{2i\theta}=-1  \nonumber \\
		\frac{d\bar{z}}{dz}=-e^{-2i\theta} =-1 \nonumber 
\ea
L'équation (\ref{d19}) devient:
$$	m^2=\frac{dS^2}{ds^2}=|\beta|^2+|\gamma|^2-\left(\beta\bar{\gamma}+\gamma \bar{\beta} \right) $$
Comme:
$$ \gamma=\mu \beta $$
 on obtient:
\be
		m^2=\frac{dS^2}{ds^2}=|\beta|^2+|\beta|^2|\mu|^2-\left(\beta\bar{\beta}\bar{\mu}+\mu \beta \bar{\beta} \right) \label{d23}
\ee
or $\mu+\bar{\mu}=2|\mu| cosarg(\mu)$, par suite l'équation (\ref{d23}) s'écrit:
\be
	m^2=\frac{dS^2}{ds^2}=|\beta|^2(1+|\mu|^2-2|\mu|cosarg(\mu)) \label{d24}
\ee
Remplaçant $\beta$ par $\displaystyle \frac{\partial f}{\partial z}(z_0)$, (\ref{d24}) devient:
\be
	\fbox{ $ m^2=\ds \frac{dS^2}{ds^2}=\left|\frac{\partial f}{\partial z}(z_0)\right|^2 \left(1+|\mu|^2-2|\mu|cosarg(\mu)\right) $}\label{d25}
\ee
\subsection*{9.9.4. Exemple de Transformation Quasi-conforme}
Lors de passage de coordonnées planes $(X,Y)_i$ d'un système géodésique $S_1$ à des coordonnées planes $(X',Y')_j$ dans un autre système géodésique $S_2$, on utilise souvent une transformation du type:
\be
\left\{\begin{array}{l}
	X'=X_0+aX+bY  \\
	Y'=Y_0+cX+dY 
	\end{array}\right.
\ee
ou encore sous forme matricielle :
\be
	\begin{pmatrix}{
	X' \cr
	Y'}
\end{pmatrix}=\begin{pmatrix}{
	X_0 \cr
	Y_0 }
\end{pmatrix}+\begin{pmatrix}{
	a & b \cr
	c & d }
\end{pmatrix}.\begin{pmatrix}{
	X \cr
	Y }
\end{pmatrix} \label{d27a}
\ee
En posant $Z=X'+iY'$ et $z=X+iY$, on obtient:
\be
	Z=(X_0+iY_0)+X(a+ic)+Y(b+id) \label{d28}
\ee        
On note par :
$$ 	Z_0=X_0+iY_0  $$
Comme $ X=(z+\bar{z})/2$ et $Y=(z-\bar{z})/2i$, alors l'équation (\ref{d28}) s'écrit:
\be
Z=Z_0+z\left(\frac{a+d}{2}+i\frac{c-b}{2}\right)+\bar{z}\left(\frac{a-d}{2}+i\frac{b+c}{2}\right) \label{d30}
\ee         
On pose :
\ba
	\beta=\frac{a+d}{2}+i\frac{c-b}{2} \nonumber  \\
	\gamma=\frac{a-d}{2}+i\frac{b+c}{2} \nonumber 
\ea
 Alors (\ref{d30}) s'écrit:
\be
	Z=Z_0+\beta z+\gamma \bar{z} \label{d33}
\ee
Pour quelles valeurs de $a,b,c,d$ la transformation (\ref{d27a}) est quasi-conforme? En comparant (\ref{d30}) avec (\ref{d14}), il faut que $|\gamma|<|\beta|$ soit:
\ba
	|\gamma|<|\beta|\Rightarrow |\gamma|^2<|\beta|^2 \Rightarrow \frac{(a-d)^2+(b+c)^2}{4}< \frac{(a+d)^2+(c-b)^2}{4} \nonumber \\
	\Rightarrow ad-bc>0 \label{d34}
\ea
C'est-à-dire que le déterminant de la matrice (\ref{d27a}) soit strictement positif.
\\

\textbf{\un{Note historique:}} \textsl{La représentation stéréographique de la sphère au plan est l'une des représentations la plus utilisée depuis l'antiquité (voir exercices n°1 et n°5 ci-dessous). Elle était connue par l'astronome et mathématicien Hipparque (190-120 avant J.C)  ainsi que Claude Ptolémée (80-168). Ce dernier connaît que la représentation stréographique transforme les cercles en cercles ou en droites, mais on ignore s'il savait que l'image de} \textit{tout cercle de la sphère est un cercle ou une droite.} \textsl{Cette propriété fut démontrée par l'astronome et ingénieur arabe Abul Abbas Al-Farghani (805-880), qui vivait entre le Caire et Baghadad au milieu du 9ème siècle. Cette représentation était employée dans la confection des astrolabes.}
\\\index{\textbf{Hipparque}}\index{\textbf{Al-Farghani A.A.}}\index{\textbf{Ptolémée C.}}

\textsl{C'était Thomas Harriot (1560-1621) mathématicien et astronome anglais qui avait montré que la représentation stéréographique était conforme et approuvée par un papier présenté par l'astronome Edmond Halley (1656-1742) à la Société Royale de Londres. }
\\ \index{\textbf{Harriot T.}}

\textsl{Le terme "projection stéréographique " fut donné par le mathématicien belge et d'origine espagnole  François d'Aiguillon (1567-1617) en 1613 dans le sixième chapitre concernant les projections de son livre d'optique "Opticorum liber extus de proiectionibus".}\index{\textbf{Aiguillon F.}}
\\

\textsl{Rappelons que l'histoire des représentations conformes était le point de départ de la géométrie différentielle moderne avec le papier de Carl Friedrich Gauss de 1827 sur la théorie générale des surfaces. Un autre apport considérable était venu du travail du mathématicien français Gaspard Monge (1746-1818) spécialement son livre sur l'application de l'analyse à la géométrie.} (H.A. Kastrup, 2008)\index{\textbf{Kastrup H.A.}}\index{\textbf{Gauss C.F.}}\index{\textbf{Monge G.}}\index{\textbf{Halley E.}}
\section{\textsc{Exercices et Problèmes}}
\bpb \lb{prob1}
 Soit $(\BbS^2)$ la sphère de rayon $R$, au point $P(\varphi, \lambda )$ on lui fait correspondre le point $p(X,Y)$ du plan $OXY$ par la représentation plane suivante définie par les formules :
	\[                         
	   p (X,Y)\,\left\{\begin{array}{ll}
	        X = 2R.tg(\ds \frac{ \pi}{4} -\frac{\varphi}{2}).sin\lambda \\ 
	                        Y = - 2R.tg(\ds \frac{ \pi}{4} -\frac{\varphi}{2}).cos\lambda 
	                        \end{array}\right.
\]
1.	Montrer que l'image d'un méridien ($\lambda$  = constante ) est une droite dont on donne l'équation.

2.	Montrer que l'image d'un parallèle ($\varphi$  = constante ) est un cercle dont on précise l'équation.

3.	En utilisant le lemme de Tissot, déterminer les directions principales.

4.	Soit $dS$ la longueur infinitésimale correspondante sur le plan, calculer $dS$.

5.  Sachant que sur la sphère $ds^2 = R^2d\varphi^2 +R^2cos^2\varphi .d\lambda ^2$, calculer le module linéaire $m$.

6.	En déduire le module linéaire $m_1$ le long du méridien.

7.	En déduire le module linéaire $m_2$ le long d'un parallèle.

8.	Comparer $m_1$  et $m_2$. Conclure sur la conformité ou la non conformité de la représentation plane.
\epb
\bpb
 Soit $(\Sigma)$ la sphère de rayon $R$, au point $P(\varphi, \lambda )$ on lui fait correspondre le point $p(X,Y)$ du plan $OXY$ par la représentation plane suivante définie par les formules :
	\[                         
	   p (X,Y)\,\left\{\begin{array}{ll}
	        X = R.\lambda  \\ 
	                        Y = R.Logtg\left(\ds \frac{ \pi}{4}+\frac{\varphi}{2}\right) 
	                        \end{array}\right.
\]
où   $Log$ désigne le logarithme népérien.

1.	Quelles sont les images des méridiens  ($\lambda$ = constante)  et des  parallèles ($\varphi$ = constante). 

2.	Soit $dS$ la longueur infinitésimale correspondante sur le plan, calculer $dS$ en fonction de $\varphi$ et  de $\lambda$ et  calculer  le module linéaire $m$.

3.	En déduire les modules linéaires $m_1$ le long du méridien et  $m_2$ le long du parallèle.

4.	Comparer $m_1$ et $m_2$  et conclure sur la conformité ou la non conformité de la représentation plane.

5.	On suppose que $P$ décrit sur la surface $(\Sigma)$ une courbe  $(\gamma)$ telle que $\varphi$ et $\lambda$ sont liées par la relation :  $tg\varphi = sin\lambda$. Pour $\varphi = 0\, gr,\, 2\, gr,\, 4\, gr,\, 6\, gr,\, 8\, gr$ et $10\, gr$,  dresser un tableau  donnant  les  valeurs  de  $\lambda$ correspondantes.

6.	Sachant que  $R= 1000\, m$, calculer les coordonnées $(X,Y)$ de la représentation  plane donnée ci-dessus  pour  les  valeurs de $\varphi$  et $\lambda$  de  la  question 5.

7.	Rapporter à l'échelle 1/100 sur le plan $OXY$, les positions $(X, Y)$ des points. Que pensez-vous de l'image de la courbe $(\gamma)$. 
\epb
\bpb
Sur une sphère de rayon unité, modèle de la terre, on désigne :

-	par  $p$  le pôle nord;

-	par $(C)$  un grand cercle qui coupe l'équateur au point  $i$  de longitude nulle;

-	par $q$ le pôle de ce grand cercle, de latitude $\varphi_0$ positive,

-	par $\omega$  et $h$ respectivement les points d'intersection de $(C)$  et du méridien de $q$  et du grand cercle issu de $q$, passant par le point  $a (\varphi,\lambda)$. 

On pose :    $   \omega h = x,\quad       ha = y$.

1.	$q$ est le pivot d'une représentation cylindrique conforme oblique tangente, dont $(C)$  est le ''pseudo-équateur''. Le plan est rapporté aux axes   $\Omega X,\Omega  Y$ images respectives de  $(C)$   et du grand cercle   $\omega pq$. Exprimer en fonction de $\varphi, \lambda$    et $\varphi_0$ les coordonnées $X,Y$ du point $A$ image de a (vérifier que pour  $\varphi_0 = 0$, on retrouve les expressions  de $X,Y$ d'une représentation transverse).

2.	Montrer que l'équation de l'image plane du parallèle de latitude $\varphi_0$  peut s'écrire :                $$e^Y cosX=tg\varphi_0$$                   

Indications : $b$ désignant un point de latitude $\varphi_0$, le triangle $pqb$ est isocèle,  décomposer ce triangle en deux triangles rectangles égaux. Etudier qualitativement les images des autres parallèles. 

3.	Montrer que l'image plan de l'équateur a pour équation :   
                                               $$ cosX + tg\varphi_0.shY = 0$$

Ecrire d'une manière analogue, l'équation de l'image du méridien $\lambda = 0$.

4.	Exprimer le gisement du méridien en fonction de $\varphi,\lambda$   et $\varphi_0$. Déterminer la valeur du module linéaire, en particulier en $p$, en un point de l'équateur, et en un point du méridien origine. 
\epb
\bpb
Etude de la représentation conforme d'une sphère de rayon unité dite représentation de Littrow\footnote{En hommage à  \textbf{Joseph Johann Littrow} (1781-1840) astronome autrichien.} définie par :
  $$ Z = sinz $$    
	avec  $z =\lambda   + iL_M$, $L_M$ la latitude de Mercator  et  $Z = X + iY$.

1. Préciser le canevas, les images des méridiens et celle de l'équateur.

2. Vérifier que les points $f$ et $f'$ $(\varphi = 0, \lambda   =\pm \pi/2)$ sont des points singuliers.

3. Etudier les images plans des cercles de diamètre $ff'$ et des petits cercles orthogonaux.

4. Soit $s$ le point $(\varphi  =\varphi_0, \lambda   =  0 )$. On appelle segment capable sphérique l'ensemble des points $b$ tels que l'angle $ \widehat{bp,bs} = \alpha $. Quelle est l'image plane de cette courbe dans cette représentation plane.
\epb
\bpb
 Soit l'application $ F(u,v):\BbR^2\longrightarrow \BbR^3 \verb|\|(0,0,1)$ définie par:
$$ OM(u,v)=F(u,v)\left\{\begin{array}{l}
x=\ds \frac{2u}{u^2+v^2+1} \\
\\
y=\ds \frac{2v}{u^2+v^2+1}\\
\\
z=\ds \frac{u^2+v^2-1}{u^2+v^2+1}
\end{array}\right. $$ 

1. Calculer la forme fondamentale $ds^2$.

2. Montrer que $OM(u,v)$ appartient à la sphère $(\BbS^2)$ d'équation $x^2+y^2+z^2=1$.

3. Calculer $u,v$ en fonction de $x,y$ et $z$.

4. Soit le point $N(0,0,1)$ de  $(\BbS^2)$, calculer les coordonnées $(X,Y)$ du point $p$ intersection de la droite $NM$ avec la plan $z=0$ en fonction de $x,y$ et $z$. 

5. Soit $\sigma$ l'application $\BbR^3 \verb|\|(0,0,1)\longrightarrow \BbR^2:(x,y,z)\longrightarrow (X,Y)=(X(x,y,z),Y(x,y,z))$. Montrer que $(\sigma \circ F)(u,v)=\sigma(F(u,v))=(u,v)$. En déduire que $F=\sigma^{-1}$.

6. Trouver le rapport de ce problème avec le problème \ref{prob1}.
\epb
\bpb
Soit un ellipsoïde de révolution $E(a,e)$ avec $a$ et $e$ respectivement le demi-grand axe de l'ellipsoïde de révolution et $e$ la première excentricité. Soit $(\BbS^2)$ une sphère de rayon $R$. On veut étudier le passage suivant: 
$$	p(\varphi,\lambda) \,\,\mbox{de l'ellipsoïde}\,\, E \Rightarrow \,P(\psi,\Lambda)\,\mbox{de la sphère }\, \,\BbS^2 $$
1. Exprimer $m$ le module linéaire de cette représentation.

2. On pose:$$ z=\m L+i\lambda,\quad Z=L+i\Lambda $$
$\m L$ est la latitude isométrique de l'ellipsoïde de révolution et $L$ la latitude de Mercator. Une transformation conforme entre $E$ et $(\BbS^2)$ est donnée par $Z=f(z)$ où $f$ est une fonction analytique. On propose le cas le plus simple à savoir:
\ba
Z=\alpha z+\beta \nonumber \\
avec \,\,\left\{\begin{array}{l}
\alpha=c_1+ic_2  \\
\beta=b_1+ib_2
\end{array} \right.\nonumber
\ea
les $c_1,c_2,b_1,b_2$ sont des constantes réelles. Donner les expressions de $L$ et $\Lambda$ en fonction de $\m L$ et $\lambda$.

3. On veut que repésentation transforme les méridiens et les parallèles de l'ellipsoïde respectivement en méridiens et parallèles de la sphère et que l'image du méridien origine $\lambda=0$ soit le méridien origine de la sphère $\Lambda=0$. Montrer que $c_2=b_2=0$ et $ L=c_1\m L+b_1, \quad 
\Lambda= c_1\lambda$.

4. Pour avoir la même orientation en longitude, on prendra $c_1>0$. On cherchera la transformation à déformation minimale autour d'un parallèle $\varphi=\varphi_0$ tel que le parallèle $\varphi=\varphi_0$ est automécoïque et le module linéaire $m$ est stationnaire pour $\varphi=\varphi_0$, c'est-à-dire $m(\varphi_0)=1$ et $\ds \left(\frac{dm}{d\varphi }\right)\biggr |_{\varphi=\varphi_0}=0$, en plus on considère aussi la condition: $$\ds \left(\frac{d^2m}{d\varphi ^2}\right)\biggr |_{\varphi=\varphi_0}=0 $$
Pour faciliter les notations, on prendra $b=b_1, c=c_1$. Montrer que la relation liant $\varphi_0$ et $\psi_0$ est:
$$  tg\psi_0=tg\varphi_0\sqrt{\frac{1-e^2}{1-e^2sin^2\varphi_0}} $$
5. Déterminer les constantes $b,c$ et $R$ en fonction de $\varphi_0$ et $\psi_0$ telles que les conditions ci-dessus soient vérifiées.

6. Montrer que l'expression du développement limité de $m(\varphi)$ de part et d'autre du parallèle $\varphi_0$ est donnée par:
$$m(\varphi)=1-\frac{2e^2(1-e^2)sin\varphi_0cos\varphi_0}{3(1-e^2sin^2\varphi_0)^2}(\varphi -\varphi_0)^3+o((\varphi -\varphi_0)^4)$$
7. On fait intervenir la deuxième excentricite $e'$, Montrer que $m(\varphi)$ s'écrit:
$$m(\varphi)=1-\frac{2e'^2sin\varphi_0cos\varphi_0}{3(1+e'^2cos^2\varphi_0)^2}(\varphi -\varphi_0)^3+o((\varphi -\varphi_0)^4)$$
\epb
\bpb
Soit $\m E(a,b)$ un ellipsoïde de référence de paramètres $a$ et $e$ respectivement le demi-grand axe et la première excentricité. On considère une représentation plane $\m P$ de $\m E$ vers le plan $(O,X,Y)$. On pose:
$$\begin{array}{l}
z=\lambda+i\m L\\
Z=X+iY=Z(z)
\end{array}$$
 avec $\m L$ la latitude isométrique.

1. Ecrire les expressions du carré des éléments infinitésimaux de longueur sur l'ellipsoïde et le plan. En déduire le module linéaire $m$.

2. On pose $\ds \zeta=\frac{\partial Z}{\partial z}$. Si $\gamma$ est le gisement de l'image du méridien passant par le point $z=(\lambda,\m L)$, montrer que $\ds  arg(\zeta)=\frac{\pi}{2}-\gamma$.

3. On cherche une représentation plane du type $Z=\alpha +\beta z+\varpi z^2$ où $\alpha,\beta$ et $\varpi$ des constantes complexes. On impose les conditions suivantes:

- pour $z=0$, $Z=0$;

- l'axe des $Y$ coïncide avec le méridien à l'origine.

Montrer que $\m Re (\beta)=0$.
 
4. En déduire que $Z$ s'écrit: $ Z=i\beta_1z+(\varpi_1+i\varpi_2)z^2$ avec $\beta_1,\varpi_1,\varpi_2$ des réels.
\epb
\bpb
L'objet de ce problème est l'étude de la comparaison de deux réseaux géodésique par l'effet d'une translation tridimensionnelle. On commence par l'étude d'un modèle sphérique, puis celui d'un modèle ellipsoïdique. Le point de départ est deux calculs $R_1$ et $R_2$ d'un réseau géodésique $\m R$ qu'on considère au voisinge de l'equateur. 

I. On considère que les coordonnées de $R_1$ et $R_2$ sont exprimées par la représentation plane Mercator directe d'un modèle sphérique (de rayon $a$). Un point $M$ a pour coordonnées:
$$ M (	X=acos\varphi cos\lambda, Y=acos\varphi sin\lambda,Z=asin\varphi)$$
1. Exprimer les coordonnées géodésiques $(\varphi, \lambda)$ en fonction des coordonnées tridimensionnelles $X,Y,Z$.

2. Une translation tridimensionnelle infinitésimale est exprimée par $(dX,dY,dZ)$ des formules ci-dessus. Montrer que $d\varphi$ et $d\lambda$ en fonction de $dX,dY$ et $dZ$ donnent: 
\ba
	d\lambda=\frac{-sin\lambda}{acos\varphi}dX+\frac{cos\lambda}{acos\varphi}dY \nonumber \\
	\frac{d\varphi}{cos\varphi}=\frac{dZ}{a}-\frac{tg\varphi cos\lambda}{a}dX-\frac{tg\varphi sin\lambda}{a}dY \nonumber 
\ea
Or l'équation ci-dessus n'est autre que la différentielle de la latitude de Mercator $L$, d'où:
$$ 	dL=\frac{d\varphi}{cos\varphi}=\frac{dZ}{a}-\frac{tg\varphi cos\lambda}{a}dX-\frac{tg\varphi sin\lambda}{a}dY$$
3. On se place au voisinage du point central $M_0(\varphi=0,\lambda=0)$ de l'équateur, on peut écrire au deuxième ordre de petitesse:
	$$ cos\varphi=1-\frac{\varphi^2}{2},\quad	sin\lambda=\lambda, \quad cos\lambda=1-\frac{\lambda^2}{2}$$
	 A partir de la formule de $L$ montrer qu'on peut écrire  $L\approx \varphi$.

4. Montrer qu'on peut écrire alors:
\ba
 dL=\frac{dZ}{a}-\frac{dX}{a}L-\frac{dY}{a}L \nonumber \\
d\lambda=\frac{-\lambda(1+\frac{L^2}{2})}{a}dX+\frac{(1-\frac{\lambda^2}{2})(1+\frac{L^2}{2})}{a}dY\nonumber 
\ea
5. En gardant les termes du 2ème ordre, simplifier l'expression de $d\lm$.

6. On pose: $z=\lm+iL$ et $\xi=d\lm+idL$. Montrer alors la formule de \textbf{Dufour-Fezzani} (\textit{H.M. Dufour}, 1979):\index{\textbf{Dufour H.M.}}\index{\textbf{Fezzani C.}}\index{Formule de Dufour-Fezzani}
$$\xi=\frac{dY+idZ}{a}-\frac{dX}{a}z-\frac{dY}{2a}z^2  $$
7. La transformation $\xi$ est-elle conforme.

II. On considère que les coordonnées de $R_1$ et $R_2$ sont exprimées par la représentation UTM d'un modèle ellipsoïdique $E(a,e)$, $a$ et $e$ sont respectivement le demi-grand axe et la première excentricité. Un point $M$ du modèle ellipsoïdique a pour coordonnées:
$$ M\,\,(X=Ncos\varphi cos\lambda, Y=Ncos\varphi sin\lambda, Z=N(1-e^2).sin\varphi)$$
 avec : $ N=a/\sqrt{1-e^2sin^2\varphi}$.
 Les coordonnées géodésiques $(\varphi, \lambda)$ sont exprimées par:
 $$ tg\lambda=\frac{Y}{X}, \quad (1-e^2)	tg\varphi=\frac{Z}{\sqrt{X^2+Y^2}} $$
8. On considère une translation tridimensionnelle infinitésimale $(dX,dY,dZ)$. On note par $\m L$ la latitude isométrique:
$$ \m L=Logtg\left(\frac{\pi}{4}+\frac{\varphi}{2}\right)-\frac{e}{2}Log\frac{1+esin\varphi}{1-esin\varphi}=L-\frac{e}{2}Log\frac{1+esin\varphi}{1-esin\varphi}$$
On opère de la même manière qu'en I., montrer qu'au deuxième ordre de petitesse, on a:
\ba
d\lambda=\frac{-dX}{a}\lambda+\frac{dY}{a}+(1+e^2)\m L^2\frac{dY}{2a} -\frac{\lambda^2}{2a}dY \nonumber \\
d\m L= \frac{dZ}{a}-\frac{\m L^2dZ}{2a} -\frac{3e^2\m L^2dZ}{2a}-\m L \frac{dX}{a}-\lambda \m L \frac{dY}{a} \nonumber
\ea
9. On pose: $z=\lambda+i\m L,  \quad \zeta =d\lambda+id\m L $. Montrer alors, la formule de \textbf{Dufour-Ben Hadj Salem} (\textit{A. Ben Hadj Salem}, 2012):\index{\textbf{Dufour H.M.}}\index{Formule de Dufour-Ben Hadj Salem}
$$ \ds  \zeta=\frac{dY+idZ}{a}-\frac{dX}{a}z-\frac{dY}{2a}z^2-\frac{e^2(dY-idZ)}{8a}(z-\bar{z})^2 $$ 
10. La transformation $\zeta$ est-elle conforme. 
\epb
\bpb
Pour une représentation plane, on dit qu'elle est \textbf{équivalente} si le produit des modules linéaires $m_1$ et $m_2$ suivant les directions principales vérifient: $$m_1.m_2=1$$ 
Soit le modèle terrestre représenté par la sphère de rayon $R$ qu'on note $\BbS^2$. Au point $P(\varphi, \lambda )$ on lui fait correspondre le point $p(X,Y)$ du plan $OXY$ par la représentation plane suivante définie par les formules :
	\be                        
	   p (X,Y)= \left( X = 2R.sin\left(\ds \frac{ \pi}{4} -\frac{\varphi}{2}\right).cos\lambda,\quad Y = 2R.sin\left(\ds \frac{ \pi}{4} -\frac{\varphi}{2}\right).sin\lambda \right)\lb{pb1z}
	  \ee
1.  Qu'elle est l'image du pôle nord $P_N$?

2.	Montrer que l'image d'un méridien ($\lambda$  = $\lambda_0$=constante ) est une droite qu'on donne l'équation.

3.	Montrer que l'image d'un parallèle ($\varphi$  = $\varphi_0$=constante ) est un cercle qu' on précise l'équation.

4.	En utilisant le lemme de Tissot, déterminer les directions principales.

5.	Soit $ds$ la longueur infinitésimale correspondante sur la sphère, donner l'expression de $ds^2$.

6.  Soit $dS$ la longueur infinitésimale correspondante sur le plan. Montrer que: $$dS^2 = R^2cos^2\left(\ds \frac{\pi}{4}-\frac{\varphi}{2}\right)d\varphi^2 +4R^2sin^2\left(\ds \frac{\pi}{4}-\frac{\varphi}{2}\right).d\lambda ^2$$
7.  En déduire le carré du module linéaire $m^2$.

8.	Calculer le module linéaire $m_1$ le long du parallèle.

9.	Calculer le module linéaire $m_2$ le long du méridien.

10.	La représentation plane définie par (\ref{pb1z}) est-elle équivalente. Justifier votre réponse.
\epb \index{Représentation \'equivalente}
\chapter{\textit{\textbf{La Représentation Plane Lambert}}}
\section{\textsc{Définition et Propriétés}}
                La représentation plane Lambert \index{Représentation Lambert}est une représentation conique\index{Représentation conique}, conforme\index{Représentation conforme} et directe\index{Représentation directe} d'un modèle ellipsoïdique :
                
- conique : on utilise les coordonnées polaires  $R$ et $\Omega$;

-	conforme : conservation des angles ou l'altération angulaire est nulle;

-	directe : les coordonnées polaires sont des fonctions de la forme :
$$ \fbox{ $ \begin{array} {l} 
	R=R(\varphi)\\
	\Omega=\Omega(\lambda)
\end{array} $}
$$
où $(\varphi , \lambda )$ sont les coordonnées d'un point sur le modèle ellipsoïdique de référence.
\\

Pour la Tunisie, on a considère le cas de la représentation tangente\index{Représentation tangente} c'est-à-dire on utilise un seul parallèle origine. Dans la suite du cours, on va étudier en détail le cas d'un seul parallèle origine.

\begin{figure}
	\centering
		\includegraphics[width=0.60\textwidth]{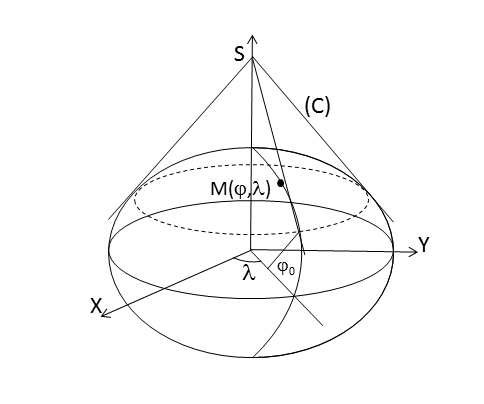}
	\caption{Interprétation géométrique}
	\label{fig:lambconique}
\end{figure}
Une interprétation de la représentation plane Lambert est comme suit:

- on considère un cône $(C)$ (\textbf{Fig. \ref{fig:lambconique}}) de sommet S tangent au parallèle origine de latitude $\varphi_0$ de l'ellipsoïde de référence $\m E$. A un point M$(\varphi,\lambda)$ de $\m E$, on lui fait correspondre son image m sur la demi-droite d'origine S tangente à la méridienne de longitude $\lambda$ et au parallèle origine.

\begin{figure}
	\centering
		\includegraphics[width=0.60\textwidth]{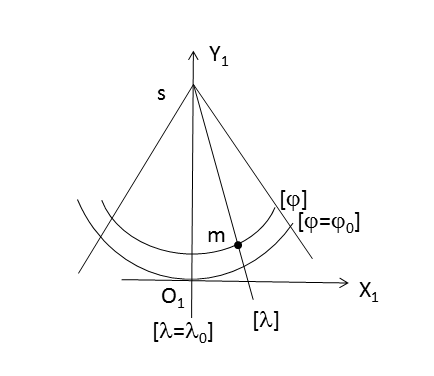}
	\caption{Images des parallèles et des méridiens }
	\label{fig:lambconique1}
\end{figure}

- on développe le cône $(C)$ sur le plan, on obtient l'image d'une portion de l'ellipsoïde (\textbf{Fig. \ref{fig:lambconique1}}).
\\

   Les images des parallèles sont des arcs de cercles concentriques de centre $\textsc{s}$ l'image du sommet du cône $(C)$, celles des méridiens sont des droites concordantes passant par $\textsc{s}$ (\textbf{Fig. \ref{fig:lambconique1}}).
\\

Les courbes coordonnées\index{Courbes coordonnées} $\varphi$  = constante et $ \lambda  =$ constante sur le modèle sont orthogonales et leurs images le sont aussi dans le plan.

\section{\textsc{Indicatrice de Tissot}}
D'après la propriété précédente des courbes coordonnées, on déduit que les directions principales \index{Directions principales} sont les tangentes au méridien et au parallèle passant par le point.
           \\
           
La représentation est conforme, par suite l'altération angulaire est nulle, l'indicatrice de Tissot\index{Indicatrice de Tissot} est un cercle et le module linéaire ne dépend pas de la direction mais seulement du point et on a l'équivalence:
$$ \fbox{$ \mbox{ Altération angulaire nulle}\Leftrightarrow      m_\varphi  = m_\lambda \Leftrightarrow \, \forall \delta  \,, m_\delta  = m $} $$
où $\delta$  désigne 'la direction'.
\section{\textsc{Calcul des modules principaux}}\index{Modules principaux}
          On commence par le calcul du module $m_\varphi $. Par définition :
$$ 	m_\varphi=\frac{dS}{ds} $$
avec $dS$ pris sur l'image de la méridienne et $ds$ sur la méridienne du modèle, or $ds = \rho d\varphi$ 
et      $dS= -dR $, le signe - provient du fait que les déplacements infinitésimaux $dR$ et $d\varphi$ sont de signe contraire. On note par $\rho$  le rayon de  courbure\index{Rayon de courbure de la méridienne}, d'où :       
$$ \fbox{ $	m_\varphi=\ds \frac{dS}{ds}=\frac{-dR}{\rho d\varphi} $} $$
Maintenant on calcule le module principal $m_\lambda$, on a :
$$ 	\fbox{ $ m_\lambda =\ds \frac{dS}{ds}=\frac{Rd\Omega}{rd\lambda} $} $$
avec $r = N.cos\varphi$     le rayon du parallèle de latitude $\varphi$.
\section{\textsc{Etablissement des Formules} $R( \varphi)$ \textsc{et} $\Omega ( \lambda)$}
Comme on a :       
	\[m_\lambda=m_\varphi
\]
d'où :
$$ \ds 	\frac{Rd\Omega}{rd\lambda}=\frac{-dR}{\rho d\varphi}\Longrightarrow \frac{d\Omega}{d\lambda}=\frac{-rdR}{\rho Rd\varphi} $$
Le terme à gauche est une fonction de $\lambda$   seulement car $\Omega$  ne dépend que de $\lambda$, le terme à droite est fonction de $\varphi$  seulement, donc l'égalité est toujours vérifiée que si les deux termes sont constants, on appelle $n$ cette constante, d'où : 
\ba
	\frac{d\Omega}{d\lambda}=n \nonumber \\
	\frac{-rdR}{\rho Rd\varphi}=n \nonumber
\ea
Par suite, en intégrant la première équation et prenant $\Omega(\lambda_0)=0$ avec $\lambda_0$ la longitude du méridien origine, on obtient:         
\be
	\fbox{ $ \Omega=n(\lambda -\lambda_0) $}
\ee
La deuxième équation différentielle s'écrit sous la forme:
\be
 \ds 	\frac{dR}{R}=\frac{-n\rho d\varphi}{r}=\frac{-n\rho d\varphi}{Ncos\varphi}=-ndL \lb{eqdif}
\ee
avec :
\be
dL=\frac{\rho d\varphi}{Ncos\varphi}
\ee
La variable $L$ appelée la latitude isométrique donnée par la formule:
\be
\fbox{ $ L(\varphi)=\ds Logtg\left(\frac{\pi}{4}+\frac{\varphi}{2}\right)-\frac{e}{2}Log\left(\frac{1+esin\varphi}{1-esin\varphi}\right) $} \lb{lisom}
\ee
En effet, calculant une primitive de $\ds \frac{\rho d\varphi}{Ncos\varphi}$ en considérant $\varphi$ n'atteint pas les valeurs extrêmes $\ds \pm \frac{\pi}{2}$, soit:
\ba
& \ds\int \frac{\rho d\varphi}{Ncos\varphi}=\int \frac{(1-e^2)d\varphi}{cos\varphi (1-e^2sin^2\varphi)}=\int \frac{(1-e^2sin^2\varphi-e^2cos^2\varphi)d\varphi}{cos\varphi (1-e^2sin^2\varphi)}=\nonumber &\\ &
\ds \int \frac{d\varphi}{cos\varphi }-\int \frac{e^2cos\varphi d\varphi}{1-e^2sin^2\varphi}=\ds Log tg\left(\frac{\pi}{4}+\frac{\varphi}{2}\right) -e\int\frac{du}{1-u^2}=\nonumber & \\
& \ds Log tg\left(\frac{\pi}{4}+\frac{\varphi}{2}\right)-\frac{e}{2}Log\frac{1+u}{1-u}=\ds Log tg\left(\frac{\pi}{4}+\frac{\varphi}{2}\right)-\frac{e}{2}Log\frac{1+esin\varphi}{1-esin\varphi}+C \nonumber &
\ea
où on a fait le changement de variables suivant $u=esin\varphi$ avec $|u|<1$. Par suite en prenant l'intégrale entre la latitude origine $\varphi=0$ et la latitude $\varphi$ et que $L(0)=0$, on retrouve la formule donnée par (\ref{lisom}) \index{Latitude isométrique} où $e$ est la première excentricité de l'ellipsoïde de référence.

On revient à l'équation (\ref{eqdif}):
 $$ \ds 	\frac{dR}{R}=\frac{-n\rho d\varphi}{r}=\frac{-n\rho d\varphi}{Ncos\varphi}=-ndL $$
En posant:
	\[L_0=L(\varphi_0)
\]
   où $\varphi_0$  est la latitude du parallèle origine, l'intégration de (\ref{eqdif}) donne:
\be
	Log\frac{R}{R_0}=-n(L-L_{0})\Longrightarrow \fbox{$ R=R_0exp{(-n(L-L_0))}=\ds R_0e^{-n(L-L_0)} $}
\ee
\section{\textsc{Détermination des Constantes} $R_0$ \textsc{et} $n$}
Pour déterminer les constantes $R_0$  et $n$, on impose que le parallèle origine soit un isomètre automécoïque et stationnaire\index{Isomètre stationnaire},  c'est-à-dire :\index{Isomètre automécoïque}
\ba
	                   m(\varphi_0) = 1 \nonumber \\               
	\mbox{et} \, \left(\frac{dm}{d\varphi}\right)_{\varphi=\varphi_0}=0 \nonumber 
\ea
soit le module linéaire admet un minimum égal à $m(\varphi_0)$. Comme: 
$$ 	m = m_\varphi  = m_\lambda  = \ds \frac{-dR}{\rho d\varphi} $$
et: 
\be
	dR=-nRdL=\ds -nR\frac{\rho d\varphi}{Ncos\varphi} \label{es18}
	\ee
D'où l'expression du module linéaire\index{Module linéaire}:                                  
\be
	m =\frac{nR}{Ncos\varphi}=\frac{nR_0e^{-n(L-L_0)}}{Ncos\varphi} \label{es19}
\ee
Pour $ \varphi = \varphi_0$, on a:
\be
	 m(\varphi_0) = 1 = \frac{nR_0}{N_{0}cos\varphi_0} \,\Rightarrow \,\fbox{ $nR_0 = N_{0}cos\varphi_0 $}   \label{es20} 
	\ee
on appelle :
$$	r=Ncos\varphi $$
Le calcul de $\ds \frac{dr}{d\varphi}$ donne :
\be
	\frac{dr}{d\varphi}= N'_{\varphi}cos\varphi - Nsin\varphi
\ee
comme:
	\[N = a(1 - e^2sin^2\varphi )^{-1/2}
\]
on obtient :
$$ 	N^{'}_{\varphi}=\ds \frac{e^2Nsin\varphi cos\varphi}{1-e^2sin^2\varphi} $$
On a alors:
	\[	\frac{dr}{d\varphi}=\frac{e^2Nsin\varphi cos^2\varphi}{1-e^2sin^2\varphi}-Nsin\varphi=Nsin\varphi \left(\frac{e^2cos^2\varphi}{1-e^2sin^2\varphi}-1\right)
\]
	\[\frac{dr}{d\varphi}=\frac{-(1-e^2)Nsin\varphi}{1-e^2sin^2\varphi}
\]
Or:
$$ \rho=\ds \frac{(1-e^2)N}{1-e^2sin^2\varphi} $$
Par suite:
\be
\fbox{ $ \ds 	\frac{dr}{d\varphi}=-\rho sin\varphi $}
\ee
On revient à l'équation (\ref{es19}) : $ m  = nR/r$ et en prenant sa différentielle logarithmique, d'où le résultat:
$$ \ds	\frac{dm}{m}=\frac{dR}{R}-\frac{dr}{r} $$
En utilisant les équations (\ref{es18}) et (\ref{es19}), on obtient :
$$ \ds	\frac{dm}{m}=\frac{-n\rho d\varphi}{r}+\frac{\rho sin\varphi d\varphi}{r}=(n-sin\varphi)\frac{\rho d\varphi}{r} $$
Soit:
\be
	\frac{dm}{d\varphi}=(n-sin\varphi)\frac{m\rho}{r} \label{es27}
\ee
Et pour $\varphi  = \varphi_0$, on a:
	\[\frac{dm}{d\varphi}_{\varphi_0}=0\Rightarrow(n-sin\varphi_0)\frac{m(\varphi_0)\rho_0}{r(\varphi_0)}=0
\]
D'où:
\be
\fbox{ $ 	n=sin\varphi_0  $}  \label{es28}
\ee
L'équation (\ref{es20}) s'écrit donc: 
\be
\fbox{ $ 	R_0=N(\varphi_0)cotg\varphi_0=N_0 cotg\varphi_0 $}
\ee
d'où les équations de la représentation plane Lambert :
\be
	\fbox{ $\begin{array}{l}
	\Omega=(\lambda-\lambda_0)sin\varphi_0   \\
  \\
  R=N_0cotg\varphi_0 e^{-sin\varphi_0(L-L_0)}   
	  	\end{array} $}
\ee
	avec:
	$$ L(\varphi)=Log tg\left(\ds \frac{\pi}{4}+\frac{\varphi}{2}\right)-\ds \frac{e}{2}Log\left(\frac{1+esin\varphi}{1-esin\varphi}\right) $$
L'expression du module linéaire est égale à:
\be
	\fbox{ $ m(\varphi)=\ds \frac{sin\varphi_0 R(\varphi)}{N(\varphi)cos\varphi} $ }
\ee
\section{\textsc{Expression des Coordonnées Cartésiennes}}
Dans ce paragraphe, on va décrire les coordonnées cartésiennes en fonction de $(\Omega,R)$. Soit un point  $M(\varphi,\lambda)$ ayant pour coordonnées polaires $(\Omega,R)$.
\\

\begin{figure}[h]
	\centering
		\includegraphics[width=0.60\textwidth]{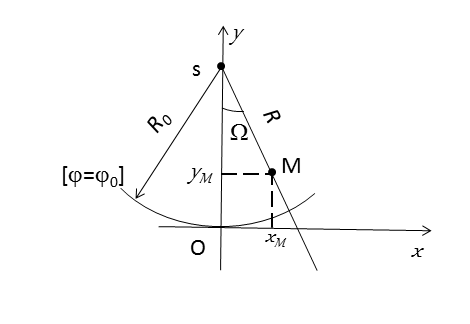}
	\caption{Le repère origine}
	\label{repereorigine}
\end{figure}

On considère un système d'axes $(O,x,y)$ qu'on nomme repère origine, tel que  l'axe $Ox$ est la tangente à l'image du parallèle origine au point $O$ dirigé vers l'Est et $Oy$ est porté par l'image du méridien origine dirigé vers le Nord (\textbf{Fig. \ref{repereorigine}}).
Soit le point $S$ de $Oy$ avec $OS=R_0$, on a alors:
\ba
		x_M=Rsin\Omega \nonumber \\
			y_M=R_0-Rcos\Omega \nonumber 
\ea
ou encore :                    
\be
\fbox{ $ \begin{array}{l}
	x_M=Rsin((\lambda-\lambda_0)sin\varphi_0) \\
	y_M=R_0 -Rcos((\lambda-\lambda_0)sin\varphi_0) 
	\end{array} $} \lb{1022}
\ee
avec $\lambda$ comptée positivement à l'Est du méridien origine des longitudes.
\section{\textsc{Passage des Coordonnées} $( R,\Omega )$ \textsc{aux Coordonnées} $(x,y)$}
   Ayant $(\varphi ,\lambda )$ et  $\varphi_0 ,\lambda_0$ , on calcule :
\ba
		\Omega=  ( \lambda - \lambda_0 )sin\varphi_0 \nonumber \\
L (\varphi ) = Logtg\left(\frac{\pi}{4}+ \frac{\varphi}{2}\right)- \frac{e}{2}Log\frac{1+esin\varphi}{1-esin\varphi}\nonumber \\
R_0 =N_0cotg\varphi_0 \nonumber \\
R = R_0exp(-sin\varphi_0(L-L_0))\nonumber \\
x =  Rsin\Omega \nonumber \\
y = R_0 - Rcos\Omega \nonumber 
\ea
\section{\textsc{Passage des Coordonnées} $(x,y)$ \textsc{aux Coordonnées} $(R,\Omega)$}
     On donne $\varphi_0$ et $\lambda_0$  et ayant  $(x,y)$, on calcule :
\ba
		R_0= N_0cotg\varphi_0 \nonumber\\
	Rcos\Omega=R_0-y \nonumber
\ea
Comme : 
	\[x = Rsin\Omega
\]
d'où : 
	\[tg\Omega   =  \frac{x}{R_0-y}
\]
Par suite :    
$$	\Omega = (\lambda  -  \lambda_0 )sin\varphi_0 = Arctg\left(\ds \frac{x}{R_0-y}\right) $$
D'où :
\be
	\fbox{ $ \lambda=  \lambda_0  + \ds \frac{1}{sin\varphi_0}Arctg\left(\ds \frac{x}{R_0-y}\right) $}
\ee
De: 
	\[y = R_0  - Rcos\Omega
\]
on obtient :
	\[R =   \frac{R_0 - y}{cos\Omega} 
\]
Et de:
	\[R = R_0 exp(-sin\varphi_0 (L - L_0)) \Rightarrow   Log\frac{R}{R_0}  = - sin\varphi_0 (L - L_0)
\]
d'où:
\be
\fbox{ $	L = L_0 + \ds \frac{1}{sin\varphi_0}Log\displaystyle \frac{R_0}{R} $ }
\ee
Le problème devient à calculer $\varphi$  à partir de la donnée de la latitude isométrique $L$. Ce calcul se fait par itérations comme suit:
\begin{enumerate}
	\item Ayant $L$, on calcule $\varphi_1$ telle que $L=Logtg\ds \left(\frac{\pi}{4}+\frac{\varphi_1}{2}\right)$.
\item On calcule $\varphi_2$ telle que $L+\ds \frac{e}{2}Log\left(\frac{1+esin\varphi_1}{1-esin\varphi_1}\right)=\ds Logtg\left(\frac{\pi}{4}+\frac{\varphi_2}{2}\right)$.
\item On réitère le processus jusqu'à ce que $|\varphi_{i+1}-\varphi_i|<\alpha$ où $\alpha$ une petite quantité fixée à l'avance. 
\end{enumerate}

\section{\textsc{Etude de l'Altération Linéaire}}
L'altération linéaire\index{Altération linéaire} est définie par:
\be
	\fbox{ $ \epsilon=m-1 $}
\ee
où $m$ est le module linéaire. Le développement limité du module linéaire au voisinage de $\varphi_0$ s'écrit:
$$ m(\varphi)=m(\varphi_0)+(\varphi-\varphi_0)\left(\ds \frac{dm}{d\varphi}\right)_{\varphi=\varphi_0}+\frac{(\varphi-\varphi_0)^2}{2}\left(\frac{d^2m}{d\varphi^2}\right)_{\varphi=\varphi_0}+o((\varphi-\varphi_0)^3) $$
Or :
$$ 	m(\varphi_0)=1 \, \, \mbox{et} \,\,	\left(\ds \frac{dm}{d\varphi}\right)_{\varphi=\varphi_0}=0 $$
car le parallèle $\varphi=\varphi_0$ est un isomètre automécoïque et stationnaire, d'où:
\be
m(\varphi)=1+\frac{(\varphi-\varphi_0)^2}{2}\left(\frac{d^2m}{d\varphi^2}\right)_{\varphi=\varphi_0}+o((\varphi-\varphi_0)^3) \label{es53}
\ee
On est amené à calculer la valeur de la dérivée seconde de $m$ pour $\varphi=\varphi_0$. Or l'équation (\ref{es27}) donne l'expression de $m^{'}_{\varphi}$. On dérive $m^{'}_{\varphi}$.
\be
\frac{d^2m}{d \varphi^2}=\frac{dm^{'}_{\varphi}}{d\varphi}=\frac{d}{d\varphi}\left[\frac{m\rho}{r}(sin\varphi-n)\right]=\frac{m\rho}{r}cos\varphi +(sin\varphi-n)\frac{d}{d\varphi}\left(\frac{m\rho}{r}\right)
\ee
d'où:
\be
	\left(\frac{d^2m}{d\varphi^2}\right)_{\varphi=\varphi_0}=\frac{cos\varphi_0 m(\varphi_0) \rho(\varphi_0)}{r(\varphi_0)}=\frac{\rho_0}{N_0}
	\ee
	car $n=sin\varphi_0$ et $m(\varphi_0)=1$. (\ref{es53}) devient:
$$ m(\varphi)=1+\frac{(\varphi-\varphi_0)^2}{2}\frac{\rho_0}{N_0}+o((\varphi-\varphi_0)^3)=1+\frac{(\varphi-\varphi_0)^2}{2}\frac{\rho^2_0}{\rho_0 N_0}+o((\varphi-\varphi_0)^3) $$
En posant :
	\[\Delta\varphi=\varphi-\varphi_0
\]
on a alors:
\be
	\fbox{ $	m(\varphi)=1+\ds \frac{1}{2N_0 \rho_0}(\rho_0 \Delta \varphi)^2
	+ o(\Delta \varphi^3) $}
\ee
Or $\rho_0 \Delta\varphi=\rho_0 (\varphi-\varphi_0)$ la distance approchée du point $M(\varphi,\lambda)$ au parallèle origine $ \varphi=\varphi_0$, d'où l'expression de l'altération linéaire:
\be
	\fbox{ $ \epsilon=m-1=\ds \frac{1}{2N_0 \rho_0}(\rho_0 \Delta\varphi)^2=\frac{\ell^2}{2N_0 \rho_0} $}
\ee
où $\ell$ est la distance du point au parallèle origine.
\subsection*{10.9.1. Calculs numériques}
On considère comme exemple numérique le cas de la représentation Lambert Nord Tunisie ayant comme parallèle origine $\varphi_0=40.0\,gr$ et l'ellipsoïde de référence est celui de Clarke Français 1880. 

On a donc les valeurs numériques du module et l'altération linéaires comme suit:
\ba
	m(\varphi_0)=1\Longrightarrow \epsilon=0 \nonumber \\
	m(42.5\, gr)=1.000\,775\,720 \Longrightarrow \epsilon=7.75720\times 10^{-4} \nonumber \\
	m(37.5\, gr)=1.000\,760\,827 \Longrightarrow \epsilon=7.60827\times 10^{-4} \nonumber
\ea
Soit une distance de $1000\, m$ sur le parallèle origine, elle se transforme à $1000\, m$ sans altération. Une distance de $1000\, m$ sur le parallèle $\varphi=42.5\, gr$ devient une distance de $ 1000.776\, m$ sur le plan, de même une distance de $1000\, m$ sur le parallèle $ \varphi=37.5\, gr$ devient une distance de $1000.761\, m$ sur le plan.
\\

Pour réduire les altérations linéaires, on multiplie le module linéaire par un coefficient $k$ dit facteur de réduction de l'échelle\index{Facteur de réduction de l'échelle}. Le module linéaire devient alors:
\be
\fbox{ $ 	m'=k.m=\ds \frac{ksin\varphi_0 R(\varphi)}{N(\varphi)cos\varphi} $}
\ee
Par suite, les modules linéaires et les altérations correspondantes deviennent (Cas de la Tunisie, le facteur $k=k_N$ est égal à 	$0.999\, 625\, 544$):
\ba
	m'(\varphi_0)=0.999\,625\,544\Rightarrow \,\epsilon=-0.000\,009\,460 \nonumber \\
	m'(42.5\, gr)=1.000\,400\,974 \Rightarrow \,\epsilon=+0.400\,974\times10^{-3} \nonumber \\
	m'(37.5\, gr)=1.000\,386\,086 \Rightarrow \,\epsilon=+0.386\,086\times10^{-3} \nonumber	
\ea
Sur le parallèle $\varphi=42.5 \,gr$, l'altération linéaire pour $1000\,m$ a passé de $+0.776\,m$ à $+0.401\,m$, d'où réduction des altérations.
\\

Avec l'introduction du facteur de réduction de l'échelle, les formules (\ref{1022}) des coordonnées rectangulaires $(x,y)$ s'écrivent:
\be 
\fbox{ $ \begin{array}{l}
	x_M=kRsin((\lambda-\lambda_0)sin\varphi_0) \\
y_M=k(R_0 -Rcos( (\lambda-\lambda_0)sin\varphi_0)) 
				\end{array} $}
\ee
Pour obtenir des coordonnées rectangulaires positives, on définit un repère $(O',X,Y)$ tels que  $O'X$ et $O'Y$ soient dirigés respectivement vers l'Est et le Nord (\textbf{Fig. \ref{reperetranslate}}) et que:
\begin{figure}[htp]
	\centering
		\includegraphics[width=0.60\textwidth]{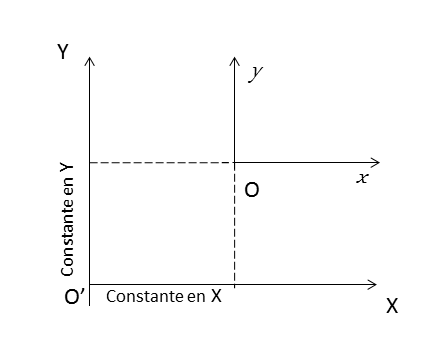}
	\caption{Le repère $(O',X,Y)$}
	\label{reperetranslate}
\end{figure}
\be
\fbox{ $ \begin{array}{l}
X=Constante\,\,X\,+x_M\\
Y=Constante\,Y\,+y_M
		\end{array} $}
\ee
Les quantités \textit{Constante $X$} et \textit{Constante $Y$} sont respectivement les constantes de translation en $X(Est)$ et en $Y(Nord)$ exprimées en mètres.
\section{\textsc{Convergence des méridiens}}\index{Convergence des méridiens}
Pour passer de l'azimut géodésique sur le modèle ellipsoïdique de la direction $m_1m_2$ au gisement $M_1M_2$ sur le plan de la représentation, on a la formule algébrique:
\be
	\fbox{ $ G= Az - \gamma + Dv $}
\ee
avec $ \gamma$ le gisement de l'image du méridien avec son signe positif ou négatif. Or l'image d'un méridien est une droite qui coupe l'axe Ox (du nord) sous l'angle $ \Omega=(\lambda-\lambda_0)sin\varphi_0$, par suite:
\be
	\fbox{ $ \gamma=\Omega=(\lambda-\lambda_0)sin\varphi_0 = \, \mbox{convergence des méridiens} $}
\ee
\section{\textsc{Calcul de la réduction de la corde}}\index{Réduction à la corde}
Sur le terrain, on observe une direction $AB$, la visée $AB$ est très voisine de la géodésique $AB$. Sa transformée sur le plan de la représentation n'est pas une droite, mais une courbe tournant sa concavité vers l'image du parallèle origine.
\\

On observe la direction $AB$ c'est-à-dire l'arc $ab$. Pour passer de l'arc $ab$ à la cordre $\overline{ab}$, on apporte une correction à la lecture de la direction $AB$ ou $\overline{ab}$. Cette correction est appelée la correction de réduction à la corde. Elle est donnée par la formule:
\be
	\fbox{ $ Dv= \ds \frac{S}{2}\Gamma\left(\frac{S}{3}\right) $}
\ee
où $S$ représente la longueur $AB$, $\Gamma\left(\frac{S}{3}\right)$ est la courbure de la transformée de la géodésique $AB$ prise au 1/3 de la distance de $A$ vers $B$.
\\

En utilisant la formule de Schols-Laborde\footnote{\textbf{Jean Laborde}: colonnel de l'armée française et géodésien cartographe. Il a défini la représentation plane qui porte son nom (représentation conforme cylindrique oblique). Celle-ci a été appliquée pour le Madagascar.} \index{Formule de Schols-Laborde} donnant la courbure de la transformée d'une géodésique, on démontre que:\index{\textbf{Laborde J.}}
\be
	\fbox{ $ Dv^{(dmgr)}=K.d\lambda $}
\ee
avec $d\lambda =$  la différence de longitude en $km$ des 2 extrémités de la visée $AB$ et $K$ vaut:
\be
	\fbox{ $ K=\ds \frac{1}{2}\frac{(R_0-R)_{1/3}}{N(\varphi_0)\rho(\varphi_0)sin1''} $}
\ee
où $(R_0 - R)$, $N(\varphi_0)$ et $\rho(\varphi_0)$ en $km$.
\section{\textsc{Exercices et Problèmes}}
\bex
 En un point $A$ de coordonnées géodésiques $\varphi = 40.9193\, gr$ et $\lambda = 11.9656\, gr$ à l'Est de Greenwich, on vise un point $B$.
 
1.	Calculer les coordonnées planes Lambert du point A, sachant que $\varphi_0=40.00\, gr$ et $\lambda_0=+11.00\,gr$.
 
2.	L'azimut géodésique de la direction $AB$ est $Azg = 55.7631\, gr$. Sachant que $Dv = 1.52\,dmgr$, calculer $G$  le gisement de la direction $AB$.

3.	La distance $AB$ réduite à l'ellipsoïde de référence est $D_e = 5421.32\, m$. Sachant que l'altération linéaire dans la région des points $A$ et $B$ vaut $- 9\, cm/km$, calculer la distance $AB$ réduite au plan.
\eex
\bex
 D'après les coordonnées de deux points $A$ et $B$ vous trouvez la distance $AB= 5427.380\, m$. Sachant que :

a -	l'altération linéaire de la représentation dans la région de $AB$ vaut $+8.10^{-5}$;

b -	les altitudes des points $A$ et $B$ sont : $H_A = 1000.00\,m$ et $H_B = 1200.00\, m$. Calculer la distance suivant la pente $D_P$ entre les points $A$ et $B$ matérialisés sur le terrain.
\eex
\bpb
On désire étudier les variations du module linéaire $m$ de la représentation Lambert tangente en fonction de la latitude $\varphi \in [0,+\ds\frac{\pi}{2}]$.

1. Donner l'expression de $m(\varphi)$.

2. Montrer que $m(0)$ est une quantité finie positive que l'on calculera.
 
3. Montrer que pour $\varphi=+\ds\frac{\pi}{2}$, $m$ devient une forme indéterminée qu'on précisera.

4. On pose $\varphi=\ds\frac{\pi}{2}-\theta$ et $m(\varphi)=M(\theta)$. Donner l'expression de $M(\theta)$.

5. Montrer qu'on peut écrire $M(\theta)=A(\theta).u$ où $A$ est une fonction de $\theta$ prenant une valeur finie positive non nulle quand $\theta\longrightarrow 0^+$ et $u$ donnée par:$$ u=\ds \frac{exp\left(sin\varphi_0Logtg\ds\frac{\theta}{2}\right)}{sin\theta}$$
avec $\varphi_0$ la latitude du parallèle origine.  

6. On pose $t=tg\ds\frac{\theta}{2}$. Montrer que $u$ s'écrit:$$ u=\ds \frac{1+t^2}{2t^{1-sin\varphi_0} }$$

et que $lim_{t\longrightarrow 0^+}u=+\infty$. 

7. En déduire la limite de $m(\varphi)$ quand $\varphi\longrightarrow \ds \frac{\pi}{2}^-$. 

8. Donner l'expression de $\ds \frac{dm}{d\varphi}$ et présenter le tableau de variation de $m(\varphi)$ pour le cas où $\varphi_0=+40\,gr$.
\epb 
\bpb
 On a mesuré une distance suivant la pente $D_P = 20\,130.858\, m$ entre deux points $A$ et $B$ avec  $H_A = 235.07\, m,\, H_B = 507.75\, m$, on prendra comme rayon terrestre  $R = 6378\, km$.

1.	Calculer la distance $D_e$ suivant l'ellipsoïde en utilisant la formule rigoureuse.

2.	Sachant que le module linéaire $m$ vaut $0.999\,850\,371$, calculer la distance $D_r$ réduite au plan de la représentation plane utilisée.

3.	Les coordonnées géodésiques du point $A$ sont : $\varphi = 10.7245\,3\, gr,\,\lambda = 41.4490\,3\, gr$. Par des observations astronomiques, on a déterminé les coordonnées astronomiques $\varphi_a = 10.7257\,4\, gr$   et $\lambda_a = 41.4505\,2\, gr$  du point $A$ et l'azimut astronomique de la direction $AB$  soit   $Aza = 89.6849\,9\, gr$.	Transformer l'azimut astronomique de la direction $AB$ en azimut géodésique en utilisant l'équation de Laplace donnée par :$$ Azg = Aza + (\lambda -\lambda_a ).sin\varphi $$
4.	Calculer le gisement $G$ de la direction $AB$ sachant que $\varphi_0=40.00\,gr, \lambda_0=+11.00\,gr$ et la correction de la corde $Dv = 0.0018\,8\, gr$.

5.	Les coordonnées Lambert de $A$ sont $X_A = 478\,022.43\, m$ et $Y_B = 444\,702.22\, m$. Déterminer alors les coordonnées Lambert de $B$.

6.	Calculer l'azimut de $B$ vers $A$ sachant qu'on néglige la correction de la corde de la direction $BA$ et que $\lambda_B =  10.9288\,4 \,gr$. 
\epb
\bpb
 On a mesuré une distance suivant la pente entre les points $A ( H_A =  1319.79 \,m)$ et  $B      ( H_B = 1025.34\, m)$ avec $D_P = 16\,483.873\, m$.

1.	Calculer la distance $D_e$ distance réduite à l'ellipsoïde de référence par la formule rigoureuse, on prendra le rayon de la Terre $R= 6378\, km$.

2.	Calculer la distance  $D_r$  réduite à la représentation plane Lambert si l'altération linéaire de la zone est de $- 14\, cm/km$.

3.	La direction $AB$  a un azimut géodésique $Azg = 297.5622\,5\, gr$. Donner l'expression du gisement $G$ de $AB$ en fonction de $Azg,\gamma$ la convergence des méridiens et $Dv$ la correction de la corde, sachant que la représentation plane Lambert utilisée a comme $\varphi_0=37\,gr,\,\lambda_0=+11.00\,gr$ que le point $A$ est au nord du  parallèle origine.   

4.	On donne $Dv = - 13.7\, dmgr$ et $\lambda    = 9.3474\, 734\, gr$ la longitude de $A$, calculer $G$.

5.	En déduire les coordonnées  $(X_B, Y_B )$ de $B$  si  $X_A = 363\,044.79\, m$  et   $Y_A  = 407\,020.09\, m$.

6.	Déterminer les coordonnées géographiques $(\varphi,\lambda)$ de $B$. 

On rappelle que: $a = 6\,378\,249.20\, m$ et $e^2 = 0.006\,803\,487\,7$.
\epb
\chapter{\textit{\textbf{La Représentation Plane UTM}}}
\section{\textsc{Définition et Propriétés}}
La représentation plane UTM (Universal Transverse Mercator) \index{La représentation UTM}est l'une des  représentations la plus utilisée dans le monde.
\\

C'est une représentation :

-	conforme d'un modèle  ellipsoïdique,

-	transverse : c'est-à-dire l'image de l'équateur (en partie) est l'axe $Ox$ (vers l'Est) et l'image d'un méridien appelé méridien central, de longitude $\lambda_0$ qu'on suppose égale à 0, est l'axe $Oy$ (vers le Nord) du plan.

Les coordonnées rectangulaires d'un point sont des fonctions de la forme :
\be
\fbox{$ \begin{array}{l}
	X=X(\varphi, \lambda) \\
	Y=Y(\varphi, \lambda)
\end{array} $}
\ee

où $(\varphi,\lambda )$ sont les coordonnées du point sur le modèle ellipsoïdique.
\\

  Soit un point $M(\varphi ,0)$ sur le méridien origine,  alors les coordonnées de $m$ son image sur le plan sont :
\ba
	 X(\varphi,0) = 0  \nonumber \\
 Y(\varphi,0) = Y(\varphi) \nonumber 
\ea
$Y(\varphi )$ sera déterminée en imposant que le long du méridien central ou origine, les longueurs sont conservées. Sur le méridien, la longueur est donnée par :
\be
	\fbox{ $ \beta(\varphi)=\ds \int^{\varphi}_{0}\rho ds $}
\ee
d'où:  
	\[              \beta(\varphi) = Y(\varphi) = Y(\varphi,0)            \]
	\\

\section{\textsc{Détermination des coordonnées UTM}}
\subsection*{11.2.1. Calcul Direct}
              Sur l'ellipsoïde, on a :
$$	      ds^2 =  \rho^2d\varphi^2 + N^2cos^2\varphi d \lambda ^2     $$
le carré de l'élément de longueur infinitésimale, avec : 
	\[                                               N(\varphi ) = a(1 - e^2sin^2\varphi )^{-1/2}                           \]
	\[                                    \rho( \varphi) = a(1-e^2)(1-e^2sin^2\varphi )^{-3/2}       \]
respectivement les rayons de courbure de la grande normale et de la méridienne, $a$ le demi-grand axe et $e$ la première excentricité de l'ellipsoïde de référence. On peut écrire que :                                            
$$  ds^2 = N^2cos^2\varphi\left( \frac{\rho^2d\varphi^2}{N^2cos^2\varphi} +  d \lambda^2\right)     $$
En posant : 
\ba
	 dL=\frac{\rho d\varphi}{Ncos\varphi} \nonumber \\
	 \mbox{ou} \,\,\, L=Logtg(\frac{\pi}{4}+\frac{\varphi}{2})-\frac{e}{2}Log\frac{1+esin\varphi}{1-esin\varphi}
\ea
avec $L$ la latitude isométrique, on a alors les coordonnées $(L,\lambda )$ symétriques et orthogonales. L'expression de $ds^2$ est égale à :
\be
	     \fbox{ $  ds^2 = N^2cos^2\varphi( dL^2 +  d \lambda^2)  $}   
\ee
Sur le plan, on a : 
$$ 	dS^2  = dX^2 + dY^2 $$
On pose : 
\be
\begin{array}{l}
	z = L +i\lambda \\
	Z=Y+iX
	\end{array}
\ee
où $i$ désigne le nombre complexe tel que $i^2= -1$. Entre $z$ et $Z$, on a la relation :
\be
	Z =  Y + iX = f(z)= f( L + i\lambda )
\ee
où $f$ est une fonction à déterminer. La représentation étant conforme, la fonction $f$ est par suite une fonction analytique dans $\BbC$ (l'ensemble des nombres complexes). La fonction $f$ est dérivable à tout ordre et développable en séries en tout point complexe. Considérons le point $z_0$ tel que $z_0 = L + i0= L$ et $z = L  + i\lambda$, ce qui donne $z - z_0 = i\lambda$.
\\

\begin{figure}
			\centering
				\includegraphics[width=0.80\textwidth]{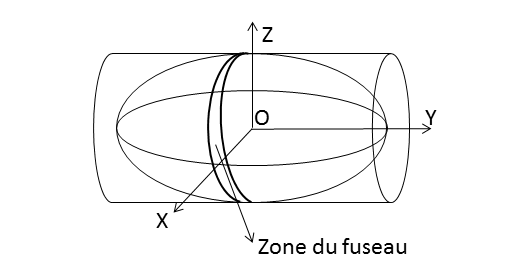}
			\caption{Interprétation géométrique de l'UTM}
			\label{fig:utmcylindre}
		\end{figure}

Dans la représentation UTM, on restreint $\lambda$ à varier dans l'intervalle $\left[-3°,+3°\right]$. Cet intervalle définit un fuseau de méridien central $\lambda_0=0°$ et d'amplitude $6°$. Ainsi, la Terre est divisé en $360°/6°=60$ fuseaux qu'on numérote de 1 à 60 ce qui explique l'utilisation mondialement de la représentation UTM. Une interprétation géométrique de la représentation UTM 	est comme suit:
	
	- on considère un cylindre ayant une base elliptique, tangent à l'ellipsoïde modèle le long de la méridienne de longitude $\lambda=\lambda_0=0°$. A un point $M(\varphi,\lambda)$ appartenant au fuseau $\left[-3°,+3°\right]$ on lui fait correspondre un point $m$ du cylindre (\textbf{Fig. \ref{fig:utmcylindre}}).
	
	- après développement du cylindre sur le plan, on obtient l'image $m(X,Y)$.
		
		On revient maintenant au développement de la fonction $f$ au voisinage de $z_0$, on a l'expression :
\be
f(z)=f(z_0)+(z-z_0)f'(z_0)+\frac{(z-z_0)^2}{2!}f''(z_0)+....\frac{(z-z_0)^n}{n!}f^{(n)}(z_0)+...
\ee
On se limite à n = 8. D'où:
\ba
& Y+iX = \ds f(L)+i\lambda f'(L)-\frac{1}{2!}\lambda ^2f"(L)-i\frac{1}{3!}\lambda ^3f^{3}(L)+ \frac{1}{4!} \lambda ^4f^{4}(L)\nonumber & \\ &\ds +i\frac{1}{5!}\lambda ^5f^{(5)}(L) -\frac{1}{6!}\lambda ^6 f^{(6)}(L) -i\frac{1}{7!}\lambda ^7f^{(7)}(L) +\frac{1}{8!}\lambda ^8 f^{(8)}(L)+...& \nonumber
\ea
Pour $\lambda=0$, on a $Y+iX=f(L)$ soit:
$$ Y=f(L)=Y(\varphi,0)=\beta(\varphi)$$
On pose :         
\be
\begin{array}{l}
		 a_0 = f(L)=\beta(\varphi)  \\
	  a_n= \ds \frac{1}{n!}f^{n}(L)=\frac{1}{n!}\frac{d^{n}\beta(\varphi)}{dL^n}  
\end{array}
\ee
Ce qui donne :
$$ Y+iX=a_0+ia_1\lambda-a_2\lambda^2-ia_3\lambda^3+a_4\lambda^4+ia_5\lambda^5-a_6\lambda^6-ia_7\lambda^7+a_8\lambda^8+... $$
D'où :
\be
\fbox{$ \begin{array}{l}
	 X = a_1\lambda  - a_3\lambda^3 + a_5\lambda^5-a_7\lambda^7+...  \\  
	Y = \beta(\varphi) - a_2\lambda^2 + a_4\lambda^4 - a_6\lambda^6 + a_8\lambda^8+... 
	\end{array} $}
\ee
avec :      
$$ 	 a_0 = f(L)=  \beta(\varphi ),   \quad 	a_1= f'(L) =\ds  \frac{d\beta}{dL}= \frac{d\beta}{d\varphi}\frac{d\varphi}{dL} $$      
En posant :
$$ 	\eta^2 = e'^2cos^2\varphi,    \quad	 e'^2 = \ds \frac{e^2}{1-e^2}, \quad  t=tg\varphi $$
	avec $e'$ la deuxième excentricité, on obtient les coefficients  :
\be
\fbox{ $ \begin{array}{l}
	 \displaystyle a_1=Ncos\varphi  \\
	\\
		\displaystyle a_2=-\frac{1}{2}Ncos\varphi sin\varphi  \\
		\\
		\displaystyle a_3=-\frac{1}{6}Ncos^3 \varphi(1+\eta^2-t^2) \\
		\\
	 \displaystyle a_4=\frac{1}{24}Ncos^3 \varphi sin\varphi (5-t^2+9\eta^2+4\eta^4) \\
	\\
		\displaystyle a_5=\frac{1}{120}Ncos^5\varphi(5-18t^2+t^4+14\eta^2-58\eta^2t^2+13\eta^4) \\
		\\
\displaystyle a_6=-\frac{1}{720}Ncos^5\varphi sin\varphi(61-58t^2+t^4+270\eta^2-330t^2\eta^2+200\eta^4-232t^2\eta^4) \\
\\
\displaystyle a_7=-\frac{1}{5040}Ncos^7\varphi(61+131t^2+179t^4+331\eta^2-3298t^2\eta^2)\\
\\
 a_8=\ds \frac{1}{40320}Ncos^7\varphi sin\varphi(165-61t^2+537t^4+9679\eta^2-23278t^2\eta^2+9244\eta^4 +\\
\\
\quad \ds 358t^4\eta^2-19788t^2\eta^4) 
\end{array} $} \lb{tabc}
\ee

Le calcul de  $\beta(\varphi) = f(L) = \beta$  se calcule à partir du développement de  $\beta(\varphi)$ en fonction de $u= e^2sin^2(\varphi )$ car  $u<1$. On exprime $sin^m \varphi$ en fonction de $sin$ de multiples de $\varphi$ soit $sinp\varphi$. En intégrant, on arrive à (voir en Annexe du présent chapitre):  
\ba
		 \beta(\varphi)= a(1-e^2).(C_0\varphi +C_2sin2\varphi +C_4sin4\varphi +C_6sin6\varphi  \nonumber \\
		+C_8sin8\varphi  +C_{10}sin10\varphi +C_{12}sin12\varphi)                                       
\ea
avec: 
\be 
\begin{array}{l}
 C_0 =\ds 1+\frac{3}{4}e^2+\frac{45}{64}e^4+\frac{175}{256}e^6+\frac{11025}{16384}e^8 +\frac{43659}{65536}e^{10}  +\frac{693693}{1048576}e^{12}\\
\\ 
C_2 =\ds-\frac{3}{8}e^2- \frac{15}{32}e^4- \frac{525}{1024}e^6- \frac{2205}{4096}e^8- \frac{72765}{131072}e^{10}-\frac{297297}{524288}e^{12}\\
\\
 	C_4 =\ds \frac{15}{256}e^4+ \frac{105}{1024}e^6 + \frac{2205}{16384}e^8 + \frac{10395}{65536}e^{10}+\frac{1486485}{8388608}e^{12} 
	\end{array}
	\ee
	\be
	\begin{array}{l}
 	 C_6 =\ds - \frac{35}{3072}e^6  -  \frac{315}{12288}e^8-\frac{31185}{786432}e^{10}-\frac{165165}{3145728}e^{12} \\
	\\
C_8 = \ds \frac{315}{131072}e^8 + \frac{3465}{524288}e^{10} + \frac{99099}{8388608}e^{12} \\
\\
	C_{10} = \ds- \frac{693}{1310720}e^{10} - \frac{9009}{5242880}e^{12} \\
	\\
	C_{12} = \ds \frac{1001}{8388608}e^{12} 
	\end{array}
\ee
On pose:
\be
	\Lambda=\lambda-\lambda_0
\ee
alors les formules définitives du calcul direct sont en s'arrêtant à l'ordre 8:
\be
\fbox{ $ \begin{array}{l}
	 X = a_1\Lambda-a_3\Lambda^3+a_5\Lambda^5-a_7\Lambda^7 \\
	 Y=\beta(\varphi)-a_2\Lambda^2+a_4\Lambda^4-a_6\Lambda^6+a_8\Lambda^8
	\end{array} $}
\ee
En général, on applique à $X,Y$ un facteur de réduction $k = 0.9996$  et une constante de translation en $X$ de $500\,000\, m$, les coordonnées obtenues sont :
\be
\fbox{ $ \begin{array}{l}
    X' = k.X + 500\,000.00\,\mbox{m}  \\                                                        
      Y' = k.Y 
\end{array} $}
\ee

\subsection*{11.2.2. Calcul Inverse}
\begin{figure}[htp]
	\centering
		\includegraphics[width=0.80\textwidth]{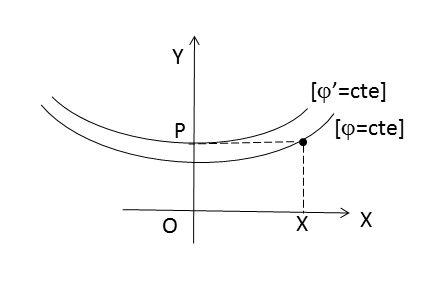}
	\caption{Passage de (X,Y) à $(\varphi,\lambda)$}
	\label{fig:utm3}
\end{figure}

         Ayant les coordonnées $(X',Y')$ en UTM et la longitude $\lambda_0$ du méridien central, comment on calcule $(\varphi,\lambda)$. On commence en revenant à: 
\be
 \begin{array}{l}
	 X= (X' -500000)/k  \\
 	 Y =Y'/k
\end{array}  
\ee
Par suite, en utilisant les variables :
     $ z = L + i(\lambda - \lambda_0)$ et  $Z = Y + iX$, on cherche à déterminer une fonction analytique $g$ telle que: 
\be
 \begin{array}{l}
	     z = g(Z)\\
	     \mbox{ou} \,\,\,\,  L  + i(\lambda - \lambda_0) =  g(Y + iX)   
    \end{array} 
\ee
Pour celà, on considère sur l'axe $OY$ le point $P(0,Y)$ (\textbf{Fig. \ref{fig:utm3}}), il lui correspond l'affixe $Z_0 = Y$, sur l'ellipsoïde il est l'image de $L'= L'(\varphi') = g(Z_0)$.
D'où le développement de $g$ au point $Z_0$: 
$$g(Z)=\ds g(Z_0)+(Z-Z_0)g'(Z_0)+\frac{(Z-Z_0)^2}{2!}g''(Z_0)+....\frac{(Z-Z_0)^n}{n!}g^{(n)}(Z_0)+...$$
Or: $ Z-Z_0=Y+iX-Y=iX$ ce qui donne :                                                                                               
\ba
& L+i(\lambda-\lambda_0)=\ds g(Z_0)+iXg'(Y)-\frac{1}{2}X^2.g"(Y)-\frac{1}{3!}iX^3g^{(3)}(Y)\nonumber & \\ & +\ds \frac{1}{4!}X^4g^{(4)}(Y) +\frac{1}{5!}iX^5g^{(5)}(Y)-  \frac{1}{6!}X^6g^{(6)}(Y)+... &\nonumber
\ea
On pose :        
$$ b_0 = g(Z_0) =g(Y) = L', \quad b_n= \ds \frac{g^{(n)}(Y)}{n!} =\frac{1}{n!}\frac{d^nL'}{dY^n} =\frac{1}{n!}\frac{d^nL'}{d\beta^n} $$
d'où :
\be
\fbox{ $ \begin{array}{l}
\lambda = \lambda_0 + b_1X - b_3X^3 + b_5X^5-b_7X^7+...     \\ 
\\
 L(\varphi)=L = L'(\varphi') - b_2X^2 + b_4X^4 - b_6X^6 +b_8X^8+...  
\end{array} $}
\ee
avec:                                                                                                   
$$ b_0 = g(Z_0)= L', \quad b_1= g'(Y) = \ds \frac{dL'}{dY} =\frac{dL'}{d\beta}  = \frac{dL'}{d\varphi'}.\frac{d\varphi'}{d\beta} $$ 
En posant  $\eta'^2 = e'^2cos^2\varphi '$ avec $e'$ la deuxième excentricité  $e'^2 = e^2/(1-e^2)$, $t'=tg\varphi' $ et $N'=N(\varphi')$, on obtient les coefficients: 
\be
\fbox{ $ \begin{array}{l}
 \displaystyle	b_1=\frac{1}{N'cos\varphi'} \\ 
\\
 \displaystyle	b_2=\frac{tg\varphi'}{2N'^2cos\varphi'} \\
\\
 \displaystyle	b_3=\frac{(1+2t'^2+\eta'^2)}{6N'^3cos\varphi'} \\
\\
 \displaystyle	b_4=\frac{tg\varphi'1(5+6t'^2+\eta'^2-4\eta'^4)}{24N'^4cos\varphi'}\\
\\
	\displaystyle	b_5=\frac{(5+28t'^2+6\eta'^2+24t'^4+8\eta'^2t'^2)}{120N'^5cos\varphi'} \\
	\\
	\displaystyle	b_6=\frac{tg\varphi'(61+180t'^2+46\eta'^2+120t'^4+48\eta'^2t'^2)}{720N'^6cos\varphi'} \\
	\\
	\displaystyle	b_7=\frac{(61+622t'^2+107\eta'^2+1320t'^4+1538\eta'^2t'^2+46\eta'^4)}{5040N'^7cos\varphi'}
	\end{array} $}
\ee
Ayant $L(\varphi)$, on calcule $ \varphi $ en utilisant la formule:
 $$L=Logtg\left( \frac{\pi}{4}+\frac{\varphi}{2}\right)-\frac{e}{2}Log\left(\frac{1+esin\varphi}{1-esin\varphi}\right) $$
Le calcul se fait par itérations.
\subsection*{11.2.3. Le Module linéaire}
Le module linéaire $m$ est tel que :
\be
	m^2  = \left(\frac{dS}{ds}\right)^2= \frac{dS^2}{ds^2}= \frac{dX^2 + dY^2}{\rho^2d \varphi^2 + N^2cos^2\varphi d \lambda^2}     
\ee
La représentation étant conforme, alors le module linéaire est indépendant de la direction, mais ne dépend que du point, on choisit de calculer $m$ le long des parallèles, soit $d\varphi =0$, ce qui donne :
$$  m^2 =  \frac{dX^2 + dY^2}{N^2cos^2\varphi d \lambda^2}                        $$
Comme on a:
 
$X = a_1\lambda  - a_3 \lambda^3 + a_5 \lambda^5+...  $

$Y=  \beta(\varphi)-a_2 \lambda^2 + a_4 \lambda^4 +...$

et que les coefficients $ a_i $ sont  des fonctions seulement de la latitude $(\varphi)$, d'où :
$$ 	   dX= a_1d\lambda  - 3a_3\lambda^ 2 d\lambda + 5a_5 \lambda^4d \lambda = (a_1 - 3a_3\lambda^ 2 + 5a_5\lambda^ 4 ) d\lambda+...$$
et: 
$$	dY = -2a_2\lambda d\lambda  + 4a_4 \lambda^3d\lambda+...$$

En gardant les termes en $\lambda$  et $\lambda^ 2$, on obtient : 
$$	 dX= Ncos\varphi  [1+ ( \lambda^2/2)cos^2\varphi (1-tg^2\varphi  + \eta^2)]d\lambda   $$
et:
$$	dY =  \lambda Ncos^2\varphi tg\varphi d\lambda $$
Par suite :  
$$	 dX^2 + dY^2 = N^2cos^2\varphi \left((1+ (\lambda^2/2)cos^2\varphi (1-tg^2\varphi  +\eta^2)\right)^2+   \lambda^2sin^2\varphi )d \lambda^2$$
En simplifiant et en négligeant les termes en  $\lambda^4$, on trouve :    
\be
	m=\sqrt{1+\lambda^2(1+\eta^2)cos^2\varphi}
\ee
Au lieu de prendre $m$ comme module linéaire, on le multiplie par un facteur $k$ appelé  facteur de réduction \index{Facteur de réduction de l'échelle}de l'échelle généralement égal à 0.9996. Le module linéaire devient :
\be
	m'=k\sqrt{1+\lambda^2(1+\eta^2)cos^2\varphi}
\ee
On remplace $\lambda$ par $\lambda-\lambda_0$, on trouve la formule du module linéaire :
\be
	\fbox{$ m'=k\sqrt{\displaystyle 1+(\lambda-\lambda_0)^2(1+e'^2cos\varphi^2)cos^2\varphi} $ }
\ee

\subsection*{11.2.4. Convergence des méridiens}
Le gisement de l'image du méridien appelé 'convergence des méridiens' et noté par $\gamma$ en un point $(\varphi,\lambda )$  est donné en première approximation  par la formule :
\be
\fbox{ $ 	tg\gamma=(\lambda-\lambda_0)sin\varphi $}
\ee

$\gamma$ est comptée dans le sens des gisements.
\section{\textsc{Annexe: Calcul de la longueur d'un arc de la méridienne d'un ellipsoïde de révolution}}
Soit $(E)$ un ellipsoïde de révolution défini par ses paramètres:

a: le demi-grand axe,

e: la première excentricité.

L'expression de la longueur de la méridienne entre l'équateur et un point $M$ de latitude géodésique $\varphi$ est donnée par:
\be
\beta=\beta({\varphi})=\int_0^{\varphi}\rho(u)du \label{ax1}
\ee
avec:
$$ \rho=\rho(u)=\ds \frac{a(1-e^2)}{(1-e^2sin^2u)^{\frac{3}{2}}} $$
$\rho$ est le rayon de courbure de la méridienne.

L'intégrale (\ref{ax1}) est une intégrale, dite elliptique,\index{Intégrale elliptique} n'est pas exprimée par une formule finie. Pour la calculer, on fait l'usage d'un développement limité de l'expression $\displaystyle (1-e^2sin^2u)^{-\frac{3}{2}}$.
\\

 On utilise la formule:
$$
(1+x)^q=\ds 1+qx+\frac{q(q-1)}{2!}x^2+\frac{q(q-1)(q-2)}{3!}x^3+...+\frac{q(q-1)...(q-1+p)}{p!}x^p+o(x^{p+1}) $$
avec $|x|< 1$, $q$ est un rationnel et $p!$ désigne factoriel $p$ soit $p(p-1)..3.2.1$. Comme $|e^2sin^2u| < 1$, on a donc à l'ordre 12:
\ba
&\ds \frac{1}{(1-e^2sin^2u)^{\frac{3}{2}}}=\ds (1-e^2sin^2u)^{-\frac{3}{2}}=\ds 1+\frac{3}{2}e^2sin^2u+\frac{15}{8}e^4sin^4u+\nonumber &\\ &\ds \frac{35}{16}e^6sin^6u+\frac{315}{128}e^8sin^8u+\frac{693}{256}e^{10}sin^{10}u+\frac{3003}{1024}e^{12}sin^{12}u \label{ax4}&
\ea
Pour pouvoir calculer les intégrales du type:
$$\int_0^{\varphi}sin^pudu $$
on va exprimer les termes $sin^pu$ en fonction des lignes trigonométriques multiples de l'argument $u$. 
Ce qui donne:
\be
\begin{array}{l}
sin^2u=\ds \frac{1}{2}-\frac{cos2u}{2} \\
\\
sin^4u=\ds \frac{3}{8}-\frac{cos2u}{2}+\frac{cos4u}{8} \\
\\
sin^6u=\ds \frac{5}{16}-\frac{15cos2u}{32}+\frac{13cos4u}{16}-\frac{cos6u}{32} \\
\\
sin^8u=\ds \frac{35}{128}-\frac{17cos2u}{16}+\frac{7cos4u}{32}-\frac{cos6u}{16}+\frac{cos8u}{128} \\
\\
sin^{10}u=\ds \frac{63}{256}-\frac{105cos2u}{256}+\frac{15cos4u}{64}-\frac{45cos6u}{512}+\frac{5cos8u}{256}-\frac{cos10u}{512} \\
\\
sin^{12}u=\ds \frac{231}{1024}-\frac{99cos2u}{256}+\frac{495cos4u}{2048}-\frac{55cos6u}{512}+\frac{33cos8u}{1024}-\frac{3cos10u}{512}+ \frac{cos12u}{2048} 
\end{array} \label{ax10}
\ee 
L'équation (\ref{ax4}) s'écrit en utilisant les expressions de droite de (\ref{ax10}):
\be
(1-e^2sin^2u)^{-\frac{3}{2}}=A_0+A_2cos2u+A_4cos4u+A_8cos8u+A_{10}cos10u+A_{12}cos12u \label{ax11}
\ee
En intégrant (\ref{ax11}) entre 0 et $\varphi$ et après multiplication par le coefficient $a(1-e^2)$, on trouve l'expression ci-dessous de la longueur de la méridienne: 
\ba
		 \beta(\varphi)= a(1-e^2).(C_0\varphi +C_2sin2\varphi +C_4sin4\varphi +C_6sin6\varphi  \nonumber \\
		+C_8sin8\varphi  +C_{10}sin10\varphi +C_{12}sin12\varphi)                                       
\ea
où les coefficient $A_k$ vérifient:
\ba
C_0=A_0\quad C_2=\frac{A_2}{2} \quad C_4=\frac{A_4}{4}\quad C_6=\frac{A_6}{6} \nonumber \\
C_8=\frac{A_8}{8} \quad C_{10}=\frac{A_{10}}{10} \quad C_{12}=\frac{A_{12}}{12}
\ea 
\newpage 
avec:
\be 
\begin{array}{l}
 C_0 =\ds 1+\frac{3}{4}e^2+\frac{45}{64}e^4+\frac{175}{256}e^6+\frac{11025}{16384}e^8 +\frac{43659}{65536}e^{10}  +\frac{693693}{1048576}e^{12}\\
\\ 
C_2 =\ds-\frac{3}{8}e^2- \frac{15}{32}e^4- \frac{525}{1024}e^6- \frac{2205}{4096}e^8- \frac{72765}{131072}e^{10}-\frac{297297}{524288}e^{12}\\
\\
 	C_4 =\ds \frac{15}{256}e^4+ \frac{105}{1024}e^6 + \frac{2205}{16384}e^8 + \frac{10395}{65536}e^{10}+\frac{1486485}{8388608}e^{12} \\
\\
 	 C_6 =\ds - \frac{35}{3072}e^6  -  \frac{315}{12288}e^8-\frac{31185}{786432}e^{10}-\frac{165165}{3145728}e^{12} \\
	\\
C_8 = \ds \frac{315}{131072}e^8 + \frac{3465}{524288}e^{10} + \frac{99099}{8388608}e^{12} \\
\\
	C_{10} = \ds- \frac{693}{1310720}e^{10} - \frac{9009}{5242880}e^{12} \\
	\\
	C_{12} = \ds \frac{1001}{8388608}e^{12} 
	\end{array}
\ee

\begin{flushright}
CQFD
\end{flushright}
\section{\textsc{Exercices et Problèmes}}
\bex
Dans cet exercice, on voudrait justifier l'arrêt à l'ordre 8 de l'expression de $Y(UTM)$ en fonction de $\lambda$. On donne: $\varphi=40.00 \,gr$ et $ a = 6\,378\,249.20\, m,\,\,\,e^2 = 0.006\,803\,4877$.

1. Calculer numériquement $e'^2,\eta^2,t^2=tg^2\varphi$ et $N(\varphi)$.

2. Calculer numériquement le coefficient $a_8$ de (\ref{tabc}).

3. On donne $\lambda= 1.235\,46\,gr$, calculer $a_8\lambda^8$ et conclure.
\eex
 
\bpb
 Soit le point $A$ de coordonnées géodésiques: $\varphi = 40.9193\, gr$ et $\lambda= 11.9656\, gr$ à l'Est de Greenwich. On considère la représentation plane UTM tronquée suivante, de méridien central $\lambda_0$ = 9° définie par les formules :
	\[  \left\{\begin{array}{ll}
	      X =  a_1.(\lambda- \lambda_0)+ a_3. (\lambda- \lambda_0)^3 \\
	        Y =   g(\varphi) + a_2.(\lambda- \lambda_0)^2  
	        \end{array} \right.
\]
où $\varphi,\,\lambda$ et $\lambda_0$ sont exprimées en $rd$, avec: 
	\[  a_1= N(\varphi).cos\varphi \]
	\[a_2 =  \frac{a_1}{2}. sin\varphi \]
	\[a_3 =  \frac{a_1cos^2\varphi}{6}(1- tg^2\varphi + e'^2.cos^2\varphi)\]
	\[N(\varphi) =  \frac{a}{\sqrt{1 - e^2.sin^2\varphi}}\]
	\[g(\varphi) = a(1 - e^2)( 1.0051353.\varphi -  0.0025731sin2\varphi )\]
	 \[ a = 6\,378\,249.20\, m,\,\,\,e^2 = 0.006\,803\,4877,\,\,\, \displaystyle e'^2=\frac{e^2}{1-e^2}   
\]

1.	Montrer que les coordonnées du point $A$ sont: $X = 157\,833.48\, m\,,Y = 4\,078\,512.97\, m$, on justifie les résultats.

2.	Soit le point $B$ de coordonnées $(X = 160\,595.98\,m ; Y = 4\,078\,564.53\, m)$. Sachant que $B$ est situé sur le même parallèle que $A$, calculer la longitude $\lambda'$ de $B$.

3.	Calculer le gisement $G$ et la distance $AB$.

4.	Sachant que la convergence des méridiens $\gamma$  est donnée par $ tg\gamma = (\lambda - \lambda_0)sin\varphi$   et qu'on néglige le $Dv$, calculer l'azimut de la direction $AB$.

5.	Calculer l'azimut de $B$ vers $A$ en négligeant le $Dv$ de $B$ vers $A$.

6.	En calculant les coordonnées UTM de $A$ et $B$, on trouve respectivement  $X_A=657\,770.34\, m,\, Y_A = 4\,076\,891.20\,m;\,  X_B = 660\,531.74\, m,\, Y_B = 4\,076\,942.76\, m$. Calculer la distance $AB$ par les coordonnées UTM. En déduire l'erreur relative sur la distance en utilisant les coordonnées de l'UTM tronquée.
\epb

 \chapter{\textit{\textbf{Les Transformations Entre Les Systèmes Géodésiques}}}\lb{transfor}
\section{\textsc{Introduction}}
        Avec l'introduction de la technologie de positionnement par GPS (Global Positioning System), laquelle fournit à l'utilisateur sa position $(X,Y,Z)$ tridimensionnelle dans le système géocentrique mondial dit $WGS84$ (World Geodetic System 1984), il est nécessaire de savoir la transformation de passage\index{Transformation de passage} du système géodésique mondial au système géodésique national ou local. On présente ci-après quelques modèles de transformations de passage entre les systèmes géodésiques.
\\
        
 On utilise par la suite les notations suivantes :
 
-	$(X_1,Y_1,Z_1)$  les coordonnées cartésiennes 3D dans le système local (système 1);

-	$(X_2,Y_2,Z_2)$  les coordonnées cartésiennes 3D dans le système géocentrique $WGS84$ (système 2);

-	$(\varphi_1, \lambda_1,he_1)$  les coordonnées géodésiques dans le système 1;

-	$(\varphi_2, \lambda_2,he_2)$ les coordonnées géodésiques dans le système 2.
\\
           
On propose d'étudier les transformations suivantes (\textit{C. Boucher}, 1979b; \textit{T. Soler}, 1998) :\index{\textbf{Boucher C.}} \index{\textbf{Soler T.}}
           
-	le modèle de Bur$\breve{s}$a\footnote{\textbf{Milan Bur}$\breve{\textbf{s}}$\textbf{a}: géodésien tchèque.}- Wolf\footnote{\textbf{Helmut Wolf} (1910-1994): géodésien allemand.} ou Helmert à 7 paramètres,\index{\textbf{Wolf H.}}\index{\textbf{Molodensky M.S.}}\index{\textbf{Bur$\breve{s}$a M.}}

-	les formules de Molodensky\footnote{\textbf{Mikhail Sergeevich Molodensky} (1909-1991): géodésien et géophysicien russe.},

-	les transformations bidimensionnelles.
\section{\textsc{Le Modèle de  BURSA - WOLF}}\index{Modèle de Bur$\breve{s}$a-Wolf}
Ce modèle s'écrit sous la forme vectorielle :
\be
	X_2 = T + ( 1+m).R(rx,ry,rz).X_1 \label{zz1}
\ee
où:

 - $X_2$  est le vecteur de composantes $(X_2,Y_2,Z_2)^T$, l'indice $T$ désigne transposée;
  
 - $T$ est le vecteur translation de composantes $(T_X,T_Y,T_Z)^T$  entre les systèmes 1 et 2;
       
 - $1+m$ est le facteur d'échelle entre les 2 systèmes;
      
 - $R(rx,ry,rz)$ est la matrice de rotation\index{Matrice de rotation} $3 \times 3$ pour passer du système 1 au système 2;
      
 -  $X_1$ est le vecteur de composantes $(X_1,Y_1,Z_1)^T$.
 \\
        
En développant (\ref{zz1}), on obtient:
      \be
			\fbox{ $ 
	\begin{pmatrix}{
	X_2 \cr  Y_2 \cr Z_2}
	\end{pmatrix}=	\begin{pmatrix}{
	T_X \cr T_Y \cr T_Z}
	\end{pmatrix}+(1+m)\begin{pmatrix}{
	1     &  rz    &   -ry  \cr
	-rz  &   1     &   rx \cr
	ry    &   -rx &  1}   
\end{pmatrix}.\begin{pmatrix}{
	X_1 \cr Y_1 \cr Z_1 }
	\end{pmatrix} $}  \label{z2}
\ee
avec $(rx,ry,rz)$ les rotations comptées positivement dans le sens contraire des aiguilles d'une montre.
Comment a-t-on obtenu cette formule? \\
\begin{figure}[htp]
	\centering
		\includegraphics[width=0.80\textwidth]{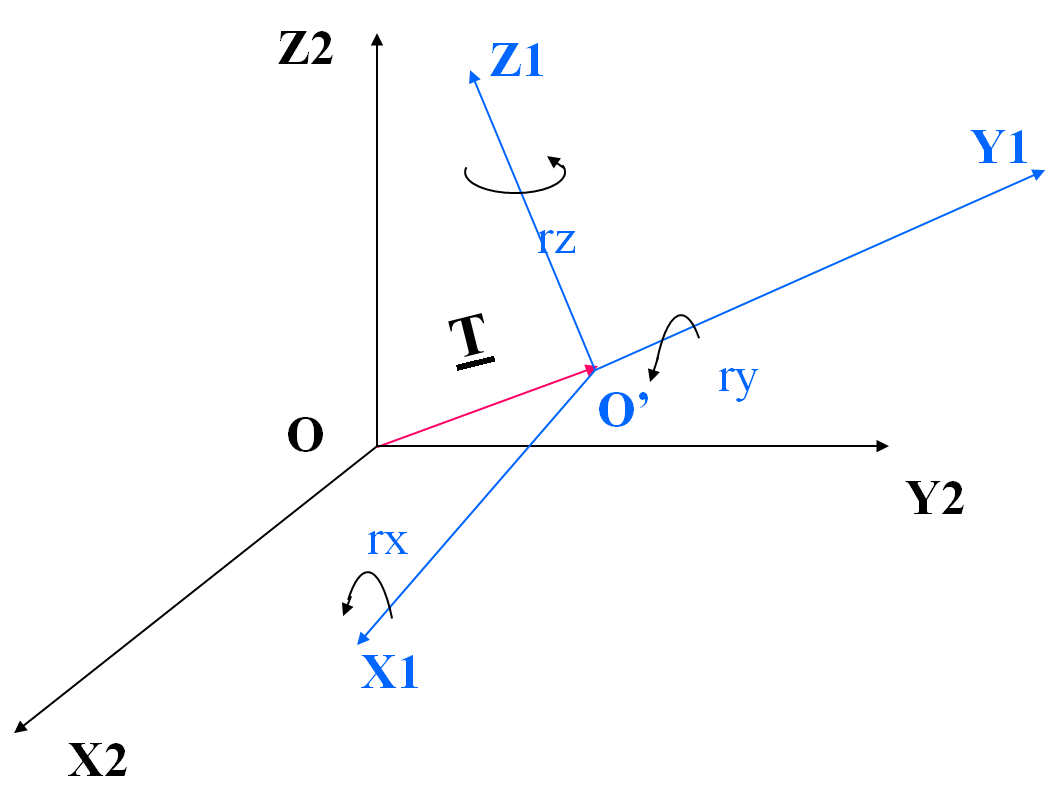}
	\caption{Le modèle de Bur$\breve{s}$a-Wolf}
	\label{fig:bursawolf}
\end{figure}
\newpage
Posons:
\ba
	\alpha=rx \label{y3} \\
	\beta=ry \label{y4} \\
	\gamma=rz \label{y5}
\ea

\subsection*{12.2.1. Matrices de Rotation}
Dans (\ref{z2}), $\alpha,\beta$ et $\gamma$ sont  les angles de rotation respectivement pour ramener les axes $O'X_1, O'Y_1$ et $O'Z_1$ parallèles aux axes $OX_2, OY_2$ et $OZ_2$. Faisant abstraction de la translation entre les systèmes 1 et 2, soit un point $M$ dans le plan $OX'_1Y'_1$ d'affixe $\xi_1=X_1+iY_1$, dans le plan $O'X_2Y_2$ le point $M$ a l'affixe $\xi_2=X_2+iY_2$. On peut écrire:
\ba
	\xi_2=X_2+iY_2=\rho e^{i\theta} \nonumber\\
	\xi_1=\rho e^{i(\theta+\gamma)}=\rho e^{i\theta}e^{i\gamma}=\xi_2e^{i\gamma} \nonumber
\ea
En passant aux coordonnées, on a:
\ba
	X_1+iY_1=e^{i\gamma}(X_2+iY_2) \nonumber \\
\Rightarrow X_2+iY_2=e^{-i\gamma}(X_1+iY_1)=(cos\gamma-isin\gamma)(X_1+iY_1) \nonumber
\ea
En séparant les parties réelles et imaginaires de l'équation précédente, on obtient:
\ba
	X_2=X_1cos\gamma+Y_1sin\gamma \nonumber \\
	Y_2=-X_1sin\gamma +Y_1cos\gamma \nonumber \\
	\mbox{et}\quad Z_2=Z_1 \nonumber
\ea
En les écrivant sous forme matricielle:
\begin{equation}
	\begin{pmatrix}{
	X_2 \cr
	Y_2 \cr
	Z_2}
	\end{pmatrix}=	\begin{pmatrix}{
	cos\gamma & sin\gamma & 0 \cr
	-sin\gamma & cos\gamma & 0  \cr
	0 & 0 & 1}
	\end{pmatrix}.\begin{pmatrix}{
	X_1 \cr
	Y_1 \cr
	Z_1}
	\end{pmatrix}=R(\gamma)	\begin{pmatrix}{
	X_1 \cr
	Y_1 \cr
	Z_1}
	\end{pmatrix} 
\end{equation}
avec :
\begin{equation}
	R(\gamma)=\begin{pmatrix}{
	cos \gamma & sin \gamma   & 0 \cr
	-sin \gamma   & cos \gamma & 0 \cr
	0 & 0 & 1}
\end{pmatrix} \label{y8}
\end{equation}
  Appelons $R(\alpha), R(\beta)$ les autres matrices de rotation. On a alors:
\ba
	R(\alpha)=\begin{pmatrix}{
	1 & 0 & 0 \cr
	0 & cos \alpha & sin \alpha \cr
	0 & -sin \alpha  & cos \alpha }
\end{pmatrix} \label{y6} \\
	R(\beta)=\begin{pmatrix}{
	cos \beta    & 0 &-sin \beta  \cr
	0 & 1      & 0     \cr
	sin \beta  & 0  & cos \beta }
\end{pmatrix}\label{y7}
\ea
Le modèle de Bur$\breve{s}$a-Wolf est obtenu comme suit:

- on fait subir une rotation autour de $O'X_1$ d'angle $\alpha $  de matrice de rotation $R(\alpha)$;

- on fait subir une rotation autour de $O'Y_1$ d'angle $\beta $ de matrice de rotation $R(\beta)$;
 
- on fait subir une rotation autour de $O'Z_1$ d'angle $\gamma$ de matrice de rotation $R(\gamma)$. 

Le résultat est la matrice :
\be
	R(\alpha,\beta,\gamma)=R(\gamma).R(\beta).R(\alpha) \label{y7a}
	\ee
Comme les angles de rotations sont petits $\leq 3^{\circ}$, on va exprimer chaque matrice $R$ en gardant seulement les termes du deuxième ordre. On utilise  les développements limités à l'ordre deux :
\ba
	sin\alpha\approx \alpha  \nonumber\\
	cos\alpha \approx 1-\frac{\alpha^2}{2} \nonumber
\ea
Alors les formules (\ref{y8}-\ref{y7}) deviennent:
\ba
	R(\alpha)=\begin{pmatrix}{
	1 & 0 & 0 \cr
		0 & 1- \frac{{\alpha}^2}{2} & \alpha \cr
		0 & -\alpha & 1- \frac{{\alpha}^2}{2} }
\end{pmatrix} \label{y12} \\
\nonumber \\
	R(\beta)=\begin{pmatrix}{
	1-\frac{{\beta}^2}{2}   & 0 &-\beta \cr
	0 & 1      & 0     \cr
	\beta & 0  & 	1-\frac{{\beta}^2}{2} }
\end{pmatrix} \label{y13} \\
\nonumber \\
	R(\gamma)=\begin{pmatrix}{
	1-\frac{{\gamma}^2}{2} & \gamma  & 0 \cr
	-\gamma  & 	1-\frac{{\gamma}^2}{2} & 0 \cr
	0 & 0 & 1 }
\end{pmatrix} \label{y14}
\ea
En revenant à la formule (\ref{y7a}), on obtient pour la matrice $R(\alpha,\beta,\gamma)$ l'expression suivante à l'ordre 2:
\begin{equation}
	R(\alpha,\beta,\gamma)=\begin{pmatrix}{
	1-\frac{{\gamma}^2}{2}-\frac{{\beta}^2}{2} & \gamma+\alpha \beta & -\beta+\alpha \gamma \cr
	\cr
	-\gamma & 1- \frac{{\gamma}^2}{2}-\frac{{\alpha}^2}{2} & \alpha+\beta \gamma \cr
	\cr
	\beta & -\alpha & 1- \frac{{\alpha}^2}{2}-\frac{{\beta}^2}{2} }
\end{pmatrix} \label{y15} 
\end{equation}
Maintenant, comme les trois angles sont petits, on va considérer que les termes du premier ordre ce qui donne pour $R(\alpha,\beta,\gamma)$:
\begin{equation}
	R(\alpha,\beta,\gamma)=\begin{pmatrix}{
	1 & \gamma & -\beta \cr
		-\gamma & 1 & \alpha \cr
		\beta & -\alpha & 1}
\end{pmatrix} \label{y16} 
\end{equation}
On revient à $(rx,ry,rz)$, on trouve:
\begin{equation}
	R(rx    ,ry   ,rz    )=\begin{pmatrix}{
	1 & rz     & -ry    \cr
		-rz     & 1 & rx     \cr
		ry    & -rx     & 1}
\end{pmatrix} \label{y17} 
\end{equation}
La formule (\ref{z2}) s'écrit:
      \begin{equation}
		\begin{pmatrix}{
	X_2 \cr Y_2 \cr Z_2 }
	\end{pmatrix}=	\begin{pmatrix}{ 
	T_X \cr T_Y \cr T_Z }
	\end{pmatrix}+(1+m)	\begin{pmatrix}{
	0     &  rz     &  -ry  \cr
 -rz  &   0     &  rx \cr
	 ry    &  -rx &  0   }
\end{pmatrix} 	\begin{pmatrix}{
	X_1 \cr Y_1 \cr Z_1 }
	\end{pmatrix} \label{y18}
\end{equation}
\subsection*{12.2.2. Calcul des Paramètres du Modèle de Bur$\breve{s}$a-Wolf par les Moindres Carrés}
 En considérant comme inconnues les paramètres $T_X,T_Y,T_Z,m,rx,ry,rz$, l'équation (\ref{z2}) s'écrit en gardant les termes du 1er ordre comme suit :
      \be
	\begin{pmatrix}{
	X_2-X_1 \cr Y_2-Y_1 \cr Z_2-Z_1}
	\end{pmatrix}=	\begin{pmatrix}{
	1 & 0 & 0 & X_1 & 0 & -Z_1 & Y_1\cr
	0 & 1 & 0 & Y_1 & Z_1 & 0 & -X_1 \cr
	0 & 0 & 1 & Z_1 & -Y_1 & X_1 & 0  }
\end{pmatrix}\begin{pmatrix}{
	T_X \cr T_Y \cr T_Z \cr m \cr rx \cr ry \cr rz }
	\end{pmatrix} \label{z3}
\ee
En utilisant l'équation (\ref{z3}) pour les n points communs dans les systèmes 1 et 2 et en posant :
\ba
	L = (X_{2i} - X_{1i})_{i =1,n}\nonumber \\
		U = ( T_X,T_Y,T_Z,m,rx,ry,rz)^T \nonumber
\ea
$A$ est la matrice $3n\times7$:
      \be
	A=\,_{3n}A_7=		\begin{pmatrix}{
	1 & 0 & 0 & X_i & 0 & -Z_i & Y_i \cr
	0 & 1 & 0 & Y_i & Z_i & 0 & -X_i \cr
	0 & 0 & 1 & Z_i & -Y_i & X_i & 0 }
\end{pmatrix}_{i=1,n}
\ee
et $V$ le vecteur des résidus de la méthode des moindres carrés, la détermination des paramètres inconnus se fait par la résolution par les moindres carrés de l'équation :
\be
	     AU= L+V
\ee
Soit :                                                                   
\be
\fbox{ $ \ov{U}  =  (A^T.A)^{-1}.A^T.L    $}
\ee
Le vecteur résidu est donné par :   
$$ 	V = A.\ov{U} - L = A.(A^T .A)^{-1}.A^T.L - L $$
Le facteur de la variance unitaire est exprimé par la formule:
\be
	\fbox{ $ \sigma^2 = \ds \frac{V^TV}{3n-7} $}
\ee
et la matrice variance-covariance du vecteur $\ov{U}$ est donnée par:
\be
\fbox{ $	\sigma_{\ov{U}}=\sigma_0 ^2(A^TA)^{-1} $}
\ee
\section{\textsc{Les Formules de MOLODENSKY}}\index{Formules de Molodensky}\index{\textbf{Molodensky M.S.}}
Un point $M$ a ses  coordonnées géodésiques cartésiennes 3D dans un référentiel donné comme suit: 
\be
\begin{array}{l}
	X = (N + he).cos\varphi .cos\lambda  \\
 Y = (N + he).cos\varphi .sin\lambda  \\   
 Z = (N(1 - e^2)+ he).sin\varphi 
\end{array} \lb{z12}
\ee
avec :  

- $ N = a.(1 - e^2.sin^2\varphi )^{-1/2}: $ le rayon de courbure de la grande normale;

- $a:$ le demi-grand axe  de l'ellipsoïde de référence;

- $e^2 = f(2-f):$ le carré de la 1ère excentricité;

- $f = (a-b)/a :$ l'aplatissement  de l'ellipsoïde de référence. 

On pose :
	\[	  w = (1-e^2.sin^2\varphi )^{-1/2}  	\]                                                               
Le rayon de courbure $\rho$ de la méridienne s'écrit:     
$$ 	 \rho = a(1-e^2)w^3   $$
Soit:
\be
	dX=\begin{pmatrix}{
dX \cr dY \cr dZ }	
\end{pmatrix}, \,\,  d\Phi =\begin{pmatrix}{
d\lambda \cr d\varphi \cr dhe }
\end{pmatrix}, \,\, dF=\begin{pmatrix}{
da \cr df }
\end{pmatrix} \label{z16}
\ee
En calculant $(dX,dY,dZ)$ des équations (\ref{z12}) en fonction de $d\varphi,d\lambda, dhe,da$ et $df$ et sachant que  $d(Ncos\varphi)= -\rho.sin\varphi .d\varphi$, on trouve :
\be
dX =  J.d\Phi +  K.dF  \label{z17}
\ee
où les matrices $J$ et $K$ sont les suivantes :
      \be
		J=\left(\begin{matrix}{
	-(N+he)cos\varphi sin\lambda & -(\rho+he)sin\varphi cos\lambda  & cos\varphi cos\lambda  \cr
	(N+he)cos\varphi cos\lambda & -(\rho+he)sin\varphi sin\lambda  & cos\varphi sin\lambda  \cr
	0       & (\rho+he)cos\varphi  & sin\varphi}  
\end{matrix}\right) \label{z18}
\ee
\be
	K=\left(	\begin{matrix}{
	wcos\varphi cos \lambda & \rho sin^2\varphi cos\varphi cos\lambda/(1-f)  \cr
	wcos\varphi sin \lambda & \rho sin^2\varphi cos\varphi sin\lambda/(1-f)  \cr
	w(1-e^2)sin\varphi  & (1-f)sin\varphi(\rho sin^2 \varphi -2N)}
\end{matrix}\right) \label{z19}
\ee
De l'équation (\ref{z17}), on tire:                                
\be
	d\Phi   = J^{-1}.dX  - J^{-1}.K.dF \label{z20}
\ee
avec:
      \be
		J^{-1}=\left(\begin{matrix}{
\displaystyle \frac{-sin\lambda}{(N+he)cos\varphi } & \displaystyle \frac{cos \lambda }{(N+he)cos\varphi } & 0  \cr
\displaystyle \frac{-sin\varphi cos\lambda }{\rho+he} & \displaystyle \frac{-sin\varphi sin\lambda }{\rho+he}  & \displaystyle \frac{cos\varphi }{\rho+he} \cr 
	cos\varphi cos \lambda  & cos\varphi sin \lambda  & sin\varphi  } 
\end{matrix}\right) \label{z21}
\ee
\be
	J^{-1}K=\left(\begin{matrix}{
	0 & 0  \cr
\displaystyle \frac{-we^2sin\varphi cos\varphi}{\rho+he} & \displaystyle \frac{-\rho sin\varphi cos\varphi (2-e^2sin^2\varphi)}{(\rho+he)(1-f)}  \cr
	\displaystyle \frac{1}{w} & -N(1-f)sin^2\varphi }   
\end{matrix}\right) \label{z22}
\ee
Or, en prenant:
\be
		dX=\left(\begin{matrix}{
X_2-X_1 \cr Y_2-Y_1 \cr Z_2-Z_1}	
\end{matrix}\right)
\ee
on a d'après l'équation (\ref{z3}) $dX = A.U$, par suite en posant:  
\ba
	          d\varphi   =   \varphi_2 -  \varphi_1 \nonumber \\
	     d\lambda  = \lambda_2 - \lambda_1  \nonumber \\
	   dhe = he_2 - he_1 \nonumber \\
	      da =  a_2 - a_1  \nonumber \\   
				df = f_2  - f_1 \nonumber
\ea
avec $(a_1, f_1)$ et $(a_2, f_2)$ sont respectivement les demi-grands axes et les aplatissements des ellipsoïdes  des systèmes 1 et 2, on a alors : 
\be
	d\Phi  = J^{-1}.A.U - J^{-1}.K.dF \label{z29}
\ee
avec $J^{-1}A $ la matrice $3\times7$ donnée ci-dessous: 
      \ba
& J^{-1}A= 	\begin{pmatrix}{
\displaystyle \frac{-sin\lambda}{(N+he)cos\varphi} & \displaystyle \frac{cos \lambda}{(N+he)cos\varphi} & 0 & 0   \cr 
\displaystyle \frac{-sin\varphi cos\lambda}{\rho+he} & \displaystyle \frac{-sin\varphi sin\lambda}{\rho+he} & \displaystyle \frac{cos\varphi}{\rho+he} & \displaystyle  \frac{-e^2N cos \varphi sin\varphi}{\rho+he}    \cr
	cos\varphi cos \lambda  & cos\varphi sin \lambda  & sin\varphi & N(1-e^2sin^2\varphi)+he    }
\end{pmatrix} \nonumber & \\
& \begin{pmatrix}{
\displaystyle \frac{(N(1-e^2)+he)tg\varphi cos\lambda}{N+he} & \displaystyle \frac{(N(1-e^2)+he)tg\varphi sin\lambda}{N+he} & -1  \cr
 \displaystyle \frac{-(N(1-e^2sin^2\varphi)+he)sin\lambda}{\rho+he} &  \displaystyle \frac{(N(1-e^2sin^2\varphi)+he)cos\lambda}{\rho+he} & 0 \cr
 -e^2Nsin\varphi cos\varphi sin\lambda & e^2Ncos\varphi sin\varphi cos\lambda & 0 }
\end{pmatrix} \nonumber&
\ea
En développant l'équation (\ref{z29}), on obtient les formules de \textbf{MOLODENSKY} de passage du système 1 au système 2:
\be
\fbox{ $ \begin{array}{l}
 \lambda_2-\lambda_1=\displaystyle -\frac{sin\lambda_1}{(N_1+he_1)cos\varphi_1}T_X+ \frac{cos \lambda_1}{(N_1+he_1)cos\varphi_1}T_Y+ \\
\\
\ds \frac{(N_1(1-e^2_1)+he_1)tg\varphi_1cos\lambda_1}{N_1+he_1}rx + \displaystyle \frac{(N_1(1-e^2_1)+he_1)tg\varphi_1sin\lambda_1}{N_1+he_1}ry-rz 
\end{array}$}\label{z32}
\ee 
\\

\be
\fbox{ $ \begin{array}{l}
\varphi_2-\varphi_1=\displaystyle -\frac{sin\varphi_1cos\lambda_1}{\rho_1+he_1}T_X-\frac{sin\varphi_1sin\lambda_1}{\rho_1+he_1}T_Y+ \frac{cos\varphi_1}{\rho_1+he_1}T_Z -\frac{e^2_1N_1cos\varphi_1sin\varphi_1}{\rho_1+he_1}m  \\
\\
-\displaystyle \frac{(N_1(1-e^2_1sin^2\varphi_1)+he_1)sin\lambda_1}{\rho_1+he_1}rx- \frac{(N_1(1-e^2_1sin^2\varphi_1)+he_1)cos\lambda_1}{\rho_1+he_1}ry \\
\\
\ds+ \frac{w_1e^2_1sin\varphi_1cos\varphi_1}{\rho_1+he_1}\Delta a+ \displaystyle \frac{\rho_1sin2\varphi_1(2-e^2_1sin^2\varphi_1)}{2(\rho_1+he_1)(1-f_1)}\Delta f 
\end{array} $} \label{z33}
\ee 
\\

\be
\fbox{ $ \begin{array}{l}
he_2-he_1= cos\varphi_1cos\lambda_1T_X+cos\varphi_1sin\lambda_1T_Y+sin\varphi_1T_Z+(N_1(1-e^2_1sin^2\varphi_1)+he_1)m  \\ 
\\
 -e^2_1N_1cos\varphi_1sin\varphi_1sin\lambda_1rx+e^2_1N_1cos\varphi_1sin\varphi_1cos\lambda_1ry -\displaystyle \frac{\Delta a}{w_1}+ N_1(1-f_1)sin^2\varphi_1\Delta f  
\end{array} $}\label{z34}
\ee
Des équations (\ref{z32})-(\ref{z34}), on remarque que :

* $ \lambda_2 - \lambda_1$ est indépendante de $ T_Z,m,a_1$  et $f_1$;

* $ \varphi_2 -  \varphi_1 $ est indépendante de $rz$;

*  $he_2 -he_1 $ est indépendante de $rz$.
\\
 
On trouve souvent dans la littérature géodésique  des formules de MOLODENSKY dites Standard et Abrégées qu'on donne ci-dessous.
\section{\textsc{Les  Formules de MOLODENSKY Standard}}\index{Formules de Molodensky standard}
Elles sont obtenues en ne tenant pas compte du facteur d'échelle et des rotations c'est-à-dire $m = 0$ et $rx = ry = rz = 0$ dans les formules (\ref{z32})-(\ref{z34}) et on obtient alors les formules suivantes en posant :
\be
\begin{array}{l}
	\Delta \varphi ''  =   \varphi_2 -  \varphi_1 \,\,\,  \mbox{en secondes sexagésimales} \\
	 \Delta \lambda ''  =   \lambda_2 -  \lambda_1\,\,\,   \mbox{en secondes sexagésimales} \\
	 \Delta he=he_2-he_1 \\
	 \Delta X=T_X \\
 \Delta Y=T_Y \\
	 \Delta Z= T_Z
	\end{array}
\ee
et en omettant les indices:
\be
\fbox{ $ \begin{array}{l}
	\Delta \varphi''= (-\Delta Xsin\varphi cos\lambda - \Delta Ysin\varphi sin\lambda + \Delta Zcos\varphi +  \\ 
 \ds Ne^2sin\varphi cos\varphi. \frac{\Delta a}{a}+\Delta f(\rho \frac{a}{b}+N\frac{b}{a}).sin\varphi cos\varphi).((\rho+he).sin1'')^{-1}  \\
\\
 \Delta \lambda''=(-\Delta Xsin\lambda+\Delta Ycos\lambda)((N+he)cos\varphi sin1'')^{-1}  \\
 \Delta he = \ds \Delta Xcos\varphi cos\lambda+\Delta Ycos\varphi sin\lambda+\Delta Zsin\varphi-a\frac{\Delta a}{N}+ \Delta f.N(1-f)sin^2\varphi 
\end{array} $}
 \ee 
avec $b$ le demi-petit axe de l'ellipsoïde 1.

\section{\textsc{Les Formules de MOLODENSKY Abrégées}}\index{Formules de Molodensky abrégées}
On fait $he=0$ et garde les termes du 1er ordre en $f$ dans les formules Standard, on trouve :
\be
\fbox{ $ \begin{array}{l} 
	\Delta \varphi'' = \displaystyle \frac{- \Delta Xsin\lambda + \Delta Ycos\lambda}{Ncos\varphi sin1''} \\
	\\
   \Delta \lambda''=\displaystyle \frac{-\Delta Xsin\varphi cos\lambda -\Delta Ysin\varphi sin\lambda +\Delta Zcos\varphi +(a\Delta f + f\Delta a)sin2\varphi}{\rho sin1''}  \\
	\\
 	\Delta he = \Delta Xcos\varphi cos\lambda + \Delta Ycos\varphi sin\lambda + \Delta Zsin\varphi - \Delta a  +(a\Delta f +f\Delta a )sin^2\varphi 
	  \end{array} $}
		\ee
\section{\textsc{La Recherche des Paramètres de Passage par les Formules de  MOLODENSKY}}
A partir de l'équation (\ref{z29}) on a :   
	\[d\Phi = J^{-1}.A.U - J^{-1}.K.dF
\]
soit :
\be
   J^{-1}.A.U=d\Phi + J^{-1}.K.dF
 \ee
où $U$ est le vecteur des inconnues $(T_X, T_Y, T_Z, m, rx, ry, rz)^T$. En écrivant l'équation précédente pour les n points communs et en posant :
\be
L = (d\Phi_i + J^{-1}_i.K_i.dF)_{i=1,n} 
\ee
le vecteur des observations $3n\times1$ et :
\be
	B = (J^{-1}_iA_i)_{i=1,n} 
\ee
la matrice des coefficients $3n\times7$ et $V$ le vecteur des résidus de la méthode des moindres carrés, la détermination des paramètres inconnus se fait par la résolution par les moindres carrés de l'équation :
\be
	 B.U =L + V 
\ee
Le vecteur solution est :
\be
\fbox{ $	\ov{U} = (B^T.B)^{-1}B^TL $}
\ee
Le vecteur résidu est :
\be
	V = B.\ov{U} - L= B.(B^T.B)^{-1}.B^T.L- L
\ee
Le facteur de la variance unitaire est donné par :
\be
	\fbox{ $ \sigma^2 = \ds \frac{V^TV}{3n-7} $} 
\ee
La matrice variance-covariance du vecteur $\ov{U}$ est donnée par :
\be
	\fbox{ $ \sigma_{\ov{U}}  = \sigma^2.(B^T.B)^{-1} $}
\ee
\section{\textsc{La Détermination des paramètres du Modèle de Bur$\breve{s}$a-Wolf}}

Dans ce paragraphe, on veut calculer manuellement les paramètres du modèle de Bur$\breve{s}$a-Wolf vu précédemment:
\be
		\textbf{\textit{X}}_2 = \textbf{\textit{T}} + (1+m).R(rx,ry,rz).\textbf{\textit{X}}_1   \label{yy1}
\ee
En développant (\ref{yy1}), on obtient:
      \be              
		\begin{pmatrix}{
	X_2 \cr
	 Y_2 \cr 
	  Z_2 }
	\end{pmatrix}=\begin{pmatrix}{
	T_X \cr
	 T_Y \cr
	  T_Z }
	\end{pmatrix}+(1+m)\begin{pmatrix}{
	1     &   rz    &  -ry  \cr
 -rz  &   1     &  rx \cr 
	 ry    &  -rx &  1 } 
\end{pmatrix}\begin{pmatrix}{
	X_1 \cr
	Y_1 \cr
	 Z_1 }
	\end{pmatrix}  \label{yy2}
\ee          
avec $(rx,ry,rz)$ les rotations comptées positivement dans le sens contraire des aiguilles d'une montre. Comment déterminer les paramètres modèle (\ref{yy1})?

\subsection*{12.7.1. Détermination de l'échelle $1+m$}
On suppose donné un ensemble de points $P_i$ pour $i=1,n$ connus dans les deux systèmes $S_1$ et $S_2$. On écrit l'équation (\ref{yy1}) pour deux points $P_j$ et $P_k$, d'où:
\ba
\textbf{\textit{X}}(P_j)_2 = \textbf{\textit{T}} + ( 1+m).R(rx,ry,rz).\textbf{\textit{X}}(P_j)_1  \label{yy3} \\
	\textbf{\textit{X}}(P_k)_2 = \textbf{\textit{T}} + ( 1+m).R(rx,ry,rz).\textbf{\textit{X}}(P_k)_1  \label{yy4} 
\ea
Par différence, on obtient :
\be
	(\textbf{P}_j\textbf{P}_k)_2=(1+m).R(rx,ry,rz).(\textbf{P}_j\textbf{P}_k)_1  \label{yy5} 
\ee
On prend la norme des deux membres de (\ref{yy5}) et que $1+m >  0$:
\be
\|(\textbf{P}_j\textbf{P}_k)_2\|=\|(1+m).R(rx,ry,rz).(\textbf{P}_j\textbf{P}_k)_1\|=(1+m)\|R(rx,ry,rz).(\textbf{P}_j\textbf{P}_k)_1 \|   \label{yy6} 
\ee
Comme $R$ est une matrice de rotation, donc son application à un vecteur est une isométrie, c'est-à-dire qu'elle laisse invariant la norme ou la longueur du vecteur, soit:
\be
	\|R.\textbf{\textit{X}}\|=\|\textbf{\textit{X}}\|, \quad \forall \,\textbf{\textit{X}} \in \BbR^3 \label{yy7}
\ee
On a donc:
\be
\|(\textbf{P}_j\textbf{P}_k)_2\|=(1+m)\|(\textbf{P}_j\textbf{P}_k)_1 \|   \label{yy8} 
\ee
Soit:
\be
\fbox{ $	1+m=\ds \frac{1}{N}\sum^N\frac{\|(\textbf{P}_j\textbf{P}_k)_2\|}{\|(\textbf{P}_j\textbf{P}_k)_1\|} $}\label{yy9}
\ee
$N$ désigne le nombre de couples de points $P_jP_k, j\neq k$.
\subsection*{12.7.2. Détermination des rotations $(rx,ry,rz)$}
Connaissant $(1+m)$, pour un couple de points $P_j,P_k$, on a :
\be
	(\textbf{P}_j\textbf{P}_k)_2=(1+m).R(rx,ry,rz).(\textbf{P}_j\textbf{P}_k)_1  \nonumber  
\ee
 On détaille la matrice $R$:
 \be           
	R=\begin{pmatrix}{
	1     &   rz    &  -ry  \cr
 -rz  &   1     &  rx \cr
	 ry    &  -rx &  1}  
\end{pmatrix}=	\begin{pmatrix}{
	1     & 0    &   0  \cr
 	0  &   1     & 0 \cr
	0    &   0 &  1    }
\end{pmatrix}+\begin{pmatrix}{
	0     &   rz    &  -ry  \cr
 -rz  &   0     &  rx \cr
	 ry    &  -rx &  0  } 
\end{pmatrix}=I_3+Q \label{yy10}
\ee             
avec $I_3$ la matrice unité et $Q$ la matrice:
\be
	Q=\begin{pmatrix} {
	0     &   rz    &  -ry  \cr
 -rz  &   0     &  rx \cr
	 ry    &  -rx &  0    }
\end{pmatrix}  \label{yy11}
\ee
Alors l'équation (\ref{yy5}) devient:
\be
	(\textbf{P}_j\textbf{P}_k)_2=(1+m).(I_3+Q(rx,ry,rz)).(\textbf{P}_j\textbf{P}_k)_1  \label{yy12} 
\ee
Comme $m \ll  1$ et $m^2\ll 1$, on a:
\be
Q(rx,ry,rz).(\textbf{P}_j\textbf{P}_k)_1=(1-m).(\textbf{P}_j\textbf{P}_k)_2-(\textbf{P}_j\textbf{P}_k)_1 \label{yy13} 
\ee
En posant:
\be 
	(\textbf{P}_j\textbf{P}_k)_2=\left|
\begin{array}{ll}
\Delta X'_{jk} \\
\Delta Y'_{jk} \\
\Delta Z'_{jk}
\end{array}\right.;\quad 	(\textbf{P}_j\textbf{P}_k)_1=\left|
\begin{array}{ll}
\Delta X_{jk} \\
\Delta Y_{jk} \\
\Delta Z_{jk}
\end{array}\right. ;\quad v=\left|
\begin{array}{ll}
v_1=(1-m)\Delta X'_{jk}-\Delta X_{jk} \\
v_2=(1-m) \Delta Y'_{jk}-\Delta Y_{jk} \\
v_3=(1-m)\Delta Z'_{jk}-\Delta Z_{jk}
\end{array}\right. \label{yy14}
\ee
Alors, on obtient l'équation:
\be               
	Q(rx,ry,rz).(\textbf{P}_j\textbf{P}_k)_1=v \label{yy15}
\ee
ou encore:
\be
	\begin{pmatrix} {
0& -\Delta Z_{jk} & \Delta Y_{jk}     \cr
	\Delta Z_{jk} & 0 &  -\Delta X_{jk} \cr
	     -\Delta Y_{jk} & \Delta X_{jk} & 0 }
\end{pmatrix}.\begin{pmatrix} {
	rx \cr
	ry \cr
	rz  }
\end{pmatrix}=\begin{pmatrix}{
	v_1 \cr
	v_2 \cr
	v_3 }
\end{pmatrix} \label{yy16}
\ee
Or le déterminant de la matrice $Q'$ :
\be
Q'=	\begin{pmatrix} {
0& -\Delta Z_{jk} & \Delta Y_{jk}     \cr
	\Delta Z_{jk} & 0 &  -\Delta X_{jk} \cr
	     -\Delta Y_{jk} & \Delta X_{jk} & 0 }
\end{pmatrix} \label{yy17}
\ee
est nul. Pour passer de cette conséquence, on utilise pour chaque ligne du système (\ref{yy16}) un couple de points $ij$ ce qui donne le système:
\be
\fbox{ $ \begin{pmatrix} {
0 & -\Delta Z_{jk} & \Delta Y_{jk}     \cr
	\Delta Z_{lm} & 0 &  -\Delta X_{lm} \cr
	     -\Delta Y_{in} & \Delta X_{in} & 0}
\end{pmatrix}.\begin{pmatrix}{
	rx \cr
	ry \cr
	rz }
\end{pmatrix}=\begin{pmatrix}{
	{v_{jk}}_1\cr
	{v_{lm}}_2\cr
	{v_{in}}_3 }
\end{pmatrix} $} \label{yy18}
\ee
Le système (\ref{yy18}) devient résolvable ce qui permet de déterminer les trois rotations $rx,ry$ et $rz$.
\subsection*{12.7.3. Détermination des composantes de la Translation $T$}
Les composantes $Tx,Ty,Tz$ du vecteur translation sont déterminées à partir des coordonnées des points $Pj$ communs dans les deux systèmes à partir de:
\ba
	Tx_j=X_{2j}-(1+m)(X_{1j}-rxY_{1j}+ryZ_{1j}) \label{yy19} \\
		Ty_j=Y_{2j}-(1+m)(rxX_{1j}+Y_{1j}-rzZ_{1j}) \label{yy20} \\
	Tz_j=Z_{2j}-(1+m)(-ryX_{1j}+rzY_{1j}+Z_{1j}) \label{yy21} 
\ea
Les composantes $Tx,Ty,Tz$ sont obtenues par une moyenne sur les $N$ points communs à savoir:
\be
\fbox{ $\begin{array}{l}
	Tx=\ds \frac{\sum^NTx_j}{N}  \\
		Ty=\ds \frac{\sum^NTy_j}{N}  \\
	Tz=\ds \frac{\sum^NTz_j}{N} 
	\end{array} $} \label{yy24}
\ee
\section{\textsc{La Transformation de HELMERT Bidimensionnelle}}\index{Transformation de Helmert}
Cette transformation  s'écrit sous la forme vectorielle :
\be
\fbox{ $ 	X_2 = T + s.R(\theta).X_1 $} \label{z60}
\ee
où:

- $X_2$  est le vecteur de composantes $(X_2,Y_2)^T$, $(X_2,Y_2)$ désignent les coordonnées planimétriques du système 2;

- $T$ est le vecteur translation de composantes $(T_x,T_y)^T$  entre les systèmes 1 et 2;

- $s$ est le facteur d'échelle entre les 2 systèmes;

- $R(\theta)$ est la matrice de rotation\index{Matrice de rotation} $2\times 2$ pour passer du système 1 au système 2;

- $X_1$ est le vecteur de composantes $(X_1,Y_1)^T$ où $X_1,Y_1$ désignent les coordonnées dans le système 1.

En développant (\ref{z60}), on obtient  :
\be
\begin{pmatrix}{
	  X_2 \cr
		Y_2 }
		\end{pmatrix}= \begin{pmatrix}{
		T_x\cr
	T_y }
		\end{pmatrix} +s \begin{pmatrix}{	
 cos\theta & - sin\theta \cr
 sin\theta  & cos\theta  }
		\end{pmatrix}.\begin{pmatrix}{	
    X_1 \cr
		Y_1 }
		\end{pmatrix} \label{z61}
		\ee
En prenant comme inconnues auxiliaires :
\ba
	    v = s.sin\theta   \\                                   
	    u = s.cos\theta   
\ea
le système (\ref{z61}) devient :
\be
	       \begin{array}{l}
												X_2 = T_x + X_1.u - Y_1.v \\
	                       Y_2 = T_y + X_1.v  + Y_1.u   
												\end{array} \lb{z64}              
\ee
Les inconnues $T_x, T_y$, $u$ et $v$ seront déterminées  par la méthode des moindres carrés en utilisant des points communs dans les deux systèmes. 
\begin{figure}[h]
	\centering
		\includegraphics{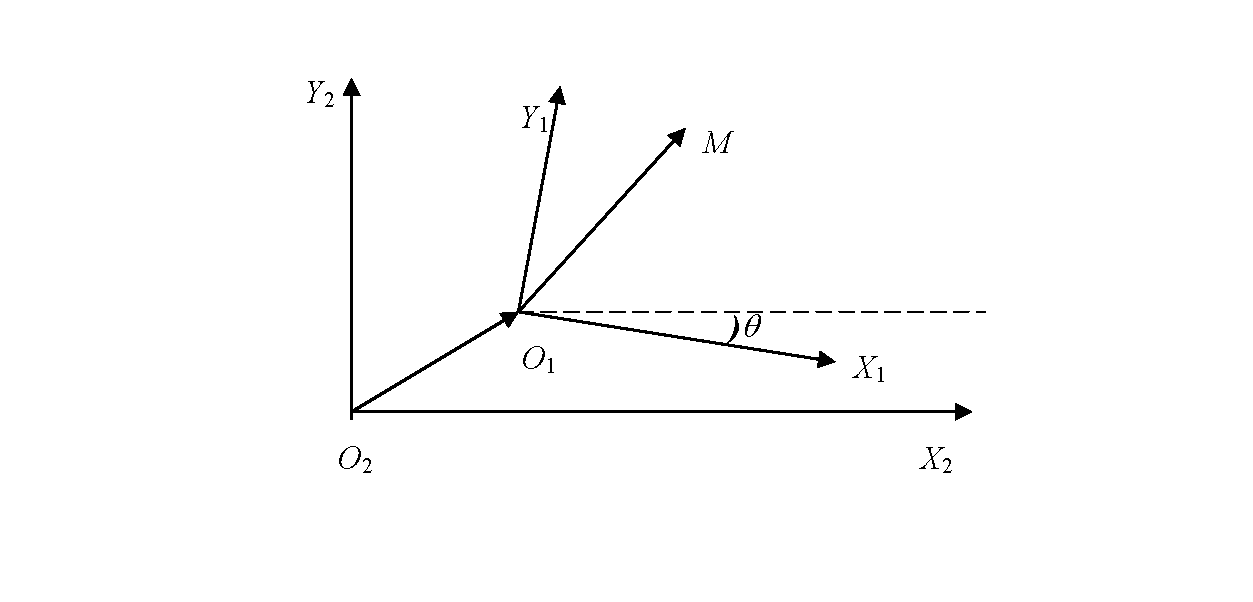}
	\caption{Modèle de Helmert}
	\label{fig:HELM}
\end{figure}
Ayant $u$ et $v$, on déduit :
\be
\fbox{ $ \begin{array}{l}
	                         s =\ds  \sqrt{u^2+v^2} \\
	                           tg \theta = \ds \frac{v}{u}
	\end{array} $}
\ee
\subsection*{12.8.1. R\'esolution par les Moindres Carr\'es}
On r\'esoud par la m\'ethode des moindres carr\'es le syst\`eme (\ref{z64}). On suppose la donn\'ee de $n$ points communs entre le syst\`eme $S_1$ et le syst\`eme $S_2$:

- $(X_i',Y_i')_{i=1,n}$ dans le syst\`eme $S_2$;

- $(X_i,Y_i)_{i=1,n}$ dans le syst\`eme $S_1$.
\\

 On pose:
\be
\ov   {X}=\ds \frac{\sum_1^nX_i}{n}, \quad \ov {Y}=\ds \frac{\sum_1^nY_i}{n}, \quad \ov {X'}=\ds \frac{\sum_1^nX_i'}{n}, \quad \ov {Y'}=\ds \frac{\sum_1^nY_i'}{n} 
\ee     
	les coordonn\'ees respectives des centres de gravit\'e, des points communs, dans $S_1$ et $S_2$.
	\\
	
	On pose de m\^eme:
\be
\begin{array}{l}             
x_i=X_i-\ov {X} \\
y_i=Y_i-\ov {Y} \\
x_i'=X_i'-\ov {X'} \\
y_i'=Y_i'-\ov {Y'}
	\end{array}
\ee                
	Dans ce cas, le syst\`eme (\ref{z64}) s'\'ecrit:
\be 
\begin{array}{l}
	 x_i' = T_x + x_i.u - y_i.v  \\
	 y_i' = T_y + x_i.v  + y_i.u     
		\end{array}\label{cq18} 
\ee     		
Soit $(T_x^0,T_y^0,u_0,v_0)$ une solution approch\'ee du syst\`eme. On note alors:
\be 
\begin{array}{l}
	T_x=T_x^0+dt_x \\
	T_y=T_y^0+dt_y \\
	u=u_0+du \\
	v=v_0+dv
\end{array}\label{cq19} 
\ee
	Alors les \'equations (\ref{cq18}) deviennent:
\be 
\begin{array}{l}          
	 x_i' =T_x^0+dt_x + x_i.(u_0+du) - y_i.(v_0+dv) \\
	 y_i' = T_y^0+dt_y + x_i.(v_0+dv)  + y_i.(u_0+du)              
\end{array}\label{cq20} 
\ee
On \'ecrit ces dernières \'equations sous la forme de l'\'equation des moindres carr\'es:
\be          
	A.X= L+ W \label{cq26}
\ee           
avec $X$ le vecteur des inconnues:
\be           
	X=\begin{pmatrix}{
	dt_x \cr
	dt_y \cr
	du \cr
	dv }
\end{pmatrix}
\ee           
$L$ le vecteur des observables:
\be             
	L=\begin{pmatrix}{
	x'_1-T_x^0-x_1u_0+y_1v_0 \cr 
	y'_1-T_y^0-x_1v_0-y_1u_0\cr  
\vdots \cr
	x'_i-T_x^0-x_iu_0+y_iv_0 \cr 
	y'_i-T_y^0-x_iv_0-y_iu_0\cr
	\vdots \cr
	x'_n-T_x^0-x_nu_0+y_nv_0 \cr
	y'_n-T_y^0-x_nv_0-y_nu_0}
\end{pmatrix}
\ee              
$W$ le vecteur des r\'esidus:
\be           
	W=\begin{pmatrix}{
	w_{x_1} \cr 
	w_{y_1} \cr 
	\vdots \cr 
	w_{x_n} \cr 
	w_{y_n}}
\end{pmatrix}
\ee           
et $A$ la matrice des coefficients:
\be           
	A=\begin{pmatrix}{
	1   &   0  &  x_1  &   -y_1  \cr
	0  &    1   & y_1  &    x_1 \cr
	\vdots &  \vdots & \vdots & \vdots \cr
		1   &   0  &  x_i  &   -y_i  \cr 
	0  &    1   & y_i  &    x_i \cr 
	\vdots &  \vdots & \vdots & \vdots \cr
		1   &   0  &  x_n  &   -y_n  \cr
	0  &    1   & y_n  &    x_n }
\end{pmatrix}
\ee              
\subsection*{12.8.2. La Solution par les Moindres Carr\'es}
La solution de (\ref{cq26}) par les moindres carr\'es donne:
\be
\overline{X}	=(A^TA)^{-1}A^TL
\ee
On pose:
\be
	N=A^TA
\ee
qu'on appelle matrice normale du syst\`eme (\ref{cq26}). On obtient alors:
\be           
	N=\begin{pmatrix}{
		n   &   0  &  \sum x_i  &   -\sum y_i  \cr
	0  &    n  & \sum y_i  &    \sum x_i \cr
		\sum x_i &   \sum y_i  &  \sum (x_i^2+y_i^2)  &   0  \cr
	-\sum y_i  &    \sum x_i   & 0  &    \sum(x_i^2+y_i^2) }
\end{pmatrix}
\ee          
Comme on travaille par rapport aux centres de gravit\'e des coordonn\'ees de $S_1$ et $S_2$, on a alors par d\'efinition:
\be
	\sum x_i=\sum y_i=\sum x'_i=\sum y'_i=0
\ee
De plus on note:
\be
	d_i^2=x_i^2+y_i^2
\ee
Alors la matrice $N$ s'\'ecrit facilement:
 \be             
	N=\begin{pmatrix}{
		n   &   0  &  0  &  0  \cr
	0  &    n  & 0  &  0 \cr
		0 &   0  &  \sum d_i^2  &   0  \cr
	0  &   0   & 0  &    \sum d_i^2 }
\end{pmatrix}
\ee           
La matrice normale est diagonale, son inverse est donn\'e par:
 \be          
	N^{-1}=\begin{pmatrix}{
		\ds \frac{1}{n}   &   0  &  0  &  0  \cr
	0  &   \ds \frac{1}{n} & 0  &  0 \cr
		0 &   0  & \ds \frac{1}{\sum d_i^2}  &   0  \cr
	0  &   0   & 0  &   \ds  \frac{1}{\sum d_i^2} }
\end{pmatrix}
\ee             
Or on sait que:
\be
	\sigma^2_{\ov{X}}=\sigma^2_0.N^{-1} \label{cq38}
\ee
o\`u $\sigma^2_0$ est le facteur de variance unitaire donn\'e par :
\be
	\sigma^2_0=\ds \frac{W^TW}{n-4}=\frac{\sum_i^nw_i^2}{n-4}
\ee
De l'\'equation (\ref{cq38}), on voit que :
\be
		\fbox{ $ \sigma^2_{dt_x}= \sigma^2_{dt_y}=\ds \frac{\sigma^2_0}{n} $}
\ee
\newpage
\textbf{Propri\'et\'e 13.1} \emph{Dans une transformation de Helmert bidimensionnelle, plus le nombre de points communs $n$ entre les deux syst\`emes est grand, plus la d\'etermination du vecteur translation $T=(T_x,T_y)^T$ est pr\'ecise.}
\\

Quant aux deux autres inconnues (facteur d'\'echelle et la rotation), on a:
\be
			\sigma^2_{u}= \sigma^2_{v}=\ds \sigma^2_0\frac{1}{\sum d_i^2}
\ee
 Soit $D=max(d_i)$. Si on veut imposer $\sigma^2_{u}$ égal \`a $\tilde{\sigma}^2_u$ donn\'e, quelles conditions doivent vérifier les $d_i$. On a:
 \be           
			\tilde{\sigma}^2_{u}= \ds \sigma^2_0\frac{1}{\sum d_i^2}\Rightarrow \sum d_i^2=\ds \frac{\sigma^2_0}{\tilde{\sigma}^2_u} \label{cqA}
\ee            
Or: $$d_i\leq D\Rightarrow \sum_{i=1}^{i=n} d_i^2\leq nD^2$$
Par suite:
\be
\fbox{ $ \ds \frac{\sigma^2_0}{\tilde{\sigma}^2_u}\leq nD^2\Rightarrow D^2\geq \ds \frac{\sigma^2_0}{n\tilde{\sigma}^2_u} $}
	\ee
	D'o\`u:
	
	\textbf{Propri\'et\'e 13.2} \textit{Dans une transformation de Helmert bidimensionnelle, en imposant un \'ecart-type donn\'e de la rotation $\tilde{\sigma}_u$, la distance maximale $D$ qu'on peut prendre vaut $\displaystyle {\frac{\sigma_0}{\tilde{\sigma}_u \sqrt{n}}}
$.}
\section{\textsc{Exercices et Problèmes}}
\bex
On donne le modèle bidimensionnel suivant, de transformation entre deux systèmes géodésiques, défini par:
	$$ \begin{pmatrix}{
	        X_2 \cr
					Y_2 }
					\end{pmatrix}=\begin{pmatrix}{
					-21.662\,m\cr
					-627.748\,m }
					\end{pmatrix}+\begin{pmatrix}{
					0.999\,988\,149 & - 0.000\,025\,928 \cr
					- 0.000\,025\,928 & 0.999\,988\,149 }
					\end{pmatrix}.\begin{pmatrix}{
					X_1\cr
					Y_1 }
					\end{pmatrix}
					$$
					
1.  S'agit-il du modèle bidimensionnel de Helmert? Justifier.

2.  Donner les valeurs numériques respectivement du facteur d'échelle et de l'angle de la rotation entre les deux systèmes.
\eex
\bex
Reprenant les tableaux données ci-dessous, calculer les coordonnées géodésiques des points $(\varphi,\lambda,he)$ sachant que l'ellipsoïde de référence est l'ellipsoïde GRS80 $(a=6\,378\,137.00\,m,e^2=0.0066\,9438\,00229)$.

1. Calculer les coordonnées UTM en prenant le fuseau adéquat.

2. En considérant le modèle de Helmert entre les systèmes $S1$ et $S2$, calculer les paramètres de la transformation.

3. Déterminer les images des points $A,B,C$ et $D$.
\eex
\bpb
 Soient les trois tableaux ci-dessous des coordonnées 3D respectivement dans les systèmes $S1$ et $S2$ et à transformer dans le système $S2$:
\begin{center}
		\[\begin{array}{cccc}
\hline
\hline
 Nom &$X(m)$& $Y(m)$& $Z(m)$\\
\hline
\hline
1& 4\,300\,244.860& 1\,062\,094.681& 4\,574\,775.629\\
2& 4\,277\,737.502& 1\,115\,558.251 &4\,582\,961.996\\
3& 4\,276\,816.431 &1\,081\,197.897 &4\,591\,886.356\\
4 &4\,315\,183.431& 1\,135\,854.241& 4\,542\,857.520\\
5& 4\,285\,934.717& 1\,110\,917.314& 4\,576\,361.689\\
6 &4\,217\,271.349& 1\,193\,915.699& 4\,618\,635.464\\
7& 4\,292\,630.700& 1\,079\,310.256& 4\,579\,117.105\\
 \hline
\hline
\end{array}\]
	\vspace{1cm}
		\[\begin{array}{cccc}
\hline
\hline
 Nom &$X(m)$& $Y(m)$& $Z(m)$\\
\hline
\hline
1 &4\,300\,245.018& 1\,062\,094.592 &4\,574\,775.510\\
2 &4\,277\,737.661& 1\,115\,558.164 &4\,582\,961.878\\
3 &4\,276\,816.590& 1\,081\,197.809& 4\,591\,886.238\\
4 &4\,315\,183.590 &1\,135\,854.153& 4\,542\,857.402\\
5 &4\,285\,934.876 &1\,110\,917.227& 4\,576\,361.571\\
6 &4\,217\,271.512 &1\,193\,915.612& 4\,618\,635.348\\
7 &4\,292\,630.858 &1\,079\,310.168& 4\,579\,116.986\\
\hline
\hline
\end{array}\]
		\[\begin{array}{cccc}
\hline
\hline
 Nom &$X(m)$& $Y(m)$& $Z(m)$\\
\hline
\hline
A &4\,351\,694.594& 1\,056\,274.819&4\,526\,994.706\\
B &4\,319\,956.455& 1\,095\,408.043& 4\,548\,544.867\\
C &4\,303\,467.472 &1\,110\,727.257 &4\,560\,823.460\\
D &4\,202\,413.995 &1\,221\,146.648 &4\,625\,014.614\\
\hline
\hline
\end{array}\]
\end{center}
1. Déterminer les paramètres du modèle de Bur$\breve{s}$a-Wolf à 7 paramètres.

2. Calculer les coordonnées 3D des points du troisième tableau dans le système $S2$. 
\epb

\chapter{\textit{\textbf{Les Systèmes des Altitudes}}}
  Le deuxième type de système géodésique est le référentiel vertical ou le datum vertical qui définit l'origine des altitudes des réseaux de nivellement. D'un pays à un autre, les altitudes sont définies autrement.
\section{\textsc{Les Systèmes d'Altitudes}}\index{Systèmes d'altitudes}
                Pour définir les altitudes, on utilise comme surface de référence, au lieu de l'ellipsoïde de référence géodésique, le géoïde (voir la définition \ref{defgeoid}). On rappelle la définition du géoïde comme suit:
								
'' \textit{Le géoïde correspond à la surface de niveau coïncidant avec le niveau moyen des mers prolongé sous les continents par la condition  d'y rester normal à toutes les lignes de forces}. ''
\\

\begin{figure}
	\centering
			\includegraphics[width=0.90\textwidth]{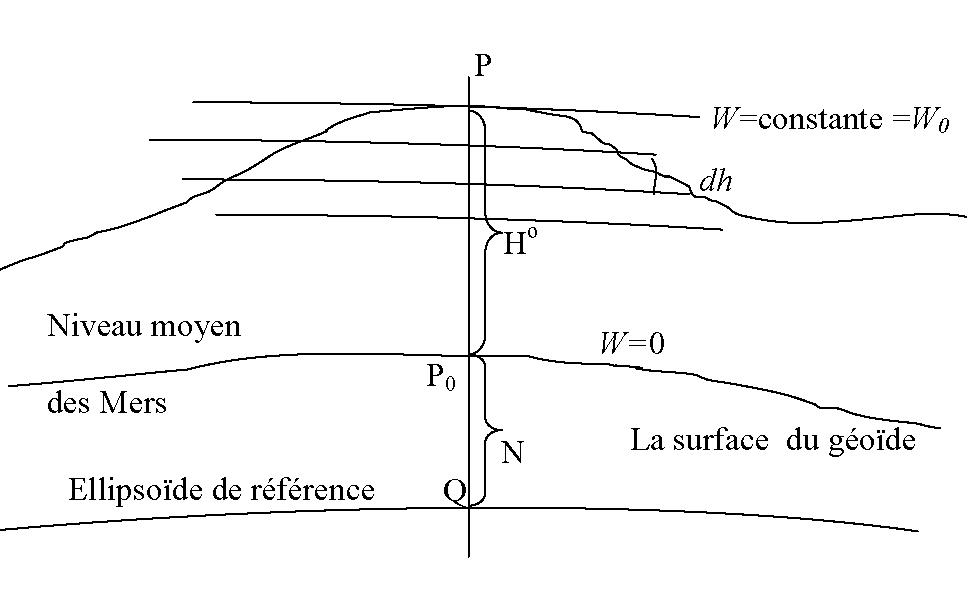}
			\caption{Le géoïde}
	\label{fig:systaltitudes}
\end{figure}
Le géoïde\index{Géoïde} est une surface équipotentielle\index{Surface équipotentielle} $W= Constante =W_o$ qu'on peut prendre égal à 0, soit $W_0=0$. Les surfaces de niveau n'étant pas parallèles et la différence de potentiel est indépendante du chemin suivi pour aller d'un point à un autre, alors on définit la cote géopotentielle\index{Cote équipotentielle} au point $P$ de la surface topographique par :
\be
	                 C(P)  = \int^{P}_{P_0}dW =W(P)-W(P_0)=W(P)-W_0=W(P)      
\ee
Or   $dW = g.dh$, d'où:                                                                                    
\be
	  \fbox{ $   C(P) = \ds \int^{P}_{P_0}gdh $} 
\ee
où $g$ est la gravité mesurée et $dh$ la denivelée mesurée (\textbf{Fig.\ref{fig:systaltitudes}}).
\\
 
Cependant, l'unité de $C(P)$ n'est pas une unité métrique, car une $ gpu = 1 kgal.m = 100\,000\,cm^2/s^2 $, 
alors pour exprimer l'altitude en unité métrique, on définit ci-après les altitudes suivantes.
\subsection*{13.1.1. L'Altitude Orthométrique}\index{Altitude orthométrique}
      L'altitude orthométrique est définie par :
\be
\fbox{ $ 	H^o=\ds \frac{C(P)}{<g>}=\frac{\ds \int^{P}_{P_0}gdh}{<g>} $}
\ee
où :                  
\be
	<g>=\ds \frac{1}{H}\int^{P}_{P_0}gdh=g'
\ee
est la moyenne de $g$ du point $P$ au point $P_0$ sur le géoïde, cette valeur est inaccessible car sa connaissance nécessite une information sur la densité du sous sol. On remplace l'intégrale par une somme finie :
\be
	C(P)=\sum gdh=\sum (g-g'+g')dh=g'\sum dh+\sum (g-g')dh
\ee
D'où :           
\be
	H^o=\frac{\sum gdh}{g'}=\frac{g'\sum dh}{g'}+\frac{\sum (g-g')dh}{g'}=\sum dh+\Delta H
\ee
$\sum dh$ est mesurée par le nivellement de précision et $\Delta H$ est une correction orthométrique donnée par :
\be
	               \fbox{ $ \Delta H = -0.0053sin2\varphi_m H_m \Delta \varphi     $}
	 \ee
avec $\varphi_m$ la latitude moyenne entre le point de départ et du point d'arrivée, $H_m$ l'altitude moyenne et $\Delta \varphi$ la différence des latitudes des points de départ et d'arrivée.     
\subsection*{13.1.2. L'Altitude Normale (ou de Molodensky)}\index{Altitude normale ou de Molodensky}
    Elle est donnée par :
\be
\fbox{ $ 	H^n=\ds \frac{1}{\gamma_m}\int^{P}_{P_0}gdh=\frac{1}{\gamma_m}\sum gdh=\frac{1}{\gamma_m}\sum _i^n g_idh_i $}
\ee
avec:                                                              
\be
	\gamma_m=\ds \frac{1}{H}\int^{H}_{0}\gamma dh
\ee
 où  $\gamma_m$  est la pesanteur normale entre le géoïde et la surface équipotentielle du point considéré. $ \gamma_m $ est  la pesanteur normale à l'ellipsoïde au point de latitude $\varphi$. 
$\gamma_ m$ est estimée par :                                                          
\be
	\gamma_m  = \gamma_0 (1 - H/R) 
\ee
  avec  $\gamma_0 $ la gravité théorique au niveau zéro et à la latitude $\varphi$, donnée par la formule de Cassini (en $gal$):\index{Formule de Cassini}
\be
\fbox{ $ \gamma_0  = 978.0490(1+0.0052884sin^2\varphi  - 0.000059sin^22\varphi )   $} 
   \ee
\subsection*{13.1.3. L'Altitude Dynamique}\index{Altitude dynamique}
  Elle est donnée par :
\be
\fbox{ $	H^d=\ds \frac{1}{\gamma_{0,45}}\int^{P}_{P_0}gdh=\frac{1}{\gamma_{0,45}}\sum _i^n g_idh_i $}
\ee
avec:
	\[	\gamma_{0,45} = \mbox{la valeur de la gravité normale à l'altitude zéro et à la latitude 45°.}
\]
\subsection*{13.1.4. L'Altitude GPS}\index{Altitude GPS}
       Le GPS fournit une altitude ellipsoïdique géométrique $h^{GPS}$. Celle-ci est en relation avec l'altitude orthométrique $H^o$ et l'ondulation du géoïde\index{Ondulation du géoïde} ou hauteur du géoïde $N$ par rapport à l'ellipsoïde du GPS, par l'équation :
\be
	              \fbox{ $   h^{GPS} = H^o + N  $}                            
\ee
\chapter{\textit{\textbf{La Géodésie Tunisienne} }}\label{geotun}

\begin{quotation}
\begin{svgraybox}

		                                  'Geodesists are amongst the men of science those that                   
	                                    operate all over the World, regardless of boundaries, 
	                                                    differences in race, religion, ideology. '
	    
	\end{svgraybox}
                                                                        
\begin{flushright}
( \textbf{Antonio Marussi}\footnote{\textbf{Antonio Marussi} (1908-1984): géodésien et géophysicien italien.}, 1974)\index{\textbf{Marussi A.}}
\end{flushright}
\end{quotation}

\section{\textsc{Introduction}}
            L'un des \'el\'ements fondamentaux de l'infrastructure d'un pays est son r\'eseau g\'eod\'esique. 
\\

\begin{figure}
	\centering
		\includegraphics[width=0.70\textwidth]{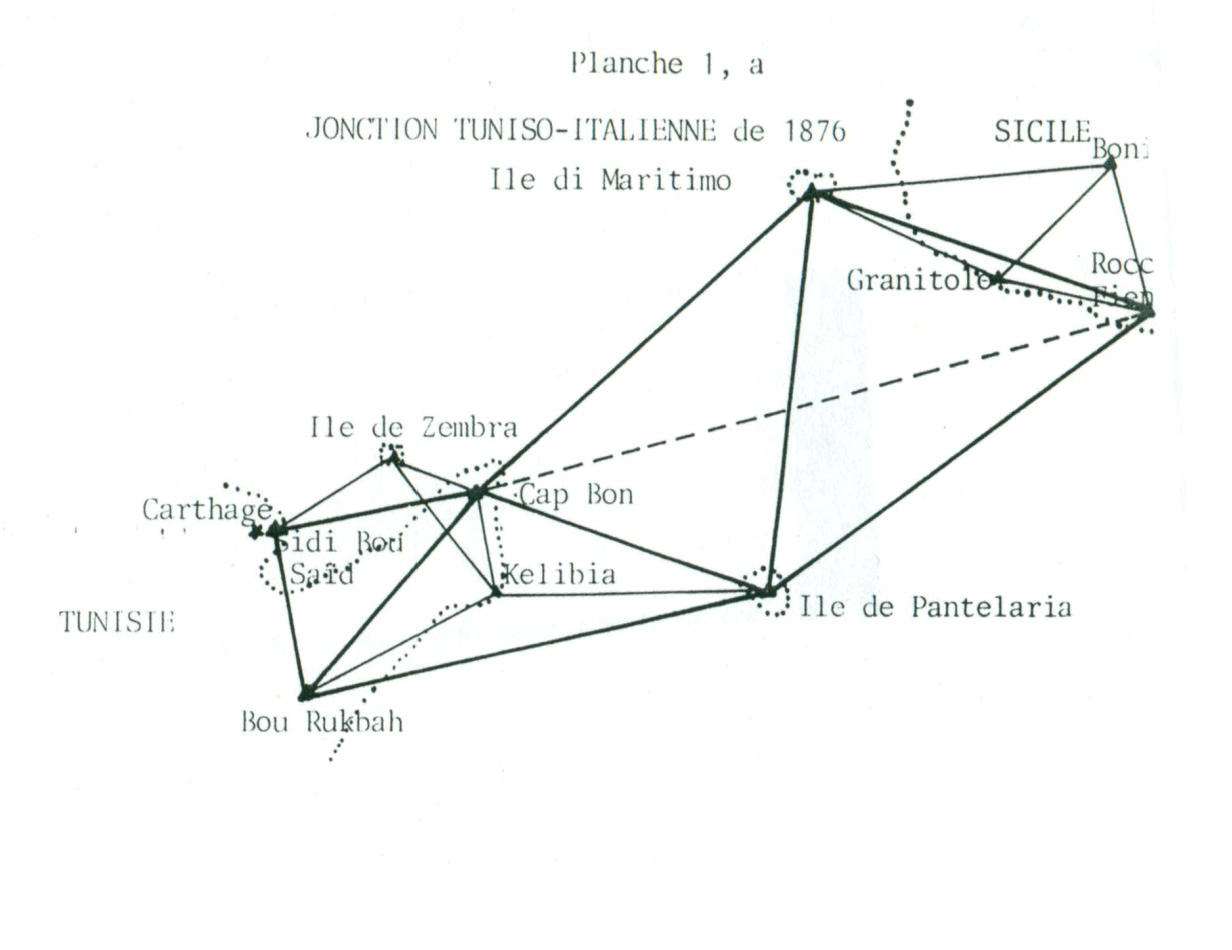}
	\caption{Jonction Tuniso-Italienne de 1876 (\textit{C. Fezzani}, 1979)}\index{\textbf{Fezzani C.}}
	\label{fig:jonction}
\end{figure}

            A ce sujet, la Tunisie a connu le d\'ebut de la cr\'eation de son premier r\'eseau g\'eod\'esique \`a partir des observations g\'eod\'esiques de la liaison entre le Cap-Bon et l'\^ile de Sicile en 1876 (\textbf{Fig. \ref{fig:jonction}}).
						\\

Sur ce r\'eseau vont \^etre rattach\'es tous les travaux cartographiques et topographiques et particuli\`erement ceux de l'Immatriculation Fonci\`ere Facultative (IFF) et du Cadastre.
\\

             Cent ans apr\`es le d\'ebut des observations du premier r\'eseau g\'eod\'esique en Tunisie, il s'\'etait av\'er\'e la n\'ecessit\'e de revaloriser et moderniser ce r\'eseau g\'eod\'esique afin de satisfaire les besoins cartographiques et topographiques du pays avec le d\'eveloppement des appareils de mesures. 
\\

             Ainsi, \`a partir de 1978, l'Office de la Topographie et de la Cartographie (OTC) a d\'ecid\'e de mener une campagne astro-g\'eod\'esique qui va durer jusqu'\`a 1983.
\\

              Le nouveau r\'eseau g\'eod\'esique obtenu fait appara\^itre un d\'ecalage entre 0 \`a 10 m\`etres avec l'ancien r\'eseau (Carthage34)\index{Réseau Carthage34}. L'adoption du nouveau r\'eseau g\'eod\'esique n'a pas \'et\'e faite.
\\

              Avec le d\'eveloppement du positionnement g\'eod\'esique avec les satellites GPS dans le syst\`eme de r\'ef\'erence mondial $WGS84$ (World Geodetic System 84), un nouveau r\'eseau de r\'ef\'erence \`a vocation spatiale vient d'\^etre observ\'e en Tunisie en 1998. Les r\'eseaux g\'eod\'esiques classiques Tunisiens peuvent -ils r\'epondre \`a la technologie GPS ? et de quelles mani\`eres ?
\\

                L'objet de ce chapitre est  de d\'ecrire l'\'evolution des systèmes g\'eod\'esiques terrestres tunisiens vers leur unification ainsi que l'adoption d'une nouvelle repr\'esentation plane ad\'equate pour toute la Tunisie. 
\\

              Dans la suite, on s'int\'eresse aux r\'eseaux g\'eod\'esiques planim\'etriques.

\section{\textsc{Historique}}
              Les premiers points g\'eod\'esiques  ont \'et\'e cr\'e\'es en 1876 par les g\'eod\'esiens italiens lors de la liaison entre le Cap-Bon et l'\^ile de Sicile (\textbf{Fig. \ref{fig:Jonction}}). Le premier r\'eseau g\'eod\'esique tunisien a \'et\'e d\'efini \`a partir de l'extension du r\'eseau Alg\'erien, par l'observation du point astronomique Carthage en 1878 et la mesure de l'azimut astronomique de la direction Carthage - Bir Bou Regba pour l'orientation du r\'eseau (\textit{C. Fezzani}, 1979).\index{\textbf{Fezzani C.}}  
         \\
            
\begin{figure}
	\centering
		\includegraphics{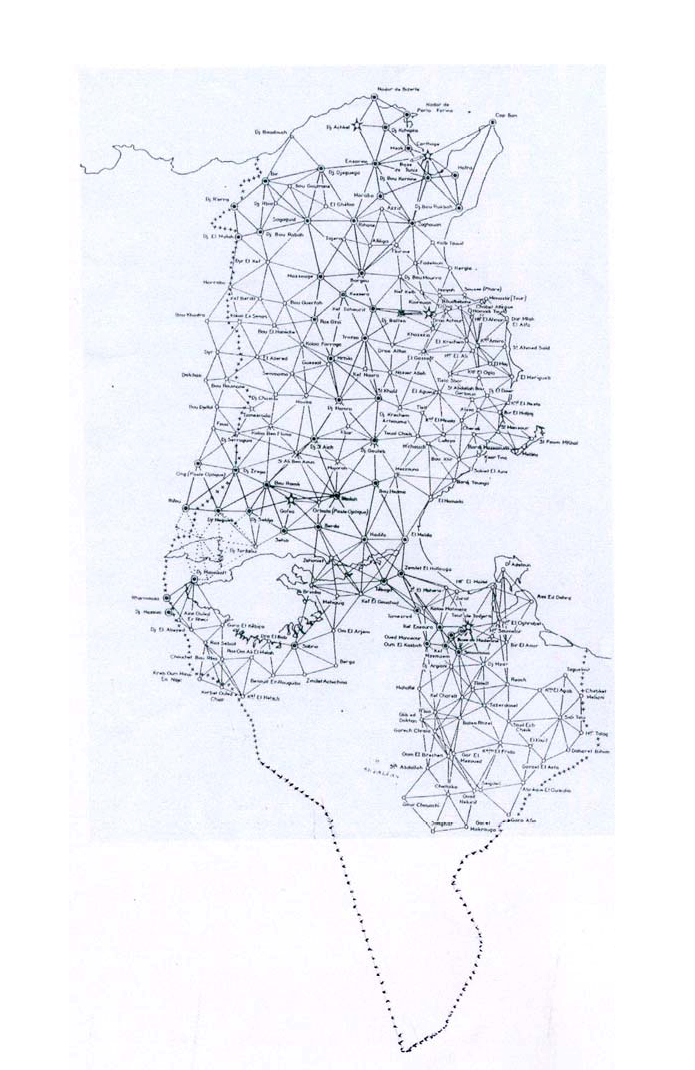}
	\caption{La structure des réseaux géodésiques avant 1978 }
	\label{reseauavant1978}
\end{figure}

 Avant 1978, la structure g\'eod\'esique tunisienne  \'etait comme suit (\textbf{Fig. \ref{reseauavant1978}}):
 
-	un r\'eseau g\'eod\'esique dit du 1er ordre form\'e par le parall\`ele de Tunis au nord et le parall\`ele de Gab\`es au sud reli\'es par la m\'eridienne dite de Gab\`es, les longueurs des c\^ot\'es varient de 30 \`a $50\, kms$;

-	un r\'eseau g\'eod\'esique du 1er ordre compl\'ementaire;

-	un r\'eseau g\'eod\'esique du 2\`eme ordre;

-	un r\'eseau g\'eod\'esique du 2\`eme ordre compl\'ementaire;

-	les r\'eseaux g\'eod\'esiques de d\'etail du 3\`eme et 4\`eme ordre;

-	un canevas de points astronomiques au sud (Sahara).
            \\
            
A un r\'eseau de points g\'eod\'esiques est associ\'e le syst\`eme g\'eod\'esique c'est-\`a-dire les \'el\'ements de r\'ef\'erence \`a partir desquels sont calcul\'ees les coordonn\'ees g\'eographiques ( latitude, longitude ) ou les coordonn\'ees planim\'etriques $(X,Y)$. On donne ci-dessous les syst\`emes les plus utilis\'es en Tunisie.

\section{\textsc{Les Syst\`emes G\'eod\'esiques en Tunisie}}
\subsection*{14.3.1. Le Syst\`eme G\'eod\'esique 'Voirol'}\index{Système Voirol}
   C'\'etait le premier syst\`eme en Tunisie caract\'eris\'e par :
             
-	le point fondamental (point de d\'epart) : Voirol pr\`es d'Alger cr\'e\'e en 1875;

-	la surface  de r\'ef\'erence c'est-\`a-dire le mod\`ele choisi pour la Terre est l'ellipsoïde de Clarke Fran\c{c}ais 1880;\index{Ellipsoïde de Clarke 1880 Français}

-	l'orientation de d\'epart est l'azimut astronomique de la direction Voirol-Meleb El Kora mesur\'e en 1874;

-	la mise \`a l'\'echelle ou la qualit\'e m\'etrique de r\'eseau : la mesure d'une distance ou base \`a Blida en Alg\'erie mesur\'ee en 1854. 
\\

Une grande partie du premier r\'eseau g\'eod\'esique terrestre tunisien \'etait calcul\'e dans ce syst\`eme.

\subsection*{14.3.2. Le Syst\`eme G\'eod\'esique 'Carthage34'}

            A la suite de la d\'etection d'une erreur dans la mise \`a l'\'echelle du syst\`eme Voirol en 1910 et vu sa qualit\'e, le Service G\'eographique de l'Arm\'ee Fran\c{c}aise (S.G.A.F) a \'etabli un nouveau syst\`eme g\'eod\'esique ind\'ependant du syst\`eme  Voirol. Les \'el\'ements de d\'efinition de ce  syst\`eme sont : 
            
-	le point fondamental : le point Carthage en Tunisie;
 	
-	l'ellipsoïde de r\'ef\'erence : l'ellipsoïde de Clarke Fran\c{c}ais 1880;
 	
-	l'azimut de l'orientation : l'azimut astronomique de la direction Carthage - Bir Bou Regba;
 	
-	la mise \`a l'\'echelle : les bases de Tunis et de Medenine.
   \\
            
  Les calculs des coordonnées des points g\'eod\'esiques de la partie nord ont \'et\'e achev\'es en 1934.
	\\
	
	Bien que l'ellipsoïde de référence est le même, les deux systèmes géodésiques Voirol et Carthage34 ont des origines différentes ce qui explique les décalages en longitude et latitude géodésiques. Ainsi au point fondamental Carthage, par exemple, on a les différences suivantes :
\ba
             \varphi_{Voirol}   -  \varphi_{Carthage} = 25.86\, dmgr \\
               \lambda_{Voirol}  -  \lambda_{Carthage} = 36.19\, dmgr
\ea
qui se traduisent par un décalage moyen planimétrique de $245\, m$  en $x(Nord)$  et de $ 280\, m$ en $y(Ouest)$.

\section{\textsc{Les  Repr\'esentations  Planes} }

                A un syst\`eme g\'eod\'esique donn\'e, on peut lui associer plusieurs types de repr\'esentations planes. On donne ci-dessous les repr\'esentations planes en usage en Tunisie.

\subsection*{14.4.1. La Repr\'esentation de Bonne}\index{Représentation de Bonne}

          La repr\'esentation de Bonne\footnote{\textbf{Rigobert Bonne} (1727-1795): ingénieur, mathématicien et cartographe français.} est une représentation \'equivalente (conserve les surfaces). Elle n'est plus en usage mais  elle \'etait  utilis\'ee  pour  le  d\'ecoupage   cartographique  des  cartes aux \'echelles 1/50 000, 1/100 000 et 1/200 000 (version ancienne).\index{\textbf{Bonne R.}}

\subsection*{14.4.2. La Repr\'esentation des Fuseaux }\index{Représentation Fuseaux }

           Elle a \'et\'e utilis\'ee dans le syst\`eme g\'eod\'esique Voirol pour le besoin de la triangulation et reste en usage dans les travaux de l'immatriculation foncière facultative pour traiter certains anciens dossiers (\textit{C. Fezzani}, 1979).
\\

\begin{figure}
	\center
		\includegraphics[width=0.80\textwidth]{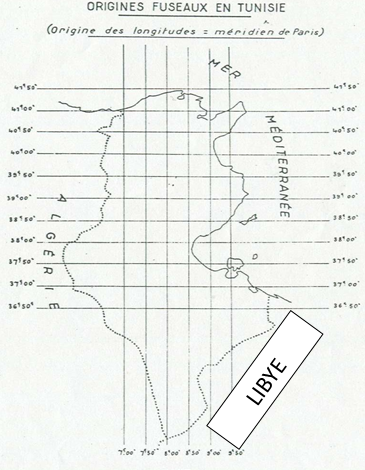}
	\caption{Le découpage des fuseaux (\textit{C. Fezzani}, 1979)}\index{\textbf{Fezzani C.}}
	\label{fuseau1}
\end{figure}

Dans cette repr\'esentation, la Tunisie était partag\'ee en six fuseaux, d'une \'etendue chacun de 0.5 grades $(gr)$ en longitude, subdivis\'es chacun en onze quadrilat\`eres curvilignes de $0.5\, gr$ de c\^ot\'e en latitude (\textbf{Fig. \ref{fuseau1}}). Cette repr\'esentation plane fut abandonn\'ee en 1922 pour \^etre remplac\'ee par la repr\'esentation plane Lambert Tunisie.\index{\textbf{Lambert J.H.}}
\\

Les formules des coordonnées Fuseaux $(x(Nord),y(Ouest))$ en un point de coordonnées $(\varphi,\lambda)$, avec $(\varphi_0,\lambda_0)$ les coordonnées du centre du quadrilatère considéré, sont comme suit (\textit{A. Ben Hadj Salem}, 2013):
\be
\fbox{$ \begin{array}{l}
\ds y=\frac{(\lambda-\lambda_0)"cos\varphi}{R} \\
\\
x=\ds \frac{(\varphi-\varphi_0)"}{P}+\frac{Q}{P}y^2
\end{array} $}
\ee
avec les coefficients:
\be
	\fbox{$ \ds P=\frac{2\times 10^6}{\pi \rho_0}, \quad Q=\ds \frac{1}{2}.\frac{tg\varphi_0}{N_0}.P, \quad \ds R=\frac{2\times 10^6}{\pi N_0}  $}
	\ee
où $(\varphi-\varphi_0),(\lambda-\lambda_0)$ sont exprimés en $dmgr(")$ en comptant les longitudes positives à l'Ouest de Greenwich et $N_0,\, \rho_0$ les 2 rayons de courbure de l'ellipsoïde de Clarke Français 1880 pour $\varphi=\varphi_0$.
\newpage
\subsection*{14.4.3. La Repr\'esentation Lambert}\index{Représentation Lambert}

             C'est une repr\'esentation conforme (conserve les angles) d'un mod\`ele ellipsoïdique. Afin d'\'eviter les d\'eformations trop importantes, la repr\'esentation Lambert Nord Tunisie a \'et\'e adopt\'ee pour la partie Nord du pays (latitude comprise entre $37.5\, gr$ et $42.5\, gr)$ et la repr\'esentation Lambert Sud Tunisie a \'et\'e adopt\'ee pour la partie sud (latitude comprise entre $34.5\, gr$ et $39.5\, gr)$. La repr\'esentation  Lambert Tunisie est nomm\'ee \`a l'OTC sous l'appellation "Origine Unique ". 
\\

\begin{figure}[h]
	\centering
		\includegraphics[width=0.60\textwidth]{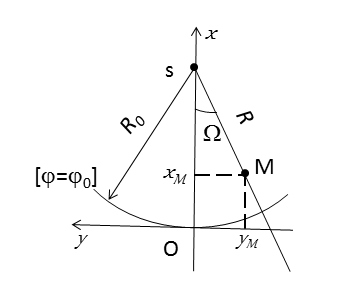}
	\caption{Le repère STT}
	\label{fig:repstt}
\end{figure}

Pour la Tunisie, on considérait un système d'axes $(O,x,y)$ tel que  l'axe $Oy$ est la tangente à l'image du parallèle origine $\varphi_0$ au point $O$ dirigé vers l'Ouest et $Ox$ est porté par l'image du méridien origine dirigé vers le Nord (\textbf{Fig. \ref{fig:repstt}}), dit repère STT\footnote{STT: Service Topographique Tunisien, premier service tunisien chargé des travaux topographiques et foncières fut créé en 1882, deviendra plus tard l'OTC.}. Soit le point $S$ de $Ox$ avec $OS=R_0$, on a alors les coordonnées planes $(x,y)$:
\ba
	x_M=R_0-Rcos\Omega \nonumber \\
		y_M=-Rsin\Omega \nonumber
\ea
avec $\lambda$ comptée positivement à l'Est du méridien origine des longitudes.
\\

\begin{figure}[htp]
	\centering
		\includegraphics[width=0.60\textwidth]{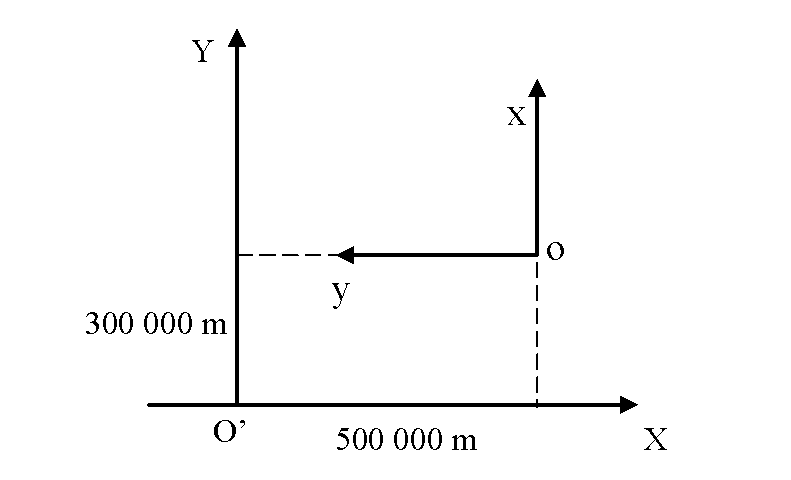}
	\caption{Le repère $(O',X,Y)$}
	\label{fig:reperexy3a}
\end{figure}

Avec l'introduction du facteur de réduction de l'échelle, les formules des coordonnées rectangulaires $(x,y)$ deviennent :
\be 
\fbox{ $ \begin{array}{l}
x_M=k(R_0 -Rcos( (\lambda-\lambda_0)sin\varphi_0)) \\
	y_M=-kRsin((\lambda-\lambda_0)sin\varphi_0) \\
		\\
		avec\quad \Omega=(\lambda - \lambda_0)sin\varphi_0
		\end{array} $}
\ee
Pour obtenir des coordonnées rectangulaires positives, on définit pour la Tunisie un repère $(O',X,Y)$ tel que  $O'X$ et $O'Y$ soient dirigés respectivement vers l'Est et le Nord (\textbf{Fig. \ref{fig:reperexy3a}}).
\\

Dans le nouveau repère $(O',X,Y)$, on a donc:
\be
\fbox{ $ \begin{array}{l}
X=500\,000.00\,m\,-y\\
Y=300\,000.00\,m\,+x
\end{array} $}
\ee
Les quantités $500\,000.00\,m\,$ et $300\,000.00\,m\,$ sont respectivement les constantes de translation en $X(Est)$ et en $Y (Nord)$.
\\

Due à l'importance de l'altération linéaire, la Tunisie est partagée en deux zones avec les caractéristiques suivantes:
\begin{center}
\[\begin{array}{cccc}
\hline
\hline
   \mbox{Représentation Plane} &\mbox{Méridien origine}&\mbox{Parallèle origine}& \mbox{k}\\ \hline \hline
      \mbox{Lambert Nord Tunisie}&  \lambda_0=11\,gr       & \varphi_0=40\, gr        & 0.999\,625\,544     \\   \hline
       \mbox{Lambert Sud Tunisie}&  \lambda_0=11\,gr       & \varphi_0=37\, gr        & 0.999\,625\,769     \\   \hline \hline
     \end{array}
\]
\end{center}

\begin{table}[h]
	\centering
	\caption{Caractéristiques de la représentation Lambert Tunisie}
	\label{tab: Caractéristiques de la Représentation Lambert Tunisie}
\end{table}
avec une amplitude de $ \pm 2.5 \,gr$ du parallèle origine.

\subsection*{14.4.4. La Repr\'esentation Plane U.T.M. }\index{Représentation UTM}

              C'est une repr\'esentation conforme. Le mod\`ele terrestre est un ellipsoïde divis\'e en 60 fuseaux de 6° d'amplitude, chacun de 3° de longitude de part et d'autre du m\'eridien central. Dans la repr\'esentation UTM, toute la Tunisie se trouve dans le fuseau 32 avec le m\'eridien de longitude 9° à l'Est de Greenwich comme m\'eridien central.
         \\
            
  Cette repr\'esentation a \'et\'e utilis\'ee dans le syst\`eme g\'eod\'esique Europe50 (usage militaire) et dans les nouvelles cartes aux \'echelles 1/200 000 et  1/50 000 \'edit\'ees par l'OTC.
	
\section{\textsc{Pourquoi Un Nouveau Syst\`eme G\'eod\'esique?}}
Suite \`a la publication du d\'ecret-loi relatif \`a l'mmatriculation fonci\`ere obligatoire en 1964 et \`a la promulgation du Code des droits r\'eels en 1965 d'une part, et le d\'ebut de la mission cartographique du Service Topographique Tunisien d'autre part, il \'etait imp\'eratif de mettre en place les infrastructures n\'ecessaires pour r\'ealiser les travaux de cartographie et de topographie dans un nouveau syst\`eme g\'eod\'esique national bas\'e sur une nouvelle g\'eod\'esie.
\\

A cet effet, un premier diagnostic a \'et\'e effectu\'e en 1969 sur l'\'etat de la g\'eod\'esie tunisienne (\textit{A. Fontaine}, 1969). 
\\

   L'analyse de l'\'etat de ces r\'eseaux a montr\'e des insuffisances aux niveaux de la qualit\'e de l'\'echelle et de l'orientation.
\\

\section{\textsc{Les Travaux  de  Modernisation des  R\'eseaux G\'eod\'esiques Tunisiens}}

              Grâce \`a la prise de conscience \`a la Direction de la Topographie et de la Cartographie (DTC)\footnote{En 1970, le Service Topographique Tunisien (STT) devenait la Direction de la Topographie et de la Cartographie au Ministère des travaux publics.} de l'importance des sciences g\'eographiques et en particulier de l'aspect g\'eod\'esique, un protocole d'accord a été conclu entre la DTC et l'Institut G\'eographique National de France (IGNF) en 1972.
\\

              Ce protocole concernait l'\'etude et l'analyse des calculs de compensation des r\'eseaux g\'eod\'esiques tunisiens du 1er et 2\`eme ordre. L'analyse de l'\'etat de ces r\'eseaux (\textit{C. Fezzani}, 1979)\index{\textbf{Fezzani C.}} a montr\'e des insuffisances aux niveaux de la qualit\'e de l'\'echelle (1/40\,000 \`a 1/30\,000) et de l'orientation $(15\,dmgr$ \`a $25\,dmgr)$. De plus, de nombreux points g\'eod\'esiques ont disparu et d'autres ont \'et\'e d\'etruits d'o\`u la n\'ecessit\'e de reprendre des travaux g\'eod\'esiques pour revaloriser les r\'eseaux g\'eod\'esiques tunisiens.

\subsection*{14.6.1. Les Travaux de la Revalorisation de la Géodésie Tunisienne   }

               A partir de 1978, l'OTC a d\'ecid\'e des travaux pour moderniser les r\'eseaux g\'eod\'esiques tunisiens afin de satisfaire les besoins cartographiques et topographiques du pays en commen\c{c}ant par le r\'eseau g\'eod\'esique de base.          
\\

Ces travaux de revalorisation des  r\'eseaux g\'eod\'esiques tunisiens ont pour objectifs:
               
   -	l'abolition d\'efinitive de l'utilisation des diff\'erents syst\`emes g\'eod\'esiques terrestres de types " isol\'es et fuseaux " en usage depuis plus d'un si\`ecle;

   -	la cr\'eation d'un r\'ef\'erentiel g\'eod\'esique terrestre unique pour la Tunisie;

   -	la mise en place d'une nouvelle repr\'esentation  plane qui convient le mieux pour la Tunisie \`a savoir l'UTM.
\\

Les travaux de revalorisation de la g\'eod\'esie Tunisienne (\textit{M. Charfi}, 1983) comprenaient :\index{\textbf{Charfi M.}}
                   
-	la r\'efection des anciens points du 1er ordre, du 1er ordre compl\'ementaire, du 2\`eme ordre et du 2\`eme ordre compl\'ementaire;

-	la construction de nouveaux points sur les sites des anciens points disparus;

-	la densification de l'ancien r\'eseau par de nouveaux points;

-	les observations angulaires azimutales et z\'enithales;

-	la d\'etermination de 8 points de Laplace;

-	la mesure des c\^ot\'es de 8 triangles g\'eod\'esiques;

-	la d\'etermination de 5 points par la m\'ethode Doppler;\index{La méthode Doppler}

-	la compensation des observations terrestres avec les donn\'ees Doppler pour obtenir les nouvelles coordonn\'ees du nouveau r\'eseau.
\\

                Les observations des 8 points de Laplace et la mesure des c\^ot\'es des 8 triangles g\'eod\'esiques, les observations et le calcul des 5 stations Doppler ainsi que la compensation du r\'eseau g\'eod\'esique ont fait l'objet d'une convention sign\'ee entre l'OTC et l'IGNF en 1982, suite à un appel d'offres international.                         
\\

                Le nouveau r\'eseau g\'eod\'esique terrestre appel\'e R\'eseau G\'eod\'esique Primordial (RGP)\index{Réseau Géodésique Primordial} est compos\'e de 312 points comme suit :
                
*	143 points anciens;

*	112 nouveaux points construits sur les sites des anciens points disparus;

*	58 nouveaux points.  

\subsection*{14.6.2. La Compensation du R\'eseau G\'eod\'esique Primordial }

                La compensation du RGP effectu\'ee par l'IGNF comprenait les compensations planim\'etrique et altim\'etrique.
\\

La compensation planim\'etrique  de 1984 effectu\'ee par l'IGNF a défini donc un nouveau système géodésique terrestre nommé OTC84\index{Système OTC84}. En comparant les coordonn\'ees issues d'OTC84 et avec celles de Carthage34, on a trouv\'e que les coordonn\'ees anciennes ont subi un d\'eplacement sous la forme d'une rotation dans le sens des gisements dont le centre se trouve dans la r\'egion de  J. Semmama et d'un angle de $20\, dmgr$ (2/1000 de grade). Les d\'eplacements planim\'etriques varient de 0 \`a $10.70\, m$ en s'\'eloignant du centre de la rotation.
\\


\subsection*{14.6.3. Le Système Carthage86}\index{Système Carthage86}
Faute d'adopter les calculs de 1984, l'OTC a effectué par ses propres moyens un calcul de compensation des nouvelles observations angulaires en fixant les coordonnées Carthage34 des points anciens existants et ce-ci en trois phases:

- la phase une : la zone du Nord;

- la phase deux : la zone du Centre;

- la phase trois: la zone du Sud. 

Cette compensation a donné naisssance au système géodésique terrestre Carthage86, ayant le même ellipsoïde Clarke Français 1880. Les coordonnées des points géodésiques obtenues dans ce système sont à $0.80 \, m$ près de celles de Carthage34. Il faut signaler que les observations astronomiques de 1980 n'ont pas été introduites dans la compensation ce qui montre que l'orientation de Carthage86 est similaire à celle de Carthage34.

   \begin{figure}[tp]
	\centering
		\includegraphics[width=0.70\textwidth]{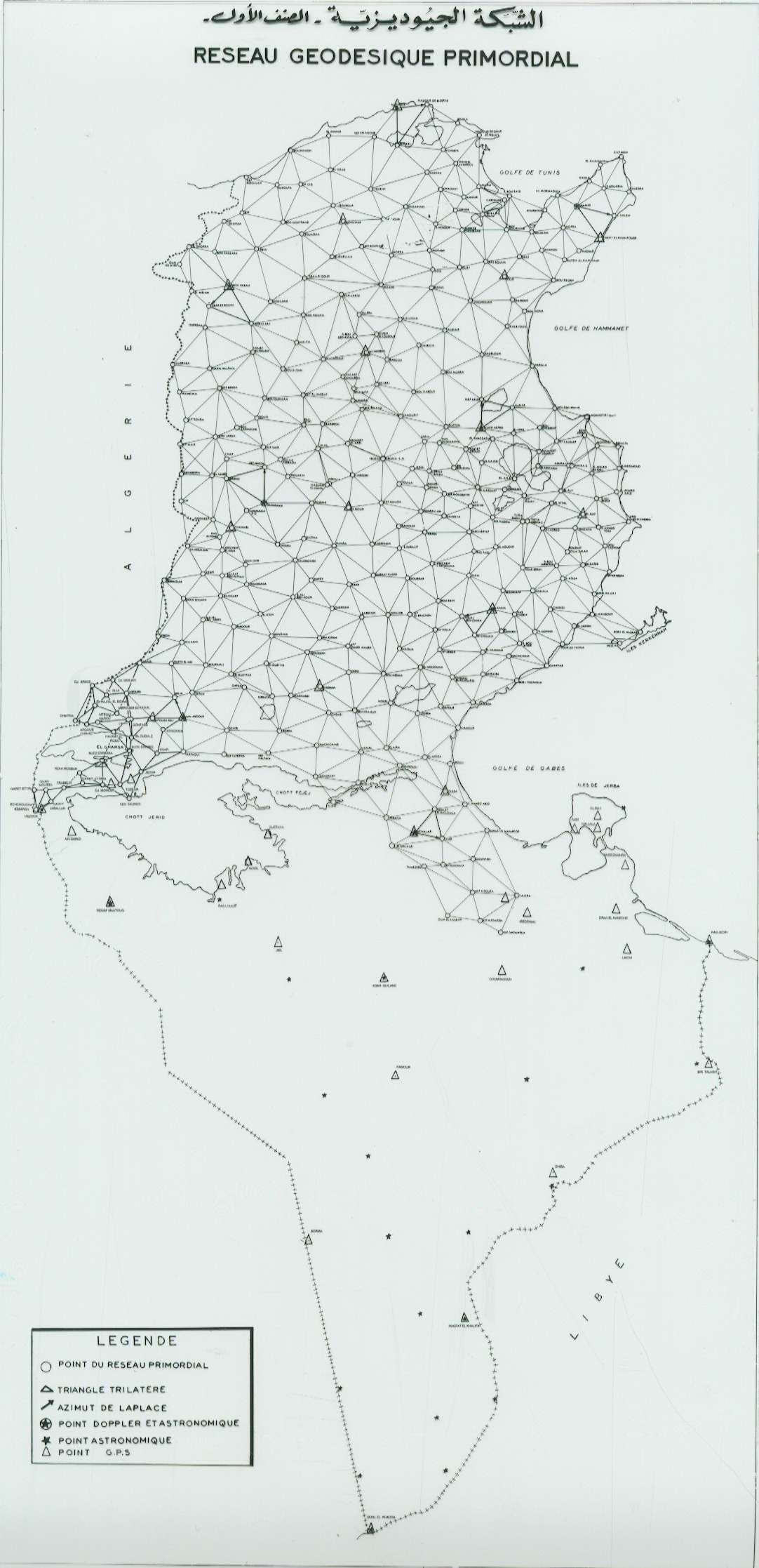}
	\caption{Le Réseau Géodésique Primordial Tunisien}
	\label{fig:geo}
\end{figure}

\subsection*{14.6.4. Les R\'esultats des Travaux de la Revalorisation de la G\'eod\'esie Tunisienne}
A l'issu des travaux de calculs en 1984, on a obtenu un nouveau syst\`eme g\'eod\'esique - dit syst\`eme g\'eod\'esique OTC84 - mat\'erialis\'e par le nouveau r\'eseau g\'eod\'esique appel\'e le R\'eseau Primordial Terrestre Tunisien constitu\'e de 312 points g\'eod\'esiques (\textbf{Fig. \ref{geo}}).   

\subsubsection*{14.6.4.1. Le Syst\`eme G\'eod\'esique  OTC84}\index{Système OTC84}
Les caract\'eristiques du syst\`eme sont:

-	ellipsoïde de r\'ef\'erence: ellipsoïde de Clarke 1880 Fran\c{c}ais;\index{Ellipsoïde de Clarke Français}

-	les nouvelles observations angulaires + les anciennes observations angulaires pour les anciens points conserv\'es  (avant 1978);

- les observations de 8 points de Laplace (latitude, longitude et azimut astronomiques) et de 24 distances dans 8 triangles;

- la fixation de 5 points anciens (observ\'es par Doppler) avec un \'ecart-type de 0.50 $m$;\index{La méthode Doppler}

-	une compensation globale par la m\'ethode des moindres carr\'es (312 points).
\\

\textbf{Le d\'ecalage entre OTC84 et Cathage34}: Il est sous la forme d'une rotation dont le centre est situ\'e dans la r\'egion de Kasserine, et d'un angle de 25 $dmgr$ dans le sens des gisements.
\\

Ce ph\'enom\`ene a \'et\'e observ\'e dans la plupart des pays qui ont chang\'e de syst\`eme g\'eod\'esique.
Voici ci-dessous un exemple de la Suisse o\`u on montre les d\'ecalages entre le syst\`eme ancien (LV03) et le nouveau (LV95). 
(\textit{Publication de l'Acad\'emie Bavaroise de G\'eod\'esie}, 1997)
\begin{figure}
	\centering
		\includegraphics[width=1.00\textwidth]{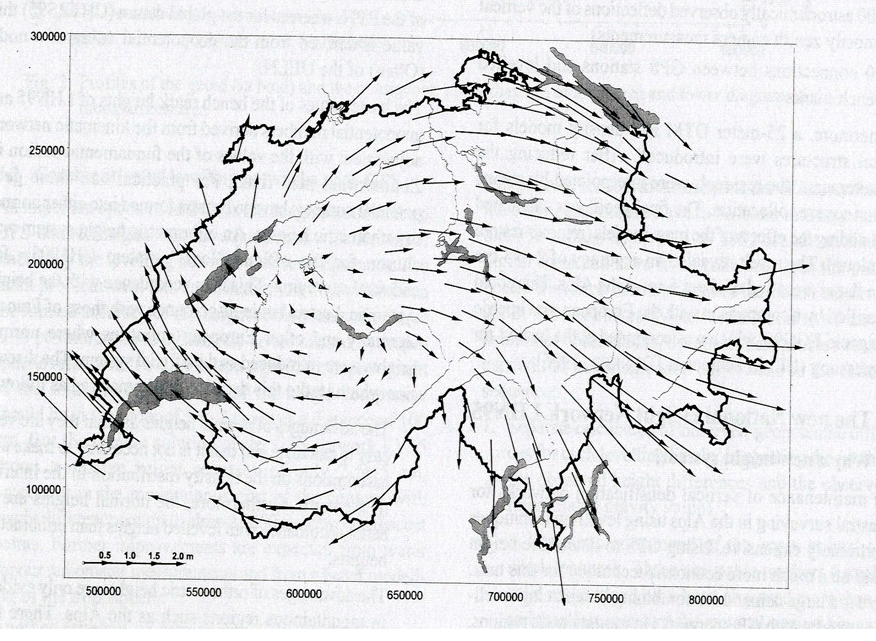}
	\caption{les d\'ecalages entre le syst\`eme ancien (LV03) et le nouveau (LV95) }
	\label{fig:lv95}
\end{figure}
\\

Le d\'ecalage entre OTC84 et Carthage34 a \'et\'e jug\'e inacceptable pour le patrimoine national en mati\`ere de cadastre. 

R\'esultat: Ce  nouveau syst\`eme géodésique n'a pas \'et\'e appliqu\'e.
\\

D'o\`u l'id\'ee de la reprise des calculs des observations.




\section{\textsc{La Mise \`a niveau de la G\'eod\'esie Tunisienne}}
En 2001, l'OTC a entrepris un programme de mise \`a niveau de la g\'eod\'esie Tunisienne (\textit{A. Ben Hadj Salem}, 1999). 

   Cette mise \`a niveau comprenait:
   
    - 1. l'unification des syst\`emes g\'eod\'esiques terrestres;
    
    - 2. la cr\'eation d'un r\'eseau géodésique de base \`a r\'ef\'erence  spatiale;
    
    - 3. la d\'etermination d'un g\'eoide pr\'ecis pour la Tunisie;
    
    - 4. la  mise  en place d'une  nouvelle  repr\'esentation cartographique pour la Tunisie.
  \\
   
    L'objectif de cette mise \`a niveau est:
    
  * d'unifier les diff\'erents syst\`emes en un seul r\'ef\'erentiel g\'eod\'esique terrestre homog\`ene et pr\'ecis, qui permettra d'une part la densification du r\'eseau et l'exploitation des nouvelles techniques spatiales de positionnement dans les meilleures conditions et d'autre part de lever et de r\'etablir les limites born\'ees sans se ref\'erer obligatoirement aux titres riverains. 
\\

Les imp\'eratifs de la mise \`a niveau \'etant:

 -	la non alt\'eration des formes et des surfaces des parcelles;
 
 -	la pr\'eservation des orientations des directions dans la limite des tol\'erances requises. 
\\

   Une commission technique a \'et\'e cr\'e\'ee en d\'ecembre 2001 pour  la mise en oeuvre de cette mise \`a niveau.
\\

\begin{figure}[htp]
	\centering
		\includegraphics[width=.700\textwidth]{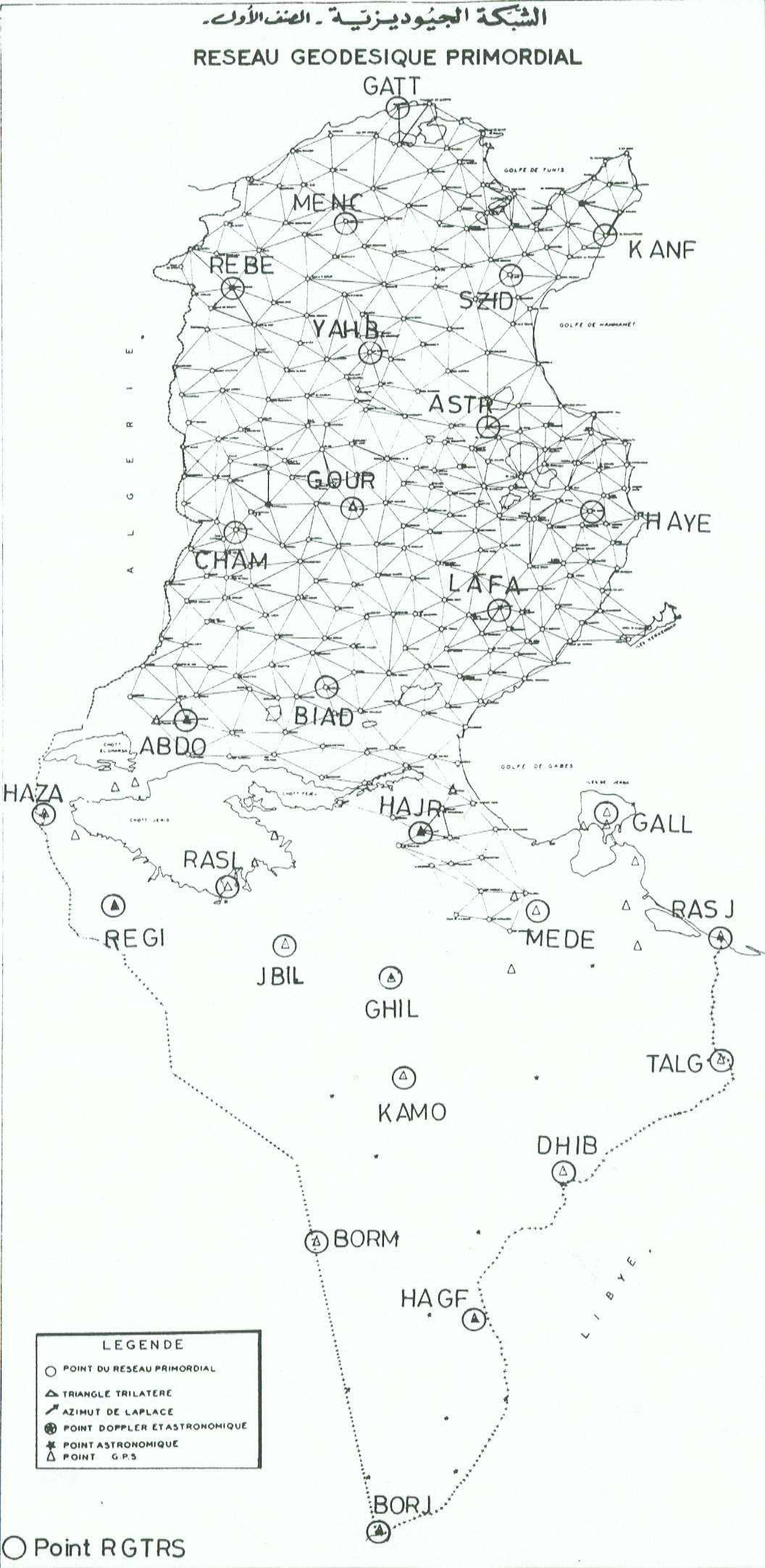}
	\caption{Le Réseau GPS Tunisien de Référence Spatiale (RGTRS)}
	\label{fig:geo1}
\end{figure}

La r\'ealisation de cette mise \`a niveau a commenc\'e par l'\'etablissement du R\'eseau G\'eod\'esique GPS Tunisien d'Ordre Z\'ero constitu\'e de 28 points r\'epartis sur tout le territoire (\textbf{Fig. \ref{fig:geo1}}).
\\
   
    A partir de ce r\'eseau, la d\'esorientation du syst\`eme Carthage34 a \'et\'e v\'erifi\'ee. 
\\

    Disposant d'un logiciel de calculs g\'eod\'esiques, la commission a proc\'ed\'e \`a un nouveau calcul des observations.

\subsection*{14.7.1. Comparaisons des coordonn\'ees des syst\`eme Carthage34 et OTC84}
La comparaison des coordonn\'ees Lambert Nord Tunisie de 48 points dans les deux syst\`emes a permis de mod\'eliser le d\'ecalage entre le syst\`eme Carthage34 et le nouveau calcul par une transformation math\'ematique conforme du type :
\be
	Z - z_0 = (z-z_0).e^{it}   \label{dq1}
\ee
avec :

   - " $z_0$ " le centre de la rotation : point fictif pr\`es du point g\'eod\'esique Gassaat Ej Jahfa, situ\'e dans la r\'egion de Kasserine;
   
   - " $t$ " l'angle de rotation = -27 $dmgr$. 
   
\subsection*{14.7.2. Les Conclusions de la Commission Technique}
Apr\`es \'etudes et tests des coordonn\'ees issues du nouveau calcul, la Commission technique  a abouti aux r\'esultats suivants:

   - le syst\`eme OTC84 laisse invariant les formes, les directions  et les surfaces dans les  tol\'erances requises;
   
   - l'homog\'en\'eit\'e de ce syst\`eme g\'eod\'esique  offre :
   
     * un meilleur passage au syst\`eme spatial $WGS84$;
     
     * une parfaite int\'egration des futurs travaux  g\'eod\'esiques et topom\'etriques dans le r\'ef\'erentiel  spatial.
\\
     
La  Commission technique a propos\'e d'adopter le système OTC84 le nouveau syst\`eme g\'eod\'esique terrestre tunisien appel\'e \textbf{NTT} (\textbf{N}ouvelle \textbf{T}riangulation \textbf{T}unisienne).\index{Système géodésique NTT}

\section{\textsc{L'Arr\^et\'e du 10 F\'evrier 2009}}\index{L'Arrêté du 10 février 2009}
L'arr\^et\'e du ministre de la D\'efense nationale du 10 f\'evrier 2009, paru dans le Journal Officiel de la R\'epublique Tunisienne n° 14 du 17 f\'evrier 2009, fixe:

   1.   le syst\`eme national de r\'ef\'erence unifi\'e de la g\'eod\'esie;

   2.   le syst\`eme national de r\'ef\'erence de la représentation cartographique;

   3.   le syst\`eme national de r\'ef\'erence du nivellement.
   
\subsection*{14.8.1. Le Syst\`eme National de R\'ef\'erence Unifi\'e de la G\'eod\'esie}
   Il est d\'efini par :
   
   - le syst\`eme national g\'eod\'esique des coordonn\'ees g\'eographiques  appel\'e - la Nouvelle Triangulation Tunisienne - (NTT);
   
   - l'ellipsoïde associ\'e: c'est l'ellipsoïde de Clarke 1880 Fran\c{c}ais ($a= 6\,378\,249.20$ m, $b= 6\,356\,515.00$ m).\index{Ellipsoïde de Clarke Français}

\subsection*{14.8.2. Le Syst\`eme National de R\'ef\'erence de la Représentation Cartographique}

Le syst\`eme national actuel de r\'ef\'erence de la représentation cartographique est d\'efini par la représentation plane Universal Transverse Mercator (UTM), fuseau 32 Nord.\index{Représentation UTM}
\\

Les caract\'eristiques de la repr\'esentation UTM :
 
\textbf{D\'efinition:} C'est une repr\'esentation:
 
- conforme c'est-\`a-dire conserve les angles;

- cylindrique ====> on utilise les coordonn\'ees rectangulaires $(X,Y)$;

- transverse ====>  $X = X(\varphi,\lambda )$ et $Y = Y(\varphi,\lambda )$

d'un mod\`ele ellipsoïdique.
\\

* Les \'el\'ements de d\'efinition sont:

- l'ellipsoïde de référence est celui de Clarke Fran\c{c}ais 1880;

- le m\'eridien origine: 9° \`a l'Est de Greenwich  ou fuseau n°32;

- le facteur d'\'echelle : $k = 0.9996$;

- la constante en $X: 500\,000.00\, m$;

- la constante en $Y : 0.00\, m$.

\subsection*{14.8.3. Le Syst\`eme National de R\'ef\'erence du Nivellement}
On utilise le syst\`eme des altitudes orthom\'etriques.\index{Altitude orthométrique}
Le point fondamental  ou r\'ef\'erence des altitude est le rep\`ere scell\'e sur le monument Porte de France \`a Tunis avec une altitude de $7.000\, m$ au dessus du niveau moyen de la mer (Nouvelle compensation du R\'eseau de Nivellement G\'en\'eral de la Tunisie - 1961). Cette référence a été déterminée à partir de mesures du niveau moyen de la mer faites au marégraphe\index{Marégraphe} de la Goulette. Le système d'altitudes est appelé le Nivellement Général de la Tunisie (NGT).\index{Nivellement général de la Tunisie}  
 \\
      
    Le Nivellement Général de la Tunisie comprend :
    
*	un réseau de nivellement de précision du 1er ordre;

*	un réseau de nivellement de précision du 2ème ordre;

*	un réseau de nivellement de précision du 3ème ordre.

\section{\textsc{Conclusions}}
L'arr\^et\'e du 10 f\'evrier 2009 a d\'efini les fondements de la g\'eomatique en Tunisie, car un pays sans les d\'efinitions pr\'ecises de ses r\'ef\'erentiels de g\'eod\'esie, de nivellement et de cartographie ne peut avancer en la mati\`ere.
\\

L'unification des syst\`emes g\'eod\'esiques terrestres va permettre: 

- d'abolir l'utilisation des syst\`emes isol\'es;

- le d\'eveloppement des bases de donn\'ees et les syst\`emes d'informations g\'eographiques sur l'ensemble du pays;

- d'assainir le syst\`eme foncier tunisien sur des bases solides;

- de d\'evelopper en harmonie l'application des technologies actuelles de positionnement spatial et celles du futur.








\chapter{\textit{\textbf{Notions sur le Mouvement d'un Satellite Artificiel de la Terre}}}
Avant de passer au chapitre sur le système GPS, il est utile pour le lecteur d'avoir des notions sur le mouvement d'un satellite artificiel autour de la Terre.
\section{\textsc{Les Equations du Mouvement}}
On considère un satellite de masse $m$ dont la position est définie par le vecteur $\textbf{\textit{OS}}=\textbf{\textit{r}}$. La Terre est considérée comme une masse ponctuelle de masse $m'$ située au point $O$ centre de la Terre ($m'=(5\,973\pm1)\times10^{21}\,kg$). 
\\

L'équation du mouvement du satellite est donnée par:
\begin{equation}
	\fbox{ $ \ds m\frac{d^2\textbf{\textit{r}}}{dt^2}=\textbf{\textit{F}}=-G\frac{mm'}{r^3}\textbf{\textit{r}} $} \label{sa18}
\end{equation}
 avec $\textbf{\textit{F}}$ la force d'attraction gravitationnelle et $G$ est la constante universelle de la gravitation de valeur égale à $(6673\pm 1)\times10^{-14}m^3s^{-2}kg^{-1}$ (\textit{H. Moritz \& I.I. Mueller}, 1988). \index{Constante universelle de la gravitation}\index{\textbf{Moritz H.}}\index{\textbf{Mueller I.I.}}
 
 L'équation (\ref{sa18}) s'écrit aussi:
\be
 m\frac{d^2\textbf{\textit{r}}}{dt^2}=\textbf{\textit{F}}=\textbf{\textit{grad}}V
\ee
On dit que $\textbf{\textit{F}}$ dérive du potentiel $V$ avec:
\be
V=G\frac{mm'}{r }
\ee
On pose:
\begin{equation}
	\fbox{$\mu=Gm' =(3986005\pm 0.5)\times10^8m^3s^{-2}$} \label{sa19}
\end{equation}
L'équation (\ref{sa18}) s'écrit:
\be
	\ddot{\textbf{\textit{r}}}=-\mu\frac{\textbf{\textit{r}}}{r^3} \label{sa20}
\ee
Comme:
$$ \textbf{\textit{r}}=\left\{
 \begin{array}{lll}
	X_C \\ Y_C \\ Z_C
\end{array}\right. $$
alors l'équation vectorielle (\ref{sa20}) s'écrit en trois équations différentielles du deuxième ordre comme suit:
\be
\fbox{ $ \begin{array}{l}
	\ds \ddot{X}_C+\frac{\mu}{r^3}X_C=0   \\
\ds	\ddot{Y}_C+\frac{\mu}{r^3}Y_C=0   \\
\ds	\ddot{Z}_C+\frac{\mu}{r^3}Z_C=0 
	\end{array} $}\label{sa22}
\ee
Après l'intégration des équations  (\ref{sa22}), on obtient six paramètres des conditions initiales qui définissent la forme et la position de l'orbite et une constante donnant la variation du mouvement du satellite avec le temps.

Les équations (\ref{sa22}) montrent qu'on a un mouvement d'un corps dans un champ central\index{Champ central}.
\subsection*{15.1.1. La 2ème Loi de Kepler}\index{Deuxième Loi de Kepler}\index{\textbf{Kepler J.}}
1. Si on applique \textit{le théorème du moment cinétique}, on obtient:
\be
	 \frac{d \mathbold{\sigma}}{dt}=\textbf{\textit{0}} \label{sa25}
\ee
car:
\begin{equation}
	\mathbold{\sigma}=\textbf{\textit{OS}}\wedge \textbf{\textit{v}}=\textbf{\textit{C}}=\mbox{constante} \label{sa26}
\end{equation}
En effet, en dérivant (\ref{sa26}) par rapport au temps, on a alors:
\begin{equation}
\frac{d\mathbold{\sigma}}{dt}=\frac{d\textbf{\textit{OS}}}{dt}\wedge \textbf{\textit{v}}+\textbf{\textit{OS}} \wedge \frac{d\textbf{\textit{v}}}{dt}=\textbf{\textit{v}}\wedge \textbf{\textit{v}}+\textbf{\textit{OS}}\wedge \ddot{\textbf{\textit{r}}}=\textbf{\textit{r}}\wedge \ddot{\textbf{\textit{r}}} \label{sa27}
\end{equation}
Or d'après (\ref{sa20}):
	\[\ddot{\textbf{\textit{r}}}=-\frac{\mu}{r^3}\textbf{\textit{r}}
\]
D'où:
\be
\frac{d\mathbold{\sigma}}{dt}=-\frac{\mu}{r^3}\textbf{\textit{r}} \wedge \textbf{\textit{r}}=\textbf{\textit{0}} \label{sa28}
\ee
On déduit donc (\ref{sa26}) et on a:
\be
	C=||\textbf{\textit{C}}||=\mbox{constante des aires} \label{sa29}
\ee
De (\ref{sa26}), le mouvement se fait dans un plan, en définissant \index{Constante des aires} $\textbf{\textit{OS}}=\textbf{\textit{r}}=\textbf{\textit{r}}(r,\upsilon)$, alors les composantes de $\textbf{\textit{v}}$ vecteur vitesse sur le rayon vecteur $\textbf{\textit{r}}$ et de la direction perpendiculaire sont:
\be
	\textbf{\textit{v}}=\left\{
\begin{array}{lll}
\displaystyle	\frac{dr}{dt} \\   \\ \displaystyle  r\frac{d\upsilon}{dt}
\end{array}\right. \label{sa30}
\ee
On a alors:
\be
	C=r.r\frac{d \upsilon}{dt}=r^2\frac{d\upsilon}{dt} \label{sa31}
\ee
et aussi:
\be 
	\dot{\Sigma}=\frac{d\Sigma}{dt}=\frac{r^2d\upsilon}{2dt} \label{sa32}
\ee
avec $\Sigma$ la surface balayée par le vecteur position. 

Des équations (\ref{sa31}) et (\ref{sa32}), on a \textit{la 2ème loi de Kepler}:
\be
	\fbox{ $ \dot{\Sigma}=\ds \frac{1}{2}C=\mbox{constante} $} \label{sa33}
\ee
\bthm
(\textbf{Deuxième loi de Kepler}, ou \textbf{loi des aires}) L'aire balayée
par le vecteur position $\textbf{\textit{r}}(t)$ varie linéairement avec le temps.
\ethm
2. Si on applique  \textit{le théorème de l'énergie cinétique} sous forme différentielle, on obtient:
\be
	d(\frac{1}{2}mv^2)=\textbf{\textit{F}}.\textbf{\textit{v}}dt \label{sa34}
\ee
où:
\be
\textbf{\textit{F}}.\textbf{\textit{v}}=-\mu \frac{m}{r^3}\textbf{\textit{r}}\frac{d\textbf{\textit{r}}}{dt}=\frac{-\mu m}{r^2}\frac{dr}{dt}=\frac{d}{dt}\left(\frac{\mu m}{r} \right) \label{sa35}
\ee
En remplaçant (\ref{sa35}) dans le second membre de (\ref{sa34}) et en intégrant, on arrive à :
$$ 	\frac{1}{2}mv^2-\frac{\mu m}{r}=\mbox{constante} $$
Soit:
\begin{equation}
	\fbox{ $ \ds 	\frac{1}{2}mv^2-\frac{\mu m}{r}=H =\mbox{constante} $} \label{sa37}
\end{equation}
où $H$ est \textit{la constante de l'énergie} ou \textit{l'énergie}.\index{Constante d'énergie}
\blm
Les fonctions $H$ et $C$ sont constantes le long des solutions: on dit que ce sont des intégrales premières du mouvement, c'est-à-dire que $H$ et $C$ sont des fonctions de la position $r$ et de la dérivée première de $r$ par rapport au temps ($t$), constantes au cours du temps. 
\elm
En effet, des équations (\ref{sa29})  et (\ref{sa37}), on a respectivement: $$H=H(r,v)=\mbox{constante}$$
et: $$C=C(r,v)=\displaystyle r^2\frac{d\upsilon}{dt}=\mbox{constante}$$
On note que $H=T-L$ s'appelle aussi le hamiltonien\footnote{En hommage à \textbf{Sir William Rowan Hamilton} (1805-1865): mathématicien, physicien et astronome irlandais.} \index{Hamiltonien} du mouvement keplérien. En effet, $T=\ds \frac{1}{2}mv^2$ est l'énergie cinétique et $U$ est le potentiel $\ds \frac{\mu m}{r}$ et on retrouve l'expression de $H$ donnée par l'équation (\ref{sa37}).
\subsection*{15.1.2. La 1ère loi de Kepler}\index{Première Loi de Kepler}\index{\textbf{Hamilton W.R.}}
En multipliant vectoriellement à droite les membres de l'équation (\ref{sa20}) par $\textbf{\textit{C}}=\textbf{\textit{r}}\wedge \textbf{\textit{v}}$, on obtient:
\ba
&	\ddot{\textbf{\textit{r}}}\wedge\textbf{\textit{C}}=\ds -\frac{\mu}{r^3}\textbf{\textit{r}}\wedge (\textbf{\textit{r}}\wedge \textbf{\textit{v}})=\ds -\frac{\mu}{r^3}\left[ \textbf{\textit{r}}\wedge (\textbf{\textit{r}} \wedge \textbf{\textit{v}})\right]=-\frac{\mu}{r^3}\left[\textbf{\textit{r}}.(\textbf{\textit{r}}.\textbf{\textit{v}})-(\textbf{\textit{r}}.\textbf{\textit{r}})\textbf{\textit{v}}\right] \nonumber& \\& =\ds -\frac{\mu}{r^3}\left(\textbf{\textit{r}}r\frac{dr}{dt}-r^2\textbf{\textit{v}}\right)=-\frac{\mu}{r^2}\left(\textbf{\textit{r}}\frac{dr}{dt}-r\textbf{\textit{v}}\right) \nonumber&\\&
\ds \ddot{\textbf{\textit{r}}}\wedge \textbf{\textit{C}}=\mu\frac{d}{dt}\left(\frac{\textbf{\textit{r}}}{r}\right)	\label{sa38}
\ea
Comme $\textbf{\textit{C}}$ est constant, l'équation (\ref{sa38}) s'écrit:
$$ 	\frac{d}{dt}\left(\dot{\textbf{\textit{r}}}\wedge \textbf{\textit{C}}\right)-\mu\frac{d}{dt}\left(\frac{\textbf{\textit{r}}}{r}\right)=0 $$
Après intégration, on a:
$$ 	\dot{\textbf{\textit{r}}}\wedge \textbf{\textit{C}}-\mu \frac{\textbf{\textit{r}}}{r}=\textbf{\textit{l}}=\,\mbox{vecteur constant}\,\,\textbf{\textit{l}}_0 $$
ou:
\be
	(\textbf{\textit{v}} \wedge \textbf{\textit{C}})-\mu \frac{\textbf{\textit{r}}}{r}=\textbf{\textit{l}}  \label{sa40}
\ee
On appelle $\textbf{\textit{l}}$ vecteur de Laplace\index{Vecteur de Laplace}. On multiplie l'équation (\ref{sa40}) par $\textbf{\textit{r}}$, on obtient :
\be
\textbf{\textit{r}}.(\textbf{\textit{v}} \wedge \textbf{\textit{C}})-\mu \frac{\textbf{\textit{r}}.\textbf{\textit{r}}}{r}=\textbf{\textit{r}}.\textbf{\textit{l}}  \label{sa41}
\ee
Comme:
$$ \textbf{\textit{r}}.(\textbf{\textit{v}}\wedge \textbf{\textit{C}})=\textbf{\textit{v}}.(\textbf{\textit{C}}\wedge \textbf{\textit{r}})=\textbf{\textit{C}}.(\textbf{\textit{r}}\wedge \textbf{\textit{v}})=\textbf{\textit{C}}.\textbf{\textit{C}}=C^2$$
et:
$$\textbf{\textit{l}}.\textbf{\textit{r}}=lrcos\upsilon $$
(\ref{sa41}) devient:
\be
	C^2=\mu r+lrcos\upsilon \label{sa42}
\ee
En posant:
\be
	p=\frac{C^2}{\mu}\,\,\,\,\,\mbox{et}\,\,\,\,\,e=\frac{l}{\mu} \label{sa43}
\ee
on déduit de (\ref{sa42}):
\be
\fbox{ $	r=\ds \frac{p}{1+ecos\upsilon} $} \label{sa44}
\ee
d'où:
\bthm 
(\textbf{1ère loi de Kepler}) L'orbite décrite par le vecteur de position $\textit{\textbf{r}}$ est une conique (ellipse) définie en coordonnées polaires $(r,\upsilon)$. L'angle $\upsilon$ compté, entre la direction du vecteur de Laplace $\textbf{\textit{l}}$ ou $\textbf{\textit{OP}}$ (périgée)\index{Périgée} et le rayon vecteur $\textbf{\textit{r}}$, s'appelle \textit{l'anomalie vraie}.\index{Anomalie vraie}\index{Vecteur de Laplace}
\ethm
On a :
\ba
	\mbox{Pour}\,\,\upsilon=0\Longrightarrow r_1=\frac{p}{1+e}\,\,\,\,\mbox{c'est la périgée}  \nonumber \\
	\mbox{Pour}\,\,\upsilon=\pi\Longrightarrow r_2=\frac{p}{1-e}\,\,\,\,\,\mbox{c'est l'apogée} \nonumber  \\
	\mbox{D'où:}\quad  r_1+r_2=2a=\frac{2p}{1-e^2}\,\,\Longrightarrow p=a(1-e^2) \nonumber
	\ea
	Par suite:
	\be
	\begin{array}{l}
	 r_1=\ds \frac{p}{1+e}=\frac{a(1-e^2)}{1+e}=a(1-e)\\
	 r_2=\ds \frac{p}{1-e}=\frac{a(1-e^2)}{1-e}=a(1+e)
	\end{array} \label{sa46b} 
	\ee
\subsection*{15.1.3. La 3ème Loi de Kepler}\index{Troisième Loi de Kepler}
D'après la 2ème loi de Kepler donnée par (\ref{sa33}):
$$ \dot{\Sigma}=\frac{d\Sigma}{dt}=\frac{C}{2}=\mbox{constante}$$ d'où:
\be
	d\Sigma=\frac{C}{2}dt \label{sa48}
\ee
En intégrant (\ref{sa48}) sur une période, on obtient:
$$ \frac{1}{2}\int_0^T Cdt=\int d\Sigma=\Sigma=\pi.a.b=\pi.a^2\sqrt{1-e^2} $$
Soit:
\be
	C=\frac{2\pi}{T}a^2\sqrt{1-e^2} \label{sa50}
\ee
Comme $C=\sqrt{p \mu}=\sqrt{a(1-e^2)\mu}\,$ et $T$ la période, on a finalement:
\be
	\fbox{ $ \ds \frac{a^3}{T^2}=\frac{\mu}{4\pi^2}=\mbox{constante} $}\label{sa51}
\ee
C'est la \textit{3ème loi de Kepler}.
\bthm
(\textbf{3ème loi de Kepler}) Le carré de la période est proportionnel au cube du demi-grand axe de l'ellipse.
\ethm
\section{\textsc{Eléments de l'orbite}}
Après l'intégration des équations du mouvement du satellite artificiel, on obtient six paramètres qui définissent la position du plan de l'orbite, ses dimensions, appelés les éléments d'orbite \index{Eléments d'orbite} et ce sont:
\begin{figure}
	\centering
		\includegraphics[width=0.90\textwidth]{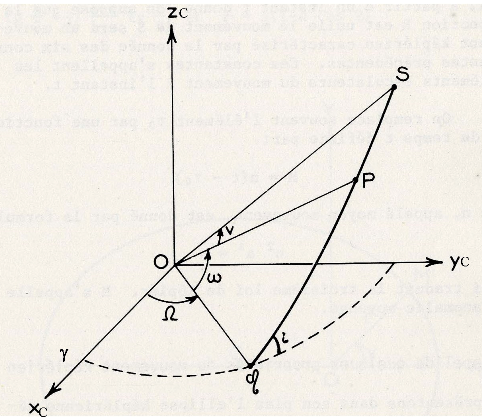}
	\caption{Le repère céleste}
	\label{fig:sat2a}
\end{figure}

- $a$ - le demi-grand axe;

- $e$ - la première excentricité;

- $i$ - l'angle d'inclinaison;

- $\Omega$ - l'ascension droite du noeud ascendant;

- $\omega$ - l'argment;

- $t_0$ - l'instant de passage au périgée. 

\subsection*{15.2.1. Les Coordonnées}
En conséquence de la 3ème loi de Kepler (\ref{sa51}), on peut écrire:
\be
	\fbox{ $ n=\ds \frac{2\pi}{T}=\sqrt{\frac{\mu}{a^3}} $} \label{sa52}
\ee
$n$ est appelé \textit{vitesse moyenne angulaire}\index{Vitesse moyenne angulaire}. A partir de (\ref{sa52}), on définit \textit{l'anomalie moyenne} $M$ \index{Anomalie moyenne} à l'instant $t$ par:
\begin{equation}
	\fbox{$M=n(t-t_0)$} \label{sa53}
\end{equation}
A l'aide de la figure (Fig. \ref{fig:sat3}), on va exprimer les coordonnées du satellite dans le plan de l'orbite:
\begin{figure}
	\centering
		\includegraphics[width=1.00\textwidth]{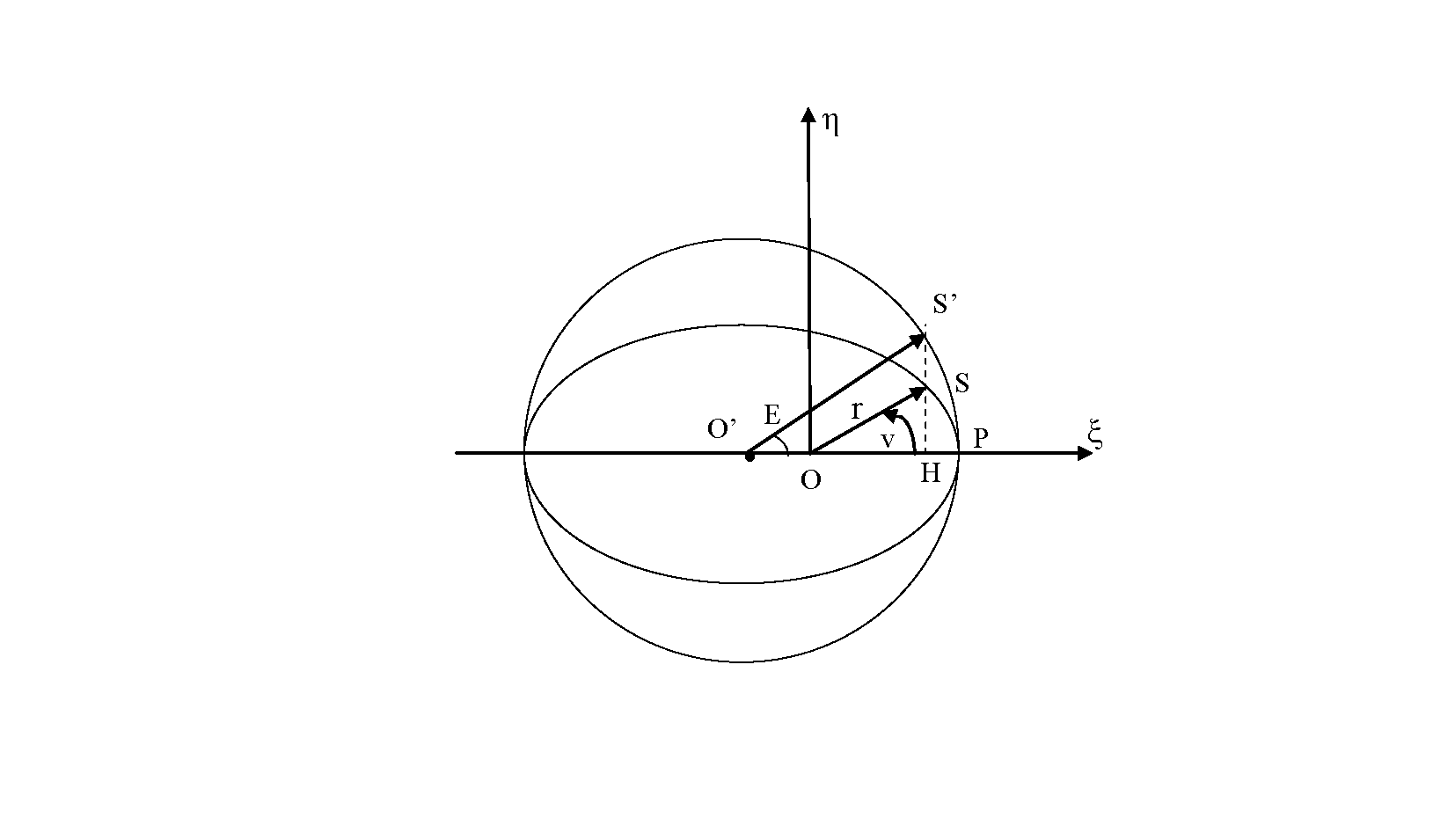}
	\caption{Plan de l'orbite}
	\label{fig:sat3}
\end{figure}

où:

- $O$ est le centre de gravité de la terre et aussi un foyer de l'ellipse;

- $S$ la position du satellite;

- $\upsilon$ l'anomalie vraie;\index{Anomalie vraie}

- $E$ l'angle $PO'S'$ est appelé l'anomalie excentrique;\index{Anomalie excentrique}

- l'axe $O\eta$ est perpendiculaire à l'axe $O\xi$ et l'axe $T\zeta$ est perpendiculaire au plan $O\xi \eta$.
\\

On sait que dans le repère $O\xi \eta$, on a:
\be 
\fbox{ $  \begin{array}{l}
	\textbf{\textit{OS}}\left\{
	\begin{array}{l}
	\xi=OScos\upsilon=rcos\upsilon \\ 
	\eta=OSsin\upsilon=rsin\upsilon
  \end{array}\right. \\
\mbox{avec}\,\,\,r=\ds \frac{p}{1+ecos\upsilon} 
\end{array} $} \label{sa55} 
\ee
Maintenant, d'après la première loi de Kepler, l'aire balayée par le vecteur position $\textbf{\textit{OS}}$ entre les instants $t_0$ et $t$ vaut:
\ba
	\Sigma=\int_{t_0}^{t}Cdt=\frac{C(t-t_0)}{2}=\frac{(t-t_0)2\pi a^2\sqrt{1-e^2}}{2T}=  \nonumber  \\ \frac{n(t-t_0)a^2\sqrt{1-e^2}}{2}=\frac{Ma^2\sqrt{1-e^2}}{2} \label{sa56}
\ea
 Comme l'ellipse de paramètres $(a,b)$ est obtenue par affinité de rapport $k=b/a=\sqrt{1-e^2}$ du cercle centré en $O'$ et de rayon $a$. 
Donc:
 $$ k=\ds \frac{\Sigma}{\sigma_1} $$
où la surface $\sigma_1$ est celle du triangle curviligne $OPS'$, elle est égale à la différence du secteur circulaire $O'PS'$ et du triangle $O'OS'$ soit:
 $$  A_1=\mbox{aire secteur O'PS'}=\ds \frac{\pi .a^2}{2\pi}.E=\frac{a^2.E}{2} $$
 et l'aire du triangle $O'OS'$ vaut:
$$	A_2=\ds \frac{O'O.HS'}{2} $$
Comme:
$$ sinE=\frac{HS'}{a}\Rightarrow \,HS'=a.sinE, \quad O'O=a-r_1=a-a(1-e)=a.e$$
d'où:
$$ \sigma_1=A_1-A_2=\ds \frac{a^2E}{2}-\frac{O'O.HS'}{2}=	\frac{a^2E}{2}-\frac{ae.asinE}{2}=\frac{a^2(E-esinE)}{2}$$
On peut écrire alors en utilisant (\ref{sa56}) que:
\ba
	\Sigma=k\sigma_1 \Longrightarrow \frac{Ma^2\sqrt{1-e^2}}{2}=\frac{b}{a}\frac{a^2(E-esinE)}{2}=\nonumber \\ \sqrt{1-e^2}\frac{a^2(E-esinE)}{2} \Longrightarrow	 E-esinE=M =n(t-t_0)\nonumber
\ea
L'équation:
\be
\fbox{$\displaystyle E-esinE=M=n(t-t_0)$} \label{sa58}
\ee
s'appelle \textit{l'équation de Kepler}.\index{Equation de Kepler}
\\

Cette relation est importante, puisqu'elle permet de calculer $E$ en fonction du temps et par suite de déterminer $\upsilon=\upsilon(t)$ voir l'équation (\ref{saaa63}) ci-dessous, et $r=r(t)$.

On peut calculer la valeur de l'anomalie excentrique $E$ par la méthode itérative. A la première itération, on prend: $$ E_1=M+esinM $$ et: $$E_2=E_1+\delta E$$ En utilisant (\ref{sa58}), on a:
$$ 	E_1+\delta E-esin(E_1+\delta E)=M $$
En faisant un développement au premier degré, on obtient:
	$$E_1+\delta E-esinE_1cos\delta E-esin\delta E cosE_1=M$$
Comme $\delta E$ est petit, on a $cos \delta E\approx1$ et $sin \delta E\approx \delta E$, on arrive à:
\be           
	\delta E=\frac{M - E_1+esinE_1}{1-ecosE_1} \label{sa61}
\ee         
On prend maintenant: $$E_1=E_1+\delta E$$ et on applique (\ref{sa61}) et ainsi de suite jusqu'à ce que $\delta E$ soit négligeable devant la précision désirée.
\\

Dans le repère $O \xi \eta \zeta$, on peut écrire les coordonnées du satellite sous la forme:
\be
\fbox{ $ \begin{array}{l}
	\xi=OH=O'H-O'O=acosE-ae=a(cosE-e)  \\
	\eta=SH= (b/a)HS'=\sqrt{1-e^2}asinE=a\sqrt{1-e^2}sinE  \\
	\zeta=0 
	\end{array} $} \label{sa64}
\ee
Ce qui donne:
$$ \ds \frac{\eta}{\xi}=\frac{\sqrt{1-e^2}sinE}{cosE-e} $$
Or d'après (\ref{sa55}), on a:
$$ \left\{
	\begin{array}{ll}
	\xi=rcos\upsilon \\
	 \eta=rsin\upsilon
\end{array}\right.\Longrightarrow tg\upsilon=\frac{\eta}{\xi} $$
d'où:
\be
\fbox{ $ tg\upsilon=\ds \frac{\sqrt{1-e^2}sinE}{cosE-e} $} \label{saaa63}
\ee
On exprime maintenant les coordonnées $(X_C,Y_C,Z_C)$ du satellite dans le référentiel céleste $X_CY_CZ_C$ à l'aide de la figure (Fig.\ref{fig:sat2a}). Il est nécessaire de faire successivement:
\begin{enumerate}
\item une rotation de $-\omega$ autour de l'axe $O\xi$;

\item une rotation de $-i$ autour de l'axe $O\Omega$;

\item une rotation de $-\Omega$ autour de l'axe $O'Z_C$.
\end{enumerate}
Les matrices de rotations sont les suivantes:
$$ R(-\Omega)=\begin{pmatrix} {
	cos\Omega & -sin\Omega & 0 \cr
	sin\Omega & cos\Omega & 0 \cr
	0 & 0 & 1 }
\end{pmatrix} $$
$$  R(-i)=\begin{pmatrix} {
	1         & 0          & 0 \cr
	0         & cosi      &-sini \cr
	0 & sini & cosi }
\end{pmatrix} $$
et:
$$ R(-\omega)=\begin{pmatrix} {
	cos\omega & -sin\omega & 0 \cr
	sin\omega & cos\omega & 0 \cr
	0 & 0 & 1 }
\end{pmatrix} $$
D'où:
\be
	\begin{pmatrix} {
	X_C  \cr
	 Y_C \cr
	  Z_C }
	\end{pmatrix}=R(-\Omega).R(-i).R(-\omega)\begin{pmatrix}{
	\xi    \cr \eta   \cr \zeta } 
	\end{pmatrix} \label{sa68}
\ee
Ce qui donne après calculs:
$$ \begin{pmatrix}{
	X_C         \cr Y_C           \cr Z_C }
\end{pmatrix}=\begin{pmatrix}{
	cos\Omega cos\omega-sin\Omega sin\omega cosi & -cos\Omega sin\omega-sin\Omega cos\omega cosi & sin\Omega sin\omega \cr
	sin\Omega cos\omega+cos\Omega sin\omega cosi & -sin\Omega sin\omega+cos\Omega cos\omega cosi & -cos\Omega sini \cr
	sini.sin\omega  & sini.cos\omega & cosi }
\end{pmatrix}\begin{pmatrix}{
	\xi \cr \eta \cr \zeta }
	\end{pmatrix} $$
En posant:
\ba
	P_X=cos\Omega cos\omega-sin\Omega sin\omega cosi  \nonumber  \\
	P_Y=-cos\Omega sin\omega-sin\Omega cos\omega cosi \nonumber \\
	Q_X=sin\Omega cos\omega+cos\Omega sin\omega cosi \nonumber \\
	Q_Y=-sin\Omega sin\omega+cos\Omega cos\omega cosi \nonumber 
\ea
On obtient comme $\zeta=0$:
\ba
	X_C=P_X \xi+P_Y \eta \label{sa74} \\
	Y_C= Q_X \xi + Q_Y \eta \label{sa75} \\
	Z_C= \xi sini.sin\omega  + \eta sini.cos\omega \label{sa76}
\ea
Si on veut calculer les coordonnées du satellite dans le référentiel terrestre $(O,X_T,Y_T,Z_T)$, on a :
\be
	\begin{pmatrix}{
	X_T \cr Y_T \cr Z_T }
\end{pmatrix}=\begin{pmatrix}{
	cos\Omega & sin\Omega & 0 \cr
	-sin\Omega & cos\Omega & 0 \cr
	0 & 0 & 1 }
\end{pmatrix}.\begin{pmatrix}{
	X_C \cr Y_C \cr Z_C }
\end{pmatrix}  \label{sa77}
	\ee
où $\Omega$ est le temps sidéral de Greenwich\index{Temps sidéral de Greenwich} au temps $t$. Il vaut:
\be            
	\Omega(\mbox{en heures})=1.002737909\times UT2 +HSG_{0TU} \label{sa77a}
\ee
avec $UT2$ le Temps Universel corrigé (en heures) et $HSG_{0TU}$ l'heure sidérale à Greenwich à 0h TU(Temps Universel).
\section{\textsc{Les Perturbations des Orbites}}
On a vu dans le chapitre précédent le mouvement d'un satellite artificiel autour de la terre sous l'action de la force gravitationnelle. Le mouvement réel du satellite est sous l'effet de la force centrale gravitationnelle et d'une force supplémentaire. On supposera que cette force est petite par rapport à la force centrale. On l'appelle force perturbatrice\index{Force perturbatrice}.  
\\

Cette force perturbatrice est la somme de forces d'origine gravitationnelle et d'autres non gravitationnelles. Dans le cas général, une force perturbatrice $\textbf{\textit{f}}$ en un point est fonction de ses coordonnées, de sa vitesse et du temps $t$. On peut écrire alors:
\ba
	f_X=f_X(X,Y,Z,\dot X,\dot Y, \dot Z,t) \nonumber \\
	f_Y=f_Y(X,Y,Z,\dot X,\dot Y, \dot Z,t)  \label{sa78} \\
	f_Z=f_Z(X,Y,Z,\dot X,\dot Y, \dot Z,t) \nonumber 
\ea
Les équations du mouvement en coordonnées rectangulaires obtenues en complètement les équations du problème de la façon suivante:
\ba
m\ddot{X}_C+\frac{\mu m}{r^3}X_C=f_{X_C}(X_C,Y_C,Z_C,\dot{X}_C,\dot{Y}_C, \dot{Z}_C,t)   \nonumber \\
	m\ddot{Y}_C+\frac{\mu m}{r^3}Y_C=f_{Y_C}(X_C,Y_C,Z_C,\dot{X}_C,\dot{Y}_C, \dot{Z}_C,t) \label{sa80} \\
	m\ddot{Z}_C+\frac{\mu m}{r^3}Z_C=f_{Z_C}(X_C,Y_C,Z_C,\dot{X}_C,\dot{Y}_C, \dot{Z}_C,t)   \nonumber
\ea
Comme une force perturbatrice est d'origine gravitationnelle, elle dérive d'un potentiel qu'on note $R$ soit:
$$ \textbf{\textit{f}}=\textit{\textbf{grad}}R $$
 Alors les équations précédentes s'écrivent:
 \ba
m\ddot{X}_C+\frac{\mu m}{r^3}X_C=\frac{\partial R}{\partial X_C} \nonumber \\
	m\ddot{Y}_C+\frac{\mu m}{r^3}Y_C=\frac{\partial R}{\partial Y_C} \label{sa84} \\
	m\ddot{Z}_C+\frac{\mu m}{r^3}Z_C=\frac{\partial R}{\partial Z_C} \nonumber 
\ea
 On fait un changement de variables tel que:
 $$ h(X_C,Y_C,Z_C,\dot{X}_C,\dot{Y}_C, \dot{Z}_C)=g(a,e,i,\Omega,\omega,M) $$
Le système (\ref{sa84}) devient un nouveau système différentiel d'ordre 2 de six inconnues de la forme:
\be
\begin{array}{l}
	\dot a=\Phi_a(a,e,i,\Omega,\omega,M,t)  \\
		\dot e=\Phi_e(a,e,i,\Omega,\omega,M,t) \\
			\dot i=\Phi_i(a,e,i,\Omega,\omega,M,t)\\
				\dot {\Omega}=\Phi_{\Omega}(a,e,i,\Omega,\omega,M,t) \\
					\dot {\omega}=\Phi_{\omega}(a,e,i,\Omega,\omega,M,t) \\	
					\dot M=\Phi_M(a,e,i,\Omega,\omega,M,t) 
					\end{array}
\ee
 Ces six nouvelles variables sont appelées les éléments osculateurs\index{Eléments osculateurs} ou instantanés.
 \\
 
 La solution des équations du mouvement est possible par une méthode analytique ou numérique. On peut dire que le satellite se mouve le long de l'orbite keplérienne, mais les éléments de l'orbite sont, dans ce cas, des fonctions du temps. On l'appelle orbite osculateur\index{Orbite osculateur}.
 \\
 
 Comme $\textbf{\textit{f}}$ a été supposée petite, la solution du système d'équations des éléments osculateurs se présentera en général sous la forme:
 
\be
\begin{array}{l}
	a=a_0+\delta a             \\
	e=e_0+\delta e             \\
	i=i_0+\delta i             \\
	\Omega=\Omega_0+\delta \Omega             \\
	\omega=\omega_0+\delta \omega             \\
	M=M_0+\delta M    
\end{array}						
\ee
 où $\delta a, \delta e,\delta i,\delta \Omega, \delta \omega, \delta M$ seront des petites quantités. Elles sont appelées les perturbations des éléments de l'orbite. L'intérêt de l'emploi des variables osculatrices est que la solution est exprimée sous la forme d'un petit complément à des quantités fixes.  
\section{\textsc{L'Influence du champ de la pesanteur sur le Mouvement du Satellite Artificiel}}
 Un point $M(X,Y,Z)$ de masse unité est soumis au potentiel $V$ de gravitation et au potentiel $\Phi $ de la force centrifuge due à la rotation de la terre.

L'expression de $V$ est :   
$$ V=G\int\!\!\!\int\!\!\!\int_{Terre}\frac{dm'}{r} $$
Malheureusement, cette expression n'est pas calculable  car on ignore la distribution des masses à l'intérieur de la Terre. Il faut appliquer un développement de $V$ en fonctions sphériques sous la forme suivante:
\be           
\fbox{ $ 	V=\ds \frac{\mu}{r}\sum_{n=0}^{+\infty}\sum_{m=0}^n\left(\frac{a}{r} \right)^n\left(C_{nm}cos m\lambda+S_{nm}sinm\lambda \right)P_{nm}(sin\varphi) $}\label{sa100}
\ee
 où:
 
 - $C_{nm},S_{nm}:$ sont les coefficients qu'on obtient par l'observation et ils sont connus;
 
 - $P_{nm}(sin\varphi):$ on les appelle les harmoniques sphériques\index{Harmoniques sphériques} ou polynômes de Legendre\index{Polynômes de Legendre} de deuxième espèce. \index{\textbf{Legendre A.M.}}
 \\
 
Les interprétations des premiers coefficients sont données par:  
\ba
\left\{\begin{array}{ll}
	C_{00}= 1  \\
	C_{10}=\displaystyle \frac{Z_C}{a},\quad C_{11}= \displaystyle \frac{X_C}{a},\quad 	S_{11}= \displaystyle \frac{Y_C}{a}
	\end{array}\right. \label{sa101}
\ea
	et:
\ba
\left\{\begin{array}{lll}
	C_{20}= \displaystyle -\frac{A+B}{2Ma^2}+ \frac{C}{Ma^2}   \\
	\\
	C_{21}=  \displaystyle\frac{E}{Ma^2},\quad C_{22}=  \displaystyle\frac{B-A}{4Ma^2}\\
	\\
		S_{21}= \displaystyle\frac{D}{Ma^2} ,\quad S_{22}= \displaystyle\frac{F}{Ma^2}
		\end{array}\right. \label{saa101}
\ea
où:

- $(X_C,Y_C,Z_C)$ les coordonnées du centre de la masse terrestre;

- $A,B,C$ les moments d'inertie principaux de la Terre;\index{Moments d'inertie principaux}

- $D,E,F$ les moments produits d'inertie;\index{Moments produits d'inertie}

- $M$ la masse de la Terre;

- $a$ le demi-grand axe de l'ellipsoïde terrestre.
\\

Dans le cas où la Terre est représentée par un ellipsoïde de révolution, $V$ est indépendant de la longitude et $V$ s'écrit:
\be
	\fbox{ $ V=\ds \frac{\mu}{r}\left[1+\sum_{n=2}^{+\infty}\left(\frac{a}{r} \right)^nC_{n0}P_n(sin\varphi)\right] $} \label{sa102}
\ee
avec:
\ba
	C_{20}=\frac{C-A}{Ma^2}=-\frac{2}{3}\left(\alpha-\frac{1}{2}q- \frac{1}{2}\alpha^2+ \frac{1}{7}q\alpha \right)\nonumber  \\
	q=\frac{\varpi^2a^3(1-\alpha)}{\mu} \nonumber 
\ea
$\alpha =(b-a)/a$ l'aplatissement de l'ellipsoïde terrestre et $\varpi$ la vitesse angulaire de la rotation de la Terre.
\\

Dans ce cas, on considère le système de coordonnées localisé au centre de la masse terrestre et l'axe $OZ$ confondu avec l'axe de rotation. En conséquence, on a:
\ba
	X_C=Y_C=Z_C=0 \nonumber \\
	A=B \,\,\,\mbox{et}\,\,\,\,D=E=0 \nonumber 
	\ea
	On a aussi tous les coefficients harmoniques tesseraux\index{Harmoniques tesseraux} et sectoriaux \index{Harmoniques sectoriaux} nuls:
	$$ C_{10}=C_{11}=S_{11}=C_{21}=S_{21}=0 $$
On peut mettre:
\be
	J_n=-C_{n0} \label{sa108}
\ee
On obtient alors la formule suivante (\textit{B. Morando}, 1974)\index{Monitor station}:\index{\textbf{Morando B.}}
\be
	\fbox{ $ V=\ds \frac{\mu}{r}\left[1-\sum_{n=2}^{+\infty}\left(\frac{a}{r} \right)^nJ_nP_n(sin\varphi)\right] $} \label{sa109}
\ee
\section{\textsc{Exercices et Problèmes}}
\bex
1. Montrer que: $r=a(1-ecosE) $.

2. Démontrer à partir des formules du cours la relation:
$$ tg\frac{\upsilon}{2}=\sqrt{\frac{1+e}{1-2}}tg\frac{E}{2}$$
Aide: exprimer $tg(\upsilon/2)$ en fonction de $tg\upsilon $.
 \eex
\bex
 A partir de l'expression de $X_C$, montrer que $X_C$ vérifie l'équation du mouvement non perturbé pour la composante $X$, soit: 
$$ \ddot{X}_C+\frac{\mu}{r^3}X_C=0 $$
\eex
\bex
Les altitudes d'un satellite artificiel au périgée et à l'apogée sont respectivement 200 km et 500 km. 

1. Donner les valeurs des paramètres de l'orbite $a,b,p$ et $e$.

2. Calculer les vitesses au périgée et à l'apogée $v_p$ et $v_a$.
\eex
\bpb
La Terre est supposée sphérique, homogène de rayon $R = 6\,371\,000\,m$. Le produit de la constante universelle de gravitation terrestre $G$ par la masse $M$ de la Terre soit $GM = 3.986\, 005\,10^{14}\,m^3s^{-2}$. Un satellite géodésique a une trajectoire telle que son altitude maximale est $1100\,km$ et son altitude minimale $800\,km$.

1. Donner la période de ce satellite.

2. Quelle est l'excentricité de sa trajectoire?

3. On mesure la distance du satellite à une station au sol de latitude $43°5$ et d'altitude nulle, lors du passage du satellite à la verticale de la station, soit $D = 812\,000\,m$.

a - Quelle est l'anomalie vraie du satellite à cet instant, sachant qu'il vient de passer au périgée.

b - Combien de temps s'est écoulé depuis le passage au périgée?
\epb
\bpb
Une comète décrit autour du Soleil une ellipse d'excentricité $e$ de demi-grand axe $a$ et de demi-petit axe $b$ où le Soleil occupe un des foyers. L'équation de l'orbite de la comète en coordonnées polaires est donnée par:
$$ r=\frac{a(1-e^2)}{1+ecos\upsilon}$$ avec $r$ la distance Soleil- comète.

1. Déterminer les distances $r_A$ et $r_P$ lorsque la comète est à l'apogée et au périgée en fonction de $a$ et $e$.

2. La comète de Halley a une orbite fortement excentrique : son apogée est à 0.53 $UA$  du Soleil et sa périgée est à 35.1 $UA$. Calculer $e$.

3. En utilisant la loi des aires et la troisième loi de Kepler, montrer que la constante des aires $C$ est exprimée par:
$$C^2=\frac{b^2}{a}G.M$$
où $G, M$ désignent respectivement la constante de la gravitation universelle et la masse du Soleil.

4. On pose : $\ds u=\frac{1}{r}$. Donner l'expression du carré de la vitesse $v^2$ de la comète en fonction de $u$ et $\ds \frac{du}{d\upsilon}$. Montrer que $v^2$ peut s'écrire sous la forme:
$$ v^2=G.M\left(\frac{2}{r}-\frac{1}{a}\right)$$

5. Déterminer l'expression du rapport des vitesses à l'apogée et au périgée $\ds \frac{v_A}{v_P}$ en fonction de $e$.

6. Calculer numériquement ce rapport pour le cas de la comète de Halley.

On donne:

- 1 $UA=149\, 597\,870 \,km$; 

- $G=6.672\times 10^{-11}\,m^3.kg^{-1}.s^{-2}$;

- $M=1.9891\times 10^{30}\,kg$.
\epb
\bpb
Soient deux corps  $\textbf{\textit{X}}_1(m_1)$ et $\textbf{\textit{X}}_2(m_2)$ soumis à l'attraction universelle dans un repère orthonormé $\m{ R}\,(O,e_1,e_2,e_3)$. On pose:
\ba
\textit{\textbf{OM}}_1=\textbf{\textit{X}}_1(m_1)\nonumber \\
\textbf{\textit{OM}}_2=\textbf{\textit{X}}_2(m_2) \nonumber
\ea
avec $\mu $ la constante universelle de la gravitation. 

1. Montrer que les équations du mouvement des points $M_1, M_2$ dans $\m R$ sont données par:
\ba
   \frac{d^2\textbf{\textit{X }}_1}{dt^2}=-\mu   m_2\frac{(\textbf{\textit{X}}_2-\textbf{\textit{X}}_1)}{\| \textbf{\textit{X}}_2-\textbf{\textit{X}}_1 \|^3}\nonumber \\
   \frac{d^2\textbf{\textit{X }}_2}{dt^2}=-\mu   m_1\frac{(\textbf{\textit{X}}_1-\textbf{\textit{X}}_2)}{\| \textbf{\textit{X}}_2-\textbf{\textit{X}}_1 \|^3}\nonumber 
\ea
2. Soit $ \m G$ le centre de gravité des deux corps $M_1,M_2$. Montrer que $\m G$ vérifie: 
$$ \textbf{\textit{O}}\m G=\frac{m_1}{m_1+m_2}\textbf{\textit{OM}}_1+\frac{m_2}{m_1+m_2}\textbf{\textit{OM}}_2 $$
En déduire que le mouvement de $\m G$ est rectiligne et uniforme dans $\m R$.

3. On pose $m=m_1+m_2$, $\alpha=\ds \frac{m_1}{m}$. En déduire que :$\ds \frac{d^2\textbf{\textit{X}}_2}{dt^2}=\frac{-\alpha}{1-\alpha}\frac{d^2\textbf{\textit{X}}_1}{dt^2}$.

4. On note :
\ba
\textbf{\textit{Y}}_1=\textbf{\textit{X}}_1-\textbf{\textit{O}}\m G \nonumber \\ 
\textbf{\textit{Y}}_2=\textbf{\textit{X}}_2-\textbf{\textit{O}} \m G \nonumber 
\ea
Montrer qu'on a:
\ba
\frac{d^2\textbf{\textit{Y }}_2}{dt^2}=\frac{-\alpha}{1-\alpha}\frac{d^2\textbf{\textit{Y}}_1}{dt^2} \nonumber \\
\frac{d^2\textbf{\textit{Y }}_1}{dt^2}=-\mu m(1-\alpha)^3\frac{\textbf{\textit{Y}}_1}{\|\textbf{\textit{Y}}_1\|^3}\nonumber
\ea
5. Quel est le mouvement de $M_2$.
\epb

\chapter{\textit{\textbf{Le Système GPS}}}\index{Système GPS}
\section{\textsc{Introduction}}
  Le système de positionnement global (Global Postioning System GPS), appelé aussi NAVSTAR/GPS (Navigation System by Timing And Ranging), est  un système de navigation par repérage du temps et mesures des distances, et c'est un système mondial de positionnement par satellites conçu et mis en service par le Département Américain de la Défense (DoD). Il détermine la position des points au sol équipés de récepteurs enregistrant des mesures d'origine satellitaire.
	
	Il a été développé en vue du remplacement du système de positionnement TRANSIT (Doppler)\index{Système Doppler} qui présentait deux défauts importants: 

		- une couverture satellitaire insuffisante;

		- une faible précision en navigation.

	Les spécifications initiales sont d'avoir accès à une position absolue dans un système de référence mondial avec une précision métrique et au temps avec la précision de la microseconde.

	Le développement du système GPS a commencé dans les années 70, la mise en place démarrant en 1978 par le lancement du premier satellite. Le système a été déclaré opérationnel par le DoD début février 1994. Depuis 52 satellites ont été lancés.
\section{\textsc{Aspects Généraux}}
\subsection*{16.2.1. Les Satellites}
Le système complet comporte 34 satellites dont 31 sont opérationnels. Ils sont répartis de façon à assurer en tout lieu une visibilité simultanée de 4 à 8 satellites  avec une élévation d'au moins 15° au dessus de l'horizon.
	
	Les satellites sont répartis sur 6 plans orbitaux, ayant tous une inclinaison d'environ 55° sur l'équateur. L'orbite des satellites est quasi-circulaire, leur altitude est d'environ 20 000 km et leur période d'environ 12 heures.
	
	Chaque satellite est muni d'un émetteur-récepteur, d'une horloge de haute précision, d'ordinateurs et d'équipements auxiliaires destinés au fonctionnement du système. On distingue cinq classes de satellites qui correspondent chacune à une étape dans la constitution du système:

*	Le \textit{Block I }: 11 satellites lancés entre 1978 et 1985. Cet ensemble a constitué la phase initiale de test du système. Aucun  satellite de ce bloc n'est utilisé à l'heure actuelle.

*	Le \textit{Block IIA} : le premier a été lancé en avril 1991. Un seul satellite est en activité.

*	Le \textit{Block IIR} : comprenant 11 satellites lancés entre 1997 et 2004 munis d'horloges en rubidium.

* Le \textit{Block IIR-M} est constitué de 7 satellites, lancés entre septembre 1997 et novembre 2009. Les derniers satellites émettent le deuxième signal civil L2C sur la fréquence L2. 

* Le \textit{Block IIF} est constitué de 12 satellites dont le premier a été lancé en 2010 en émettant en plus le deuxième signal civil L2C et le troisième signal civil L5. Le dernier a été lancé en février 2016. La durée de vie de ces  satellites est de 12 ans.

* Trois satellites sont en réserve.
\\

Le premier satellite de GPS III, la nouvelle génération des satellites GPS, serait lancé probablement en 2017. 

	Un secteur de contrôle\index{Secteur de contrôle} composé de 5 stations au sol (Hawaii, Colorado Springs, Ascencion, Diago Garcia et Kwajalein) qui enregistrent en continu les signaux GPS sur les deux fréquences L1, L2 et sur la fréquence L2C (pour les satellites du  \textit{Block IIR-M} et du \textit{Block IIF}) et la fréquence L5 (pour les satellites du \textit{Block IIF}). Les tâches du secteur de contrôle sont:

	- capter les satellites GPS;

	- analyser les orbites et prédire;
 
  - mettre à jour les messages de navigation.

    Le secteur de contrôle se compose de :

	- \textit{Monitor station} (MS) : on observe les éphémérides et l'horloge,\index{Monitor station}

	- \textit{Master control station} (MCS) située à Colorado Springs où on effectue les opérations suivantes : \index{Master control station}

  *	les calculs des erreurs;

  *	les corrections de l'orbite et la fréquence d'horloge;
 
  *	la création de nouveaux messages de navigation;

  - \textit{Upload station} (antennes sur terre) : on y envoie vers les satellites les messages de navigation.
\subsection*{16.2.2. Le Message GPS}
  Les satellites émettent en permanence deux signaux ultra-stables sur les deux fréquences spécifiques du système L1 et L2. Sur ces deux fréquences, est modulé un code particulier dont le décodage fournit le message GPS. 

	Ce message se compose des éléments suivants :

-	prévisions de l'orbite des satellites, provenant des stations de contrôle, afin de permettre aux utilisateurs d'effectuer des calculs de navigation en temps réel;

-	l'information concernant la qualité des horloges des différents satellites et le modèle de développement polynomial du comportement de ces horloges;

-	information d'ordre général sur le système.
 \subsection*{16.2.3. Protection du Système GPS}
Le système GPS était pourvu de deux processus qui permettent de limiter son utilisation civile et de protéger son utilisation militaire :

*	L'accès sélectif (SA : Selective Availability) comporte deux dégradations :\index{Selective Availability}

			- dégradation de la fréquence de l'horloge des satellites par ajout d'un déphasage;
			
			- dégradation des éléments du message radiodiffusé (éphémérides des satellites, paramètres d'horloge des satellites).
\\

	Ces dégradations étaient connues et diffusées uniquement à des utilisateurs autorisés.
	
*	L'antibrouillage (AS : Anti-Spoofing) est une protection destinée à éviter le brouillage volontaire du système GPS par un utilisateur qui enverrait un signal proche de celui de GPS, créant ainsi la confusion et des erreurs de positionnement pour les autres utilisateurs.  \index{Anti-Spoofing}

Ces deux protections ont été levées.
\section{\textsc{Les Instruments de mesures GPS}}
Pour optimiser l'utilisation du GPS, une gamme très variée de récepteurs a été mise en \oe uvre. Chaque type de récepteurs a ses spécifications :

 *  GPS portatif : dans ce type d'instrument, il existe plusieurs modèles;
	
 *	poids varie de 397 g à 850 g;

 * 	précision :
 
 - mesure absolue : 15 à 100 m;
 
 - mesure différentielle : 2 à 3m;

 *	détermination de la position en deux dimensions avec trois satellites ou en mode de mesures à 3 dimensions (position et altitude avec 4 satellites);

 *  GPS stationnaire : ce type d'instrument permet d'obtenir une meilleure précision que les instruments portatifs. On distingue suivant le mode de fixation du récepteur :

  -	GPS fixé sur trépied : il permet de faire des mesures satellitaires en mode statique et statique rapide;

  -	GPS mobile : l'antenne est fixée sur une tige. Il permet de réaliser le mode STOP and GO. Mais pour réaliser le temps réel, il nécessite un équipement radio supplémentaire qui permet à un poste fixe d'envoyer sa position à l'autre utilisateur.
	\subsection*{16.3.1. Les Récepteurs  Géodésiques}\index{Récepteur géodésique}
  Les récepteurs géodésiques GPS sont des appareils qui enregistrent le message des satellites par l'intermédiaire d'une antenne stationnée sur un point dont on veut connaître la position. Par un processus basé essentiellement sur la connaissance du code qui module le signal reçu, les récepteurs effectuent les mesures GPS.
 
Ils décodent par eux-mêmes les messages provenant des satellites. 
\\

	Il existe essentiellement deux types de récepteurs:

 -	les récepteurs mono-fréquences qui n'enregistrent que les signaux de la fréquence L1;
 
 -	les récepteurs bi-fréquences qui enregistrent les signaux des deux fréquences L1 et L2.

\section{\textsc{Principes de mesures GPS}}
Les satellites émettent des signaux qui sont reçus, interprétés par des récepteurs au sol. A partir d'une fréquence fondamentale (10.23 $MHz$), l'émetteur génère deux ondes L1 et L2 de fréquence respectivement 1575.42 $MHz$ et 1227.60 $MHz$.\index{Fréquence fondamentale}
\\

\begin{figure}[b]
	\centering
		\includegraphics[width=0.90\textwidth]{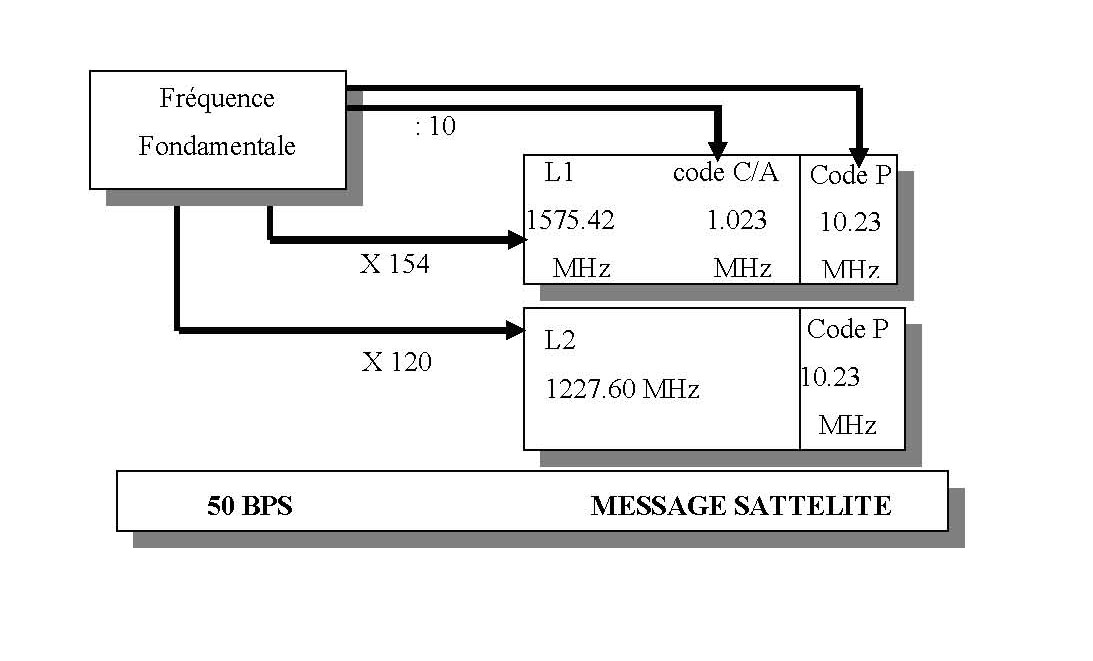}
	\caption{Les fréquences du GPS}
	\label{fig:gps1}
\end{figure}

	L1 et L2 sont modulées par des codes aléatoires\index{Codes aléatoires} (code C/A et code P) ainsi que par un message qui contient en particulier les données sur les éphémérides du satellite. Sur l'onde L1 on trouve le code C/A et le code P, tandis que sur L2, il n'y a que le code P. Le message existe sur les deux ondes.

	Le récepteur génère une ou deux ondes qui sont une réplique des ondes L1 et L2 émises par le satellite, ainsi que le code C/A et le code P. Les mesures consistent à comparer cette réplique synthétisée avec le signal reçu soit en mesurant le décalage en temps sur les codes, soit la différence de phases sur les ondes porteuses.

	Les mesures sur les codes donnent directement la distance satellite - récepteur en connaissant la durée de propagation du signal émis par le satellite jusqu'à son arrivée à un récepteur sur la terre.

	\section{\textsc{Les Equations Fondamentales d'Observations}}

\subsection*{16.5.1. L'Equation d'Observations en Pseudo-Distances }\index{Observations en pseudo-distances}
La mesure  de la pseudo-distance est l'observation GPS de base pour tout récepteur. Elle utilise les codes C/A et P. La mesure enregistre le temps apparent $\Delta t$  mis par le signal pour arriver du satellite au récepteur :
$$ \Delta t =  t_R - t^S = (t_{R(GPS)} - \delta_R) - (t^{S(GPS)} -  \delta^S) =  \Delta t_0 + \Delta \delta   $$
où $\Delta t_0$  désigne le temps vrai et :
$$ 	\Delta \delta = \delta ^S - \delta_R $$
la différence entre les corrections aux temps du récepteur et du satellite. L'intervalle  $\Delta t$ ci-dessus est multiplié par la vitesse de la lumière $c$ pour avoir une distance :
$$    R = c\Delta t  = c.\Delta t_0 + c\Delta \delta   = \rho + c\delta                            $$
Si la position du satellite est représentée par le vecteur $\textbf{\textit{r}}$ et celle du récepteur par $\textbf{\textit{R}}$, la vraie distance  $\rho $  peut être calculée à partir de :
$$  \rho     = \|\textbf{\textit{r}} - \textbf{\textit{R}} \|                                    $$
Si  l'erreur sur l'orbite du satellite est notée par $d\rho$  et si on prend compte des retards dûs au passage du signal dans la troposphère et l'ionosphère, alors la pseudo-distance peut être calculée par l'équation:\index{Troposphère}\index{Ionosphère}
\be
\fbox{ $ 	R=\rho +d\rho+c\Delta \delta + \Delta\rho_{ion} + \Delta\rho_{trop} $} \lb{eqpd}
\ee
\subsection*{16.5.2. L'Equation d'Observation de Phases}\index{Observations de phases}
Notons $\varphi^{S}(t)$ la phase de l'onde reçue avec la fréquence $f^S$ et $\varphi_{R}(t)$  la phase générée par le récepteur avec la fréquence $f_R$. Ici le temps $t$ est une époque dans le temps $GPS$ compté à partir d'un instant initial $t_0 = 0$, alors on peut écrire que:
\ba
	\varphi^S(t)=f^S(t-\frac{\rho}{c})-\varphi^S_0\nonumber \\
	\varphi_R(t)=f_R .t-\varphi_{0R} \nonumber
\ea
Les phases initiales $\varphi^{S}_{0}$ et $\varphi_{0R}$ sont causées par les erreurs des horloges du satellite et du récepteur et elles sont égales à:
\ba
\varphi^S_0=	f^S.\delta^S \nonumber \\
\varphi_{0R}=f_R.\delta_R \nonumber
\ea
Par suite, la phase de battement est égale à :\index{Phase de battement}
\be
	\varphi^{S}_{R}(t)=\varphi^{S}(t)-\varphi_{R}(t)=-f^S.\frac{\rho}{c}-f^S.\delta^S+f_R.\delta_R+(f^S-f_R).t \label{gps10}
\ee
Comme  les fréquences $f^S$   et $f_R$   sont presque les mêmes, l'équation (\ref{gps10}) peut s'écrire :
\be
	\varphi^{S}_{R}(t)=-f.\frac{\rho}{c}-f.\Delta \delta
\ee

 avec $\Delta\delta   =   \delta^S -\delta_R$ .
\\

A l'époque $t$, la phase de battement est :
\be
	\varphi^{S}_{R}(t)=\left[\Delta\varphi^{S}_{R}(t)\right]^{t}_{t_0}=\Delta\varphi^{S}_{R}(t)+N 
\ee
où $N$ désigne le nombre de cycles entre le satellite et le récepteur. $N$ est appelé aussi  ambiguïté\index{Ambiguïté entière} entière. Tant que le récepteur reçoit le signal du satellite, $N$ reste le même  et $\Delta\varphi^{S}_{R}(t)$ est la fraction de phase mesurée depuis le début de la réception du signal.  

On pose :
\be
	\Phi=-\Delta\varphi^{S}_{R}(t)\label{gps13}
\ee
Alors l'équation (\ref{gps13}) devient:
\be
	\Phi=f.\frac{\rho}{c}+f.\Delta\delta+N \label{gps14}
\ee
En introduisant la longueur d'onde $\lambda$, l'équation  (\ref{gps14}) devient :
\be
\fbox{  $ 	\Phi=\ds \frac{\rho}{\lambda} + c.\frac{\Delta \delta}{\lambda}+ N  $} \label{gps15}
\ee
L'équation (\ref{gps15}) représente l'équation d'observations de phases. $\rho$  désigne la distance entre le récepteur et le satellite. La phase peut être mesurée mieux que 0.01 cycles ce qui correspond à une précision millimétrique.  En effet, de (\ref{gps15}), on peut écrire que:
   \be
	\lambda d\Phi=24\,cm\times 0.01=0.24\, cm=2.4\,mm \label{gps15a}
\ee
où on a considéré que la longueur d'onde $\lambda=24\,cm$.                                                                    

\subsection*{16.5.3. Le Référentiel $WGS84$}\index{Référentiel WGS84}
Le GPS fournit les coordonnées géodésiques ($\phi$,$\lambda$,$he$) dans le système géocentrique $WGS84$(World Geodetic System 84)\index{Système géocentrique $WGS84$} ($O,X, Y, Z$).         
  
   Il est défini par : 
  
-	l'origine : $O$ centre des masses de la Terre;

-	l'axe $OZ$ : parallèle à la direction de l'origine vers le Pôle Conventionnel Terrestre (PCT);\index{Pôle Conventionnel terrestre}

-	l'axe $OX$ : le plan $OZX$ est parallèle au méridien de longitude zéro défini par le Bureau International de l'Heure (BIH);\index{Bureau International de l'Heure}

-	l'axe $OY$ : dans l'équateur du PCT et perpendiculaire au plan $OZX$ dans le sens direct.

	Au référentiel $WGS84$, on associe l'ellipsoïde $WGS84$ \index{Ellipsoïde $WGS84$} dont les paramètres sont :
	
-	$a$ le demi-grand axe : $6\,378\,137.00\, m$;

-	l'aplatissement $f = 1/298.257\,223\,563$.  
   
Au lieu d'utiliser l'ellipsoïde $WGS84$, on utilise l'ellipsoïde $GRS80$ adopté par l'Association Internationale de Géodésie (AIG) ayant le même demi-grand axe que le $WGS84$ et un aplatissement de $1/298.257\,222\,101$. 
\subsection*{16.5.4. Les Coefficients de précision des dispositions des satellites GPS}
Comme en levés topographiques classiques où le géomètre opère des observations, le cas d'un relèvement, en visant des points connus bien disposés, les observations convenablement exécutées des observations GPS obéissent à cette règle.

En effet, on définit les coefficients de précision relatifs à la disposition des différents satellites par rapport au point de l'observation. Ces coefficients permettent de donner à l'opérateur le choix de la période des observations où les satellites GPS sont bien disposés dans le ciel de façon à obtenir des précisions meilleures. Ces coefficients sont notés comme suit:
\be
\fbox{$ \begin{array}{l}
GDOP:\quad \mbox{le coefficient de précision géométrique} \\ \index{Le GDOP} 
PDOP:\quad \mbox{le coefficient de précision de position} \\ \index{Le PDOP}
HDOP:\quad \mbox{le coefficient de précision horizonthale} \\ \index{Le HDOP}
VDOP:\quad \mbox{le coefficient de précision verticale} \\ \index{Le VDOP}
TDOP:\quad \mbox{le coefficient de précision du temps} \index{Le TDOP}
\end{array} $}
\ee 
Comment sont calculés ou estimés ces coefficients?

A partir des équations fondamentales d'observations GPS (\ref{eqpd}) ou (\ref{gps15}), on arrive au système suivant, pour les observations GPS  en un point $A$ de coordonnées approchées tridimensionnelles $(X_0,Y_0,Z_0)$ ou géodésiques $(\varphi_0,\lambda_0,h_0)$:
\be
 	A.U=B+V \lb{eqpd1}
\ee
où:

- $U=(dX,dY,dZ,\Delta \delta )^T$le vecteur des corrections à la position approchée $(X_0,Y_0,Z_0)$ et $\Delta \delta $ la correction du temps;

- $A$ la matrice des coefficients;

- $B$ le vecteur des observations;

- et $V$ le vecteur des résidus.

La solution de (\ref{eqpd1}) par la méthode des moindres carrés (voir (\ref{c49})) est donnée par :
\be
\ov{U}=(A^TA)^{-1}A^TB \lb{eqpd2}
\ee
La matrice $ Q=(A^TA)^{-1}$, qu'on appelle la matrice de variance des inconnues, s'écrit sous la forme:
\be
Q=\begin{pmatrix}{
q _{XX} & q_{XY} & q_{XZ} & q_{Xt} \cr
 q       _{YX} & q       _{YY} & q       _{YZ} & q       _{Yt} \cr
q       _{ZX} & q       _{ZY} & q       _{ZZ} & q       _{Zt} \cr
q       _{tX} & q       _{tY} & q       _{tZ} & q       _{tt} }
\end{pmatrix} \lb{eqpd3}
\ee
D'où les expressions de :
\be
\fbox{ $ \begin{array}{l}
GDOP=\sqrt{q       _{XX} + q       _{YY} + q       _{ZZ} + q       _{tt}} \\
\\
PDOP=\sqrt{q       _{XX} + q       _{YY} +q        _{ZZ}                } \\
\\
TDOP=\sqrt{q       _{tt}} 
\end{array} $}\lb{eqpd4}
\ee
Ces coefficients sont définis car la matrice $Q$ est définie positive. Pour donner les expressions de $HDOP$ et $VDOP$, on doit se ramener au repère local du point $A$.

\subsubsection*{16.5.4.1. Le Repère local}
  Soit un point $ A (\varphi ,\lambda,he)$ relatif à un ellipsoïde de r\'evolution associ\'e \`a un r\'ef\'erentiel g\'eocentrique donn\'e $ \mathcal{R}$. 
  
\begin{figure}[hbp]
	\centering
		\includegraphics[width=0.80\textwidth]{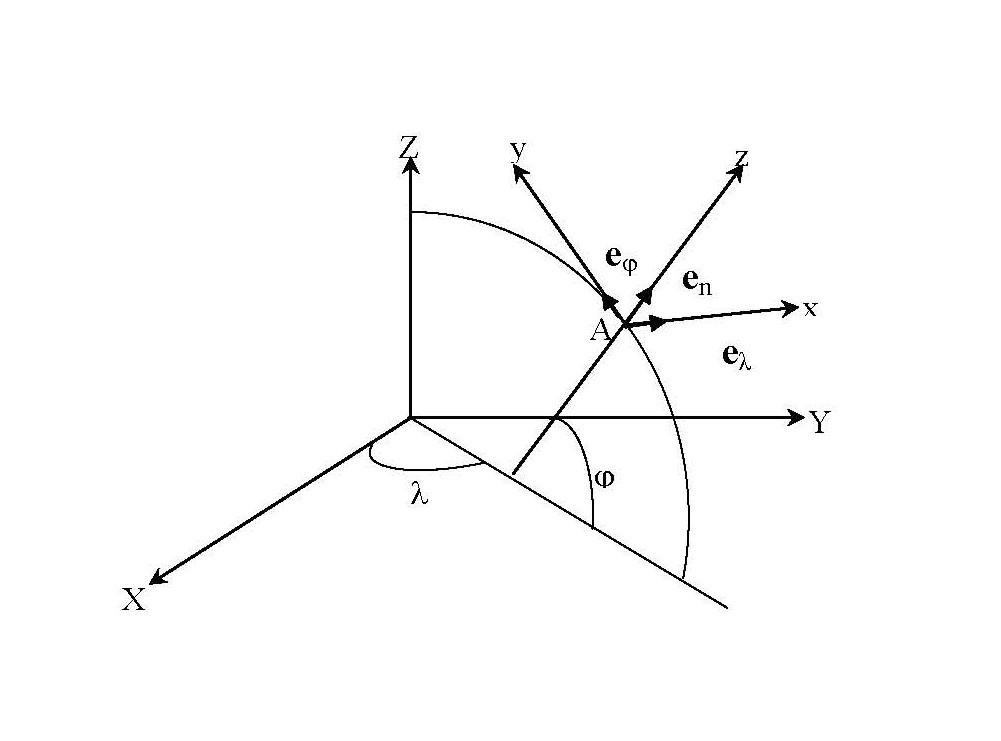}
	\caption{Le rep\`ere local}
	\label{fig:Normale1a}
\end{figure}
  
  On consid\`ere le rep\`ere orthonorm\'e local en A $(e_{\lambda},e_{\varphi},e_n)$ d\'efini dans la base orthonorm\'ee $(i,j,k)$ de $ \mathcal{R}$ (\textbf{Fig. \ref{fig:normale1a}}) par:
\be 
	e_{\lambda}=\left|
\begin{array}{lll}
	-sin\lambda \\
	cos\lambda \\
	0
\end{array}\right. ;\quad	e_{\varphi}=\left|
\begin{array}{lll}
	-sin\varphi cos\lambda \\
	-sin\varphi sin\lambda \\
	cos\varphi
\end{array}\right. ;\quad 	e_n=\left|
\begin{array}{ll}
	cos\varphi cos\lambda \\
	cos\varphi sin\lambda \\
	sin\varphi
\end{array} \right. \label{e1}
\ee
Matriciellement, on peut \'ecrire (\ref{e1}) sous la forme:
\be
	\begin{pmatrix}{
	e_{\lambda} \cr 
	e_{\varphi} \cr
	e_n }
\end{pmatrix}=\begin{pmatrix} {
		-sin\lambda & cos\lambda & 	0 \cr
		-sin\varphi cos\lambda &	-sin\varphi sin\lambda & 	cos\varphi \cr
		cos\varphi cos\lambda & cos\varphi sin\lambda & sin\varphi }
\end{pmatrix}\begin{pmatrix}{
	i \cr
	j \cr
	k }
\end{pmatrix}= R\begin{pmatrix}{
	i \cr
	j \cr
	k }
\end{pmatrix}  \label{e2}
\ee
avec $R$ la matrice:
\be
	\fbox{ $ R=\begin{pmatrix}{
		-sin\lambda & cos\lambda & 	0 \cr
		-sin\varphi cos\lambda &	-sin\varphi sin\lambda & 	cos\varphi \cr
		cos\varphi cos\lambda & cos\varphi sin\lambda & sin\varphi }
\end{pmatrix} $} \label{e3}
\ee
Soit un vecteur  de composantes $W=(X,Y,Z)^T$ et $w=(x,y,z)^T$ respectivement par rapport aux bases $(i,j,k)$ et $(e_{\lambda},e_{\varphi},e_n)$, on a alors la relation:
$$ 	\begin{pmatrix}{
	x              \cr
	y           \cr
	z   }
\end{pmatrix}=\begin{pmatrix} {
		-sin\lambda & cos\lambda & 	0 \cr
		-sin\varphi cos\lambda &	-sin\varphi sin\lambda & 	cos\varphi \cr
		cos\varphi cos\lambda & cos\varphi sin\lambda & sin\varphi }
\end{pmatrix}\begin{pmatrix}{
	X \cr
	Y \cr
	Z }
\end{pmatrix}= R\begin{pmatrix}{
	X \cr
	Y \cr
	Z }
\end{pmatrix} $$
ou encore:
$$   w=R.W $$
Or d'après la formule (\ref{c60}), la matrice de variance $Q_w$ de $w$ s'écrit:
\be
Q_w=RQ_W.R^T \lb{eqpd5a}
\ee
\subsubsection*{16.5.4.2. Calcul de $HDOP$ et $VDOP$ }
Pour calculer $HDOP$ et $VDOP$, on exprime la matrice $Q'$ extraite de la matrice $Q$ donnée par (\ref{eqpd3}) soit:
\be
Q'=\begin{pmatrix}{
q _{XX} & q_{XY} & q_{XZ}  \cr
 q       _{YX} & q       _{YY} & q       _{YZ}  \cr
q       _{ZX} & q       _{ZY} & q       _{ZZ} }
\end{pmatrix} \lb{eqpd5}
\ee
On calcule alors la matrice transformée de $Q'$ par (\ref{eqpd5a}), d'où:
\be
Q''=R.Q'.R^T= \begin{pmatrix}{
q'' _{xx} & q''_{xy} & q''_{xz}  \cr
 q''       _{yx} & q''       _{yy} & q''       _{yz}  \cr
q''       _{zx} & q''       _{zy} & q''       _{zz} }
\end{pmatrix} \lb{eqpd6}
\ee
Par suite:
\be
\fbox{ $ \begin{array}{l}
HDOP=\ds \sqrt{q'' _{xx} + q''_{yy}}  \\
\\
VDOP=\ds \sqrt{q'' _{zz}}
\end{array} $}
\ee
Là aussi, les coefficients $HDOP$ et $VDOP$ sont définis car la matrice $Q''$ est aussi définie positive.
\section{\textsc{Les Différents types de Positionnement Par GPS}}
\subsection*{16.6.1. Le Positionnement Absolu}\index{Positionnement absolu}
Le positionnement absolu d'un point est le processus de collection de données à partir de plusieurs satellites en un endroit donné, constamment avec une éphéméride, pour déterminer la position d'une station indépendante.

L'utilisateur muni d'un récepteur mesure la distance entre sa station et plusieurs satellites à l'aide du code C/A ou P. L'intérêt de cette méthode est l'obtention d'une position instantanée en temps réel.

 -	Précision :  - 8 mètres avec le code C/A et S/A inactivé (théoriquement);
 
  					    - 40 mètres avec le code C/A et S/A activé;
  					    
				         - 8 mètres avec le code P.
				          
 -	Les besoins : 
 
 * 4 heures d'observations par session;
 
 * 2 à 3 sessions par emplacement;
 
 * avec un récepteur bi-fréquence (L1/L2);
 
 * éphémérides précises;
 
 * un logiciel de traitement.
\begin{figure}[ht]
	\centering
		\includegraphics[width=0.80\textwidth]{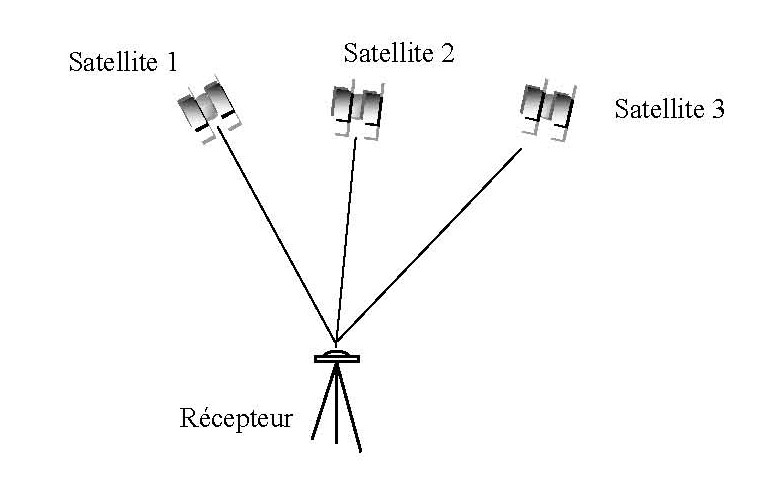}
	\caption{Positionnement absolu}
	\label{fig:gps4}
\end{figure}

\subsection*{16.6.2. Le Positionnement Relatif}\index{Positionnement relatif}
	C'est un processus de collection de données simultanément en deux ou plusieurs stations, à partir d'au moins quatre satellites, pour déterminer les positions des points relativement à d'autres points connus. Le positionnement des points est relatif d'un point à un autre.
	
	La distance est calculée soit à partir du code P ou C/A, soit à partir de mesures de phases. Ce positionnement permet de diminuer l'influence des erreurs d'orbites, des corrections de propagation, des décalages des horloges satellites et récepteurs.

	Le positionnement peut être statique ou dynamique (en général un poste fixe, un mobile). Le temps réel peut être réalisé par ajout d'un équipement de transmission radio qui permet à un poste d'envoyer ses observations à l'autre qui peut alors calculer la position et la vitesse du mobile.

	Le champ d'utilisation est vaste allant de la navigation à l'établissement de réseaux de précision (géodésie, topographie,...etc).
\\	
	
Procédures :
  * double différence (deux satellites), les stations observent les satellites simultanément à un instant $t$.

\begin{figure} [bh]
\centering
		\includegraphics[width=0.90\textwidth]{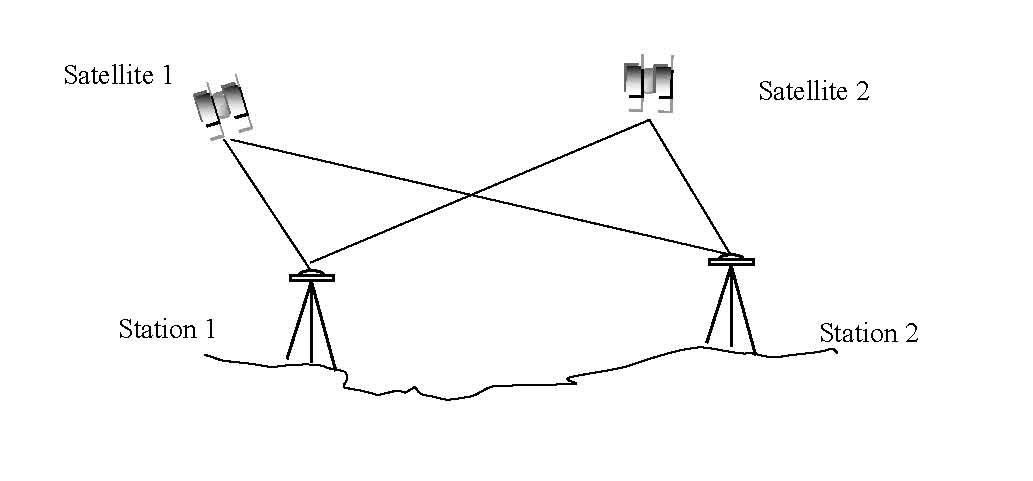}
\caption{Positionnement relatif(double différence)}
\label{fig:gps5}
\end{figure}

	* triple différence ($2\times$double différence), les stations observent les satellites à un instant $t$ puis à un instant $t+1$.

\begin{figure}[b]
	\centering
		\includegraphics[width=0.90\textwidth]{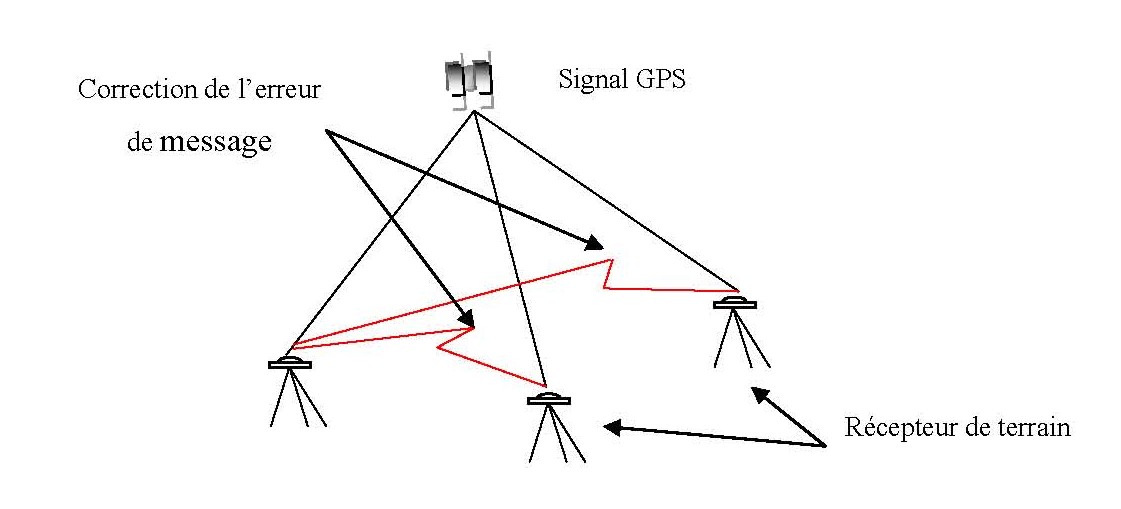}
	\caption{Postionnement relatif (Triple différence)}
	\label{fig:gps6}
\end{figure}

Précision : 

 - observations en mode statique 2 à $3\,mm$;
 
 - observations en mode cinématique 2 à $3\,cm$.
%
\subsection*{16.6.3. La Navigation Différentielle DGPS}\index{Navigation Différentielle}
La navigation différentielle ou le GPS différentiel (DGPS) consiste à faire des observations de pseudo distances en plusieurs stations simultanément.

	Le premier type d'application du DGPS est dynamique, c'est à dire une station fixe appelée station de référence et l'autre est mobile. La station fixe dont les coordonnées sont connues envoie, à l'aide d'un équipement radio (UHF, VHF, HF, MF), des corrections à la station mobile qui calcule alors sa position. Ce type de positionnement est très utilisé en navigation maritime à l'approche des côtes ainsi que pour l'hydrographie. Le deuxième type d'application est statique et profite de la bonne précision du DGPS pour l'établissement de réseaux de points dont la précision n'a pas besoin d'être meilleure qu'un mètre.
\subsection*{16.6.4. Le Mode Statique Absolu}\index{Mode statique absolu}
Le statique est le mode le plus utilisé par la technique GPS pour l'établissement de canevas ou réseaux en géodésie, topométrie et photogrammétrie. Il s'agit d'observer sur chaque station les phases sur au moins 4 satellites, le résultat est un positionnement relatif, c'est-à-dire un positionnement par rapport à une station de référence en connaissant $Dx, Dy$. Le temps d'observations doit être suffisant pour résoudre les ambiguïtés et dépend de la longueur de la ligne de base, de la géométrie de la constellation des satellites et des conditions atmosphériques. Par expérience, la durée d'observations est d'au moins une heure pour des lignes de bases courtes (15 km) et elle est de plusieurs heures pour les lignes de longues bases.

	La précision dépend du type de récepteurs (mono ou bi-fréquence), du nombre de récepteurs et du type d'orbites utilisées :

*	GPS statique mono fréquence :

         - mesures de phases sur L1;
         
		     - durée d'observations : une heure;
		     
 		     - distances jusqu'à 20 km en zone tempérée, 10 km en zone équatoriale;
 		     
         - orbites radio- diffusées;
         
         - précision : $2\,mm + 2\times10^{-6}$ $D$;
         
où $D$ représente la distance entre les points.

*	GPS statique bi-fréquence :

 		   - mesures de phases sur L1 et L2;
 		  
       - durée d'observations : de 1 à 4 heures;
                
       - distances : jusqu'à 300 km;
                
       - orbite radio diffusée;
                
      - précision : $2\, mm + 10^{-6}$ $D$.   

*	GPS statique ultra précis :

               - mesures de phases sur L1 et L2;
               
               - durée d'observations : de 1 à plusieurs jours;
               
               - distances de l'ordre de 5000 km;
               
               - calcul par traitement semi-dynamique ou avec orbites précises;
                
               - précision : $10^{-8}D$ à $10^{-9}D$.
               
\subsection*{16.6.5. La Statique Rapide}\index{Mode statique rapide}

	On choisit une station de référence sur laquelle un récepteur mesure en continu, tandis que des récepteurs mobiles se déplacent sur les autres stations en observant quelques minutes sur chacune d'elles.

	Cette technique très productive convient aux petits réseaux pour le cadastre, la topométrie, la densification de réseaux géodésiques et l'équipement photogrammétrique.
	
\subsection*{16.6.6. La Réoccupation ou pseudo-cinématique }
La réoccupation est une autre variante du statique. Le principe des calculs est de grouper toutes les observations faites sur un site, lorsque le site est stationné plusieurs fois. Si on observe trois satellites à la première occupation, et trois satellites à la deuxième  occupation, le logiciel fera comme si  six satellites avaient été observés.

\subsection*{16.6.7. Le Mode Opératoire }
Le récepteur de référence (base) est stationné sur un point connu, il collecte les données, calcule la correction de 'pseudo-range' en se référant à sa position connue et envoie ces corrections aux utilisateurs (itinérants). Les récepteurs de terrain reçoivent ces corrections et les utilisent pour corriger leurs positions relativement à la position connue de la station de référence.

\section{\textsc{Les Applications du GPS}}
Les applications du système GPS sont multiples:

*	canevas de détails;

*	navigation de précision;

*	levés hydrographiques de grande précision;

*	contrôle d'altitudes (plan d'eau, marée,...);

*	relevés de type sismique;

*	surveillance de position et de mouvement (micro-géodésie);

*	levés de détails, cadastre et topographie;

*	implantation;

*	localisation des points;

*	mines, prospection;

*	stéréo-préparation.

\subsection*{16.7.1. Les Avantages de système GPS}
Les avantages de système GPS sont :

	- précision centimétrique;
	
	- indépendance des conditions atmosphériques;
	
	- nécessite un seul opérateur;
	
	- productivité active;
	
	- résultats exploitables en divers domaines:
	
*	photogrammétrie;

*	système d'information géographique (SIG);

*	génie civil (collection des eaux usées, drainage);

	- résultats en coordonnées $WGS84$ et locales.
	\section{\textsc{Almanach}}\index{Almanach}
L'almanach est l'ensemble de paramètres radiodiffusés par chaque satellite GPS,  permettant d'estimer la trajectoire et le comportement des horloges du satellite. Il est utilisé pour des prévisions à moyen terme. 

	Il est intéressant de regarder le nombre de satellites et le $PDOP$ précédemment défini avant de faire des observations pour éviter les périodes défavorables ($GDOP$ > 4, $PDOP$ >  3).
	
	La détermination des almanachs nécessite une observation d'environ deux heures pour les récupérer en leur totalité, et un logiciel de traitement. L'almanach peut être valable pour l'observation des autres zones mais dans un délai ne dépassant pas les 60 jours.
	
	Pour obtenir des résultats propres à chaque zone d'études, on doit fixer sa latitude, sa longitude et son altitude ainsi que la date du jour de l'observation. 
	
	Les produits sont les graphiques:
	
-	de la visibilité des satellites;

-	des coefficients $PDOP$ et $GDOP$;

-	de la configuration des satellites.

	




	
\clearpage
\chapter{\textit{\textbf{Bibliographie I}}}\label{biblio1}
\begin{enumerate}
 
\item \textbf{F.R. Helmert}. 1884. \textit{Die Mathematischen und Physikalischen Theorien der Höheren Geodäsie}; Vol 2, Leibzig, B.G Teubner(reprinted 1962).\index{\textbf{Helmert F.R.}}

\item \textbf{H. Poincaré}. 1905. Sur les Lignes g\'eod\'esiques des surfaces convexes. Transactions of the American Mathematical Society. n°6, pp. 237-274; \OE uvres 6, pp. 38-84.\index{\textbf{Poincaré H.}}

\item \textbf{F. Tisserand \& H. Andoyer}. 1912. \textit{Leçons de Cosmographie}. 6ème édition. Librairie Armand Colin. 396p.\index{\textbf{Tisserand F.}}\index{\textbf{Andoyer H.}}

\item \textbf{G. Julia}. 1955. Cours de l'Ecole Polytechnique. \textit{Cours de Géométrie Infinitésimale}. Cinquième Fascicule, Deuxième Partie: Théorie des Surfaces. Deuxième édition entièrement refondue. Editeur Gauthier-Villars. 141p.\index{\textbf{Julia G.}}
 
\item \textbf{W.A. Heiskanen \& H. Moritz}. 1967. \textit{Physical Geodesy}. Freeman, San Francisco. Reprint, 1979. Institute of Physical Geodesy, Technical University, Graz, Austria. 364p. \index{\textbf{Moritz H.}}\index{\textbf{Heiskanen W.A.}}  

\item \textbf{J. Dieudonné}. 1968. \textit{Calcul Infinitésimal}. 1ère édition. Collection Les Méthodes. Hermann,
Paris. 479p. \index{\textbf{Dieudonné J.}}

\item \textbf{A. Fontaine}. 1969. Rapport sur la Géodésie de la Tunisie. OTC.\index{\textbf{Fontaine A.}}

\item \textbf{C.A.C.G.G. (le Comité Associé Canadien de Géodésie et de Géophysique)}. 1973. Canadian Surveyor, Vol 27, n°3. 

\item \textbf{A. Marussi}. 1974. Africa and Modern Geodesy. Proceedings of the First Symposium in Geodesy in Africa, 14-19 January. Khartoum. Soudan.\index{\textbf{Marussi A.}} 

\item \textbf{B. Morando}. 1974. \textit{Mouvement d'un Satellite Artificiel de la Terre}. Gordon \& Breach, Paris, London et New York. 255p.\index{\textbf{Morando B.}}

\item \textbf{L. Bers}. 1977. Quasiconformal mappings, with applications to differential equations, function theory and topology. Bulletin of the American Mathematical Society, vol 83, n°6, pp. 1083-1100, 1977.\index{\textbf{Bers L.}}
    
\item \textbf{C. Fezzani}.\index{\textbf{Fezzani C.}} 1979. Analyse de la structure des réseaux astro-géodésiques tunisiens. Thèse de Docteur Ingénieur. Ecole Nationale des Sciences Géographiques. IGN France. 314p.\index{\textbf{Fezzani C.}}

\item \textbf{C. Boucher}. 1979a. Systèmes géodésiques de référence et autres problèmes géodésiques liés à la localisation en mer. Colloque national sur la localisation en mer. Brest, 1-5 octobre 1979. IGN France.\index{\textbf{Boucher C.}}

\item \textbf{C. Boucher}. 1979b. Les Transformations géométriques entre systèmes géodésiques. Rapport Technique RT/G n°3, SGNM, IGN France.

\item \textbf{P. Dombrowski}. 1979. 150 Years after Gauss " disquisitiones generales circa superficies curvas". Astérisque n°62. Publication de la Société Mathématique de France. 153p.\index{\textbf{Dombrowski P.}}

\item \textbf{H.M. Dufour}. 1979. Systèmes de références: Systèmes Projectifs. Conférence présentée au Colloque national sur la Localisation en mer. Brest, 1-5 octobre 1979. 27p.\index{\textbf{Dufour H.M.}}

\item \textbf{H. Monge}. 1979. L'Ellipsoïde de Clarke 1880. Bulletin d'Information de l'IGN France n°39-1979/3, pp. 37-38. \index{\textbf{Monge H.}}

\item \textbf{J. Lemenestrel}. 1980. \textit{Cours de Géodésie Elémentaire}, ENSG, IGN France.\index{\textbf{Lemenestrel J.}}

\item \textbf{J. Commiot}. 1980. \textit{Cours de Cartographie Mathématique}, ENSG, IGN France.\index{\textbf{Commiot J.}}

\item \textbf{A. Danjon}. 1980. \textit{Astronomie Générale: Astronomie sphérique et Eléments de mécanique céleste}. Seconde édition, revue et corrigée. Librairie Scientifique et Technique Albert Blanchard. 454p.\index{\textbf{Danjon A.}}

\item \textbf{B. Doubrovine, S. Novikov et A. Fomenko}. 1982. \textit{Géométrie Contemporaine: Méthodes et Application}. Première Partie: Géométrie des surfaces, des groupes de transformations et des champs. Edition Mir, Moscou. 438p. \index{\textbf{Doubrovine B.}} \index{\textbf{Novikov S.}} \index{\textbf{Fomenko A.}}

\item \textbf{M. Charfi}. 1983. Les Travaux de revalorisation de la Géodésie Tunisienne. OTC.\index{\textbf{Charfi M.}}

\item \textbf{P. Vani$\vec{\breve{c}}$ek \& E.J. Krakiwsky}. 1986. \textit{Geodesy: the Concepts}. North Holland Compagny. 2ème Edition. 697p.\index{\textbf{Vani$\vec{\breve{c}}$ek P.}}\index{\textbf{Krakiwsky E.J.}}

\item \textbf{Defense Mapping Agency(DMA)}. 1987a. DMA Technical Report 8350.2. Dept of Defense, World Geodetic System 1984, Its definition and relationships with local geodetic systems. 121p.\index{Defense Mapping Agency}
  
\item \textbf{Defense Mapping Agency}. 1987b. DMA Technical Report 8350.2-A. Part I: Methods, Techniques, and Data used in WGS84 development. Supplement to Department of Defense World Geodetic System 1984 Technical Report. 412p.

\item \textbf{Defense Mapping Agency}. 1987c. DMA Technical Report 8350.2-A. Part II: Parameters, formulas, and graphics for the practical application of WGS84. Supplement to Department of Defense World Geodetic System 1984 Technical Report. 775p.

\item \textbf{H. Moritz \& I.I. Mueller}. 1988. \textit{Earth Rotation :Theory and Observation}. Ungar Publishing Compagny. New York. 617p.\index{\textbf{Mueller I.I.}}

\item \textbf{B. Hofmann-Wellenhof, H. Lichtenegger et J. Collins}. 1992. \textit{Global Positioning System, Theory and Practice}. Springer-Verlag Wien New York. 382p.\index{\textbf{Hofmann-Wellenhof B.}}\index{\textbf{Lichtenegger H.}}\index{\textbf{Collins J.}}

\item \textbf{M. Bur$\vec{\breve{s}}$a \& K. P$\vec{\breve{e}\breve{c}}$}. 1993. \textit{Gravity Field and Dynamics of the Earth}. Springer-Verlag. 385p.\index{\textbf{Bur$\breve{s}$a M.}}\index{\textbf{P$\breve{e}\breve{c}$ K.}}

\item \textbf{Publications de l'Acad\'emie Bavaroise de G\'eod\'esie, n°58}. 1997. Symposium de la Commission EUREF, Sofia, 4-7 juin 1997. Extrait du rapport présenté par la Suisse. pp. 212-218.

\item \textbf{T. Soler}. 1998. A Compendium of transformation formulas useful in GPS work. Journal of Geodesy, Vol.72, n°7/8, pp. 482-490.\index{\textbf{Soler T.}}

\item \textbf{E.W. Grafarend}. 1998. Helmut Wolf: Das Wissenschaftliche Werk/ The Scientific Work. Publication de la Deutsche Geodätische Kommission bei der Bayerischen Akademie der Wissenschaften, Reihe A, Heft n°115, München 1998. 97p.	\index{\textbf{Grafarend E.W.}}\index{\textbf{Wolf H.}}

\item \textbf{P. Petersen}. 1998. \textit{Riemannian Geometry}. Graduate Texts in Mathematics, n°171. Springer-Verlag. 435p.\index{\textbf{Petersen P.}}

\item \textbf{A. Ben Hadj Salem}. 1999. La Mise \`a niveau de la G\'eod\'esie Tunisienne par l'unification Carthage2000 de ses r\'eseaux, pr\'esent\'ee au Premier Atelier Maghr\'ebin de G\'eod\'esie. Tunis, 18-20 mai 2000. Publi\'ee dans la revue G\'eo-Top de l'OTC num\'ero sp\'ecial mai 2000, pp. 6-16.

\item	\textbf{A. Ben Hadj Salem}. 1999. Note sur les comparaisons des azimuts g\'eod\'esiques du r\'eseau g\'eod\'esique primordial Tunisien. 5p.

\item \textbf{H. Moritz}. 2000. Geodetic Reference System 1980. Journal of Geodesy, vol. 74 n°1, pp. 128-134.\index{\textbf{Moritz H.}}

\item	\textbf{A. Ben Hadj Salem}. 2001. Contr\^ole de l'azimut d'orientation du syst\`eme g\'eod\'esique tunisien Carthage34.  Publi\'e dans la revue G\'eo-Top de l'OTC, n°4 f\'evrier 2002, pp. 63-68.

\item \textbf{H.A. Kastrup}. 2008. On the Advancements of Conformal Transformations and their Associated Symmetries in Geometry and Theoretical Physics. arxiv:http://arxiv.org/physics.hist-ph/0808.2730v1. www.arxiv.org.\index{\textbf{Kastrup H.A.}} 82p.

\item \textbf{Arr\^et\'e du ministre de la D\'efense nationale du 10 f\'evrier 2009}. 2009. Journal Officiel de la R\'epublique Tunisienne n°14 du 17 f\'evrier 2009.

\item \textbf{A.N. Pressley}. 2010. \textit{Elementary Differential Geometry}. Springer-Verlag Heidelberg. 395p.\index{\textbf{Pressley A.N.}}

\item \textbf{M. Lemmens}. 2011. \textit{Geo-information: Technologies, Applications and the Environment}. Series Geotechnologies and the Environment. Volume 5. Springer Netherlands. 349p. \index{\textbf{Lemmens M.}}
   
\item \textbf{E. Zeidler}. 2011. \textit{Quantum Field Theory III: Gauge Theory A Bridge between Mathematicians and Physicists}. Springer-Verlag Heidelberg. 1158p.\index{\textbf{Zeidler E.}}

\item \textbf{A. Ben Hadj Salem}. 2012. Selected Papers. Tome II. pp 165-172.

\item \textbf{A. Ben Hadj Salem}. 2013. Histoire de la Topographie en Tunisie : Les Coordonnées Origines Fuseaux. 11p.

\item \textbf{United Nations Commitee of Experts on Global Geospatial Information Management (UN-GGIM)}. 2015. La Résolution adoptée par l'Assemblée générale le 26 février 2015. www.ggim.un.org/docs/A\_RES\_69\_266\_F.pdf. 3p.

\end{enumerate}

\part{\textit{Eléments de la Théorie des Moindres Carrés}}

\chapter{\textit{\textbf{El\'ements Math\'ematiques pour La M\'ethode des Moindres Carr\'es }}}\lb{ELMMC}




	Dans ce chapitre, on pr\'esente les th\'eor\`emes math\'ematiques relatifs \`a l'application de la th\'eorie des moindres carr\'es.


\vspace{2mm}
\section{\textsc{D\'efinitions}}
\bdf
Soient $E$ et $F$ deux espaces vectoriels norm\'es de dimension finie, et $f$ une application de $E$ dans $F$. $f$ est d\'erivable en $a \in E$ s'il existe une application lin\'eaire $L$ de $E$ dans $F$ qui v\'erifie:$$ \forall \epsilon >\,0\,\,\exists \alpha>\,0\,\,/\, || u|| \leq \alpha \Rightarrow ||f(a+u)-f(a)-L(u)||\leq \epsilon$$
ou: $$ f(a+u)=f(a)+L(u)+||u||\epsilon(u),\quad \lim_{u\rightarrow 0}\epsilon(u)=0$$
\edf
On note :
\be
	L=f'(a)
\ee
$f'(a)\in \mathcal{L}(E,F)$ l'espace vectoriel des applications lin\'eaires de $E$ dans $F$.
\\

Si $f'$ est d\'erivable en $a$, on pose:
\be
	f"(a)=(f')'(a)
\ee
$f"(a)$ est une application lin\'eaire de $E$ dans $\mathcal{L}(E,F)$, donc $f"(a)\in \mathcal{L}(E,\mathcal{L}(E,F))$ qui est isomorphe \`a l'espace vectoriel $\mathcal{L}_2(E,F)$ des applications bilin\'eaires de $E$ dans $F$, c'est-\`a-dire des applications $B$ telles que:
\ba
	\forall (x,x_1,x_2)\in E^3,\forall (\alpha_1,\alpha_2)\in \BbR^2,\nonumber \\
	\left\{
\begin{array}{ll}
	B(\alpha_1x_1+\alpha_2x_2,x)=\alpha_1B(x_1,x)+\alpha_2B(x_2,x) \\
		B(x,\alpha_1x_1+\alpha_2x_2)=\alpha_1B(x,x_1)+\alpha_2B(x,x_2) 
\end{array}\right.
\ea
$f"(a)$ est une application bilin\'eaire sym\'etrique si:
\be
	\fbox{ $ f"(a)(h,k)=f"(a)(k,h) $}
\ee
\subsection*{18.1.1. Composition des d\'eriv\'ees premi\`eres}
\bdf
Soient $E,F$ et $G$ trois espaces vectoriels norm\'es de dimension finie, et $f$ une application de $E$ dans $F$ d\'erivable en $a$ et $g$ une application de $F$ dans $G$ d\'erivable en $b=f(a)$; alors $h=gof$ est d\'erivable en $a$ et :
\be
	h'(a)=(gof)'(a)=g'(b)of'(a)
\ee  
\edf       
\section{\textsc{Condition n\'ecessaire d'extremum local}}\index{Extremum local}
\bthm
Si une fonction $f$ d'un espace vectoriel norm\'e $E$ dans $\BbR$ admet un extremum local en un point $a$ o\`u elle est d\'erivable, alors:
\be
	f'(a)=0
\ee        
\ethm
\textbf{D\'emonstration:} Soit $I\subset \BbR$ un intervalle ouvert contenant 0. On d\'efinit la fonction $\varphi:I\rightarrow \BbR,\,\varphi(t)=f(a+tu)$, on a alors $\varphi'(t)=f'(a+tu)u$. En particulier $\varphi'(0)=f'(a)u$.

Si $a$ est un minimum relatif, soit:
\ba           
	&f(a)\leq f(a+tu) \,\forall t\in I\Rightarrow \varphi(0)\leq \varphi(t)\Rightarrow 0\geq lim_{t\rightarrow 0^-}\frac{\varphi(t)-\varphi(0)}{t}=\nonumber& \\&
	\varphi'(0)=lim_{t\rightarrow 0^+}\frac{\varphi(t)-\varphi(0)}{t}\geq 0&
\ea
car $f$ est d\'erivable en $a$ ou encore $\varphi$ est d\'erivable en $0 \in I$, on a donc:
\be          
	\varphi'(0)=0\Rightarrow f'(a)u=0 \,\,\forall u \Rightarrow f'(a)=0
\ee
\section{\textsc{Formulation du probl\`eme des Moindres Carr\'es - le Cas Lin\'eaire}}
Soient $B$ une matrice (m,n) r\'eelle de rang $n$ et $L$ un vecteur de $\BbR^m$. On cherche un vecteur $X$ de $\BbR^n$ tel que:
\be          
	||B.X-L||=\min_{\zeta \in R^n}||B.\zeta-L||
\ee          
C'est-\`a-dire approchant au mieux la solution du syst\`eme:
\be
	B.X=L
\ee          
 Si l'on pose:
\be
	V=B.X-L
\ee          
$V$ est appel\'e vecteur r\'esidu. On cherche donc un vecteur $X$ qui minimise $||V||$ ou $||V||^2$ c'est-\`a-dire qui minimise la somme des carr\'es des r\'esidus d'o\`u le terme " $X$ est solution des moindres carr\'es (Least-Squares Solution)" .
\\

Pour cela, on introduit la forme quadratique $f$ d\'efinie par:
\be
	f(U)=\frac{1}{2}||B.U-L||^2-\frac{1}{2}||L||^2 \label{mq1}
\ee        
avec $U$ un vecteur de $\BbR^n$. $f$ est donc une application de $\BbR^n\longrightarrow \,\BbR$.
On introduit un produit scalaire:
\be
	(\zeta,\xi)=\zeta^T.\xi=\sum^{i=n}_{i=1}\zeta_i.\xi_i \label{m2}
\ee        
 D'o\`u la norme du vecteur $\zeta=(\zeta_1,\zeta_2,...,\zeta_n)^T$:
 \be
	\|\zeta \|^2=(\zeta,\zeta)=\zeta^T.\zeta =\sum_{i=1}^n.\zeta_i^2 \label{m3}
\ee          
dans l'espace vectoriel $\BbR^n$.

D\'eveloppons l'\'equation (\ref{mq1}), on obtient:
\ba
&f(U)=\frac{1}{2}||B.U-L||^2-\frac{1}{2}||L||^2= \frac{1}{2}((B.U-L,B.U-L)-(L,L))=\nonumber &\\&	\frac{1}{2}.\left((BU-L)^T.(BU-L)-L^T.L\right)=\nonumber &\\&\frac{1}{2}.(B^T.B.U,U)-(B^T.L,U)=\frac{1}{2}U^T.B^T.B.U-L^T.B.U \label{mq2}&
\ea          
 Le probl\`eme consiste \`a minimiser $f(U)$ sur $\BbR^n$, donc chercher un minimum absolu ( et non pas seulement relatif). On a \'etabli une condition n\'ecessaire de minimum relatif: $f'(U)=0$. La recherche de la solution compl\`ete fait appel aux propri\'et\'es des \underline{fonctions convexes}.\index{Fonctions convexes}
 
\subsection*{18.4.1. Calcul des d\'eriv\'ees premi\`ere et seconde de $f$}
On commence par le calcul de $f'$. D'o\`u:
\ba
&	f(U+h)=\frac{1}{2}.((U+h)^T.B^T.B.(U+h)-L^T.B.(U+h)\nonumber &\\&
	= \frac{1}{2}.(U^T.B^T.B.U+2h^T.B^T.B.U +h^T.B^T.B.h)\nonumber & \\ &-L^T.B.U - L^T.B.h&\label{mq3} 
	 \ea          
	 Par suite, en posant:
\ba
	h^T.B^T.B.h=2||h||.\epsilon(h) \label{mq4} \\
	et \quad 	\lim_{h \rightarrow 0}\epsilon(h)=0
\ea           
Alors:
\be
	f(U+h)-f(U)=h^T.B^T.B.U-h^T.B^T.L+||h||\epsilon(h) \label{mq5}
\ee          
Or l'\'equation pr\'ec\'edente (\ref{mq5}) s'\'ecrit:
\be
		f(U+h)-f(U)=(B^T.B.U-B^T.L,h)+||h||\epsilon(h) \label{mq6}
\ee           
Or en posant :
\be
	Q(h)=(B^T.B.U-B^T.L,h)
\ee             
Il est imm\'ediat que $Q$ est une application lin\'eaire de $\BbR^n$ dans $\BbR$ qui \`a $h$ associe le r\'eel $(B^T.B.U-B^T.L,h)$. D'apr\`es la d\'efinition 18.1, on a donc:
\ba
&	f(U+h)-f(U)=(B^T.B.U-B^T.L,h)+||h||\epsilon(h)=(f'(U),h)+||h||\epsilon(h) \nonumber & \\
 &\mbox{avec}\quad f'(U)=B^T.B.U-B^T.L \quad \mbox{et}\quad \lim_{h \rightarrow 0}\epsilon(h)=0  & \label{mq6}
\ea         
Il s'ensuit que:
\be 
	f"(U)(h,k)=f"(U)(k,h)=h^T.B^T.B.k \label{mq7}
\ee             
\bthm
(\textbf{de Taylor\footnote{\textbf{Brook Taylor} (1685-1731): mathématicien anglais.}-Young\footnote{\textbf{William Henry Young} (1863-1942): mathématicien anglais.}})\index{Théorème de Taylor-Young}\index{\textbf{Taylor B.}} Soit $f$ une fonction de $\BbR^n$ dans $\BbR$. Si $f$ est d\'erivable en $a \in \BbR^n$, alors:\index{\textbf{Taylor B.}}\index{\textbf{Young W.H.}}
\be
	f(a+h)=f(a)+f'(a)h+||h||\epsilon(h),\quad \lim_{h\rightarrow 0}\epsilon(h)=0
\ee         
Si $f$ est deux fois d\'erivable en $a$, alors:
\be
	f(a+h)=f(a)+f'(a)h+\frac{1}{2}f"(a)(h,h)+||h||^2\epsilon(h),\quad \lim_{h\rightarrow 0}\epsilon(h)=0
\ee        
\ethm
\bthm
(\textbf{de Taylor-Maclaurin\footnote{\textbf{Colin Maclaurin} (1698-1746): mathématicien écossais.}})\index{Théorème de Taylor-Maclaurin}\index{\textbf{Maclaurin C.}} Soit $f$ une fonction de $ \BbR^n$ dans $ \BbR$. Si $f$ est d\'erivable en $ ]a,a+h[$, alors:
\be          
	f(a+h)=f(a)+f'(a+\theta h)h,\quad 0<\theta<1
\ee
Si $f$ est deux fois d\'erivable sur $]a,a+h[$, alors:
\be           
	f(a+h)=f(a)+f'(a)h+\frac{1}{2}f"(a+\theta h)(h,h)+||h||^2\epsilon(h),\quad 0<\theta<1
\ee
\ethm
\section{\textsc{Convexit\'e}}
\bdf
Une partie $U$ de $   \BbR^n$ est convexe si:
\be
	\fbox{ $ \forall (u,v)\in U^2,\,\,\forall \theta \in [0,1],\,\,\theta u+(1-\theta)v \in U $}
\ee        
\edf
\bdf
Une fonction $f$: $U$ convexe $\subset \BbR^n \longrightarrow \BbR$ est convexe si:
\be
	\fbox{ $ \forall (u,v)\in U^2,\,\,\forall \theta \in [0,1],\,\,f[\theta u+(1-\theta)v]\leq \theta f(u)+(1-\theta)f(v) $ }
\ee         
\edf
\bthm
(\textbf{Condition n\'ec\'essaire de minimum relatif sur un ensemble convexe}) Soient $f$: $\Omega$ ouvert $\subset    \BbR^n \longrightarrow  \BbR$ et $U$ convexe $\subset \Omega$. Si $f$ est d\'erivable en $u \in U$ et si elle admet en $u$ un minimum relatif par rapport \`a $U$, alors:
\be
	f'(u)(v-u)\geq 0\quad \forall v \in U
\ee          
 \ethm
\textbf{D\'emonstration:} Soient $v \in U$  et $w=v-u$. $U$ est convexe $\Rightarrow \forall \theta \in[0,1],\,u+\theta w=(1-\theta)u+\theta v \in U$. $f(u+\theta w)=f(u)+f'(u)\theta w+||\theta w|| \epsilon_1(\theta w)$ d'apr\`es la formule de Taylor-Young, soit : $f(u+\theta w)=f(u)+\theta[f'(u)w+\epsilon(\theta w)]$, $\lim_{\theta \rightarrow 0}\epsilon(\theta)=0$. 

$u$ minimum relatif $\Rightarrow 0\leq f(u+\theta w)-f(u)=\theta [f'(u)w+\epsilon (\theta)] \Rightarrow f'(u)w +\epsilon(\theta) \geq 0 \Rightarrow f'(u)w\geq 0$ d\`es que $\epsilon(\theta)$ est suffisamment petit. D'o\`u: $$ f'(u)(v-u) \geq 0$$  
\bthm
(\textbf{Convexit\'e et d\'erivabilit\'e premi\`ere}) Soient $f$: $\Omega$ ouvert $\subset \BbR^n \longrightarrow    \BbR$, d\'erivable dans $\Omega$ et $U$ convexe $\subset \Omega$. La fonction $f$ est convexe si et seulement si:
\be
	f(v)\geq f(u)+f'(u)(v-u)\quad \forall (u,v)\in U^2 \label{el31}
\ee       
\ethm  
\textbf{D\'emonstration:} Soient $(u,v) \in  U^2, u\neq v, \theta \in ]0,1[$. $$f \,\mbox{est convexe} \Rightarrow f(u+\theta(v-u))\leq (1-\theta)f(u)+\theta f(v)$$
$$\quad \Rightarrow \frac{f(u+\theta(v-u))-f(u)}{\theta}\leq f(v)-f(u) $$
Soit $h=\theta (v-u)\Rightarrow f(u+\theta(v-u))=f(u+h) = f(u)+f'(u)h+||h||\epsilon(h)=f(u)+\theta f'(u)(v-u)+||h||\epsilon(h)$
$$\Rightarrow\lim_{\theta \rightarrow 0}\frac{f(u+\theta (v-u))-f(u)}{\theta}=f'(u)(v-u) \Rightarrow f'(u)(v-u)\leq f(v)-f(u)$$
C'est l'\'equation (\ref{el31}).

Reciproquement: on a:
\be
	f(v)\geq f(u)+f'(u)(v-u)\quad \forall (u,v)\in U^2 \nonumber
\ee        
On peut \'ecrire:
\be
	\left\{ \begin{array}{lll}
	f(v)\geq f(v+\theta(u-v))+f'(v+\theta(u-v))(v-v-\theta(u-v)) \\
\mbox{on a pris}\,\, u=v+\theta(u-v) \\
f(u)\geq f(v+\theta(u-v))+f'(v+\theta(u-v))(u-v-\theta(u-v)) \\ 
\mbox{on a pris}\,\, v=v+\theta(u-v)\quad
		\mbox{avec}\quad u\neq v,0<\theta<1
\end{array}  \right. 
\ee
En utilisant que $f'(a)$ est une application lin\'eaire, on a:
\be
	\left\{
\begin{array}{ll}
	f(v)\geq f(v+\theta(u-v))-\theta f'(v+\theta(u-v))(u-v) \\
		f(u)\geq f(v+\theta(u-v))+(1-\theta)f'(v+\theta(u-v))(u-v) 
		\end{array} \right.
\ee
En multipliant la premi\`ere \'equation par $(1-\theta)$ et la deuxi\`eme par $\theta$ et faisant la somme on obtient:
$$(1-\theta)f(v)+\theta f(u)\geq f(v+\theta(u-v))=f((1-\theta)v+\theta u)\Rightarrow f \mbox{ est convexe}$$
\bthm
(\textbf{Convexit\'e et d\'erivabilit\'e seconde}) Soient $f$: $\Omega$ ouvert $\subset \BbR^n\rightarrow \BbR$, deux fois d\'erivable dans $\Omega$ et $U$ convexe $\subset \Omega$. La fonction $f$ est convexe sur $U$  si et seulement si:
\be
	f"(u)(v-u,v-u) \geq 0\quad \forall (u,v)\in U^2 \label{el32}
\ee 
\ethm          
\textbf{D\'emonstration:} La formule de Taylor-Mac Laurin s'\'ecrit:
$$	f(a+h)=f(a)+f'(a)h+\frac{1}{2}f"(a+\theta h)(h,h),\quad 0<\theta<1 $$
En l'appliquant \`a $a=u$ et $h=v-u$, on obtient: 
$$	f(v)-f(u)-f'(u)(v-u)=\frac{1}{2}f"((1-\theta)u+\theta v)(v-u,v-u)$$
Si $w=(1-\theta)u+\theta v \in U$ convexe, on a : $v-u=(1-\theta)(v-u)$, soit $\rho=1/(1-\theta) >0$, il vient:
$$  	f(v)-f(u)-f'(u)(v-u)=\frac{\rho^2}{2}f"(w)(v-w,v-w),\quad w\in ]u,v[$$
Si $f"(w)(v-w,v-w)\geq 0 ,\forall (w,v) \in U^2$, alors $f(v)-f(u)-f'(u)(v-u)\geq 0$ donc $f$ est convexe d'apr\`es le th\'eor\`eme 18.5 'convexit\'e et d\'erivabilit\'e premi\`ere'.

Reciproquement: Soit $g: \Omega \rightarrow \BbR,\, g(v)=f(v)-f'(u)v$, on a alors:
$$g(v)-g(u)=f(v)-f'(u)v-f(u)+f'(u)u=f(v)-f(u)-f'(u)(v-u)\geq 0$$ 
car f est convexe, donc $g$ admet un minimum sur $U$ en $u$. On a:
$$g'(v)=f'(v)-f'(u)\Rightarrow g'(u)=f'(u)-f'(u)=0 \quad \mbox{et} \,\,g"(v)=f"(v)$$
La formule de Taylor-Young s'écrit:
$$	\varphi(a+h)=\varphi(a)+\varphi'(a)h+\frac{1}{2}\varphi"(a)(h,h)+||h||^2\epsilon(h),\quad \lim_{h\rightarrow 0}\epsilon(h)=0 $$
s'\'ecrit pour $\varphi=g, a=u,h=tw$ avec $w=v-u$ d'o\`u $a+h=u+t(v-u) \in U$ pour $t \in [0,1]$:
\ba
& g(u+tw)=\ds g(u)+g'(u)tw+\frac{1}{2}g"(u)(tw,tw)+||tw||^2\epsilon_1(tw)=\nonumber &\\
&\ds g(u)+\frac{t^2}{2}f"(u)(w,w)+ \frac{t^2}{2}\epsilon(t) \nonumber &
\ea
Comme $g$ admet un minimum en $u$, on a alors:
$$g(u+tw)-g(u)\geq 0 \Rightarrow \frac{t^2}{2}\left(f"(u)(w,w)+\epsilon(t)\right)\geq 0 $$
comme $\lim_{t\rightarrow 0}\epsilon(t)=0$ et pour $t$ suffisamment petit, on a $f"(u)(w,w)\geq 0$
donc:
$$f"(u)(v-u,v-u)\geq 0$$
\bthm
(\textbf{Convexit\'e, d\'erivabilit\'e premi\`ere et minimum global}) Soient $U$ convexe $\subset    \BbR^n$ et $f$ une fonction convexe de $U$ dans $   \BbR$.

1. Si $f$ admet un minimum relatif en un point de $U$, elle y admet un minimum global.

2. Si $f$ est d\'efinie sur un ouvert $\Omega$ contenant $U$ et d\'erivable en $u \in U$, alors $f$ admet en $u$ un minimum global sur $U$ si et seulement si:
$$	\forall v \in U,\,\,f'(v)(v-u)\geq 0 $$ 
3. Si $U$ est ouvert, la condition pr\'ec\'edente \'equivaut \`a $	f'(u)=0 $.
\ethm

\textbf{D\'emonstration}

1. Soient $v\in U$ et $w=v-u$. $f$ est convexe $\Rightarrow \forall \theta \in[0,1] \quad f(v+\theta w)\leq (1-\theta)f(u)+\theta f(v) \Rightarrow f(v+\theta w)-f(u)\leq\theta (f(v)-f(u))$.

 $u$ minimum relatif $\Rightarrow \exists \theta_0 > 0\,\, /\,\, 0\leq f(v+\theta_0 w)-f(u) \Rightarrow 0\leq \theta_0 (f(v)-f(u))\Rightarrow f(v)\geq f(u)$.
\\

2. La condition est n\'ecessaire (minimum relatif sur un convexe). Reciproquement: soit $f$ convexe sur $U$ convexe, $f$ d\'erivable en $u$, telle que $f'(u)(v-u) \geq 0, \,\,\,\forall v\in U$.

 $f$ convexe $\Rightarrow f(v)\geq f(u)+f'(u)(v-u)\quad \forall (u,v)\in U^2$ (convexit\'e et d\'erivabilit\'e premi\`ere) donc $\forall v, f(v)-f(u)\geq f'(u)(v-u)\geq 0\Rightarrow u$ est un minimum. 
\\

3. $f'(u)=0\Rightarrow f'(u)(v-u)=0\,\, \forall v\in U\Rightarrow f'(u)(v-u)\geq 0, \forall v\in U$.

Reciproquement: $U$ ouvert, $\exists$ boule ouverte de $U$ contenant $u$ tel que :

$\forall h\in    \BbR^n,\,\, \exists \alpha \in    \BbR / \epsilon < \alpha \Rightarrow (u+\epsilon h,u-\epsilon h)\in U^2 \longrightarrow:$
$$ \left\{
\begin{array}{ll}
	\epsilon f'(u)h \geq 0 \\
	-\epsilon f'(u)\geq 0
\end{array}\right. \Rightarrow f'(u)h=0\quad \forall h \Rightarrow f'(u)=0$$
\section{\textsc{Application au probl\`eme des moindres carr\'es}}
  On a:
\ba             
	f(U)=\frac{1}{2}(B^T.BU,U)-(B^T.L)=\frac{1}{2}||BU-L||^2-\frac{1}{2}||L||^2 \\
	f'(U)=B^T.BU-B^T.L \\
	f"(U)(h,h)=h^T.B^T.B.h>0  
\ea
car $B^T.B$ est une matrice sym\'etrique définie positive.
Donc:
\ba
&	\forall (U,W)\in    \BbR^n \times    \BbR^n,\quad f"(U)(W-U,W-U)=B^T.B(W-U,W-U)\geq 0 \nonumber & \\
&	\Rightarrow \mbox{(convexit\'e et d\'erivabilit\'e seconde) f est convexe sur}\,\,   \BbR^n &
\ea           
D'apr\`es le th\'eor\`eme 18.7 'convexit\'e, d\'erivabilit\'e premi\`ere et minimum global', la recherche du minimum de $f$ coïncide alors avec l'ensemble des solutions de l'\'equation $ f'(U)=0 $. Soit:
\be
	B^T.B.U=B^T.L
\ee   
D'où respectivement la solution des moindres carr\'es et le vecteur résidu:
\ba
&\fbox{ $ \overline{U}=(B^TB)^{-1}.B^T.L $}&\\
&\fbox{ $V=B.\overline{U}-L=(B.(B^TB)^{-1}.B^T-I).L $}&
\ea
\section{\textsc{Exercices et Problèmes}}
\bex
 Soient $U$ un ouvert convexe d'un espace de Banach\footnote{\textbf{Stefan Banach} (1892-1945): mathématicien polonais.} réel $E$ c'est-à-dire un espace vectoriel normé complet sur $\BbR$ et $f$ une fonction à valeurs réelles, différentiable et convexe dans $U$. Montrer que si $f'(x_0)=0$ en un point $x_0\in U$, alors $f$ a un minimum absolu en $x_0$.\index{\textbf{Banach S.}}
\eex
\bex
 Montrer que dans un espace de Banach réel $E$, la fonction $f=\|\,. \|^2$ est strictement convexe, c'est-à-dire, $\forall \, \alpha \in ]0,1[,\,\,f(\alpha x+(1-\alpha) y)< \alpha f(x)+(1-\alpha)f(y), $ pour tout couple $(x,y) \in E^2$.

Aide: utiliser l'identité remarquable:
$$\|\alpha x + (1-\alpha )y\|^2 = \alpha \|x\|^2 + (1-\alpha )\|y\|^2-\alpha (1-\alpha )\|x- y\|^2$$
\eex
\bex
 On note $F$ une surface de $\BbR^3$ définie par la représentation paramétrique:
$$\textbf{\textit{OM}}=(a_1(u,v),a_2(u,v),a_3(u,v))^T $$
où $u$ et $v$ sont deux paramètres réels. On se donne un point $P(x,y,z)\in \BbR^3$.

1. Donner une condition géométrique portant sur le plan tangent à $F$ au point $M_0(u_0,v_0)$ pour que la différentielle de la fonction $(u,v)\longrightarrow \varphi(u,v)=\left\|\textbf{\textit{OP}}-\textbf{\textit{OM}}(u,v)\right\|^2$ soit nulle en $M_0(u_0,v_0)$.
\eex
\bpb
Soient $U$ un ouvert convexe d'un espace de Banach réel $E$ et $f$ une application différentielle de $U$ dans $\BbR$.

1. Montrer que $f$ est convexe dans $U$ si et seulement si:
$$ f(x)\geq f(x_0)+f'(x_0)(x-x_0)$$
pour tout couple de points $x,x_0\in U$.

2. On suppose $E=\BbR^n$ et $f$ de classe $C^2$ soit deux fois différentiable et $f"$ continue; pour $x\in U$, soit $\varphi_x$ la forme quadratique  définie par :
$$\varphi_x(h)= \sum_{i,j=1}^n\frac{\partial ^2 f}{\partial x_i \partial x_j}(x)h_ih_j,\quad \,\,\,\,h=(h_1,h_2,...,h_n)\in \BbR^n $$
Montrer que $f$ est convexe dans $U$ si et seulement si $\varphi_x$ est positive pour tout $x\in U$ soit $\varphi_x(h)\geq 0$ pour $x\in U$ et $h\in \BbR^n$. 
\epb

\chapter{\textit{\textbf{El\'ements de la M\'ethode des Moindres Carr\'es}}}
\section{\textsc{Introduction}}
\subsection*{19.1.1. D\'efinition du Probl\`eme}
On veut d\'eterminer $r$ grandeurs scalaires inconnues $G_1, G_2,..., G_r$ \`a l'aide de $n$ grandeurs observ\'ees distinctes ou non des pr\'ec\'edentes mais qui leur sont li\'ees g\'eom\'etriquement. Les grandeurs $G_i$ peuvent \^etre des coordonn\'ees de points, des altitudes ou des constantes d'orientation.
\\

On dispose de $n$ mesures ou observations: $l_1, l_2, ..., l_n$. Le nombre des mesures est g\'en\'eralement surabondant par rapport au nombre $r$ des inconnues \`a d\'eterminer. Le but est de faire la compensation ou ajustement des mesures (en anglais : adjustment) et d'estimer au mieux les inconnues.
\\

Les grandeurs $G_i$ \`a d\'eterminer peuvent \^etre ou non li\'ees entre elles; s'il y'a $p$ relations entre les grandeurs, on dit que le nombre de degr\'es de libert\'e de l'ensemble est $r-p$. Par suite, il y aura \`a d\'eterminer $r-p$ inconuues $G_k$. Si le nombre des observations $n$ est sup\'erieur \`a $r-p$, c'est une condition n\'ecessaire pour d\'eterminer les $r-p$ inconnues, mais elle n'est pas suffisante. Il faut trouver parmi les $n$ observations un groupe de mesures qui permet de calculer les inconnues.
\section{\textsc{Les Mesures}}
Les mesures $l_1, l_2, ..., l_n$ \'ecrites sous forme matricielle:
\begin{equation}  
l=	\begin{pmatrix}{
		l_1 \cr  
	l_2 \cr 
	\vdots  \cr 
	l_n }
\end{pmatrix}
\end{equation} 
- peuvent \^etre directes, dans ce cas, ind\'ependantes, outre leur valeur, on connaît leur exactitude donn\'ee par leur moyenne quadratique $\sigma _i$, dans ce cas, la matrice variance du vecteur $\textbf{\textit{l}}$ est alors diagonale.
\be
\Gamma _ l =\begin{pmatrix}{
	 \sigma^2_1 &   0         & 0 & 0          \cr
	         0  &  \sigma^2_2 & 0 & 0          \cr
	         \cdots  & \cdots     & \ddots & \cdots      \cr
	         0  &  0          & 0 & \sigma^2_n }
\end{pmatrix}
\ee
et les $l_i$  $\in$ $\m N(E(l), \sigma): \m N$ loi normale $l$ $\,$ $\in$ $\m N(\dot{L}_i, \sqrt{\Gamma _l})$, où $\dot{L}$ désigne la valeur réelle inconnue avec $E(l)=\dot{L}$.

\begin{figure}
	\centering
		\includegraphics{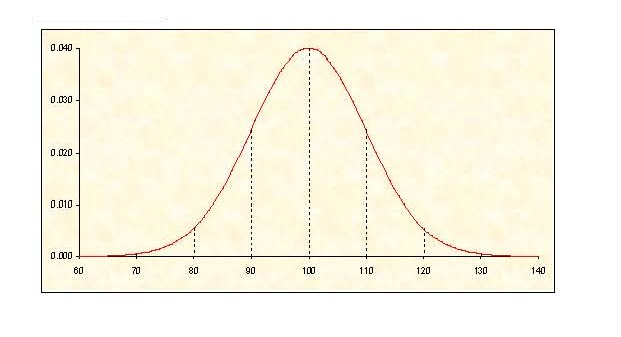}
	\caption{La courbe de Gauss\index{Courbe de Gauss}}
	\label{fig:gauss}
\end{figure}

- peuvent \^etre indirectes, alors en g\'en\'eral corr\'el\'ees donc la matrice $ \Gamma _ l$ n'est pas diagonale:

\begin{equation}
\Gamma _ l =\begin{pmatrix}{
	 \sigma^2_1 &   \sigma_{12}         & \cdots & \sigma_{1n}         \cr
	         \sigma_{21}  &  \sigma^2_2 & \cdots & \sigma_{2n}          \cr
	         \cdots  & \cdots     & \ddots & \cdots      \cr
	         \sigma_{n1}  &  \cdots   & 0 & \sigma^2_n }
\end{pmatrix}
\end{equation}
Le vecteur $l$ sera consid\'er\'e comme normal et sa densit\'e de probabilit\'e \index{Densité de probabilité} sera la fonction:
\be
\m P=\m P(l)=\frac{1}{\sqrt{(2\pi)^n\left|\mbox{D\'et}(\Gamma_l)\right|}}e^{\ds -\frac{1}{2}(l-\dot{L})^T \Gamma^{-1}_l(l-\dot{L})}
\ee          
où Dét$(\Gamma)$ désigne le déterminant de la matrice $\Gamma$. Quand $l$ est une variable scalaire, le graphique de $\m P (l)$ est donné par la courbe de Gauss (\textbf{Fig. \ref{fig:gauss}}). \index{\textbf{Gauss C.F.}}
\\

Le but qu'on se propose est alors:

* d'estimer au mieux \`a partir de la donn\'ee de $n$ mesures la valeur des $r$ grandeurs inconnues.

 On appelle :
\be
\overline{X}=	\begin{pmatrix}{
	\overline{X}_1 \cr 
	\overline{X}_2 \cr 
	\vdots  \cr 
	\overline{X}_r }
\end{pmatrix}
\ee
le vecteur des estimateurs correspondants.

* de trouver le meilleur estimateur de la matrice variance $\Gamma_{\overline{X}}$ de $ \overline{X} $ qui chiffrera la pr\'ecision de $\overline{X}$ et d'\'evaluer l'exactitude de $\overline{X}$. 
\newpage
\section{\textsc{Poids - Matrice de Poids - Variance de la Mesure de Poids Unitaire}}
Deux cas se posent:

- on conna\^it $\Gamma_l$ : mesures directes et ind\'ependantes;

- on conna\^it $\Gamma_l$ \`a un facteur pr\`es que la th\'eorie permettra de d\'eterminer.

 On pose alors:
\be
\Gamma_l	=\sigma_0^2 P^{-1} \label{mat6}
\ee               
avec:

- $P$ matrice de poids (fix\'ee avant les calculs);

- $\sigma^2_0$ variance unitaire scalaire inconnue \`a d\'eterminer.\index{Variance unitaire}
\\

La matrice des poids doit \^etre choisie comme inversement proportionnelle \`a la matrice variance des erreurs de mesures. En effet de (\ref{mat6}), on a:
\be
\fbox{ $ P=\sigma^2_0\Gamma_l^{-1}=\left(\ds \frac{\Gamma_l}{\sigma_0^2}\right)^{-1} $}
	\ee           
En particulier, si les mesures sont ind\'ependantes, les poids doivent \^etre choisis comme inversement proportionnels aux carr\'es des erreurs moyennes quadratiques.\index{Erreurs moyennes quadratiques}
\\

G\'en\'eralement si $\Gamma_l $ est connue, alors $ P = \Gamma_l^{-1}$, donc on consid\`ere que $\sigma^2_0 =1$. Dans le cas de mesures directes, la matrice $P=(p_i)$ est diagonale avec:
\be
\fbox{ $	p_i=\ds \frac{1}{\sigma^2_i} $}
\ee               
Une comparaison de la variance unitaire \textit{a priori} $\sigma^2_0 $ et la valeur estim\'ee $ s^2_0$ permettra de d\'eceler les incoh\'erences dans les donn\'ees relatives \`a l'exactitude des mesures.


\section{\textsc{Mod\`eles Fonctionnel et Stochastique}}
Ces mod\`eles d\'ecrivent les ph\'enom\`enes physiques.
\subsection*{19.4.1. Le Mod\`ele Fonctionnel}
Le mod\`ele fonctionnel d\'ecrit la relation entre les observables, les param\`etres et les relations \`a utiliser. Qu'est-ce qu'il faut choisir comme mod\`ele fonctionnel pour d\'ecrire un ph\'enom\`ene physique?
\\

\textbf{Hypoth\`ese}: \`a une grandeur, on peut associer sa valeur maximale vraie. Soit: 
$$ 	l\Longrightarrow \dot{l} $$
Soit $_{n}\dot{L}_1$ le vecteur des valeurs nominales des observables. S'il y a $n-r$ relations ind\'ependantes donn\'ees par l'exp\'erience, alors $r$ est appel\'e le nombre de \un{degr\'es de libert\'e} \index{Degré de liberté} du mod\`ele fonctionnel. 
\\
 
\textbf{Exemple 1}: d\'eterminer les angles $A, B, C$ d'un triangle plan. On a:
$$ \dot{A}+\dot{B}+\dot{C}=\pi $$
Donc $r=3-1=2$.
\\

\textbf{Exemple 2}: d\'eterminer un c\^ot\'e d'un triangle. On a 4 observations (3 angles et un c\^ot\'e) et une relation ind\'ependante (la somme des angles d'un triangle plan vaut $\pi$) d'o\`u $r=4-1=3$. 
\\

\textbf{Exemple 3}: d\'eterminer les 3 angles et les 3 c\^ot\'es d'un triangle plan. On a 3 relations ind\'ependantes:
\ba
	\dot{A}+\dot{B}+\dot{C}=\pi \nonumber \\
	\dot{a}^2=\dot{b}^2+\dot{c}^2-2\dot{b}\dot{c}cos\dot{A} \nonumber \\
	\dot{b}^2=\dot{a}^2+\dot{c}^2-2\dot{a}\dot{c}cos\dot{B} \nonumber
\ea
Le vecteur des observables est:
\be
	_6\dot{L}_1=\begin{pmatrix}{
	\dot{A} \cr
	\dot{B} \cr
	\dot{C} \cr
	\dot{a} \cr
	\dot{b} \cr
	\dot{c} }
\end{pmatrix} \label{aq8}
\ee
La d\'etermination de l'\'equation (\ref{aq8}) se fait par la m\'ethode des \'equations de condition c'est-\`a-dire que le vecteur $\dot{L}$ v\'erifie:
\begin{equation}
	\Phi(\dot{L})=0
\end{equation}
avec $n-r$ relations ind\'ependantes et $\Phi$ une certaine fonction.
\subsection*{19.4.2. La M\'ethode des Equations d'observations}
Avec cette m\'ethode, on d\'ecrit le mod\`ele fonctionnel par des param\`etres commodes. Soit $_{p}\dot{X}_1$ le vecteur des param\`etres avec $p\geq r$ o\`u $r$ est le degr\'e de libert\'e du mod\`ele. On cherche alors \`a exprimer les observables $\dot{L}$ en fonction des param\`etres du mod\`ele, soit:
\begin{equation}
	_{n}\dot{L}_1=\Phi(\dot{X})\Longleftrightarrow \Psi(\dot{X},\dot{L})=0
\end{equation}
Les \'equations d'observations sont si $p \geq r$ :
\ba
	_{n}\dot{L}_1=\Phi(\dot{X}) \\
	\Xi(\dot{X})=0 \label{aq12}
\ea
L'\'equation (\ref{aq12}) repr\'esente les $(p-r)$ relations ind\'ependantes dites \textbf{\'equations de condition}.
\\

\textbf{Exemple}: D\'etermination des coordonn\'ees des sommets d'un triangle plan $A_1A_2A_3$ \`a partir de l'observation des trois angles aux sommets $A_1, A_2, A_3$. On choisit comme param\`etres $\dot{x}_1, \dot{y}_1, \dot{x}_2, \dot{y}_2, \dot{x}_3, \dot{y}_3,$ donc $\dot{X}=\,_6\dot{X}_1$ s'\'ecrit:
\be
	\dot{X}=\,_6\dot{X}_1=\begin{pmatrix}{
	\dot{x}_1 \cr
	\dot{y}_1 \cr
	\dot{x}_2 \cr
	\dot{y}_2 \cr
	\dot{x}_3 \cr
	\dot{y}_3 }
\end{pmatrix} \label{aq10}
\ee
Les observables sont les angles $\dot{A}_1, \dot{A}_2, \dot{A}_3$, on a alors:
\ba
	\dot{A}_1=\varphi_1(\dot{X}) \nonumber \\
	\dot{A}_2=\varphi_2(\dot{X}) \nonumber \\
	\dot{A}_3=\varphi_3(\dot{X}) \nonumber
\ea
et :
\be
	_3\dot{L}_1=\dot{L}=\begin{pmatrix}{
	\dot{A}_1 \cr
	\dot{A}_2 \cr
	\dot{A}_3 }
\end{pmatrix}=\Phi(\dot{X})=\begin{pmatrix}{
	\varphi_1 \cr
	\varphi_2 \cr
	\varphi_3 }
\end{pmatrix} \label{aq17}
\ee 
Les \'equations de condition sont quatre \`a savoir:
\ba
	\dot{x}_1=0 \nonumber \\
	\dot{y}_1=0 \nonumber \\
	\dot{y}_3=0 \nonumber \\
	\dot{x}_3=100\, m \nonumber
\ea
Il reste \`a d\'eterminer $(\dot{x}_2,\dot{y}_2)$. Comme $n=6$, on a $6-r=4\Longrightarrow r=2$ degr\'es de libert\'e. Alors:
\be
	\dot{X}=\begin{pmatrix}{
		\dot{x}_2 \cr
			\dot{y}_2 }
\end{pmatrix}
\ee
On laisse \`a titre d'exercice la d\'etermination de la fonctionnelle $\Phi$.
\subsection*{19.4.3. Le Modèle Stochastique}
 C'est le modèle qui décrit  les lois qui réagissent les erreurs des mesures des observables.
\bdf
On appelle erreur la quantité $e$ telle que:
 \begin{equation}
\fbox{ $ e=l-\dot{l}= \,\,\mbox{"Faux" - "Vrai"}	$}
\end{equation}
et correction la quantité $-e$.
\edf
Donc $\dot{l}=l+(-e)$. On corrige la valeur observée pour avoir la valeur réelle.
\\

On utilise un modèle des erreurs centrées avec une variance. En répétant les mesures, les erreurs de ces mesures sont la réalisation d'une variable aléatoire centrée $(e)$. Soit: $e=(e_1, e_2,..., e_n)^T$. On a alors:
\begin{equation}
E(e)=0\,\,\,\mbox{ou encore}\,\,\,lim_{n\rightarrow +\infty} \frac{\sum_i^n e_i}{n}= 0	
\end{equation}
Le vecteur erreur est:
\begin{equation}
	e=\,  _ne_1=L - \dot{L}
\end{equation}
La matrice variance est définie par:
\begin{equation}
	_{n}\sigma^2_n=\sigma^2=E(ee^T)
\end{equation}
C'est une matrice symétrique $n\times n$ et elle est positive c'est-à-dire:
\begin{equation}
	\forall \,U,\,\,\,U^T\sigma^2U \geq 0
\end{equation}
où $U$ est un vecteur à $n$ composantes. En effet:
$$U^T\sigma^2U=U^TE(ee^T)U=E(U^Tee^TU)=E((U^T.e).(e^T.U))$$
Or $U^T.e=w \in \BbR$ et $(e^T.U)=U^T.e=w$, par suite:
\be 
\forall \,U,\,\,\,U^T\sigma^2U=E((U^T.e).(e^T.U))=E(w.w)=E(w^2)\geq 0
\ee
Si on pose $S=\sigma^2$, on a $U^TSU\geq 0$. $S$ est positive et on démontre que les valeurs propres de $S$ sont positives et qu'on peut écrire $S$ sous la forme:
\begin{equation}
	S=V \begin{pmatrix}{
	\lambda_1^2  &   0 & \cdots & 0 \cr
  0  &  \lambda_2^2 & \cdots & 0 \cr
  \cdots  & \cdots & \ddots & \cdots \cr
	0 & \cdots & \cdots & \lambda_n^2 }
\end{pmatrix}V^T=V\Lambda^2V^T\Longrightarrow \sqrt{S}=V \Lambda V^T
\end{equation}
avec $V$ une matrice orthogonale normale c'est-à-dire : $V^{-1}=V^T$. Alors:
\begin{equation}
	\sigma=\sqrt{S}=\mbox{la matrice écart-type du vecteur  erreur }\,\,e
\end{equation}
En multipliant les mesures, on a un accroissement mécanique de la précision. Seulement la première hypothèse est très fragile: le centrage des erreurs  $E(e)=0$ est-il toujours vérifié ?

Les erreurs systématiques subsistent encore  vue l'insuffisance de définition des valeurs nominales $\dot{l}$ ainsi que les erreurs dues aux instruments. Ces effets montrent les limites de la méthode des moindres carrés.
\\

Le modèle des erreurs des mesures est la loi normale centrée, c'est-à-dire que $e \in \mathcal{N}(0,\sigma)$ où $\sigma$ est la matrice écart-type.
En utilisant le modèle stochastique $\mathcal{N}(0,\sigma)$ avec :
\begin{equation}
	\sigma^2=\begin{pmatrix}{
	\sigma_1^2  &   0 &\cdots & 0 \cr
  0  &  \sigma_2^2 & \cdots & 0 \cr
  \cdots  & \cdots & \ddots & \cdots \cr
	0 & \cdots & \cdots & \sigma_n^2 }
\end{pmatrix}
\end{equation}
et soit la ième équation d'observation:
\begin{equation}
	\sum_{j=1}^{j=n}A_{ij}x_j=l_i-e_i
\end{equation}
La règle qu'il faut utiliser c'est qu'il faut diviser chaque équation par l'écart-type correspondant:
 \begin{equation}
	\fbox{ $ \ds \frac{\sum_{j=1}^{j=n}A_{ij}x_j}{\sigma_i}=\ds \frac{l_i}{\sigma_i}-\frac{e_i}{\sigma_i} $}
\end{equation}
On obtient les estimations optimales, soit $\sigma$ petit, l'observation est précise et le poids $\ds p_i=\frac{1}{\sigma_i}$ est grand. 

\section{\textsc{Présentation de la Méthode des Equations d'Observations}}
Pour la méthode des équations d'observations, on a :
\begin{equation}
	\dot{L}=\Phi(\dot{X}) \label{eq111}
\end{equation}
On pose:
\be
\left\{\begin{array}{lllll}
	\dot{X}=X_0+\dot{x} \\
	\dot{L}=L_0+\dot{l} \\
	L=L_0+l= \mbox{"observé"} \\
	L_0=\Phi(X_0)=\mbox{"calculé"} \\
	e=L-\dot{L}=l-\dot{l}=\mbox{"l'erreur"} 
	\end{array}\right.
\ee             
En développant l'équation (\ref{eq111}), on obtient:
\ba
	L_0+\dot{l}=\Phi(\dot{X})=\Phi(X_0+\dot{x})=\Phi(X_0)+\left(\frac{\partial \Phi}{\partial X}\right)_{X_0}.\dot{x}+\epsilon(\dot{x})||\dot{x}|| \\
	\mbox{avec}\quad lim_{\dot{x} \rightarrow 0}\epsilon(\dot{x})\rightarrow 0 \nonumber
\ea                
On appelle :
\begin{equation}
	A=\left(\frac{\partial \Phi}{\partial X}\right)_{X_0}
\end{equation}
$A$ est une matrice $m\times n$ où $m$ est le nombre des relations d'observations et $n$ le nombre des paramètres inconnus. On obtient alors:
\begin{equation}
	A\dot{x}=\dot{l}-\epsilon(\dot{x}).\dot{x}
\end{equation}
On substitue à l'équation précédente l'équation:
\begin{equation}
		A\dot{x}=\dot{l} \label{eq88}
\end{equation}
sachant que $lim\,\epsilon(\dot{x})\rightarrow 0$ quand $\dot{x}\rightarrow 0 \Longrightarrow |\epsilon(\dot{x})| \ll 1$ et $ |\epsilon(\dot{x})| \|\dot{x}\| \leq |l_i |\Rightarrow |\epsilon(\dot{x})|\|\dot{x}\| \leq \|e\|$. On utilise une procédure d'itération en prenant $X_0$ une valeur estimée de l'inconnue.
\\

\textbf{Exemple}: Observation d'une distance entre deux points 1 et 2. On a les notations suivantes:
\ba
	\dot{x}_1=x_{10}+\Delta \dot{x}_1 \nonumber \\
	\dot{y}_1=y_{10}+\Delta \dot{y}_1 \nonumber \\
	\dot{x}_2=x_{20}+\Delta \dot{x}_2 \nonumber \\
	\dot{y}_2=y_{20}+\Delta \dot{y}_2 \label{eq2} \\
	\dot{d}=d_0+\Delta \dot{d} \nonumber  \\
	d=d_0+\Delta d \label{eq4 } \\
	 d^2_0=(x_{10}-x_{20})^2+(y_{10}-y_{20})^2 \nonumber 
\ea             
A partir de (\ref{eq2}), on a:
\ba
	&	 \dot{d}^2=(\dot{x}_1-\dot{x}_2)^2+(\dot{y}_1-\dot{y}_2)^2=(d_0+ \Delta \dot{d})^2  \nonumber  & \\ &
		 (d_0+ \Delta \dot{d})^2=(x_{10}+\Delta \dot{x}_1-x_{20}-\Delta \dot{x}_2)^2+(y_{10}+\Delta \dot{y}_1-y_{20}-\Delta \dot{y}_2)^2 &\nonumber \\ &
		 d_0^2+2d_0\Delta \dot{d}+\Delta \dot{d}^2=(x_{10}-x_{20})^2+(\Delta \dot{x}_1-\Delta \dot{x}_2)^2+2(x_{10}-x_{20})(\Delta \dot{x}_1-\Delta \dot{x}_2)+ \nonumber & \\ &
		 (y_{10}-y_{20})^2+(\Delta \dot{y}_1-\Delta \dot{y}_2)^2+2(y_{10}-y_{20})(\Delta \dot{y}_1-\Delta \dot{y}_2) &\nonumber
\ea               
Soit:
\begin{equation}
	\fbox{ $ \Delta \dot{d}=\ds \frac{x_{10}-x_{20}}{d_0}(\Delta \dot{x}_1-\Delta \dot{x}_2)+\frac{y_{10}-y_{20}}{d_0}(\Delta \dot{y}_1-\Delta \dot{y}_2)=d-d_0-e_d $}\label{eq66}
\end{equation}
L'équation (\ref{eq66}) représente l'équation d'observation d'une distance. On montre de même que l'équation d'un tour d'horizon s'écrit sous la forme:
\begin{equation}
	\fbox{ $ \ds 	\frac{y_{10}-y_{20}}{d_0^2}(\Delta \dot{x}_1-\Delta \dot{x}_2)-\frac{x_{10}-x_{20}}{d_0^2}(\Delta \dot{y}_1-\Delta \dot{y}_2)+dv_0=\alpha-\alpha_0-e_{\alpha} $} \label{eq7}
\end{equation}

\section{\textsc{La Solution des Moindres Carrés}}
Le système (\ref{eq88}) s'écrit en prenant en compte de la matrice des poids $P=\sigma^{-2}$:
\begin{equation}
	\sigma^{-1}A\dot{x}=\sigma^{-1}\dot{l} \Longleftrightarrow \sqrt{P}A\dot{x}=\sqrt{P}\dot{l} 
\end{equation}
La solution du système précédent vient en minimisant en norme euclidienne l'expression $\| \sqrt{P}A\dot{x}-\sqrt{P}\dot{l}\|$ ou encore:
\be
	U(x)=(\sqrt{P}Ax-\sqrt{P}l)^T(\sqrt{P}Ax-\sqrt{P}l)=x^T(A^TPA)x-2x^TA^TPL+l^TPl 
\ee
On obtient:
$$ \ds 	\frac{\partial U}{\partial x}=grad U(x)=2A^TPAx-2A^TPl $$
Si le minimum existe, $x$ est solution si:
	$$ 2A^TPAx-2A^TPl=0\Rightarrow A^TPAx=A^TPl $$
$A^TPA$ est symétrique, régulière et inversible car elle est carrée de rang $n$. Donc la solution des moindres carrés est:
\begin{equation}
	\fbox{ $ \overline{X}=(A^TPA)^{-1}A^TPl  \Rightarrow \overline{X}=X_0+\overline{x} $}
\end{equation}
Les observations compensées:
\begin{equation}
		\fbox{ $ \overline{L}=L_0+\bar{l}=L_0+A\overline{x} $}
	\end{equation}

\section{\textsc{Propriétés des Estimateurs}}
\textbf{Propriété 19.1} \textit{L'estimateur des moindres carrés $\ov{X}$ est centré  et sans biais c'est-à-dire:}
\begin{equation}
	\fbox{ $ E(\overline{X})=\dot{X}\quad \mbox{ou}\quad E(\overline{X})=\dot{x} $}
\end{equation}
Ainsi à des jeux différents d'observations, on fait correspondre des valeurs $\overline{X}$ différentes dont on peut estimer l'espérance mathématique ou moyenne. Cet estimateur est un vecteur aléatoire. Comme:
$$l=\dot{l}+e \Longrightarrow E(l)=E(\dot{l})+E(e)=\dot{l}+0=\dot{l}$$
alors:
\ba
	E(\overline{X})=E\left( (A^TPA)^{-1}A^TPl  \right)=(A^TPA)^{-1}A^TPE(l)=\nonumber \\
	(A^TPA)^{-1}A^TP\dot{l}=(A^TPA)^{-1}A^TPA.\dot{x} \nonumber 
	\ea
Soit:
\begin{equation}
E(\overline{X})=\dot{x}
\end{equation}
Parmi les estimateurs sans biais qui est le plus précis est $\overline{X}$  (celui dont l'écart-type est le plus petit), soit:
\begin{equation}
	\overline{X}=\begin{pmatrix}{
	\overline{X}_1 \cr
	\overline{X}_2\cr
	\vdots \cr
	\overline{X}_n }
\end{pmatrix}
\end{equation}
 $\overline{X}$ est le plus précis pour chaque composante. On obtient l'estimation la plus précise (la matrice $\sigma_{\ov X}$ est minimale). On détermine cette matrice:
 \be
 \sigma^2_{\overline{X}}=E\left\{\left(\overline{X}-E(\overline{X})\right)\left(\overline{X}-E(\overline{X})\right)^T\right\}=E\left\{\left(\overline{X}-\dot{x}\right)\left(\overline{X}-\dot{x}\right)^T\right\} 
\ee               
Or:
\begin{equation}
	\overline{X}-\dot{x}=(A^TPA)^{-1}A^TP(l-\dot{l})=(A^TPA)^{-1}A^TPe
\end{equation}
Par suite:
\ba
& \sigma^2_{\overline{X}}=(A^TPA)^{-1}A^TPE(ee^T)PA(A^TPA)^{-1}=\nonumber &\\ & (A^TPA)^{-1}A^TP.\sigma^2_0P^{-1}PA(A^TPA)^{-1}\nonumber & 
\ea 
\be
\fbox{$ \sigma^2_{\overline{X}}=\sigma^2_0(A^TPA)^{-1} $}\label{eq100} 
\ee                 
car $E(ee^T)=\sigma^2_0P^{-1}$ avec $\sigma^2_0$ un facteur scalaire.  
\\

Un estimateur est d\'efini en fonction d'un \'echantillon $l_1,l_2,...,l_n$ al\'eatoire, soit:
\begin{equation}
	\tilde{x}=\varphi(l_1,l_2,...,l_n)
\end{equation}
$\tilde{x}$ est proche de $\dot{x}$.
On cherche les estimateurs sans biais, c'est-\`a-dire:
\begin{equation}
	E(\tilde{x})=\dot{x}
\end{equation}
\textbf{Exemple:} soient une quantit\'e $\dot{x}$ et $l_1,l_2,...,l_n$ les observables ind\'ependants d'une m\^eme variable al\'eatoire de $x$, alors un estimateur de $\dot{x}$ est donn\'e par:
\begin{equation}
	\ov m=\frac{1}{n}\sum_{i=1}^nl_i
\end{equation}
Alors cet estimateur est sym\'etrique, sans biais\index{Estimateur sans biais} (les erreurs de mesures sont centr\'ees), de pr\'ecision maximale parmi les estimateurs lin\'eaires et sans biais, estimateur des moindres carr\'es et asymptotiquement normal (si $n \longrightarrow+\infty$), la fonction de r\'epartion de $\ov m$ comme variable al\'eatoire converge uniform\'ement vers une fonction de r\'epartition de la loi normale.
\section{\textsc{Les Résidus}}
Le vecteur résidu est défini par:
\begin{equation}
	\fbox{ $ V=\overline{L}-L=\bar{l}-l=\mbox{"compensé" - "observé"} $}
\end{equation}
Pourquoi calculer les résidus $v$:

1- Les composantes de $V$ sont de l'ordre des erreurs de mesures. La redondance des mesures est donnée par $n-r$, avec taux de redondance est égale à $(n-r)/n$. L'accroissement des mesures implique la précision mais ce n'est pas toujours vrai. 

2- Possibilité de vérifier la normalisation de la résolution:
\begin{equation}
	\fbox{ $ A^TPV=0 $}
\end{equation}
On a en effet:
	$$ A^TPV=A^TP(\bar{l}-l)=A^TP(A\overline{X}-l)=A^TPA\overline{X}-A^TPl=0$$
Soit $\tilde{x}$ une solution, $V'=A\tilde{x}-l\Longrightarrow A^TPV'=0\Longrightarrow \quad \tilde{x}=\overline{X}$.  
\subsection*{19.8.1. Expression du vecteur résidu}
Le vecteur résidu s'exprime par:
\begin{equation}
	V=(A(A^TPA)^{-1}A^TP-I)(l-\dot{l})=(A(A^TPA)^{-1}A^TP-I)e=B.e
\end{equation}
où $B=(A(A^TPA)^{-1}A^TP-I)$ est une matrice singulière (non inversible). Utilisant l'équation (\ref{eq100}), $V$ s'écrit:
\begin{equation}
V=(\ds \frac{1}{\sigma^2_0}A\sigma^2_{\overline{X}}A^TP-I)e \label{eq101}
\end{equation}
Comme $\bar{l}=A\overline{X}\Longrightarrow \bar{l}-\dot{l}=A(\overline{X}-\dot{x})$, d'où la variance des observations compensées:
\begin{equation}
\fbox{ $	\sigma^2_{\bar{l}}=E\left\{(l-\dot{l})(l-\dot{l})^T\right\}=A\sigma^2_{\overline{X}}A^T $}
\end{equation}
Utilisant l'équation (\ref{eq101}), on obtient :
\begin{equation}
	V=\left(\frac{\sigma^2_{\bar{l}}}{\sigma^2_0}P-I\right)e
\end{equation}
\section{\textsc{La Variance des Mesures}}
Par d\'efinition, la variance des mesures $l_i$ est :
\begin{equation}
	\tilde{\sigma}^2=E[(l_i-E(l_i))^2]
\end{equation}
or:
	$$ E(l_i)=\frac{1}{n}\sum_{i=1}^nl_i=\ov{m} $$
Par suite:
\ba
\tilde{\sigma}^2=\frac{1}{n}\sum_{i=1}^n(l_i-E(l_i))^2=\frac{1}{n}\sum_{i=1}^n(l_i^2-2l_iE(l_i)+E(l_i)^2)= \nonumber \\
\frac{1}{n}\sum_{i=1}^n(l_i^2-2l_i\ov m  +\ov m^2)=\frac{1}{n}\sum_{i=1}^nl_i^2-2\ov m  \frac{1}{n}\sum_i^nl_i+\ov  m ^2=\frac{1}{n}\sum_i^nl_i^2-\ov  m ^2 \nonumber
\ea
Soit:
\be
\fbox{ $\left\{\begin{array}{l}
\ov {m}=E(l_i)=\ds \frac{1}{n}\sum_{i=1}^nl_i \\
\tilde{\sigma}^2=\ds \frac{1}{n}\sum_i^nl_i^2-\ov {m}^2 
\end{array}\right. $}
\ee
Pour une variable centr\'ee, on a :
\begin{equation}
	\sigma^2=E(l^2)-E(l)^2=E(l^2)-\dot{x}^2\Longrightarrow E(l^2)=\sigma^2 +\dot{x}^2
\end{equation}
Comme:$$\ov {m}=\ds \frac{1}{n}\sum_il_i\longrightarrow \sigma_{\ov m}=\frac{1}{\sqrt{n}}\sigma_l=\frac{\sigma}{\sqrt{n}}$$
et: $$\sigma^2_{\ov m}=E(\ov {m}^2)-E(\ov {m})^2=  E(\ov {m}^2)-\dot{x}^2\Rightarrow E(\ov {m}^2)=\dot{x}^2+\sigma^2_{\ov m}=\dot{x}^2+\frac{\sigma^2}{n}$$
Donc:
\ba
&E(\tilde{\sigma}^2)=\ds E\left(\frac{1}{n}\sum_i^nl_i^2-\ov {m}^2\right)=\frac{1}{n}\sum_i^nE(l_i^2)-E(\ov  m ^2)=E(l_i^2)-E(\ov {m}^2)\nonumber &\\&
=\ds \sigma^2+\dot{x}^2-\left(\frac{\sigma^2}{n}+\dot{x}^2\right)=\frac{n-1}{n}\sigma^2\neq\sigma^2 &
\ea
Alors, l'estimateur $\tilde{\sigma}^2$ est biais\'e. Si on adopte comme formule de l'estimateur de la variance:
\begin{equation}
\fbox{ $ \tilde{\sigma}^2=\ds \frac{1}{n-1}\sum_{i=1}^n(l_i-\ov {m})^2 \label{eq132} $}
\end{equation}
on obtient:
\begin{equation}
	E(\tilde{\sigma}^2)=\sigma^2\Longrightarrow \,\,\,\,\mbox{l'estimateur est sans biais}
\end{equation}
(A vérifier à titre d'exercice!). De m\^eme, pour: $$\tilde{\sigma}'=\sqrt{\tilde{\sigma}^2}\Longrightarrow E(\tilde{\sigma}')\neq \sigma $$
\section{\textsc{La matrice d'Information de Fisher}}\index{Matrice d'information de Fisher}\index{\textbf{Fisher R.A.}}
\bdf
On appelle matrice d'information de Fisher\footnote{\textbf{Sir Ronald Aylmer Fisher} (1890-1962): staticien anglais. } relative \`a la fonction densit\'e de probabilit\'e $\m P(l_1,l_2,...,l_n,\dot{x}_1,\dot{x}_2,...,\dot{x}_r)=\m P(l_1,l_2,...,l_n,\dot{x})$ la matrice :
\be
	\fbox{ $ F=(F_{ij})=\left(E\left(-\ds \frac{\partial^2Log\m P}{\partial \dot{x}_i\partial \dot{x}_j}\right)\right) $}
\ee   
\edf           
\subsection*{19.10.1. La Quantit\'e d'Information de Fisher}
\bdf
Si $x$ est un scalaire, on appelle la quantit\'e d'information de Fisher la valeur:\index{Quantité d'information de Fisher}\index{\textbf{Fisher R.A.}}
\begin{equation}
	\fbox{ $ E\left(-\ds \frac{\partial^2Log\m P}{\partial \dot{x}^2}\right) $}
\end{equation}
\edf
Existe-t-il des estimateurs (sans biais) meilleurs que d'autres? La r\'eponse est donn\'ee par l'in\'egalit\'e de Rao\footnote{\textbf{Calyampudi Radhakrishna Rao} (1920 -\quad): mathématicien indien, éminent spécialiste en statistiques.}-Cramér\footnote{\textbf{Harald Cramér} (1893-1985): mathématicien et statisticien suédois.}.\index{\textbf{Rao C.R.}}\index{\textbf{Cramér H.}}


\section{\textsc{L'In\'egalit\'e de Rao-Cramér}}\index{Inégalité de Rao-Cramér}
L'erreur moyenne quadratique d'un estimateur sans bais $\tilde x$ peut s'écrire en utilisant la matrice variance-covariance $\sigma^2_x(\tilde x)=(\sigma_{ij}) \,\, i,j=1,2,...,n$ où :
$$ \sigma_{ij}\stackrel{def}{=}E[(\tilde x_i(X)-x_i)(\tilde x_j(X)-x_j)^T] $$
\bthm 
\textbf{(L'In\'egalit\'e de Rao-Cramér)}(S. Amari et H. Nagaoka, 2000) La matrice variance-covariance $\sigma^2_x(\tilde x)$ d'un estimateur sans biais $\tilde x$ vérifie:\index{\textbf{Amari S.}}\index{\textbf{Nagaoka H.}}
\be
\sigma^2_x(\tilde x) \geq F(\tilde x)^{-1}
\ee
dans le sens que $\sigma^2_x(\tilde x)-F(\tilde x)^{-1}$ est une matrice définie semi-positive, soit:
\be
\forall u \in \BbR^r \,\,\,u^T.\sigma^2_{\tilde{x}}.u\geq u^T.F^{-1}(\tilde x).u
\ee
\ethm
\bdf
$\tilde{x}$ est un estimateur efficace de $\dot{x}$, si\index{Estimateur efficace}: 
\be
	\forall u \in \BbR^r \,\,\,u^T.\sigma^2_{\tilde{x}}.u = u^T.F^{-1}(\tilde x).u
\ee
\edf
En appliquant le principe de maximum de vraisemblance \`a la fonction de densit\'e $$ \m P(l_1,l_2,...,l_n,\dot{x}_1,\dot{x}_2,...,\dot{x}_r)=\m P(l_1,l_2,...,l_n,\dot{x})$$
on obtient:
\ba
&\m P(l_1,l_2,...,l_n,\tilde{x}_1,\tilde{x}_2,...,\tilde{x}_r)=
\m P(l_1,l_2,...,l_n,\dot{x}_1,\dot{x}_2,...,\dot{x}_r)+\sum_i^r(\tilde{x}_i-\dot{x}_i)\frac{\partial \m P}{\partial \dot{x}_i}& \nonumber \\& +(\mbox{termes 2\`eme ordre}) \Longrightarrow \frac{\partial \m P}{\partial \dot{x}_i}=0\Longleftrightarrow \frac{1}{\m P}\frac{\partial \m P}{\partial \dot{x}_i}=0&	\nonumber \\&
\Longrightarrow \ds \frac{\partial Log \m P}{\partial \dot{x}_i}=0\,\, \mbox{pour}\,\,i=1,2,...,r &
\ea
Donc à l'ordre deux, on a : $\m P(l_1,l_2,...,l_n,\tilde{x}_1,\tilde{x}_2,...,\tilde{x}_r)=
\m P(l_1,l_2,...,l_n,\dot{x}_1,\dot{x}_2,...,\dot{x}_r)+(\mbox{termes 2\`eme ordre})$. Alors $\tilde{x}_1,\tilde{x}_2,...,\tilde{x}_r$ sont dits asymptotiquement sans biais, asymptotiquement de pr\'ecision maximale et asymptotiquement normaux et $\m P$ est maximale.

\paragraph{\textbf{Exemple de l'application de l'inégalité de Rao-Cramér:}}
Soit le cas d'un seul param\`etre inconnu $\dot{x}$, d'\'echantillon $l_1,l_2,...,l_n$ de densit\'e de probabilit\'e:
\begin{equation}
\m P(l,\dot{x})=\frac{1}{(\sqrt{2\pi})^n}\frac{1}{\sigma^n}e^{-\ds \frac{1}{2}\frac{\sum_{i=1}^{i=n}(l_i-\dot{x})^2}{\sigma^2}}
\end{equation}
On a alors:
\ba
		Log\m P=Constante -\frac{1}{2}\frac{\sum_i^n(l_i-\dot{x})^2}{\sigma^2} \nonumber \\
	\frac{\partial Log\m P}{\partial \dot{x}}=\frac{1}{\sigma^2}\sum_i^n(l_i-\dot{x})\Rightarrow \frac{\partial^2Log\m P}{\partial \dot{x}^2}=-\frac{n}{\sigma^2} \nonumber
\ea           
Par suite:
\be
	E\left(-\ds \frac{\partial^2Log\m P}{\partial \dot{x}^2}\right)=E\left(\frac{n}{\sigma^2}\right)=nE\left(\frac{1}{\sigma^2}\right)=\frac{n}{\sigma^2}
\ee           
L'inégalité de Rao-Cramér est bien vérifiée, en effet $\sigma^2_{\tilde x}=\sigma^2$, on a bien sûr:
\begin{equation}
		\sigma_{\tilde{x}}^2	\geq \frac{1}{E\left(-\ds \frac{\partial^2Log\m P}{\partial \dot{x}^2}\right)}=\frac{1}{\frac{n}{\sigma^2}}=\frac{\sigma^2}{n }
\end{equation}
\\

Soit un triangle plan o\`u on mesure les deux c\^ot\'es $b, c$ et l'angle $A$ et on veut estimer le c\^ot\'e $a$, on a:
\begin{equation}
	\dot{a}=\sqrt{\dot{b}^2+\dot{c}^2-2\dot{b}\dot{c}cos\dot{A}}
\end{equation}
Un estimateur de $a$ est:
\begin{equation}
	\tilde{a}=\sqrt{b^2+c^2-2bc.cosA}
\end{equation}
Cet estimateur est biais\'e.
\\

Comme $\dot{x}=f(\dot{l})$	et $\tilde{x}=f(l)$, donc:
\be
	E(\tilde{x})\stackrel{?}{=}\dot{x}\Longrightarrow E(f(l))\stackrel{?}{=}f(\dot{l})
\ee
Or:
$	\tilde{x}=\dot{x}+\underbrace{(l-\dot{l})^Tf'}_{e}+\mbox{termes ordre 2}=\dot{x}+e+\mbox{termes ordre 2} $
\be
	\Rightarrow E(\tilde{x})=E(\dot{x})+E(e)+E(\mbox{termes ordre 2})=\dot{x}+E(\mbox{termes ordre 2})
\ee          
si $e$ est petite par rapport aux valeurs. Donc $\tilde{x}$ est quasiment sans biais.  
\\


\section{\textsc{L'Ecriture Matricielle des Equations d'Observations}}
Les \'equations d'observations peuvent s'\'ecrire matriciellement:
\be
	AX+K=V  \label{c34}
\ee                 
o\`u:

- $A$ est la matrice des coefficients;
 
 - $V$ le vecteur des r\'esidus;
 
-  $K$ est le vecteur: $ K=  \textit{"calcul\'e"} - \textit{"observ\'e"} $;

- $X$ est le vecteur des inconnues (les corrections $dx_i, dy_i, ...$  aux valeurs approch\'ees $x^0_i, y^0_i, ...$).
\\

 Pour $n$ \'equations d'observations et $r$ le nombre des inconnues, on a :
\be
	_{n}A_{r}._{r}X_{1}+\,_{n}K_{1}=\, _{n}V_{1}, \ \ \ n>r  \label{c35}
\ee              
La solution de (\ref{c35}) par la m\'ethode des moindres carr\'es est obtenue en minimisant la forme quadratique:
\be
	\sum_{i=1}^{i=n}V^2_i=V^TV=(AX+K)^T(AX+K)  \label{c36}
\ee               
o\`u $V^T$ d\'esigne la matrice transpos\'ee de $V$.

Sachant que:
\ba
	(A+B)^T=A^T+B^T \nonumber  \\
	(AB)^T=B^TA^T \nonumber 
\ea                
L'\'equation (\ref{c36}) s'\'ecrit:
$$ 	V^TV=X^T(A^TA)X + X^TA^TK+K^TAX+K^TK  $$
Comme $X^TA^TK$ et $K^TAX$ sont des scalaires, alors:
	\[X^TA^TK=K^TAX
\]
Par suite:
$$ 	V^TV=X^T(A^TA)X + 2X^TA^TK+K^TK $$
On pose:
\be
	N=A^TA  \label{c41}
\ee               
Alors la matrice N est sym\'etrique $(N^T=N)$ et inversible c'est-\`a-dire r\'eguli\`ere et de rang $r$. En effet la matrice $
_{n}A_{r}$ est de rang $r$ c'est-\`a-dire on peut extraire une sous matrice $_{r}A^{'}_{r}$ de $A$ telle que son d\'eterminant est diff\'erent de z\'ero:
$$ 	\mbox{\textit{Dét}}(_{r}A'_{r}) \neq 0 $$
La matrice $N$ est d\'efinie positive, on entend par l\`a que :
\be
	\forall X\neq0 \ \ X^TNX > 0  \label{c43}
\ee             
car : 	$$ X^TNX=X^TA^TAX=(AX)^T(AX)=\|AX.AX\|=\|AX\|^2 \geq 0  $$ 
le carré de la norme du vecteur $A.X$ et pour $ X\neq 0$ on a $A.X\neq 0$ sinon la matrice $A$ serait de rang < $r$, par suite (\ref{c43}) est v\'erifi\'ee.

La matrice $N$ est dite \textit{la matrice normale}.\index{Matrice normale}

En posant:
\ba
	F(X)=X^T(A^TA)X +2X^TA^TK+K^TK \nonumber   \\
\mbox{ou} \ F(X)=V^TV=X^TNX +2X^TA^TK+K^TK \label{c46}
\ea          
c'est une fonction scalaire du vecteur X. La solution de (\ref{c46}) avec $F(X)$ minimum est telle que:
$$ 	\frac{\partial F}{\partial X}=0  $$
La d\'eriv\'ee de (\ref{c46}) par rapport au vecteur des inconnues est telle que:
\be
	\frac{\partial F}{\partial X}=2N.X+2A^T.K = 0 \Rightarrow N.X+A^T.K=0 	\label{c48}
\ee           
Soit:
\be
\fbox{ $ 	\overline{X}=-N^{-1}A^TK $} \label{c49} 
\ee               
o\`u $\overline{X}$ est le vecteur d\'etermin\'e des inconnues. Si $X_0$ est le vecteur approch\'e des inconnues, on a :
\be
	\tilde{X}=X_0+\overline{X} \label{c50}
\ee               
o\`u $ \tilde{X} $ est le vecteur des valeurs d\'efinitives des inconnues.

Le vecteur $\tilde{V}$ se d\'etermine par:
\be
	\tilde{V}=A \overline{X}+K=-AN^{-1}A^TK+K=(I-AN^{-1}A^T)K \label{c51}
\ee             
Le vecteur des observations compens\'ees est donn\'e par:
\be
	\tilde{L}=l+\tilde{V} \label{c52}
\ee
Notons que:
\be
	\fbox{ $ A^T \tilde{V}=0   $} \label{c53}
\ee
En effet, $A^T \tilde{V}=A^T(I-AN^{-1}A^T)K=A^T(K-AN^{-1}A^TK)=A^TK-A^TAN^{-1}A^TK$, soit $A^T \tilde{V}=0 $ en tenant compte que $N=A^TA$. La condition (\ref{c53}) est appel\'e renormalisation\index{Condition de renormalisation}. Elle est importante car elle garantit que le r\'esultat obtenu $\tilde{X}$ est bien celui des moindres carr\'es. 
\\

Dans ce chapitre, on a consid\'er\'e le syst\`eme:
	\[AX+K=V
\]
sans parler de la matrice $P$. Si on consid\`ere la matrice de poids $P$ et la matrice de variance des observations:
$$ \Gamma_l	=\sigma_0^2 P^{-1}  $$
La solution de (\ref{c34}) des moindres carr\'es est obtenue \`a partir du syst\`eme normal:
$$ 	(A^TPA)\overline{X}+A^TPK=0  $$
en minimisant $V^TPV$ ou encore si la matrice $P$ est diagonale $P=(p_{ii}=p_i)$: 
\be
	\fbox{ $\ds  \sum_{i=1}^{i=n}p_iv^2_i \,\,\,\,minimum $} \label{c56}
\ee               
Posant encore $N=A^TPA$ \textit{la matrice normale}\index{Matrice normale}, la solution de (\ref{c56}) est:
\be
	\fbox{ $ \overline{X}=-(A^TPA)^{-1}A^TPK=-N^{-1}A^TPK $}  \label{c57}
\ee          
La condition de la renormalisation devient:
\be
	\fbox{ $ A^TP\tilde{V}=0 $} \label{c58}
\ee          
Pour l'estimation du facteur de la variance unitaire, on admet que l'estimateur de $\sigma^2_0$ est donn\'e par :
\be
	\fbox{ $ s^2=\ds \frac{\tilde{V}^TP\tilde{V}}{n-r} $}\label{c59}
\ee              
Pour l'estimateur de la variance de $\overline{X}$, on va utiliser la propri\'et\'e suivante:
\be
	\fbox{ $ \mbox{si} \ Y=A.Z \Longrightarrow \ \Gamma_Y=A.\Gamma_Z.A^T $}  \label{c60}
\ee             
avec $\Gamma_Y$ la matrice variance de $Y$ et $\Gamma_Z$ celle de $Z$.
En effet:
\ba
	\Gamma_Y=E\left\{(E(Y)-Y)(E(Y)-Y)^T\right\} \nonumber  \\
	E(Y)=E(AZ)=AE(Z) \nonumber
\ea           
avec $E$ l'op\'erateur esp\'erance math\'ematique. D'o\`u:
\ba
	&\Gamma_Y=E\left\{(AE(Z)-AZ)(AE(Z)-AZ)^T\right\}= E\left\{A(E(Z)-Z)(A(E(Z)-Z))^T\right\}& \nonumber \\
	&=AE\left\{(E(Z)-Z)(E(Z)-Z)^TA^T\right\} & \nonumber \\
	&= AE\left\{(E(Z)-Z)(E(Z)-Z)^T\right\}A^T & \label{c62}
\ea	
Or :
	\[\Gamma_Z=E\left\{(E(Z)-Z)(E(Z)-Z)^T\right\}
\]
Par suite:
\be
	\Gamma_Y= A\Gamma_ZA^T  \label{c63}
\ee               
Comme: $$\overline{X}=-N^{-1}A^TPK $$
 d'o\`u:
\be
	Var(\overline{X})=-N^{-1}A^TP.Var(K).(-N^{-1}A^TP)^T=N^{-1}.A^T.P.Var(K).P.A.N^{-1} \label{c64}
\ee           
Comme $ K=L_0-l \Rightarrow Var(K)=Var(l)=\Gamma_l=\sigma_0^2 P^{-1}$, d'o\`u:
\ba
&\Gamma_{\overline{X}}=N^{-1}A^TP.\sigma_0^2.P^{-1}PAN^{-1}=\sigma_0^2.N^{-1}A^TPAN^{-1}=\sigma_0^2.N^{-1} & \nonumber  \\ &
\fbox{ $ \Gamma_{\overline{X}}=\sigma_0^2.N^{-1} $}  & \label{c65}
\ea           
D'o\`u l'estimateur de la variance de $\overline{X}$:
\be
	\fbox{ $ \overline{\Gamma}_{\overline{X}}=s^2.N^{-1} $}  \label{c66}
\ee             
La matrice $\overline{\Gamma}_{\overline{X}}$ permet de d\'eterminer la pr\'ecision des inconnues d\'etermin\'ees par la m\'ethode des moindres carr\'es.


\section{\textsc{M\'ethode des Equations d'Observations avec Equations de condition}}\index{Equations de condition}
\subsection*{19.13.1. Cas o\`u les $r$ grandeurs \`a d\'eterminer sont li\'es par $p$ relations}
On suppose qu'il existe parmi les n grandeurs observ\'ees $l_i$ un ensemble au moins de $r$ grandeurs observ\'ees permettant de d\'eterminer les $r$ grandeurs inconnues, c'est-\`a-dire que la matrice $A$ de $AX+K=V$ est de rang $r$.

Apr\`es lin\'earisation (si n\'ecessaire), les liaisons entre les $r$ grandeurs inconnues donneront bien \`a un syst\`eme de $p$ \'equations de condition:
\be
\left. \begin{array}{ll}
	BX=M \\
	\mbox{ou} \, _{p}B_{r}._{r}X_{1}=\,_{p}M_{1}
	\end{array} \right\}
	\ee
avec rang $B$ = $p$ et $B$ et $M$ ne d\'ependent que de $X$.

Pour estimer $X$, on a consid\'er\'e l'ensemble des relations:
\be
\left\{\begin{array}{ll}
	AX+K=V \\
	BX=M
	\end{array}\right.
\ee
avec rang $A=r$ et rang $B=p$. Apr\`es compensation, on aura le syst\`eme:
\be
\left\{\begin{array}{ll}
	A\overline{X}+K=\overline{V} \\
B\overline{X}=M  
	\end{array}\right. \lb{eq141}
\ee             
On remarque qu'une \'equation de condition est \'equivalente \`a une \'equation d'observations dont le poids est infini: en effet $M$ peut \^etre consid\'er\'e comme un vecteur d'observations fictives certaines dont les erreurs moyennes quadratiques $(emq)$ sont nulles et par suite les poids infinis. 
\\

On se ram\`ene au cas plus haut, \`a partir de $BX=M$ exprimer $p$ inconnues en fonction des $r-p$ autres et les reporter dans $AX+K=V$.
\section{\textsc{Application de la M\'ethode des Moindres Carr\'es}}
 On a le mod\`ele fonctionnel lin\'eaire (\ref{eq141}) qu'on écrit sous la forme:
\be
	\left\{\begin{array}{ll}
	A\dot{x}=l-\dot{e} \\
	B\dot{x}= k
	\end{array}\right.
\ee
avec:
$$       A=\begin{pmatrix}{
	a_{11}     & a_{12} & \cdots & a_{1p} \cr
	a_{21}     & a_{22} & \cdots & a_{2p} \cr
      	\vdots & \vdots       & \ddots & \vdots      \cr
	a_{n1}     & a_{n2} & \cdots & a_{np} } 
\end{pmatrix} $$ 
$$ \dot{x}=\begin{pmatrix} {
	\dot{x}_1 \cr
	\dot{x}_2 \cr
	\vdots \cr
	\dot{x}_p  }
	\end{pmatrix};\,\,\,\,l=\begin{pmatrix}{
	l_1 \cr
	l_2 \cr
	\vdots \cr
	l_n }
	\end{pmatrix};\,\,\,\,\dot{e}=\begin{pmatrix}{
	\dot{e}_1 \cr
	\dot{e}_2 \cr
	\vdots \cr
	\dot{e}_n }
	\end{pmatrix}\,\,\,\mbox{avec}\,\,n\geq p $$
$$ B=\begin{pmatrix}{
	b_{11}     & b_{12} & \cdots & b_{1p} \cr
	b_{21}     & b_{22} & \cdots & b_{2p} \cr
      	\vdots & \vdots       & \ddots & \vdots      \cr
	b_{p1}     & b_{p2} & \cdots & b_{pp} }
\end{pmatrix},\quad et \quad k =\begin{pmatrix}{
	k_1 \cr
	k_2 \cr
	\vdots \cr
	k_p }
\end{pmatrix} $$
On suppose que les erreurs $\dot{e}$ suivent la loi normale  $ \Longleftrightarrow   \dot{e} \in \m N (0,\sqrt{P^{-1}}\sigma_0)$, $P$ est le poids et $\sigma^2=P^{-1}\sigma^2_0$. Ecrivons la fonction de densit\'e de probabilit\'e, soit:
$$ \m P(l_1,l_2,...,l_n,\dot{x}_1,\dot{x}_2,...,\dot{x}_p)=\frac{1}{(\sqrt{2\pi})^n}\frac{1}{(\sigma_0^2)^{n/2}det\sqrt{P^{-1}}}e^{\ds -\frac{1}{2}(l-A\dot{x})^T\frac{P}{\sigma^2_0}(l-A\dot{x})}	$$
Les inconnues sont $\dot{x}$ et $\sigma^2_0$. On consid\`ere la fonction scalaire:
$$ 	U=Log\m P+\Lambda^T.(B\dot{x}-k) $$
avec $\Lambda$ le vecteur des multiplicateurs de Lagrange\index{Multiplicateurs de Lagrange}:\index{\textbf{Lagrange J.L.}}
	\be 
\Lambda      =\begin{pmatrix}{
	\lambda_1 \cr
	\lambda_2 \cr
	\vdots \cr
	\lambda_p }
\end{pmatrix}	
\ee
Les inconnues sont obtenues en cherchant les extr\'emums de $U$. $U$ s'\'ecrit:
\be
\fbox{ $ U=\ds Cte-\frac{n}{2}Log\sigma^2_0-\frac{1}{2}(l-A\dot{x})^T\frac{P}{\sigma^2_0}(l-A\dot{x})+\Lambda^T.(B\dot{x}-k) $}
\ee
 Les extr\'emums de $U$ sont obtenus par la solution de:
\ba             
	\frac{\partial U}{\partial \sigma^2_0}=-\frac{n}{2\sigma^2_0}+\frac{1}{2\sigma_0^4}(l-A\dot{x})^TP(l-A\dot{x})=0\label{aq51} \\
	\frac{\partial U}{\partial \dot{x}}=-\frac{A^TPA}{\sigma^2_0}\dot{x}+\frac{A^TPl}{\sigma_0^2}+B^T\Lambda=0 \label{aq52}\\
	\frac{\partial U}{\partial \Lambda}=B\dot{x}-k=0 \label{aq53}
	\ea
En notant que :
$$ 	A\dot{x}-l=v=\,\,\mbox{le vecteur r\'esidu} $$
L'\'equation (\ref{aq51}) donne:
\be
	\tilde{\sigma}_0^2=\frac{v^TPv}{n}
\ee
 (\ref{aq52}) et (\ref{aq53}) donnent le syst\`eme :
\be
	\fbox{ $ \begin{pmatrix}{
	A^TPA & B^T \cr
	B & 0 }
\end{pmatrix}.\begin{pmatrix}{
	\tilde{x} \cr
	-\Lambda \sigma_0^2 }
\end{pmatrix}=\begin{pmatrix}{
	A^TPl \cr
	k }
\end{pmatrix} $}
\ee
La solution du syst\`eme pr\'ec\'edent donne la solution des moindres carr\'es. L'estimateur $\tilde{x}$ est un estimateur sans biais, asymptotiquement normal, alors que l'estimateur :
\be
	\fbox{ $ \tilde{\sigma}_0^2=\ds \frac{v^TPv}{n} $}
\ee
est biais\'e mais asymptotiquement sans biais. Par contre, si on prend:
\be
\fbox{ $	\tilde{\sigma}_0^{'2}=\ds \frac{v^TPv}{n-r} $}
\ee
$\tilde{\sigma}_0^{'2}$ est sans biais.
\\

Une question se pose: le choix des param\`etres influence-t-il sur les solutions. Pour la m\'ethode des moindres carr\'es, la solution est ind\'ependante. En effet, soit le syst\`eme:
$$ A\dot{x}=l-\dot{e}\Longrightarrow \mbox{la solution}\,\,\,\tilde{x}$$
Soit le syst\`eme avec d'autres param\`etres:
\be
	\dot{y}=C\dot{x}+k\Longrightarrow \tilde{y}=C\tilde{x}+k
\ee
Soit $\tilde{y}=C\tilde{x}'+k$ une autre solution, on a alors:
\be
	\sigma^2_{\tilde{y}}=C\sigma^2_{\tilde{x}'}C^T
\ee
En utilisant l'in\'egalit\'e de Rao-Cramér:
\be
u^T\sigma^2_{\tilde{y}}u=u^TC\sigma^2_{\tilde{x}'}C^Tu=(C^Tu)^T\sigma^2_{\tilde{x}'}(C^Tu)\,\,\,\mbox{minimum si}\,\,\,\sigma^2_{\tilde{x}'}\,\,\mbox{est minimum}\Rightarrow \tilde{x}'=\tilde{x}
\ee 

\section{\textsc{Exemples de Pose d'Equations d'Observations}}
On a les \'el\'ements de base:
\be
\begin{array}{l}
	\dot{L}=l+\dot{v}    \\
	\sum v^2=\mbox{minimum} \\
	\Phi (\dot{X},\dot{L})=0 
	\end{array}
\ee
o\`u $\dot{X},\dot{L} $ sont respectivement les vecteurs des valeurs nominales des param\`etres inconnus et de l'observable. $\dot{v}$ le vecteur des résidus réels et $\Phi$ une fonction liant les inconnues et l'observable.

A partir de la relation:
\be
	dL+L_0 - l=v \label{c4}
	\ee              
soit:
\begin{center}
\fbox{ $ \textbf{Compens\'e + Calcul\'e - Observ\'e = R\'esidu} $}
\end{center}
	
On va \'ecrire la relation (\ref{c4}) pour diff\'erents types d'observations.
\subsection*{19.15.1. La G\'eod\'esie Bidimensionnelle}
\subsubsection*{Observation d'une distance}
Soient deux points $M_1$ et $M_2$ de coordonn\'ees approximatives $(x_1,y_1)$ et $(x_2, y_2$). La distance calcul\'ee est :
\be
	D_0=\sqrt{(x_1-x_2)^2+(y_1-y_2)^2}  \label{c5}
\ee                 
La fonction $\Phi$ est:
\be
	\Phi=\Phi(x_1, x_2, y_1, y_2)=D- \sqrt{(x_1-x_2)^2+(y_1-y_2)^2}=0  \label{c6}
\ee              
Soit: 
\be
	D= \sqrt{(x_1-x_2)^2+(y_1-y_2)^2} \label{c7}
\ee               
D'o\`u:
\be
	dD=\frac{x_1-x_2}{D}(dx_1-dx_2)+\frac{y_1-y_2}{D}(dy_1-dy_2)  \label{c8}
\ee             
En posant: 
\be
\begin{array}{l}
	\Delta X =x_1 - x_2     \\
	\Delta Y = y_1 -y_2  
	\end{array}
\ee
On a alors la relation d'observations d'une distance entre les points 1 et 2:
\be
	\fbox{ $ \ds \frac{\Delta X}{D_0}dx_1+\frac{\Delta Y}{D_0}dy_1 -\frac{\Delta X}{D_0}dx_2-\frac{\Delta Y}{D_0}dy_2+D_0 - D_{obs}=v_{12} $} \label{c11}
\ee               
Si le point 1 est connu, on a alors:
\be
-\frac{\Delta X}{D_0}dx_2-\frac{\Delta Y}{D_0}dy_2+D_0 - D_{obs}=v_{12}  \label{c12}
\ee                
\subsubsection*{Observation Angulaire}
Ecrivons la relation (\ref{c4}) pour la direction $M_1M_2$:
\be
g=\lambda +\vartheta_A       \label{c13}   
\ee              
\be
\mbox{\textbf{Gisement réel = Lecture réelle + Constante de désorientation du limbe}}  \label{c14}
\ee               
Cette relation est remplac\'ee par une relation d'estimation:
\be
	G=L + V_A  \label{c15}
\ee             
Comme $ L= l+v=$ lecture faite + correction de compensation, d'o\`u:
\be
	G=L + V_A=l+v+V_1   \label{c16}
\ee             
On introduit un r\'eseau g\'eom\'etrique approch\'e $(x_1, y_1,x_2, y_2, G_0)$ et une valeur approch\'ee du calage du limbe $V_{01}$ au point $M_1$ d'o\`u:
\ba
	G_0+ dG = l+v+V_{01} +dV_1  \nonumber    \\
	\mbox{ou encore} \ dG - dV_1 + G_0 - l - V_{01} = v \nonumber
\ea
$dG$ est la diff\'erentielle de $G$:
\ba
	tg G = \frac{x_2 - x_1}{y_2 - y_1}   \nonumber \\
	(1+tg^2G)dG=\frac{(y_2-y_1)(dx_2-dx_1)-(x_2-x_1)(dy_2-dy_1)}{(y_2-y_1)^2} \nonumber 
\ea
$dG$ est calcul\'ee \`a partir du r\'eseau approch\'e et on a :
\ba
	x_2 - x_1=DsinG_0 \nonumber   \\
	y_2 - y_1=DcosG_0 \nonumber
\ea
Il vient:
$$ 	dG = \frac{cosG_0(dx_2-dx_1)}{D}-\frac{sinG_0(dy_2-dy_1)}{D}  $$
D'o\`u la relation d'observations angulaires au point $M_1$ vers le point $M_2$:
\be
\fbox{ $	\ds -\frac{cosG_0}{D}dx_1+\frac{sinG_0}{D}dy_1+\frac{cosG_0}{D}dx_2-\frac{sinG_0}{D}dy_2-dV_1 +G_0-l-V_{01}=v_{12} $} \label{c24}
\ee              
Si le point $M_1$ est connu et qu'il y'a une inconnue d'orientation, alors l'\'equation pr\'ec\'edente devient:
\be
-dV_1+\frac{cosG_0}{D}dx_2-\frac{sinG_0}{D}dy_2 +G_0-l-V_{01}=v_{12} \label{c25}
\ee              
Dans le cas d'un rel\`evement sur un point $M_2$ connu, alors on a $dx_2=dy_2=0$ et:
\be
-dV_1-\frac{cosG_0}{D}dx_1+\frac{sinG_0)}{D}dy_1+(G_0-l-V_{01})=v_{12} \label{c26}
\ee                
Pour les \'equations (\ref{c24}), (\ref{c25}), et (\ref{c26}) les r\'esidus $v_{12}$ sont exprim\'es en radians. Pour avoir la m\^eme unit\'e (c-a-d le $m$) que les r\'esidus des \'equations (\ref{c11}) et (\ref{c12}), on multiplie les \'equations (\ref{c24}), (\ref{c25}) et (\ref{c26}) par la distance $D$. On aura pour l'\'equation (\ref{c26}), l'\'equation:
\be
	\fbox{ $ -DdV_1-cosG_0dx_1+sinG_0dy_1+D(G_0-l-V_{01})=v_{12}^{'} $} \label{c27}
\ee    
\subsection*{19.15.2. La G\'eod\'esie Tridimensionnelle}
\subsubsection*{Observation d'une distance}
Soient $(X_1, Y_1, Z_1)$ et $(X_2, Y_2, Z_2)$ les coordonn\'ees tridimensionnelles approch\'ees des points $M_1$ et $M_2$. La distance spatiale $M_1M_2$ est calcul\'ee par:
$$ 	D^{cal}=\sqrt{(X_2-X_1)^2+(Y_2-Y_1)^2+(Z_2-Z_1)^2} $$
D'o\`u l'\'equation d'observations d'une distance spatiale en g\'eod\'esie 3D:
\be
	dD + D^{cal} - D^{obs}=v \label{c29}
\ee                
avec:
$$ 	dD=\frac{(X_2 -X_1)(dX_2-dX_1) + (Y_2 -Y_1)(dY_2-dY_1)+(Z_2 -Z_1)(dZ_2-dZ_1)}{D^{cal}} $$
On pose:
\ba
	\Delta X= X_2 - X_1  \nonumber \\
\Delta Y= Y_2 - Y_1 \label{c31}  \\
\Delta Z= Z_2 - Z_1 \nonumber
\ea
La relation (\ref{c29}) devient:
\be
\fbox{ $ \ds -\frac{\Delta X}{D^{cal}}dX_1-\frac{\Delta Y}{D^{cal}}dY_1-\frac{\Delta Z}{D^{cal}}dZ_1+\frac{\Delta  X}{D^{cal}}dX_2+\frac{\Delta Y}{D^{cal}}dY_2+\frac{\Delta Z}{D^{cal}}dZ_2+D^{cal} - D^{obs}=v_{12} $} \label{c32}
\ee                
Si le point $M_1$ est connu, alors $ dX_1 = dY_1 = dZ_1 = 0$ et :
\be
\frac{\Delta  X}{D^{cal}}dX_2+\frac{\Delta Y}{D^{cal}}dY_2+\frac{\Delta Z}{D^{cal}}dZ_2+D^{cal} - D^{obs}=v_{12} \label{c33}
\ee               
\section{\textsc{Les Inverses Généralisées}}\index{Inverses généralisées}
Dans le probl\`eme des compensations, on est amen\'e \`a inverser des matrices.
\\

Alors si $A$ est une matrice carr\'ee inversible, alors son inverse $A^{-1}$ est unique et v\'erifie:
\be
AA^{-1}  = A^{-1}A = I
\ee
o\`u $I$ est la matrice unité.
\\

Si $A$ est une matrice rectangulaire:
\begin{equation}
	A^{-1}\,\,\mbox{est l'inverse de}\,\,A\,\,\mbox{si}\,\,AA^{-1}A = A 
\end{equation}
$A^{-1}$ est solution de :
\begin{equation}
	AXA=A
\end{equation}
Ce système a une infinit\'e de solutions.
\\

Les inverses g\'en\'eralis\'ees servent \`a r\'esoudre les \'equations lin\'eaires impossibles. Par exemple:
\be
\begin{array}{llll}
	x_1+x_2=k_1  \\
	x_1+2x_2=k_2 \\
	x_1+3x_2=k_3 \\
	2x_1+x_2=k_4
	\end{array}
\ee
syst\`eme qu'on \'ecrit sous la forme:
\be
	Ax\cong k\Rightarrow x=A^{-1}k,\,\,\,\,A^{-1}\,\mbox{inverse g\'en\'eralis\'ee de}\,A\Rightarrow \mbox{donne toutes les solutions correctes} 
\ee
\subsection*{19.16.1. Propri\'et\'es des matrices}
\textbf{Propriété 19.2} \textit{Soient deux matrices $A$ et $B$, alors:}
\be
		rang (AB)\leq \left\{  
	\begin{array}{ll}
	  	rang(A) \\
		rang(B)
	\end{array}\right.
\ee
\textbf{Propriété 19.3} (\textbf{Th\'eor\`eme des pivots}): \textit{Soit la matrice $A=\,_nA_p$, s'il existe une sous matrice d'ordre $r$ $_{r}A_r$$=A_{11}$ r\'eguli\`ere de rang $r$, on peut \'ecrire la matrice $A$ sous la forme:}
\ba
	A=\begin{pmatrix}{
	A_{11} & A_{12} \cr
	A_{21} & A_{22}}
\end{pmatrix} \label{ac0}\\
\mbox{\textit{avec}}\,\,_{r}A_r=A_{11};\,\,_{r}A_{p-r}=A_{12};\,\,_{n-r}A_r=A_{21},\,\,_{n-r}A_{p-r}=A_{22} \\
\mbox{\textit{et}}\,\,A_{11}\,\,\mbox{\textit{inversible et }}\,\, A_{22}=A_{21}A_{11}^{-1}A_{12}
\ea
Les $r$ premi\`eres colonnes de $A$ sont ind\'ependantes $\Rightarrow$ les colonnes de $r+1$ \`a $p$ s'expriment en fonction des $r$ premi\`eres colonnes. On a donc:
 \be
	\begin{pmatrix}{
	A_{12} \cr
		A_{22} }
\end{pmatrix}=\begin{pmatrix}{
	A_{11} \cr
A_{21} }
\end{pmatrix}H\Longrightarrow \left\{
\begin{array}{ll}
	A_{12}=A_{11}H  \\
	A_{22}=A_{21}H 
\end{array} \right.\Rightarrow \left\{
\begin{array}{ll}
	H=A_{11}^{-1}A_{12} \\
	A_{22}=A_{21}A_{11}^{-1}A_{12}   
\end{array}\right.
\ee
\bdf
Soit une matrice $A$ quelconque, $A^{-1}$ inverse g\'en\'eralis\'ee de $A$ si $AA^{-1}A=A$. (\textit{A. Bjerhammer}, 1973).\index{\textbf{Bjerhammer A.}} 
\edf
\subsection*{19.16.2. Existance}
Soit une matrice $A$ de rang $r$, d'apr\`es le th\'eor\`eme des pivots, on peut \'ecrire $A$ sous la forme donn\'ee par (\ref{ac0}), alors une inverse g\'en\'eralis\'ee de $A$ est donn\'ee par: 
\be
	A^{-1}_0=\begin{pmatrix}{
	A_{11}^{-1} & 0 \cr
	0 & 0 }
\end{pmatrix}
\ee
\textbf{V\'erification}. On a:
\ba
&	AA^{-1}_0A=\begin{pmatrix}{
	A_{11} & A_{12} \cr
	A_{21} & A_{22} }
\end{pmatrix}\begin{pmatrix} {
	A_{11}^{-1} & 0 \cr
	0 & 0 }
\end{pmatrix}\begin{pmatrix}{
	A_{11} & A_{12} \cr
	A_{21} & A_{22} }
\end{pmatrix}=\begin{pmatrix}{
	A_{11} & A_{12} \cr
	A_{21} & A_{22} }
\end{pmatrix}\begin{pmatrix}{
	I & A_{11}^{-1}A_{12} \cr
	0 & 0 }
\end{pmatrix} \nonumber & \\
&=\begin{pmatrix} {
	A_{11} & A_{11}A_{11}^{-1}A_{12} \cr
	A_{21} & A_{21}A_{11}^{-1}A_{12} }
\end{pmatrix}= \begin{pmatrix}{
	A_{11} & A_{12} \cr
	A_{21} & A_{22} }
\end{pmatrix}=A &
\ea
\textbf{Propriété 19.4} \textit{Si $A_0^{-1}$ est une inverse g\'en\'eralis\'ee, les autres matrices inverses g\'en\'eralis\'ees s'expriment comme suit:}
\be
	A^{-1}=A_0^{-1}+(I-A_0^{-1}A)M+N(I-AA_0^{-1})
\ee
On v\'erifie ais\'ement que $AA^{-1}A=A$.
\\

\textbf{Propriété 19.5} \textit{Les matrices $AA^{-1}$ et $A^{-1}A$, $I-AA^{-1}$ et $I-A^{-1}A$ sont des matrices carr\'ees et nilpotentes c'est-\`a-dire $M.M=M$.}
\subsection*{19.17.3. Les Syst\`emes Lin\'eaires}
Soit le syst\`eme lin\'eaire suivant:
\be
	AX=K    \label{ss1}
\ee
avec $A$=$_{n}A_p$, $X$=$_{p}X_q$ et $K$=$_{n}K_q$.
\\

Une condition n\'ecessaire et suffisante pour l'existance de solutions de (\ref{ss1}) est:
\be
	AA^{-1}K=K
\ee
Les solutions de (\ref{ss1}) sont donn\'ees par:
\be
	X=A_0^{-1}K+(I-A_0^{-1}A)M
\ee
o\`u $A_0^{-1}$ est une inverse g\'en\'eralis\'ee particuli\`ere et $M$ une matrice arbitraire.
\\

%
\textbf{\un{Note Historique:}} \textsl{La méthode des moindres carrés fut publiée par la première fois par Adrien-Marie Legendre\footnote{\textbf{Adrien-Marie Legendre} (1752-1833): mathématicien et géodésien français.} en 1809. La justification comme procédure statistique de la méthode des moindres carrés fut donnée par Carl Friedrich Gauss en 1809, puis en 1810 dans son mémoire sur l'astéroïde Pallas découvert par Heinrich Wilhelm Olbers\footnote{\textbf{Heinrich Wilhelm Olbers} (1758-1840): astronome et physicien  allemand.} le 28 mars 1802. Selon Gauss, la méthode des moindres carrés conduit à la meilleure combinaison possible des observations quelle que soit la loi de probabilité  des erreurs. Elle fut immédiatement reconnue comme une contribution majeure. Gauss affirma l'avoir déjà utilisée dès 1795. Il est certain qu'il s'en servit en 1801 pour déterminer l'orbite de la comète Cérès découverte par Giuseppe Piazzi\footnote{\textbf{Giuseppe Piazzi} (1746-1826): astronome et mathématicien italien.} le 1er janvier 1801. }\index{\textbf{Piazzi G.}}\index{\textbf{Legendre A.M.}}\index{\textbf{Gauss C.F.}}\index{\textbf{Olbers H.W.}}
\\

\textsl{Le mathématicien américain d'origine irlandaise Robert Adrain\footnote{\textbf{Robert Adrian} (1775-1843): mathématicien américain.}  avait, à l'occasion d'une question de topographie, publié  un article daté  de 1808 (mais paru en 1809) dans lequel il exposait également la méthode des moindres carrés. Ce travail passa totalement inaperçu en Europe. En 1818, Adrain appliqua encore cette méthode à  la détermination de l'aplatissement de la Terre à  partir de mesures du méridien et en tira une estimation des axes de
l'ellipsoïde terrestre.}\index{\textbf{Adrain R.}}
\\ 

\textsl{Parmi les méthodes de calcul de l'inversion de la matrice normale, on cite la méthode dite de Cholesky\footnote{\textbf{André-Louis Cholesky} (1875-1918): ingénieur polytechnicien et géodésien militaire français.}. Ce dernier a été le Chef du Service Topographique Tunisien entre mai 1913 et août 1914.}  (C. Brezinski, 2005, § \ref{biblio2})\index{\textbf{Cholesky A.L.}}\index{\textbf{Brezinski C.}}
\section{\textsc{Exercices et Problèmes}}
\bex
Soit un triangle $ABC$, on observe les angles $\hat{A},\,\hat{B},\,\hat{C}$ et les côtés $BC=a,\, AC=b$ et $AB=c$:
	\[
\left\{\begin{array}{l}
\hat{A}=43.7716\,0\,gr\,\quad \sigma_{\hat{A}}=3.1\,dmgr\\
\hat{B}=98.3904\,3\,gr\,\quad \sigma_{\hat{B}}=3.1\,dmgr\\
\hat{C}=57.8385\,8\,gr\,\quad \sigma_{\hat{C}}=3.1\,dmgr\\
a=333.841\,m,\quad \sigma_a=0.005\,m \\
b=525.847\,m,\quad \sigma_b=0.010\,m\\
c=414.815\,m,\quad \sigma_c=0.005\,m 
\end{array}\right. \nonumber
\]
1. Calculer les angles et les côtés compensés.

2. Calculer les poids de l'angle $\hat{A}$ et du côté $a$.

3. Déterminer une estimation du facteur de variance unitaire.
\eex
\bpb
Les directions suivantes sont observées respectivement aux stations $A,B,C$ et $D$ d'un quadrilatère $ABDC$ comme suit:
 \ba
Station\,\,A=\left\{\begin{array}{l}
vers\,\,B: \,\,0.0000\,0\,gr \\
vers\,\,C: 74.1666\,7\,gr  
\end{array}\right. \nonumber \\
Station\,\,B=\left\{\begin{array}{l}
vers\,\,D: \,\,0.0000\,0\,gr \\
vers\,\,C: 82.4608\,0 \,gr \\
vers\,\,A: 170.6253\,1\,gr 
\end{array}\right. \nonumber \\
Station\,\,C=\left\{\begin{array}{l}
vers\,\,A: \,\,0.0000\,0\,gr \\
vers\,\,B: 37.6709\,9\,gr \\
vers\,\,D:85.0830\,2\,gr 
\end{array}\right. \nonumber \\
Station\,\,D=\left\{\begin{array}{l}
vers\,\,C: \,\,0.0000\,0\,gr \\
vers\,\,B: 70.1280\,9\,gr 
\end{array}\right. \nonumber
 \ea
Les observations sont non corrélées. l'écart quadratique moyen de ces observations est identique et vaut $\sigma_d=6.2\,dmgr$. 

1. Compenser les directions et calculer leurs poids et celui de l'angle $CBA$.

2. Calculer l'estimateur $s^2$ du facteur de variance unitaire et celui de $\ds \frac{s^2}{\sigma^2}$.

3. Des observations de nivellement ont été effectuées sur les lignes $ABC$ et $BCD$. Les différences d'altitudes observées sont les suivantes:
\ba
H_A-H_B=0.509\,m \nonumber \\ 
H_B-H_D=1.058\,m \nonumber \\ 
H_A-H_C=3.362\,m \nonumber \\ 
H_D-H_C=1.783\,m \nonumber \\ 
H_B-H_C=2.829\,m \nonumber  
\ea 
Les observations sont non corrélées et de précision identique. Compenser les observations ci-dessus et calculer un estimateur du facteur de variance unitaire.
\epb
\bpb
 1. Montrer que dans un cheminement altimétrique de précision, le poids de l'observation entre deux repères est inversement proportionnel de leur distance en supposant l'égalité des portées et que les observations sont non corrélées. 

\begin{figure}
	\centering
		\includegraphics{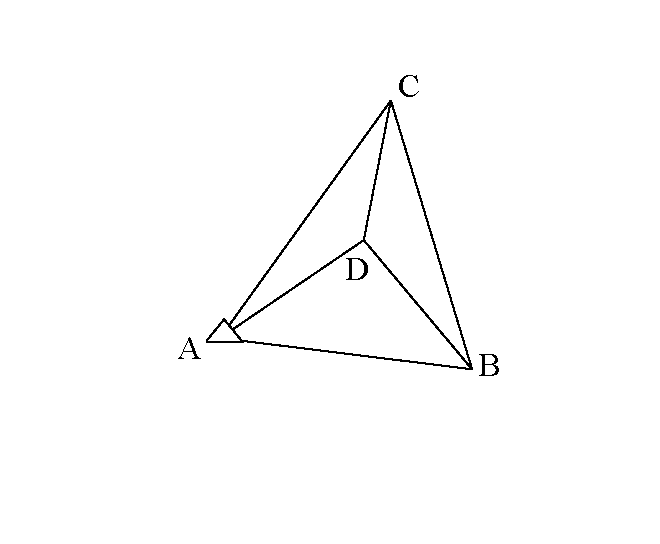}
	\caption{La polygonale observée}
	\label{fig:triangleniv}
\end{figure}

2. Une polygonale $ABCD$ (voir \textbf{Fig. \ref{fig:triangleniv}}) a été observée par le nivellement de précision. L'instrument utilisé a une précision de $2\,mm$ par $km$. Les observations considérées non corrélées sont les suivantes:
\ba
H_C-H_A=1.878\,m, \quad 	AC=6.44\,km \nonumber \\ 
H_D-H_A=3.831\,m, \quad 	AD=3.22\,km \nonumber \\ 
H_D-H_C=1.954\,m, \quad 	CD=3.22\,km \nonumber \\ 
H_B-H_A=0.332\,m, \quad 	AB=6.44\,km \nonumber \\ 
H_D-H_B=3.530\,m, \quad 	BD=3.22\,km \nonumber \\ 
H_C-H_B=1.545\,m, \quad 	BC=6.44\,km \nonumber  
\ea 

L'altitude du repère $A$ est de $3.048\,m$ et non entachée d'erreurs. Calculer par compensation des observations les altitudes des repères $B,C$ et $D$ et leurs écarts-types.

3. Calculer l'écart-type de la différence d'altitudes entre les repères $C$ et $D$.

4. Donner une estimation de la précision par $km$ du nivellement effectué.
\epb
\bpb
		On veut étalonner un anéroide, appareil donnant la pression de l'air, par la formule:
	\[ D=d+\alpha t+\gamma
\]
où $\alpha,\gamma$ sont deux constantes, $t$ la température en degrés centigrades. Les paramètres $d$ et $D$ sont lus respectivement de l'anéroide et à partir d'un baromètre en mercure, et exprimés en $mm$ $Hg$. 
\\

Pour déterminer $\alpha$ et $\gamma$, des lectures ont été prises à différentes températures (voir tableau \ref{tab: Tableau des observations}). 
\begin{table}[ht]
	\centering
\begin{center}
\[\begin{array}{lll}
\hline
 \quad  \quad	t		 & d   &      D   \\ \hline
\mbox{° Centigrade} \,\,      &    mm\,\,Hg \,\,   &              mm\,\,Hg     \\   \hline
 \quad  \quad6.0         &       761.3                 &       762.3              \\   \hline
 \quad  \quad10.0        &      759.1                 &        759.5              \\  \hline 
  \quad  \quad14.0       &    758.4                    &    758.7                  \\ \hline
 \quad  \quad 18.0       &    763.1                    &    763.0                  \\ \hline
\end{array} 
\]
\caption{Table des observations}
	\label{tab: Tableau des observations}
\end{center}
\end{table}
Ces observations sont non corrélées. L'écart-type de la lecture de $d$ est de $0.14\,mm\,\,Hg$; $t$ et $D$ sont supposées sans erreurs.

1. Calculer par la méthode des moindres carrés les constantes $\alpha$ et $\gamma$.

2. Estimer le facteur de variance unitaire.

3. Déterminer la variance et la covariance de $\alpha$ et $\gamma$.
\epb
\bpb
On considère le réseau suivant constitué des points $P_1,P_2,P_3$. Le point $P_1$ est considéré fixe.
\begin{table}[htp]
	\centering
\begin{center}
\[\begin{array}{lll}
\hline
 \quad  \mbox{Nom du point} 		 & \varphi \,(gr)   &    \lambda\,(gr)   \\ \hline
 \quad  P_1         &       37.0830\,6094                 &       11.5451\,6843             \\   \hline
 \quad  P_2          &      37.1229\,0536                 &       11.2861\,5241              \\  \hline 
  \quad  P_3         &    37.0542\,4612                 &       11.4288\,7620                \\ \hline
\end{array} 
\]
\caption{Table des coordonnées}
	\label{tab: Tableau des coordonnées}
\end{center}
\end{table}
Les coordonnées des points sont dans le système Carthage34 dont l'ellipsoïde de référence est celui de Clarke Français 1880 ($a= 6378\,249.20\,m; e^2=0.0068\,034877)$. On veut calculer la compensation des observations dans le plan en utilisant les coordonnées Lambert Sud Tunisie. On prendra comme coefficient d'échelle nominale $k_S=1.0000$. On rappelle que pour la représentation Lambert Sud Tunisie: $\varphi_0=37\,gr$ Nord; $\lambda_0=11\,gr$ Est de Greenwich.

1. Calculer les coordonnées planes $(X_1,Y_1)$ du point fixe $P_1$.

2. Calculer les coordonnées approchées $(X,Y)$ des autres points.

3. On donne ci-dessous les tours d'horizon effectués aux points du réseau:
 \ba
Station\,\,P_1=\left\{\begin{array}{l}
vers\,\,P_2: \,\,91.0566\,34\,gr \\
vers\,\,P_3: 61.1946\,66\,gr \\
\end{array}\right. ;\quad Station\,\,P_2=\left\{\begin{array}{l}
vers\,\,P_1: \,329.2207\,96\,gr \\
vers\,\,P_3:350.8497\,39\,gr\\
\end{array}\right. \nonumber \\
Station\,\,P_3=\left\{\begin{array}{l}
vers\,\,P_1: \, 236.9308\,57\,gr \\
vers\,\,P_2: 88.4215\,38\,gr \\
\end{array}\right. \nonumber 
 \ea
4. On mesure la distance suivant l'ellipsoïde $P_1P_3=10\,156.963\,m$. Calculer $P_1P_3$ réduite au plan de la représentation plane.

5. Donner la liste des inconnues à déterminer.

6. Ecrire l'équation matricielle du système des moindres carrés et déterminer la matrice des coefficients et le vecteur des observables. On prendra comme matrice de poids la matrice unité $I$.

7. Calculer la matrice normale.

8. Résoudre le système de la solution des moindres carrés.
 
9. Exprimer les coordonnées planes définitives et en déduire leurs coordonnées géographiques définitives.
\epb
\bpb
En statistiques, la loi normale est une famille de distributions de probabilités caractérisées par la fonction de densité:
$$p(x,\mu, \sigma)=\frac{1}{\sqrt{2\pi}\sigma}e^{-\frac{(x-\mu)^2}{2\sigma^2}}$$
où $\mu$ est la moyenne et $\sigma^2$ la variance. On note par $l(x,\mu,\sigma)=Logp(x,\mu,\sigma)$, soit: $$ l(x,\mu,\sigma)=-Log\sigma-\frac{(x-\mu)^2}{2\sigma^2}$$
Soit $X$ une variable aléatoire ayant comme fonction de densité $p(x,\mu,\sigma)$. On rappelle les opérateurs suivants espérance mathématique ou moyenne et variance:
\ba
E[f(X)]=\int_{-\infty}^{+\infty}f(x)p(x,\mu,\sigma)dx\nonumber \\
V(f(X))=E[(E[f(X)]-f(X))^2]\nonumber
\ea
On donne la formule: $\ds \int_0^{+\infty}e^{-u^2}du=\frac{\sqrt{\pi}}{2}$.

1. Montrer que:
\ba
E(X)=\ds \int^{+\infty}_{-\infty}p(x,\mu,\sigma)dx=\mu \nonumber \\
\sigma^2(X)=Var(X)=Cov(X,X)=\ds \int^{+\infty}_{-\infty}(x-\mu)^2p(x,\mu,\sigma)dx=\sigma^2 \nonumber
\ea 
2. Montrer que:
$$\ds \int_{-\infty}^{+\infty}u^4e^{-u^2}du=\ds \frac{3\sqrt{\pi}}{4}$$

\noindent 3. Calculer $\ds \frac{\partial l}{\partial \mu}\,\,\frac{\partial l}{\partial \sigma}$.
\\

\noindent 4. On pose $\theta=(\mu,\sigma)$. Soit $T_\theta$ l'espace engendré par $(\ds \frac{\partial l}{\partial \mu},\,\,\frac{\partial l}{\partial \sigma})$. On définit sur $ T_\theta$ l'opérateur $<.,.>: T_\theta \times T_\theta\longrightarrow \BbR$ à $A,B$ deux variables aléatoires $\in T_\theta$:
$$<A,B>=E[A(x)B(x)]$$
Justifier qu'on peut écrire:
$$E[A(x)B(x)]=Cov(A(x),B(x))=E[(E[A(x)]-A(x))(E[B(x)]-B(x))]$$
5. Montrer que $<.,.>$ définit un produit scalaire sur $T_\theta$.
\\

\noindent 6. On pose: $e_1=\ds \frac{\partial l}{\partial \mu}$ et $e_2=\ds \frac{\partial l}{\partial \sigma}$. On définit le tenseur métrique sur $T_\theta$ par:
$$g_{ij}=<e_i,e_j>$$
Montrer que la matrice $g=(g_{ij})$ est donnée par: 
$$g=\ds \frac{1}{\sigma^2}\begin{pmatrix}{
1 & 0 \cr
0 & 2 }
\end{pmatrix}$$
et que la première forme fondamentale sur $T_\theta$ s'écrit:
$$ ds^2=\ds \frac{1}{\sigma^2}(d\mu^2+2d\sigma^2)$$
\epb
\bpb
Soit un triangle de côtés $a,b,c$ et d'angles $A,B$ et $C$. On se propose:

- d'estimer $\dot{a},\dot{b}$ et $\dot{c}$,  et les variances de ces déterminations. Les observations sont:
\be
\left\{\begin{array}{l}
   a=96.48\,mm \\
	b=115.50\,mm \\
 	A=63.042 \,gr  \\
  B=99.802\,gr     \\
   C=37.008\,gr    
	\end{array}\right. \lb{AF21}
 \ee
 On choisit ici comme \un{unités normalisées} le décimillimètre $(0.1\,mm)$ pour les mesures de distances, et le décimilligrade $(0.1\,gr)$ pour les angles. 
\\

On prend les poids égaux aux inverses des carrés des $emq$ de chaque observation. On donne la matrice des poids $P$:
$$ P=\begin{pmatrix}{
0.277 & 0&0&0&0 \cr
0&0.160&0&0&0 \cr
0&0&1.524&0&0 \cr
0&0&0&1.524 &0 \cr
0&0&0&0&1.524 }
\end{pmatrix}$$
On prendra comme valeurs approchées des inconnues $ a_0=a;\quad b_0=b;\quad c_0=\ds a\frac{sinC}{sinA}$.
\\

1. Ecrire les paramètres observées et les valeurs observées des inconnues dans les nouvelles unités. 

2. Soit $X=(a,b,c)$ le vecteur des inconnues. On adopte le système suivant liant les inconnues aux observables:
\be
\left\{\begin{array}{l}
\dot{a}=\dot{a} \\
\dot{b}=\dot{b} \\
Arccos\ds \frac{\dot{b}^2+\dot{c}^2-\dot{a}^2}{2\dot{b}\dot{c}}=\dot{A} \\
\\
 Arccos\ds \frac{\dot{c}^2+\dot{a}^2-\dot{b}^2}{2\dot{c}\dot{a}}=\dot{B} \\
\\
Arccos\ds \frac{\dot{a}^2+\dot{b}^2-\dot{c}^2}{2\dot{a}\dot{b}}=\dot{C} 
\end{array}\right. \lb{BF21-13}
\ee
Ceci étant, on posera pour les grandeurs à déterminer:
\ba
\dot{a}=a_0+da=a+da \nonumber \\
\dot{b}=b_0+db=b+db \nonumber \\
\dot{c}=c_0+dc \nonumber
\ea
et pour les grandeurs observées:
\ba
\dot{a}=a+v_a;\quad \dot{b}=b+v_b;\quad \dot{A}=A+v_A\nonumber \\
\dot{B}=B+v_B; \quad \dot{C}=C+v_C \nonumber 
\ea
En linearisant la troisième équation de (\ref{BF21-13}), montrer que l'équation d'observation s'écrit:
$$\frac{1}{sinA}\frac{a_0}{b_0c_0}\frac{2000}{\pi}da-\frac{1}{sinA}\frac{a_0^2+b_0^2-c_0^2}{2b_0^2c_0}\frac{2000}{\pi}db-\frac{1}{sinA}\frac{a_0^2+c_0^2-b_0^2}{2b_0c_0^2}\frac{2000}{\pi}dc=-k_A\frac{2000}{\pi}+v_A$$
où :
$$k_A=\frac{b^2_0+c^2_0-a^2_0-2b_0c_0cosA}{2b_0c_0sinA}$$
(étant entendu qu'on exprime $v_A$ en $dcgr$).
\\

3. Montrer que le système des moindres carrés $ AX=L+V$ s'écrit:
$$
\begin{pmatrix}{
1. & 0.& 0. \cr
0. & 1. & 0.& \cr
1.00375 & -0.83924 & 0.00143 \cr
-1.00571 & 1.20285 & -0.66128 \cr
0.00094 & -0.36239 & 0.65918 }
\end{pmatrix}.\begin{pmatrix} {
da \cr
db \cr
dc }
\end{pmatrix}=\begin{pmatrix}{
0. \cr
0. \cr
0.97981 \cr
-2.88449 \cr
0.42396}
\end{pmatrix}+\begin{pmatrix}{
v_a \cr
v_b \cr
v_A \cr
v_B \cr
v_C }
\end{pmatrix} $$
4. Résoudre le système précédent par la méthode des moindres carrés et montrer que la matrice normale $N=A^TPA$ est donnée par:
$$ N=\begin{pmatrix} {
3.35605 & -3.13044 & 1.01750 \cr 
-  & 3.64132 & -1.57937 \cr
- & - & 1.32971 }
\end{pmatrix}$$
5. Montrer que:
$$ X=\begin{pmatrix}{
+0.62971\cr
-0.90962\cr
 0.94782 }
\end{pmatrix}$$
6. Déterminer les variances des inconnues $\sigma^2_a,\,\sigma^2_b$ et $\sigma^2_c$.
\epb

\chapter{\textit{\textbf{Pr\'esentation des Aspects Théoriques de la Géométrie de la Compensation non-Linéaire par les Moindres Carrés} }}



	L'objet de ce chapitre est de pr\'esenter les aspects th\'eoriques de la g\'eom\'etrie de la compensation non-lin\'eaire par les moindres carr\'es. Il est inspir\'e du travail r\'ealis\'e par P.J.G. Teunissen (\textit{P.J.G. Teunissen}, 1985).\index{\textbf{Teunissen P.J.G.}}
\section{\textsc{Introduction}}
La th\'eorie des moindres carr\'es telle qu'elle a \'et\'e d\'efinie par Gauss repose sur la lin\'earisation. On arrive au mod\`ele de Gauss-Markov\footnote{\textbf{Andreï Markov} (1856-1922): mathématicien russe.}\index{Modèle de Gauss-Markov}:\index{\textbf{Markov A.}}
\begin{equation}
	\fbox{ $ A.X=L+V $}\label{tu1}
\end{equation}
avec:

- $X$ le vecteur des inconnues;

- $A$ la matrice des coefficients;

- $L$ le vecteur des observables;

- $V$ le vecteur des r\'esidus.
\\

La compensation des mod\`eles lin\'eaires a \'et\'e l'objet de plusieurs recherches au cours des derni\`eres ann\'ees. Concernant la compensation non-lin\'eaire, la recherche n'a pas \'et\'e aussi d\'evelopp\'ee et la bibliographie afférente a ce sujet est limit\'ee.

\section{\textsc{Pr\'esentation du Probl\`eme}}
On consid\`ere une application non lin\'eaire: 
\be
y:M \rightarrow N \,\,\mbox{qui \`a}\,\, x  \in M  \rightarrow  y=y(x) \in N \label{tu2}
\ee            
o\`u $M$ et $N$ deux vari\'et\'es g\'eom\'etriques lisses munies respectivement d'un produit scalaire ou encore d'une m\'etrique $<.,.>_M$ et $<.,.>_N$. On appelle :
\begin{equation}
	 \tilde{N}=Im(M)=\left\{y \in N\,/ \exists x \in M\, \mbox{unique avec}\quad y=y(x)\right\}  \label{tu3}
\end{equation}
 
Le probl\`eme de la compensation non lin\'eaire peut \^etre divis\'e en deux sous probl\`emes:

1. le premier probl\`eme est de trouver les estimateurs $\hat{x}$ et $\hat{y}$ respectivement de $x$ et $y$ v\'erifiant la condition (\ref{tu5}) ci-dessous;

2. le deuxi\`eme sous probl\`eme est d'obtenir les propri\'et\'es statistiques des deux estimateurs pr\'ec\'edemment cit\'es.

Afin de pouvoir r\'esoudre ces deux questions importantes, on se limite au cas de la r\'esolution du probl\`eme par les estimateurs des moindres carr\'es. Pour cel\`a, on suppose que l'application $y$ est injective c'est-\`a-dire:
\be
	\forall x_1,\,x_2 \in M \,\mbox{avec}\,x_1\neq x_2\Rightarrow y(x_1)\neq y(x_2) \label{tu4}
\ee         
ou encore, l'ant\'ec\'edent d'une image est un seul point.
\\
 
La condition des moindres carr\'es est :
\be
	\min_{x \in M}2E(x)=\min_{y\in \tilde{N}=y(M)} ||y_s-y||^2_N =\min_{y\in \tilde{N}=y(M)} <y_s-y(x),y_s-y(x)>_N \label{tu5}
\ee           
o\`u on a not\'e  $y_s$ le vecteur des observables.
\\

On appelle $P:x\in M\,\longmapsto y(x) \,\in N$ et  $Q:y\in \tilde{N}\,\longmapsto x \,\in M$ les applications non lin\'eaires avec: 
\be
	\hat{y}=P(\hat{x}) \quad \mbox{et}\quad \hat{x}=Q(\hat{y})\quad \mbox{avec} \quad P\circ Q=I \label{tu7}
\ee
Due \`a la non lin\'earit\'e des applications $P$ et $Q$, il est difficile de trouver une formule ferm\'ee en $y$. La solution est d'appliquer des m\'ethodes it\'eratives:
\vspace{4mm}

- on part d'un point $x_0 \in M$ et on essaye de g\'en\'erer une suite $x_0,x_1,x_2,...$ qui convergera vers le point solution $\hat{x}$. Les m\'ethodes trouv\'ees dans la litt\'erature proc\`edent conform\'ement aux étapes suivantes (\textit{J.M. Ortega \& W.C. Rheinboldt}, 1970):\index{\textbf{Ortega J.M.}}\index{\textbf{Rheinboldt W.C.}}
\begin{equation}
	x_{q+1}^{(\beta)}=x_q^{\beta}+t_q\Delta x_q^{\beta};\quad \beta=1,2,...,n \label{tu7} 
\end{equation}

1. on prend $q=0$ et on donne une premi\`ere valeur $x_0$;
\\

2. on d\'etermine le pas $\Delta x_q$;
\\

3. on d\'etermine un scalaire $t_q$ tel que :
\begin{equation}
	||y_s-y(x_{q+1})||_N \leq ||y_s-y(x_q)||_N \label{tu8}
\end{equation}
cette relation assure la convergence;
\\

4. voir si la convergence est assur\'ee, dans ce cas $x_{q+1}=\hat{x}_1$, sinon, $q\rightarrow q+1$ et on passe \`a l'\'etape 2.
\\

La r\'esolution de (\ref{tu8}) d\'epend du choix de $\Delta x_q $ et $t_q$. 
\\

La m\'ethode it\'erative comprend deux cat\'egories:

- la premi\`ere est une m\'ethode it\'erative directe qui ne n\'ecessite pas l'utilisation des donn\'ees de la fonction $E(x)$ (\ref{tu5});

- la deuxi\`eme cat\'egorie concerne la m\'ethode it\'erative du gradient qui fait appel aux calculs des d\'eriv\'ees de l'expression $E(x)$.
\newpage
Concernant la d\'etermination des propri\'et\'es statistiques, elle n'est pas aussi facile comme pour le cas lin\'eaire. On a peu d'information sur une m\'ethodologie pour traiter cette question.
\section{\textsc{Les Eléments Mathématiques de l'Optimisation des Problèmes Non-Linéaires }}
Les problèmes non-linéaires sont rencontrés dans divers phénomènes physiques tels que le mouvement des fluides, l'élasticité, la relativité, la géodésie, et autres. On présente ci-après quelques définitions et théorèmes relatives à ces problèmes. On s'intéresse à l'étude des fonctions d'un domaine de $\BbR^n \longrightarrow \BbR$. 
\\

Soit $\Omega \subset \BbR^n$ et $E$ une fonction: $x\in \Omega \longrightarrow E(x) \in \BbR$ avec $x=(x_1,...,x_n)^T$. D'où les définitions:

\bdf
Soit $E$ la fonction à optimiser définie sur $\Omega \subset \BbR^n$, à $x\in \Omega \longrightarrow E(x) \in \BbR$. On appelle $E$ la fonction \un{objectif} ou \un{énergie}.\index{Fonction objectif ou énergie}
\edf
\bdf
Un point $x^* \in \Omega$ est dit un \un{minimum global} de la fonction énergie $E$ si:\index{Minimum global}
\be
\fbox{ $ E(x^*)\leq E(x)\quad \forall x \in \Omega $}
\ee 
\edf
\bdf
Un point $x^* \in \Omega$ est dit un \un{minimum global strict} de la fonction énergie $E$ si:\index{Minimum global strict}
\be
\fbox{ $ E(x^*)< E(x)\quad \forall x \in \Omega\,\, x\neq x^* $}
\ee 
\edf
Si le domaine $\Omega$ est un compact de $\BbR^n$ d'où le théorème suivant qu'on énonce sans démonstration: 
\bthm
Si $E:\Omega \longrightarrow \BbR$ est continue où $\Omega \subset \BbR^n$ est un compact ( c'est-à-dire $\Omega$ est un fermé borné), alors la fonction $E$ admet au moins un minimum global $x^* \in \Omega$.
\ethm
Pour définir le cas du maximum relatif ou global, on considère l'inégalité ($E(x)\leq E(x^*)$ ou $E(x)<E(x^*)$).
\\

\subsection*{20.3.1. Dérivée Directionnelle}
Soit $E:\Omega \subset \BbR^n\longrightarrow \BbR$, la définition standard du \textit{gradient} de $E$ est le vecteur noté:
\be
\nabla E(x)=\begin{pmatrix}{
\ds \frac{\partial E}{\partial x_1} \cr
\vdots \cr
\ds \frac{\partial E}{\partial x_n} }
\end{pmatrix}=\left(\frac{\partial E}{\partial x_1},\frac{\partial E}{\partial x_2},\cdots,  \frac{\partial E}{\partial x_n} \right)^T
\ee
\bdf
Soit $V$ un espace vectoriel réel muni d'un produit scalaire qu'on note $<\textbf{.},\textbf{.}>$. Le gradient d'une fonction $E: V\longrightarrow \BbR$ en un point $x \in V$ est le vecteur $\nabla E(x)\in V$ qui vérifie:
\be
\fbox{ $ <\nabla E(x),v>=\ds \frac{d E}{dt}(x+tv) \ds \biggr |_{t=0} \quad \forall v \in V $}
\ee
\edf
La définition telle qu'elle est donnée ci-dessus est appelée \textit{la dérivée directionnelle de $E$ suivant la direction du vecteur $v$}.
\\

On suppose qu'on est dans le cas d'un espace vectoriel $V$ muni du produit scalaire euclidien. Soient $E(x)=E(x_1,...,x_n)$ et $v=(v_1,...,v_n)^T$, on a donc:
\be
 \frac{d E}{dt}(x+tv)=\frac{d E}{dt}(x_1+tv_1,x_2+tv_2,...,x_n+tv_n)=\frac{\partial E}{\partial x_1}(x+tv).v_1+...+ \frac{\partial E}{\partial x_n}(x+tv).v_n
\ee
Pour $t=0$; on obtient:
\be
\frac{d E}{dt}(x+tv) \ds \biggr |_{t=0}=\frac{\partial E}{\partial x_1}(x).v_1+...+ \frac{\partial E}{\partial x_n}(x).v_n=\nabla E(x).v=<\nabla E(x),v>
\ee
C'est le produit scalaire du vecteur $v$ avec le gradient défini usuellement dans $\BbR^n$.

Dans la suite de ce chapitre, on considère que les fonctions \textit{objectifs} ou \textit{énergies} sont des fonctions de classe $C^2$ deux fois continûement différentiable dans leurs domaines de définition. 
\subsection*{20.3.2. Propriété du gradient}
On considère que le point $x$ se mouve le long d'une courbe $\Gamma$ paramétrée par une variable $t$, soit $x(t)$. Alors $E$ devient une fonction de $t$:
\be
 \m E(t)=E(x(t))=E(x_1(t),x_2(t),...,x_n(t))
\ee
En différentiant l'équation  précédente, on obtient:
\ba
&d\m E(t)=dE(x)=gradE.dx=gradE.\ds \frac{dx(t)}{dt}dt=<\nabla E(x),\ds \frac{dx}{dt}>dt \Longrightarrow \nonumber& \\&
\fbox{ $ \ds \frac{d}{dt}E(x(t))=<\nabla E(x),\ds \frac{dx}{dt}> $}\lb{eq2017}&
 \ea
C'est la dérivée directionnelle suivant le vecteur $v=\dot{x}=\ds \frac{dx}{dt}$ tangent à la courbe $\Gamma$ en $x(t)$. Si $\ds \frac{dE}{dt} >0$, on est dans le sens où la fonction $E$ croît, par contre si $\ds \frac{dE}{dt} <0$, on est dans le sens où la fonction $E$ décroît. Alors si $dE=0\Longrightarrow \nabla E(x)=0$, on dit que $x$ est un point critique d'où la définition:
\bdf
Un point $x^*$ est un point critique de la fonction objectif $E(x)$ si:
$$ \fbox{ $ \nabla E(x^*)=0 $} $$\index{Point critique}
\edf
Par suite, on a le théorème suivant:
\bthm
Tout minimum relatif $x^*$ à l'intérieur d'un domaine $\Omega$ de la fonction $E(x)$ est aussi un point critique. \index{Minimum relatif}
\ethm
\textbf{Démonstration :} Soient $x^*$ un minimum relatif $\in \Omega \subset \BbR^n$ et un vecteur $v\neq0 \in \BbR^n$, et on considère la fonction  scalaire :
\be
g(t)=E(x^*+tv)=E(x_1^*+tv_1,...,x_n^*+tv_n)
\ee
où $t \in \BbR$ tel que $x=x^*+tv \in \Omega$. Comme $x^*$ est minimum relatif, donc $ E(x^*)\leq E(x^*+tv)$ ce qui implique $g(0)\leq g(t)$, soit $g(0)$ est un minimum pour la fonction scalaire, dans ce cas, on a nécessairement $g'(t)=0$ pour $t=0$. Or d'après la formule (\ref{eq2017}): 
$$ g'(t)=\frac{dg}{dt}=<\nabla E(x),\ds \frac{dx}{dt}> $$
Pour $t=0$, on obtient:
\be
g'(0)=0=\ds \frac{d E}{dt}(x^*+tv) \ds \biggr |_{t=0}=<\nabla E(x),\ds \frac{dx}{dt}>\biggr |_{t=0}=<\nabla E(x^*),v>=0 \Rightarrow  \nabla E(x^*)\perp v
\ee
Comme $v$ est arbitraire, par suite on a:
\be
 \nabla E(x^*)=0
\ee
Donc $x^*$ est un point critique. On note que la reciproque n'est pas toujours vraie. On peut avoir $\nabla E(x^*)=0$ sans que $x^*$ soit un point minimum relatif ou global.
\subsection*{20.3.3. Etude de la Dérivée Seconde}
Dans le cas d'une fonction réelle $f$ d'une seule variable, telle que $f'(x_0)=0$ au point $x_0\in \BbR$, l'étude de $f"(x_0)$ peut donner une information sur le statut du point critique $x_0$, d'où la proposition (\textit{P.J. Olver}, 2013)\index{\textbf{Olver P.J.}}:
\bpr \lb{prop201}
Soit $g(t)$ une fonction scalaire réelle de classe $C^2$ et supposons que $t^*$ est un point critique de $g$ c'est-à-dire $g'(t^*)=0$:

- si $t^*$ est un minimum relatif, alors $g"(t^*)\geq 0$;

- si $g"(t^*)>0$, alors $t^*$ est un minimum relatif strict.\index{Minimum relatif strict}
\epr
La  démonstration utilise l'expression de $g(t)$ en utilisant le développement de Taylor au voisinage de $t^*$ soit :
\ba 
&g(t)=g(t^*)+(t-t^*)g'(t^*)+\frac{1}{2}(t-t^*)^2g"(t^*)  &\\
&t^*\,\,point \,\,critique\Longrightarrow g'(t^*)=0\Longrightarrow g(t)=g(t^*)+\frac{1}{2}(t-t^*)^2g"(t^*)\Longrightarrow g"(t)=g"(t^*)\nonumber &
\ea
 Si considère maintenant le cas de $n$ variables pour une fonction $E(x)=E(x_1,...,x_n)$, la dérivée seconde est représentée par une matrice $n\times n$ dite le Hessien \footnote{Nommé en hommage au mathématicien allemand \textbf{Ludwig Otto Hesse} (1811-1874).}:\index{\textbf{Hesse L.O.}}
\be
\fbox{ $ \nabla^2E(x)=\ds \begin{pmatrix}{
\ds \frac{\partial^2 E}{\partial x_1^2}&\ds \frac{\partial^2 E}{\partial x_1 \partial x_2} & \cdots & \ds \frac{\partial^2 E}{\partial x_1 \partial x_n} \cr
 \ds \frac{\partial^2 E}{\partial x_2 \partial x_1}&\ds \frac{\partial^2 E}{\partial x_2^2} & \cdots &\ds  \frac{\partial^2 E}{\partial x_2 \partial x_n} \cr
\vdots & \vdots & \ddots & \vdots \cr
 \ds \frac{\partial^2 E}{\partial x_n \partial x_1}&\ds \frac{\partial^2 E}{\partial x_n \partial x_2} & \cdots & \ds \frac{\partial^2 E}{\partial x_n^2} }
 \end{pmatrix} $}
\ee
Comme $E$ est supposée de classe $C^2\Longrightarrow \ds \frac{\partial^2 E}{\partial x_i \partial x_j}=\frac{\partial^2 E}{\partial x_j \partial x_i}$, par suite la matrice hessienne ci-dessus est symétrique:
$$ \nabla^2E(x)=(\nabla^2E(x))^T $$
On a vu d'après la proposition (\ref{prop201}), que pour les fonctions scalaires à une seule variable de classe $C^2$, que si la dérivée seconde est strictement positive en un point $t$, alors ce point est un minimum relatif strict donc un minimum relatif. Cependant pour les fonctions scalaires de classe $C^2$ à plusieurs variables, on a le théorème suivant (\textit{H. Cartan}\footnote{\textbf{Henri Cartan} (1904-2008): un des grands mathématiciens français du XXème siècle, fils du célèbre mathématicien Elie Cartan.}, 1979)\index{\textbf{Cartan H.}} :
\bthm
Soit $E(x)=E(x_1,x_2,...,x_n)$ de classe $C^2(\Omega)$ où $\Omega$ est un ouvert de $\BbR^n$:

- si $x^*$ est un minimum relatif de $E$, alors il est nécessaire un point critique, c'est-à-dire $\nabla E(x^*)=0$;

- de plus, la matrice hessienne $\nabla^2E(x)$ est semi-définie positive au point $x^*$, c'est-à-dire $(\forall X \in \Omega, \,\, X^T.\nabla^2E(x^*).X \geq 0)$; \index{Matrice semi-définie positive}

- inversement, si $x^*$ est un point critique avec la matrice hessienne $\nabla^2E(x^*)$ définie positive $(\forall X\neq 0 \in \Omega, \,\, X^T.\nabla^2E(x^*).X > 0)$, alors $x^*$ est un minimum relatif strict de $E$. \index{Matrice définie positive}
\ethm
\section{\textsc{La Méthode de Newton}}\index{La méthode de Newton}
En revenant à la condition des moindres carrés donnée par l'équation (\ref{tu5}):
$$ 	\min_{x \in M}2E(x)=\min_{y\in \tilde{N}=y(M)} ||y_s-y||^2_N =\min_{y\in \tilde{N}=y(M)} <y_s-y(x),y_s-y(x)>_N $$
on prend :
\be 
\begin{array}{lllll}
M=\Omega \subset \BbR^n \,\,\mbox{ouvert de}\,\,\BbR^n; \\
N=\BbR; \\
X\in \BbR^n\,\,||X||^2=X^T.P.X,\quad P \,\,\mbox{matrice $m\times m$ définie positive;}\\
y=A(x)\quad \mbox{fonction vectorielle de classe $C^2$ non linéaire de }\,\,\BbR^n\longrightarrow \BbR^m; \\
y_s=\,\mbox{vecteur des observables}\,\,\in \BbR^m.
\end{array} 
\ee 
\un{Le problème posé est la recherche de}:
\be
	\min_{x \in \Omega}E(x)=\min_{x\in \Omega} \frac{1}{2}||y_s-A(x)||^2 =\min_{x\in \Omega}\frac{1}{2}(y_s-A(x))^T.P(y_s-A(x)) \lb{eq**}
	\ee
On notera le minimum de l'équation précédente par $x^*$.

\subsection*{20.4.1. Les Méthodes Itératives de Descente}\index{Les méthodes itératives de descente}
Parmi les méthodes de résolution de l'équation (\ref{eq**}), on cite les méthodes itératives de descente. Celà signifie que la fonction $E(x)$ prend des valeurs en décroissant quand l'indice de l'itération croît.

On part d'une valeur $x_0$, on calculera $x_{k+1}$ en fonction de $x_k$ tel que :
\be
\fbox{ $ x_{k+1}=x_k+a_kv_k\quad k=0,1,2,.. $} \lb{eq**1}
\ee
où $x_{k+1},x_k,v_k$ sont des vecteurs de $\BbR^n$ et $a_k$ un réel. On définit alors les étapes suivantes: 
\begin{enumerate}
	\item pour k=0, on part d'une valeur approchée $x_0$ donnée;
	\item on fait le choix d'une direction $v_k$ telle que la fonction $E$ décroît;
	\item on détermine le coefficient $a_k$ tel que : $ E(x_{k+1}) \leq E(x_k)$;
	\item on teste si la convergence est obtenue soit $||x_{k+1}-x_k||< \lambda$ où $\lambda$ une quantité qui sera fixée \textit{a priori} pour la convergence de la méthode d'itération. 
\end{enumerate}
 On note que les méthodes des itérations de descente se distinguent en deux groupes:

- le premier où on n'utilise pas les dérivées partielles de la fonction \textit{énergie} $E$;

- par contre pour le deuxième groupe, on fait usage des dérivées partielles de la fonction $E$, pour le choix de la direction du vecteur $v_k$, c'est la technique dite des gradients. 
\subsection*{20.4.2. Préliminaires}
D'après l'équation (\ref{eq**1}), liant les vecteurs $x_{k+1},x_k,v_k$ et le coefficient $a_k$, on a:
$$x_{k+1}=x_k+a_kv_k$$
\newpage
On parlera de la méthode de descente s'il existe $a_k$ tel que:\index{La méthode des gradients} \index{La méthode de descente}
\be
E(x_{k+1})< E(x_k) \lb{eq**2}
\ee
soit:
\be
E(x_k+a_kv_k)< E(x_k)\lb{eq**3}
\ee
On écrit alors un développement de Taylor de $E(x)$ au point $x_k$:
\be
E(x_{k+1})=E(x_k+a_kv_k)=E(x_k)+a_k\nabla E(x_k)^T.v_k+|a_k|\epsilon(x_k), \quad lim_{x_k\rightarrow 0}\epsilon(x_k)=0 \lb{eq**4}
\ee
ou encore:
\be
E(x_{k+1})-E(x_k)=a_k\nabla E(x_k)^T.v_k+|a_k|\epsilon(x_k)<0 \lb{eq**5}
\ee
Comme $a_k>0$, on peut avoir ce coefficient si:
\be
\nabla E(x_k)^T.v_k< 0  \lb{eq**6}
\ee
Donc on peut choisir un vecteur $v_k$ tel que l'inégalité précédente soit vérifiée par:
\be
v_k=-Q(x_k).\nabla E(x_k) \lb{eq**7}
\ee
où $Q(x_k)$ est une matrice arbitraire $n\times n$ définie positive dépendante de $x_k$. En effet:
\be
\nabla E(x_k)^T.v_k=\nabla E(x_k)^T.(-Q(x_k).\nabla E(x_k))=-\nabla E(x_k)^T.Q(x_k).\nabla E(x_k) < 0 \lb{eq**8}
\ee
car la matrice $Q(x_k)$ est définie positive. 
\\

On revient maintenant à l'équation (\ref{eq**1}), elle s'écrit:
\be
\fbox{ $ x_{k+1}=x_k-a_kQ(x_k).\nabla E(x_k)  $} \lb{eq**9}
\ee
Or le coefficient $a_k$ dépend aussi de $x_k$, alors on peut écrire l'équation (\ref{eq**9}) en posant:
\be
\fbox{ $ \Phi(x)=x-a(x)Q(x)\nabla E(x) $} \lb{eq**10}
\ee
sous la forme:
\be
\fbox{ $ x_{k+1}=\Phi(x_k) $}  \lb{eq**11}
\ee
où $\Phi:\Omega \subset \BbR^n\longrightarrow \BbR^n$. 

Au point de la convergence, on obtient:
\be
\fbox{ $ x^*=\Phi(x^*)  $}\lb{eq**12}
\ee
Il revient donc ce qu'on appelle à la recherche d'un point fixe de la fonction $\Phi$. La méthode itérative concernant la fonction $\Phi$ est dite la méthode itérative du point fixe. Le théorème suivant précise les conditions de convergence vers l'unique solution de l'équation (\ref{eq**12}):
\bthm
(\textbf{Itération du point fixe}): Soit $\Omega \subset \BbR^n$ tel que :

i)- $\forall x\in \Omega ,\, \Phi(x) \in \Omega$;

ii)- $\forall x \in \Omega,\,\, \Phi(x)$ est une fonction continue;

iii)- $\forall x_1,x_2 \in \Omega \,\,\exists c\in \BbR\,\,$ tel que $0\leq c <1$ avec :
\be
||\Phi(x_2)-\Phi(x_1)||\leq c||x_2-x_1|| \lb{eq**13}
\ee
Alors:
\begin{enumerate}
	\item une solution $x^*$ de $x=\Phi(x)$ existe dans $\Omega$;
	\item la solution $x^*$ est unique;
	\item la méthode itérative converge vers $x^*$ c'est-à-dire $lim_{k\rightarrow +\infty}\,x_k=x^*$.
\end{enumerate}
\ethm
Pour une démonstration de ce théorème voir (\textit{P.J.G. Teunissen}, 1990).\index{\textbf{Teunissen P.J.G.}}

Si on peut vérifier la condition iii), on peut appliquer dans ce cas le théorème cité ci-dessus. Si on considère que $\Omega \subset \BbR^n$ est un ouvert convexe, dans ce cas l'application du théorème 18.3 de Taylor- Maclaurin du chapitre § [\ref{ELMMC}] donne:
\be
||\Phi(x_2)-\Phi(x_1)||=\left\|\frac{\partial \Phi}{\partial x}(\overline{x})(x_2-x_1)\right\|\leq \left\|\frac{\partial \Phi}{\partial x}(\overline{x})\right\|.||x_2-x_1|| \lb{eq**14}
\ee
où $\overline{x}=x_1+\mu(x_2-x_1)$ avec $\mu \in [0,1]$. Si on prend :
\be
\fbox{ $ c=max_{x\in \Omega}\left\|\ds \frac{\partial \Phi}{\partial x}(x)\right\|<1 $}
\ee
alors la condition iii) est vérifiée.
\\

On va voir maintenant ce que $E(x)$ et $Q(x)$ doivent vérifier. Partant de l'équation (\ref{eq**10}):
$$\Phi(x)=x-a(x)Q(x)\nabla E(x) $$
on pose:
\be
q(x)=a(x)Q(x)
\ee
$q(x)$ c'est une matrice définie positive, alors l'équation (\ref{eq**10}) s'écrit :
$$ \Phi(x)=x-q(x).\nabla E(x) $$
On prend sa dérivée partielle par rapport à $x$, ce qui donne:
\ba
\frac{\partial \Phi(x)}{\partial x}=\frac{\partial}{\partial x}\left(x-q(x).\nabla E(x)\right)=\frac{\partial}{\partial x}\left(I.x-q(x).\nabla E(x)\right) \nonumber \\
=I-\sum_{j=1}^{j=n}\frac{\partial q_j(x)}{\partial x}.\nabla E(x)-q(x).\frac{\partial ^2 E(x)}{\partial x ^2} \nonumber
\ea
avec $q_j(x)$ les vecteurs colonnes de la matrice $q(x)$ et $I$ la matrice unité $n\times n$. Passant à la norme, d'où:
\be
\left\|\frac{\partial \Phi(x)}{\partial x}\right\|\leq \left\| I -q(x).\frac{\partial ^2E(x)}{\partial x^2}\right\|+\left\|\sum_{j=1}^n\frac{\partial q_j(x)}{\partial x}.\nabla E(x)\right\|
\ee
Donc:
\be
c=max_{x\in \Omega}\left\|\frac{\partial \Phi}{\partial x}(x)\right\|\leq \left\| I -q(x).\frac{\partial ^2E(x)}{\partial x^2}\right\|+\left\|\sum_{j=1}^n\frac{\partial q_j(x)}{\partial x}.\nabla E(x)\right\|
\ee
Si on veut que $c<1$, on prendra :
\be
\left\| I -q(x).\frac{\partial ^2E(x)}{\partial x^2}\right\|+\left\|\sum_{j=1}^n\frac{\partial q_j(x)}{\partial x}.\nabla E(x)\right\| <1
\ee
Sachant que pour la solution $x=x^*$ le vecteur $\nabla E(x^*)=0$ implique $\exists \delta >0$ tel que si $||x-x^*||<\delta$, par continuité de $\nabla E(x)$, le terme:
\be
 \left\|\sum_{j=1}^n\frac{\partial q_j(x^*)}{\partial x}.\nabla E(x^*)\right\|\ll 1
\ee
est infiniment petit devant 1. Revenons à $Q(x)$, on obtient pour que la condition iii) soit vérifiée :
\be
\fbox{ $ \left\|I-a(x)Q(x).\ds \frac{\partial^2 E(x)}{\partial x^2}\right\|<1 $}
\ee
On verra par la suite comment appliquer cette condition.
\subsection*{20.4.3. La Méthode de Newton: Cas d'une fonction réelle}\index{La méthode de Newton}
Supposons qu'on ait une fonction $G:\Omega\subset \BbR\rightarrow \BbR$ de $C^1(\Omega)$ et on cherche la solution de :
\be
G(x)=0\lb{eq**19}
\ee
L'équation de la tangente à la courbe de la fonction $y=G(x)$ en un point $M_k(x_k,G(x_k))$ est donnée par:
$$y-G(x_k)=\left(\frac{dy}{dx}\right)_{x=x_k}(x-x_k)$$
soit:
$$ y-G(x_k)=G'(x_k)(x-x_k) $$
Cette tangente coupe l'axe des $x$ au point $x_{k+1}$ tel que:
\be
0-G(x_k)=G'(x_k)(x_{k+1}-x_k)\Longrightarrow \fbox{ $ x_{k+1}=\ds x_k-\left(\frac{dG(x)}{dx}\right)^{-1}_{x_k}G(x_k) $} \lb{eq**20}
\ee
On suppose évidemment que $G'(x_k)\neq0$. L'équation (\ref{eq**20}) représente la formule de récurrence pour déterminer la solution $x^*$ de $G(x)=0$.
\subsection*{20.4.4. La Méthode de Newton: Cas de la fonction \textit{énergie} $E(x)$}
On rappelle que le problème posé est la recherche de la solution du minimum de $E(x)=\frac{1}{2}||y_s-A(x)||^2$, soit:
\be
\partial_x E(x)=\nabla E(x)=0
\ee
On remplace dans (\ref{eq**19}) $G(x)$ par $\partial_x E(x)$. Alors l'équation de récurrence (\ref{eq**20}) devient pour la fonction $E(x)$ à plusieurs variables:
$$ x_{k+1}=x_k-\left(\frac{\partial }{\partial x}\left(\frac{\partial E(x)}{\partial x}\right)\right)^{-1}_{x=x_k}.\frac{\partial E(x)}{\partial x}(x_k)=x_k-\left(\frac{\partial }{\partial x}\left(\frac{\partial E(x)}{\partial x}\right)\right)^{-1}_{x=x_k}.\nabla E(x_k) $$
qu'on écrit aussi sous la forme:
\be
\fbox{ $ x_{k+1}=x_k-\left(\nabla^2E(x_k) \right)^{-1}.\nabla E(x_k) $}\lb{eq**23}
\ee
En comparant l'équation (\ref{eq**23}) avec l'équation (\ref{eq**9}), on déduit que le choix de $a_kQ(x_k)$ est le suivant:
\be
a_kQ(x_k)=\left(\nabla^2E(x_k) \right)^{-1} \lb{eq**24}
\ee
La méthode de Newton permet de vérifier la condition:
\be
 \left\|I-a_k(x)Q(x).\ds \frac{\partial ^2E(x)}{\partial x ^2} \right\|=\left\|I-\left(\nabla^2E(x_k) \right)^{-1}.\ds \nabla^2E(x)\right\|=||I-I||< 1 
\ee
ce qui permet d'obtenir la convergence de l'itération.
\subsubsection*{\textit{20.4.4.1. La Vitesse de Convergence}}
La vitesse de convergence de l'itération $k+1$ est obtenu en exprimant $||x_{k+1}-x^*||$ en fonction de $||x_k-x^*||$. On commence par développer l'équation (\ref{eq**23}) au voisinage de la solution $x^*$. Comme:
$$ x_{k+1}=\Phi(x_k)=\Phi(x^*)+\frac{\partial \Phi(x^*)}{\partial x}.(x_k-x^*)+\frac{1}{2}\partial ^2_{xx}\Phi(x^*)((x_k-x^*),(x_k-x^*))+o(||x_k-x^*||^2)$$  
soit:
\be
x_{k+1}=x^*+\frac{\partial \Phi(x^*)}{\partial x}.(x_k-x^*)+\frac{1}{2}\partial ^2_{xx}\Phi(x^*)((x_k-x^*),(x_k-x^*))+o(||x_k-x^*||^2) \lb{eq**25}
\ee
Or d'après (\ref{eq**23}):
\ba
 &\Phi(x)=x-\left(\nabla^2E(x) \right)^{-1}.\nabla E(x) \Rightarrow \ds \frac{\partial \Phi(x)}{\partial x}=I -\partial_x[(\nabla^2E(x))^{-1}\nabla E(x) ] &\nonumber \\
&= I-\partial_x(\nabla^2E(x))^{-1}\nabla E(x)-(\nabla^2E(x))^{-1}.(\nabla^2E(x))= &\nonumber \\
&I-\partial_x(\nabla^2E(x))^{-1}\nabla E(x)-I\Rightarrow & \nonumber \\
& \ds \frac{\partial \Phi(x)}{\partial x}=-\partial_x(\nabla^2E(x))^{-1}\nabla E(x) &
\ea
Par suite pour $x=x^*$, on a :
\be
\partial_x\Phi(x^*)=-\partial_x(\nabla^2E(x^*))^{-1}\nabla E(x^*)=0
\ee
L'équation (\ref{eq**25}) devient:
\be
x_{k+1}=x^*+\frac{1}{2}\partial ^2_{xx}\Phi(x^*)((x_k-x^*),(x_k-x^*))+o(||x_k-x^*||^2) \lb{eq**28}
\ee
On calcule maintenant la dérivée seconde de $\Phi(x)$ soit $\nabla^2\Phi(x)$ ou $\partial^2_{xx}\Phi(x)$:
\ba
\nabla^2\Phi(x)=\partial^2_{xx}\Phi(x)=\partial_x[-\partial_x(\nabla^2E(x))^{-1}\nabla E(x)]\Longrightarrow \nonumber \\
\partial^2_{xx}\Phi(x)=-\partial^2_{xx}(\nabla^2E(x))^{-1}\nabla E(x)-\partial_x(\nabla^2E(x))^{-1}\nabla^2E(x) \nonumber 
\ea
Comme la matrice $\nabla^2E(x)$ est une matrice carrée définie positive, donc elle inversible, on a alors:
$$ (\nabla^2E(x))^{-1}.\nabla^2E(x)=I\Rightarrow \partial_x(\nabla^2E(x))^{-1}.\nabla^2E(x)=-(\nabla^2E(x))^{-1}.\partial_x\nabla^2E(x) $$
Par suite:
\ba
\partial^2_{xx}\Phi(x)=-\partial^2_{xx}(\nabla^2E(x))^{-1}\nabla E(x)-\partial_x(\nabla^2E(x))^{-1}\nabla^2E(x)\nonumber \\
=-\partial^2_{xx}(\nabla^2E(x))^{-1}\nabla E(x)+(\nabla^2E(x))^{-1}.\partial_x\nabla^2E(x) \nonumber 
\ea
Pour $x=x^*$, on a en tenant compte de $\nabla E(x^*)=0$:
\be
\fbox{ $ \partial^2_{xx}\Phi(x^*)=(\nabla^2E(x^*))^{-1}.\partial_x\nabla^2E(x^*) $}
\ee
L'équation (\ref{eq**28}) s'écrit:
$$ x_{k+1}=x^*+\frac{1}{2}(\nabla^2E(x^*))^{-1}.\partial_x\nabla^2E(x^*)((x_k-x^*),(x_k-x^*))+o(||x_k-x^*||^2) $$
ce qui donne:
$ ||x_{k+1}-x^*||\leq||\ds \frac{1}{2}(\nabla^2E(x^*))^{-1}.\partial_x\nabla^2E(x^*)((x_k-x^*),(x_k-x^*))||  $
soit:
\be 
\fbox{$ ||x_{k+1}-x^*||\leq \ds \frac{1}{2}||(\nabla^2E(x^*))^{-1}.\partial_x\nabla^2E(x^*)||.||x_k-x^*||^2=\ds \frac{M}{2}||x_k-x^*||^2$} 
\ee
où $M$ est une constante qui dépend de $x^*$. Donc la vitesse de convergence de l'itération $k+1$ est majorée par une fonction quadratique de la vitesse de convergence de l'itération $k$. D'où le théorème:\index{Vitesse de convergence}
\bthm
Si $\nabla^2E$ est continue et inversible au voisinage d'une solution $x^*$, alors la vitesse de la convergence de l'itération $k+1$ de la méthode de Newton  est quadratique de la vitesse de convergence de l'itération $k$.
\ethm
En d'autres termes, l'avantage de la méthode de Newton est que la convergence est très rapide. Mais comme inconvéniants de la méthode de Newton c'est le calcul et le stockage des éléments de la matrice $\nabla^2E(x)$ laquelle il faut s'assurer qu'elle est définie positive et qu'elle est inversible.

\section{\textsc{La Méthode de Gauss-Newton}}
Dans ce paragraphe, on va présenter la méthode de Gauss-Newton pour minimiser la fonction:\index{La méthode de Gauss-Newton}
\be
E(x)=\frac{1}{2}||y_s-A(x)||^2 \lb{eq**35}
\ee
avec $||.||=(.)^TP(.)$ et $P$ une matrice carré définie positive et $A:\Omega\subset \BbR^n\longrightarrow \BbR^m$. De plus, on suppose que $m > n$. Quand $x$ varie dans $\Omega$, le point $A(x)$ décrit une variété $\m A$ de dimension $n$ plongée dans $\BbR^m$. $||y_s-A(x)||^2$ représente le carré de la distance  du point $y_s$ au point $A(x)$. La recherche de $minE(x)$ revient donc à trouver le point $A(x^*)$ de $\m A$ tel que sa distance au point $y_s$ soit minimale.

On pose :
\be
\fbox{ $ e(x)=y_s-A(x) $}
\ee
Comme:
$$E(x)=\frac{1}{2}(y_s-A(x))^T.P.(y_s-A(x))=\frac{1}{2}(A(x)^T.P.A(x)-2y_s^T.P.A(x)+y_s^T.P.y_s)$$
Ce qui donne:
\be
\fbox{$ \partial_xE(x)=\ds \frac{\partial E(x)}{\partial x}=\nabla E(x)=\frac{\partial A}{\partial x}^T.P.(A(x)-y_s)=-\frac{\partial A(x)}{\partial x}^T.P.e(x)$}
\ee
avec $\ds \frac{\partial A(x)}{\partial x}$ une matrice $m\times n$. Pour appliquer la méthode de Newton, on a besoin de calculer $\nabla^2E(x)$, soit:
$$ \nabla^2E(x)=\frac{\partial}{\partial x} \left(-\frac{\partial A(x)}{\partial x}^T.P.e(x)\right) $$
d'où:
$$ \nabla^2E(x)=-\nabla^2A(x)^T.P.e(x)-\frac{\partial A(x)}{\partial x}^T.P.\frac{\partial e(x)}{\partial x} $$
or $ \ds \frac{\partial e(x)}{\partial x}=-\frac{\partial A(x)}{\partial x}$, on obtient finalement:
\be
\nabla^2E(x)=-\nabla^2A(x)^T.P.e(x)+\frac{\partial A(x)}{\partial x}^T.P.\frac{\partial A(x)}{\partial x}
\ee
Par la formule (\ref{eq**23}), on a la relation entre $x_{k+1}$ et $x_k$:
$$x_{k+1}=x_k-\left(\nabla^2E(x_k) \right)^{-1}.\frac{\partial E(x_k)}{\partial x}$$
soit:
\be
\fbox{$ x_{k+1}=x_k+\left(\partial_x A(x_k)^T.P.\partial_x A(x_k)-\nabla^2A(x_k)^T.P.e(x_k)\right)^{-1}.(\partial_xA(x_k)^T.P.e(x_k)) $}
\ee
Seulement, la formule d'itération ci-dessus de la méthode de Newton ne prend pas en considération de l'avantage de l'expression de la fonction $E(x)$ donnée par (\ref{eq**35}). Au lieu de prendre $ a_k(x_k)Q(x_k)=(\nabla^2E(x_k))^{-1}$, on prendra :\index{La méthode de Newton}
\be
\fbox{ $ Q(x_k)=\left( \ds \frac{\partial A(x_k)}{\partial x}^T.P.\frac{\partial A(x_k)}{\partial x}\right)^{-1} $ }
\ee
et on fait abstraction du terme $-\nabla^2A(x_k)^T.P.e(x_k)$ dont le calcul est plus compliqué.
  
La matrice $Q(x_k)$ est une matrice $n\times n$ définie positive. En effet, soit un vecteur $X\neq 0 \in \Omega \subset \BbR^n$, comme $P$ est une matrice $m\times m$ définie positive, d'où:
$$ X^T.Q^{-1}.X=X^T.(\partial_xA^T.P.\partial_xA).X=(\partial_xA.X)^T.P.(\partial_xA.X) >0 $$
car $P$ est définie positive et que $\partial_xA(x).X\neq0$. Donc $Q(x_k)^{-1}$ est définie positive par suite $Q(x_k)$ est définie positive.
\\

On obtient pour la méthode dite de Gauss-Newton la formule:
\be
\fbox{$ x_{k+1}=x_k+\left(\partial_x A(x_k)^T.P.\partial_x A(x_k)\right)^{-1}.(\partial_x A(x_k)^T.P.e(x_k)) $}
\ee
La solution $x^*$ vérifie:
\be
\fbox{ $ \left\{\begin{array}{l}
\ds \frac{\partial E(x^*)}{\partial x}=0 \\
\\
\nabla^2E(x^*)\,\, \mbox{est définie positive}
\end{array}\right. $} \lb{eq**76}
\ee
L'interprétation de la première équation est que le vecteur résidu $v=-e(x^*)=A(x^*)-y_s$ est orthogonal au vecteur $\ds \frac{\partial A(x^*)}{\partial x}$ soit $e(x^*)$ est perpendiculaire au plan tangent au point $A(x^*)$ de la variété $\m A$, donc $||y_s-A(x^*)||$ représente la distance minimale du point $y_s$ à la variété $\m A$.
\\

La deuxième équation de (\ref{eq**76}) a une autre interprétation géométrique moins évidente qu'on découvrira dans le chapitre suivant. 

\section{\textsc{Exercices et Problèmes}}
\bex
 On considère $(u,v)\in \BbR^2$ et on définit la fonction par :
$$ f(u,v)=u^4+6uv+1.5v^2+36v+405 $$
1. Chercher les points critiques réels de $f$.

2. Montrer que le point $x^*=(u,v)=(3,-18)$ est un point minimum de $f$.

3. Montrer que le Hessien de $f$ est une matrice définie positive si $u^2>1$ et indéfinie si $u^2< 1$.

4. Montrer que la formule de récurrence de Newton s'écrit avec $ J=1.5(u_k^2-1)$:
	\[u_{k+1}=\frac{u_{k}^3+9}{J},\quad v_{k+1}=-\frac{2u_{k}^3+18u_k^2}{J}
	\]
	\eex
\bpb
Soient le plan $(P)$ et la sphère $(\BbS^2)$ d'équations respectivement: $x+y+z=1$ et $x^2+y^2+z^2=1$. On veut chercher le point $M\in (\BbS^2)$  tel que sa distance au plan $(P)$ soit maximale.

1. Montrer que la distance d'un point $M(X,Y,Z)$ au plan $(P)$ est donnée par : 
	\[d=|X+Y+Z-1|/\sqrt{3} 	\]
2. Pour répondre à la question posée ci-dessus, on considère la fonction: $E(x,y,z,\lambda)=-(x+y+z-1)^2$$-\lambda(x^2+y^2+z^2-1)$. Ecrire le système d'équations donnant les points critiques de $E$ qu'on note par (1).

3. Montrer que si $\lambda=-1$, on arrive à une contradiction. On suppose que $\lambda \neq-1$. Que représente le cas $\lambda=0$.

4. On suppose que $\lambda \notin \{-3,-1,0\}$. Résoudre le système (1). Soit le point $M_2$ tel que ses coordonnées sont négatives.

5. Montrer que la matrice hessienne de $E$ pour $M_2$ s'écrit sous la forme:
	\[ H=\begin{pmatrix}{
\mu^2 & -2 & -2 \cr
-2& \mu^2 & -2 \cr
-2 & -2 & \mu^2 }
\end{pmatrix}\quad avec\,\,\,\mu=1+\sqrt{3}\]
6. Si on pose $U=(X,Y,Z)^T \in (\BbS^2)$. Montrer que $U^T.H.U=2\left[3+\sqrt{3}-(X+Y+Z)^2\right]$. En déduire que $U^T.H.U>0$ pour tout $U\neq 0 \in (\BbS^2)$. 

7. Montrer que pour le point $M_2$, on obtient un minimum strict de $E$. A-t-on répondu à la question du problème. 
\epb

\chapter{\textit{\textbf{Interpr\'etation G\'eom\'etrique de la Compensation Non-Lin\'eaire}}}
\section{\textsc{Introduction}}

 E. Grafarend et B. Schaffrin ont étudié dans un article (\textit{E.W. Grafarend \& B. Schaffrin}, 1989), la g\'eom\'etrie de la compensation ou l'ajustement non-lin\'eaire et ont pr\'esent\'e le cas du probl\`eme d'intersection plane en utilisant le mod\`ele de Gauss-Markov,\index{\textbf{Gauss C.F.}}\index{\textbf{Markov A.}} par les moindres carr\'es. Dans ce chapitre, on pr\'esente les principes de la g\'eom\'etrie de la  compensation non-lin\'eaire par la m\'ethode des moindres carr\'es en s'appuyant sur le lemme de P\'azman (\textit{A. P\'azman}, 1984).\index{\textbf{Grafarend E.W.}}\index{\textbf{Schaffrin B.}}\index{\textbf{P\'azman A.}}

\vspace{4mm}
\section{\textsc{La G\'eom\'etrie Non Lin\'eaire du Mod\`ele de Gauss-Markov}}
Le mod\`ele non lin\'eaire de Gauss-Markov est d\'efini par:\index{Modèle de Gauss-Markov}
\be
	\fbox{ $ \zeta(X)=L-e, \quad e \in \mathcal{N}(0,\Gamma) $} \label{m1}
\ee
avec:

- $L$: le vecteur des observations $(n\times1)=(L_1,L_2,...,L_n)^T$;

- $X$: le vecteur des inconnues $(m\times1)=(X_1,X_2,...,X_m)^T$;

- $e$: le vecteur des erreurs $(n\times1)=(e_1,e_2,...,e_n)^T$ suit la loi normale $\mathcal{N}(0,\Gamma)$ avec $E(e)=0$ et $\Gamma=E(ee^T)$ la matrice de dispersion ou variance, on prendra $\Gamma=\sigma^2_0.P^{-1}$. $P$ est la matrice des poids et $\sigma_0$ une constante positive;

- $\zeta$: est une fonction donn\'ee injective d'un ouvert $U\subset \BbR^m \rightarrow \BbR^n$ et  $m<n$.
\\

Remarque: dans le cas d'un mod\`ele lin\'eaire, la fonction $\zeta=A.X$ o\`u $A$ est une matrice $n\times m$.
\\

On note $Im \zeta=\left\{\zeta(X) / X \in U\right\}$ l'image de $U$ par la fonction $\zeta$. $Im\zeta$ est une vari\'et\'e de dimension $m$ v\'erifiant les conditions:
\vspace{4mm}

(i): les vecteurs $\displaystyle \frac{\partial \zeta}{\partial X_1},\frac{\partial \zeta}{\partial X_2},...,\frac{\partial \zeta}{\partial X_m}$ sont lin\'eairement ind\'ependants en chaque point $X \in U$;
\vspace{4mm}

(ii): les fonctions $\displaystyle \frac{\partial^2 \zeta}{\partial X_i \partial X_j}$ sont continues sur $U$ pour $i,j \in \left\{1,2,...,m\right\}$.
\\

On introduit un produit scalaire:
\be
	<\zeta_1,\zeta_2>=\zeta_1^T.P.\zeta_2 \label{m2}
\ee
 D'o\`u la norme du vecteur $\zeta=(\zeta_1,\zeta_2,...,\zeta_n)^T$:
 \be
	\|\zeta \|^2=<\zeta,\zeta>=\zeta^T.P.\zeta =\sum_{i=1}^np_i.\zeta_i^2 \label{m3}
\ee
dans l'espace vectoriel $\BbR^n$ en prenant la matrice de poids $P$ une matrice diagonale. 
\\

Alors la solution par les moindres carr\'es $\overline{X}$ sera d\'efinie par:
\be
		\fbox{ $ \|L-\bar{\zeta}(\overline{X}) \|=min\left\{\|L-\zeta(X) \| \,/ \,X\in U \right\}  $}\label{m4}
\ee
Cette condition est exprim\'ee par les \'equations suivantes:
\be
\fbox{ $ \ds 	\frac{\partial}{\partial X_i} \|L-\zeta(X) \|^2=0 \quad \mbox{pour}\,\,i\in \left\{1,2,...,m\right\} $}\label{m5}
\ee
En effet, on veut minimiser la fonction:
\be
	F(X)=F(X_1,X_2,...,X_m)=\|L-\zeta(X) \|=\|L-\zeta(X_1,X_2,...,X_m) \| \label{m6}
\ee
Comme $F$ est une fonction positive, minimiser $F$ c'est aussi minimiser $F^2$, soit $J(X)=F^2(X)$ d'où:
$$	\frac{\partial J(X)}{\partial X_i}=0 $$ 
soit:
\be 
	\frac{\partial}{\partial X_i} \|L-\zeta(X_1,X_2,...,X_m) \|^2=0 \quad \mbox{pour}\,\,i \in \left\{1,2,...,m\right\} \label{m7}
\ee
or:
\ba
	& \|L-\zeta(X_1,X_2,...,X_m) \|^2=(L-\zeta(X_1,X_2,...,X_m))^T.P.(L-\zeta(X_1,X_2,...,X_m))=\nonumber & \\ &
	 \zeta(X)^T.P.\zeta(X) -2L^T.P.\zeta(X)+L^T.P.L \label{m8} &
	\ea
Par suite:
\be 
		\frac{\partial  J(X)}{\partial X_i}=2\zeta(X)^T.P.	\frac{\partial \zeta(X)}{\partial X_i}-2L^T.P.\frac{\partial \zeta(X)}{\partial X_i} \quad \mbox{pour}\,\,i\in \left\{1,2,...,m\right\} \label{m9}
\ee
ou encore :
$$ 		\frac{\partial  J(X)}{\partial X_i}=2(\zeta(X)-L)^T.P.\frac{\partial \zeta(X)}{\partial X_i} \quad \mbox{pour}\,\,i\in \left\{1,2,...,m\right\} $$

ce qui donne en utilisant (\ref{m5}):
\ba
<L-\zeta(X),\frac{\partial \zeta(X)}{\partial X_i}>\,=0 \quad \mbox{pour}\,\,i\in \left\{1,2,...,m\right\} \nonumber  \\
\fbox{ $ \mbox{ou}\quad <e,\ds \frac{\partial \zeta(X)}{\partial X_i}>\,=0 \quad \mbox{pour}\,\,i\in \left\{1,2,...,m\right\} $}\label{m11}
\ea
G\'eom\'etriquement, cela veut dire que le vecteur erreur $e =L-\zeta(X)$ est perpendiculaire (produit scalaire nul) au plan tangent de la vari\'et\'e $Im \zeta$ au point $\bar{\zeta}(\overline{X})$ (s'il existe).
\\

Pour le cas non-lin\'eaire, la condition (\ref{m11}) est n\'ecessaire mais non suffisante. Pour obtenir le minimum, il faut que la matrice $\left(\displaystyle \frac{\partial^2 J}{\partial X_i \partial X_j} \right)\, i,j\in \left\{1,2,...,m\right\} $ soit d\'efinie positive.

\section{\textsc{Interpr\'etation G\'eom\'etrique}}
Dans cette section, on va clarifier l'interpr\'etation g\'eom\'etrique de fa\c{c}on que la solution de (\ref{m11}) soit localement unique.
\\

On considère la matrice $m\times m$ d\'efinie par:
\be
	G(X)=(G_{\alpha\,\beta}) \quad \mbox{avec}\,\,G_{\alpha\,\beta}=\,<\frac{\partial \zeta(X)}{\partial X_{\alpha}},\frac{\partial \zeta(X)}{\partial X_{\beta}}> \quad \left\{
\begin{array}{ll}
	\alpha=1,2,...,m \\
	\beta=1,2,...,m
\end{array} \right. \label{m12}
\ee
Or:
\be 
	 ds^2=G_{\alpha\,\beta}dX_{\alpha}dX_{\beta} \label{m13}
	 \ee
repr\'esente la m\'etrique de la vari\'et\'e $Im \zeta$. La matrice $G(X)= (G_{\alpha\,\beta})$ est appel\'ee en terme statistique la matrice d'information de Fisher. \index{Matrice d'information de Fisher}\index{\textbf{Fisher R.A.}}
\\

On introduit la matrice $B$ d\'efinie par:
\be
	B(X,L)=(B_{\alpha\,\beta}) \quad \mbox{avec}\,\,\,B_{\alpha\,\beta}=\frac{1}{2}\frac{\partial^2 }{\partial X_{\alpha}\partial X_{\beta}}\|L-\zeta(X)\|^2 \quad \left\{
\begin{array}{ll}
	\alpha=1,2,...,m \\
	\beta=1,2,...,m
\end{array} \right. \label{m14}
\ee
 or: 
\ba
&	 \|L-\zeta(X)\|^2=(L-\zeta(X))^TP(L-\zeta(X))=(L^T-\zeta^T(X))(PL-P\zeta(X))=\nonumber& \\ & L^TPL-2L^TP\zeta(X)+\zeta^T(X)P\zeta(X)\label{m15} &
\ea
 D'o\`u:
 \ba
&	\ds \frac{\partial}{\partial X_{\alpha}}(L^TPL-2L^TP\zeta(X)+\zeta^T(X)P\zeta(X))=\ds -2L^TP\frac{\partial \zeta}{\partial X_{\alpha}}+2\zeta^T(X)P\ds \frac{\partial \zeta(X)}{\partial X_{\alpha}}= \nonumber& \\ & -2(L^T-\zeta^T(X))P\ds \frac{\partial \zeta(X)}{\partial X_{\alpha}} \label{m16}&
\ea
Donc:
\ba
&	\ds \frac{1}{2}\frac{\partial^2 }{\partial X_{\alpha}\partial X_{\beta}}\|L-\zeta(X)\|^2=-L^TP\frac{\partial ^2 \zeta(X)}{\partial X_{\alpha}\partial X_{\beta} }+\frac{\partial \zeta(X)}{\partial X_{\beta}}P\frac{\partial \zeta(X)}{\partial X_{\alpha}}+\zeta^T(X)P\frac{\partial^2 \zeta(X)}{\partial X_{\alpha} \partial X_{\beta}}= \nonumber & \\  & \ds \frac{\partial \zeta^T(X)}{\partial X_{\beta}}P\frac{\partial \zeta(X)}{\partial X_{\alpha}}-(L^T-\zeta^T(X))P\frac{\partial^2 \zeta(X)}{\partial X_{\alpha} \partial X_{\beta}}= \nonumber & \\ & \ds <\frac{\partial \zeta(X)}{\partial X_{\beta}},\frac{\partial \zeta(X)}{\partial X_{\alpha}}> -<L-\zeta(X),\frac{\partial ^2 \zeta(X)}{\partial X_{\beta}\partial X_{\alpha}}> \label{m17} &
\ea
Soit:
\be
	B_{\alpha\,\beta}=G_{\alpha\,\beta}-<L-\zeta(X),\frac{\partial ^2 \zeta(X)}{\partial X_{\beta}\partial X_{\alpha}}> \label{m18}
\ee
On pose:
\be
	H(L,X)=(h_{\alpha \beta})=\left(<L-\zeta(X),\frac{\partial ^2 \zeta(X)}{\partial X_{\beta}\partial X_{\alpha}}>\right)\quad \left\{\begin{array}{ll}
	\alpha=1,2,...,m \\
	\beta=1,2,...,m
\end{array} \right. \label{m19} 
\ee
c'est-\`a-dire:
\be
	B=G-H \label{m20}
\ee
On revient \`a $Im\zeta(X)=\left\{ \zeta(X)\quad / X\in U\right\}$. Soit une ligne g\'eod\'esique de $Im \zeta(X)$ passant par un point $\zeta=\zeta(X)$ param\'etr\'ee par son abscisse curviligne $s$, on a alors:
\be
  \fbox{ $ \chi(s)=\zeta(X(s)),\,\,\,s\in [s_1,s_2] $} \label{m21}
\ee
o\`u $X(s)$ d\'ecrit une certaine courbe dans le domaine $U\subset \BbR^m$.
\\

Le vecteur :
\be
	\chi'(s)=\frac{ d\chi(s)}{ds} \label{m22}
\ee
repr\'esente le vecteur tangent \`a la ligne g\'eod\'esique au point $\zeta(X(s))$ de $Im \zeta(X)$. Ce vecteur v\'erifie:
\be
	\|\chi'(s)\|^2=1  \label{m23}
\ee
Par suite, la d\'eriv\'ee de ce vecteur par rapport \`a $s$ est un vecteur orthogonal \`a $\chi'(s)$ donc orthogonal \`a $Im \zeta(X(s))$ au point $\zeta(X(s))$:
\be
	\fbox{ $ 	\chi"(s)=\ds \frac{ d\chi'(s)}{ds}\, \bot \,\,\chi'(s)  $} \label{m24}
\ee
 c'est-\`a-dire parall\`ele au vecteur normal \`a la surface ou la vari\'et\'e $Im \zeta(X)$ et on retrouve la propri\'et\'e que $\chi(s)$ est une g\'eod\'esique.
 
 On remarque que pour une ligne g\'eod\'esique, la courbure g\'eod\'esique est nulle et la courbure normale coïncide avec la courbure de la courbe $\chi(s)$ soit:
\be
	\fbox{ $ \rho(s)=\ds \frac{1}{\|\chi"(s)\|} $} \label{m25}
\ee          
le rayon de courbure. On appelle:
\be
	n(s)=\frac{\chi"(s)}{\|\chi"(s)\|}=\chi"(s).\rho(s) \label{m26}
\ee             
C'est un vecteur unitaire perpendiculaire au plan tangent \`a la surface $Im \zeta(X)$.

D'apr\`es l'\'equation (\ref{m11}), au point $\zeta(\overline{X})$, le vecteur $e$ est perpendiculaire \`a $L-\zeta(\overline{X})$. On note aussi : 
  \be
	K(\zeta(\overline{X}))=\left\{ Z /\,Z\in \BbR^n, \mbox{ avec} <Z,\frac{\partial \zeta(\overline{X})}{\partial X_{\alpha}}>\, =0\,\,,\alpha=1,2,...,m \right\} \label{m27} 
\ee              
On a donc $e=L-\zeta(\overline{X}) \in K$. Ce dernier est un espace vectoriel de dimension $(n-m)$ orthogonal \`a $Im\zeta(X)$ au point $\zeta(\overline{X})$. On a aussi $\chi"(s) \in K$.
\subsection*{21.3.1. Lemme de P\'azman}
 On peut maintenant \'enoncer le lemme de P\'azman (\textit{A. P\'azman}, 1984) comme suit:\index{Lemme de P\'azman} \index{\textbf{P\'azman A.}}
 
 \blm
(\textbf{de P\'azman:}) Pour tout vecteur d'observation $L \in \BbR^m $, et toute solution appropri\'ee $\overline{X}$ des \'equations:
  $$ \fbox{ $ <L-\zeta(\overline{X}),\ds \frac{\partial \zeta(\overline{X})}{\partial X_{\alpha}}>\, =0,\,\,\, \alpha=1,2,...,m  $} $$
 les conditions suivantes sont \'equivalentes:
  
 1 - La matrice:
$$ \fbox{ $ B(\overline{X},L)=\ds G(\overline{X})-(<L-\zeta(\overline{X}),\frac{\partial ^2\zeta(\overline{X})}{\partial X_{\alpha} \partial X_{\beta}}>) $} $$ est d\'efinie positive.
 
 2 - Pour toute ligne g\'eod\'esique $\chi(X(s))$ v\'erifiant :
 $$ \chi(\bar{s})=\zeta(\overline{X}(\bar{s})) $$
on a l'in\'egalit\'e:
\be
	\fbox{ $  <L-\zeta(\overline{X}),n(\bar{s})>\,\, < \, \rho(\bar{s}) $} \label{m28}
\ee 
\elm         
\subsubsection*{$\un{\textbf{1} \Longrightarrow \textbf{2}}$}
 En effet, on suppose que la matrice $B(\overline{X},L)$ est d\'efinie positive c'est-\`a-dire:
 \be
	\forall\, Y \in   \BbR^m,\,\, Y\neq 0 \Longrightarrow Y^T.B.Y > 0 \label{m29}
	\ee             
	On prend alors: $Y=\chi'(\bar{s}) $. On a:
	\be
	\chi'^T(\bar{s}).B(L,\overline{X}).\chi'(\bar{s}) > \,0  \label{m30}
\ee             
Comme $B=G(\overline{X}) - H(L,\overline{X})$, on obtient:
$$\chi'^T(\bar{s}).(G(\overline{X})- H(L,\overline{X})).\chi'(\bar{s}) > \,0 $$
soit:
\be
	\chi'^T(\bar{s}).G(\overline{X}).\chi'(\bar{s})- \chi'^T(\bar{s}).H(L,\overline{X}).\chi'(\bar{s}) > \,0 \label{m31}
\ee          
Or pour $s$:
\be
 Y=\chi'(s)=\ds \frac{d\chi(s)}{ds}=\sum_{i=1}^{i=m}\frac{\partial \zeta(X(s))}{\partial X_i}\frac{d X_i(s)}{ds}=\sum_{i=1}^{i=m}\chi'_i(s)\frac{\partial \zeta(X(s))}{\partial X_i} \label{m32}
\ee          
en notant $\chi'_i(s)=\ds \frac{d X_i(s)}{ds}$ les composantes de $\chi'(s)$ dans le plan tangent \`a $Im \zeta$ au point $\zeta (X(s))$. Comme:
\ba
& \| \chi'(s)\|^2=1=\chi'^T(s).\chi'(s)= & \nonumber \\ & \ds <\sum_{i=1}^{i=m}\chi'_i(s)\frac{\partial \zeta(X(s))}{\partial X_i},\sum_{j=1}^{j=m}\chi'_j(s)\frac{\partial \zeta(X(s))}{\partial X_j}>=\nonumber  & \\ &\ds \sum_{i=1}^{m}\sum_{j=1}^m\chi'_i(s)\left( <\ds \frac{\partial \zeta(X(s))}{\partial X_i},\frac{\partial \zeta(X(s))}{\partial X_j}> \right)\chi'_j(s) =\nonumber  & \\ &\chi'(s)^T.G.\chi'(s)=1 \label{m33} &
\ea
On prend $s=\bar{s}$, alors (\ref{m31}) devient:
\be
	 \chi'^T(\bar{s}).H(L,\overline{X}).\chi'(\bar{s}) < \,1 \label{m34}
\ee             
Comme $B=G-H$  donc la matrice $H$ est exprim\'ee dans la base de $B$ soit $\left(\ds \frac{\partial \zeta(X(s))}{\partial X_j}\right)$.
En utilisant (\ref{m32}), le nombre r\'eel $\chi'^T(\bar{s}).H(L,\overline{X}).\chi'(\bar{s})$ s'\'ecrit:
\be
\chi'^T(\bar{s}).H(L,\overline{X}).\chi'(\bar{s})=\sum_{i=1}^{i=m}\chi'_i(\bar{s})\left(\sum_{j=1}^{j=m}h_{ij}.\chi'_j(\bar{s})\right)=\sum_{i,j=1}^{m}\chi'_i(\bar{s}).\chi'_j(\bar{s}).h_{ij} \label{m35}
\ee           
On remplace $h_{ij}$ par:
 $$<L-\zeta(\overline{X}),\frac{\partial ^2\zeta(\overline{X})}{\partial X_i \partial X_j}> $$ Par un calcul simple, l'\'equation (\ref{m35}) devient:
\be
\chi'^T(\bar{s}).H(L,\overline{X}).\chi'(\bar{s})= <L-\zeta(\overline{X}),\sum_{i,j=1}^{m}\chi'_i(\bar{s}).\chi'_j(\bar{s})\frac{\partial ^2\zeta(\overline{X})}{\partial X_i \partial X_j}> \label{m36}
\ee           
Maintenant, on va s'int\'eresser au membre à droite du produit scalaire de l'\'equation (\ref{m36}). En diff\'erentiant l'\'equation (\ref{m32}) par rapport \`a $s$, on obtient:
\be
	\chi"(s)=\sum_i \frac{d\chi'_i(s)}{ds}.\frac{\partial \zeta(X(s))}{\partial X_i}+\sum_i\chi'_i(s)\sum_j\frac{\partial^2 \zeta(X(s))}{\partial X_i \partial X_j}.\frac{dX_j}{ds} \label{m37}
	\ee             

	Alors on a pour $s=\bar{s}$:
	\ba
& <L-\zeta(\overline{X}),\chi"(\bar{s})>=<L-\zeta(\overline{X}),\ds \sum_i \frac{d\chi'_i(s)}{ds}.\frac{\partial \zeta(X(s))}{\partial X_i}+\nonumber &\\
& \sum_i \chi'_i(s)\sum_j \ds \frac{\partial^2 \zeta(X(s))}{\partial X_i \partial X_j}.\frac{dX_j}{ds}>= \ds \sum_i \frac{d\chi'_i(s)}{ds}<L-\zeta(\overline{X}),\frac{\partial \zeta(X(s))}{\partial X_i}>+\nonumber & \\ 
& \sum_i \sum_j \chi'_i(s) <L-\zeta(\overline{X}),\ds \frac{\partial^2 \zeta(X(s))}{\partial X_i \partial X_j}.\frac{dX_j}{ds}> \label{m38} &
\ea             
Or en utilisant l'\'equation (\ref{m11}), le premier terme de la deuxi\`eme ligne de l'\'equation pr\'ec\'edente est nul:
$$ \sum_i \frac{d\chi'_i(s)}{ds}.<L-\zeta(\overline{X}),\frac{\partial \zeta(X(\bar{s}))}{\partial X_i}>\, =0 $$
 et comme :
$$ \chi'_j(s)=\ds \frac{dX_j}{ds}$$
Alors l'\'equation (\ref{m38}) devient:
	\ba
&<L-\zeta(\overline{X}),\chi"(\bar{s})>=\ds \sum_i \chi'_i(s)\sum_j <L-\zeta(\overline{X}),\frac{\partial^2 \zeta(X(s))}{\partial X_i \partial X_j}.\chi'_j(s)> = \nonumber &\\ &\ds \sum_i \sum_j \chi'_j(\bar{s})\chi'_i(\bar{s}) <L-\zeta(\overline{X}),\frac{\partial^2 \zeta(X(\bar{s}))}{\partial X_i \partial X_j}> \label{m39}&
\ea         
Or le deuxi\`eme membre n'est autre que l'\'equation (\ref{m36}). En utilisant (\ref{m34}), on obtient:
\be
	<L-\zeta(\overline{X}),\chi"(\bar{s})> \,\,\,< 1 \label{m40}
\ee           
Or:
$$ n(\bar{s})=\rho(\bar{s}).\chi"(\bar{s})$$
D'o\`u:
\be
	<L-\zeta(\overline{X}),n(\bar{s})> \quad < \, \rho(\bar{s}) \label{m41}
\ee              
\subsubsection*{$\un{\textbf{2} \Longrightarrow \textbf{1}}$}
On suppose que reciproquement, on a pour toute g\'eod\'esique $\chi(s)=\zeta(X(s))$ de $Im\zeta$ passant par le point $ \chi(\bar{s})=\zeta(\overline{X}(\bar{s}))$ v\'erifiant:
\be
	<L-\zeta(\overline{X}),n(\bar{s})>\quad< \, \rho(\bar{s}) \label{m42}
\ee           
et :
$$ <L-\zeta(\overline{X}),\frac{\partial \zeta(X)}{\partial X_i}>\,=0 \quad \mbox{pour}\,\,i\in \left\{1,2,...,m\right\} $$
tel que $\chi'(s)$ v\'erifiant :$$ \|\chi'(s) \|^2=1 $$
Comme: $$n(s)=\frac{\chi"(s)}{\|\chi"(s)\|}=\rho(s)\chi"(s)$$ 
le rempla\c{c}ant dans l'\'equation (\ref{m42}), on obtient:
\be
	<L-\zeta(\overline{X}),\rho(s)\chi"(s)>\quad <\, \rho(\bar{s}) \label{m43}
\ee              
et en simplifiant par $\rho\neq0$, soit:
\be
	<L-\zeta(\overline{X}),\chi"(s)>\quad < \,1 \label{m44}
\ee              
Comme:
$$ \chi(s)=\zeta(X(s))\Rightarrow \chi'(s)=\frac{d\chi(s)}{ds}=\sum_{i=1}^m\frac{\partial \zeta}{\partial X_i}\frac{dX_i(s)}{ds}=\sum^m_{i=1}\chi'_i(s)\frac{\partial \zeta}{\partial X_i} $$
D'o\`u en d\'erivant une deuxi\`eme fois par rapport \`a $s$:
\ba
		\chi"(s)=\sum_i \frac{d\chi'_i(s)}{ds}.\frac{\partial \zeta(X(s))}{\partial X_i}+\sum_i\chi'_i(s)\sum_j\frac{\partial^2 \zeta(X(s))}{\partial X_i \partial X_j}.\frac{dX_j}{ds} \nonumber \\
	 = \sum_i \frac{d\chi'_i(s)}{ds}.\frac{\partial \zeta(X(s))}{\partial X_i}+\sum_i\sum_j \chi'_i(s)\chi'_j(s)\frac{\partial^2 \zeta(X(s))}{\partial X_i \partial X_j} \label{m45}
\ea            
En rempla\c{c}ant $\chi"(s)$ dans (\ref{m44}), on obtient:
\be
	<L-\zeta(\overline{X}),\sum_i \frac{d\chi'_i(s)}{ds}.\frac{\partial \zeta(X(s))}{\partial X_i}+\sum_i\sum_j \chi'_i(s)\chi'_j(s)\frac{\partial^2 \zeta(X(s))}{\partial X_i \partial X_j}>\quad < \,1 \label{m46}
\ee   
ou encore:
\be
		\sum_i\chi_i"(s)<L-\zeta(\overline{X}),\frac{\partial \zeta(X(s))}{\partial X_i}>+\sum_{i,j}\chi'_i(s)\chi'_j(s)<L-\zeta(\overline{X}),\frac{\partial^2 \zeta(X(s))}{\partial X_i \partial X_j}>\quad < \,1 \label{m47}
\ee           
Or la premi\`ere somme est nulle en vertu de l'\'equation (\ref{m11}). Il reste:
\be
	\sum_i\sum_j \chi'_i(s)\chi'_j(s)<L-\zeta(\overline{X}),\frac{\partial^2 \zeta(X(s))}{\partial X_i \partial X_j}>\quad < \,1 \label{m48}
\ee          
Comme:
$$ \chi'(s)=\frac{d\chi(s)}{ds}=\sum^m_{i=1}\frac{\partial \zeta}{\partial X_i}\frac{dX_i(s)}{ds}=\sum^m_{i=1}\chi'_i(s)\frac{\partial \zeta}{\partial X_i} $$
et $\chi'(s)$ v\'erifie $\|\chi'(s)\|=1$ car c'est un vecteur unitaire tangent \`a la g\'eod\'esique $\chi(s)$. Donc:
\be
		<\chi'(s),\chi'(s)>\,=1 \Rightarrow	<\sum^m_{i=1}\chi'_i(s)\frac{\partial \zeta}{\partial X_i},\sum^m_{j=1}\chi'_j(s)\frac{\partial \zeta}{\partial X_j}>\,=1 \label{m49}
\ee          
soit:
\be
	\sum^m_{i=1}\sum^m_{j=1}\chi'_i(s)	<\frac{\partial \zeta}{\partial X_i},\frac{\partial \zeta}{\partial X_j}>\chi'_j(s)=	\sum^m_{i=1}\sum^m_{j=1}\chi'_i(s)G_{ij}\chi'_j(s)=1 \label{m50}
\ee             
et ce-ci n'est autre que :
\be
	\chi'^T(s).G.\chi'(s)=1 \label{m51}
\ee          
On utilise l'\'equation (\ref{m48}), on a:
\be
	\sum_i\sum_j \chi'_i(s)\chi'_j(s)<L-\zeta(\overline{X}),\frac{\partial^2 \zeta(X(s))}{\partial X_i \partial X_j}>\quad < \,\chi'^T(s).G.\chi'(s) \label{m52}
\ee                
ou encore:
\be
0< \chi'^T(s).G.\chi'(s)-	\sum_i\sum_j \chi'_i(s).h_{ij}.\chi'_j(s) \label{m53}
\ee           
Finalement, on obtient:
\be
0< \chi'^T(s).G.\chi'(s)-	\chi'^T(s).H.\chi'(s) \label{m54}
\ee             
C'est-\`a-dire pour tout vecteur $Y=\chi'(s)\neq 0$ du plan tangent de $Im \zeta$:
\be
	Y^T.B.Y > \,0\, \Longrightarrow \quad \mbox{ la matrice $B$ est d\'efinie positive} \label{m55}
\ee          
\begin{flushright}
C.Q.F.D
\end{flushright}
Maintenant, on peut dire quand $\overline{X}$ solution de (\ref{m11}) est solution des moindres carr\'es en r\'ef\'erence \`a l'\'equation (\ref{m4}) et ce \`a partir du corollaire suivant:
\\

\textbf{Corollaire 21.1}\textit{ Si $\overline{X}$ est solution de:}
\be
\fbox{ $ <L-\zeta(X),\ds \frac{\partial \zeta(X)}{\partial X_i}>\,=0 \quad \mbox{pour}\,i\in \left\{1,2,...,m\right\} $} \label{m56} 
\ee         
\textit{ avec :}
\be
	\|L-\zeta(\overline{X}\| < r \label{m57}
\ee           
\textit{o\`u $r$ d\'esigne le rayon de courbure minimum de la  vari\'et\'e $Im \zeta$,  d\'efini par:}
\be
	\fbox{ $ r=min\left\{ \rho_{\chi}(s)\,\,\mbox{avec}\,\,\chi(s)=\zeta(X(s)) \,\,\mbox{une g\'eod\'esique passant par}\,\zeta(X),\forall X\in U \right\}$}  \label{m58}
\ee           
\textit{Alors $\overline{X}$ coïncide avec la solution des moindres carr\'es $\overline{X}=\hat{X}(L)$.}
\\

\textbf{D\'emonstration:}

Comme les deux vecteurs $e=L-\zeta(\overline{X})$ et $n(\bar{s})$ sont orthogonaux \`a $Im \zeta$, ils sont colin\'eaires et comme $n(\bar{s})$ est un vecteur unitaire alors leur produit scalaire $<L-\zeta(\overline{X}),n(\bar{s})>$  est plus petit ou \'egal \`a $\|L-\zeta(\overline{X}\|$. Or ce terme est plus petit que $r$ d'apr\`es (\ref{m57}) et comme $r$ est le plus petit rayon de courbure, on a $r \leq \rho_{\chi}(\bar{s})$  pour toute g\'eod\'esique $\chi(\bar{s})$ passant par $\zeta(\overline{X})$. Ce-ci est traduit par l'\'equation:
\be
	<L-\zeta(\overline{X}),n(\bar{s})>\quad \leq\|L-\zeta(\overline{X}\| < r \leq \rho_{\chi}(\bar{s}) \label{m59}
\ee         
 De cette derni\`ere \'equation, on a:
 \be
	<L-\zeta(\overline{X}),n(\bar{s})> \quad \leq \,\rho_{\chi}(\bar{s}) \label{m50}
\ee           
On utilise le lemme de P\'azman cit\'e ci-dessus, la matrice $B$ est d\'efinie positive donc $\zeta(\overline{X})$ est un minimum strict (\textit{H. Cartan}, 1979). Or on a suppos\'e que l'application $\zeta$ est injective (si $\zeta(X_1)=\zeta(X_2)\Longrightarrow X_1=X_2$), alors $\overline{X}$ coïncide avec la solution des moindres carr\'es $\overline{X}=\hat{X}(L)$.
\section{\textsc{Exercices et Problèmes}}
\bpb
 Dans le plan affine $\mathcal P$, on a mesuré trois distances planes entre un point inconnu $P(X_1,X_2)$ vers trois points connus $P_i(a_i,b_i)_{i=1,3}$ dans trois directions différentes. On considère le mod\`ele non lin\'eaire de Gauss-Markov d\'efini par:
$$	\zeta(X)=L-e, \quad e \in \mathcal{N}(0,\Gamma) $$
avec:

- $L$: le vecteur des observations $(3\times1)=(L_1,L_2,L_3)^T$;

- $X$: le vecteur des inconnues $(2\times1)=(X_1,X_2)^T$;

- $e$: le vecteur des erreurs $(3\times1)=(e_1,e_2,e_3)^T$ suit la loi normale $\mathcal{N}(0,\Gamma)$ avec $E(e)=0$ et $\Gamma=E(ee^T)$ la matrice de dispersion ou variance, on prendra $\Gamma=\sigma^2_0.P^{-1}$, $P$ est la matrice des poids égale à la matrice unité $I_3$,  $\sigma_0$ une constante positive;

- $\zeta$: est une fonction donn\'ee injective d'un ouvert $U\subset \BbR^2 \rightarrow \BbR^3$ définie par:
$$\zeta(X)=\zeta(X_1,X_2)=\begin{pmatrix} {
\frac{1}{2}\left[(X_1-a_1)^2+(X_2-b_1)^2\right] \cr
\frac{1}{2}\left[(X_1-a_2)^2+(X_2-b_2)^2\right] \cr
\frac{1}{2}\left[(X_1-a_3)^2+(X_2-b_3)^2\right] }
\end{pmatrix} $$
On prendra comme composante $L_i$ du vecteur  observation la quantité $L_i=\ds \frac{D^2_{i\,\mbox{observée}}}{2}$

1. Montrer que les vecteurs $\displaystyle \frac{\partial \zeta}{\partial X_1},\frac{\partial \zeta}{\partial X_2}$ sont lin\'eairement ind\'ependants en chaque point $X \in U$.

2. Montrer que les fonctions $\displaystyle \frac{\partial^2 \zeta}{\partial X_i \partial X_j}$ sont continues sur $U$ pour $i,j \in \left\{1,2\right\}$.

3. Posons: $ J=\|L-\zeta(X)\|^2$

Calculer les coefficients de la matrice $(\displaystyle \frac{\partial^2 J}{\partial X_i \partial X_j} ),i,j\in \left\{1,2\right\} $.

4. Soit la matrice carrée d\'efinie par:
$$	g(X)=(g_{ij}) \quad \mbox{avec}\,g_{ij}=<\frac{\partial \zeta(X)}{\partial X_i},\frac{\partial \zeta(X)}{\partial X_j}> \quad \left\{
\begin{array}{ll}
	i=1,2 \\
	j=1,2
\end{array} \right.$$
Calculer les coefficients $g_{ij}$.

5. Introduisons la matrice $B$ d\'efinie par:
$$	B(X,L)=(B_{ij}) \quad \mbox{avec}\,\,\,B_{ij}= g_{ij}-<L-\zeta(X),\frac{\partial^2 \zeta}{\partial X_i\partial X_j}> \quad \left\{
\begin{array}{ll}
	i=1,2 \\
	j=1,2
\end{array} \right. $$
Calculer les éléments de la matrice $B$ et montrer qu'elle est définie positive.
\epb
\chapter{\textit{\textbf{Bibliographie II}}}\label{biblio2}

1. \label{bibb5} \textbf{J.M. Ortega \& W.C. Rheinboldt}. 1970. \textit{Iterative Solution of Nonlinear Equations in Several Variables}. Academic Press. 527p. \index{\textbf{Ortega J.M.}}\index{\textbf{Rheinboldt W.C.}}

2. \textbf{A. Bjerhammer}. 1973. \textit{Theory of Errors and Generalized Matrix Inverses}. Elsevier Scientific Publishing Compagny. Amsterdam. 420p.\index{\textbf{Bjerhammer A.}}

3. \label{bibb3} \textbf{H. Cartan}. 1979. \textit{Cours de Calcul Diff\'erentiel}. Collection Les M\'ethodes. Hermann, Paris. 362p.\index{\textbf{Cartan H.}}

4. \textbf{P. Hottier}. 1980. \textit{Th\'eorie des Erreurs}. Ecole Nationale des Sciences Géographiques. IGN France. 81p.\index{\textbf{Hottier P.}}

5. \label{bibb2} \textbf{A. P\'azman}. 1984. Probability distribution of the multivariate nonlinear least-squares estimates. Kybernetika n°20, pp. 209-230.\index{\textbf{P\'azman A.}}

6. \label{bibb4} \textbf{P.J.G. Teunissen}. 1985. \textit{The Geometry of Geodetic Inverse Linear Mapping and Non-Linear Adjustment}. Publications on Geodesy, n°1, Volume 8, Netherlands Geodetic Commission. 177p. \index{\textbf{Teunissen P.J.G.}}

7. \label{bibb1} \textbf{E.W. Grafarend \& B. Schaffrin}. 1989. The Geometry of non-linear adjustment - the planar trisection problem. \textit{FESTCHRIFT to TORBEN KRARUP} edited by E. Kejlo, K. Poder and C.C. Tscherning. Geod\ae tisk Institut, Meddelelse n°58, pp. 149-172. K\o benhavn, Danmark. \index{\textbf{Grafarend E.W.}}\index{\textbf{Schaffrin B.}} 

8. \label{bibb4a} \textbf{ P.J.G. Teunissen}. 1990. Nonlinear least squares. Manuscripta Geodaetica, Vol 15, n°2, pp. 137-150. 

9. \textbf{S. Amari \& H. Nagaoka}. 2000. \textit{Methods of Information Geometry}. Traduit du japonais par Daishi Harada. Translations of Mathematical Monographs, Vol 191. Oxford University Press. AMS. 206p.

10. \textbf{C. Brezinski}. 2005. La Méthode de Cholesky. Revue d'Histoire des Mathématiques, publication de la Société Mathématique de France. Vol 11 (2005), pp. 205-238.\index{\textbf{Cholesky A.L.}}\index{\textbf{Brezinski C.}}

11. \textbf{P.J. Olver}. 2013. \textit{Nonlinear Systems}. 58p. Université de Minnesota.







\clearpage
\addcontentsline{toc}{chapter}{\textsc{\textbf{Liste des Figures}}}
\addcontentsline{toc}{section}{Liste des Figures}

\listoffigures
\addcontentsline{toc}{section}{Liste des Tables}

\listoftables
\printindex{\textsc{Index}}
\newpage
\clearpage
\section*{\textsc{L'Auteur} }
\begin{figure}
\begin{flushleft}
		\includegraphics[width=0.25\textwidth]{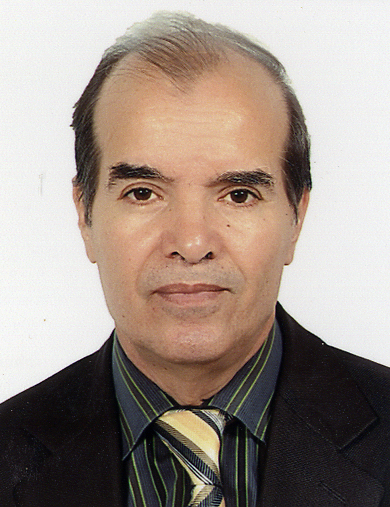}
	\label{fig:abdelmajid}
\end{flushleft}
\end{figure}
L'auteur \textbf{Abdelmajid Ben Hadj Salem }est Ingénieur Géographe Général, retraité de l'Office de la Topographie et du Cadastre (OTC), ancien élève de l'Ecole Nationale des Sciences Géographiques de l'IGN France. Spécialiste en géodésie, il avait participé en 1982 aux travaux de terrain de la revalorisation de la géodésie tunisienne. 
\\

Membre de la commission technique de géodésie de l'OTC, il avait étudié en détail les systèmes et les réseaux géodésiques tunisiens et avait participé à la mise à niveau de la géodésie tunisienne. Il est aussi enseignant, depuis une vingtaine d'années en matière de géodésie à l'université tunisienne et également un formateur. Il a rédigé plus de 100 notes et rapports en la matière collectées en trois tomes de Selected Papers.
\\

A. Ben Hadj Salem était membre de l'Association Internationale de Géodésie (AIG) et il avait participé aux calculs du projet ADOS (African Doppler Survey) initié par l'AIG,  dans le cadre du Consortium Africain de Traitement des Données Doppler. Il  avait représenté l'OTC à l'Organisation Africaine de Cartographie et de Télédétection en qualité d'expert technique. Il était aussi membre du comité de rédaction de \textit{Géo-Top} la publication scientifique de l'OTC.

\end{document}